\documentclass[11pt,twoside,oneandahalfspace,openany]{mitthesis}
\usepackage{lgrind}
\usepackage{myMacros}
\usepackage{subcaption}
\usepackage{amsrefs}
\usepackage[disable]{todonotes}
\usepackage{chngcntr}
\usepackage{setspace}
\usepackage{verbatim}
\renewcommand{\sp}{\mathrm{sp}}
\newcommand{\spT}{{^\sp T \mathcal M}}

\newcommand{\phgd}{\mathcal A_{\mathrm{phg}}^{\mathcal E}}
\newcommand{\phg}[1]{\mathcal A_{\mathrm{phg}}^{#1}}
\newcommand{\phgi}[2]{\mathcal A_{\mathrm{phg},#1}^{#2}}
\newcommand{\Mp}{\widebar{\mathcal M'}}
\DeclareMathOperator{\sob}{sob}
\DeclareMathOperator{\reg}{reg}

\hypersetup{
      colorlinks=true, 
      linkcolor=[RGB]{0, 0, 0}, 
      citecolor=[RGB]{0, 0, 0}, 
      filecolor=magenta, 
      urlcolor=[RGB]{0,0,0}, 
}

\newcommand{\loc}{\mathrm{loc}}
\allowdisplaybreaks

\begin{document}
\pagestyle{plain}

\title{Asymptotic description of the formation of black holes from short-pulse data}

\author{Ethan Yale Jaffe}

\department{Department of Mathematics}

\degree{Doctor of Philosophy}
\degreemonth{May}
\degreeyear{2020}
\thesisdate{March 11, 2020}

\supervisor{Richard B. Melrose}{Professor of Mathematics}

\chairman{Davesh Maulik}{Chairman, Department Committee on Graduate Studies}

\maketitle

\cleardoublepage
\setcounter{savepage}{\thepage}
\begin{abstractpage}
In this thesis we present partial progress towards the dynamic formation of black holes in the four-dimensional Einstein vacuum equations from Christodoulou's short-pulse ansatz. We identify natural scaling in a putative solution metric and use the technique of real blowup to propose a desingularized manifold and an associated rescaled tangent bundle (which we call the ``short-pulse tangent bundle'') on which the putative solution remains regular. We prove the existence of a solution solving the vacuum Einstein equations formally at each boundary face of the blown-up manifold and show that for an open set of restricted short-pulse data, the formal solution exhibits curvature blowup at a hypersurface in one of the boundary hypersurfaces of the desingularized manifold.

This thesis is intended to be partially expository. In particular, this thesis presents an exposition of double-null gauges and the solution of the characteristic initial value problem for the Einstein equations, as well as an exposition of a new perspective of Christodoulou's monumental result on the dynamic formation of trapped surfaces \cite{ChrForm}.

\end{abstractpage}


\cleardoublepage

\section*{Acknowledgements}
\begin{singlespace}
I would like to thank my advisor, Richard Melrose, for his support and guidance, and for sharing with me many ideas that were essential to make progress in this project.

I would also like to thank Semyon Dyatlov, Peter Hintz and Jared Speck for their interest in this project and many useful conversations. I would like to thank David Jerison, Andrew Lawrie, Casey Rodriguez, and Gigliola Staffilani for their support both during this project and throughout my doctoral studies.

I would like to single out Prof. Melrose, Prof. Hintz, and Prof. Jerison for additional thanks for serving on my thesis committee.

I thank Andrew Ahn, Vishal Arul, Thao Do, Ricardo Grande Izquierdo, Vishesh Jain, and Kaavya Valiveti for many useful conversations and late nights listening to me blabber about whatever was on my mind. I thank Andrew and Vishal additionally for making sure I never suffocated under a bar while bench pressing. I thank Kaavya and Thao, as well as Robert Burklund, Sahana Vasudevan, Hong Wang and Jane Wang for making 2-231D the best graduate-student office in the department while we were all still there. I also thank my other friends in the department: Hood Chatham, Robin Elliott, Paul Gallagher, Lucas Mason-Brown, Chris Ryba and Nick Strehlke for always being willing to go the common room or get tea and then talk about anything. I thank the rest of the graduate students in the MIT department of mathematics -- both current and former -- for making these past five years so enjoyable.

I thank Ricardo, as well as Alex van Grootel and Jon Paul Janet for their friendship and livening up a dreary apartment with their high spirits.

I thank Assaf Bar-Natan, Adam Jaffe and Jo Ann Peritz for reminding me that old friends and family were never so far away from Boston, and Joe Kupferman, Janet Freilich, and Jacob Kupferman for their warm hospitality and being my family \emph{in} Boston.

I need to thank Kelly McColough for being a constant source of support and companionship, and moreover putting up with all my shenanigans during the completion of this thesis. I would not have had the peace of mind to focus on this project without you.

Lastly, I thank my parents, Brandon Jaffe and Elaine Peritz, for their unconditional support of all my endeavours.\bottomnote{Much of the research undertaken in the development of this thesis was carried out with support from an NSERC award PGSD3-487557-2016. Some of the thesis-writing was undertaken with support from D. Jerison's NSF grant 1500771.}
\end{singlespace}

\tableofcontents
\listoffigures

\counterwithout{equation}{section}
\counterwithout{thm}{section}
\numberwithin{footnote}{chapter}

\chapter{Introduction}
\section{Background and main results}
In this thesis, we study singularity formation in the four-dimensional Einstein vacuum equations (EVEs). Specifically, we study metrics solving the EVEs with data given by Christodoulou's short-pulse ansatz, first defined in \cite{ChrForm}.

Recall that the EVEs are
\begin{equation}\label{eq:C1:EVE}\Ric(g) = 0,\end{equation}
where $g$ is a metric of Lorentzian signature $(-,+,+,+)$ on a manifold $M$. We refer to the pair $(M,g)$ as a \emph{spacetime}. The fundamental example of a spacetime is the Minkowski spacetime where $M = \R_t \times \R^3_y$ and
\[g = \mathring{g} := -dt^2 + dy^2.\]

Equation \eqref{eq:C1:EVE} is a nonlinear second-order partial differential equation for the metric $g$. One key feature of the equation is diffeomorphism, or \emph{gauge}, invariance: if $\phi: M \to M$ is any diffeomorphism, and $g$ is a Lorentzian metric satisfying $\Ric(g) = 0$, then $\Ric(\phi^\ast g) = 0$, too. Thus, when solving \eqref{eq:C1:EVE}, one must impose some fixed \emph{gauge} which breaks the diffeomorphism invariance. 

Perhaps the simplest example of a gauge is the DeTurck/wave gauge, introduced in \cite{DeTExis}. One fixes a background metric $k$ on $M$ (not necessarily Lorentzian) and requires in addition to $g$ solving \eqref{eq:C1:EVE}, that the identity map $\mathrm{id}:(M,g) \to (M,k)$ is a wave map. One says that ``$g$ is in wave gauge with respect to $k$.'' In the simple case that $M = \R_t \times \R^3_y$ and $k = \mathring{g}$, this just means that the coordinate functions $t$ and $y^i$ ($i = 1,2,3$) solve the linear wave equation, i.e.\
\[\Box_g t = \Box_g y^i = 0,\]
and says that the coordinates are ``harmonic.'' The advantage of this is that one may write
\[\Ric(g) = P(g) + G(g),\] where $G(g)$ vanishes if $g$ is in wave gauge with respect to $k$, and $P(g)$ is a quasilinear hyperbolic operator with leading part $g^{\alpha\beta}\nabla_{\alpha}\nabla_{\beta}g_{\gamma\delta}$, $\nabla$ being the Levi-Civita connection of $k$, and we have adopted abstract index notation and the Einstein summation convention to express tensor contraction. Thus, in a wave gauge, the EVEs manifest themselves as a (nonlinear) wave equation. We refer the reader to Chapter XI of \cite{ChoGene} for a more complete description (although see also \cref{C:C2:gettingsolution}).

Although attractive, a wave gauge will not be sufficient for our purposes. We will use instead what is called a \emph{double-null gauge}. Double-null gauges have their origins in the work of Christodoulou \cite{ChrNonl}, and were subsequently used by Klainerman--Nicol\`o \cite{KlaNicEvel} to obtain a new proof of the stability of Minkowski spacetime. Roughly speaking, one chooses transverse functions $u$, $\und{u}$ defined over $M$, and requires that $u$, $\und{u}$ are optical functions for the spacetime metric $g$, i.e., $g(\grad u,\grad u) = g(\grad \und{u},\grad \und{u}) = 0$ (so the level sets of $u$, $\und{u}$ are null hypersurfaces). An important feature of the transversality of $u$ and $\und{u}$ is that the map $(u,\und{u})$ is a submersion, and we require its fibres to be compact. Thus, one may view the spacetime as fibred over a base, the range of $(u,\und{u})$, with fibres diffeomorphic to some compact manifold. To fix the remaining gauge freedom and complete the gauge, one also needs to specify one null geodesic congruence (either $\gradt_g u$ or $\gradt_g \und{u}$). See \cref{C:C2} for details. The level sets of $u$ and $\und{u}$ (without the specification of a congruence) form what we call a ``double-null foliation.'' 

As one may expect from the hyperbolic character of the Einstein equations in a wave gauge, the EVEs admit a well-posed Cauchy problem. This was first shown by Choquet-Bruhat\footnote{Choquet-Bruhat's work predates the complete development of the wave gauge, so in fact she only used harmonic coordinates.} \cite{ChoTheo} (see also Chapter XI of \cite{ChoGene} for an exposition of the proof), who later with Geroch \cite{ChoGerGlob} showed the existence of a ``maximal solution,'' in the sense of a maximal \emph{globally hyperbolic solution}. Global hyperbolicity is a statement that the causal structure of the spacetime is ``reasonable'' (we refer the reader to Chapter XII of \cite{ChoGene} for a precise definition of globally hyperbolic).

Given the Cauchy problem, it makes sense to ask whether one may prescribe data which are free of singularities, such that the solution of the EVEs with these data develops singularities in the future, and in particular whether a ``black hole'' can form dynamically. In other words, can the Einstein equations cause a black hole to form from regular initial data? The history of solutions to the EVEs with singularities goes back to the discovery of the famous Schwarzschild black hole spacetime \cite{SchUber}. However, it is a much harder problem to determine if a black hole can form dynamically. This question was answered in the affirmative by Christodoulou in the case of the Einstein scalar field system (where the equation for the metric $g$ is coupled to a scalar field, which itself evolves via the wave equation for $g$) under spherical symmetry \cites{ChrForm91,ChrBoun}. Singularity formation (albeit not necessarily ``black hole'' formation) in the higher-dimensional ($1+n$ for $n \geq 4$) EVEs is also known (see for instance \cite{XinXueExam}).

In a monumental work \cite{ChrForm}, Christodoulou took the first step to answering the question for the ($1+3$)-dimensional Einstein \emph{vacuum} equations, without symmetry. In it, he proved that a codimension-two compact \emph{trapped surface} can form dynamically from data he dubbed ``the short pulse ansatz.''

Trapped surfaces were first introduced by Penrose \cite{PenGrav}. To define a trapped surface, recall that at any point of a codimension-two spacelike submanifold $\Sigma$ of a spacetime $(M,g)$ (with a time orientation), there are two independent, future-directed null directions which are orthogonal to $\Sigma$. The surface $\Sigma$ is called \emph{trapped} if the mean curvature in both these directions is negative. By contrast, the spheres $\{t = t_0\}\n \{|y| = r_0\} \subseteq \R_t \times \R^3_y$, $t_0 \in \R$, $r_0 > 0$, in the Minkowski spacetime have mean curvature increasing in one direction, and decreasing in the other.

The importance of trapped surfaces is because Penrose showed: \begin{thm}[Penrose Incompleteness Theorem \cite{PenGrav}]\label{thm:C1:Penrose}A maximal globally hyperbolic spacetime\footnote{Since possessing a Cauchy hypersurface implies global hyperbolicity, this first assumption is, in fact, redundant.} with a non-compact Cauchy hypersurface satisfying the Einstein vacuum equations \eqref{eq:C1:EVE} and containing a compact trapped surface is future causally geodesic incomplete, i.e.\ there is at least one future-directed timelike or null geodesic which does not exist for all values of its affine parameter.\end{thm}
The failure of the existence of this causal geodesic means that, in a suitable sense, the metric $g$ cannot have a global solution. 

Hawking-Ellis \cite{HawEllLarg} extended \cref{thm:C1:Penrose}, showing that if a spacetime containing a closed trapped surface also possesses a complete \emph{future null infinity}, then the spacetime must contain an \emph{event horizon}. This means that the existence of a trapped surface may be taken as evidence that a ``black-hole'' region will form. Since Christodoulou showed that a trapped surface could form dynamically, using \cref{thm:C1:Penrose} he could conclude that a black hole must be present in the future development of the solution he found, although he was unable to discuss the exact nature of the singularity.

The formation of a black hole in the EVEs cannot result from the collapse of matter, but rather must be the result of strong incoming gravitational radiation. Recently, experiments have been able to measure gravitational radiation in the real\footnote{i.e.\ non-mathematical} universe (in this case form the merger of black holes) \cite{AbbObse}. Christodoulou's ``short pulse'' is supposed to describe the requisite gravitational radiation.

The short-pulse data are posed for what is known as the \emph{characteristic initial value problem}. In contrast to the traditional Cauchy problem, for which data are posed on a codimension-one spacelike submanifold, initial data for the characteristic initial value problem are posed on a pair of null hypersurfaces $\mathcal H_1$ and $\mathcal H_2$ intersecting at a compact codimension-two manifold. The precise description the initial data set for the characteristic initial value problem is more complicated than the initial data for the Cauchy problem (see \cref{thm:C2:CIVP}), but roughly speaking allows one to specify the conformal class of the metric on the initial hypersurfaces $\mathcal H_1$ and $\mathcal H_2$, as well as some data at their intersection, and obtain a solution in a small region containing $\mathcal H_1 \n \mathcal H_2$. The well-posedness of the characteristic initial value problem is due to Rendall \cite{RenRedu}, while a good exposition of his original approach can be found in Luk \cite{LukLoca} (although see also our viewpoint in \cref{C:C2}).

To introduce the short-pulse data for a metric $g$, let us first define on $\R_t \times \R^3_y$ the functions
\[u = 1-\frac{|y|-t}{2},\ \und{u} = \frac{|y|+t}{2}.\] The level sets of $u$ and $\und{u}$ are ``outgoing'' and ``incoming'' cones, respectively. The functions $u$ and $\und{u}$ are generically smooth and transverse, and provide the optical functions of a double-null gauge. A fibre of $(u,\und{u})$ is generically a sphere $S^2$, although it degenerates to a point at the tips of cones. The short-pulse ansatz is really a family of initial data, parametrized by $\delta > 0$ small. For the short-pulse ansatz, $\mathcal H_1 = \{\und{u} = 0, 0 \leq u < 1\}$, $\mathcal H_2 = \{u = 1, 0 \leq \und{u} \leq \delta\}$. Thus, $\mathcal H_1 \n \mathcal H_2$ is a sphere $S^2$. Also notice that $\mathcal H_1$ and $\mathcal H_2$ are both fibred by $S^2$. Let $\slash{g}$ denote the (Riemannian) metrics $g$ induces on the fibred spheres. The data along $\mathcal H_1$ are to coincide with the restriction of $\mathring{g}$ (recall that $\mathring{g}$ denotes the Minkowski metric) to $\mathcal H_1$. The data on $\mathcal H_2$ are non-trivial, and correspond to a ``pulse'' of gravitational radiation. One fixes a tracefree, $\mathring{\slash{g}}$-symmetric tensor $\T \in C^{\infty}([0,1]\times S^2;TS^2\otimes T^\ast S^2)$ (recall that $\mathring{\slash{g}}$ is the metric on the round sphere), with $\T(0) = 0$, which we call the \emph{short-pulse tensor}, and specifies the conformal class of data on $\mathcal H_2$ to be given by \[\slash{g} = \mathring{\slash{g}}\cdot \exp\left(\sqrt{\delta}\T(\und{u}/\delta)\right),\] where $\exp(\bullet)$ is the exponential map (which is well-defined from $C^{\infty}([0,1]\times S^2;TS^2\otimes T^\ast S^2)$ to itself via power series), and the $\cdot$ indicates the contraction of a symmetric type $(0,2)$ tensor and type $(1,1)$ tensor to a type $(0,2)$ tensor. The reason that the data are prescribed to be trivial along $\mathcal H_1$ is so that one may extend the solution for $\und{u} \leq 0$ to be the Minkowski spacetime, which Christodoulou requires in order to extend his solution to all of $\{u = 1\}$ (the extension to $\und{u} > \delta$ being arbitrary) for his spacetime to possess a non-compact Cauchy hypersurface and apply \cref{thm:C1:Penrose}.

Not all tensors $\T$ necessarily result in the formation of a trapped surface. Define the energy $\mathbf{E}(v,\theta) = \frac{1}{2}\int_0^v |\pa_s\T|^2(s,\theta)\ ds$. Then Christodoulou requires $\mathbf{E}(1,\theta)$ to be nonzero in every direction in $\theta \in S^2$. Furthermore, the energy cannot be too large, otherwise the initial data will possess a trapped surface (and thus the trapped surface does not form \emph{dynamically}, i.e.\ from data free of trapped surfaces).

\begin{figure}[htbp]
\centering
\includegraphics{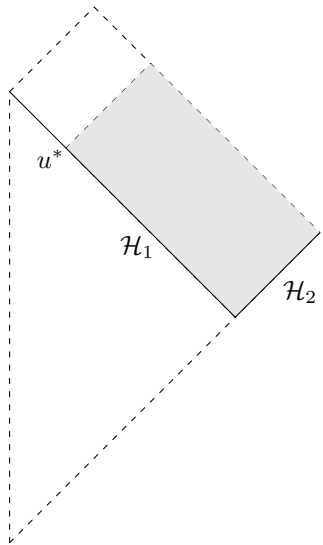}
\caption{A slice at a certain angle of the region of existence of the metric of Christodoulou. The shaded region is the region of existence, while the white triangle enclosed by dotted lines is Minkowski spacetime. $\mathcal H_1$ is a backwards cone in Minkowski spacetime. The upper vertex of the left triangle corresponds to the tip of the cone
$(u,\protect\und{u}) = (1,0)$.
}
\end{figure}
We may vaguely state Christodoulou's theorem as:
\begin{thm}[Christodoulou \cite{ChrForm}]\label{thm:C1:ChrVague}Fix $0 < u^\ast < 1$ and short-pulse data given on $\mathcal H_1 \un \mathcal H_2$. Then if $\delta > 0$ is small enough, there exists a smooth solution $g$ to the EVEs with data given by the short pulse ansatz, which exists in a )double-null gauge for $0 \leq \und{u} \leq \delta$ and $0 \leq u \leq u^\ast$. For an open set of short-pulse tensors $\T$ (those with $\inf_{\theta} \mathbf{E}(1,\theta) > 0$ but $\sup_{\theta} \mathbf{E}(1,\theta)$ not too large), $g$ possesses a trapped sphere $S^2$, and the initial data are free of trapped surfaces. \end{thm}
See \cref{C:C3:ansatz} for a precise statement of the theorem and description of the initial data, and \cite{DafForm} for a more detailed exposition of Christodoulou's proof.
\begin{rk}Christodoulou actually proves something slightly stronger: instead of $1$, he allows himself to pose data on a backwards light cone whose base has radius $R$, and then take $R \to \infty$. This allows him to prove a existence theorem from past null infinity. For this thesis, we will content ourselves with the problem posed in a finite region.\end{rk}

The primary goal of this thesis is to continue Christodoulou's program of black-hole formation in the (1+3)-dimensional EVEs, by taking the first steps to solving for Christodoulou's metric $g$ all the way up until its future boundary, where incompleteness happens, where one hopes to see the formation of a singularity which one can call a ``black hole.'' 

In a double-null gauge, the future domain of dependence of the initial surface is the region \[\mathcal S = \{0 \leq u \leq 1, \ 0 \leq \und{u} \leq \delta, \ \delta \geq 0\},\]
which is singular in a number of places. It is both singular at $\{\und{u} = \delta = 0\}$ and the ``cone points'' $\{\und{u} = 1-u = 0\}$. Using the technique of real blowup,\footnote{See \cref{C:A4:blowup} for a primer on blowups.} we desingularize this space to a manifold with corners (mwc), $\widebar{\mathcal M}$. See \cref{C:C4:blowup} for details. The desingularized space $\widebar{\mathcal M}$ has several faces: two initial faces corresponding to $\mathcal H_1$ and $\mathcal H_2$, $\mathbf{lf}$ and $\mathbf{rf}$, respectively, a bottom face corresponding to $\{\delta = 0\}$, $\mathbf{bf}$, an ``artificial'' far face corresponding to $\und{u} = \delta$, $\mathbf{ff}$, two faces corresponding to the singular point $\und{u} = \delta = 0,\ u = 1$, $\mathbf{if}$ and $\mathbf{sf}$, and a face corresponding to $u = 1, \mathrm{u} = 0$, $\mathbf{cf}$. The faces $\mathbf{sf}$ and $\mathbf{cf}$ should be thought of being in the chronological future of $\mathbf{if}$, at least according to the time orientation induced by the vector field $\pa_u + \pa_{\und{u}}$.

\begin{figure}[htbp]
\centering
\includegraphics{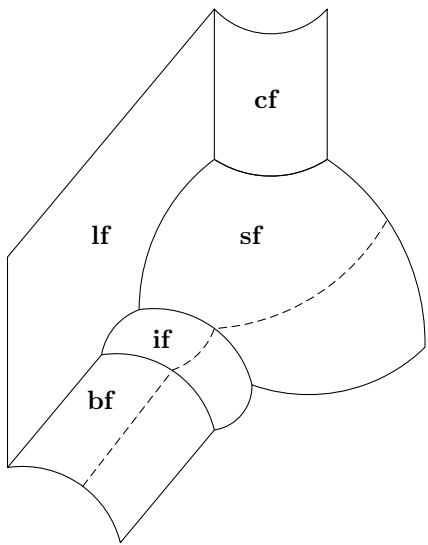}
\caption{A stylized view of $\protect\widebar{\mathcal M}$. The dashed line is where $\mathbf{ff}$ intersects $\mathbf{bf}$, $\mathbf{if}$ and $\mathbf{sf}$. $\mathbf{rf}$ is not drawn, and is a face meeting $\mathbf{lf}$, $\mathbf{bf}$ transversely.}
\label{fig:C1:blowup}
\end{figure}

The trade-off for desingularizing $\mathcal S$ is that the putative solution metric, $g$, is not expected to be a section of $\Sym^2(T^\ast\widebar{\mathcal M})$. However, there is a rescaled bundle, which we call ${}^\sp T\widebar{\mathcal M}$ ($\mathrm{sp}$ for ``short-pulse''), such that we expect $g$ to be a section of $\bm{w}\Sym^2({}^\sp T^\ast\widebar{\mathcal M})$, for some smooth weight function $\bm{w}$ which is a product of various boundary-defining functions (bdfs) of the faces of $\widebar{\mathcal M}$. The notion of a metric on ${}^\sp T\widebar{\mathcal M}$ being in a double-null gauge with $u$, $\und{u}$ extends from the interior $\widebar{\mathcal{M}}^{\circ}$ (which is identified with the regular region $\{0 \leq u < 1, \ 0 \leq \und{u} < 1, \delta > 0\}$) to give a well-defined notion on all of $\widebar{\mathcal M}$. We present the details of this in \cref{C:C4:DNG}.

For given short-pulse data, we introduce a smooth function $\digamma: \mathbf{if}\setminus \mathbf{sf} \to \R$, which depends only on the energy $\mathbf{E}$. We remark that the hypersurface $\{\digamma = 0\}$ extends to a smooth hypersurface of $\mathbf{if}$ (although it only intersects $\mathbf{lf}$ and $\mathbf{sf}$ at their common intersection, i.e.\ at a corner), and the region $\{\digamma > 0\}$ contains $\mathbf{bf}\n \mathbf{if}$ and $(\mathbf{lf}\setminus \mathbf{sf})\n \mathbf{if}$. The main goal of this thesis is to present evidence that $\widebar{\mathcal M}$ is the ``correct'' space on which to look for a continuation of the solution of Christodoulou. We will show the existence of a solution $g$ ``in asymptotic series,'' i.e.\ a polyhomogeneous\footnote{See \cref{C:A4:phg} for a primer on polyhomogeneity.} solution which solves $\Ric(g) = 0$ up to a rapidly vanishing error. More precisely, we show: 
\begin{thm}\label{thm:C1:bigthm}
For any short-pulse data, there exists an open subset $W \subseteq \widebar{\mathcal M}$ containing \hspace{0.25pt} $\mathbf{lf} \n \{u < 1\}$, $\mathbf{rf} \n \{\und{u} \leq \delta\}$, $\mathbf{bf}$ and $\mathbf{if} \n \{\digamma > 0\}$, and a family of Lorenztian metrics $g$, which together constitute a polyhomogeneous section of $\bm{w}\Sym^2(^\sp T^\ast \widebar{\mathcal M})$ over $W$, which is in a double-null gauge and induces the given short-pulse data, such that $\Ric(g)$ is rapidly vanishing at all boundary faces of $\widebar{\mathcal M}$ contained in $W$. Such a metric is unique up to a section rapidly vanishing on the same boundary faces. Moreover, the restriction of $g$ to the $S^2$ fibres over $\mathbf{if}$ degenerates as $\digamma \to 0$.

Furthermore, for an open subset of restricted short-pulse data (data which we call ''commutative'' and have a ``very'' short pulse), any such metric cannot be continued beyond $\{\digamma = 0\}$.
\end{thm}
\begin{rk}Since the ultimate goal is to provide some form of singularity development, it is important that the solution in asymptotic series cannot be continued beyond $\{\digamma = 0\}$.\end{rk}
Since the solution of Christodoulou solves $\Ric(g) = 0$ exactly on $W \n\{ u < u^\ast \}$ for any $u^\ast < 1$ the existence portion of \cref{thm:C1:bigthm} may be understood as a continuation theorem in formal series for the solution $g$. 

Notice, however, the \cref{thm:C1:bigthm} is unable to say whether, in any open set of the boundary of $\widebar{\mathcal M}$, there is a family of Lorentzian metrics $g$ satisfying $\Ric(g) = 0$ exactly. The theorem only says that given such a family of metrics, which is suitably smooth and a section of $\bm{w}\Sym^2({}^\sp T^\ast \widebar{\mathcal M})$, then its asymptotic development towards the boundary of $\widebar{\mathcal M}$ is determined. In particular, we do not prove a true existence theorem for solutions of the EVEs beyond what Christodoulou is able to prove. While we suspect it is possible to push some of the techniques developed in our reproof of \cref{thm:C1:ChrVague} (see below) to obtain a solution on a larger region than Christodoulou's (see \cref{thm:C1:conjj} and \cref{thm:C1:conjjj}), the existence of a solution right up until blow up happens (i.e.\ until the formation of the black hole) is beyond the scope of this thesis. Additionally, we are only able to prove a non-continuation theorem for a small class of data. What we hope is true is:
\begin{conj}\label{thm:C1:biggconj}Let $W$ be as in the statement of \cref{thm:C1:bigthm}. Then there is an extension of $\digamma$ to $W$ and a polyhomogeneous section $g$ of $\bm{w}\Sym^2(^\sp T^\ast \widebar{\mathcal M})$ over $W$, which is in a double-null gauge and induces the given short-pulse data, such that $\Ric(g) \equiv 0$. Moreover, the curvature of such a metric blows up at $\{\digamma = 0\}$.\end{conj} We discuss some of the difficulties in proving \cref{thm:C1:biggconj} in the proof outline, below.

\section{Aspects of the proof}
This thesis begins by reexamining the proof of \cref{thm:C1:ChrVague}, both to begin analyzing the solution metric, but also as a proving ground for the techniques of real blowup and the use of the short-pulse tangent bundle. The main advantage of our perspective on the proof is that in contrast to Christodoulou, who needs to prove a ``large data'' theorem, our main analytic work will be in proving a much easier ``small data'' theorem.

Let us begin by considering the space $\mathcal S$ on which we seek to prove the existence of solutions to $\Ric(g) = 0$ given by \[\mathcal S = \{0 \leq u < u^\ast, \ 0 \leq \und{u} \leq \delta, \ \delta \geq 0\}.\] Observe that $\mathcal S$ is singular at $\delta = 0$. We may desingularize it by performing a version parabolic blowup, introducing $v = \und{u}/\delta$, and $x = \sqrt{\delta}$. Let $\mathcal M$ denote the blown-up space, which is a mwc. The map $(u,v)$ remains a submersion onto its range with fibres diffeomorphic to $S^2$.

Although the space is now smooth, we expect the metric itself to be singular. For example, the Minkowski metric takes the form
\[x^{-2}\mathring{g} = -2(du\otimes dv + dv\otimes du) + (1-u+x^2v)^2(x^{-2}\slash{\mathring{g}}),\]
$\mathring{\slash{g}}$ denoting the metric on the round sphere. The singularities of $\mathring{g}$ (or rather $x^{-2}\mathring{g}$) can be described very precisely: rather than be a smooth combination of $du,dv$ and 1-forms $d\theta^i$ on the sphere, it is a smooth combination of $du,dv$ and \emph{singular} 1-forms $d\theta^i/x$. This suggests introducing a rescaled cotangent bundle $\mathcal M$ whose smooth sections are smooth linear combinations of $du,dv,d\theta^i/x$ (but not $dx$). We will call this bundle the short-pulse (co)tangent bundle $^\sp T^\ast\mathcal M$, and nondegenerate sections of $\Sym^2(^\sp T^\ast \mathcal M)$ \emph{short-pulse metrics}. Of course, for $x = x_0 > 0$, the space $\mathcal M \n \{x = x_0\}$ is canonically identified with \[M_{x_0^2} := \{0 \leq u < u^\ast, \ 0 \leq \und{u} \leq \delta, \ \delta = x_0^2\},\]
and the pullback of $^\sp T\mathcal M$ to $M_{x_0^2}$ is the same as the ordinary tangent bundle $TM_{x_0^2}$.

Now, notice that the short-pulse data prescribed along $\mathcal H_2$ become non-singular after blowing up (at least on $^\sp T\mathcal M)$, namely
\[x^{-2}\slash{g} = (x^{-2}\mathring{\slash{g}})\exp(x\T(v)),\] so short-pulse data should fit into this geometric framework.

The point of introducing $^\sp T\mathcal M$ is that for short-pulse metrics there is a well-posed characteristic initial value problem for the EVEs, which we describe in \cref{C:C3:SCIVP}. To prove well-posedness, we proceed in two steps. After this, the formation of trapped surfaces becomes an easy, albeit lengthy, computation.

The EVEs in a double-null gauge are well-adapted to computations in formal series, so we begin by establishing the existence of a solution satisfying the EVEs only modulo a rapidly vanishing error (both at the initial hypersurfaces $\mathcal H_1 \un \mathcal H_2$ and at the new face $\{x = 0\}$). Indeed, over $\mathcal H_1 \un\mathcal H_2$, the EVEs in a double-null gauge become transport equations for components of the metric and all their derivatives (including the derivatives transverse to the initial surfaces). Over $\{x = 0\}$, as evidenced by the appearance of $x\pa_{\theta^i}$ ($i = 1,2$),\footnote{The ``dual'' vector fields to the basic short-pulse one-forms $d\theta^i/x$.} the angular derivatives drop out, leaving only a $(1+1)$-dimensional wave equation to solve. The miracle of the short-pulse ansatz is that the only non-linear wave equation left in the limit is \emph{explicitly solvable}, and all equations for higher-order terms in the formal series are of course linear. 

After solving for a solution modulo a rapidly vanishing error, we may try to perturb off of it, correcting it by a perturbation rapidly vanishing at $\mathcal H_1 \un \mathcal H_2$ and $\{x = 0\}$. Using techniques introduced by Rendall in the solution of the characteristic initial-value problem for non-linear wave equation \cite{RenRedu}, we reduce the characteristic initial value problem to the Cauchy problem, which is easier to handle. Thus, the main analytic portion in our presentation of the proof of \cref{thm:C1:ChrVague} is to prove a long-time existence theorem for small data solutions to a non-linear wave equation (albeit for short-pulse metrics). While the non-linearity interacts unfavourably with the singularity in $x$ of short-pulse metric, because we are only seeking to solve in a space of rapidly vanishing tensors, we may safely absorb any loss in powers of $x$. To prove the small data theorem, we are thus able to leverage classical results in the well-posedness theory of high-regularity solutions to non-linear wave equations, albeit carefully transferred over the short-pulse setting.

We undertake the reexamination of the proof in \cref{C:C3}. Our ``formal-series first'' perspective on the characteristic initial value problem appears new, even in the case of ordinary (i.e.\ \emph{not} short-pulse) metrics. We therefore present an exposition of it in \cref{C:C2}. One advantage of our approach is that unlike the proof of Rendall, who must use coordinate charts to cover the fibres of what we call $(u,\und{u})$, we are able to work \emph{globally} on the fibres.

As a side effect of our proof, we will be able to slightly extend \cref{thm:C1:ChrVague}:
\begin{thm}\label{thm:C1:Smoothness}Let $g_\delta$ denote the metric on $M_\delta$ provided by \cref{thm:C1:ChrVague} for $\delta$ small. After pulling back to $\mathcal M$, the family $g_\delta = g_{x^2}$ is smooth in $x$, all the way down to $x=0$. In other words, the family $g_\delta$ is smooth in $\sqrt{\delta}$ and $\und{u}/\delta$.\end{thm}

After reproving the theorem of Christodoulou, we extend the short-pulse bundle to $\widebar{\mathcal M}$ start to prove the existence/uniqueness portion of \cref{thm:C1:bigthm}, and the degeneracy of the fibre-metric. See \cref{C:C4:existence} and specifically \cref{thm:C4:continuation} for a precise statement. 


The degeneracy of the fibre metric does not imply that the full metric itself becomes degenerate or singular as $\digamma \to 0$, since the degeneracy may be a gauge artefact. For instance, in the coordinates $(u,\und{u},\theta)$ on $(\R_t \times \R^3_y)\setminus \{y = 0\}$, the Minkowski metric takes the form
\[\mathring{g} = -2(du\otimes d\und{u} + d\und{u}\otimes du) + (1-u+\und{u})^2\mathring{\slash{g}}.\]
The fibre metric, $(1-u+\und{u})^2\mathring{\slash{g}}$, degenerates as $u \to 1$ and $\und{u} \to 0$. However, the total metric of course remains smooth, and the apparent degeneracy is only an artefact of the coordinate system.

To show that the metric does not continue, we use the Kretschmann scalar, which is a scalar curvature invariant defined by
\[\mathcal K = |\Riem(g)|_g^2 = R^{\alpha\beta\gamma\delta}R_{\alpha\beta\gamma\delta}.\] Its blow up is manifestly coordinate invariant.

We expect the $\mathcal K$ associated to a solution $g$ to be polyhomogeneous at $\mathbf{if}$, so we will need to rescale it before its restriction to $\mathbf{if}$ makes sense. Let $\tilde{\mathcal K}$ denote this rescaled version. In fact we may take $\tilde{\mathcal K} = \xi^{12}\mathcal K$, where $\xi$ is a boundary-defining function of $\mathbf{if}$. Then $\tilde{\mathcal K}$ makes sense on $(\mathbf{if} \n \{\digamma > 0\})\setminus \mathbf{bf}$ (since $\{\digamma = 0\} \n \mathbf{bf} = \emptyset$, this rescaling is sufficient to see blow up, although we could further rescale by a boundary defining of $\mathbf{bf}$ to make the restriction well-defined on all of $\mathbf{if}$). The blow up of the scalar $\tilde{\mathcal K}$ at $\{\digamma = 0\}$ will show that the metric cannot continue past $\{\digamma = 0\}$.

We conjecture:
\begin{conj}\label{thm:C1:bigconj}For arbitrary (non-trivial) short-pulse data, the rescaled Kretschmann scalar $\tilde{\mathcal K}$, when restricted to $\mathbf{if} \n \{\digamma > 0\}$ is unbounded in a neighbourhood of least one point of $\{\digamma = 0\}$.\end{conj}

Proving this conjecture is unfortunately beyond the scope of this thesis. However, we can prove blow up in the case of ``commutative'' data which is data for which the short-pulse tensor $\T$ factors into $\psi(v)\T_0$, where $\psi$ is smooth in $v$ and doesn't depend on $\theta$ and $\T_0$ does not depend on $v$. The reason for calling this data commutative will become clear in \cref{C:C5}.

\begin{thm}\label{thm:C1:blowupvague}There exists an open set of commutative short-pulse data such that the rescaled Kretschmann scalar $\tilde{\mathcal K}$, when restricted to $\mathbf{if}$, blows up in a neighbourhood of an open subset of $\{\digamma = 0\}$. More precisely, $\tilde{\mathcal K} \gtrsim \digamma^{-3}$ on this subset.\end{thm}
We state a more precise version of this theorem in \cref{C:C6:prelims}, specifically \cref{thm:C6:blowup}. The open subset of $\{\digamma = 0\}$ should be thought of as ``large,'' which will be quantified in the precise statement.

In order to prove this, we will first need to analyze the behaviour of the components of the metric at $\{\digamma = 0\}$. We undertake this analysis in \cref{C:C5}. In fact, we give a precise description of what happens for ``generic'' (in the sense of open and dense) choices of commutative data.

We suspect that it is possible to obtain an actual solution just before the singularity occurs $\{\digamma > \epsilon\}$ for any $\epsilon > 0$ using a modification of the analytic techniques introduced in \cref{C:C3}. More specifically, it should be possible to convert to the Cauchy problem and use the idea of Rendall in conjunction with the solution in asymptotic series given by \cref{thm:C1:bigthm} to prove:
\begin{conj}\label{thm:C1:conjj}Fix $\epsilon > 0$. Then there exists an open subset $V \subseteq \widebar{\mathcal M}$ containing $\mathbf{lf} \n \{u < 1\}$, $\mathbf{rf} \n \{\und{u} \leq \delta\}$, $\mathbf{bf}$ and $\mathbf{if} \n \{\digamma \geq \epsilon\}$, on which there is a unique section of $\bm{w}\Sym^2({}^\sp T^\ast \widebar{\mathcal M})$ which is in a double-null gauge and induces the given short-pulse data, such that $\Ric(g) \equiv 0$.\end{conj}

\begin{rk}While we conjecture this holds for all $\epsilon > 0$, the neighbourhood $V$ may shrink with $\epsilon$, so a priori as $\epsilon \to 0$ we recover no solution in an open subset:\ only the top-order at the boundary of $\widetilde{\mathcal M}$, which we already know exists by \cref{thm:C1:bigthm}.\end{rk}

If we assume the data are (generic) commutative, our analysis in \cref{C:C5} should be able to be furthered to obtain more precise information of higher-order terms of the asymptotic expansion at $\mathbf{if}$ near the singular point $\{\digamma = 0\}\n \{v = 0\} \subseteq \mathbf{if}$. This leads to:
\begin{conj}\label{thm:C1:conjjj}Fix $\epsilon > 0$ and an extension of $\digamma$ to all of $\widebar{\mathcal M}$. Then on the set $\{\digamma \geq v\epsilon\}$, there is a unique section of $\bm{w}\Sym^2({}^\sp T^\ast \widebar{\mathcal M})$ which is in a double-null gauge and induces the given (generic commutative) short-pulse data, such that $\Ric(g) \equiv 0$.\end{conj}

Let us end the discussion of the proof by examining some of the difficulties in attempting to prove \cref{thm:C1:biggconj}. One issue is the fibre anisotropy of function $\digamma$: because $\digamma$ is not constant on the fibres, it is not true that $\mathbf{if} \n \{\digamma = 0\}$ is the base of an $S^2$-fibre bundle, and the notion of a double-null gauge breaks down. 

A more serious problem, however, is that even for the ``commutative'' data for which \cref{C:C5} provides precise asymptotics for restriction to $\mathbf{if}$ of the metric near $\{\digamma = 0\}$, the linearized EVEs involved in computing the higher-order terms become very difficult to work with. In particular, we are not able to even guess the behaviour near $\{\digamma = 0\}$ of the higher-order terms at $\mathbf{if}$ in the asymptotic series for the metric.

To see why, we may deduce from \cref{C:C5} that some metric components behave roughly like $\digamma^A = \exp(A\log \digamma)$, where $A$ is a smooth $2\times 2$ tracefree-matrix-valued function. The assumption of ``commutative data'' implies that the values of $A$ are simultaneously diagonalizable, so essentially $\digamma^A = \left(\begin{smallmatrix} \digamma^{\lambda} & 0\\0 & \digamma^{-\lambda}\end{smallmatrix}\right)$, for some real-valued function $\lambda$. This means that it becomes difficult to even evaluate the error terms in the iterative process, let alone solve the linear equations which should remove the error. For instance, one quantity (among dozens) which appears is
\[\Ric(\mathring{\slash{g}}\digamma^A).\] Without a good a priori guess of the form the solution should take, the number of terms explodes.

However, for data that are not commutative, we do not even have available results which are as strong as those in \cref{C:C5} (the strongest result available in this case is \cref{thm:C4:nonblowup}, which guarantees very weak control on the blowup of certain metric components). We conjecture that certain metric components still behave in some sense like $\digamma^A$. However, now the values of $A$ should no longer be simultaneously diagonalizable. Indeed, we conjecture that the eigenvectors are ``rotating rapidly.'' From a technical perspective, this means that differentiating $\digamma^A$ is problematic, as we do not expect the derivative of $A$ to commute with $A$. We are unable to prove any precise results in this case.

\section{Related work}
The initial work of Christodoulou \cite{ChrForm} has been extended in several different ways. Using the existence theorem as a black box, Klainerman--Luk--Rodnianski \cite{KlaLukRodFull} undertook a more careful analysis of the structure of the spacetime and found a trapped sphere for $g$ provided the anisotropic condition $\sup_{\theta \in S^2} \mathbf{E}(1,\theta) > 0$. Separately, Klainerman--Rodnianski \cites{KlaRodForm, KlaRodEmer} greatly simplified the proof of \cref{thm:C1:ChrVague} by keeping better track of scaling, in particular by noting a certain ``parabolic scaling'' analogous to our use of parabolic blowup. An--Luk \cite{AnLukTrap} extend these ideas and prove an existence theorem for a broader class of short-pulse data (which we do not consider in this thesis) which lie in a ``scaling-critical'' space. An points out in \cite{AnEmer} that the existence theorem may be iterated and obtains, for a variant of short-pulse data which are not smooth (even in $\und{u}$) existence in the region $u \leq 1-\und{u}/\delta$, and extends the ideas in \cite{KlaLukRodFull} to describe the boundary of the trapped region. An has also been able to use the simplified framework of Klainerman--Rodnianski to give a new proof of Christodoulou's result on the formation of trapped surfaces from past null infinity \cite{AnForm}, and has recently \cite{AnScale} been able to further simplify the proof of \cref{thm:C1:ChrVague}. In a separate direction, Li--Yu \cite{LiYuCons} managed to find Cauchy (rather than characteristic) data exhibiting trapped surface formation. Yu has also managed to extend Christodoulou's result to the Einstein-Maxwell equations where the Einstein equations are coupled to an electromagnetic field \cite{YuDyn}.

The work of An \cite{AnEmer} is the first to give an indication of the black hole region. Despite this, there has been so far no adequate description of the maximal domain of existence of the spacetime, nor of the future boundary where the incompleteness predicted by \cref{thm:C1:Penrose} occurs.

Blowups and rescaled (co)tangent bundles have a long history, especially in the elliptic setting, for instance the b-tangent bundle $^{\mathrm{b}}TM$ in Melrose's b-calculus \cite{MelTapsit}, the zero tangent bundle $^{0}TM$ in the zero calculus of Mazzeo--Melrose \cite{MazMelMero}, the more general edge bundle $^{\mathrm{e}}TM$ in Mazzeo's edge calculus \cite{MazElli}, or the scattering tangent bundle $^\mathrm{sc}TM$ in Melrose's scattering calculus \cite{MelSpec}. Perhaps most similar to our bundle are the ``adiabatic'' bundles considered by Mazzeo--Melrose in \cite{MazMelAdia}. In terms of applications to physical theories about gravity, Zhu \cite{ZhuElev} used the edge calculus to study the eleven-dimensional supergravity equations (although this was still done in the elliptic setting).

Blowups and rescaled tangent bundles have more recently also found use in the hyperbolic setting, and indeed in the study of the Einstein equations, although not in the study of singularity formation. We mention the articles of Baskin--Vasy--Wunsch \cite{BasVasWunAsym} in which it was shown that blowing up the future null infinity of spacetimes modelled on $(\R^4,\mathring{g})$ allows for an asymptotic expansion of the forward fundamental solution to the wave equation, and the new proof of Hintz--Vasy of the nonlinear stability of Minkowski space \cite{HinVasGlobA}, where a more complicated blowup is used. We also mention the articles by Hintz--Vasy and Hintz \cites{HinVasGlobN,HinNonL}, where the b-calculus (and in particular the scaling afforded by the b-tangent bundle) was used with great effect to establish the nonlinear stability of the Kerr-de Sitter family of rotating black holes and the Kerr-Newman-de Sitter family of charged, rotating black holes, respectively.

\section{Notation}
If $\phi$ is a tensor, a section of $(TM)^{\otimes p}\otimes (T^\ast M)^{\otimes q}$, for some manifold $M$, then we say that $\phi$ is a type $(p,q)$ tensor. We often use the shorthand $T^p_q M := (TM)^{\otimes p}\otimes (T^\ast M)^{\otimes q}$. We will prefer as much as possible to use coordinate-free notation to describe our tensors, but occasionally will resort to using abstract index notation, especially to denote tensor contraction. We will always keep track of covariant (i.e.\ covector) and contravariant (i.e.\ vector) indices. In abstract index notation, covariant indices will be written ``downstairs,'' and contravariant indices will be written ``upstairs.'' For instance the symbol $\phi^{a}_{b}$ will refer to a type $(1,1)$ tensor, an element (or section) of $T^\ast M \otimes TM$, and not a particular component, and likewise $\phi^{ab}$, $\psi_{ab}$ will denote type $(2,0)$ and type $(0,2)$ tensors, respectively. We will employ the Einstein summation convention for repeated summation, especially in conjunction with abstract index notation where it expresses contraction. For instance, $\phi^a_a$ will denote the trace, $\Tr(\phi)$. If $\phi$ and $\psi$ are type $(1,1)$ tensors, then they are naturally linear maps $TM \to TM$, and the notation $\phi\psi = \phi^a_c\psi^c_b$
traditionally and unambiguously corresponds to multiplication of linear maps.\footnote{\label{fn:C1:wut}Of course one could also interpret $\phi$ and $\psi$ as linear maps $T^\ast M \to T^\ast M$, in which case $\phi\psi = \phi_b^c\psi_c^a$ corresponds to the \emph{reverse} contraction as the one described previously. Let us interpret this in coordinate-free language. Let $F:T^{1}_1M \to \Hom(TM,TM)$ and $G:T^{1}_1 M \to \Hom(T^\ast M, T^\ast M)$ be the natural isomorphisms. Then it is easy to check that $F(\phi)^\ast = G(\phi)$, where $\bullet^\ast$ denotes the natural adjoint of a map $TM \to TM$ defined by $\langle \beta, Tv\rangle = \langle T^\ast \beta, v\rangle$ for all $v \in TM$ and $\beta \in T^\ast M$ (over the same point), and $\langle \bullet,\bullet \rangle$ denotes the natural pairing $T^\ast M \otimes TM \to M\times\R$. Thus 
\[G(\phi)G(\psi) = F(\phi)^\ast F(\psi)^\ast = (F(\psi)F(\phi))^\ast = G(F^{-1}(F(\psi)F(\phi))),\]
i.e.\ multiplication of $\phi$ an $\psi$ interpreted one way is the reverse of multiplication interpreted the other way.

Since our definition is purely notational, it does not matter that the symbol $\phi\psi$ is not ``canonically'' defined, so long as the choice is consistent.} We extend this by an abuse of notation if $\phi$ and $\psi$ are type $(i,j)$ tensors for $i+j = 2$, where the notation $\phi\psi$, indicates contracting the ``inner'' indices of $\phi$ and $\psi$, as would be done if the tensors were written as matrices without regards to upper and lower indices. For instance, if $\phi$ is a type $(2,0)$ tensor and $\psi$ is a type $(1,1)$ tensor, then $\phi\psi = \phi_{ac}\psi^{c}_b$, while $\psi\phi = \psi^{a}_c\phi_{cb}$, and if $\phi$ is a type $(2,0)$ tensor and $\psi$ is a type $(0,2)$ tensor then $\phi\psi = \phi_{ac}\psi^{cb}$. If $\phi$ is a metric, with associated Levi-Civita connection $\nabla$, and $\psi$ is a type $(p,q)$ tensor, then we will denote by $\phi^{-1}\nabla^2 \psi$ the contraction of the two upper indices of $\phi^{-1}$ with the two lower indices of $\nabla^2 \psi$ appearing from taking the covariant derivative, i.e.\ \[\phi^{-1}\nabla^2 \psi = \phi^{ab}\nabla^2_{ab}\psi_{\alpha}^{\beta} = \phi^{ab}(\nabla^2 \psi)_{ab\alpha}^{\beta},\] where $\alpha$ and $\beta$ are (disjoint) strings of abstract indices of length $q$ and $p$, respectively.
 
If $\phi$, $\psi$ are both symmetric type $(0,2)$ tensors, then we may form their Kulkarni-Nomizu product
\[(\phi \owedge \psi)_{abcd} := \phi_{ac}\psi_{bd} + \phi_{bd}\psi_{ac}-\phi_{ad}\psi_{bc} - \phi_{bc}\psi_{ad}.\]
The tensor $\phi \owedge \psi$ has the same algebraic symmetries as the curvature tensor of a metric.

If $g$ is a metric, i.e.\ a nondegenerate section of $\Sym^2(T^\ast M)$, then $g^{-1}$ will denote the dual cometric on $\Sym^2(T M)$, and vice-versa if $g$ is a cometric. This is due to the fact that when expressed in abstract index notation $(g^{-1})^{ab}g_{bc} = \delta^a_c$ is the identity $(1,1)$ tensor. When using abstract index notation, we will just denote $(g^{-1})^{ab} = g^{ab}$. We will use the metric to raise and lower indices. For example, if $\phi^{ab}$ is a $(2,0)$ tensor, then $\phi^a_{\ b} = g_{bc}\phi^{ca}$ denotes the lowered $(1,1)$ tensor. In the case of $g$-symmetric tensors, we omit the space and write $\phi^a_b$, since which index is lowered does not matter. As mentioned above, if $\phi^a_b$ is a $(1,1)$ tensor, then $\Tr(\phi) = \phi^a_a$ is well-defined. In the presence of a metric, one may also define the traces of symmetric $(0,2)$ and $(2,0)$ tensors via contraction with the metric. We will denote these with $\tr_g$ or just $\tr$ if the metric is understood from context. For instance, if $\phi_{ab}$ is a symmetric $(0,2)$ tensor, then $\tr \phi = \Tr(g^{-1}\phi) = g^{ab}\phi_{ba}$.

We will use the convention for curvature:
\begin{align*}
  \Riem(g)(X,Y,Z,T) &= R(X,Y,Z,T) = g(R(X,Y,Z),T)\\
  &= g(\nabla_Y\nabla_XZ-\nabla_X\nabla_YZ + \nabla_{[X,Y]}Z,T).
\end{align*}

Will also use notation, common in the general relativity literature, which we call ``slash notation.'' Very often, we will augment a spacetime $(M,g)$ with a submersion $\pi: M \to N$, so that the fibres of $\pi$ are Riemannian submanifolds (for this thesis, these fibres will mainly have codimension two). The restriction of $g$ to the fibre-tangent directions is thus a fibre-metric of Riemannian signature. 
Thus there will often be two different notions of metrics to consider: a metric on the entire manifold, $M$, and a fibre metric. To keep the two notions distinct, we will always use slashes, such as $\slash{g}$, $\slash{h}$, etc., to denote fibre metrics. If there exists a pair of metrics with the same symbol, one with a slash, and one without, for instance $g$ and $\slash{g}$, then $\slash{g}$ will denote the restriction of $g$ to the fibre-tangent directions.
Associated to metrics are operators and quantities like the Levi-Civita connection, the extrinsic trace, curvature, etc. To keep these distinct, the operations with respect to the fibre metrics will also be denoted with slashes. For instance, if $(M,g)$ is a spacetime as above and $\slash{g}$ the associated fibre metric, then $\nabla$, $\slash{\nabla}$ could denote their respective associated Levi-Civita connections, while $\tr$ and $\slash{\tr}$ could denote their respective trace operators.

The metrics $\mathring{g}$ and $\mathring{\slash{g}}$ will always denote the Minkowski metric and the metric on the round sphere $S^2$, respectively. This convention is a minor abuse of the previous convention. If we take $(u,\und{u},\theta)$ coordinates on $(\R_t \times \R^3_y)\setminus \{y = 0\}$, then $(u,\und{u})$ becomes the projection of an $S^2$-bundle, for which the fibre metric associated to $\mathring{g}$ is $(1-u+\und{u})^2\mathring{\slash{g}}$, rather than just $\mathring{\slash{g}}$ itself.

If $X$ is a mwc, we denote by $X^\circ$ its interior, i.e.\ $X\setminus \bigcup_{F \in \mathcal F} F$, where $\mathcal F$ denotes the collection of boundary faces of $X$.

For our notation for blowups and polyhomogeneous functions, see \cref{C:A4}. If $x$ is a bdf of a face $F$ of a mwc $X$, and $t \in \R$, then the notation $x^{t^{\pm}}$ is shorthand for all possible functions $x^{t \pm \epsilon}$, where $\epsilon >0 $ is sufficiently small. For example, the space $x^{t^{-}}L^{\infty}(X)$ consists of all functions $f$ on a space $X$ for which $f \in x^{t-\epsilon}L^{\infty}(X)$ for all $\epsilon > 0$ (sufficiently small).

We use the notation $A \lesssim B$ to indicate there is a constant $C$, not depending on $A$, $B$, or other implicit quantities, such that $A \leq CB$. If not clear from context, we will mention what these implicit quantities are.

\section{Outline of the thesis}
We give an outline of the thesis.

In \cref{C:C2} we rigorously define double-null gauges and provide a proof of local well-posedness of the characteristic initial value problem for the EVEs in the spirit of Rendall.

In \cref{C:C3} we provide a proof of Christodoulou's theorem, \cref{thm:C1:ChrVague}.

In \cref{C:C4} we construct the blown-up manifold $\overline{\mathcal M'}$ and prove the existence/uniqueness and degeneracy statement of \cref{thm:C1:bigthm}.

In \cref{C:C5} we restrict ourselves to commutative data and analyze the behaviour of the metric components at the boundary $\{\digamma = 0\}$.

In \cref{C:C6} we state and prove a precise version of \cref{thm:C1:blowupvague}.

\numberwithin{equation}{section}
\numberwithin{thm}{section}
\chapter{Double-null gauges and the characteristic initial value problem for the Einstein vacuum equations}
\label{C:C2}

\section{Double-null gauges}
\label{C:C2:DNG}\begin{defn}\label{def:C2:doublefoliated}Let $\mathcal R \subseteq \R^2$ be a codimension-zero embedded manifold with corners. We consider $\R^2$ as $\R_u \times \R_v$, and equip $\mathcal R$ with the coordinates $(u,v)$. Let $M$ be a smooth manifold with corners of dimension at least $3$ and $\pi: M \to \mathcal R$ a proper smooth surjective submersion with fibres diffeomorphic to the same closed (and connected) manifold $S$. Using the coordinates $(u,v)$ on $\mathcal R$, identify the components of $\pi = (u,v)$.\footnote{We will conflate points $(u,v) \in \mathcal R$ with the coordinate functions $u,v$ on $\mathcal R$, and also with the components $u,v$ of $\pi$.} Make the assumptions that, separately, the fibres of $u$ and $v$ are themselves connected. A \emph{doubly-foliated manifold} consists of the data $(M,\mathcal R, \pi = (u,v),S)$ specified above.\end{defn}
We will write $(M,u,v)$ as shorthand for the entire data of a doubly-foliated manifold. If $R \subseteq \mathcal R$ is a subset, we will use the notation $M(R) = \{p \in M\: (u,v)(p) \in R\}$. We will also write $S_{u,v} = \pi^{-1}(u,v)$ for the fibre of $(u,v) \in \mathcal R$ to emphasize that it is a diffeomorphic copy of $S$.

\begin{exam}The most important example for us of a doubly-foliated manifold consists of the subset $M \subseteq \R^4 = \R_t\times \R^3_y$ between the cones defined by
\[0 \leq 1- \frac{|y|-t}{2} < 1, \ 0 \leq \frac{t+|y|}{2} < \infty,\] with $u := 1- \frac{|y|-t}{2}, v = \und{u} := \frac{t+|y|}{2}$, and $\mathcal R = [0,1)\times [0,\infty)$.
\end{exam}

The connectedness of the fibres of $u$ and $v$ together with the compactness of $S$ imply:
\begin{lem}\label{thm:C2:complete}Let $X, \ Y \in C^{\infty}(M;TM)$ be lifts to $M$ of the coordinate vector fields $\pa_u$ and $\pa_v$ on $\mathcal R$. Then it is always possible to flow along $X$ from a point $p \in S_{u_1,v_0}$ to a point $q \in S_{u_2,v_0}$, and similarly for $Y$.\end{lem}

If $(M,u,v)$ is a doubly-foliated manifold, we may look at the bundle of fibre-tangent vectors, which we will denote (abusing notation) by $TS \subseteq TM$. We also consider its dual bundle, $T^\ast S$. Sections of $T^\ast S$ may be interpreted as a smoothly-varying family of $1$-forms on the fibres $S_{u,v}$. We will often adopt this point of view. One has a canonical inclusion operator $\iota: TS \to TM$ and a dual operator $\iota^\ast: T^\ast M \to T^\ast S$. These have left and right inverses, respectively, but the inverses are not canonical. This remark extends to type $(0,q)$ and type $(p,0)$ tensors, respectively. One can also form tensor products of these bundles, leading to the notion of general type $(p,q)$ tensors. We will call such a tensor a fibre tensor. The following lemma defines the Lie derivative of fibre tensors along vector fields on $M$ which are $\pi$-related to vector fields on $\mathcal R$.

\begin{lem}[\cite{ChrForm}*{Lemma~1.1}]\label{thm:C2:LieDerI}Let $X,\ Y \in C^{\infty}(M;TM)$ be $\pi$-related to vector fields $X',\ Y' \in C^{\infty}(\mathcal R; T\mathcal R)$. Then $[X,Y]$ is $\pi$-related to $[X',Y']$. In particular, if $X$ is a lift of a coordinate vector field $\pa_u$ or $\pa_v$, and observing that the fibre-tangent vector fields are those $\pi$-related to $0$, if $V$ is a fibre-tangent vector field, then $\Lie_X V = [X,V]$ is as well. 
\end{lem}
\begin{proof}This is a standard property of submersions.\end{proof}

Using this, we may define the Lie derivative of any fibre tensor along any vector field on $M$ which are $\pi$-related to one on $\mathcal R$. We have just treated the case of vector fields. By tensoring, it suffices to now treat the case for $1$-forms.

\begin{lem}[\cite{ChrForm}*{Lemma~1.3}]\label{thm:C2:LieDerII}Let $X \in C^{\infty}(M;TM)$ be $\pi$-related to some vector field $Y \in C^{\infty}(\mathcal R; T\mathcal R)$. Let $\omega \in C^{\infty}(M;T^\ast S)$ be a fibre one-form. Extend $\omega$ arbitrarily to a one form $\tilde{\omega} \in C^{\infty}(M;T^\ast M)$. Then $\iota^\ast \Lie_X \tilde{\omega}$ does not depend on the choice of extension, and hence $\Lie_X \omega := \iota^\ast \Lie_X \tilde{\omega}$ is well-defined.\end{lem}
\begin{proof}If $V$ is a fibre vector field, then
\[\Lie_X \tilde{\omega}(V) = X\tilde{\omega}(V) - \tilde{\omega}(\Lie_X V).\] The first term is equal to $X\omega(V)$. Using \cref{thm:C2:LieDerI}, $\Lie_X V$ is fibre-tangent, and thus the second term is $\omega(\Lie_X V)$ which does not depend on the choice of extension.\end{proof}

We will make the distinction below between double-null foliation and a double-null gauge for a Lorentzian metric $g$, since in our definition the former is not enough to kill the diffeomorphism group. A double-null gauge will be a double-null foliation with extra data.
\begin{defn}Let $(M,u,v)$ be a doubly-foliated manifold, and $g$ a Lorentzian metric on $M$. We say that $(M,g,u,v)$ is a \emph{double-null foliation} if $u,v$ are optical functions for $g$, i.e.\ 
\[g^{-1}(du,du) = g^{-1}(dv,dv) = 0.\] We additionally suppose that $(M,g)$ is time orientable and that $u,v$ are increasing towards the future.\end{defn}
The level sets of $u,v$ are the ``null foliations'' of the double-null foliation.

If $(M,g,u,v)$ is a double-null foliation, we define $\iota^\ast g = \slash{g}$ to be the restriction of $g$ to $TS$ and $N' = -2\grad v$ and $L' = -2\grad u$ to be the future directed null geodesic generators of the levels sets of $v$ and $u$ respectively, and putatively define $0 < \Omega$ by $-2\Omega^{-2} = g(N',L')$. The function $\Omega$ is not yet well-defined, since it requires $g(N',L') < 0$. That it is in fact well-defined, and more, follow from the following:
\begin{prop}[Geometry of double-null foliations]\label{thm:C2:DNGeo} Let $(M,g,u,v)$ be a double-null foliation. Then $\slash{g}$ is Riemannian (i.e.\ the fibres of $(u,v)$ are spacelike), $g(N',L') < 0$, and $\vspan\{N',L'\}$ is transverse to $TS$. \end{prop}
\begin{rk}In fact, the proposition is true more generally. Let $(V,g)$ be a vector space equipped with a Lorentzian metric, and let $T$ be a timelike vector providing a time orientation for $V$. Suppose $\mu$ and $\nu$ are linearly independent null covectors satisfying $\mu(T), \nu(T) > 0$. Set $L' = g^{-1}(-2du,\bullet)$ and $N' = g^{-1}(-2dv,\bullet)$, and let $W = \vspan\{L',N'\}^{\bot} = \ker \mu \n \ker \nu$. Then the proposition remains true with $(\mu,\nu)$ replacing $(du,dv)$, and with $\vspan\{L',N'\}^{\bot}$ replacing $TS$ (since in a double-null foliation, $TS = \vspan\{L',N'\}$ by definition). As we do not need the general statement, we only prove the proposition using the language of double-null gauges.\end{rk}
\begin{proof}
Since $du$ and $dv$ are linearly independent null covectors, $L'$ and $N'$ are linearly independent and $g^{-1}(du,dv) \neq 0$. To show this, first diagonalize $g^{-1}$ over a point $p \in M$ to find a timelike vector $w$ with $g^{-1}(w,w) = -1$ and an orthogonal spacelike hyperplane $\Sigma$. We may find $\alpha,\beta,\gamma,\delta \in \R$ and $\sigma,\tau \in \Sigma$ of length $1$ such that $du = \alpha w + \beta\sigma$ and $dv = \gamma w + \delta\tau$. Since $du$ and $dv$ are null, $|\alpha| = |\beta|$ and $|\gamma| = |\delta|$ and none of $\alpha,\beta,\gamma,\delta$ are $0$. Suppose for contradiction that $g^{-1}(du,dv) = 0$.  Then also $\alpha\gamma = \beta\delta g^{-1}(\sigma,\tau)$. Since $g^{-1}|_{\Sigma\times\Sigma}$ is a (positive-definite) inner product, we may apply the Cauchy-Schwarz inequality to conclude that $|\alpha\gamma| \leq |\beta\delta|$, with equality holding if and only if $\sigma$ and $\tau$ are scalar multiples of one another. Since $|\alpha\gamma| = |\beta\delta|$, equality holds and $\sigma$ and $\tau$ are scalar multiples of each other, and so $\sigma = \pm\tau$. Thus $g^{-1}(\sigma,\tau) = \pm 1$, and so $\alpha\gamma = \beta \beta\delta$ (with the same sign as in $\sigma = \pm \tau)$ and hence $\beta\gamma/\delta = \pm \beta^2/\alpha = \pm \alpha$, and so \[\beta/\delta dv = \pm \alpha w + \beta\tau = \pm \alpha w \pm \beta\sigma = \pm du,\] which contradicts linear independence.

We next claim that a nontrivial linear combination of $N',L'$ is linearly independent of $TS$, and thus $\vspan(N',L')$ and $TS$ are transverse. Indeed if $aN' + bL'$ ($a,b \in \R$) were in $TS$, then $0 = g(L',aN'+bL') = ag(L',N')$, which means $a = 0$ (since $g^{-1}(dv,du) \neq 0$), and likewise $0 = g(N',aN'+bL') = bg(N',L')$ implies $b = 0$. 

Now we show that $\slash{g}$ is Riemannian.

We may express $g$ acting between $N',L',TS$ as a block-diagonal matrix
\[g = \begin{pmatrix}0 & g(N',L') & 0\\
g(N',L') & 0 & 0\\
0 & 0 & \slash{g}\end{pmatrix}.\] Partially diagonalizing this matrix as
\[\begin{pmatrix} -g(N',L') & 0 & 0\\
0 & g(N',L') & 0\\
0 & 0 & \slash{g}\end{pmatrix}\] and using that $g$ is Lorentzian, it follows that $\slash{g}$ must be Riemannian.

Finally, we show that $g(N',L') < 0$. If $T$ is timelike vector field for which $du\cdot T, dv\cdot T > 0$\footnote{Such a vector exists because $u,v$ are increasing towards the future.} then we may write $T = aN' + bL' + \Theta$, for some $a,b \in \R$ and $\Theta \in TS$ (since now we know that $\vspan(N',L')$ is transverse to $TS$). So $0 > g(T,L') = ag(N',L')$ and $0 > g(T,N') = bg(L',N')$. Thus $a$ and $b$ have the same sign. Hence
\[0 > g(T,T) = abg(N',L') + g(\Theta,\Theta),\] and so $g(N',L') < 0$.
\end{proof}
\begin{rk}Notice that the proposition implies immediately that $u+v$ is a time function for $(M,g)$.\end{rk}

It will be more convenient to work with rescaled versions of $N',L'$.
Set $N = \Omega^2 N'$ and $L = \Omega^2 L'$ to be the normalized null vector fields. Observe that $N$ and $L$ are lifts of $\pa_u$ and $\pa_v$ respectively, i.e.\
\begin{align*}
Nu = 1, \ Nv = 0\\
Lu = 0, \ Lv = 1.\end{align*}
Since $[N,L]$ is $\pi$-related to $[\pa_u,\pa_v] = 0$, $[N,L]$ is a fibre vector field.

We will use slash notation, such as $\slash{\grad}$, $\slash{\div}$ to denote operations with respect to $\slash{g}$ on the fibres $S_{u,v}$, rather than with respect to $g$ on $M$. We also define the torsion $\zeta$, a fibre 1-form by $\zeta(V) = \frac{1}{2}g(\nabla_V (L/\Omega),(N/\Omega))$ (where $\nabla$ denotes the Levi-Civita connection of $g$). We remark that $\zeta = \frac{1}{4\Omega^2}\slash{g}([N,L],\bullet)$. This follows from the computation of the Levi-Civita connection carried out in Appendix~\ref{C:A1:connection}.

 These components uniquely determine the metric:
\begin{prop}\label{thm:C2:Reconstruct}Let $(M,u,v)$ be a doubly-foliated manifold. Let $N$, $L$ be vector fields on $M$ which are lifts of $\pa_u,\ \pa_v$, respectively, let $\slash{g}$ be a fibre Riemannian metric, and $0 < \Omega \in C^\infty(M;\R)$. Then there is a unique Lorentzian metric $g$ such that $(M,g,u,v)$ is a double-null foliation and $g$ induces the given data, i.e.\ $\iota^\ast g = \slash{g}$, $-2\Omega^{-2} = g^{-1}(-2du,-2dv)$, $N = -2\Omega^2 \grad v$, $L = -2\Omega^2 \grad u$.\end{prop}
\begin{proof}It is clear how we should define $g$. \Cref{thm:C2:DNGeo} shows that at a point $p$, $N,L,T_pS$ together span $TM$ so we just need to specify how $g$ should act between all of them. We should have, for $V,W \in T_pS$,
\begin{itemize}
\item $g(N,N) = g(L,L) = 0$;
\item $g(N,L) = -2\Omega^2$;
\item $g(N,V) = g(L,V) = 0$;
\item $g(V,W) = \slash{g}(V,W)$.
\end{itemize}
This proves uniqueness. For existence, we just need to check that this $g$ is actually in a double null gauge with $u,v$. First of all, $g$ is certainly Lorentzian. Second, the vector field $N+L$ provides a time orientation along which $u,v$ are increasing. We need to show $g^{-1}(du,du) = g^{-1}(dv,dv) = 0$, $g^{-1}(-2du,-2dv) = -2\Omega^{-2}$ and $N = -2\Omega^2\grad v$, $L = -2\Omega^2\grad u$. All follow from the last two. By non-degeneracy, we just need to check that $N$ and $-2\Omega^2\grad v$ have the same values when tested again $N,L,V$, and same for $L$. Indeed,
\begin{itemize}
\item $g(-2\Omega^2\grad v,N) = -2\Omega^2Nv = 0$;
\item $g(-2\Omega^2\grad v,L) = -2\Omega^2Lv = -2\Omega^2$;
\item $g(-2\Omega^2\grad v,V) = -2\Omega^2Vv = 0$,
\end{itemize}
and likewise for $L$.\end{proof}

A double-null foliation $(M,g,u,v)$ is not enough to break the diffeomorphism invariance of the Einstein equations. Indeed, all diffeomorphisms are $\phi:M \to M$ are allowed, provided $\phi^{\ast}u = u$ and $\phi^{\ast}v = v$. To remedy this, we need to add in some extra data and introduce a double-null \emph{gauge}.
\begin{defn}\label{def:C2:doublenullgauge}Let $(M,g,u,v)$ be a double-null foliation. Let $N \in C^\infty(M;TM)$ be a lift of $\pa_u$. We call the quintuple $(M,g,u,v,N)$ \emph{double-null gauge} if $N = -2\Omega^2 \grad v$ (with $\Omega$ as above). In other words, $N$ must equal the $N$ introduced above.\end{defn}
 If $(M,g,u,v,N)$ is a double-null gauge, we will often say that $g$ is in a double null gauge with $u,v,N$ (the background manifold $M$ being clear from context).

Versions of double-null gauges have been used before in the literature (see for instance \cites{ChrForm, LukLoca}) under the name ``canonical coordinates,'' but the definition given here is more coordinate invariant. The relation is as follows. Choose some $(u_0,v_0) \in \mathcal R$ and set $\mathcal H = u^{-1}(u_0)$. Choose arbitrary local coordinates $\theta^i$ on $S_{u_0,v_0}$. This induces coordinates $(u,v,\theta^i)$ on $M$ by flowing out $\theta^i$ via $L$ along $\mathcal H$, and and then via $N$ to all of $M$. In these coordinates, $g$ looks like
\[g = -2\Omega^2(du\otimes dv + dv\otimes du) + \slash{g}_{ij}(d\theta^i-f^idv)\otimes(d\theta^j-f^jdv),\]
where
$\partial_v f^i = N f^i = [N,L]^i$ and $f^i = 0$ on $\mathcal H$. The above coordinates are what are referred to as ``canonical.'' This procedure also allows one to obtain a trivialization $M \iso \mathcal R \times S_{u_0,v_0}$.

A double-null gauge fixes the diffeomorphism invariance of the Einstein equations, up to a choice of diffeomorphism fixing some data on the initial surfaces. Fix a doubly-foliated manifold $(M,u,v)$, and let $\mathcal H_1 = v^{-1}(v_0)$ and $\mathcal H_2 = u^{-1}(u_0)$ denote fibres of $v$, $u$ respectively at $v_0$ and $u_0$. Also fix $N$, a lift of $\pa_u$ to $M$. Then:
\begin{prop}\label{thm:C2:diffeo}Let $\phi: M \to M$ be a diffeomorphism mapping $\mathcal H_i$ into $\mathcal H_i$ ($i = 1,2$) such that
\begin{romanumerate}
\item $\phi|_{\mathcal H_2} = \id$;
\item $\phi^\ast u = u$ on $\mathcal H_1$;
\item $du = d(\phi^\ast u)$ and $dv = d(\phi^\ast v)$ on $S_{u_0,v_0}$.
\end{romanumerate}
Suppose $g$ and $\phi^\ast g$ are both in a double null gauge with $u,v,N$. Then $\phi = \id$ on all of $M$.\end{prop}
\begin{proof}By assumption
\[g^{-1}(\phi^\ast du,\phi^\ast du) = g^{-1}(\phi^\ast dv,\phi^\ast dv) = 0.\] 

Also by assumption, $u = \phi^\ast u$ on $\mathcal H_1$ and $du = d(\phi^\ast u)$ on $S_{u_0,v_0}$. Now observe that $L'$, which is the projection of the bicharacteristic flow of $g^{-1}(\xi,\xi)$ with initial data $-2du$ over $\mathcal H_1$, is transverse to $\mathcal H_1$, since $\mathcal H_1$ is spanned by $N'$ and $TS$. Hamilton-Jacobi theory applied to the equation \[g^{-1}(du,du) = g^{-1}(d(\phi^\ast u),d(\phi^\ast u)) = 0\] now implies that $u = \phi^\ast u$ on all of $M$. The same argument implies that $v = \phi^\ast v$ on $M$. In particular, $\Omega^2_{g} = \Omega^2_{\phi^\ast g}$, and so $N = \phi_\ast N$.

The action of $\phi$ may thus be described as follows. For a point $p \in S_{u,v}$, flow backwards along $N$ for time $u-u_0$ until the flow reaches $\mathcal H_2$, then apply $\phi = \id$, and flow forwards via $\phi_\ast N = N$ for time $u\circ \phi(p) -u_0 = u-u_0$. Thus $\phi = \id$.\end{proof}

Double-null gauges are not rare. Any metric (avoiding the obvious obstruction that $u,v$ need to be future oriented) may locally be put into a double-null gauge. Let $(M,u,v)$, $\mathcal H_i$ ($i = 1,2$) and $N$ be as above. Suppose $g$ is a metric on $M$ and further suppose that $\mathcal H_i$ ($i = 1,2$) are null for $g$ and that $g$ has a time orientation for which $u$ and $v$ are increasing towards the future. For the moment let us also suppose that $(u_0,v_0) \in \mathcal R^\circ$ (although we will reduce to this case in \cref{thm:C2:onbd}).
\begin{prop}\label{thm:C2:badconversion}There exists a neighbourhood $U$ of $S_{u_0,v_0}$, and a smooth map $\phi:U \to M$, a diffeomorphism onto its image, $M'$, mapping $\mathcal H_i$ into $\mathcal H_i$ ($i = 1,2$) and satisfying
\begin{romanumerate}
\item $\phi|_{\mathcal U\n H_2} = \id$;
\item $\phi^\ast u = u$ on $U\n \mathcal H_1$;
\item $du = d(\phi^\ast u)$ and $dv = d(\phi^\ast v)$ on $S_{u_0,v_0}$;
\end{romanumerate}
such that $(M',\phi_\ast g,u,v,N)$ is a double-null gauge. Two such diffeomorphisms agree on their common domain of definition.\end{prop}
\begin{proof}\Cref{thm:C2:diffeo} shows uniqueness.
Using Hamilton-Jacobi theory, let us first show the existence of functions $u',\ v'$ in a neighbourhood $U$ of $S_{u_0,v_0}$ satisfying
\[g^{-1}(du',du') = g^{-1}(dv',dv') = 0\]
and $u' = u$ on $\mathcal H_1$, $v' = v$ on $\mathcal H_2$, $du' = du$, $dv' = dv$ on $S_{u_0,v_0}$. Let $N'$ and $L'$ denote the vector fields generating the future-directed null-geodesic congruences tangent to $\mathcal H_1$ and $\mathcal H_2$, respectively. Since $v$ and $u$ are constant along the null hyperplanes $\mathcal H_1$ and $\mathcal H_2$, respectively, $\grad v$ and $\grad u$ are perpendicular to them, and hence are multiples of $N'$ and $L'$, respectively. We may normalize $N'$ and $L'$ by the initial condition $N' = -2\grad v$ and $L' = -2\grad u$ on $S_{u_0,v_0}$.

In particular $g^{-1}(du,du) = d^{-1}(dv,dv) = 0$ on $S_{u_0,v_0}$, $L'$ is transverse to $\mathcal H_1$ and $N'$ is transverse to $\mathcal H_2$ at $S_{u_0,v_0}$. Writing $(x,\xi)$ for a generic cotangent variable, the Hamiltonian vector field of $-\frac{1}{2}g^{-1}(\xi,\xi)$ is
\[H = -g^{-1}(\xi,\bullet) + \frac{1}{2}(\pa_{x}g^{-1})(\xi,\xi)\pa_{\xi}.\] For $x \in S_{u_0,v_0}$ and $\xi = du$, it follows that $H = \alpha L' + V$, for some nonvanishing $\alpha$ and vector field tangent to the fibres of the projection $T^\ast M \to M$. In particular, the projection of $H$ evaluated at $(x,du)$ to $M$ is transverse to $\mathcal H_1$ at $S_{u_0,v_0}$. Likewise, the projection of $H$ evaluated at $(x,dv)$ to $M$ is transverse to $\mathcal H_2$ at $S_{u_0,v_0}$ projection of $H$ to over $S_{u_0,v_0}$. Therefore, Hamilton-Jacobi theory provides the existence of $U$, $u'$ and $v'$.

Now, let us prove that $v'|_{U\n\mathcal H_1} = v_0$ and $u'|_{U\n \mathcal H_2} = u_0$. Using the correspondence between the bicharacteristic flow of $-\frac{1}{2}g^{-1}(\xi,\xi)$ and the cogeodesic flow of $g$, it follows from the Hamilton-Jacobi theoretic construction of $v'$ that it is constant along the null geodesic starting at $S_{u_0,v_0}$ with initial direction $N'/2$. Indeed, If $\gamma(t)$ is an integral curve of $N'/2$\footnote{i.e.\ $\gamma$ is a null geodesic tangent to $\mathcal H_1$.} starting at a point in $S_{u_0,v_0}$, then
\[(x(t),\xi(t)) := (\gamma(t),-g(\gamma'(t),\bullet)) = (\gamma(t),g(N'/2,\bullet))\] is a bicharacteristic of $-\frac{1}{2}g^{-1}(\xi,\xi)$, which at time $t=0$ is equal to $(x,dv')$. So $dv' = -g^{-1}(N'/2,\bullet)$ over $\gamma(t)$, and hence 
\[N'/2v' = dv'\cdot N'/2 = g(N'/2,-N'/2) = 0.\] Since $v' = v = v_0$ on $S_{u_0,v_0}$, $v' = v_0$ on all of $\mathcal H_1 \n U$. The same argument works for $u'$.

Shrinking $U$, we will show that $(U,u',v')$ is a doubly-foliated manifold. Shrinking $U$, we may assume that $du'$ and $dv'$ are independent on $T^\ast U$, and are future directed for $g$, since $du' = du$ and $dv' = dv$ on $S_{u_0,v_0}$. In particular, $(u',v')$ is a submersion onto its image in $\R^2$. 

Let $V \subseteq U$ be a compactly contained neighbourhood of $S_{u_0,v_0}$ (this is possible since $S_{u_0,v_0}$ is compact). Since $(u',v')$ is a submersion, if $\epsilon > 0$ is sufficiently small then $(u',v')(V)$ contains $(u_0-\epsilon,u_0+\epsilon)\times(v_0-\epsilon,v_0+\epsilon)$, and let us replace $U$ by the preimage of this set. With this domain, $(u',v')$ is a proper map, too. Indeed, $(u',v')$ was defined over $V$, so is proper given domain $U$. Therefore, \[(u',v'):U \to (u',v')(U) = (u_0-\epsilon,u_0+\epsilon)\times(v_0-\epsilon,v_0+\epsilon)\] is also proper.

From general properties of proper submersions, it follows that each fibre is diffeomorphic to $S_{u_0,v_0} \iso S$. Thus we have proven that $(U,u',v')$ is a doubly-foliated manifold. Denote by $S'_{u',v'}$ a generic fibre of $(u',v')$. We remark that $\{u' = u_0\} = \mathcal H_2 \n U$ and $\{v' = v_0\} = \mathcal H_1 \n U$. Indeed, $\{u' = u_0\}$ is a connected hyperplane containing $\mathcal H_2 \n U$, but $\mathcal H_2 \n U$ has no boundary in $U$, and similarly for $\{v' = v_0\}$.

Let $\tilde{N} = -2\Omega^2 \grad v'$, where $-2\Omega^{-2} = g(du',dv')$. Then $(U,g,u',v',\tilde{N})$ is a double-null gauge. Since $U\n\mathcal H_1 \subseteq \{v' = v_0\}$, it follows that $\tilde{N}$ is tangent to $U\n \mathcal H_1$. Recall that $\tilde{N} u' = 1$, and $\tilde{N} v' = 0$.

 We now need to construct $\phi$. If $p \in U$, we may flow back from $p$ along $\tilde{N}$ for time $u'(p)-u_0$ to $\{u' = u_0\} = \mathcal H_2 \un U$. Then flow forwards along $N$ for time $u'(p)-u_0$. Denote the flow of $N$ by $\Phi$ and the flow of $\tilde{N}$ by $\Psi$. In symbols,
 \[\phi(p) = \Phi(u'(p)-u_0,\Psi(-(u'(p)-u_0),p)).\] 
 Since $Nu = 1$, $u(\Phi(s,q)) = u(q) + s$ for any point $q$ and time $s$ for which $\Phi(s,q)$ is defined. Thus, since $u = u_0$ on $\mathcal H_2$,
 \[u(\phi(p)) = u(\Psi(-(u'(p)-u_0)),p) + u'(p)-u_0 = u_0 + u'(p)-u_0 = u'(p),\]
 and so $\phi^\ast u = u'$. Since $Nv = 0$, $v(\Phi(s,q)) = u(q)$ for any point $q$ and time $s$ for which $\Phi(s,q)$ is defined, and since $\tilde{N}v' = 0$, $v'(\Psi(s,q)) = v'(q)$ for any point $q$ and time $s$ for which $\Psi(s,q)$ is defined. Thus, since $v = v'$ on $\mathcal H_2$,
 \[v(\phi(p)) = v(\Psi(-(u'(p)-u_0)),p) = v'(\Psi(-(u'(p)-u_0)),p) = v'(p),\]
 and so $\phi^\ast v = v'$.
 
 Since $\tilde{N}u' = 1$, $u'(\Psi(s,q)) = u'(q)+s$ for any point $q$ and time $s$ for which $\Psi(s,q)$ is defined. Thus
 \begin{align*}
d\phi \tilde{N}_p &= \left.\frac{d}{ds}\right|_{s = 0} \phi(\Psi(s,p))\\
&= \left.\frac{d}{ds}\right|_{s = 0} \Phi(u'(p)+s-u_0,\Psi(-(u'(p)+s-u_0),\Psi(s,p)))\\
&= \left.\frac{d}{ds}\right|_{s = 0} \Phi(u'(p)+s-u_0,\Psi(-(u'(p)-u_0),p))\\
&= N_{\phi(p)},
 \end{align*}
 and so $\phi_\ast \tilde{N} = N$. 

Let us now show that $\phi$ is a diffeomorphism onto its image. From the above, $\phi:S'_{u',v'} \to S_{u',v'}$, and is in fact a diffeomorphism, since it is given by the composition of flows. In particular $\phi$ defined on $U$ is injective and surjective onto its image. However, $d\phi$ is injective, and thus $\phi$ is a global diffeomorphism. Indeed, if $d\phi \cdot x = 0$ for some $x \in TU$, then 
\[du' \cdot x = d(\phi^\ast u)\cdot x = 0 = d(\phi^\ast v)\cdot x = dv' \cdot x,\]
and thus $x \in TS'_{u',v'}$, and so $d\phi \cdot x = 0$ implies $x=0$ since $d\phi|_{TS'_{u',v'}}$ is injective.

Now let us check that $\phi$ satisfies (i), (ii) and (iii). It is clear that (i) holds. Property (ii) follows from observing that $u = u'$ on $U \n \mathcal H_1$. Property (iii) is true since $du' = du$ and $dv' = dv$ on $S_{u_0,v_0}$.

Finally, it is clear by definition that $(M'= \phi(U),\phi_\ast g,u,v,N)$ is a double-null gauge.\end{proof}

In the previous proposition, we assumed that $(u_0,v_0)$ lay in the interior of $\mathcal R$. However, it is possible that at least one of $u_0$ or $v_0$ lies in the boundary of $\mathcal R$. The conclusion still holds assuming a mild non-degeneracy assumption on the position of $(u_0,v_o) \in \pa \mathcal R$. We will assume that there is a neighbourhood $V$ of $(u_0,v_0)$ in which $\pa \mathcal R$ is locally $u^{-1}(u_0) \n V$, $v^{-1}(v_0) \n V$, or $(u^{-1}(u_0) \un v^{-1}(v_0)) \n V$). We show that we can reduce this case to the case that $(u_0,v_0)$ lies in the interior of $\mathcal R^\circ$.
\begin{cor}\label{thm:C2:onbd}The conclusions of \cref{thm:C2:badconversion} continue to hold.\end{cor}
\begin{proof}The hardest case is when $(u_0,v_0)$ is a corner of $\mathcal R$, and we treat this case. Let us also assume for simplicity that $u \geq u_0, \ v \geq v_0$ in a neighbourhood of $(u_0,v_0)$ in $\mathcal R$ (the other cases are similar). Let us replace $\mathcal R$ with $[u_0,u_0+\epsilon]\times[v_0,v_0+\epsilon]$ for $\epsilon$ sufficiently small.
We may trivialize $M \iso \mathcal R \times S_{u_0,v_0}$, and for $\epsilon' > 0$ small extend $\mathcal R$ to $\widetilde{\mathcal R} = [u_0-\epsilon',u_0+\epsilon]\times[v_0-\epsilon',v_0+\epsilon]$ so that $(u_0,v_0)$ lies in the interior of $\widetilde{\mathcal R}$, and extend the projection to a map $(u,v)$ from $\widetilde{M} := \widetilde{\mathcal R}\times S_{u_0,v_0} \to \widetilde{\mathcal R}$ so that $(\widetilde{M},u,v)$ is a doubly-foliated manifold. Extend $g$ and $N$ to $\widetilde{M}$, and also extend the timelike vector field $T$ giving $g$ its time orientation to $\widetilde{M}$. Notice that $M = \{u \geq u_0, \ v \geq v_0\} \n \tilde{M}$. Now we apply the proposition to $\tilde{M}$ obtaining a diffeomorphism $\phi$ from a neighbourhood $U$ of $S_{u_0,v_0}$ inside $\tilde{M}$ onto its image. By construction, in the notation of the proof, $(U,g,u',v',\tilde{N})$ is a double-null gauge. The key observation is that the vector fields $\tilde{N}$ and $\tilde{L} := -2\Omega^2 \grad u'$ are lifts of $\pa_{u'}$ and $\pa_{v'}$, respectively, from the base of the underlying doubly-foliated manifold. Shrinking $U$, we may also assume that $\tilde{L}v > 0$ and $\tilde{N}u > 0$ on $\mathcal H_1\n U$ and $\mathcal H_2\n U$, respectively (since this is certainly true on $S_{u_0,v_0}$).

To show the conclusions for $M$, rather than the extension $\tilde{M}$, we only need to check that $\phi(U
\n M) = \phi(U) \n M$. Since $\phi^\ast u = u'$, $\phi^\ast v = v'$, it suffices to show that for $p \in U$ 
\begin{equation}\label{eq:C2:idk}u'(p) \geq u_0, \ v'(p) \geq v_0\text{ if and only if }p \in U \n M \subsetneq U \n \tilde{M}.\end{equation} We already know that $u'(p) = u_0$ or $v'(p) = v_0$ if and only if $p \in U \n \mathcal H_2$ or $U \n \mathcal H_1$, respectively, so we just need to focus on $U \setminus(\mathcal H_1 \un \mathcal H_2)$. Using $\tilde{L}$ or $\tilde{N}$ to flow from $p \in U\setminus(\mathcal H_1 \un \mathcal H_2)$ to $\mathcal H_1 \un \mathcal H_2$, \eqref{eq:C2:idk} will be true provided that the forwards and backwards flow of $\tilde{L}$ for nonzero time from $\mathcal H_1$ stays inside $\{v > v_0\}$ and $\{v < v_0\}$, respectively, and the forwards and backwards flow of $\tilde{N}$ for nonzero time from $\mathcal H_2$ stays inside $\{u < u_0\}$ and $\{u > u_0\}$, respectively. Without loss of generality let us prove the second of these. 

Since $du N > 0$ on $\mathcal H_2$, the claim is true for small nonzero times. Suppose it is not true for all times. Then there exists $p \in \{u > u_0\}$ and $q \in \{u < u_0\}$ and a forwards-directed integral curve $\gamma$ of $\tilde{N}$ with $\gamma(0) = p$, $\gamma(1) = q$. Let $t^\ast$ be the first time for which $u(\gamma(t^\ast)) \in \mathcal H_2$. Then \[0 < du \tilde{N}_{\gamma(t^\ast)} = (u\circ \gamma)'(t^\ast) \leq 0,\] since $u\circ \gamma$ is non-increasing at the first time of intersection. This is a clear contradiction.\end{proof}

It is inconvenient that in general one may not choose the diffeomorphism $\phi$ to be the identity on $\mathcal H_1$. However, if we assume by fiat that the vector field $N$ is a null (rescaled) geodesic for $g$ along $\mathcal H_1$, then we may arrange this.
\begin{cor}\label{thm:C2:conversion}Suppose that $N$ is a null (rescaled) geodesic for $g$ tangent to $\mathcal H_1$. Then, there exists a neighbourhood $U$ of $S_{u_0,v_0}$, and a smooth map $\phi:U \to M$, a diffeomorphism onto its image, $M'$, which is the identity on $(\mathcal H_1 \un \mathcal H_2)\n U$, such that $(M',\phi_{\ast}g,u,v,\tilde{N})$ is a double-null gauge. Two such diffeomorphisms are unique on their common domain of definition.\end{cor}
\begin{proof}
The assumptions are stronger than those of \cref{thm:C2:badconversion}, so the conclusion still holds, providing a unique diffeomorphism $\phi$ which we need to show is additionally the identity on $U\n \mathcal H_1$. Let $\tilde{N}$ be the vector field from the proof of \cref{thm:C2:badconversion}. It suffices to show that $\tilde{N}= N$ on $\mathcal H_1$. By construction, $\tilde{N}$ is a null vector field tangent to $\mathcal H_1$ which satisfies $\tilde{N}u' = 1$. Since $u' = u$ on $\mathcal H_1$, $\tilde{N}u = 1$. Thus $\tilde{N}$ and $N$ are both null vector fields tangent to $\mathcal H_1$ and which satisfy $\tilde{N}u = 1 = Nu$. It follows that $\tilde{N} = N$.
\end{proof}

\section{The characteristic initial value problem}
\label{C:C2:CIVPSectionI}
In this section, we state the characteristic initial value problem for the Einstein vacuum equations and outline our proof.
Let $(M,u,v)$ be a doubly-foliated manifold, let $N$ be a lift of $\pa_u$, and assume that $M = M([0,a)\times [0,b))$ for some $a,b > 0$ (and possibly $a, b = \infty$). Set $\mathcal H_1 = \{v=0\}$ and $\mathcal H_2 = \{u = 0\}$. In our context, the local well-posedness of the characteristic initial value problem to the Einstein vacuum equations states:

\begin{thm}[Rendall \cite{RenRedu}, Luk \cite{LukLoca}]\label{thm:C2:CIVP}Let $(M,u,v)$ and $N$ be as above. Suppose the following data are given:
\begin{enumerate}[label = (\roman*)]
\item a vector field $L$, a lift of $\pa_v$ to $T\mathcal H_2$;
\item a Riemannian metric $\slash{g}$ on $S_{0,0}$;
\item a section $\hat{\slash{g}}$ of $\Sym^2(T^\ast S)$ over $\mathcal H_1\un \mathcal H_2$ which has Riemannian signature and satisfies $\hat{\slash{g}} = \slash{g}$ at $S_{0,0}$;
\item a pair of smooth functions $f_1$, $f_2$ on $S_{0,0}$;
\item a fibre 1-form $\zeta$ on $S_{0,0}$ (or a fibre vector field $W$ on $S_{0,0}$).
\end{enumerate}Then there exist $0 < u^\ast \leq a$, $0 < v^\ast \leq b$ and a unique Lorentzian metric $g$ on $M' = M([0,u^\ast)\times[0,v^\ast))$ satisfying $\Ric(g) = 0$ such that $(M',g,u,v,N)$ is a double-null gauge and
\begin{enumerate}[label = (\roman*)]
\item $L|_{\mathcal{H}_2} = -2\Omega^2\grad u$;
\item $\iota^\ast g = \slash{g}$ over $S_{0,0}$;
\item $\iota^\ast g$ lies in the conformal class of $\hat{\slash{g}}$ over $\mathcal H_1 \un \mathcal H_2$;
\item $f_1$ is the mean curvature of $\slash{g}$ on $S_{0,0}$ in the direction of $N$, and $f_2$ is the mean curvature of $\slash{g}$ on $S_{0,0}$ in the direction of $L$;
\item $\zeta$ is the torsion of $g$ at $S_{0,0}$ (or $W = [N,L]$ on $S_{0,0}$);
\item $\Omega = 1$ on $\mathcal H_1 \un \mathcal H_2$.
\end{enumerate}
\end{thm}

\begin{figure}[htbp]
\centering
\includegraphics{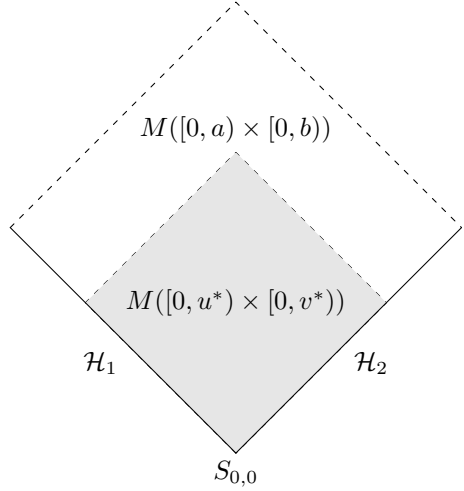}
\caption{The doubly-foliated manifold in \cref{thm:C2:CIVP}. The region of existence is shaded. Fibre directions not pictured.}
\end{figure}

Our proof of \cref{thm:C2:CIVP} will proceed similarly to the original proof of Rendall, that is it will proceed by reduction to the Cauchy problem. Our approach differs from both his approach and that of Luk by not requiring the metric to be in harmonic gauge with a Minkowski background. Rather, we first find a metric that solves the Einstein equations in Taylor series at the boundary $\mathcal H_1 \un \mathcal H_2$. Then we perturb off the solution in Taylor series and construct the solution metric by requiring it to be in a wave gauge with respect to the solution in Taylor series. Subsequently, we use \cref{thm:C2:conversion} to switch gauges back to a double-null gauge. The advantage of a double-null gauge is that the Einstein equations in a double-null gauge are well-adapted to Taylor series computations.

More precisely, we first prove:
\begin{thm}\label{thm:C2:TSI}Let initial data be given as in \cref{thm:C2:CIVP}. Then there exist $u', v' > 0$ and a smooth Lorentzian metric $g$ on $M([0,u')\times[0,v'))$ satisfying $\Ric(g) = 0$ in Taylor series at $\mathcal H_1\un \mathcal H_2$ such that $(M',g,u,v,N)$ is a double-null gauge and induces the initial data in the sense of \cref{thm:C2:CIVP}. If $\widetilde{g}$ is any other such metric satisfying the same conditions then $g=\widetilde{g}$ in Taylor series. Moreover, there is an iterative algorithm to compute the Taylor series expansions at $\mathcal H_1$ and $\mathcal H_2$.\end{thm}
\begin{rk}Our method may be seen as reversing the method of Rendall. Rendall is able to solve in Taylor series any second-order hyperbolic equation, which he uses as his first step to solving the equation exactly. Then, he finds harmonic coordinates which allow him to reduce to this case. Our method is to find the Taylor series for the Einstein equations \emph{first} and use this solution as our gauge. The advantage is that our method is inherently global on the fibres $S$, and we will modify this argument to prove the \cref{thm:C1:ChrVague}.\end{rk}

Then, we will try to find a tensor $h$, zero in Taylor series at $\mathcal H_1 \un \mathcal H_2$ so that $g+h$ is in wave gauge with respect to $g$ and $g+h$ solves the reduced Einstein equations $\widetilde{\Ric}(g+h) = 0$.
The operator $Q(h) := \widetilde{\Ric}(g+h)$ is quasilinear hyperbolic. We are then in a position to use a theorem of Rendall, which we now describe. Let $V$ be a vector bundle over $M$, with connection $D$, and $P$ a second-order quasilinear hyperbolic operator given by
\[Ph = \psi(h)D^2 h + Lh,\] where $L$ is a first-order non-linear differential operator and $\psi(h)$ is a Lorentzian metric depending on $h$ but not its derivatives. Then:
\begin{thm}[Rendall \cite{RenRedu}]\label{thm:C2:Rend}Let $P$ and $V$ be as above, and let $(M,u,v)$ be the doubly-foliated manifold as above. Let $\phi_i$ be smooth sections of $V$ over $\mathcal H_i$ which agree on $S_{0,0}$. Suppose that with $\phi_i$ prescribed, $\mathcal H_i$ are null for $\psi(\phi_i)$. Then the equation $P(\phi) = 0$ with $\phi|_{\mathcal H_i} = \phi_i$ has a unique solution $\phi$ in a neighbourhood of $S_{0,0}$. Furthermore, any two solutions in Taylor series at $\mathcal H_1 \un \mathcal H_2$ are equal in Taylor series. \end{thm}
We will prove a version of \cref{thm:C2:Rend} for the special case of $P = Q$ and $V$ is the bundle of symmetric $(0,2)$ tensors because our proof of \cref{thm:C1:ChrVague} depends on a modified version, and it is instructive to see the easier case done first.

Once $h$ is found, it will be clear that $h$ is actually zero in Taylor series at $\mathcal H_1 \un \mathcal H_2$, and thus the wave gauge constraints are satisfied initially, and so by standard arguments are satisfied on the entire domain of existence. See \cref{thm:C2:wavegauge}. Since $h$ is $0$ in Taylor series, we may apply \cref{thm:C2:conversion} to change to a double-null gauge.

Uniqueness follows from similar arguments, but we will need to work slightly harder. See \cref{C:C2:Uniqueness}.

As is known (see for instance \cite{LukLoca}), the equations $\Ric(g)(L,L) = 0$ and $\Ric(g)(N,N) = 0$ provide a pair of obstructions or \emph{constraints} to solving for a metric even in Taylor series. Let $e^{\Phi}$ denote the conformal factor on $\mathcal H_1 \un \mathcal H_2$, i.e.\ $\slash{g} = e^{\Phi}\hat{\slash{g}}$. We will show below in \cref{C:C2:SolvingTaylor} that these equations take the form:
\begin{subequations}
\label{eq:C2:con}
\begin{align}
\label{eq:C2:con1}2N^2\Phi + (N\Phi)^2 + \Tr(\hat{\slash{g}}^{-1}\Lie_N\hat{\slash{g}})N\Phi + \frac{1}{2}|\Lie_N\hat{\slash{g}}|^2_{\hat{\slash{g}}} + N\Tr(\hat{\slash{g}}^{-1}\Lie_N\hat{\slash{g}})&= 0\\
\label{eq:C2:con2}
2L^2\Phi + (L\Phi)^2 + \Tr(\hat{\slash{g}}^{-1}\Lie_L\hat{\slash{g}})L\Phi + \frac{1}{2}|\Lie_L\hat{\slash{g}}|^2_{\hat{\slash{g}}} + L \Tr(\hat{\slash{g}}^{-1}\Lie_L\hat{\slash{g}})&= 0\end{align}
\end{subequations}
on $\mathcal H_1,\mathcal H_2$, respectively, with initial data $\Phi|_{S_{0,0}} = 0$ and $N\Phi|_{S_{0,0}}$ and $L\Phi|_{S_{0,0}}$ coming from the mean curvature specification. It is clear that these are nonlinear ODEs, and hence do not necessarily have solutions on all of $\mathcal H_1$ and $\mathcal H_2$, respectively. This is the \emph{only} obstruction to solving in Taylor series, however. Thus we will in fact prove:
\begin{prop}\label{thm:C2:TSbetter}Suppose that the constraints are satisfied on $\mathcal H_1$ and $\mathcal H_2$, i.e.\ \eqref{eq:C2:con1} and \eqref{eq:C2:con2} have solutions on all of $\mathcal H_1$ and $\mathcal H_2$, respectively. Then the conclusion of \cref{thm:C2:TSI} holds on all of $M([0,a)\times[0,b))$.\end{prop}

\section{Ricci curvature in a double-null gauge}
\label{C:C2:RicciDNG}
Let $(M,u,v)$ be a doubly-foliated manifold, and $(M,g,u,v,N)$ a double-null gauge. Let $\Omega, \slash{g},L,\zeta$ be the quantities associated to the double-null foliation $(M,g,u,v)$ as defined in \cref{C:C2:DNG}.

For symmetric tensors $\theta,\psi \in \Sym^2(T^\ast S)$, we will denote
\[(\theta \times \psi)_{ab} = \theta_{ac}\slash{g}^{cd}\psi_{db},\] where the metric $\slash{g}$ is the one induced by $g$. Let us also denote by $\slash{\tr}\theta$ the extrinsic trace $\slash{\tr}\theta = \slash{g}^{ab}\theta_{ba}$. We also define the divergence
\[(\slash{\div}\theta)_a = \slash{g}^{bc}\slash{\nabla}_{b}\theta_{ca},\]
 the $\slash{g}$-norm of a tensor by $|\bullet|$ and the inner product
\[\slash{g}(\theta,\psi) = \slash{g}^{ac}\slash{g}^{bd}\theta_{ab}\psi_{cd}.\]

Denote by $\slash{d}$ the exterior derivative mapping functions on $TS$ to fibre one forms,
\[(\slash{d}f)_a = \pa_{\theta^a}f.\]
The Ricci curvature in a double null gauge naturally splits into $6$ different quantities, comprising $\frac{1}{2}(\dim M)(\dim M+1)$ components: a symmetric tensor, a section of $\Sym^2(T^\ast S\otimes T^\ast S)$ over $M$, two fibre $1$-forms, and 3 smooth functions. Let us denote by $-$ an argument taken by a vector in $TS$. If $X \in TS$, let $Z = \slash{g}([N,L],\bullet)$ denote the $\slash{g}$ dual-vector field. The components of the Ricci curvature are:
\begin{subequations}
\label{eq:C2:Ric}
\begin{align}
\begin{split}
\label{eq:C2:one}2\Omega^2\Ric_{-,-} &= \Lie_N\Lie_L\slash{g} - \frac{1}{2}\Lie_{[N,L]} \slash{g} + \frac{1}{4}(\slash{\tr}\Lie_N \slash{g})\Lie_L\slash{g} + \frac{1}{4}(\slash{\tr}\Lie_L\slash{g})\Lie_N \slash{g}\\
&-\frac{1}{2}(\Lie_L\slash{g} \times \Lie_N \slash{g} + \Lie_N\slash{g}\times \Lie_L\slash{g})\\
&+2\Omega^2\slash{\Ric} - 2\slash{\Hess}(\Omega^2) + 2\Omega^2(\slash{d}\log \Omega)\otimes(\slash{d}\log\Omega)\\
&- \frac{1}{4\Omega^2}g([N,L],\bullet)\otimes g([N,L],\bullet)\end{split}\\
\begin{split}
\label{eq:C2:two}4\Omega^2\Ric_{N,-} &= \Lie_N Z + 2\Omega^2\slash{\divg}(\Lie_N \slash{g})-2\Omega^2\slash{d}\slash{\tr}\Lie_N\slash{g}\\
&+ \frac{1}{2}(\slash{\tr}\Lie_N\slash{g})Z + \slash{\tr}\Lie_N\slash{g}\slash{d}\Omega^2 - 4\Omega^2\slash{d}N\log\Omega - 2(N\log\Omega)Z
\end{split}\\
\begin{split}
\label{eq:C2:three}-4\Omega^2\Ric_{L,-} &= \Lie_L Z-2\Omega^2\slash{\divg}(\Lie_L\slash{g})+2\Omega^2\slash{d}\slash{\tr}\Lie_L\slash{g}\\
&+\frac{1}{2}(\slash{\tr}\Lie_L\slash{g})Z- \slash{\tr}\Lie_L\slash{g}\slash{d}\Omega^2 + 4\Omega^2\slash{d}L\log\Omega - 2(L\log\Omega)Z
\end{split}\\
\begin{split}
\label{eq:C2:four}-2\Ric_{NN} &= N\slash{\tr}\Lie_N \slash{g} + \frac{1}{2}|\Lie_N \slash{g}|^2 - 2(N\log\Omega)\slash{\tr}\Lie_N\slash{g}\end{split}\\
\begin{split}
\label{eq:C2:five}-2\Ric_{LL} &= L\slash{\tr}\Lie_L\slash{g} + \frac{1}{2}|\Lie_L\slash{g}|^2 - 2(L\log\Omega)\slash{\tr}\Lie_L\slash{g}\end{split}\\
\begin{split}
\label{eq:C2:six}-\frac{1}{2}\Ric_{NL} &=NL\log\Omega + \frac{1}{4}N\slash{\tr}\Lie_L\slash{g} - \frac{1}{8}\slash{\tr}\Lie_{[N,L]} \slash{g} - \frac{1}{2}[N,L]\log\Omega- \frac{1}{2}\slash{\Delta}\Omega^2\\
&+\frac{1}{8}\slash{g}(\Lie_L\slash{g},\Lie_N \slash{g}) + \frac{1}{8\Omega^2}|[N,L]|^2.
\end{split}
\end{align}
\end{subequations}
We will derive these formulae in \cref{C:A1:ComputeRicci} (see also \S 3.1.6 of \cite{KlaNicEvel} for the same equations, albeit written very differently).

It is important to observe that, given $N$, the Ricci curvature may be computed entirely from $L,\slash{g},\Omega$, which we treat as unknowns. Thus, by \cref{thm:C2:Reconstruct}, if we can find $L,\slash{g},\Omega$ with each component of $\Ric(g)$ being equal to $0$ in Taylor series, then the associated metric $g$ will have $\Ric(g) = 0$ in Taylor series. This is the approach we take.

The system \eqref{eq:C2:one}-\eqref{eq:C2:six} appears overdetermined. $L$ has only $\dim M -2$ components, since $2$ are fixed by $Lu = 0, Lv = 1$. Thus, counting components, these are really a system of $\frac{1}{2}(\dim M)(\dim M+1)$ equations in $\frac{1}{2}(\dim M)(\dim M - 1)$ unknowns. However, the contracted Bianchi identities,
\[g^{\beta\gamma}\nabla_\beta R_{\gamma\alpha} = \frac{1}{2}\nabla_\alpha R\]
(where $\nabla$ is the $g$-covariant derivative, $R_{\gamma\alpha}$ is the Ricci curvature, and $R = g^{\alpha\beta}R_{\alpha\beta}$ is the scalar curvature), implies certain integrability conditions on the system. Denoting $\chi^N,\chi^L$ the second fundamental forms of $g$ in the directions $N,L$, respectively, and $H^N,H^L$ the mean curvature (i.e.\ their traces), these conditions are as follows:
\begin{subequations}
\begin{align}
\label{eq:C2:onei}
\begin{split}
\Lie_L\Ric_{N,-}+\Lie_N \Ric_{L,-} &+H^L\Ric_{N,-} + H^N\Ric_{L,-} = 2\Omega^2\slash{\div}\Ric_{-,-}-\Omega^2\slash{d}\slash{\tr}\Ric_{-,-}\\
&\hspace{4.7cm}+2\Ric_{-,-}\cdot \slash{\grad}\Omega^2+\slash{d}\Ric_{NL}.\end{split}
\\
\label{eq:C2:twoi}
\begin{split}
L\Ric_{NN} +H^L\Ric_{NN}
&=-\Omega^2N\slash{\tr}\Ric_{-,-}-2\Omega^2\slash{g}(\chi^N, \Ric_{-,-})-H^N\Ric_{NL}\\
&+2\Omega^2\slash{\div}\Ric_{N,-}+2\Ric_{N,-}\cdot (\slash{\grad}\Omega^2) - \Ric_{N,-}\cdot [N,L].
\end{split}
\\
\label{eq:C2:threei}
\begin{split}
N\Ric_{LL} + H^N\Ric_{LL} &=-\Omega^2L\slash{\tr}\Ric_{-,-}-2\Omega^2\slash{g}(\chi^L, \Ric_{-,-})-H^L\Ric_{LN}\\
&+2\Omega^2\slash{\div}\Ric_{L,-}+2\Ric_{L,-}\cdot (\slash{\grad}\Omega^2) + \Ric_{L,-}\cdot [N,L].
\end{split}
\end{align}
\end{subequations}
Here, a $\cdot$ denotes a single contraction. We will derive these in \cref{C:A1:computebianchi}.

The first equation shows that only $\Ric_{-,-}$, $\Ric_{NL}$ and one of either $\Ric_{N,-}$ or $\Ric_{L,-}$ are sufficient to know the other, provided that the other is known on the corresponding initial surface. The second and third show that $\Ric_{-,-}$, $\Ric_{NL}$, $\Ric_{N,-}$, $\Ric_{L,-}$ are sufficient to know $\Ric_{NN}$ and $\Ric_{LL}$ everywhere, provided that each is known on the appropriate initial surface. This brings down the number of equations required to be solved everywhere to $\frac{1}{2}(\dim M)(\dim M - 1)$, which means the system is formally determined.

\section{Solving in Taylor series}
\label{C:C2:SolvingTaylor}
In this section, we prove \cref{thm:C2:TSI} and \cref{thm:C2:TSbetter}. This is at its heart a fairly standard Taylor series computation. We will solve to top order, and then solve iteratively for the higher order terms, which involves solving only the linearized equations, which we will compute. Let $\mathring{L}$ be the vector field satisfying $[N,\mathring{L}] = 0$ with $\mathring{L} = L$ on $\mathcal H_2$. Then $\mathring{L}u = 0, \mathring{L}v = 1$, and so $L = \mathring{L} + \overline{L}$, for some $\overline{L} \in C^\infty(M;TS)$. Observe that $[N,L] = [N,\overline{L}]$. Since $N$ is fixed, by \cref{thm:C2:Reconstruct}, it suffices to solve for $\slash{g},\overline{L},\Omega$. It will actually be more convenient to use $\omega = \log\Omega$, so we consider $\phi = (\slash{g},\overline{L},\omega)$, a section of $\Sym^2 T^\ast S \oplus TS\oplus \R$. We need to find such a $\phi$ which solves \eqref{eq:C2:one}--\eqref{eq:C2:six} in Taylor series, since then we may find the solution $g$ in Taylor series using \cref{thm:C2:Reconstruct}. For the uniqueness of the series, it suffices to establish the uniqueness of $\phi$ since every solution $g$ of $\Ric(g) = 0$ gives a solution $\phi$ of \eqref{eq:C2:one}--\eqref{eq:C2:six} and vice-versa.

Let us use $N,\mathring{L}$ to fix a product decomposition of $M$, turning it into $[0,a) \times [0,b) \times S^2$. We will thus conflate $(u,v)$ with coordinates on $[0,a)\times[0,b)$, and $N$ with $\pa_u$, $\mathring{L}$ with $\pa_v$. 

We use the notation $O(u^j)$ and $O(v^j)$ to denote a quantity equal to $u^j$ and $v^j$, respectively, times a smooth section of an appropriate bundle. The notation $O(u^\infty)$ or $O(v^\infty)$ indicates that for all $j$ the quantity is in $O(u^j)$ or $O(v^j)$, respectively.

Let us first deal with the overdetermined nature of the equations:
\begin{lem}\label{thm:C2:integrability}Suppose that $\phi$ is given and \eqref{eq:C2:one}, \eqref{eq:C2:three}, \eqref{eq:C2:six} are $0$ in Taylor series, and \eqref{eq:C2:two}, \eqref{eq:C2:four} are identically $0$ on $\mathcal H_1$ and \eqref{eq:C2:five} is identically $0$ on $\mathcal H_2$. Then all equations \eqref{eq:C2:one}-\eqref{eq:C2:six} are $0$ in Taylor series.\end{lem}
\begin{proof}
Since $\Omega = \exp(\omega) > 0$, we may divide by it. From \eqref{eq:C2:onei}, $\Ric_{N,-}$ solves an equation of the form
\begin{equation}\label{eq:integrability:integral}\Lie_{\mathring{L}}\Ric_{N,-} + (\Lie_{\overline{L}}+H^L)\Ric_{N,-}= f,\end{equation} and $f$ is $0$ in Taylor series at $\mathcal H_1 \un \mathcal H_2$, and $\Ric_{N,-} = 0$ on $\mathcal H_1$. Since $\overline{L}$ is a fibre vector field, $\Lie_{\overline{L}}\Ric_{N,-} = 0$ on $\mathcal H_1$, and thus the equation above reduces to $\Lie_{\mathring{L}}\Ric_{N,-} = 0$ on $\mathcal H_1$. This shows that $\Ric_{N,-} = 0$ to order one. To derive it is $0$ to order $k > 1$, apply $\Lie_{\mathring{L}}^{k-1}$ to the previous equation to obtain
\[\Lie_{\mathring{L}}^k \Ric_{N,-} + (\Lie_{\overline{L}}+H^L)\Lie_{\mathring{L}}^{k-1}\Ric_{N,-} = \Lie_{\mathring{L}}^{k-1}f + \cdots,\]
where the $\cdots$ involve only products of known quantities with fibre-derivatives of $\Lie_{\mathring{L}}^j\Ric_{N,-}$ for $j < k$. Inductively, we may show that $\Lie_{\mathring{L}}^k \Ric_{N,-} = 0$ on $\mathcal H_1$ for all $k$, and thus $\Ric_{N,-} = 0$ in Taylor series there. Now \eqref{eq:integrability:integral} also implies that $\Ric_{N,-} = 0$ on $\mathcal H_2$, since here it solves a homogeneous linear transport equation with trivial data (recall $\overline{L} \equiv 0$ on $\mathcal H_2$). We may apply $\Lie_N^{k}$ to the equation to obtain
\begin{equation}\label{eq:C2:transport}\Lie_{\mathring{L}}\Lie_N^k \Ric_{N,-} + (\Lie_{\overline{L}}+H^L)\Lie_N^k \Ric_{N,-} = \Lie_N^k f + \cdots,\end{equation}
where the $\cdots$ involve only products of known quantities with fibre-derivatives of $\Lie_N^j \Ric_{N,-}$, for $j < k$. Equation \eqref{eq:C2:transport} is a transport equation for $\Lie_N^k \Ric_{N,-}$. Since $\Lie_N^k \Ric_{N,-} = 0$ on $\mathcal H_1$ by assumption, the equation implies that the same is true on $\mathcal H_2$, we well.

With $\Ric_{N,-} = 0$ in Taylor series, the same argument applied to \eqref{eq:C2:twoi} and \eqref{eq:C2:threei} works for the other two quantities.\end{proof}

We now begin the process of solving in Taylor series. We will keep track of two different Taylor series, one a series in $v$ and at $\mathcal H_1$, and one a series in $u$ and at $\mathcal H_2$. The functions and terms involving the series at $\mathcal H_1$ will be indicated with a superscript $1$, while those at $\mathcal H_2$ will be indicated with a superscript $2$.

Let us first solve to top order:
\begin{lem}\label{thm:C2:toporder}There exists $0 < u' \leq a$ and $0 < v' \leq b$ and a unique choice of sections $\phi^1_0$, $\phi^1_1$ defined on $\mathcal H_1\n\{u \leq u'\}$ and $\phi^2_0$, $\phi^2_1$ defined on $\mathcal H_2\n\{v \leq v'\}$, satisfying $\phi^1_0|_{S_{0,0}} = \phi^2_0|_{S_{0,0}}$, $\pa_u \phi^1_0|_{S_{0,0}} = \phi^2_1|_{S_{0,0}}$, $\phi^1_1|_{S_{0,0}} = \pa_v \phi^2_0|_{S_{0,0}}$ and such that any $\phi$ smooth with $\pa^i_v\phi|_{\mathcal H_1} = \phi^1_i$ ($i=0,1$) and $\pa^i_u \phi|_{\mathcal H_2} = \phi^2_i$ ($i=0,1$) satisfies \eqref{eq:C2:one}, \eqref{eq:C2:three}, \eqref{eq:C2:six} $=0$ on $\mathcal H_1 \un \mathcal H_2$, \eqref{eq:C2:two},\eqref{eq:C2:four} $=0$ on $\mathcal H_1$, \eqref{eq:C2:five} $=0$ on $\mathcal H_2$, and has the correct initial data, i.e.\ $\slash{g}$ belongs to the correct conformal class along $\mathcal H_1 \un \mathcal H_2$, $\slash{g}$ has the correct value at $S_{0,0}$ and the trace of its derivatives correspond to the correct mean curvatures, $\overline{L} = 0$ on $\mathcal H_2$, $\omega=0$ on $\mathcal H_1 \un \mathcal H_2$, and $\zeta$ (or equivalently $[N,L]$) has the correct value at $S_{0,0}$. If the constraints are satisfied, we may choose $u' = a$, $v' = b$.\end{lem}
\begin{proof}Let $\phi$ satisfy the equations. We first solve for the conformal class of the metric $\slash{g}$ over $\mathcal H_1 \un \mathcal H_2$. Let us set $\slash{g} = e^\Phi\hat{\slash{g}}$ for some conformal factor $e^\Phi$ to be solved for. Since we are assuming $\omega = 0$ on $\mathcal H_1 \un \mathcal H_2$ by assumption, then \eqref{eq:C2:four}, \eqref{eq:C2:five} $=0$ on $\mathcal H_1$ and $\mathcal H_2$, respectively, are equivalent to the following ODEs along the flow of $N,L$ respectively:
\begin{align*}
2N^2\Phi + (N\Phi)^2 + \Tr(\hat{\slash{g}}^{-1}\Lie_N\hat{\slash{g}})N\Phi + \frac{1}{2}|\Lie_N\hat{\slash{g}}|^2_{\hat{\slash{g}}} + N\Tr(\hat{\slash{g}}^{-1}\Lie_N\hat{\slash{g}})&= 0\\
2L^2\Phi + (L\Phi)^2 + \Tr(\hat{\slash{g}}^{-1}\Lie_L\hat{\slash{g}})L\Phi + \frac{1}{2}|\Lie_L\hat{\slash{g}}|^2_{\hat{\slash{g}}} + L \Tr(\hat{\slash{g}}^{-1}\Lie_L\hat{\slash{g}})&= 0,\end{align*}
on $\mathcal H_1,\mathcal H_2$, respectively. By assumption, $\Phi|_{S_{0,0}}= 0$, and $N\Phi|_{S_{0,0}}$ and $L\Phi|_{S_{0,0}}$ are determined by the specification of the mean curvature on $S_{0,0}$:
\begin{align*}
f_1 &= \frac{1}{2}\slash{\tr}_{\slash{g}}\Lie_N \slash{g}|_{S_{0,0}} = N\Phi|_{S_{0,0}} + \frac{1}{2}\slash{\tr}_{\hat{\slash{g}}}\Lie_N \hat{\slash{g}}|_{S_{0,0}}\\
f_2 &= \frac{1}{2}\slash{\tr}_{\slash{g}}\Lie_L \slash{g}|_{S_{0,0}} = L\Phi|_{S_{0,0}} + \frac{1}{2}\slash{\tr}_{\hat{\slash{g}}}\Lie_L \hat{\slash{g}}|_{S_{0,0}}.
\end{align*}
Since the fibres are compact, there is some finite time along the flow of $N$ and $L$ on which $\Phi$ does not blow up, i.e.\ there is some $u' \leq a$ and $v' \leq b$ such that solutions exist for $0 \leq u < u'$ and $0 \leq v < v'$ on $\mathcal H_1$ and $\mathcal H_2$, respectively. If we are operating under the hypotheses of \cref{thm:C2:TSbetter}, then we have by assumption $u' = a$ and $v' = b$. This gives us the top order behaviour of $\slash{g}$.
Now let us deal with \eqref{eq:C2:two} $=0$ along $\mathcal H_1$. Here, we know all quantities except for $Z$ and $\Lie_N Z$. Thus it gives an linear transport equation for $Z$, whose initial value we know in terms of $\zeta$ (or equivalently we know $Z$ initially). Knowing $Z$ and $\slash{g}$, we know $[N,L] = [N,\overline{L}]$, and thus $\overline{L}$ and $L$ along $\mathcal H_1$. This gives us the zeroth order behaviour of $\overline{L}$ at $\mathcal H_1$. Similarly, by setting \eqref{eq:C2:three} $=0$ we may find $[N,L]$ along $\mathcal H_2$, which gives us the first-order behaviour of $\overline{L}$ ($\overline{L} = 0$ initially on $\mathcal H_2$ is by definition).
Now \eqref{eq:C2:one} $=0$ can be written as a linear transport equation for $\Lie_{\mathring{L}} \slash{g}$ along $\mathcal H_1$ in terms of known quantities and, commuting the derivatives, $\Lie_N \slash{g}$ along $\mathcal H_2$ in terms of known quantities, so we may solve it, obtaining the first-order behaviour of $\slash{g}$. Recalling that $\omega = 0$ on $\mathcal H_1 \un \mathcal H_2$, we do something similar with \eqref{eq:C2:six} to get $N\omega$ on $\mathcal H_2$ and $\mathring{L}\omega$ on $\mathcal H_1$, which gives us the first-order behaviour. Finally, \eqref{eq:C2:three} $ =0$ on $\mathcal H_1$ takes the form of an \emph{algebraic} equation for $\Lie_L Z$ on $\mathcal H_1$ in terms of known quantities. By definition,
\begin{align*}\Lie_L Z &= \slash{g}([\mathring{L},[N,\overline{L}],\bullet) + (\Lie_{\mathring{L}}\slash{g})([N,\overline{L}],\cdot)\\
&+ \slash{g}([\overline{L},[N,\overline{L}],\overline{L}],\bullet) + (\Lie_{\overline{L}}\slash{g})([N,\overline{L}],\bullet).\end{align*}
Since $\overline{L}$ is tangent to the fibres, we know all quantities except the first. Thus we also know
\[[\mathring{L},[N,\overline{L}]] = [N,[\mathring{L},\overline{L}]] + [[\mathring{L},N],\overline{L}] = [N,[\mathring{L},\overline{L}]] .\]
Letting $K$ denote a known-quantity, we thus have
\[\Lie_N [\mathring{L},\overline{L}] = K.\] We know the initial data at $S_{0,0}$ is $0$ (since $L = \mathring{L}$ there). This is a transport equation for $[\mathring{L},\overline{L}]$ which we then solve to find the first-order behaviour.

Thus we know what $\phi$ and its $N,\mathring{L}$ derivatives must be, and conversely if $\phi$ and its first derivatives are as prescribed, then $\phi$ satisfies the equations. This proves the lemma.
\end{proof} 

Now we may iteratively now solve for the remainder of the Taylor series in $u$ and $v$ at the same time. Suppose we are given for $j \geq 2$ and $0 \leq \ell \leq j$ and $i = 1,2$ sections $\phi_{\ell}^i$ defined over $\mathcal H_i\n C_i$ (where $C_1 = \{u \leq u'\}$ and $C_2 = \{v \leq v'\}$), and $\phi_0^i, \phi_1^i$ agree with the sections found in \cref{thm:C2:toporder}. Extend each $\phi_\ell^1$ to be constant in $v$ and each $\phi_\ell^2$ to be constant in $u$. Write
\[\phi^1_{(j)} = \phi_0^1 + u\phi_1^1 + \cdots + u^j\phi_1^j,\] and similarly for $\phi^2_{(j)}$. Assume that $\phi^1_{(j)}$ solves \eqref{eq:C2:one}, \eqref{eq:C2:three}, \eqref{eq:C2:six} $=0$ mod $O(u^{j})$, and similarly for $\phi^2_{(j)}$, and that for all $\ell$, \[\phi^1_{\ell}|_{S_{0,0}} = \frac{1}{\ell!}\pa_v^{\ell}\phi_0^2|_{S_{0,0}}, \ \phi^2_{\ell}|_{S_{0,0}} = \frac{1}{\ell!}\pa_u^{\ell}\phi_0^2|_{S_{0,0}}\]
(observe that this is true for $\ell = 0,1$ by assumption). We seek to find $\phi^i_{j+1}$ for $i = 1,2$.

\begin{lem}\label{thm:C2:higherorder}There exists a unique section $\phi^1_{j+1}$ defined on $\mathcal H_1 \n C_1$ and $\phi^2_{j+1}$ defined on $\mathcal H_2 \n C_2$ with $\phi^1_{j+1}|_{S_{0,0}} = \frac{1}{(j+1)!}\pa_v^{j+1}\phi^2_0|_{S_{0,0}}$ and $\phi^2_{j+1}|_{S_{0,0}} = \frac{1}{(j+1)!}\pa_u^{j+1}\phi^1_0|_{S_{0,0}}$ such that $\phi^1_{(j)} + v^{j+1}\phi^1_{j+1}$ solves the equations mod $O(v^{j+1})$ and $\phi^2_{(j)}+u^{j+1}\phi^2_{j+1}$ solves the equations mod $O(u^{j+1})$.\end{lem}
\begin{proof}Let us first show how to solve for $\phi^1_{j+1}$. Let $P = (P^1,P^2,P^3)$ denote the operator taking in a section $\phi$ and returning the right-hand side of \eqref{eq:C2:one}, \eqref{eq:C2:three}, \eqref{eq:C2:six}, respectively. We may write
\[P(\phi^1_{(j)} + v^{j+1}\phi^1_{j+1}) = P(\phi^1_{(j)}) + v^j T^1_{j+1,\phi^1_{(j)}}(\phi^1_{j+1}) + v^{j+1}M_{j+1,\phi^1_{(j)}}(\phi^1_{j+1}),\] where $T^1_{j+1,\phi^1_{(j)}}$ is linear, $M^1_{j+1,\phi^1_{(j)}}$ is nonlinear, and both depend on $j$ and $\phi^1_{(j)}$. By assumption, $P(\phi^1_{(j)}) \in O(v^{j})$, and so this sets up a linear equation for $\phi^1_{j+1}$ 
\[T^1_{j+1,\phi^1_{(0)}}(\phi^1_{j+1}) = -(v^{-j-1}P(\phi^1_{(j)}))|_{v=0},\] (since we ignore terms of higher order, and the dependence of $T^1_{j+1,\phi^1_{(0)}}$ to top order depends only on $j+1$ and $\phi^1_{(1)}$). We will show that $T_{j+1,\phi^1_{(0)}}$ takes the form (essentially) of a transport operator along $\mathcal H_1$, so will have a unique solution with the give initial data. The same argument will allow us to solve for $\phi^2_{j+1}$ (with a different, but analogous, linear operator $T^2_{j+1,\phi^2_{(0)}}$, instead).

All that's left to do is compute the $T$'s and verify that we can reduce the linear PDE above to transport equations. Fix some arbitrary $\phi = (\slash{g},\overline{L},\omega)$ smooth (with $\Omega := e^\omega$) and $\psi = (\slash{g}_{j+1},\overline{L}_{j+1},\omega_{j+1})$ not depending on $v$, and we will compute $T^1_{j+1,\phi}(\psi)$. We have:
\begin{subequations}
\begin{align}
\label{eq:C2:TS1}
\begin{split}
2\Omega^2\Ric_{-,-}&(\phi+v^{j+1}\psi) = 2\Omega^2\Ric_{-,-}(\phi)\\
&+ (j+1)v^{j}\left(\pa_u \slash{g}_{j+1} + \frac{1}{4}\slash{\tr}(\Lie_N \slash{g})\slash{g}_{j+1} + \frac{1}{4}\slash{\tr}(\slash{g}_{j+1})\Lie_N \slash{g}\right.\\
&\left. -\frac{1}{2}\left(\slash{g}_{j+1}\times \Lie_N \slash{g} + \Lie_N \slash{g} \times \slash{g}_{j+1}\right)\right) + O(v^{j+1})\end{split}\\
\label{eq:C2:TS2}
\begin{split}
-4\Omega^2\Ric_{L,-}&(\phi+v^{j+1}\psi) = -4\Omega^2\Ric_{L,-}(\phi)\\
&+ (j+1)v^{j}\left(\slash{g}(\pa_u \overline{L}_{j+1},\bullet)+\slash{g}_{j+1}([N,\overline{L}],\bullet) -2\Omega^2\slash{\div}(\slash{g}_{j+1}) \right.\\
&\left.+ 2\Omega^2\slash{d}\slash{\tr} \slash{g}_{j+1}+\frac{1}{2}(\slash{\tr}\slash{g}_{j+1})\slash{g}([N,\overline{L}],\bullet) - (\slash{\tr}\slash{g}_{j+1})\slash{d}\Omega^2\right.\\
&\left. + 4\Omega^2\slash{d}\omega_{j+1} - 2\omega_{j+1}\slash{g}([N,\overline{L}],\bullet)\right) + O(v^{j+1})\end{split}\\
\label{eq:C2:TS3}
\begin{split}
-\frac{1}{2}\Ric_{N,L}&(\phi+v^{j+1}\psi) = -\frac{1}{2}\Ric_{N,L}(\phi)\\
&+ (j+1)v^{j}\left(\pa_u \omega_{j+1} + \frac{1}{4}\pa_u \slash{\tr}\slash{g}_{j+1}\right.\\
& \left. +\frac{1}{8}\slash{g}(\slash{g}_{j+1},\Lie_N \slash{g})\right) + O(v^{j+1}).
\end{split}
\end{align}\end{subequations}
All slash-operations are performed with respect to $\slash{g}$.

$T^1_{j+1,\phi}$ is thus $j+1$ times the terms in brackets. Setting this equal to the appropriate right-hand side is essentially a transport equation, since the first component of $T^1_{j+1,\phi}$ is certainly a transport equation for $\slash{g}_{j+1}$, and then with this known, the last component is a transport equation for $\omega_{j+1}$, and then with both these known, the second component is, after raising via $\slash{g}^{-1}$, is a transport equation for $\overline{L}_{j+1}$.

Now let $\psi$ not depend on $u$, and let us compute the linearization at $\mathcal H_2$. Assume now that the $\overline{L}$ of $\phi$ is $0$ at $\mathcal H_2$ (which it is for the case we care about, since by our assumption $\phi$ has the correct initial data on $\mathcal H_2$, and so the associated $L = \mathring{L}+\overline{L}$ is just $\mathring{L}$). We will compute $T^2_{j+1,\phi}(\psi)$. We have:
\begin{subequations}
\begin{align}
\begin{split}
2\Omega^2\Ric_{-,-}&(\phi+u^{j+1}\psi) = 2\Omega^2\Ric_{-,-}(\phi)\\
&+ (j+1)u^j\left(\pa_v \slash{g}_{j+1} + \frac{1}{2}\Lie_{\overline{L}_{j+1}}\slash{g} + \frac{1}{4}(\slash{\tr}\pa_v \slash{g})\slash{g}_{j+1}\right.\\
&\left. + \frac{1}{4}(\slash{\tr}\slash{g}_{j+1})\pa_v \slash{g}- \frac{1}{2}\left(\slash{g}_{j+1}\times \pa_v\slash{g} + \pa_v \slash{g} \times \slash{g}_{j+1}\right)\right.\\
&\left. - \frac{1}{4\Omega^2}\left(\slash{g}(\overline{L}_{j+1},\cdot)\otimes \slash{g}([N,\overline{L}],\bullet) + \slash{g}([N,\overline{L}],\bullet) \otimes (\slash{g}(\overline{L}_{j+1},\bullet)\right)\right)\\
&+ O(u^{j+1})\end{split}\\
\begin{split}
-4\Omega^2\Ric_{L,-}&(\phi+u^{j+1}\psi) = -4\Omega^2\Ric_{L,-}(\phi)\\
&+ (j+1)u^{j}\left(\slash{g}(\pa_v\overline{L}_{j+1},\bullet)\right) + O(u^{j+1})\end{split}\\
\begin{split}
-\frac{1}{2}\Ric_{N,L}&(\phi+u^{j+1}\psi) = -\frac{1}{2}\Ric_{N,L}(\phi)\\
&+ (j+1)u^{j}\left(\pa_v \omega_{j+1} + \frac{1}{4}\pa_v\slash{\tr}\slash{g}_{j+1} + \frac{1}{8}\slash{\tr}\Lie_{\overline{L}_{j+1}}\slash{g}\right.\\
&\left. + \frac{1}{2}\overline{L}_{j+1}\omega + \frac{1}{8}\slash{g}(\pa_v\slash{g},\slash{g}_{j+1}) +\frac{1}{4\Omega^2}\slash{g}(\overline{L}_{j+1},[N,\ovelrine{L}])\right) + O(u^{j+1}).
\end{split}
\end{align}
\end{subequations}$T^2_{j+1,\phi}$ is thus $j+1$ times the terms in brackets. Setting this equal to the appropriate right-hand side, the second component is, after raising via $\slash{g}^{-1}$, a transport equation along the flow of $\mathring{L}=\pa_u$ on $\mathcal H_2$, for $\overline{L}_{j+1}$, in terms of known quantities. Thus the first component is a transport equation for $\slash{g}_{j+1}$ involving known quantities, and the third component becomes a transport equation for $\omega_{j+1}$ in terms of known quantities.
\end{proof}

From \cref{thm:C2:toporder}, \cref{thm:C2:higherorder}, and Borel's lemma, we may find $\phi^1$, $\phi^2$ solving \eqref{eq:C2:one}, \eqref{eq:C2:three}, \eqref{eq:C2:six}$=0$ in Taylor series at $\mathcal H_1$ and $\mathcal H_2$, respectively, and which are unique in Taylor series. We wish to find a $\phi$ whose Taylor series at $\mathcal H_1$ agrees with that of $\phi^1$ and whose Taylor series at $\mathcal H_2$ agrees with that of $\phi^2$. To do this, it suffices by Borel's lemma for $\phi^1,\phi^2$ to satisfy the compatibility condition $\pa_u^i \pa_v^j \phi^1|_{S_{0,0}} = \pa_u^i \pa_v^j \phi^2|_{S_{0,0}}$ for all $i,j$.
\begin{lem}The compatibility condition holds.\end{lem}
\begin{proof}By assumption, \[\pa_u^i \phi^1|_{S_{0,0}}= \pa_u^i \phi^1_0|_{S_{0,0}} = i!\phi^2_i|_{S_{0,0}} = \pa_u^i \phi^2|_{S_{0,0}},\] and similarly $\pa_u^i \phi^1|_{S_{0,0}} = \pa_v^i \phi^2|_{S_{0,0}}$. Since $\phi^1,\phi^2$ satisfy \eqref{eq:C2:one}, \eqref{eq:C2:three}, \eqref{eq:C2:six} $=0$ in Taylor series at $\mathcal H_1$, $\mathcal H_2$, respectively, there is a non-linear differential operator, $F$, involving only derivatives tangent to the fibres $S$, and not depending on $i=1,2$ such that
\begin{align}
\label{eq:C2:eqqq}
\begin{split}
\pa_u\pa_v \phi^1 = F(\phi^1,\pa_u \phi^1,\pa_v\phi^1) + O(v^\infty)\\
\pa_u\pa_v \phi^2 = F(\phi^2,\pa_u \phi^2,\pa_v\phi^2) + O(u^\infty).\end{split}\end{align}
Indeed, for equations \eqref{eq:C2:one} and \eqref{eq:C2:six}, just extract the top order $\pa_u\pa_v$ and move everything else to the other side, and for \eqref{eq:C2:three}, first raise using $\slash{g}$ and do the same.
Starting from the fact that we know $\phi^1|_{S_{0,0}}= \phi^2|_{S_{0,0}}$, $\pa_u\phi^1|_{S_{0,0}} = \pa_u \phi^2|_{S_{0,0}}$, $\pa_v\phi^1|_{S_{0,0}} = \pa_v \phi^2|_{S_{0,0}}$, and using that by assumption $\pa_u^i \phi^1|_{S_{0,0}}= \pa_u^i \phi^2|_{S_{0,0}}$, $\pa_u^i \phi^1|_{S_{0,0}} = \pa_v^i \phi^2|_{S_{0,0}}$ for all $i \geq 2$, we may apply derivatives $\pa_u,\pa_v$ to \eqref{eq:C2:eqqq} to inductively show for all $\ell \geq 2$, $\pa_u^{i}\pa_v^j \phi^1|_{S_{0,0}} = \pa_u^{i}\pa_v^j\phi^2|_{S_{0,0}}$ whenever $i+j = \ell$.
\end{proof}

The section $\phi$ we have constructed satisfies the hypotheses of \cref{thm:C2:integrability}, since by definition it solves equation \eqref{eq:C2:one}, \eqref{eq:C2:three}, \eqref{eq:C2:six} in Taylor series, and solves the other equations identically on the correct hypersurfaces since it has the correct top-order behaviour. Thus by \cref{thm:C2:integrability}, we have found our solution $\phi$ in Taylor series. The solution $\phi$ is unique in Taylor series since any other solution has expansions at $\mathcal H_1$, and $\mathcal H_2$, and we already know that these are unique.

The only place in the proof where we replaced $a$ with $u' \leq a$ and $b$ with $v' \leq b$ was in \cref{thm:C2:toporder} when we solved the constraint equations. Thus we observe that \cref{thm:C2:TSbetter} has also been proven.

\section{Obtaining a solution in wave gauge}
\label{C:C2:gettingsolution}
Denote by $g$ the Taylor series solution provided by \cref{thm:C2:TSI}, which exists on $M([0,u')\times[0,v'))$. We wish to find a symmetric $(0,2)$ tensor $h$ such that $\Ric(g+h) = 0$ and $g+h$ is in wave gauge with respect to $g$. Let $F$ be the vector-field-valued first-order differential operator defined for tensors close to $0$ by
\[F^\alpha(h) = (g+h)^{\beta\gamma}(\tilde{\Gamma}^{\alpha}_{\beta\gamma}-\Gamma^\alpha_{\beta\gamma}),\] where $\tilde{\Gamma},\Gamma$ are the Christoffel symbols of $g+h$,$g$, respectively.\footnote{$F$ is a vector field because the difference of two connections is a tensor.} Observe that $F$ vanishes precisely when the identity map is a wave map if the domain is given the metric $g+h$ and the codomain is given the metric $g$.

Write $\nabla$ for the Levi-Civita connection of $g$. Introduce $G$, a second-order differential operator, taking values in symmetric tensors, defined by
\[G(h)_{\alpha\beta} = \frac{1}{2}((g+h)_{\alpha\lambda}\nabla_\beta F^{\lambda} + (g+h)_{\beta\lambda}\nabla_\alpha F^{\lambda}).\]
Then, as is well known,
\[Q(h) := \tilde{\Ric}(g+h) := \Ric(g+h) + G(h)\] is a quasilinear hyperbolic partial differential operator with leading part
\[-\frac{1}{2}(g+h)^{\alpha\beta}\nabla_\alpha \nabla_\beta h_{\mu\nu}.\] The equation $\tilde{\Ric}(g+h) = 0$ is called the \emph{reduced Einstein equation}.
In fact, $Q(h)$ has the form (see \cite{ChoGene}, \S VI.7.1-4)
\begin{equation}\label{eq:C2:hastheform}Q(h) = -\frac{1}{2}(g+h)^{\alpha\beta}\nabla_\alpha \nabla_\beta h_{\mu\nu} + P((g+h)^{-1})(\nabla h,\nabla h) - (g+h)^{-1}\cdot \Riem(g)\cdot(g+h),\end{equation}
where $P((g+h)^{-1})$ is a $(2,6)$ tensor on $M$ depending on $(g+h)^{-1}$, but not any of its derivatives, and which we interpret here as acting on a pair of $(0,3)$ tensors, and is symmetric interpreted this way, and the $\cdot$ in the last term denote specific contractions. More explicitly, for $H$ a type $(2,0)$ tensor, $K$ a type $(0,2)$ tensor and $\xi$ a type $(0,3)$ tensor 
\begin{align*}
P(H)(\xi,\xi)_{\alpha\beta} &= -\frac{1}{2}(\xi_{\beta\rho\sigma}H^{\rho\lambda}H^{\sigma\mu}\xi_{\lambda\alpha\mu} + \xi_{\alpha\rho\sigma}H^{\rho\lambda}H^{\sigma\mu}\xi_{\lambda\beta\mu})\\
&- \frac{1}{4}H^{\mu\rho}(\xi_{\alpha\rho\mu} + \xi_{\mu\alpha\rho} - \xi_{\rho\alpha\mu})H^{\lambda\sigma}(\xi_{\beta\sigma\lambda}+\xi_{\lambda\beta\sigma}-\xi_{\sigma\beta\lambda}).
\end{align*}
and
\[(H\cdot \Riem (g_0) \cdot K)_{\alpha\beta}= -\frac{1}{2}H^{\lambda\mu}(K_{\alpha\rho}R^{\ \rho}_{\lambda \ \beta\,\mu} + K_{\beta\rho}R^{\ \rho}_{\lambda \ \alpha\mu}),\]
where $R^{\ \rho}_{\lambda \ \alpha\mu}$ are the components of the curvature tensor associated to $g$.

We may appeal to \cref{thm:C2:Rend} directly to find the solution $h$ to $Q(h) = 0$ which vanishes in Taylor series at $\mathcal H_1$ and $\mathcal H_2$, but we will prove it in this special case, to convince the reader that the same argument can go through almost unchanged to prove \cref{thm:C1:ChrVague}. We will be quite pedantic in doing so. 

The idea of Rendall is to reduce to the Cauchy problem, so we will need to embed our setup into an appropriate setting for the Cauchy problem, i.e.\ a sliced Lorentzian manifold.

\begin{lem}\label{thm:C2:Extension}Fix $0 < u'' < u'$ and $0 < v'' < v'$. Then there exists an embedding of $M([0,u'']\times [0,v''])$ into $\tilde{M} = \R_t \times \R_y \times S_{0,0}$ such that $\mathcal H_1 \n M([0,u'']\times [0,v''])$ is a subset of $\{t = -y\}$, $\mathcal H_2 \n M([0,u'']\times [0,v''])$ is a subset of $\{t = y\}$, and there exists an extension of $g$ to $\tilde{M}$ for which the hyperplanes $\{t=\pm y\}$ are null, the hyperplanes $\{t = \text{const.}\}$ are spacelike, $\pa_t$ is timelike and future-oriented, and $g = -dt^2 + dy^2 + \mathring{\slash{g}}$ for $t$ or $y$ large enough (here $\mathring{\slash{g}}$ is any fixed Riemannian metric on $S_{0,0}$).\end{lem}
\begin{proof}For $p \in M([0,u'']\times [0,v''])$, defined its image in $\tilde{M}$ as follows. Define $t,y$ by $u = \frac{1}{2}(t-y)$, $v = \frac{1}{2}(t+y)$. The first two coordinates are $(t(p),y(p))$. The point in $S_{0,0}$ is obtained by flowing via $N$ back to $\mathcal H_2$, and then via $\mathring{L}$ to $S_{0,0}$.

It is clear that $N = \pa_t - \pa_y + X$, $L = \pa_t+\pa_y + Y$, where $X,\ Y$ are vector fields tangent to $S_{0,0}$. Extend $X,\ Y$ to be $0$ near infinity, which then extends $N,\ L$. Now extend $g$ by extending its action between $N,\ L$ and $TS_{0,0}$ as follows. Extend $g(N,L)$ to be negative always, and $-2$ near infinity, keep $g(N,N)=g(L,L) = g(N,\Theta) = 0$ (for $\Theta$ in $TS_{0,0}$ arbitrary) and extend $g|_{S_{0,0}}$ arbitrarily provided it is equal to $\mathring{\slash{g}}$ near infinity (we may extend this way since all quantities are actually defined up until $u = u'$ and $v= v'$, and we are allowed to change the quantities in the area region $\{u'' < u < u', \ v' < v < v''\}$).\end{proof}
 
Now we prove existence. The proof below is a special case of \cref{thm:C2:Rend} proved in \cite{RenRedu}, but proved in our language for purposes of unity and later modification.
\begin{prop}\label{thm:C2:existence}There exist $0 < u^\ast < u''$ and $0 < v^\ast < v''$ and a tensor $h$ solving $Q(h) = 0$ on $M([0,u^\ast)\times [0,v^\ast))$ with $h=0$ in Taylor series at $(\mathcal H_1\un \mathcal H_2) \n M([0,u^\ast)\times [0,v^\ast))$.\end{prop}
\begin{proof}
Embed $M([0,u'']\times [0,v''])$ into $\tilde{M}$ and extend $g$, as provided by \cref{thm:C2:Extension}. Let us set $P = Q-\Ric(g)$, so that $P(0) = 0$. We are interested in finding an $h$ with $P(h) = -\Ric(g)$. Notice that extending $g$ also extends $P$.

The idea of Rendall is to define a section $\rho_1$ by
\[\rho_1 = \begin{cases}-\Ric(g) & u,v \geq 0\\
0 & \text{otherwise}.\end{cases}\]
Notice that $\rho_1$ is smooth provided $u \leq u''$ and $v \leq v''$, since $-\Ric(g)$ is $0$ in Taylor series there. Fix $0 < u_3 < u''$ and $0 < v_3 < v''$, and let $\rho$ be any smooth section with $\rho = \rho_1$ for $u \leq u_3$, $v \leq v_3$, and $\rho$ compactly supported. Then $\rho$ is a compactly supported smooth section.

Now we may solve the Cauchy problem $P(h) = \rho$ with Cauchy data \[(h|_{t=0},\pa_t h|_{t=0})= (0,0),\] at least for time $0 \leq t < \epsilon$. Shrinking $\epsilon$, we may assume that if $t < \epsilon$ and $u,v \geq 0$ then $u < u_3, v < v_3$. Let $u^\ast,v^\ast$ be largest values such that if $0 \leq u < u^\ast$ and $0 \leq v < v^\ast$, then $t < \epsilon$ . Let us now show that $h$ is $0$ in Taylor series at the required hypersurfaces. Set
\[R = \{u < u^\ast, v < v^\ast, u < 0 \text{ or } v < 0, t \geq 0\}.\] The boundary of $R$ consists of $(\mathcal H_1\un \mathcal H_2) \n M([0,u^\ast)\times [0,v^\ast))$ and $\{t = 0\}$, the former of which is null for $g$, and the latter of which is spacelike. It follows that the causal past of $p \in R$ according to $g$ does not exit $R$. In particular, $P(0)=0 = \rho$ on the causal past of $p$, and has the correct Cauchy data at $t=0$. Thus by uniqueness of solutions to quasilinear hyperbolic equations in causal sets, $h \equiv 0$ in all of $R$. But this of course means that $h = 0$ in Taylor series at $(\mathcal H_1\un \mathcal H_2) \n M([0,u^\ast)\times [0,v^\ast))$, as desired.\end{proof}

\begin{figure}[htbp]
\centering
\includegraphics{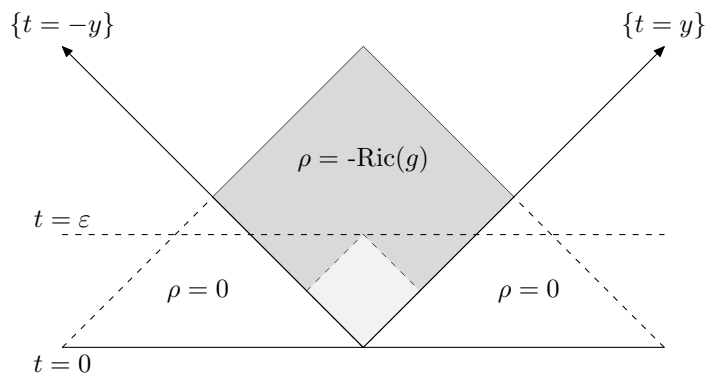}
\caption{The proof of \cref{thm:C2:existence}. The dark shaded region corresponds to $M([0,u_3]\times[0,v_3))$. The light shaded region corresponds to the domain of existence of the solution $h$ of the characteristic initial value problem.}
\label{fig:C2:shorttime}
\end{figure}

Shrinking $u^\ast$ and $v^\ast$, we may assume that $g+h$ is a Lorentzian metric on $M([0,u^\ast)\times[0,v^\ast))$. To complete existence, we just need:
\begin{lem}\label{thm:C2:wavegauge}The solution $h$ of $Q(h) = 0$ with $h$ vanishing in Taylor series at the initial surfaces satisfies that $g+h$ is in wave gauge with respect to $g$ on its entire domain of definition, and so $\Ric(g+h) = 0$.\end{lem}
\begin{proof}
Since $h$ is $0$ in Taylor series, $F(h) = 0$ initially. As is well-known (see Chapter VI of \cite{ChoGene}), the contracted Bianchi identities imply that $F(h)$ satisfies a linear hyperbolic system, and so by \cref{thm:C2:Rend}, $F(h) \equiv 0$ identically on $M([0,u^\ast)\times [0,v^\ast))$, and so $g+h$ is in wave gauge with $g$, and thus $G(h) = 0$ and $\Ric(g+h) = 0$.\end{proof}

To convert back to double-null gauge, we simply apply \cref{thm:C2:conversion}, perhaps shrinking $u^\ast$ and $v^\ast$. This shows existence.

\section{Uniqueness}
\label{C:C2:Uniqueness}
We have a found a solution $g$ to $\Ric(g) = 0$ in double-null gauge and having the correct initial data.

Now we address uniqueness. We start with a local result.
\begin{prop} Let $g_1$, $g_2$ be two solutions to the characteristic initial value problem\footnote{We continue to use the notation in \cref{thm:C2:CIVP}.} in the same double-null gauge with the same initial data, and assume that they have the same domain of existence $M'$. Then $g_1 = g_2$ in a small neighbourhood of $S_{0,0}$ in $M'$.\end{prop}
\begin{proof}
We first seek to find a neighbourhood $U$ of $S_{0,0}$ in $M'$, and a map $\phi:U \to M'$, a diffeomorphism onto its image and the identity along $(\mathcal H_1 \un \mathcal H_2)\n U$ such that $\phi^{\ast}g_2$ is in wave gauge with $g_1$. This will be possible provided we can solve the wave maps equation for $\phi:U \to M'$ with the identity as characteristic initial data, and $\phi$ is a diffeomorphism near $S_{0,0}$. We can solve the wave maps equation, since it is semilinear hyperbolic, using \cref{thm:C2:Rend}, after extending $g_1$ and $g_2$ to a larger manifold. Details are similar to those in the proof of \cref{thm:C2:existence} and are omitted. We just need to check, perhaps after shrinking $U$, that $d\phi$ is non-singular in $U$. Indeed, at $S_{0,0}$, $\phi$ being the identity on $\mathcal H_1 \un \mathcal H_2$ implies that $d\phi = \id$ on $S_{0,0}$.\footnote{This means that the full differential is the identity, not the weaker condition the fibre-differential $\slash{d}\phi$ is the identity.}

Now $\phi^{\ast}g_2$ is in wave gauge with respect to $g_1$ and solves $\Ric(\phi^{\ast}g_2) = 0$. Since $\phi$ was the identity on $\mathcal H_1\un \mathcal H_2$, $\phi^{\ast}g_2$ has the same characteristic initial data as $g_1$ and $g_2$. Set $h = \phi^\ast g_2-g_1$. Then $\Ric(g_1+0) = \Ric(g_1+h) = 0$, and both $g_1$, $g_1+h$ are in wave gauge with respect to $g_1$, and so by \cref{thm:C2:Rend}, $h \equiv 0$. 

We now wish to show that $\phi$ is the identity. From the above, we know that $\phi^{\ast}g_2 = g_1$, and both $g_2,g_1$ are in a double-null gauge. Since $\phi$ fixes $\mathcal H_1 \un \mathcal H_2$, \cref{thm:C2:diffeo} implies that $\phi$ is the identity. Thus, $g_1 = g_2$ in a small neighbourhood of $S_{0,0}$ in $M'$.\end{proof}

We need to extend this local result to a result over all of $M'$. We do this by using an open-closed argument, and a theorem on uniqueness to the Einstein equations in ``causal sets.''
\begin{thm}\label{thm:C2:uniqueness}Let $g_1$, $g_2$ be two solutions to the characteristic initial value problem in the same double-null gauge with the same initial data, and assume that $g_1,g_2$ both exist on $M([0,u^\ast)\times[0,v^\ast))$ for $u^\ast \leq a$, $v^\ast \leq b$. Then $g_1 = g_2$ on $M([0,u^\ast)\times[0,v^\ast))$ .\end{thm}
\begin{proof}We first prove the special case $u^\ast = v^\ast$.

Let $\mathcal S \subseteq \R$ be the set of all $0 \leq c < u^\ast = v^\ast$ such that $g_1 = g_2$ for $0 \leq u \leq c, 0 \leq v \leq c$. It suffices to show that $\mathcal S$ is open, closed, and non-empty. Closed is clear, and non-empty follows from the previous proposition. Suppose $a \in \mathcal S$. Consider the fibre $S_{c,0}$. We may set up a characteristic initial value problem with initial hypersurfaces $\mathcal H_1 \n \{u \geq c\}$ and $\{u = c\}$ with initial data the same as on $\mathcal H_1$, and on the second hyperplane induced by $g_1 = g_2$ (since by assumption we know they are equal there). Choose on $S_{a,0}$ the data induced by $g_1=g_2$. The previous proposition now gives some small $\epsilon > 0$ such that $g_1 = g_2$ for $0 \leq v \leq \epsilon$ and $c \leq u \leq c+\epsilon$. Call this region $M(\mathcal R_1)$ for $\mathcal R_1 \subseteq \R^2$. Examining the analogous problem around $S_{0,c}$ (and perhaps shrinking $\epsilon$) shows that $g_1 = g_2$ for $0 \leq u \leq \epsilon$ and $a \leq v \leq c+\epsilon$. Call this region $M(\mathcal R_2)$ for $\mathcal R_2 \subseteq \R^2$. Now consider the region $\mathcal R_3$ in $\R^2$ bounded by the three curves $u+v = c+\epsilon/2$, $u = c+\epsilon/2$ and $v = c + \epsilon/2$. See figure \ref{fig:C2:uniquenessfigure}. Shrinking $\epsilon$, we may ensure this region is contained inside $[0,u^\ast)\times [0,v^\ast)$. If $u \leq c+\epsilon/2$ and $v \leq c+\epsilon/2$, then either $(u,v) \in \mathcal R_i$ for some $i$, or else $u,v \leq c$. Thus $c+\epsilon/2 \in \mathcal S$ provided that we can show $g_1 = g_2$ on $M(\mathcal R_3)$. We will use an open/closed argument to show this.

\begin{figure}[htbp]
\centering
\includegraphics{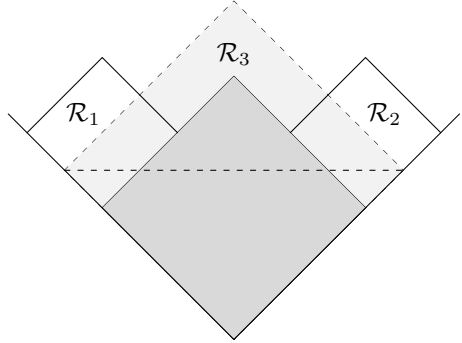}
\caption{The regions in the proof of the special case in \cref{thm:C2:uniqueness}. We know $g_1 = g_2$ on the dark shaded region, and we would like to know $g_1 = g_2$ on the light shaded region.}
\label{fig:C2:uniquenessfigure}
\end{figure}

For $c+\epsilon/2 \leq t \leq 2a+\epsilon$, consider the sets 
\begin{align*}
P_t& := M(\{ u < c+\epsilon/2, v < \epsilon/2, c+\epsilon/2 \leq u+v \leq t)\})\\
F_t&:= M(\{ u < c+\epsilon/2, v < \epsilon/2, t \leq u+v \leq 2c+\epsilon\})\\
H_t &:= M(\{ u < c+\epsilon/2, v < \epsilon/2, u+v = t\}).\end{align*}
$H_t$ is a spacelike hyperplane in $M(\mathcal R_3)$, and $P_t$, $F_t$ are its causal past and future, according to both metrics $g_1\ g_2$ (since they are in the same double-null gauge).

It suffices to show that $g_1 = g_2$ on $P_t$ for all $c+\epsilon/2 \leq t \leq 2c+\epsilon$, so let us consider $\mathcal T$, the set of all $c+\epsilon/2 \leq t \leq 2c+\epsilon$ for which $g_1 = g_2$ on $P_t$. $\mathcal T$ is certainly closed. $\mathcal T$ is also non-empty since $g_1 = g_2$ on $P_{c+2\epsilon/3} \subseteq \mathcal R_1 \un \mathcal R_2 \un M([0,c]\times[0,c])$.

Now we just need to show that $\mathcal T$ is open. Fix $t \in \mathcal T$. We may assume $t < 2c+\epsilon$. $F_t$ is a causal set for both $g_1$ and $g_2$, and so by uniqueness of the Einstein equations in causal sets (see for instance \cite{ChoGene}, \S VI, Theorem 8.8) there exists $s_1,s_2 > t$ and a diffeomorphism $\phi:F_t\n P_{s_1} \to F_t\n P_{s_2}$ which agrees with the identity on $H_t$ to first order such that $\phi^{\ast}g_2 = g_1$. $g_1$, $g_2$ being in double-null gauge with $u,v,N$ and $g_1 = \phi^\ast g_2$ means that $g_2(\phi^\ast du,\phi^\ast du) = g_2(\phi^\ast dv,\phi^\ast dv) = 0$. Now, $\phi^\ast u$ and $u$ agree to first order at $H_t$, as do $\phi^\ast v$ and $v$. By Hamilton-Jacobi theory, this means that $\phi^\ast u = u$ and $\phi^\ast v = v$ everywhere both pairs are defined. It follows that $\phi_\ast N = N$, since $g_1$ and $g_2$ are already in a double null gauge. The action of $\phi$ is the same as flowing backwards along $N$ from a point $p$ to $H_t$ for some time $\tau$, then flowing forwards along $\phi_{\ast} N = N$ for the same time. In particular, the above implies that $s_1 = s_2$ and $\phi$ is the identity, and hence $g_1 = g_2$.\footnote{Notice that we cannot apply \cref{thm:C2:diffeo} since $\phi$ is assumed to be the identity to first order on a spacelike hypersurface, not the identity on a pair of intersecting null hypersurfaces.} Thus there is some $s > t$ such that $g_1 = g_2$ on $P_s$, which shows that $\mathcal T$ is open.

Thus, $g_1 = g_2$ on $M(R_3)$, which completes the proof of the special case.

Now suppose $u^\ast \neq v^\ast$. It suffices to show that the set of points $(u',v') \in [0,u^\ast)\times[0,v^\ast)$ for which $g_1 = g_2$ on $M([0,u']\times[0,v'])$ is open, closed and nonempty. Closed is clear and non-empty follows from the special case. Let us now shown open. Suppose $g_1 = g_2$ on $M([0,u']\times[0,v'])$. If we can show that $g_1 = g_2$ on $M([u',u'+\epsilon)\times [0,v'))$ and $M([0,u')\times[0,v'+\epsilon))$, for some $\epsilon > 0$, then we can use the special case, with initial surfaces $\{u' \leq u \leq u'+\epsilon, \ v = v'\}$ and $\{u = u', \ v' \leq v \leq v'+\epsilon\}$, to show that $g_1= g_2$ on $M([0,u'+\epsilon]\times[0,v'+\epsilon])$. Without loss of generality, let us prove $g_1 = g_2$ on the first region.

\begin{figure}[htbp]
\centering
\includegraphics{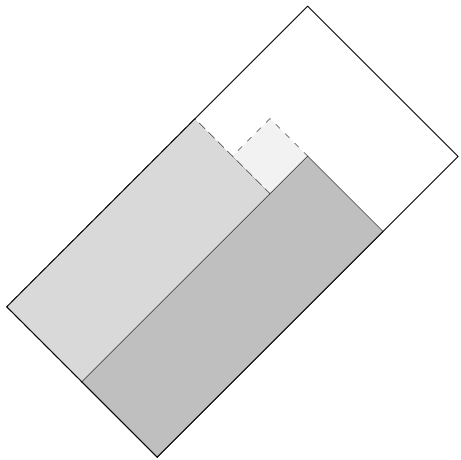}
\caption{The regions in the proof of the general case of \cref{thm:C2:uniqueness}. The dark shaded region is $M([0,u')\times[0,v'))$, the medium shaded region is $M([0,u^\ast)\times[0,k(u^\ast-u'))$ and the light shaded region is $M([u',u'+(v'-k(u^\ast-u')))\times [k(u^\ast-u'),v'))$.}
\end{figure}
With initial surfaces $\{u' \leq u < u^\ast, \ v = 0\}$ and $\{u = u', \ 0 \leq v \leq v'\}$, we may use the special case to obtain $g_1 = g_2$ on \[M([u',u' + \min(u^\ast-u',v')),v')\times[0, \min(u^\ast-u',v'))).\] If $\min(u^\ast-u',v') = v'$, then we're done and $\epsilon = v'$. Otherwise, we may iterate, obtaining $g_1 = g_2$ on $M([u',u^\ast)\times[0, k (u^\ast-u')))$, where $k$ is the largest integer such that $k(u^\ast-u') \leq v'$. If $k(u^\ast-u') = v'$, then we're done again, with $\epsilon = u^\ast - u'$. Otherwise, $v'> k(u^\ast-u')$ but $v' \leq (k+1)(u^\ast-u')$, and so we may use the special case again with initial surfaces $\{u' \leq u \leq u'+(v'-k(u^\ast-u')), \ v = k(u^\ast-u')\}$ and $\{u = u', \ k(u^\ast-u') \leq v \leq v'\}$ to obtain $g_1 = g_2$ on \[M([u',u'+(v'-k(u^\ast-u')))\times [k(u^\ast-u'),v')),\] and hence on \[M([u',u'+(v'-k(u^\ast-u')))\times[0,v')),\] as well.
\end{proof}
\chapter{The theorem of Christodoulou}
\label{C:C3}

\section{\hspace{-2.55pt}The short-pulse ansatz}
\label{C:C3:ansatz}

The goal of this chapter is to prove \cref{thm:C1:ChrVague}, as well as its extension \cref{thm:C1:Smoothness}. Let us begin by reviewing the short-pulse ansatz of Christodoulou, and stating the precise theorem we wish to prove.

Consider $\R^4 = \R_t \times \R^3_y$, and set $u =1-\frac{|y|-t}{2}$, $\und{u} = \frac{t+|y|}{2}$. For $\delta > 0$ consider the region $M_\delta \subseteq \R^4$ picked out by
\[M_\delta = \{0 \leq u < 1, 0 \leq \und{u} \leq \delta\}.\]

Then $(M_\delta,u,\und{u})$ is a doubly-foliated manifold, in the sense of \cref{def:C2:doublefoliated}, over $[0,1)\times[0,\delta]$ and with fibre $S^2$. We will fix a vector field $N = \pa_u = \pa_t - \pa_r$ on $M_\delta$. Here $r\pa_r := y^i\pa_{y_i}$. Let us set up a characteristic initial value problem with initial data on $\mathcal H_1\un \mathcal H_2 = \{u = 0\} \un \{\und{u} = 0\}$. We call the initial data the \emph{short-pulse ansatz} or \emph{short-pulse data}. Following \cref{thm:C2:CIVP}, we need to specify a number of things:
\begin{romanumerate}
\item a vector field $L$ on $\mathcal H_2$ with $L\und{u} = 1$;
\item a Riemannian metric $\slash{g}$ on $S_{0,0}$;
\item a section $\hat{\slash{g}}$ of $\Sym^2(T^\ast S)$ over $\mathcal H_1\un \mathcal H_2$ which has Riemannian signature and satisfies $\hat{\slash{g}} = \slash{g}$ at $S_{0,0}$;
\item a pair of smooth functions $f_i$ on $S_{0,0}$;
\item a vector field $W$ on $S_{0,0}$.
\end{romanumerate}

For us $L = \pa_{\und{u}} = \pa_t + \pa_r$ on $\mathcal H_2$. We will use $N$, $L$ to fix a product decomposition $M_\delta = [0,1)\times[0,1]\times S^2$ obtained by flowing from $S_{0,0}$ to $\mathcal H_2$ via $L$, and then to all of $\mathcal M$ via $N$.

Let is now specify the rest of the data. Let $\mathring{\slash{g}}$ denote the metric on the round sphere $S^2$ of radius $1$. Let us use $N,L$ to extend $\mathring{\slash{g}}$ to a section of $\Sym^2(T^\ast S)$ over $\mathcal H_1\un \mathcal H_2$.

Consider for $t \in [0,1]$ a smooth one-parameter family $\T(t)$ of type $(1,1)$ tensors on $S^2$, i.e.\ $\T \in C^\infty([0,1]\times S^2; TS^2\otimes T^\ast S^2)$. Let us assume that $\T(0) = 0$, $\Tr(\T) = 0$, and $\T$ is $\mathring{\slash{g}}$-symmetric\footnote{Interpreted either as a map $TS^2 \to TS^2$ or $T^\ast S^2 \to T^\ast S^2$; symmetry as a map in one interpretation is equivalent to the other.} The choice of the tensor $\T$ determines all of the interesting behaviour of the short-pulse ansatz, and for this reason we call it \emph{the short-pulse tensor}. Let us recall that we may define\footnote{The symbol $\mathbf{T}^k$ is well-defined interpreting $\mathbf{T}$ as a map $TS^2 \to TS^2$ or $T^\ast S^2 \to T^\ast S^2$, or indeed just by a $k$-fold contradiction; in this case all interpretations give the same result.  Cf. \cref{fn:C1:wut}.}
\[\exp(\T) = \sum_{k=0}^\infty \frac{\T^k}{k!},\] which is positive and $\mathring{\slash{g}}$-symmetric. Using $L$ to fix the product decomposition on $\mathcal H_2$, for $\delta > 0$ the tensor $\T_\delta$ defined in coordinates by $\T_\delta(\und{u}) = \T(\und{u}/\delta)$ is a well-defined section of $TS\otimes T^\ast S$ over $\mathcal H_2$.

Let us define the energy $\mathbf{E}\: [0,1]\times S^2 \to \R$ of $\T$ by 
\[\mathbf{E}(t,\theta) = \frac{1}{2}\int_0^t |\pa_t \T|_{\mathring{\slash{g}}}^2(s,\theta) \ ds.\]

Following Christodoulou \cite{ChrForm}, we may give our initial data, which we call \emph{short-pulse data}:
\begin{romanumerate}
\item $L = \pa_{\und{u}}$;
\item $\slash{g} = \mathring{\slash{g}}$;
\item $\hat{\slash{g}} = (1-u)^2\mathring{\slash{g}}$ over $\mathcal H_1$, and $\hat{\slash{g}} = \mathring{\slash{g}}\exp(\sqrt{\delta}\T_\delta)$ over $\mathcal H_2$;
\item $f_1 = 2$, $f_2 = -2$;
\item $W = 0$.
\end{romanumerate}

The precise version of the Theorem of Christodoulou we seek to prove is:
\begin{thm}[Christodoulou \cite{ChrForm}]\label{thm:C3:ChristMain}Consider the characteristic initial value problem over $M_\delta$ with initial data as given above. Then, for any $u^\ast < 1$, if $\delta > 0$ is small enough, there exists a solution $g_\delta$ of $\Ric(g_\delta) = 0$ on $M_\delta([0,u^\ast]\times [0,\delta])$ in double-null gauge with $M_\delta,u,\und{u},N$ and $g$ induces the initial data. 

Moreover, if $K$ is any compact subset of
\[U = \{(u,t,\theta) \in [0,1)\times [0,1]\times S^2 \: 1-\frac{1}{4}\mathbf{E}(t,\theta) < u \leq u^\ast\},\]
then for $\delta$ perhaps smaller, the mean curvatures $\frac{1}{2}\slash{\tr}\Lie_L \slash{g}_\delta, \ \frac{1}{2}\slash{\tr}\Lie_N \slash{g}_\delta$ are both strictly negative at any point $(u,\und{u},\theta) \in M_\delta$ whenever $(u,\und{u}/\delta,\theta) \in K$. 

If $\sup_{\theta \in S^2} \mathbf{E}(1,\theta) < 4$, then choosing $\delta$ smaller, there are no trapped surfaces on $\mathcal H_1 \un \mathcal H_2$.
\end{thm}

The formation of trapped surfaces follows easily, as we now indicate. Suppose \[\inf_{\theta \in S^2} \mathbf{E}(1,\theta) > 4(1-u^\ast).\] Set\[t^\ast = \inf\left\{t\: \inf_{\theta \in S^2} \mathbf{E}(t,\theta) > 4(1-u^\ast)\right\}.\]

\begin{figure}[htbp]
\centering
\includegraphics{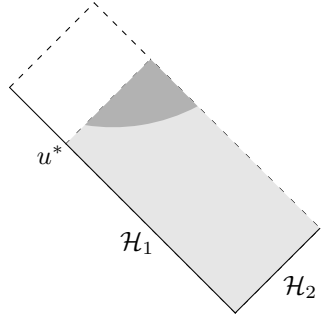}
\caption{The base of the submersion of the double-null foliation $(M_\delta,u,\protect\und{u})$. The light-shaded region is the domain of existence, and the dark-shaded region is the set $\{(u,\protect\und{u})\: (u,\protect\und{u}/\delta) \in \epi F\}$. For $\delta$ small enough, the fibre of any point in the dark-shaded region (away from the boundary) is trapped.}
\end{figure}

Then on $(t^\ast,1)$,
\[F(t) := 1-\inf_{\theta \in S^2}\frac{1}{4}\mathbf{E}(t,\theta) < u^\ast,\]
and so the epigraph of $F$ on $(t^\ast,1)\times [0,u^\ast)$,
\[\epi F := \{(u,t) \: 0 \leq u < u^\ast,\ t^\ast < t < 1,\ u > F(t)\}\] is nonempty, and $(u,t) \in \epi F$ implies that $(u,t,\theta) \in U$. Thus, for any $(u,\und{u})$ for which $(u,\und{u}/\delta) \in \epi F$, (for $\delta$ small enough), the entire fibred sphere $S_{u,\und{u}}$ is trapped. The condition $\sup_{\theta \in S^2} < 4$ ensures that the mean curvature of $\slash{g}_{\delta}$ in the direction of $L$ are positive on $\mathcal H_2$. We will show in \cref{C:C3:FormTrapped} that this means that there are \emph{no} trapped surfaces in $\mathcal H_1 \un \mathcal H_2$.
It is not hard to find $\T$ so that $0 < \inf_{\theta \in S^2} \mathbf{E}(1,\theta) < \sup_{\theta \in S^2} \mathbf{E}(1,\theta) < 4$, so taking $u^\ast$ sufficiently close to $1$, we conclude that trapped surfaces can form dynamically.

Now let us consider the region
\[\mathcal R = \{(p,\delta) \in \R^4\times [0,1)\: p \in M_\delta\},\] and extend the definitions of $u$, $\und{u}$ to $\mathcal R$. Considering $(u,\und{u})$ as a projection, we may consider $\mathcal R$ as a (topological) fibre bundle over
\[[0,1)\times [0,\delta]\times [0,1)_\delta\]
with fibre $S^2$.
We may interpret the collection of all short-pulse data for $\delta > 0$ as initial data for a problem on $\mathcal R$. However, $\mathcal R$ is singular at $\delta = 0$. To resolve this, we perform a parabolic blowup at $\delta = 0$, $\und{u}=0$, i.e.\ we blow up the cylinder $\{\delta = 0, \ \und{u} = 0\}$, and then redefine the smooth structure by requiring a square root of a bdf of the introduced face to be smooth. Call the resulting manifold $\mathcal M$. 

 Since the fibres of $\mathcal R$ are also fibred over $\{\delta = 0, \und{u} = 0\}$, $\mathcal M$ is also a fibre bundle. To see this in coordinates, let us set $v = \und{u}/\delta$, $x = \sqrt{\delta}$. Then $\mathcal M$ is diffeomorphic to $[0,1)_u \times[0,1]_v \times S^2_\theta\times [0,1)_\delta$, and the projection is given by $(u,v,\theta,x) \mapsto (u,v,x)$. There is also a natural blowdown map $\beta:\mathcal M \to \R$ sending $(u,v,\theta,x) \mapsto (u,x^2v,\theta,x^2)$.

\begin{figure}[htbp]
\centering
\includegraphics{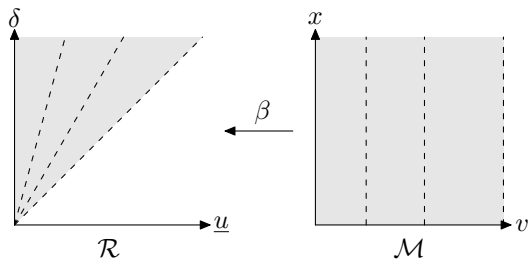}
\caption{A side view of $\mathcal M$ and $\mathcal R$. Also shown are some lines in $\mathcal R$ and their preimages under $\beta$.}
\end{figure}

One sees that on $\mathcal M$, $\hat{\slash{g}}$ pulls back to a smooth section. However, $L =x^2 \pa_v$ is now degenerate at $\delta = 0$. This is related to the fact, pointed out in the introduction, that the Minkowski metric takes the form
\[x^{-2}\mathring{g} = -2(du\otimes dv + dv\otimes du) + (1-u+x^2v)^2(x^{-2}\slash{\mathring{g}}),\]
so cannot be a smooth metric on $T\mathcal M$. To resolve this difficulty, we instead consider a rescaled bundle $^\sp T\mathcal M$, which we call the short-pulse bundle, and which is defined by specifying a local basis of smooth sections (such a bundle exists by Swan's theorem). These vector fields are $\pa_u$, $\pa_v$, and $x\Theta_1$, $x\Theta_2$, where $\Theta_i$ ($i = 1,2$) is a local basis of the vertical bundle of the fibred spheres. The Minkowski metric (after rescaling it by $x^2$) is easily seen to be a metric on ${}^\sp T\mathcal M$ (observe that we do not take $\pa_x$ to be a smooth section). Furthermore, denoting by $g_\delta$ the solution to the characteristic initial value problem on $M_\delta$ with short-pulse data, we expect that the family $x^{-2}\beta^\ast g_\delta$ is smooth, \emph{as a metric on ${}^\sp T\mathcal M$}.

We will derive the existence statement of \cref{thm:C3:ChristMain} from \cref{thm:C3:SPCIVP}, which is a more general long-time existence theorem for certain classes of data (of which the short-pulse ansatz is an example) for the characteristic initial value problem for metrics on the short-pulse bundle. After this, the statement about the sign of the mean curvatures (and hence the formation of trapped surfaces) will become an easy computation.

In the following three sections, we will abandon the specific set up of the short-pulse ansatz, and consider more general short-pulse bundles and metrics on them. This shows that the long-time existence portion of \cref{thm:C3:ChristMain} is not particular to the exact setup of the short-pulse ansatz. Only in \cref{C:C3:FormTrapped}, the final section of this chapter, will we return to the short-pulse ansatz, where we will apply our general results to establish long-time existence, and then prove the formation of trapped surfaces.

\section{The characteristic initial value problem for the Einstein equations for short-pulse metrics}
\label{C:C3:SPCharValue}
\subsection{The short pulse bundle}
\label{C:C3:shortpulsegeometry}
We start off with the analogue of a doubly-foliated manifold. 
\begin{defn}\label{def:C3:doublefoliated}Let $\mathcal R \subseteq \R^2$ be a codimension-zero embedded manifold with corners. We consider $\R^2$ as $\R_u \times \R_v$, and equip $\mathcal R$ with the coordinates $(u,v)$. Let $M$ be a smooth manifold with corners of dimension at least $3$. Fix $\epsilon > 0$ and define the space $\mathcal M = [0,\epsilon)_x \times M$ and let $x:\mathcal M \to [0,\epsilon)$ be the map taking $(x,p) \mapsto x$. Also denote by $x: [0,\epsilon) \times \mathcal R \to [0,\epsilon)$ the map taking $(x,q) \mapsto x$. Let $\pi: \mathcal M \to [0,\epsilon) \times \mathcal R$ be a smooth proper surjective submersion with fibres diffeomorphic to the same closed (and connected) manifold $S$. Using the coordinates $(u,v)$ on $\mathcal R$ and $x$ on $[0,\epsilon)$, identify the components of $\pi = (x,u,v)$. Make the assumptions that, separately, the fibres of $(x,u)$ and $(x,v)$ are themselves connected. A \emph{parametrized doubly-foliated manifold} consists of the data $(\mathcal M = [0,\epsilon)_x \times M,\mathcal R, \pi = (x,u,v),S)$.\end{defn}
We will make the convention that if $(\mathcal M,x,u,v)$ is a parametrized doubly-foliated manifold, with the total space given by a caligraphic Roman letter, then the associated manifold $M$ for which $\mathcal M = [0,\epsilon)_x\times M$ is denoted by the same Roman letter in majuscule form. We will denote $M_{x_0} = \{p \in \mathcal M: x(p) = x_0\}$, and make the usual conflation between the function $x$ and its values in $[0,1)$. We will drop the adjective ``parametrized'' if there is no risk of confusion with a doubly-foliated manifold of \cref{def:C2:doublefoliated}, and will call the doubly-foliated manifolds of \cref{def:C2:doublefoliated} ``ordinary'' if instead we wish to emphasize the distinction.

We will primarily be interested in the behaviour of objects on $\mathcal M$ near $x = 0$. Thus, we will often (and sometimes implicitly) replace $\mathcal M$ with $\mathcal M \n \{0 \leq x < \delta\}$, for $0 < \delta \ll \epsilon$, as needed, and say that the desired behaviour holds ``for small $x$.'' In particular, we will never care about the exact value of $\epsilon$ in the definition.

Notice that the definition implies that for each $x$, $(M_x,u,v)$ is an ordinary doubly-foliated manifold. We will denote by $S_{u,v,x}$ the fibre of $\pi$, and $S_{u,v} \iso [0,\epsilon)_x\times S$ the fibre of $(u,v)$. As for an ordinary doubly-foliated manifold, if $R \subseteq \mathcal R$, set $\mathcal M(R) := \{p \in \mathcal M\: (u,v)(p) \in R\}$.

The manifold $M_0$ is special and constitutes a boundary face of $\mathcal M$. Let us suppose that the projection $\pi$ is canonically trivial when restricted to $M_0$, i.e.\ $M_0$ is equipped with a diffeomorphism $\tau: M_0 \to \mathcal R \times S$ which we regard as fixed, so that $\pi = (0,u,v):M_0 \to \mathcal R$ is just the projection onto the first factor of $M_0 \iso \mathcal R \times S$. If $\tau$ is given, we will abuse notation and simply treat $\tau$ as the identity. Henceforth, all of our parametrized doubly-foliated manifolds will come equipped with such a diffeomorphism $\tau$.

\begin{lem}\label{thm:C3:globaltriv}Suppose $(\mathcal M,x,u,v)$ is a parametrized doubly-foliated manifold (with $M_0$ canonically trivial). Then there is a diffeomorphism $\phi:\mathcal M \to [0,\epsilon)\times M\times \mathcal R$ such that $\phi$ is the identity on $M_0$.\end{lem}
\begin{proof}The key is that it is always possible to lift vector fields from the codomain of a proper submersion to the domain, and by the assumptions that the fibres are connected, various flows are well-defined. The vector field $\pa_x$ on $\mathcal M$ coming from the initial product structure is well-defined. Via the trivialization $M_0 = \mathcal R \times S$, we may lift the coordinate vector fields $\pa_u, \pa_v$ from $\mathcal R$ to $M_0$. Extend the vector field $\pa_u$ to a vector field $U$ on all of $T\mathcal M$ by requiring $[\pa_x,U] = 0$. Notice that $U$ is tangent to each $M_x$, and $Uu = 1$, $Uv = 0$, $Ux = 0$. Extend $\pa_v$ to a vector field $V$ in the same way. Since $[U,V] = 0$ on $M_0$, it follows that $[U,V] = 0$ on all of $\mathcal M$.

Fix $(u_0,v_0) \in \mathcal R$. \Cref{thm:C2:complete} then implies that we may pick any $p \in \mathcal M$, and flow it to $S_{u_0,v_0,0}$ by flowing along $U$, $V$, $\pa_x$ in any order (since they commute). If $(x,u,v)(p) = (x,u,v)$, define the $\phi:[0,\epsilon)\times\mathcal M \to [0,\epsilon)\times \mathcal R\times S_{u_0,v_0,0}$ by requiring its first coordinate to be $(x,u,v)$, and its last coordinate is the image of $p$ under the aforementioned flow. Then $\phi$ is a diffeomorphism. It is also a global trivialization by definition, and by construction is the identity on $M_0$.
\end{proof}
\begin{rk}It appears that in the definition of a doubly-foliated manifold, the level sets of $(x,u)$ and $(x,v)$ should be null hyperplanes for some metric to appear. However, this is not necessarily the case, and we may also foliate a manifold by space-like hyperplanes in one direction, and the product of a time-like direction and a codimension-two hyperlane in the other This will be important when we discuss the Cauchy problem in \cref{C:C3:Cauchy}, where $(u,v)$ will be replaced by functions $(t,y)$ on a space of the form $\R_t\times \R_y \times S$, and $\R_y\times S$ will be a spacelike hypersurface fibred over by fibres $S$.\end{rk}

We will be interested in a rescaled version of $T\mathcal M$.
\begin{defn}\label{def:C3:spbundle}Suppose $(\mathcal M,x,u,v)$ is a parametrized doubly-foliated manifold. Let ${^\sp T\mathcal M}$ denote the vector bundle over $\mathcal M$ whose smooth sections are smooth sections of $T\mathcal M$ which are annihilated by $dx$ and which at $M_0 = \mathcal R \times S$ are tangent to the fibres of the projections onto the second factor. We call ${^\sp T\mathcal M}$ the \emph{short-pulse tangent bundle}.\footnote{The reader familiar with the semiclassical calculus and the edge calculus of Mazzeo \cite{MazElli} will recognize $\spT$ as a sort of \emph{semiclassical edge bundle}, with the semiclassical parameter $x$. However, the edge structure is mixed up with the semiclassical structure. This will be a useful point of view when discussing the Cauchy problem, when we will consider sections of $^\sp T\mathcal M$ as $x$-parametrized sections of $T M_x$ which scale according to the scaling in $^\sp T\mathcal M$. However, we prefer to use the term ``short-pulse'' to describe the bundle in order to avoid piling on the adjectives. The bundle is also of a similar form to the \emph{adiabatic} bundles considered by Mazzeo-Melrose in \cite{MazMelAdia}.}\end{defn}
The choice of ${^\sp T\mathcal M}$ of course depends on the choice of trivialization of $M_0$, although we suppress this dependence since the interaction of different trivializations will not play a role in our analysis.

Using \cref{thm:C3:globaltriv}, we may choose local coordinates $(x,u,v,\theta)$, where $\theta = (\theta^1,\cdots,\theta^n)$ are local coordinates on $S$. In these coordinates, a basis of smooth sections of ${^\sp T\mathcal M}$ consists of $\pa_u,\pa_v,x\pa_{\theta^1},\ldots,x\pa_{\theta^n}$ ($n = \dim S$). Swan's theorem provides for the existence of such a vector bundle, although it is not hard to construct it by hand.

We also consider the bundle $^\sp TS \subseteq {^\sp T\mathcal M}$, the subbundle of all vectors which remain tangent to the fibres $S$, and write $\iota:{}^\sp TS \to {^\sp T\mathcal M}$ for the inclusion map. Notice that the putative subbundle $^\sp T\mathcal M$ of vectors tangent to the fibres of $x$ is already all of $\spT$, so we will not have cause to adopt specific notation for it.

Associated to ${^\sp T\mathcal M}$ is its dual bundle ${{^\sp T^\ast\mathcal M}}$, the \emph{short pulse cotangent bundle.} We may also consider ${^\sp T^\ast S}$, the dual bundle to $^\sp TS$.

In the previous coordinates, sections of ${{^\sp T^\ast \mathcal M}}$ are smooth combinations of \[du,dv, d\theta^1/x,\ldots d\theta^n/x.\] 

Observe that ${^\sp T\mathcal M} \subseteq \ker dx \subseteq T\mathcal M$ and $^\sp TS \subseteq TS$. Furthermore, over $x > 0$, the inclusions ${}^\sp T\mathcal M \embeds \ker dx$ and $\iota: {}^\sp TS \embeds TS$ are isomorphisms, since we are allowed to multiply/divide by $x$. Dually, $T^\ast \mathcal M/ \vspan\{dx\} \subseteq {}^\sp T^\ast \mathcal M$ and $T^\ast S \subseteq TS$ and the inclusions are isomorphisms off of $x > 0$.

One may also form the tensor products ${^\sp T^p_q \mathcal M} := {^\sp T\mathcal M}^{\otimes p}\otimes {{^\sp T^\ast \mathcal M}}^{\otimes q}$, and the tensor products ${^\sp T^p_q S} := {^\sp TS}^{\otimes p}\otimes {{^\sp T^\ast S}}^{\otimes q}$.

\begin{defn}We call sections of ${^\sp T^p_q \mathcal M}$ \emph{type $(p,q)$ short-pulse tensors}.\footnote{Unfortunately this conflicts with the notion of \emph{the} short-pulse tensor $\mathbf{T}$ of the initial data in the short-pulse ansatz. Since we will always use $\mathbf{T}$ to denote this tensor, and will never use the terms in the same context, this will not be a cause for confusion.} We call type $(1,0)$ short-pulse tensors \emph{short-pulse vector fields} and type $(0,1)$ short-pulse tensors \emph{short-pulse one forms}. We similarly call sections of the tensor product ${^\sp T^p_q S} $ \emph{type $(p,q)$ short-pulse fibre tensors}. We call type $(1,0)$ short-pulse tensors \emph{short-pulse fibre vector fields} and type $(0,1)$ short-pulse tensors \emph{short-pulse fibre one forms}.\end{defn} If $h$ is a short-pulse tensor, we will denote by $h_x$ the restriction of $h$ to $M_x$, which, for $x > 0$, may be thought of as an ordinary tensor.

\begin{defn}Let $g$ be a symmetric type $(0,2)$ short-pulse tensor, i.e.\ $g$ is a section of $\Sym^2({{^\sp T^\ast\mathcal M}})$. Suppose that $g$ is non-degenerate and is of Riemannian (resp. Lorentzian) signature. We will call $g$ a \emph{short-pulse Riemannian (resp. Lorentzian) metric}. We similarly call a non-degenerate section (no signature assumptions) of $\Sym^2({{^\sp T^\ast S}})$ a \emph{short-pulse fibre metric}, or a \emph{short-pulse pseudo-Riemannian metric}. \end{defn}

We observe that ${^\sp T\mathcal M}$ is closed under Lie bracket, i.e.:
\begin{lem}\label{thm:C3:algebroid}Suppose $X_1,\ X_2$ are short-pulse vector fields. Then $[X_1,X_2]$, with the commutator interpreted as between sections of $T\mathcal M$, in fact is a short-pulse vector field.\end{lem}
\begin{proof}This is easiest to see in the local coordinates given above, since we may write for $i=1,2$
\[X_i = f_i\pa_u + g_i\pa_v + xh_i^1\pa_{\theta^1} + \cdots + xh_i^n\pa_{\theta^n},\]
where $f_i,g_i,h^j_i$ are smooth functions. Evaluating the commutator, it is clear that $[X_1,X_2]$ is a linear combination (over $C^{\infty}(\mathcal M)$) of \[\pa_u,\pa_v,x\pa_\theta^j, x^2\pa_\theta^j,\] with coefficients given in terms of $f_i,g_i,h_i^j$ and their derivatives. Such a vector field is clearly a section of ${^\sp T\mathcal M}$.\end{proof}

Using \cref{thm:C3:algebroid}, we may extend \cref{thm:C2:LieDerI} and \cref{thm:C2:LieDerII} to this setting. Namely:
\begin{cor}\label{thm:C3:LieDer}Let $X$ be a short-pulse vector field which is $\pi$-related to a vector field on the base $[0,\epsilon)\times \mathcal R$. Then $\Lie_X$ is a well-defined differential operator between sections of $^\sp T_q^p S^2$.\end{cor}
\begin{proof}We know that if $\alpha$ is a section of $T_q^p S^2$ over $x> 0$, then so is $\Lie_X \alpha$. So it suffices to show that if $\alpha$ is a section of $^\sp T_q^p S^2$ then so is $\Lie_X \alpha$. \Cref{thm:C3:algebroid} shows that this is true if $\alpha$ is a short-pulse fibre vector field, since $\Lie_X \alpha$ is both fibre-tangent and a short-pulse vector field. If $\alpha$ is a short-pulse fibre one-form, then for any short-pulse fibre vector field, $V$,
\[(\Lie_X \omega) (V) = X\omega(V) - \omega(\Lie_X V),\] which then shows that $\Lie_X \omega$ is a short-pulse fibre one-form, since it pairs with short-pulse fibre vector fields. The general result follows form tensoring.
\end{proof}

Let $g$ be a short-pulse pseudo-Riemannian metric, and $\nabla$ its Levi-Civita connection (defined in the usual way on $TM_x$ for $x>0$). Certainly for $x > 0$, $\nabla:C^\infty(M_x;TM_x) \to C^\infty(M_x;T^\ast M_x \otimes T M_x)$. As a consequence of \cref{thm:C3:algebroid}, we have the following corollary on how it interacts with $\spT$:
\begin{cor}\label{thm:C3:connectionisnice}The connection $\nabla$ extends to a map \[\nabla: C^\infty(\mathcal M;\spT) \to C^\infty(\mathcal M;{}^\sp T^1_1 \mathcal M).\]\end{cor}
\begin{proof}For short-pulse vector fields $X,Y,Z$ we may express the connection via the Koszul formula:
\begin{align*}
g(\nabla_X Y,Z) &= \frac{1}{2}\left(Xg(Y,Z) + Yg(Z,X) - Zg(X,Y)\right.\\
&\left. + g([X,Y],Z)-g([Y,Z],X) - g([X,Z],Y)\right).\end{align*}
The inner products in the first line are all smooth since $g$ is a short-pulse metric. Since $X,Y,Z$ are in particular smooth vector fields, it follows that the fist line is smooth. The second line is smooth since all commutators are short-pulse vector fields by \cref{thm:C3:algebroid}. Thus $\nabla_X Y$ is a smooth short-pulse vector even down to $x=0$. This implies that we may define $\nabla_v Y(p)$ at a fixed vector $v \in {}^\sp T \mathcal M$ at a point $p \in \mathcal M$ by extending $v$ arbitrarily to a short-pulse vector field $X$ in a neighbourhood of $p$ and setting $\nabla_v Y(p) = \nabla_X Y(p)$.\end{proof}

We obtain immediately:
\begin{cor}\label{thm:C3:curvatureisnice}If $g$ is a short-pulse metric, then $\Riem(g)$ and $\Ric(g)$ are both short-pulse tensors.\end{cor}
\begin{proof}Let $X,Y,Z,W$ be short-pulse vector fields. From the definition of curvature
\[\Riem(g)(X,Y,Z,W) = R(X,Y,Z,W) = g(\nabla_Y\nabla_X Z - \nabla_X\nabla_Y Z -\nabla_{[Y,X]},Z,W),\]
it follows that $R(X,Y,Z,W)$ is smooth. However, it is linear, and therefore is a short-pulse tensor. It follows that $\Ric(g)$, being a contraction of $\Riem(g)$ against $g^{-1}$, is also a short-pulse tensor.\end{proof}

Similarly, we have:
\begin{lem}\label{thm:C3:paxisnice}The vector field $x\pa_x$ determines a well-defined differential operator $\Lie_{x\pa_x}$ mapping $C^\infty(\mathcal M;\spT)$ to itself.\end{lem}
\begin{proof}The operator $\Lie_{\pa_x}$ is a well-defined operator between ordinary vector fields of $T\mathcal M$. To show that short-pulse vector fields are mapped back to themselves, pick a local basis $\pa_u,\pa_v, x\pa_\theta^j$ and apply $\Lie_{x\pa_x}$.\end{proof}
We abuse notation and write $x\pa_x$ for $\Lie_{x\pa_x}$ because of the product decomposition of $\mathcal M$ into $[0,\epsilon)_x \times M$.

\begin{defn}\label{def:C3:DNG}Let $(\mathcal M ,x,u,v)$ be a parametrized doubly-foliated manifold, and $g$ a Lorentzian short-pulse metric on $\mathcal M$. We say that $(\mathcal M,g,x,u,v)$ is a \emph{short-pulse double-null foliation}, if $u$, $v$ are optical functions for $g$, i.e.\ 
\[g^{-1}(du,du) = g^{-1}(dv,dv) = 0,\]
where we treat $d$ as mapping sections of $C^{\infty}(\mathcal M)$ to $C^\infty(\mathcal M; {}^\sp T^\ast \mathcal M)$ (i.e.\ quotient out by $\vspan\{dx\}$). We additionally suppose that $(\mathcal M,g)$ is time-orientable, and that $u,v$ are increasing towards the future, in the sense that there is a globally defined short-pulse vector field $T$ with $g(T,T) < 0$, and $Tu > 0$, $Tv > 0$. \end{defn}
We will often omit the adjective ``short-pulse'' if it is clear from context.

We also have the analogue of \cref{thm:C2:DNGeo}. If $(\mathcal M,g,x,u,v)$ is a double-null foliation, we set $\slash{g} = \iota^\ast g$ to be the restriction of $g$ to fibre-tangent vectors, and $N' = -2\grad v$, $L' = -2\grad u$ (which are short-pulse vector fields), and putatively define $0 < \Omega$ by $-2\Omega^{-2} = g(N',L')$. Then:
\begin{prop}Let $(\mathcal M,g,x,u,v)$ be a double-null foliation. Then $\slash{g}$ is Riemannian, $\Omega > 0$ defined by $-2\Omega^{-2} = g(N',L')$ is well-defined, and $\vspan\{N',L'\}$ is transverse to $^\sp TS$.\end{prop}
\begin{proof}Observe that $\gradt u = g^{-1}(du,\bullet)$, so $\gradt u$ is a short-pulse vector field, and likewise for $\gradt v$. The rest of the proof is identical to \cref{thm:C2:DNGeo}.\end{proof}

We will also rescale and set $N = -2\Omega^2\grad u$, $L = -2\Omega^2\grad v$, so that $Nu = 1, Nv = 0$, $Lu = 0,Lv = 1$.
Observe $[N,L]$ is a short-pulse vector field, and hence a short-pulse fibre vector field since $[N,L] u = [N,L]v = 0$.

We may of course reconstruct a short-pulse metric from $N,L$, $\Omega$ and $\slash{g}$, as in \cref{thm:C2:Reconstruct}.
\begin{prop}\label{thm:C3:ReconstructSP}Let $(\mathcal M,x,u,v)$ be a parametrized doubly-foliated manifold. Let $N,\ L$ be short-pulse vector fields which are lifts of $\pa_u$ and $\pa_v$, respectively, let $\slash{g}$ be a short-pulse fibre Riemannian metric, and $0 < \Omega \in C^\infty(M;\R)$. Then there is a unique short-pulse Lorentzian metric $g$ such that $(\mathcal M,g,x,u,v)$ is a double-null foliation and $g$ induces this data, i.e.\ $\iota^\ast g = \slash{g}$, $-2\Omega^{-2} = g^{-1}(-2du,-2dv)$, $N = -2\Omega^2 \grad v$, $L = -2\Omega^2 \grad u$.\end{prop}

Finally, we introduce the notion of a short-pulse double-null gauge.
\begin{defn}Let $(\mathcal M,g,x,u,v)$ be a short-pulse double-null foliation. Let $N \in C^\infty(\mathcal M;{^\sp T\mathcal M})$ be a lift of $\pa_u$. We call the $6$-tuple $(M,g,x,u,v,N)$ \emph{short-pulse double-null gauge} if $N = -2\Omega^2 \grad v$ (with $\Omega$ as above). In other words, $N$, must equal with the $N$ introduced above.\end{defn}
As above, we will omit the adjective ``short-pulse'' if it is clear from context, and refer to the double-null gauges of \cref{def:C2:doublenullgauge} as ``ordinary'' double-null gauges if we with to emphasize the distinction. If $(\mathcal M,g,x,u,v,N)$ is a double-null gauge, we will often say that $g$ is in a double null gauge with $u,v,N$ (the background manifold $\mathcal M$ and function $x$ being clear from context).

It is clear that for any $x > 0$, over $M_x$ a short-pulse double-null foliation/gauge is the same thing as a ordinary double-null foliation/gauge, respectively, in the sense of \cref{def:C2:doublefoliated}/\cref{def:C2:doublenullgauge}, respectively.

Let us now see what this looks like in canonical coordinates. Choose some $(u_0,v_0) \in \mathcal R$ and set $\mathcal H = u^{-1}(u_0)$. Choose arbitrary local coordinates $\theta^i$ on $S_{u_0,v_0,0}$. This induces coordinates $(x,u,v,\theta^i)$ on $\mathcal M$ by flowing out $\theta^i$ via $\pa_x$ to $S_{u_0,v_0}$, then via $L$ to $\mathcal H$, and then via $N$ to all of $\mathcal M$. Set $\slash{k} = x^2\slash{g}$, an ordinary fibre-metric, and let $\widetilde{f}^i$ be the smooth functions defined by $\pa_v \widetilde{f}^i = [N,L]^i$. Since $[N,L]$ is a short-pulse fibre vector field, $\widetilde{f}^i = xf^i$ for a smooth function $f$.

In these coordinates, $g$ looks like
\[g = -2\Omega^2(du\otimes dv + dv\otimes du) + x^{-2}\slash{k}_{ij}(d\theta^i-xf^idv)\otimes(d\theta^j-xf^jdv).\]
This is of the same form as the rescaled Minkowski metric after performing the parabolic blowup in \cref{C:C3:ansatz}.

A short-pulse double-null gauge breaks the diffeomorphism invariance for the same reason as an ordinary double-null gauge does.

Unlike \cref{thm:C2:conversion}, it may not be possible to being a short-pulse metric into double-null gauge, even with the correct assumptions on the initial surfaces. This is because in our definition of a short-pulse metric, we did not specify the ``Riemannian part'' to be on the fibres $S$ while the ``Lorentzian part'' should be on the fibres $\mathcal R$ at $M_0$. While one may restrict the class further, we will not have a need for this. However, we will need the following lemma, which guarantees that perturbations of metrics in a double-null gauge may always be brought back into a double-null gauge over almost the entire manifold we started off with, at least for small $x$. We prove it in the case that $\mathcal R = [0,a)\times[0,b)$ for some $a,b > 0$ (if $\mathcal R$ does not have boundary, the proof is similar but slightly easier), and set $\mathcal H_1 = u^{-1}(0)$, $\mathcal H_2 = v^{-1}(0)$. Let $h \in \Sym^2(\spT)$ be any symmetric tensor which is $0$ on $\mathcal H_1 \un \mathcal H_2$. We will consider $xh$ as a perturbation. Then:
\begin{lem}\label{thm:C3:SPconversion}Let $(\mathcal M,g,x,u,v,N)$ be a short-pulse double-null gauge, and set $g_x = g+xh$. Fix $a' < a$ and $b' < b$. Then there exists an open set $U \subseteq \mathcal M$ containing $S_{0,0,0}$, $\epsilon > 0$ and a unique smooth map $\phi:U \to \mathcal M$, a diffeomorphism onto its image, the identity on $(\mathcal H_1 \un \mathcal H_2)\n U$ and $\{x = 0\} \n U$, with $\mathcal M([0,a')\times[0,b'))\n \{0 \leq x <\epsilon\}$ in the range of $\phi$, and satisfying $\phi^\ast x = x$, such that $\phi_{\ast}g_x$ is in a double-null gauge with $x,u,v,N$.\end{lem}
The proof of this lemma is quite technical and involved, and the arguments are somewhat orthogonal to the main body of this work, so we will prove it in \cref{C:A2:dumblemma}.

\subsection{The characteristic initial value problem for short pulse metrics}
\label{C:C3:SCIVP}

Let $(\mathcal M,x,u,v)$ be a doubly-foliated manifold, let $N$ be a short-pulse vector field which is a lift of $\pa_u$, and assume that $\mathcal M = \mathcal M([0,a)\times [0,b))$ for some $a,b > 0$. Set $\mathcal H_1 = \{v=0\}$ and $\mathcal H_2 = \{u = 0\}$.
We want to prove a result analogous to \cref{thm:C2:CIVP}. A preliminary version of the theorem we want to prove is:

\begin{thm}\label{thm:C3:SPCIVPbad}Let $(\mathcal M,x,u,v)$ and $N$ as above. Suppose the following data are given:
\begin{enumerate}[label = (\roman*)]
\item a short-pulse vector field $L$, a lift of $\pa_v$ to ${}^\sp T\mathcal H_2$;
\item a short-pulse fibre Riemannian metric $\slash{g}$ on $S_{0,0}$;
\item a section $\hat{\slash{g}}$ of $\Sym^2({}^{\sp}T^\ast S)$ over $\mathcal H_1\un \mathcal H_2$ which has Riemannian signature and satisfies $\hat{\slash{g}} = \slash{g}$ at $S_{0,0}$ such that $\Lie_L \hat{\slash{g}}$ vanishes at $x=0$ along $\mathcal H_2$;
\item a pair of smooth functions $f_1$, $f_2$ on $S_{0,0}$ with $f_2$ vanishing at $x=0$;
\item a short-pulse fibre vector field $W$ on $S_{0,0}$ which vanishes at $x=0$.
\end{enumerate}Then there exist $0 < u^\ast \leq a$, $0 < v^\ast \leq b$ and a unique short-pulse Lorentzian metric $g$ on $\mathcal M' = \mathcal M([0,u^\ast)\times[0,v^\ast))$ satisfying $\Ric(g) = 0$ such that $(\mathcal M',g,x,u,v,N)$ is a double-null gauge and:
\begin{enumerate}[label = (\roman*)]
\item $L|_{\mathcal{H}_2} = -2\Omega^2\grad u$;
\item $\iota^\ast g = \slash{g}$ over $S_{0,0}$;
\item $\iota^\ast g$ lies in the conformal class of $\hat{\slash{g}}$ over $\mathcal H_1 \un \mathcal H_2$;
\item $f_1$ is the mean curvature of $g$ on $S_{0,0}$ in the direction of $N$, and $f_2$ is the mean curvature of $g$ on $S_{0,0}$ in the direction of $L$\footnote{The mean curvature is well-defined down to $x=0$ since the contraction of a type $(0,2)$ short-pulse fibre tensor (the second fundamental form) with a type $(2,0)$ short-pulse fibre tensor (the metric) is a smooth function.};
\item $W = [N,L]$ on $S_{0,0}$;
\item $\Omega = 1$ on $\mathcal H_1 \un \mathcal H_2$.
\end{enumerate}
\end{thm}
This theorem is not what we want, since it is still only a local existence theorem. As in the ordinary characteristic initial value problem, the equations $\Ric_{NN}(g) = 0$ and $\Ric_{LL}(g) = 0$ provide obstructions to solving in Taylor series on $\mathcal H_1 \un \mathcal H_2$, since they are nonlinear ODEs for the conformal factor between $\iota^\ast g = \slash{g}$ and $\hat{\slash{g}}$ on $\mathcal H_1$ and $\mathcal H_2$, respectively. This is not different in the short-pulse case (as we will see in the next subsection). We will assume by fiat that these equations are solvable, so we make the following definition.
\begin{defn}\label{def:C3:regular}We will call initial data as stated in \cref{thm:C3:SPCIVPbad} \emph{regular} if we may solve for the conformal factor on all of $\mathcal H_1$ and $\mathcal H_2$, i.e.\ there exists a smooth function $\Phi$ defined on all of $\mathcal H_1$ and $\mathcal H_2$, such that, with $\slash{g} = e^{\Phi}\hat{\slash{g}}$, $\Phi$ satisfies
\begin{align*}
2N^2\Phi + (N\Phi)^2 + \Tr(\hat{\slash{g}}^{-1}\Lie_N\hat{\slash{g}})N\Phi + \frac{1}{2}|\Lie_N\hat{\slash{g}}|^2_{\hat{\slash{g}}} + N\Tr(\hat{\slash{g}}^{-1}\Lie_N\hat{\slash{g}})&= 0\\
2L^2\Phi + (L\Phi)^2 + \Tr(\hat{\slash{g}}^{-1}\Lie_L\hat{\slash{g}})L\Phi + \frac{1}{2}|\Lie_L\hat{\slash{g}}|^2_{\hat{\slash{g}}} + L \Tr(\hat{\slash{g}}^{-1}\Lie_L\hat{\slash{g}})&= 0,\end{align*}
on $\mathcal H_1, \mathcal H_2$, respectively, with initial data for $\Phi|_{S_{0,0}} = 0$ and initial data for $N\Phi|_{S_{0,0}}$ and $L\Phi|_{S_{0,0}}$ chosen to ensure that
\begin{align*}f_1 &= \frac{1}{2}\slash{\tr}_{\slash{g}}\Lie_N \slash{g}|_{S_{0,0}} = N\Phi|_{S_{0,0}} + \frac{1}{2}\slash{\tr}_{\hat{\slash{g}}}\Lie_N \hat{\slash{g}}|_{S_{0,0}}\\
f_2 &= \frac{1}{2}\slash{\tr}_{\slash{g}}\Lie_L \slash{g}|_{S_{0,0}} = L\Phi|_{S_{0,0}} + \frac{1}{2}\slash{\tr}_{\hat{\slash{g}}}\Lie_L \hat{\slash{g}}|_{S_{0,0}}.\end{align*}
\end{defn}
This is the only obstruction. Thus we have:
\begin{thm}\label{thm:C3:SPCIVP}Provide data as in \cref{thm:C3:SPCIVPbad}, and assume they are regular. Then \emph{for any} $u' < a$ and $v' < b$, there exists $\delta > 0$ and a unique short-pulse Lorentzian metric $g$ on $\mathcal M' = \mathcal M([0,u')\times[0,v'))\n\{x \leq \delta\}$ satisfying $\Ric(g) = 0$ such that $(\mathcal M',g,x,u,v,N)$ is a double-null gauge and has the correct initial data in the sense of \cref{thm:C3:SPCIVPbad}.\end{thm}

Uniqueness follows immediately from \cref{thm:C2:CIVP}, since for fixed $x$, $g|_{M_x}$ solves the usual characteristic initial value problem, so we need only treat existence.

The rest of the proof of \cref{thm:C3:SPCIVP} follows the same outline as the proof of \cref{thm:C2:CIVP}. Namely, we first find a solution $g$ in Taylor series, both at $\mathcal H_1 \un \mathcal H_2$ and also at $x = 0$. Then we perturb, looking for a perturbation $h$ such that $g+h$ is in wave gauge with respect to $g$ and solves $\Ric(g+h) = 0$. Afterwards, we may switch back to a double-null gauge. 

The solution in wave gauge will be obtained similarly to as in the proof of \cref{thm:C2:existence}. The main tool used in the proof was the solution of the Cauchy problem for quasilinear hyperbolic PDE. We will need to examine the Cauchy problem for wave equations if the background Lorentzian metric is a short-pulse metric. The idea is that since $h$ vanishes to infinite order at $x=0$, the non-linear terms in the equation will be effectively negligible, which guarantees a solution for long enough time, provided that $x$ is small enough.

More precisely, we will first prove:
\begin{thm}\label{thm:C3:SPTSI}Let regular initial data be given as in \cref{thm:C3:SPCIVP}. Then there exists a short-pulse Lorentzian metric $g$ on $M([0,a)\times[0,b))$ satisfying $\Ric(g) = 0$ in Taylor series at $\mathcal H_1\un \mathcal H_2$ and at $\{x=0\}$ such that $(\mathcal M,g,u,v,N)$ is a double-null gauge and induces the initial data in the sense of \cref{thm:C2:CIVP}. If $\widetilde{g}$ is any other such metric satisfying the same then $g=\widetilde{g}$ in Taylor series at all boundaries. Moreover, there is an iterative algorithm to compute the Taylor series expansions at $\mathcal H_1$, $\mathcal H_2$ and $\{x = 0\}$.\end{thm}

Now we will try to find a perturbation $h$ of $g$ so $g+h$ solves the reduced Einstein equations $\widetilde{\Ric}(g+h) = 0$. Let us denote $Q(h) = \widetilde{\Ric}(g+h)$. Then the major existence result we need to prove is:
\begin{thm}\label{thm:C3:SPexistence}For any $u' < a$ and $v' < b$, there exists $\delta > 0$ and a unique smooth, symmetric, type $(0,2)$ short-pulse tensor $h$ which is $0$ in Taylor series at $\mathcal H_1 \un \mathcal H_2$ and $x=0$ and which solves $Q(h) = 0$ on $\mathcal M([0,u']\times[0,v'])\n\{x \leq \delta\}$.\end{thm}

Choose $\delta$ smaller, $g+h$ is a short-pulse Lorentzian metric. For fixed $x > 0$, $(g+h)|_{M_x}$ is an ordinary Lorentzian metric solving the reduced Einstein equations, and the wave-gauge condition is satisfied on the initial hypersurfaces. Thus by \cref{thm:C2:wavegauge}, $(g+h)|_{M_x}$ is in wave gauge on its entire domain of existence, and thus $\Ric(g+h) = 0$ over $x > 0$, and thus by continuity, $\Ric(g+h) = 0$ for all $0 \leq x \leq \delta$. Now using \cref{thm:C3:SPconversion} we may convert back to double-null gauge and conclude \cref{thm:C3:SPCIVP}.

In the remaining subsections of \cref{C:C3:SPCharValue} we prove \cref{thm:C3:SPTSI}. Since the proof of \cref{thm:C3:SPexistence} will require an understanding of the short-pulse Cauchy problem, we will devote \cref{C:C3:Cauchy} to this undertaking, and complete the proof in \cref{C:C3:completetheproof} after the necessary machinery has been developed. 

\subsection{Ricci curvature in a double null gauge}
Let $(\mathcal M,g,x,u,v,N)$ be a double-null gauge. Let $\Omega$, $\slash{g}$, $L$ be the quantities associated to $\slash{g}$ in a double-null foliation, i.e.\ $-2\Omega^{-2} = g^{-1}(2du,2dv)$, $L = -2\Omega^2 \grad v$, $g = \iota^\ast \slash{g}$. While the equations in \cref{C:C2:RicciDNG} still hold for short-pulse metrics in a parametrized double-null gauge, the scaling is not as evident, so we rescale to better see the scaling. Define $\mathring{L}$ by setting $\mathring{L} = L$ on $\mathcal H_2$ and extending it to all of $\mathcal M$ by requiring $[N,\mathring{L}] = 0$. We may write $L = \mathring{L} + x\overline{L}$, where $\overline{L}$ is an ordinary fibre vector field. We may also write $\slash{g} = x^{-2}\slash{h}$, where $\slash{h}$ is an ordinary fibre metric. Let us also introduce the shorthand $Z = x^{-1}\slash{h}([N,L],\bullet) = \slash{h}([N,\overline{L}],\bullet)$, which is an ordinary fibre one-form.

Keep the notation as in \cref{C:C2:RicciDNG}, except take all slash tensorial operations with respect to $\slash{h}$. Observe that since $T^\ast S \subseteq {}^\sp T^\ast S$, $\slash{d}$ maps functions to short-pulse one-forms. Rescaling \eqref{eq:C2:one}-\eqref{eq:C2:six} we have that
\begin{subequations}
\begin{align}
\begin{split}
\label{eq:C3:one}2\Omega^2x^2\Ric_{-,-} &= \Lie_N\Lie_L\slash{h} - \frac{x}{2}\Lie_{[N,\overline{L}]} \slash{h} + \frac{1}{4}(\slash{\tr}\Lie_N \slash{h})\Lie_L\slash{h} + \frac{1}{4}(\slash{\tr}\Lie_L\slash{h})\Lie_N \slash{h}\\
&-\frac{1}{2}(\Lie_L\slash{h} \times \Lie_N \slash{h} + \Lie_N\slash{h}\times \Lie_L\slash{h})\\
&\hspace{-4pt}+2x^2\Omega^2\slash{\Ric} - 2x^2\slash{\Hess}(\Omega^2) + 2x^2\Omega^2(\slash{d}\log \Omega)\otimes(\slash{d}\log\Omega)\\
&- \frac{1}{4\Omega^2}Z\otimes Z\end{split}\\
\begin{split}
\label{eq:C3:two}4x\Omega^2\Ric_{N,-} &= \Lie_N Z + 2x\Omega^2\slash{\divg}(\Lie_N \slash{h})-2x\Omega^2\slash{d}\slash{\tr}\Lie_N\slash{h}\\
&+ \frac{1}{2}(\slash{\tr}\Lie_N\slash{h})Z + x\slash{\tr}\Lie_N\slash{h}\slash{d}\Omega^2 - x4\Omega^2\slash{d}N\log\Omega - 2(N\log\Omega)Z
\end{split}\\
\begin{split}
\label{eq:C3:three}-4x\Omega^2\Ric_{L,-} &= \Lie_L Z-2x\Omega^2\slash{\divg}(\Lie_L\slash{h})+2x\Omega^2\slash{d}\slash{\tr}\Lie_L\slash{h}\\
&+\frac{1}{2}(\slash{\tr}\Lie_L\slash{h})Z- x\slash{\tr}\Lie_L\slash{h}\slash{d}\Omega^2 + 4x\Omega^2\slash{d}L\log\Omega - 2(L\log\Omega)Z
\end{split}\\
\begin{split}
\label{eq:C3:four}-2\Ric_{NN} &= N\slash{\tr}\Lie_N \slash{h} + \frac{1}{2}|\Lie_N \slash{h}|^2 - 2(N\log\Omega)\slash{\tr}\Lie_N\slash{h}\end{split}\\
\begin{split}
\label{eq:C3:five}-2\Ric_{LL} &= L\slash{\tr}\Lie_L\slash{h} + \frac{1}{2}|\Lie_L\slash{h}|^2 - 2(L\log\Omega)\slash{\tr}\Lie_L\slash{h}\end{split}\\
\begin{split}
\label{eq:C3:six}-\frac{1}{2}\Ric_{NL} &=NL\log\Omega + \frac{1}{4}N\slash{\tr}\Lie_L\slash{h} - \frac{x}{8}\slash{\tr}\Lie_{[N,\overline{L}]} \slash{h} - \frac{x}{2}[N,\overline{L}]\log\Omega- \frac{x^2}{2}\slash{\Delta}\Omega^2\\
&+\frac{1}{8}\slash{h}(\Lie_L\slash{h},\Lie_N \slash{h}) + \frac{1}{8\Omega^2}|[N,\overline{L}]|^2.
\end{split}
\end{align}
\end{subequations}

Observe that we are plugging in vectors in $TS$ for the $-$, not vectors in ${}^\sp TS$. Thus the factors of $x$ out front are consistent with $\Ric(g)$ being a smooth short-pulse tensor.

We may also rescale the integrability conditions \eqref{eq:C2:onei}-\eqref{eq:C2:threei}. Denote by $\chi^N$, $H^N$, $\chi^L$, $H^L$ the second fundamental forms and mean curvature of $\slash{h}$ in the directions of $N$ and $L$, respectively. Then (all slash operations still taken with respect to $\slash{h}$):
\begin{subequations}
\begin{align}
\label{eq:C3:onei}
\begin{split}
x\Lie_L\Ric_{N,-}+x\Lie_N \Ric_{L,-} &+xH^L\Ric_{N,-} + xH^N\Ric_{L,-}+2\Omega^2x^2\slash{\div}\Ric_{-,-}-\Omega^2x^2\slash{d}\slash{\tr}\Ric_{-,-}\\
&\hspace{4.7cm}+2x^2\Ric_{-,-}\cdot \slash{\grad}\Omega^2+\slash{d}\Ric_{NL}.
\end{split}\\
\label{eq:C3:twoi}
\begin{split}
L\Ric_{NN} + H^L\Ric_{NN}&=-\Omega^2x^2N\slash{\tr}\Ric_{-,-}-2\Omega^2x^4\slash{h}(\chi^N, \Ric_{-,-})-H^N\Ric_{NL}\\
&+2\Omega^2x^2\slash{\div}\Ric_{N,-}+2x^2\Ric_{N,-}\cdot (\slash{\grad}\Omega^2) - \Ric_{N,-}\cdot [N,L].\end{split}\\
\label{eq:C3:threei}
\begin{split}
N\Ric_{LL} + H^N\Ric_{LL}&=-\Omega^2x^2L\slash{\tr}\Ric_{-,-}-2\Omega^2x^4\slash{h}(\chi^L, \Ric_{-,-})-H^L\Ric_{LN}\\
&+2\Omega^2x^2\slash{\div}\Ric_{L,-}+2x^2\Ric_{L,-}\cdot (\slash{\grad}\Omega^2) + \Ric_{L,-}\cdot [N,L].\end{split}
\end{align}
\end{subequations}

We use the notation $O(x^j)$ to denote a quantity equal to $x^j$ times a smooth section of an appropriate bundle. The notation $O(x^\infty)$ indicates that for all $j$ quantity is in $O(x^j)$. Keep the notation $O(u^j)$, $O(v^j)$, $O(u^\infty)$, and $O(v^{\infty})$ from \cref{C:C2:RicciDNG}.

It is not difficult to see that if $g$ is a metric whose associated $\Omega$, $\slash{g}$, $L$ have the right-hand side of each of the above being in $O(x^k)$, then $\Ric(g)$, considered as a tensor on the rescaled bundle, is also in $O(x^k)$. The converse also holds.

\subsection{Solving in Taylor series}
\label{C:C3:SPTaylor}
In this section we will prove \cref{thm:C3:SPTSI}. Set $\omega = \log \Omega$. Instead of solving for $g$, we will solve for $\phi = (\slash{h},\overline{L},\omega)$, a section of $\Sym^2 T^\ast S \oplus TS \oplus \R$. We need to find such a $\phi$ which solves \eqref{eq:C3:one}--\eqref{eq:C3:six} in Taylor series, since then we may find the solution $g$ in Taylor series using \cref{thm:C3:ReconstructSP}. For the uniqueness of the series, it suffices to establish the uniqueness of $\phi$ since every solution $g$ of $\Ric(g) = 0$ gives a solution $\phi$ of \eqref{eq:C3:one}--\eqref{eq:C3:six} and vice-versa.

First, let us deal with the integrability conditions.
\begin{lem}\label{thm:C3:integrability}Suppose $\phi$ is given such that \eqref{eq:C3:one}, \eqref{eq:C3:three}, \eqref{eq:C3:six} are $0$ in Taylor series at $\mathcal H_1 \un \mathcal H_2$ and $\{x = 0\}$, and \eqref{eq:C3:two},\eqref{eq:C3:four} are identically $0$ on $\mathcal H_1$ and \eqref{eq:C3:five} is identically $0$ on $\mathcal H_2$. Then all equations \eqref{eq:C3:one}-\eqref{eq:C3:six} are $0$ in Taylor series at $\mathcal H_1 \un \mathcal H_2$ and at $\{x = 0\}$.\end{lem}
\begin{proof}
As in \cref{thm:C2:integrability}, the conditions are sufficient to be $0$ in Taylor series at $\mathcal H_1 \un \mathcal H_2$. Observe that $x\Ric_{N,-}$ is a smooth (ordinary!) one-form, since $\Ric_{N,-}$ is a short-pulse one-form. From \eqref{eq:C3:onei}, $x\Ric_{N,-}$ satisfies an equation of the form
\begin{equation}\label{eq:C3:integrabilityeq}\Lie_{\mathring{L}}x\Ric_{N,-} + H^Lx\Ric_{N,-} + \Lie_{\overline{L}}x\Ric_{N,-} \in O(x^\infty).\end{equation}

Equation \eqref{eq:C3:integrabilityeq} implies that $x\Ric_{N,-} = 0$ at $x=0$, since it is $0$ on $\mathcal H_1$ by assumption. Now, we may differentiate \eqref{eq:C3:integrabilityeq} $k$ times in $x$ to obtain
\[\Lie_{\mathring{L}}\pa_x^k x\Ric_{N,-} + H^L \pa_x^k x\Ric_{N,-} + \Lie_{\overline{L}} \pa_x^k x\Ric_{N,-} = \ldots + O(x^\infty),\]
where the $\ldots$ only products of known quantities with $\pa_x^j x\Ric_{N,-}$ and its derivatives in fibre directions. Thus proceeding inductively, we deduce that $\pa_x^k x\Ric_{N,-} = 0$ for all $k$. With $x\Ric_{N,-} \in O(x^\infty)$ known, the same argument applied to \eqref{eq:C3:twoi} and \eqref{eq:C3:threei} works for the other two quantities.
\end{proof}

Arguing exactly as in the proof of \cref{thm:C2:TSI} (and using the \cref{thm:C2:TSbetter}), but treating $x$ as a parameter, we know that there is a unique solution in Taylor series $\psi$ which solves \eqref{eq:C3:one}-\eqref{eq:C3:six} in Taylor series on $\mathcal H_1\un \mathcal H_2$. So we just need to construct $\phi$, the Taylor series at $x=0$ and prove the uniqueness of the Taylor series there.

\begin{lem}\label{thm:C3:toporder}There exist a unique section $\phi_0$ defined on $M_0([0,a)\times[0,b))$ such that any $\phi$ smooth with $\phi|_{x=0} = \phi_0$ satisfies \eqref{eq:C3:one}, \eqref{eq:C3:three}, \eqref{eq:C3:six} $=0$ at $x=0$, \eqref{eq:C3:two}, \eqref{eq:C3:four} $=0$ on $\mathcal H_1 \n \{x=0\}$, and \eqref{eq:C3:five} $=0$ on $\mathcal H_2 \un \{x=0\}$ and has the correct initial data. In fact, if $\phi_0 = (\slash{h}_0,\overline{L}_0,\omega_0)$, then $\overline{L}_0 = 0$, $\omega_0 = 0$.\end{lem}
\begin{proof}To obtain $\phi_0$, we propagate from $S_{0,0,0}$ to $(\mathcal H_1 \un \mathcal H_2)\n\{x=0\}$ and finally to all of $\{x=0\}$. By assumption the data are regular, so \eqref{eq:C3:four}, \eqref{eq:C3:five}$=0$ are satisfied on $\mathcal H_1$ and $\mathcal H_2$, respectively. \eqref{eq:C3:two} ,\eqref{eq:C3:three}$=0$ at $x=0$ reduce to homogeneous linear ODEs for $Z$ on $\mathcal H_1$ and $\mathcal H_2$, respectively. Since $[N,L]$ vanishes as a section of $^\sp TS$ at $S_{0,0}$ by assumption, $[N,\overline{L}]$ vanishes at $x=0$ as a section of $TS$, and hence $Z$ does too. It follows that $Z$ is $0$ along $\mathcal H_1 \un \mathcal H_2$. Since we know $Z$ and $\slash{h}$, we know $[N,\overline{L}]$. In particular, since $\overline{L} = 0$ on $S_{0,0}$ by assumption, $\overline{L} = 0$ on $\mathcal H_1$.

Now at $x=0$, $Z=0$ certainly solves \eqref{eq:C3:two} $= 0$ on all of $M_0$ and has the correct initial data (i.e.\ $Z=0$ on $\mathcal H_2$, which is true by assumption). If $\phi$ is any smooth solution, then this equation is linear in $Z$ (since everything else is known), and so $Z=0$ is the unique such solution. 

If $\Phi$ is the function such that $\slash{g} = e^{\Phi}\hat{\slash{g}}$, then along $\mathcal H_2$ \eqref{eq:C3:five} $ =0$ becomes
\[2L^2\Phi + (L\Phi)^2 + \Tr(\hat{\slash{g}}^{-1}\Lie_L\hat{\slash{g}})L\Phi + \frac{1}{2}|\Lie_L\hat{\slash{g}}|^2_{\hat{\slash{g}}} + L \Tr(\hat{\slash{g}}^{-1}\Lie_L\hat{\slash{g}}) .\]
By assumption $\Lie_L \hat{\slash{g}}$ vanishes at $x=0$ on $\mathcal H_2$, and therefore the equation reduces to
\[2L^2\Phi + (L\Phi)^2 = 0.\] By assumption $\Phi = 0$ at $S_{0,0,0}$.

Also, \[\Lie_L \slash{g} = e^{\Phi}(L\Phi\hat{\slash{g}} + \Lie_L\hat{\slash{g}}),\] and by assumption the second term vanishes at $x=0$. Thus $L\Phi$ is the mean curvature at $S_{0,0,0}$, which is $0$ by assumption. Since $\Phi = 0$ at $S_{0,0,0}$ by assumption, it follows that $\Phi \equiv 0$ is the only solution, and hence $\Lie_L \slash{g} = 0$ on $\mathcal H_2$ at $x=0$ by the previous formula.

Now let us look at \eqref{eq:C3:one} $=0$. One easily see that any $\slash{h}$ with $\Lie_L \slash{h} = 0$ solves the equation at $x=0$. As for uniqueness, if there is a smooth solution $\slash{h}$, then treating $\slash{h}$ and $\Lie_N \slash{h}$ as known, the equation is linear in $\Lie_L \slash{h}$, and thus must be $0$. Since $\Lie_L\slash{h} = 0$, we may find $\slash{h}$ from integrating out from $\mathcal H_1$ where we have found initial data.

Since $Z=0$, and $\overline{L} = 0$ on $\mathcal H_2$, $\overline{L} = 0$ on $M_0$. Thus \eqref{eq:C3:six}$=0$ has unique solution $\omega = 0$ with the correct initial data (i.e.\ $\omega = 0$ initially).\end{proof}

Let $P = (P^1,P^2,P^3)$ denote the operator taking in a section $\phi$ and given the right-hand side of \eqref{eq:C3:one}, \eqref{eq:C3:three}, \eqref{eq:C3:six}, respectively. As in the proof of \cref{thm:C2:TSI}, we will need to linearize. Suppose $\phi_0$ is a section of $\Sym^2 T^\ast S \oplus TS \oplus \R$, and write $D_{\phi} P$ for the linearization of $P$ at $\phi$. Then:
\begin{lem}\label{thm:C3:linearization}Suppose a section $\phi$ of $\Sym^2(T^\ast S)\otimes TS \otimes \R$ satisfies $\phi|_{x=0} = (\slash{h}_0,0,0)$, with $\Lie_{\mathring{L}}\slash{h}_0 = 0$. Let $\xi = (\slash{h},\overline{L},\omega)$ be any other section. Then the components of $DP_{\phi|_{x=0}}\xi$ are (with all tensorial operations taken with respect to $\slash{h}_0$):
\begin{subequations}
\begin{align}
\begin{split}
&\Lie_N\Lie_{\mathring{L}}\slash{h} + \frac{1}{4}(\slash{\tr}\Lie_N \slash{h}_0)\Lie_{\mathring{L}} \slash{h} + \frac{1}{4}(\slash{\tr}\Lie_{\mathring{L}} \slash{h})\Lie_N \slash{h}_{0}\\
&-\frac{1}{2}(\Lie_{\mathring{L}} \slash{h} \times \Lie_N \slash{h}_0 + \Lie_N \slash{h}_0 \times \Lie_{\mathring{L}} \slash{h}))
\end{split}\\
\begin{split}
&\Lie_{\mathring{L}}\slash{h}_0([N,\overline{L}],\bullet)
\end{split}\\
\begin{split}
&N\mathring{L}\omega + N\frac{1}{4}\slash{\tr}\Lie_{\mathring{L}}\slash{h} + \frac{1}{8}\slash{h}_0(\Lie_{\mathring{L}}\slash{h},\Lie_N \slash{h}_0).\end{split}\end{align}
\end{subequations}
\end{lem}

Now we may iteratively solve for the complete Taylor series at $x=0$. Let us first mention that if the compatibility conditions are to hold, then \[\phi|_{\mathcal H_1 \un \mathcal H_2} = \psi|_{\mathcal H_1 \un \mathcal H_2} + O(x^\infty),\] and so the top-order behaviour of $\psi$ tells us the initial data for each term in the expansion of $\phi$.

Suppose we are given for $j \geq 1$ and $0 \leq \ell \leq j$ sections $\phi_{\ell}$ which are defined over $M_0$, and extend them to all of $\mathcal M$ by requiring $\pa_x \phi_{\ell} = 0$. Write
\[\phi_{(j)} = \phi_0 + \cdots + x^{j}\phi_{j},\]
and assume $\phi_{(j)}$ solves \eqref{eq:C3:one}, \eqref{eq:C3:three}, \eqref{eq:C3:six} $=0$ mod $O(x^{j+1})$ and induces the correct initial data mod $O(x^{j+1})$. We seek to find $\phi_{j+1}$:

\begin{lem}\label{thm:C3:higherorder}
There exists a unique section $\phi_{j+1}$ defined over $M_0$ such that, after being to $\mathcal M$ via $\pa_x \phi_{j+1} = 0$, $\phi_{(j)} + x^{j+1}\phi_{j+1}$ induces the correct initial data mod $O(x^{j+2})$ and satisfies the equations mod $O(x^{j+2})$.\end{lem}
\begin{proof}
With $P$ as above,
\[P(\phi_{(j)} + x^{j+1}\phi_{j+1}) = P(\phi_{(j)}) + x^{j+1}DP_{\phi_{(j)}}\phi_{j+1} + x^{j+2}M_{\phi_{(j)}}(\phi_{j+1}),\] where $DP$ denotes the linearization of $P$ and $M$ is some non-linear operator. By assumption, $P(\phi_{(j)}) \in O(x^{j+1})$, so this sets up a linear equation for $\phi_{j+1}$:
\[DP_{\phi_{(j)}}\phi_{j+1} = -(x^{-(j+1)}P(\phi_{(j)}))_{x=0},\]

\Cref{thm:C3:toporder} tells us that our $\phi_{(j)}$ fits into the assumptions of the previous lemma, so we know the linearizations $DP_{\phi_{(j)}}|_{x=0}$.

Now let us mention what the data for $\phi_{j+1}$ must be. Expand $\psi$ as a Taylor series at $\{x = 0\}$ on $\mathcal H_1 \un \mathcal H_2$:
\[\psi \sim \psi_0 + x\psi_1 + \cdots.\]
By assumption, we know that $\phi_{(j)}$ agrees with $\psi$ (recall that $\psi$ is the Taylor series solution at $\mathcal H_1 \un \mathcal H_2$) to order $O(x^{j+1})$, and thus we must set $\phi_{j+1} = \psi_{j+1}$ on $\mathcal H_1 \un \mathcal H_2$ in order to induce the correct initial data.

With this correct initial data, we proceed to solve
\[DP_{\phi_0}\phi_{j+1} = -(x^{-(j+1)}P(\phi_{(j)}))_{x=0}.\]
From the form of the linearization, the equation for the first component is a linear transport equation for $\Lie_{\mathring{L}} \slash{h}_{j+1}$, which we solve, using the initial data for $\slash{h}_{j+1}$ on $\mathcal H_2$. Using this, we can then solve for $\slash{h}_{j+1}$ itself. We may also use the equation for the second component to solve for $\slash{h}_0([N,\overline{L}_{j+1}],\bullet)$, using that we know it on $\mathcal H_1$. Since we know $\slash{h}_0$ we know $[N,\overline{L}_{j+1}]$, and using the data for $\overline{L}_{j+1}$ on $\mathcal H_2$, $\overline{L}_{j+1}$ itself. Lastly, we may solve for $\omega_{j+1}$ using the same technique.
\end{proof}

Using \cref{thm:C3:toporder}, \cref{thm:C3:higherorder} and Borel's lemma, we have a section $\tilde{\phi}$ which is a solution to \eqref{eq:C3:one}, \eqref{eq:C3:three}, \eqref{eq:C3:six} $=0$ in Taylor series at $\{x = 0\}$ and which is unique in Taylor series. We wish to find a function $\phi$ whose series expansion at $\{x = 0\}$ coincides with that of $\tilde{\phi}$ and whose series expansions at $\mathcal H_1$ and $\mathcal H_2$, respectively, coincides with those of $\psi$. We only we need the compatibility condition to be satisfied i.e.\ the $\pa_x$, $\pa_u$, and $\pa_v$ derivatives of $\tilde{\phi}$, and $\psi$ must coincide at $(\mathcal H_1 \un \mathcal H_2)\n \{x=0\}$.

\begin{lem}\label{thm:C3:compat}The compatibility conditions are satisfied.\end{lem}
\begin{proof}Expand $\psi$ in Taylor series
\[\psi \sim \psi_0 + x\psi_1 + x^2\psi_2 + \cdots\]
and $\tilde{\phi}$ in Taylor series
\[\tilde{\phi} \sim \phi_0 + x\phi_1 + x^2\phi_2 + \cdots.\]
Let us also write $\psi_{(j)}$ and $\phi_{(j)}$ for the sum of the first $j$ terms.

The compatibility condition is equivalent to $\psi_j = \phi_j$ in Taylor series at $\mathcal H_1 \un \mathcal H_2$. We will use induction on $j$.

Let us treat the base case $j = 0$. By construction, at $\{x=0\}$, $P\phi_0 \equiv 0$ and $P\psi_0 = 0$ in Taylor series at $\mathcal H_1 \un \mathcal H_2$, and both $\phi_0$ and $\psi_0$ satisfy \eqref{eq:C3:two}, \eqref{eq:C3:three}$=0$ on $\mathcal H_1\n \{x=0\}$ and \eqref{eq:C3:five}$=0$ on $\mathcal H_2\n \{x=0\}$. Since both $\psi_0$ and $\phi_0$ induce the same initial data, either by using versions of \cref{thm:C2:toporder} and \cref{thm:C2:higherorder}, or carrying through a similar argument to that of \cref{thm:C3:toporder}, it follows that $\phi_0 = \psi_0$ in Taylor series.

Now let us use induction, assuming that $\phi_{i}-\psi_i$ vanishes in Taylor series at $\mathcal H_1 \un \mathcal H_2$ for $i \leq j$, $j > 0$. We may write $\psi-\psi_{(j)} = x^{j+1}\psi_{j+1} + x^{j+2}\tilde{\psi}$, for some smooth $\tilde{\psi}$.

Let us write
\[P(\psi) = P(\psi_{(j)} + (\psi-\psi_{(j)})) = P(\psi_{(j)} + x^{j+1}\psi_{j+1} + x^{j+2}\tilde{\psi}).\]
Thus
\begin{equation}\label{eq:compat:eq1}P(\psi) = P(\psi_{(j)}) + x^{j+1}DP_{\psi_{(j)}}(\psi_{j+1}) + O(x^{j+2}).\end{equation}
Since each term in $\psi_{(j)}$ is equal to each term in $\phi_{(j)}$ up to a quantity vanishing in Taylor series at $\mathcal H_1 \un \mathcal H_2$,
\begin{equation}\label{eq:compat:eq2}P(\psi_{(j)}) = P(\phi_{(j)}) + \sum_{i=0}^{j+1} x^i \kappa_i + O(x^{j+2}),\end{equation} where $\kappa_i$ does not depend on $x$ and vanishes in Taylor series at $\mathcal H_1 \un \mathcal H_2$. Also, by assumption $P(\psi)$ is rapidly vanishing in Taylor series at $\mathcal H_1 \un \mathcal H_2$. In particular
\begin{equation}\label{eq:compat:eq3}P(\psi) = \sum_{i=0}^{j+1} x^i \lambda_i + O(x^{j+2}),\end{equation} where $\lambda_i$ does not depend on $x$ and vanishes in Taylor series at $\mathcal H_1 \un \mathcal H_2$. Substituting \eqref{eq:compat:eq2} and \eqref{eq:compat:eq3} into \eqref{eq:compat:eq1} and equating the coefficients of $x^{j+1}$, and and taking $x \to 0$, it follows that
\[\lambda_{j+1} = (x^{-j+1}P(\phi_{(j)})|_{\{x=0\}}) + DP_{\psi_0}\psi_{j+1} + \kappa_{j+1}.\]
But $\psi_0 = \phi_0$ in Taylor series, and $(x^{-j+1}P(\phi_{(j)})|_{\{x=0\}}) + DP_{\phi_0}\phi_{j+1} = 0$ by construction, and so we conclude that \[DP_{\phi_0}(\psi_{j+1}-\phi_{j+1})\] vanishes in Taylor series at $\mathcal H_1 \un \mathcal H_2$. Moreover, $\phi_{j+1} - \psi_{j+1} = 0$ on $\mathcal H_1 \un \mathcal H_2$ (since $\tilde{\phi}$ has the correct initial data).

Arguing as in the proof of the previous lemma using the form of $DP_{\phi_0}$, it follows that $\psi_{j+1}-\phi_{j+1}$ also vanishes in Taylor series at $\mathcal H_1 \un \mathcal H_2$. This completes the proof.
\end{proof}
The $\phi$ we have constructed satisfies the hypotheses of \cref{thm:C3:integrability}, since by definition it solves equations \eqref{eq:C3:one}, \eqref{eq:C3:three}, \eqref{eq:C3:six} $=0$ in Taylor series at the correct faces, and solves the other equations identically on the correct hypersurfaces (being equal to $\phi$ at $x=0$ and $\psi$ on $\mathcal H_1 \un \mathcal H_2$). Thus, by \cref{thm:C3:integrability}, we have found our solution $\phi$ in Taylor series. The solution $\phi$ is unique in Taylor series, since any other solution has the same expansions at $\{x = 0\}$, $\mathcal H_1$, and $\mathcal H_2$, and we already know these are unique.

\section{The Cauchy problem for short-pulse metrics}
\label{C:C3:Cauchy}
\subsection{Outline and main result}
\label{C:C3:Cauchyintro}
In this section we study the Cauchy problem for short-pulse tensors, with the main goal to derive a long-time existence result which we can use to prove \cref{thm:C3:SPexistence}. In this section, we will be working with a very special type of doubly-foliated manifold. Let $S$ be a closed manifold, and set $Y = \R_y \times S$, $M = [0,\infty)_t\times Y = [0,\infty)_t\times\R_y\times S$, and as usual $\mathcal M = [0,1)_x\times M$. The coordinate projection $(t,y)$ turns $\mathcal M$ into a doubly-foliated manifold with parameter $x$ over $[0,\infty)\times \R$ with fibre $S$. We choose the diffeomoprhism at $\{x = 0\}$ the identity map. We may consider the subbundle $^\sp TY \subseteq {}^\sp T\mathcal M$ consisting of those vector fields tangent to $Y$, and also its tensor powers
\[^\sp T^p_q Y := ({}^\sp TY)^{\otimes p} \otimes ({}^\sp T^\ast Y)^{\otimes q}.\]

Let us fix $\slash{k}$ a Riemannian metric on $S$, and use it to construct the short-pulse Riemannian metric $\mathbf{e} = dt^2 + dy^2 + x^{-2}\slash{k}$. We may also construct metric $e_x = dy^2 + x^{-2}\slash{k}$ on on $^\sp TY$. Let $D$ denote the Levi-Civita connection of $\mathbf{e}$, which by \cref{thm:C3:connectionisnice} maps sections of $^\sp T^p_q \mathcal M$ to sections of $^\sp T^p_{q+1}\mathcal M$. Also let $\nabla$ denote the Levi-Civita connection of $e$. Since the action of $\nabla$ on sections of $^\sp TY$ is given by the orthogonal projection onto $^\sp TY$ of the the action of $D$, $\nabla$ maps sections of $^\sp T^p_q Y$ to sections of $^\sp T^p_{q+1} Y$.

We wish to study hyperbolic PDE of the form
\begin{equation}\label{eq:C3:type}\phi(h)D^2 h + a(h,Dh)D h + b(h,Dh) h = f,\end{equation}
where $\phi:U_1 \subseteq {}^\sp T^p_q \mathcal M \to \Sym^2({}^\sp T \mathcal M)$ takes values in Lorentzian cometrics, $a,b$ are smooth bundle maps
\begin{align*}
a:U_2 \subseteq {}^\sp T^p_q \mathcal M \oplus{}^\sp T^p_{q+1} \mathcal M \to {}^\sp T^{q+1+p}_{p+q} \mathcal M = \Hom( {}^\sp T^p_{q+1} \mathcal M, {}^\sp T^p_q \mathcal M)\\
b:U_3 \subseteq{}^\sp T^p_q \mathcal M \oplus{}^\sp T^p_{q+1} \mathcal M \to {}^\sp T^{q+p}_{p+q} \mathcal M = \Hom({}^\sp T^p_q \mathcal M, {}^\sp T^p_q \mathcal M),
\end{align*}
where $U_i$ ($i=1,2,3$) are open sets containing the zero section, and $f$ is a smooth type $(p,q)$ short-pulse tensor. Here and throughout this section, we will not explicitly write the variables in the base as arguments of our bundle maps. We wish to emphasize that $a, \ b$ are not linear in the fibre variables; any smooth map on the fibres is permissible (for now).

Equation~\eqref{eq:C3:type} takes the form of a singular family of hyperbolic PDE for $x>0$. With this in mind, we will \emph{not} immediately be looking for solutions $h$ which are smooth short-pulse tensors. Rather we will be looking for each $x>0$ a distinct solution $h_x$, the collection of which for all $x > 0$ form a short-pulse tensor. In particular for each $x$, a section of $T^p_q M_x$ is the same thing as a section of $^{\sp} T^p_q \mathcal M$ defined over $M_x$ (i.e.\ in the pullback bundle) This approach is keeping in line with thinking of ${}^\sp T\mathcal M$ as a sort of semiclassical bundle. In the semiclassical calculus, with $h$ the semiclassical parameter, one is not often concerned with regularity in the $h$ variable, and we will at first adopt this point of view here, too.

We will make the following regularity assumptions on our coefficients $(\phi,a,b)$ which we will call ``admissible.''
\begin{defn}\label{def:C3:admissible}We will call $(\phi,a,b)$ \emph{admissible} if
\begin{enumerate}[label=(\roman*)]
\item $\phi$ is a smooth bundle map from ${}^\sp T^p_q \mathcal M$ to $\Sym^2({}^\sp T\mathcal M)$, which takes values in Lorentzian cometrics for which the surfaces $\{x\}\times \{t\} \times Y$ are spacelike and $\pa_t$ is time-like and future-oriented;
\item $\phi$ is the constant cometric $-\pa_t^2+e^{-1}_x$ outside a compact subset in the base $\mathcal M$;
\item $a,b$ do not depend on the fibre-variables outside a compact subset of the base $\mathcal M$.
\end{enumerate}
\end{defn}

The main theorem we wish to prove in this section is:
\begin{thm}\label{thm:C3:bigtheoremI}Suppose $(\phi,a,b)$ is admissible, and $f$ is $0$ in Taylor series at $x=0$ and is compactly supported. Then for any $T > 0$ there exists $\delta > 0$ and a unique smooth short-pulse tensor $h$ solving \eqref{eq:C3:type} on $\mathcal M\n\{x \leq \delta\} \n\{t \in [0,T]\}$ with $0$ Cauchy data. The solution $h$ is $0$ in Taylor series at $x=0$. \end{thm}
We remark immediately that the usual theory of quasilinear hyperbolic equations gives us uniqueness, and smoothness away from $x > 0$. The content of \cref{thm:C3:bigtheoremI} is that we establish a long-time existence theorem, and prove smoothness down to $0$.

\begin{rk}If we were working over $T\mathcal M/\vspan \ker dx$ rather than $^\sp T\mathcal M$, then the family of equations \eqref{eq:C3:type} would be, rather than a degenerate family of hyperbolic equations, a smooth family of hyperbolic equations. Thus, if the coefficients were ``admissible'' (in an appropriate analogous sense) and the Cauchy data vanishing, the proof of the analogue of \cref{thm:C3:bigtheoremI} would be a routine exercise using perturbative arguments. The power of the short-pulse bundle framework is that it provides a framework that allows us to use almost \emph{exactly the same approach} to prove \cref{thm:C3:bigtheoremI}. In fact, the only thing that needs changing is that instead of using Sobolev spaces constant in $x$, we will need to use varying Sobolev spaces which incorporate a scaling to deal with the degeneracy of $^\sp T\mathcal M$ relative to $TM$. This will, of course, mean that when we use the Sobolev inequalities to control the nonlinearity, our estimates will lose an inverse power $x^{-c}$, $c > 0$, relative to the usual Sobolev inequalities, which depends only on the coefficients and the amount of regularity our Sobolev spaces encode. However, since we expect our solution to be vanishing rapidly in $x$, we may trade in some of this decay to offset the loss of a power of $x$. More explicitly, we will divide equation~\eqref{eq:C3:type} by some large power of $x$, $x^{\ell}$, and instead solve for $h^\ell = x^{-\ell}h$. For instance, in controlling $a(h,Dh)$ we will use an inequality roughly of the form
\[\norm{a(h,Dh)} = \norm{a(x^{\ell}h_\ell,x^\ell Dh_\ell)} \lesssim x^{\ell-c} \leq x,\]
which better than being bounded even decays. It is this decay that lets us prove long-time existence.
\end{rk}

\subsection{Classical energy estimates}
\label{C:C3:energyestimates}
The main estimate in our proof will be classical energy estimates for hyperbolic PDE. The goal of this section will be to carefully track down the constants we need and state the theorem in a language useful to us. For this section, we let $(X,e)$ be a Riemannian manifold and will consider the Cauchy problem on $[0,\infty)_t \times X$. For this subsection, set $\mathbf{e} = dt^2 + e$, a Riemannian metric on $[0,\infty) \times X$. Write $D$ for the Levi-Civita connection of $\mathbf{e}$ and $\nabla$ for the Levi-Civita connection of $e$.

 We will suppose that $(X,e)$ has \emph{bounded geometry}, which for us will mean that $\Riem(e)$, together with all of its $e$-covariant derivatives, is bounded on $X$, and $X$ has a positive injectivity radius.

For non-negative integer $k$, we let $H^k(X)$ be the Sobolev spaces of tensors, whose norm is defined using the metric $e$ via
\[\norm{h}_{H^k(X)} = \sum_{i=0}^k \int_X |\nabla h|^2_e\ d\vol_e,\]
(we will omit the dependence of the space $H^k(X)$ on the type of the tensor from our notation). We will also let $C^k(X)$ be the space of $k$-times continuously differentiable tensors, whose norm is defined via
\[\norm{h}_{C^k(X)} = \sup_{i \leq k}\sup_X |\nabla^i h|_e.\]

We assume that $(X,e)$ satisfies the Sobolev embedding and multiplication inequalities for $u_1$ a type $(p_1,q_1)$ tensor, $u_2,u_3$ type $(p_2,q_2)$ and $(p_3,q_3)$ tensors, respectively.
\begin{enumerate}
\item $\norm{u_1}_{C^i(X)} \leq C_1(i,j,p_1,q_1)\norm{u}_{H^j(X)}$ for $j > \dim X/2+i$.
\item $\norm{u_2 \otimes u_3}_{H^k(X)} \leq C_2(k,k_1,k_2,p_2,p_3,q_2,q_3)\norm{u_2}_{H^{k_1}(X)}\norm{u_3}_{H^{k_2(X)}}$ for $k_1+k_2 > k+\dim X/2$, $k_1,k_2 \geq k$.
\end{enumerate}

We will also assume that $C_c^\infty(X)$ is dense in each $H^k(X)$. For brevity write $X_T = [0,T]\times X$. The we may also consider the space $C^k(X_T)$ and $H^k(X_T)$ for $k \in \N$, with the norm again defined via $\mathbf{e}$. Since $(X,e)$ satisfies the Sobolev inequalities, so does $(X_T,\mathbf{e})$, i.e.\
\begin{enumerate}
\item $\norm{u_1}_{C^i(X_T)} \leq C_1(i,j,p_1,q_1)\norm{u}_{H^j(X_T)}$ for $j > \dim X/2+i+1/2$.
\item $\norm{u_2 \otimes u_3}_{H^k(X_T)} \leq C_2(k,k_1,k_2,p_2,p_3,q_2,q_3)\norm{u_2}_{H^{k_1}(X_T)}\norm{u_3}_{H^{k_2(X+T)}}$ for $k_1+k_2 > k+\dim X/2+1/2$, $k_1,k_2 \geq k$,
\end{enumerate}
and the constants do not depend on $T$.
Let us set $\sob_e(k',r,s)$ to be the worst $C_1,C_2$ on the right-hand side of the Sobolev inequality, over all possible $i,j,k_i,k \leq k'$ ($i = 1,2$) $p_\ell \leq r, q_\ell \leq s$ ($\ell = 1,2,3$).

Now, let $g$ be a Lorentzian metric on $\R \times X$.
\begin{defn}We say $g$ is \emph{regular} (cf.\ the ``fundamental regularity hypothesis'' of \cite{ChoGene}*{appendix III}) with respect to $\mathbf{e}$ over $X_T$ if
\begin{enumerate}[label = (\roman*)]
\item $X_t := \{t\}\times X$ is spacelike for $g$ for $t \in [0,T]$, $\pa_t$ is timelike and future-oriented;
\item the Riemannian metrics $g_t = g|_{X_t}$ are uniformly equivalent to $e$, i.e.\ there exists a constant $B$, not depending on $t$, such that
\[B^{-1}e(V,V) \leq g_t(V,V) \leq Be(V,V)\] for all tangent vectors $V$ in $TX_t$;
\item the lapse $\Lambda$ (i.e.\ the projection of $ \partial_t$ onto the past-directed $g$-unit normal to $X_t$) is uniformly equivalent to $1$, i.e.\ there a constant $A$ such that
\[A^{-1} \leq \Lambda \leq A;\]
\item the shift $\beta$ (i.e.\ the projection of $\partial_t$ onto $M_t$), is uniformly bounded in the $e$ norm, i.e.\ there exists a constant $c$ such that
\[e(\beta,\beta) \leq c;\]
\item $g \in C^1(X_T)$.
\end{enumerate}
\end{defn}

Let us set $\reg_{\mathbf{e},T}(g) = \max(A^{-1},A,B^{-1},B,c,\norm{g}_{C^1(X_T)})$. Let us consider the hyperbolic equation for a type $(r,s)$ tensor $h$:
\begin{equation}\label{eq:C3:hyper}g^{-1}D^2 h + (a_1 + a_2)Dh + (b_1+b_2)h = f\end{equation} with Cauchy data $(h,\pa_t h)|_{t=0} = (\phi,\psi)$. The reason for the splitting $a_1+a_2$ and $b_1+b_2$ will be used when we state the estimate. In applications, $a_1,\ b_1$ will be the ``fixed'' coefficients, and $a_2,\ b_2$ will depend non-linearly on the solution. We also split $g^{-1} = g_1 + g_2$ into two symmetric cotensors. For $k \in \N$, and $T > 0$, let us assume that
\[Dg_1, a_1, b_1 \in C^k(X_T)\] and
\[Dg_2,a_2,b_2 \in H^k(X_T).\]

The following theorem is classical (see for instance appendix III of \cite{ChoGene}).
\begin{thm}[Energy Estimates]\label{thm:C3:energyestimate}Fix $T > 0$ and $k > \dim X/2 + 3$, and assume that $g$ is regular with respect to $\mathbf{e}$ on $X_T$. Then there exists natural numbers $n,r,s$ depending only on $k$, and a constant $C$, depending only on $T$, $k$, the $C^k(X_T)$ norms of $Dg_1, a_1,b_1$ and $\Riem(\mathbf{e})$, $(1+\sob_e(k,r,s))^n$ times the $H^k(X_T)$ norms of $Dg_2, a_2,b_2$, and $\reg_{\mathbf{e},T}(g)$, i.e.\
\begin{align*}
C = C(T,&k,\norm{D g_1}_{C^k(T)},\norm{ a_1}_{C^k(X_T)}, \norm{ b_1}_{C^k(X_T)},\\
&(1+\sob_e(k,r,s))^n\norm{D g_2}_{H^k(X_T)},(1+\sob_e(k,r,s))^n\norm{a_2}_{H^k(X_T)},\\
&(1+\sob_e(k,r,s))^n\norm{b_2}_{H^k(X_T)}, \reg_{\mathbf{e},T}(g),\norm{\Riem(e)}_{C^k(X_T)}),
\end{align*}
such that any solution $h \in L^2(X_T) \n C^1(X_T)$\footnote{one could take a number of different spaces here instead.} of \eqref{eq:C3:hyper} with Cauchy data $(h|_{t=0},\pa_t h|_{t=0}) = (\phi,\psi) \in H^{k+1}(X)\times H^k(X)$ is in fact in $H^{k+1}(X_T)$ and satisfies:
\[\norm{h}_{H^{k+1}(X_T)} \leq C(\norm{\phi}_{H^{k+1}(X)} + \norm{\psi}_{H^k(X)} + \norm{f}_{H^k(X_T)}).\]
Moreover $C$ may be taken to be increasing in all of its arguments.
\end{thm}
\begin{proof}[Proof outline.]This is very classical except for perhaps tracking why it is not the $H^k$ norms which appears in the dependence of the coefficient but some Sobolev constant times the $H^k$ norms. This follows from the usual proof in the following way. The usual energy estimates for $k=0$ depend on the $C^0$ norm of the coefficients, which in turn then depend on the $H^k$ norm of the coefficients via the Sobolev inequality. So in fact they depend on $\sob_e(k,r,s)$ times the $H^k$ norm. 

When deriving higher-order energy estimates, one differentiates the equation to get a similar equation for the derivatives of the solution. The coefficients remain essentially the same (with the addition of replacing $a_i$ by $a_i + Dg_i$, ($i = 1,2$) and $b_2$ with $b_2 + \Riem(\mathbf{e})$) but one obtains as sources products of derivatives of the coefficients, the Riemannian curvature tensor, and fewer derivatives of the solution. We want to put the sources in $L^2$. In fact, if we differentiate the equation $s \leq k$ times, the contributions to the source contain terms only of the form
\begin{align}
  \label{eq::energyestimates:eq1}
\begin{split}
&D^{s-j}g_1D^{j+2}h, \ D^{s-j}a_1D^{j+1}h, \ D^{s-j}b_1D^{j}h,\\
&D^{s-j}g_2D^{j+2} h, D^{s-j}a_2D^{j+1}h, \ D^{s-j}b_2D^{j}h,\ D^{s-j}\Riem(\mathbf{e})D^j h,\end{split}\end{align}
for $0 \leq j \leq s-1$, except the first term of each line, in which case $0 \leq j \leq s-2$ (see \cite{ChoGene}*{appendix II, section 3.6}). Arguing inductively, we may assume that $\norm{h}_{H^s(X_T)} < \infty$ with a bound of the form
\[\norm{h}_{H^s(X_T)} \leq C(\norm{\phi}_{H^{k+1}(X)} + \norm{\psi}_{H^k(X)} + \norm{f}_{H^k(X_T)}).\]and we wish to derive an estimate for $\norm{h}_{H^{s+1}(X_T)}$ of the same form.

We need to put each of the terms in \eqref{eq::energyestimates:eq1} in $L^2$ (with the appropriate bounds). For the terms in the second line, the first factor is in $C^0$ and the second is in $L^2$. For the first line, the first factor in in $H^{k+j-s}$ and the second factor is in $H^{s-j-2}$ for $0 \leq j \leq s-2$ and in $L^2$ for $j = s-1$ (the first term is not present in the case $j=s-1$). Sobolev multiplication therefore bounds all the terms in $L^2$, but the bound using Sobolev multiplication is of the form \[\norm{AB}_{L^2(X_T)} \leq (\sob_e(k)\norm{A}_{H^{k_1}(X_T)})\norm{B}_{H^{k_2}(X_T)},\] (for appropriate $k_1,k_2$), so at each iteration, it is not the $H^k$ norms of $Dg_1,a_1,b_1$ which enter, but $\sob_e(k)$ times the norms. \end{proof}

\subsection{Function spaces}
In this subsection we introduce the Banach spaces we will be using and prove some elementary properties. Recall the metrics $\mathbf{e}$ and $e$ from \cref{C:C3:Cauchyintro}. For $x > 0$, we denote by $\mathbf{e}_x$ and $e_x$ metrics on ${}^\sp TM$ and ${}^\sp TY$ obtained from $\mathbf{e}$ and $e$ by identifying on $M_x$ with $M$ and $\{x\}\times Y$ with $Y$, respectively. Let $D_x$ denote the Levi-Civita connection of $\mathbf{e}_x$ and let $\nabla_x$ denotes the Levi-Civita connection of $e_x$.

For the following two definitions, all our tensors will be ordinary tensors over $Y$.
\begin{defn}Fix $x>0$, and $k \geq 0$ an integer. Suppose $h$ is a tensor-distribution on $Y$ with $k$ measurable weak derivatives. We define
\[\norm{h}^2_{H^k_x(Y)} = \sum_{i=0}^k \int_{Y} |\nabla_x^i h|^2_{e_x}\ x^2d\vol_{e_x},\]
and let $H^k_x(Y)$ be the space of all such tensors with finite $H^k_x(Y)$ norm.\end{defn}

Observe that by definition $x^{2}d\vol_{e_x} = x^{2} dy\wedge x^{-2} d\vol_{\slash{k}} = d\vol_{e_1}$. The point of the $x^2$ factor is so that smooth and compactly supported short-pulse tensors have bounded $H^n_x$ norms down to $x=0$. We may similarly define:
\begin{defn}Fix $x> 0$ and fix $k \geq 0$. For $h$ a $k$-times continuously differentiable tensor distribution, set
\[\norm{h}_{C^k_x(Y)} = \sup_{i \leq k} \sup_Y |\nabla_x^i h|_{e_x},\]
and let $C^k_x(Y)$ be the space of all such tensors with finite $C^k_x(Y)$ norms.\end{defn}

We have the following elementary relationship between the spaces:
\begin{lem}Suppose $h$ is a compactly supported and $k$ times continuously differentiable. Then \[\norm{h}_{H^k_x}(Y) \leq C\norm{h}_{Ck_x}(Y),\] where $C$ depends only on the $e_1$-volume of $\supp h$.\end{lem}

For $T > 0$ and set $Y_T = [0,T]\times Y$. 

We define two more spaces:

\begin{defn}Fix $x > 0$, and fix $k \geq 0$. For $h$ a tensor distribution over $Y_T$ with $k$ measurable weak derivatives, set
\[\norm{h}_{H^k_x(Y_T)} =\sum_{i=0}^k \int_{Y} |D_x^i h|^2_{\mathbf{e}_x}\ x^{2}d\vol_{\mathbf{e}_x},\]
and let $H^k_x(Y_T)$ be the space of all such tensors with finite $H^k_x(Y_T)$ norms.\end{defn}

\begin{defn}Fix $x > 0$, and fix $k \geq 0$. For $h$ a $k$-times continuously differentiable tensor distribution over $Y_T$, set
\[\norm{h}_{C^k_x(Y_T)} = \sup_{i \leq k} \sup_{[0,T]\times Y} |D^i h|_{\mathbf{e}_x},\]
and let $C^k_x(Y_T)$ be the space of all such tensors with finite $C^k_x(Y_T)$ norms.\end{defn}

The following is a consequence of \cref{thm:C3:connectionisnice}.
\begin{lem}\label{thm:C3:tensoringisbounded}Suppose $h$ is a smooth type $(p,q)$ short-pulse tensor defined over $\mathcal M$, and $\psi$ is a smooth cutoff function of a compact set in $M$. Fix $T > 0$. Then $\sup_{x > 0} \norm{\psi h}_{C^k_x(Y_T)}$ is finite for all $k$, and tensoring with $\psi h$ is a bounded map from $H^k_x(Y_T)$ to itself.
\end{lem}

Fox fixed $x > 0$, $e_x$ is a metric of bounded geometry. Thus:
\begin{prop}For fixed $x >0$, $H^k_x(Y)$, $C^k_x$, $C^k_x(Y_T)$, $H^k_x(Y_T)$ are Banach spaces, and compactly supported smooth sections are dense in $H^k_x(Y)$, and $H^k_x(Y_T)$. 
\end{prop}

At $x=1$, recall that $e_1 = dy^2 + \slash{k}$, and let $\nabla_1$ be its Levi-Civita connection.

We have the following:
\begin{lem}\label{thm:C3:compare}Let $u$ be a type $(p,q)$ tensors defined on $Y$. Then pointwise in $Y$, for any $x > 0$:
\begin{subequations}
\label{eq:C3:sobprelim}
\begin{align}
|\nabla_x^ku|^2_{e_x} &\leq x^{-2p} |\nabla_1^k u|^2_{e_{1}}\\
|\nabla_1^ku|^2_{e_{1}} &\leq x^{-2(q+k)}|\nabla_x^k u|^2_{e_x}.\end{align}
\end{subequations}

Furthermore,
\begin{subequations}
\label{eq:C3:sobact}
\begin{align}
x^{q+k}\norm{u}_{H^k_1(Y)} &\leq \norm{u}_{H^k_x(Y)} \leq x^{-p}\norm{u}_{H^k_1(Y)}\\
x^{q+k}\norm{u}_{C^k_1(Y)} &\leq \norm{u}_{C^k_x(Y)} \leq x^{-p}\norm{u}_{C^k_1(Y)}\\
x^{q+k}\norm{u}_{C^k_1(Y_T)} &\leq \norm{u}_{C^k_x(Y_T)} \leq x^{-p}\norm{u}_{C^k_1(Y_T)}\\
x^{q+k}\norm{u}_{H^k_1(Y_T)} &\leq \norm{u}_{H^k_x(Y_T)} \leq x^{-p}\norm{u}_{H^k_1(Y_T)}
\end{align}
\end{subequations}
\end{lem}
\begin{proof}
Equations \eqref{eq:C3:sobact} follow immediately from \eqref{eq:C3:sobprelim}, so we show the latter. Choose an orthonormal frame $\{E_a\}$ of $S$, and denote by roman letters components in this basis. The index $1$ will denote plugging in $\pa_y$.
Then by definition
\begin{align*}
(e_x)_{ab} &= x^{-2}\slash{k}_{ab}\\
(e_x)^{ab} &= x^2\slash{k}^{ab}\\
(e_x)_{11} &= 1\\
(e_x)^{11} &= 1\\
(e_x)_{1a} &= 0\\
(e_x)^{1a} & = 0.\end{align*}
Suppose $h$ is a type $(p,q)$ tensor. Thinking of $u$ as a linear map taking $p$ covectors and $q$ vectors as arguments, label the $p$ covector arguments with $0,\ldots, p-1$ and and the $q$ vector arguments with $0,\ldots, q-1$. Let $I \subseteq \{0,\ldots, p-1\}$, $J \subseteq \{0, \ldots, q-1\}$, of size $i$ and $j$, respectively, and denote by $u^I_J$ the type $(p-i,q-j)$ tensor on $S$, the result of plugging in the vector field $dy$ in arguments labelled by $I$ and $\pa_y$ in arguments labelled by $J$. Then
\[|u|_{e_x}^2 = \sum_{I,J} |u^I_J|^2_{x^{-2}\slash{k}} = \sum_{I,J} x^{2(q-j) - 2(p-i)}|u^I_J|^2_{\slash{k}} \leq x^{-2p}\sum_{I,J}|u^I_J|_{\slash{k}} = |u|_{e_1}^2,\] where the sum runs over all $I,J$ of all lengths. This proves the first inequality in the case $k=0$. The second is proved in the same fashion.

Observe that $\nabla_x = \nabla_1$, since $e_x$ is given as a product metric and the connection induced by $x^{-2}\slash{k}$ does not depend on $x$.

Since $\nabla_x$ maps $(p,q)$ tensors to $(p,q+1)$ tensors, the $k=0$ case implies the full result.\end{proof}

As a corollary, we have the Sobolev inequalities:
\begin{prop}\label{thm:C3:Sobolev}Let $u$ be a type $(p,q)$ tensor defined on $Y$. Then
\begin{align*}
\norm{u}_{C^k_x(Y)} &\lesssim x^{-q-p-n}\norm{u}_{H^n_x(Y)}, \ \ n > k + \frac{1}{2}\dim Y\\
\norm{u_1\otimes u_2}_{H^n_x(Y)} &\lesssim x^{-p-2q-n_1-n_2}\norm{u_1}_{H^{n_1}_x(Y)}\norm{u_2}_{H^{n_2}_x(Y)},\\
& n_1+n_2 > n + \frac{1}{2}\dim Y, \ n_1,n_2 \geq n.\end{align*}\end{prop}
The analogous inequalities hold for $C^k_x(Y_T)$ and $H^k_x(Y_T)$ provided we use $\frac{1}{2}\dim Y + 1/2$ instead of $\frac{1}{2}\dim Y$.
\begin{proof}The Sobolev inequalities are certainly true for $e_1$. Now use \cref{thm:C3:compare}.\end{proof}

\begin{rk}The only important thing to remember about the exponents is that they are negative and depend only on the type of the tensor to which they are being applied and the amount of regularity measured by the spaces.\end{rk}

\subsection{Existence}
With our function spaces defined, we may now state a version of \cref{thm:C3:bigtheoremI} which we will prove first in this section. Suppose $(\phi,a,b)$ is admissible and $f$ is $0$ in Taylor series and compactly supported, and consider the equation
\begin{equation}\label{eq:C3:typeII}\phi(x^{\ell}h_x)D^2 h_x + a(x^{\ell}h_x,x^{\ell}Dh_x)Dh_x + b(x^{\ell}h_x,x^{\ell}Dh_x)h_x = x^{-\ell}f.\end{equation} Then:
\begin{thm}\label{thm:C3:bigtheoremII}For $T > 0$ and $\ell$ large enough, there exists $\delta > 0$ such that if $0 < x < \delta$ there exists a unique smooth section $h_x$ of $T^p_q(Y_T)$ solving \eqref{eq:C3:typeII} with $0$ Cauchy data. Furthermore, if $h_x$ is considered as a tensor over $(0,\delta)\times Y_T$, then $h_x$ is smooth (also as a function of $x$) and for any $i,j \geq 0$ there exists $c(i,j) > 0$ (not depending on $\ell$ or $h_x$) so that \[\sup_{0 < x < \delta} x^{c(i,j)}\norm{(x\pa_x)^ih_x}_{H^{j}_x(Y_T)} < \infty.\]\end{thm}

Of course, Sobolev embedding means the last statement implies that for all $i,j \geq 0$, there is some other constant $c(i,j) > 0$ so that \[\sup_{0 < x < \delta} x^{c(i,j)}\norm{(x\pa_x)^i h_x}_{C^{j}_x(Y_T)} < \infty.\]

\Cref{thm:C3:bigtheoremI} follows easily from \cref{thm:C3:bigtheoremII}.
\begin{proof}[Proof of \cref{thm:C3:bigtheoremI}]
If $h_x^{\ell}$ solves \eqref{eq:C3:typeII}, then $x^{\ell}h_x^{\ell}$ solves \eqref{eq:C3:type}. Fix $T > 0$. \Cref{thm:C3:bigtheoremII} provides us with $\delta > 0$ and at least one $\ell$ for which $x^{\ell}h_x^{\ell}$ solves \eqref{eq:C3:type} for $0 < x < \delta$ and $0 \leq t \leq T$. Set $h = x^{\ell}h^{\ell}$. Then $h$ is a smooth short-pulse tensor over $\{x > 0\}$ (since over $x > 0$ there is no difference between a short-pulse and ordinary tensor). In order to complete the proof we just need to show that $h$ is $0$ in Taylor series at $x=0$. A sufficient condition for this is for all $i,j$
\[\limsup_{x \to 0} \norm{\pa_x^i h_x}_{C^j_x(Y_T)} = 0.\]

By the usual uniqueness theorems, if $h_x^{\ell'}$ also solves \eqref{eq:C3:typeII}, then $h = x^{\ell'}h_x^{\ell'}$, too. In particular, \cref{thm:C3:bigtheoremII} says that, for all $i,j$ there is $c = c(i,j) > 0$ such that for all $\ell \geq i$ sufficiently large
\[\sup_{0 < x < \delta} x^{c}\norm{(x\pa_x)^i x^{-\ell}h_x}_{C^{j}_x(Y_T)} < \infty.\]

Let us compute $\pa_x^i h_x$ in terms of $(x\pa_x)^m x^{-\ell}h_x$, for $0 \leq m \leq i$:
\[\pa_x^i h_x = \sum_{k=0}^i \binom{i}{k}\frac{\ell!}{(\ell-k)!}x^{\ell-k}\pa_x^{i-k} (x^{-\ell}h_x) = \sum_{k=0}^i \binom{i}{k}\frac{\ell!}{(\ell-k)!}x^{\ell-i}x^{i-k}\pa_x^{i-k} (x^{-\ell}h_x).\]

The operator $x^{j-k}\pa_x^{i-k}$ is just a linear combination of $(x\pa_x)^m$ for $m \leq i-k$. Thus,
\[\norm{\pa_x^i h_x}_{C^j_x(Y_T)} \lesssim x^{\ell-i-\sup_{m \leq i} c(m,j)}\sum_{m=0}^i x^{c(m,i)}\norm{(x\pa_x)^m x^{-\ell}h_x}_{C^j_x(Y_T)}.\]
for $c(m,j)$ as above. In particular, the sum on the right-hand side is uniformly bounded for $0 < x < \delta$. We are always free to choose $\ell$ large so that $\ell-i-\sup_m c(m,j) > 0$, which means the right-hand side converges to $0$ as $x \to 0$. This completes the proof of the theorem.\end{proof}
We focus now on proving \cref{thm:C3:bigtheoremII}.

A type of smooth bundle map will appear often, so we give it a name.
\begin{defn}Let $X$ be a smooth manifold, and $E$, $F$ vector bundles over $[0,1)\times X$. We call a smooth bundle map $\psi:E \to F$ \emph{tame} if 
\begin{enumerate}[label=(\roman*)]
\item $\psi(0) = 0$ (i.e.\ $\psi$ maps the zero section to the zero section);
\item $\psi$ vanishes outside a large compact neighbourhood of the base;
\item $\psi$ is constant on each fibre outside a compact neighbourhood of the zero section, i.e.\ there exists a smooth section $f$ of $F$ and a compact neighbourhood $K$ of the $0$ section in $E$ such that off of $K$, $\psi$ is identically $f$.
\end{enumerate}
\end{defn}

If $(\phi,a,b)$ is admissible, then we can almost write $\phi = \phi_1 + \phi_2$, $a = a_1+a_2$ and $b = b_1+b_2$, where $\phi_1,a_1,b_1$ are tame and $\phi_2,a_2,b_2$ are just fixed sections (i.e.\ constant in the fibres). Indeed, just set $\phi_2 = \phi(0)$, $a_2 = a(0)$ and $b_2 = b(0)$ ($0$ denoting the zero section). However, we did not assume anything about the behaviour of $\phi,a,b$ at fibre infinity (for instance they are not necessarily even defined!) so $\phi_1,a_1,b_1$ need not be tame. 

However, we have:
\begin{lem}Given $(\phi,a,b)$ admissible, there is another admissible triple $(\tilde{\phi},\tilde{a},\tilde{b})$ which is equal to $(\phi,a,b)$ on neighbourhood of the $0$ section, but is constant on the fibres outside a large compact set. In particular, we may decompose $\tilde{\phi} = \phi_1 + \phi_2$, $\tilde{a} = a_1+a_2$ and $\tilde{b} = b_1+b_2$ as above into a tame bundle map and a fixed smooth section.\end{lem}
\begin{proof}Let $a_0 = a(0)$, $b_0 = b(0)$ where $0$ denotes the zero section. Let $\chi$ be a smooth bundle map from ${}^\sp T^p_q \mathcal M\oplus {}^\sp T^p_{q+1} \mathcal M$ to itself, which is the identity near $a_0$, and $0$ away from it. Set $\tilde{a} = \chi(a)$. Since $a$, is assumed to be constant outside a large compact set, $\tilde{a}$ satisfies the correct properties. The same argument works for $b$.

We cannot quite do this for $\phi$ because it is required to be Lorentzian among other things. Let $\phi_0 = \phi(0)$. Now let $\rho$ be a smooth bundle map from $^\sp T^p_q \mathcal M$ to $^\sp T^0_0 \mathcal M = M\times \R$ which is $1$ in some sufficiently small neighbourhood of the zero section and $0$ outside a slightly larger neighbourhood. Set $\tilde{\phi} = \rho \phi + (1-\rho)(\phi_0)$. Certainly $\tilde{\phi}$ is equal to $\phi$ near the zero section and is constant on the fibres far away from it. We need to check it is Lorentzian, considers $\{t\}\times Y_x$ spacelike for all $y$, and has $\pa_t$ timelike and future oriented. By assumption $\phi = -\pa_t^2+e^{-1}_x$ outside some large compact subset $K$ of the base, so we just need to check what happens over $K$. By definition
\[\tilde{\phi}-\phi_0 = \rho(\phi-\phi_0).\] Since $\rho$ is bounded, this tensor may be made arbitrarily small (in the sense that it takes values arbitrarily close to the zero section), provided that $\phi$ and $\phi_0$ are close on $\supp \rho$. This is sufficient to ensure all the conditions.

Now we just decompose by setting $a_2 = \tilde{a}(0)$, $b_2 = \tilde{b}(0)$, $\phi_2 = \tilde{\phi}(0)$.
\end{proof}

The advantage is that we may replace $(\phi,a,b)$ with $(\tilde{\phi}, \tilde{a},\tilde{b})$ without affecting solutions to \eqref{eq:C3:type} or \eqref{eq:C3:typeII}. In fact:
\begin{prop}\label{thm:C3:mightaswell}For $\delta$ small, a solution to \eqref{eq:C3:type} for $x < \delta$ with $0$ Cauchy data and with coefficients $(\tilde{\phi}, \tilde{a},\tilde{b})$ is the same as a solution to \eqref{eq:C3:type} with coefficients $(\phi,a,b)$. For $\delta$ small and for $\ell$ large enough, a solution to \eqref{eq:C3:typeII} for $0 < x < \delta$ with $0$ Cauchy data and with coefficients $(\tilde{\phi}, \tilde{a},\tilde{b})$ is the same as a solution to \eqref{eq:C3:type} with coefficients $(\phi,a,b)$. \end{prop}
\begin{proof}Assume we have a solution to \eqref{eq:C3:type} with coefficients $(\tilde{\phi},\tilde{a},\tilde{b})$. By \cref{thm:C3:bigtheoremI}, the solution $h$ vanishes in Taylor series. Thus $\phi(h) = \tilde{\phi}(h)$ for $x$ sufficiently small, and similarly for $a,b$.

Now assume we have a solution to \eqref{eq:C3:typeII} with coefficients $(\tilde{\phi},\tilde{a},\tilde{b})$. By \cref{thm:C3:bigtheoremII}, the solutions $h_x$ satisfy $\norm{h_x}_{C^1(Y_T)} \lesssim x^{-c}$ for some $c > 0$ not depending on $h$ or $x$. Thus provided $\ell > c$, $\phi(x^{\ell}h) = \tilde{\phi}(x^{\ell}h)$ for $x$ sufficiently small, and similarly for $a,b$.\end{proof}

In light of the proposition, we will make an additional assumption on an admissible triple $(\phi,a,b)$: that is constant at fibre-infinity, and call the new assumptions admissible*.

The composition of a tame map with a section in $H^k_x(Y_T)$ satisfies a good estimate. More precisely,
\begin{lem}\label{thm:C3:composition}Let $(r,s)$ and $(p_i,q_i)$ $i = 1,\ldots, n$ be pairs of integer, and suppose
\[\psi:\bigoplus_{i=1}^n{}^\sp T^{p_i}_{q_i}\mathcal M \to {}^\sp T^r_s \mathcal M\] is a tame bundle map. Fix $T > 0$ and $k > \dim Y/2+1$ and suppose $h_i \in H_x^k(Y_T)$ are sections of ${}^\sp T^{p_i}_{q_i}\mathcal M$, for $k > \dim Y/2+1$. Then there exists $c>0$, $N \geq 0$, depending only on $(p_i,q_i)$, $n$, $\psi$ and $k$ such that, setting
\[B = \left(\max_{1 \leq i \leq n} \norm{h_i}_{H^k_x(Y_T)}\right),\]
\[\norm{\psi(h_1,\ldots, h_n)}_{H^k_x(Y_T)} \lesssim x^{-c}B(1+B^N).\]
\end{lem}

The last subsection of this section will be devoted to proving this.
\begin{rk}The point of tameness is to obtain the factor $B$ in front of the estimate, rather than just $1$.\end{rk}

As a corollary we have:

\begin{cor}\label{thm:C3:bounds}Suppose $(\phi,a,b)$ is admissible$^{\ast}$, and suppose $k > \dim Y/2+1$. Then there exist $c> 0, N \geq 0$, depending only on $k$ and the type of tensors involved such that
\begin{enumerate}[label=(\roman*)]
\item $\phi(h)^{-1}$ is regular and $\reg_{\mathbf{e},T}(\phi(h)^{-1}) \lesssim 1+x^{-c}\norm{h}_{H^{k+1}(Y_T)}(1+\norm{h}_{H^{k+1}_x(Y_T)}^N)$;
\item $\phi = \phi_1 + \phi_2$, where $\phi_1$ is a fixed section satisfying $\norm{D \phi_1}_{C_x^k(Y_T)} < \infty$,\\
and $\norm{D \phi_2(h)}_{H_x^{k}(Y_T)} \lesssim x^{-c}\norm{h}_{H_x^{k+1}(Y_T)}(1+\norm{h}_{H^{k+1}_x(Y_T)}^N)$;
\item $a = a_1 + a_2$, where $a_1$ is a fixed section satisfying $\norm{a_1}_{C_x^k(Y_T)} < \infty$,\\
 and $\norm{a_2(h,Dh)}_{H_x^{k}(Y_T)} \lesssim x^{-c}\norm{h}_{H_x^{k+1}(Y_T)}(1+\norm{h}_{H^{k+1}_x(Y_T)}^N)$;
\item $b = b_1 + b_2$, where $b_1$ is a fixed section satisfying $\norm{b_1}_{C_x^k(Y_T)} < \infty$,\\
 and $\norm{b_2(h,Dh)}_{H_x^{k}(Y_T)} \lesssim x^{-c}\norm{h}_{H_x^{k+1}(Y_T)}(1+\norm{h}_{H^{k+1}_x(Y_T)}^N)$.
\end{enumerate}
The implied bounds depend only on $\phi, \ a, \ b$, and not on $h,\ x$.
\end{cor}
\begin{proof}In light of the decomposition of $\phi$, $a$, $b$, into a fixed section and a tame map, the last two statements follow from \cref{thm:C3:composition}. The second does, too, since $\norm{D\phi_2(h)}_{H^{k}_x(Y_T)} \leq \norm{\phi_2(h)}_{H^{k+1}_x(Y_T)}$.

Let us show (i). Condition (v) of regularity is equivalent to $\phi(h) \in C^1_x(Y_T)$ with bounds. Certainly $\phi_1 \in C^1_x(Y_T)$. We know that $\phi_2(h) \in H^{k+1}_x(Y_T)$ since $\phi_2$ is tame, and so Sobolev embedding gives us (v).

Conditions (i)-(iv) of regularity are not properties of the regularity of the function $\phi(h)$, but rather only properties of the range of $\phi$. Let $r$ denote a variable in the base $M$, and $\xi$ a variable in the fibre of $TM$, and for $x \geq 0$ set $g = \phi(x,r,\xi)^{-1}$. By assumption (i) is true for $g$. Now write
\begin{align*}
g &= -adt\otimes dt + bdy\otimes dy + 2c(dy\otimes dt + dt\otimes dy) + x^{-2}\slash{g}\\
&+ (dt\otimes \omega/x + \omega/x \otimes dt) + (dy\otimes \alpha/x + \alpha/x \otimes dy)\\
g^{-1} &= -a'\pa_t\otimes \pa_t + b'\pa_y\otimes \pa_y + c'(\pa_y\otimes\pa_t + \pa_t \otimes \pa_y) + x^2\slash{g}'\\
&+ (\pa_t \otimes x\omega' + x\omega'\otimes \pa_t) + (\pa_y \otimes x\alpha' + x\alpha' \otimes \pa_y),\end{align*}
where $a,b,c,a',b',c'$ are scalars, $a,b,a',b' > 0$, $\slash{g} \hspace{-1.4pt}\in\hspace{-1.4pt} \Sym^2(T_{(x,r)}^\ast S)$, $\slash{g}'\hspace{-1.4pt} \in\hspace{-1.4pt} \Sym^2(T_{(x,r)} S)$, $\omega,\alpha \in T_{(x,r)}^\ast S$, $\omega',\alpha' \in T_{(x,r)} S$, and $bdy^2 + x^{-2}\slash{g}$ is Riemannian.

By admissibility, $\phi$ attains all of its values in a compact set of $(r,\xi)$. Thus we may suppose that $r$ and $\xi$ lie in a compact set. Now let us prove property (ii) of regularity. Write $V = \Upsilon\pa_y + x\Theta$, where $\Upsilon \in \R$, $\Theta \in TS$, and $\Upsilon^2 + |\Theta|^2_{\slash{k}} = 1$. Then \[(x,r,\xi,V) \mapsto \phi(x,r,\xi)^{-1}(V,V)\] is a nonzero continuous function of $x,r,\xi$ and $\Upsilon$, $\Theta$ lying in a compact set, and hence is bounded above by some uniform $B$ and below by $B^{-1}$. This shows (ii).

Now let us show (iii) and (iv). Let us find a function $\tilde{\Upsilon}$ and a vector field $\tilde{\Theta}$ so that $W = -\pa_t + \tilde{\Upsilon}\pa_y + x\tilde{\Theta}$ is perpendicular to $X_t$. In fact, $\tilde{\Upsilon},\tilde{\Theta}$ are given by
\begin{gather*}
\tilde{\Upsilon} = b'c + \alpha'\cdot \omega\\
\tilde{\Theta} = c\alpha' + \slash{g}'\cdot \omega.\end{gather*}
Now notice that because of compactness, $\tilde{\Upsilon}$ and $\tilde{\Theta}$ are bounded above in the $\slash{k}$ norm. Also, $-|W|^2_g$ is a continuous, positive, function of $x,r,\xi$, and thus is bounded above and below because of compactness. Write $n = -\frac{1}{|W|^2_g}W$. Then $n = \frac{1}{|W|^2_g}\pa_t + \Upsilon\pa_y + x\Theta$ is the past-directed unit normal to $X_t$, where $\Upsilon = (-|W|^2_g)\tilde{\Upsilon}$ and $\Theta = (-|W|^2_g)\tilde{\Theta}$. Properties (iii) and (iv) follow.\end{proof}

We are now in a position to prove \cref{thm:C3:bigtheoremII}. We will prove the theorem in two steps. In step one, we will establish most of the theorem, and prove:
\begin{prop}\label{thm:C3:Ihatebigtheorem}Suppose $(\phi,a,b)$ is admissible and $f$ is $0$ in Taylor series and compactly supported. Fix $k > \frac{1}{2}\dim Y + 10$.\footnote{We choose $10$ to be safe. One could get away with a slightly smaller number.} Then for $T > 0$, and $\ell$ large enough, there exists $\delta > 0$ such that if $0 < x < \delta$ there exists a unique smooth section $h_x$ of $T^p_q(Y_T)$ in $H^k_x(Y_T)$ solving the equation
\eqref{eq:C3:typeII}
with $0$ Cauchy data. Furthermore \[\sup_{0 < x < \delta} \norm{h_x}_{H^{k}_x(Y_T)} < \infty.\]\end{prop}
In step two, we establish the higher regularity, i.e.\ we show:
\begin{prop}\label{thm:C3:IhatebigtheoremII}Let $k$, $T$, $\delta$, $h_x = h_x(\ell)$ be as in the previous proposition. Then for $\ell$ large enough, and for any $i,j \geq 0$ there exists $c(i,j) > 0$ (not depending on $\ell$ or $h_x$, but depending on the fixed $k$) so that \[\sup_{0 < x < \delta} x^{c(i,j)}\norm{(x\pa_x)^ih_x}_{H^{j}_x(Y_T)} < \infty.\] 
\end{prop}
\begin{rk}The $\ell$ in \cref{thm:C3:IhatebigtheoremII}, although larger than the $\ell$ in \cref{thm:C3:Ihatebigtheorem}, still does not depend on $i,j$.\end{rk}

\begin{proof}[Proof of \cref{thm:C3:Ihatebigtheorem}]
Notice that for all $j \geq 0$,
The spaces $H^{j}_x(Y_T)$ are just the spaces $H^{j}(Y_T)$ of \cref{C:C3:energyestimates} associated to the metric $\mathbf{e}_x$, except the norms differ by a factor of $x$, i.e.\ $\norm{\bullet}_{H^{j}_x(Y_T)} = x \norm{\bullet}_{H^{j}(Y_T)}$.

By \cref{thm:C3:Sobolev}, the associated $\sob_{\mathbf{e}_x}(m,p,q)$ is controlled by $x^{-d(m,p,q)}$, for some constant $d > 0$. Thus for any $m,p,q$, $(1+\sob_{\mathbf{e}_x}(m,p,q))^n \leq 1+x^{-c(m,p,q)}$ for some $c$.

Observe by \cref{thm:C3:curvatureisnice} that $\Riem(\mathbf{e})$ is a smooth short-pulse tensor.

Thus, by the usual theorem quasilinear existence theorems, for any $x$ and $\ell$, there is a smooth solution which exists for small time. Furthermore, if a solution exists on $Y_{S}$ for some $0 < S \leq T$, then the solution is unique and smooth on $Y_S$, and for $j \geq k$ the estimate
\[\norm{h_x}_{H^j_x(Y_S)} \lesssim_{j,x} \norm{h_x}_{H^k_x(Y_S)} + \norm{x^{-\ell}f}_{H^j_x(Y_S)}\] is valid, where the implied constant depends on $j$, $x$, but not on $h_x$ or $S \leq T$.\footnote{In the course of proving \cref{thm:C3:bigtheoremII}, we will prove a more refined version of this estimate (provided $\ell$ is sufficiently large) which shows that the constant blows up like a negative power of $x$. We do not need this precision at this point.} Thus, if we can show that for $x$ small and $\ell$ large, there exists a constant $M$, such that any solution $h_x$ on $H^k_x(Y_S)$, for any $S \leq T$, satisfies $\norm{h_x}_{H^k_x(Y_S)} \leq M$, then we can iterate local well-posedness and obtain a solution on all of $Y_T$, the key point being $M$ does not depend on $S$.\footnote{The argument we use here loses some regularity; indeed, in order to iterate well-posedness, we need to use Sobolev restriction to control $h_x|_{\{t = S\}}$ and $\pa_t h_x|_{\{t = S\}}$ in $H^{k}_x(Y)$ and $H^{k-1}_x(Y)$ respectively, which necessitates controlling $h_x$ in $H^{k+1}_x(Y_S)$. Instead, one can immediately use the estimate on the $H^k_x(Y_s)$ norm to control $h_x$ in $C^0_t H^{k-1}(Y_s) \n C^1_t H^{k-2}(Y_s)$. One would like to replace $k-1$ and $k-2$ respectively with $k$ and $k-1$ to iterate well-posedness. To do this, differentiate the equation $k-1$ times to obtain a linear equation with coefficients in $C^1(Y_s)$ and a source in $L^2(Y_s)$. The linear theory now applies yielding a solution which is in fact in $L^{\infty}_t L^2(Y_s)$, which one eventually can upgrade to $C^0_t H^1(Y_s) \n C^1_t L^2(Y_s)$.

Since we work with very regular data and coefficients, we do not need the more precise approach.}

Using \cref{thm:C3:energyestimate} and \cref{thm:C3:bounds}, we have that there is some $c > 0$ such that
\begin{equation}\label{eq:bigtheoremII:eq3}\norm{h_x}_{H^{k}_x(Y_S)} \leq C(S,x^{-c}\norm{x^\ell h_x}_{H^{k}_x(Y_S)}(1+\norm{x^\ell h_x}_{H^{k}_x(Y_S)}^N))\norm{x^{-\ell} f}_{H^{k-1}_x(Y_S)}.\end{equation}
Here, we suppressed the explicit dependence of the constant $C$ on $a,b,\phi,k,\mathbf{e}$, and only keep track of the contributions to $C$ which depend on $S$ and $h$.
We may assume $C$ increasing in both arguments.

Notice that neither $c$ nor the function $C$ depends explicitly on $\ell$. Thus we are free to take $\ell > c + 1$. For shorthand write $A(S) = \norm{h_x}_{H^{k}_x(Y_S)}$. Since $f$ vanishes in Taylor series and is compactly supported, $\sup_x \norm{x^{-\ell} f}_{H^{k-1}_x(Y_T)} \leq D < \infty$. Factoring out an $x^{\ell}$ from $\norm{x^{\ell}h^{\ell}_x}_{H^k_x}$ to cancel with $x^{-c}$, as long as $S \leq T$, we may continue from \eqref{eq:bigtheoremII:eq3}
\[A(S) \leq C(T,xA(S)(1+A(S)^N))D := F(xA(S)(1+A(S)^N)),\]
where we have set $F(y) = C(T,y)D$.

$F$ is increasing, so 
\[F(0^+) = \lim_{\hspace{5pt} y \to 0^+} \hspace{-3pt}F(z)\] exists. Set $M = 1+F(0^+)$ (which only depends on $k$, $\ell$). Thus, we may choose $\delta$ small enough so that $F(xz(1+z^N)) \leq M$ whenever $x < \delta$ and $z \leq 2M$.
Thus, if $x< \delta$ and $A(S) \leq 2M$ for $S \leq T$, then
\[A(S) \leq F(x2M(1+(2M)^N)) \leq M.\]

Recall that the function $A(s)$ is continuous on $(0,S]$,\footnote{it is not defined as $s = 0$.} and $A(s) \to 0$ as $s \to 0$. Thus, we may use a bootstrap argument to conclude $A(s) \leq M$ for all $s \in (0,S]$, and in particular $A(S) \leq M$. This completes the proof.
\end{proof}

Before proceeding, we will need one definition, primarily used for notational convenience.
\begin{defn}A \emph{constant-coefficient polynomial} is a smooth map
\[P:\bigoplus {}^\sp T^{p_i}_{q_i}M \to {}^\sp T^r_s M\] which tensors its arguments with itself several times, and then contracts them in some fashion. In other words, $P$ is an iterated composition of, the identity map, the maps $\xi \mapsto \xi\otimes \xi$, and a contraction map $\xi\otimes \eta \mapsto \xi\cdot \eta$.\end{defn}

\begin{proof}[Proof of \cref{thm:C3:IhatebigtheoremII}]
Choosing $\ell$ large enough, the case $i=0$, $j\leq k$ was proved in the previous proposition. Note that we will need to pick $\ell$ even larger once in the course of the proof.

Let us start with the special case $i = 0$, and show the proposition for all $j$ (i.e.\ there are no $x\pa_x$ derivatives). We have already mentioned that the subcase $j \leq k$ comes from the previous proposition. So let us start with $j=k+1$. We will differentiate \eqref{eq:C3:typeII} with $D$ to establish an equation for $u = Dh$. Due to the quasilinear nature of the problem and the fact that $A$ and $B$ depend on $Dh$ in addition to $h$, when taking the first derivative it will not be sufficient to treat as source terms all terms which result when $D$ hits a coefficient in \eqref{eq:C3:typeII}. Let us start by examining what happens when $D$ hits a coefficient, using the coefficient $\phi$ as an example. Let us decompose $\phi = \phi_1 + \phi_2$, where $\phi_1$ is tame, and $\phi_2$ is a fixed smooth section. Notice that $D \phi_2$ is a smooth section. Using a version of the chain rule (which can be proved by choosing a basis of short-pulse vector fields and using the Christoffel symbols of $\mathbf{e}$ to express the connection $D$),
\[D\phi_1(x^\ell h_x) = \psi_1(x^{\ell} h_x) + \psi_2(x^{\ell} h_x)(x^{\ell}Dh_x),\]
where $\psi_1$ is tame, and $\psi_2$ is the sum of a tame map and a fixed smooth section. Indeed, the term $\psi_2$ comes from the chain rule and the term $\psi_1$ comes from the connection not being an exact coordinate-wise derivative. Thus one term we obtain after differentiating \eqref{eq:C3:typeII} is
\[(\psi_1(x^\ell h_x) + \psi_2(x^{\ell} h)(x^{\ell} Dh_x)  + D\phi_2)Du,\]
which comes if $D$ hits the first factor of $\phi(x^{\ell} h_x)D^2 h_x$.

A similar result is true when differentiating the terms with $a$ or $b$. 

When applying $D$ to \eqref{eq:C3:typeII}, there are two sorts of terms. Those that come from hitting a coefficient and those that come from hitting a a derivative of $h_x$. Those that come from hitting a derivative of $h_x$ may be replaced by the corresponding derivative of $u$. For those that come from hitting a coefficient, there are two types of terms: those which only contain derivatives of order at most $1$, and those which contain a derivative of order $2$. Those that contain a derivative of order at most $1$ may be moved to the right-hand side and treated as a source. Those that containing a derivative of order $2$ must be replaced with an order-one derivative of $u$ and treating as part of the coefficient of $Du$.

Thus, $u$ satisfies an equation of the form
\begin{equation}\label{eq:C3:Ihatethis}\phi(x^{\ell}h)D^2 u + a'(x^\ell h_x,x^\ell Dh_x)Du + b(x^\ell h_x, x^\ell Dh_x)u = F(h_x,Dh_x),\end{equation}
where $F(h_x,Dh_x)$ takes the form
\[F = D(x^{-\ell}f) + \sum_i F_i(x^\ell h_x, x^{\ell}Dh_x)P_i(h_x,Dh_x)\] where $F_i$ is the sum of a tame map and a fixed smooth section, and $P_i$ is a constant-coefficient polynomial (which is not just the trivial polynomial that is just a constant), and where
\[a'(\xi,\eta) = a'_0 + a'_1(\xi,\eta) + a'_2(\xi,\eta) \xi + a'_3(\xi,\eta) \eta,\]
and $a'_i$, for $i =1,2,3$ is the sum of a  tame map and a fixed smooth section, and $a'_0$ is a fixed smooth section.

Next, we show that $u$ has trivial Cauchy data. Writing $u = Dh_x = (\pa_t h_x, \nabla_x h_x)$, it is then clear that $u = 0$, and $\pa_t u = (\pa_t^2 h_x, \pa_t \nabla_x h_x)$. Observe $\pa_t\nabla_x h_x = \nabla_x \pa_t h_x = 0$ (exchanging the derivatives is valid since $[\pa_t,\nabla] = 0$ on $C^2$ sections and $h_x \in C^2$ follows from $h_x \in H^k_x(Y_T)$). To obtain $\pa_t^2 h_x$, write
\[\phi(0)|_{\{t = 0\}} = -a\pa_t\otimes \pa_t + (\pa_t \otimes \omega + \omega\otimes \pa_t) + \slash{g},\]
where $a >0 $ is a smooth scalar, $\omega \in C^{\infty}(\{t = 0\};{}^\sp TY)$, and $\slash{g} \in C^{\infty}(\{t = 0\};\Sym^2({}^\sp TY))$.
Then
\[\phi(0)D^2 h = -a\pa_t^2 h + 2\pa_t\nabla h + \slash{k}\nabla^2 h.\] Using \eqref{eq:C3:typeII} and that $x^{\ell}h|_{\{t= 0\}},\pa_t x^{\ell}h_{\{t=0\}} = (0,0)$ by assumption, we obtain that $a\pa_t^2 h = 0$, and thus since $a > 0$, $\pa_t^2 h = 0$.
Although the triple $(\phi,a',b)$ is not admissible*, because of the terms $a'_2, a'_3$, \cref{thm:C3:bounds} still applies to the triple since we may use Sobolev multiplication to control the products by the polynomial factors in $x^\ell h_x$ and $x^{\ell}Dh_x$. Thus, like in the proof of \cref{thm:C3:Ihatebigtheorem}, we may use \cref{thm:C3:energyestimate} to obtain, for some $c > 0$, the estimate
\begin{equation}\label{eq:C3:Iwanttobedone}\norm{u}_{H^k_x(Y_T)} \lesssim C(x^{\ell-c}\norm{h}_{H^k_x(Y_T)})\norm{F}_{H^{k-1}_x(Y_T)},\end{equation} were we have omitted the exact dependence of $C$ on quantities which are certainly bounded. Choosing $\ell$ large (this is the only point in the proof where we will need to increase $\ell$), we may assume $\ell-c \geq 0$. Thus
\[\norm{u}_{H^k_x(Y_T)} \lesssim\norm{F}_{H^{k-1}_x(Y_T)}.\]

Now again using Sobolev multiplication, and \cref{thm:C3:composition}, we obtain the bound
\[\norm{u}_{H^k_x(Y_T)} \lesssim \norm{F}_{H^{k-1}_x(Y_T)} \lesssim x^{-c'}(1+\norm{h_x}_{H^{k}(Y_T)})^N \lesssim x^{-c'},\] which is the desired result for $i=0, j=k+1$.

Now let us treat the case $i=0$ and $j > k+1$. Write $j = j'+k$. We will use induction on $j'$. The coefficients depend only on $h_x$ and $Dh_x$, so all terms resulting from hitting $a'$ or $b$ with $D$ may safely be treated as source terms. The same is not true of the first term, because we cannot safely control $D^2 u_{j'}$ to obtain an equation for $u_{j'+1}$. Thus, we will need to modify the coefficient of $D^2u_{j'}$ as we did to obtain an equation for $u_{1} = u$. We may differentiate the sources inductively using the variant of the chain rule described above.

Differentiating \eqref{eq:C3:Ihatethis}, we obtain an equation of the form
\begin{equation}\label{eq:C3:higherders}\phi(x^{\ell}h)D^2 u_{j'} + a^{j'}(x^\ell h_x,x^\ell Dh_x)Du_{j'} + b(x^\ell h_x, x^\ell Dh_x)u_{j'} = F^{j'}(u_0,u_1,\ldots, u_{j'}),\end{equation}
where $F^{j'}(h_x,Dh_x)$ takes the form
\[F^{j'} = D^{j'}(x^{-\ell} f)\sum_i F^j_i(x^\ell u_0,\ldots, x^\ell u_{j'})P_i(u_0,\ldots,u_{j'}),\] where $F_i^j$ is the sum of a tame map and a fixed smooth section, and $P_i$ is a constant-coefficient polynomial (which is not just the trivial polynomial that is just a constant), and $a^{j'}$ is given inductively by
\[a^{j'}(\xi,\eta) = a^{j'-1}(\xi,\eta) + a^{j'}_0 + a_1^{j'}(\xi,\eta) + a_2^{j'}(\xi,\eta)\xi + a_3^{j'}(\xi,\eta)\eta\]
where $a^{j'}_i$, for $i =1,2,3$ is the sum of a  tame map and a fixed smooth section, and $a^{j'}_0$ is a fixed smooth section. Notice that even though we changed the coefficient $a^{j'}$, the type of tensor did not change.

If we can show that the data for $u_{j'}$ are smooth short-pulse tensors and compactly supported, then we can conclude by induction on $j'$ using the same argument as the case $j'=1$.\footnote{The appearance of $u_{j'}$ as an argument of $F$ is not harmful, since $\norm{u_{j'}}_{H^{k-1}_x(Y_T)} \lesssim \norm{h_x}_{H^{k+j'-1}_x(Y_T)}$, which is controlled by induction.} Since the type of tensors appearing as coefficients does change, the loss $x^{-c}$ when applying \cref{thm:C3:bounds} is the same for each $j'$, and hence $\ell$ does not need to be taken larger. 

We will show that the data has the desired properties by induction on $j'$. The case $j'=1$ has already been established, and the data are trivial. Now suppose it is true for $j'-1$. We show it for $j'$, too. Write $u_{j'} = Du_{j'-1} = (\pa_t u_{j'-1},\nabla_x u_{j'-1})$. Thus, $u_{j'}$ is smooth of compact support, and $\pa_t u_{j'} = (\pa_t^2 u_{j'-1}, \pa_t \nabla_x u_{j'-1})$. Like the case $j'=1$, the second component is smooth and compact support, being equal to $\nabla_x \pa_t u_{j'-1}$. For the first component, using the same argument as in the case $j' = 1$, except applied to \eqref{eq:C3:higherders} for the case $j'-1$, we obtain that $a\pa_t^2 u_{j'-1}$ is smooth and of compact support. Since $a > 0$, $1/a$ is smooth on the compact support of $u_{j'-1}$, and thus $\pa_t^2 u_{j'-1}$ is smooth and of compact support.

To prove the full proposition, i.e.\ establish the bounds for any $i,j$, one may simply differentiate \eqref{eq:C3:typeII} using $x\pa_x$ multiple times, and use a similar argument to the above. The key difference is that we do not need to change the coefficients, since all terms may be treated as sources, at the expense of allowing the sources to depend on $D$-derivatives up to order $2$. For instance, the equation for $v = x\pa_x h_x$ is
\begin{equation}\label{eq:C3:Ihatehighders}\phi(x^\ell h_x)D^2v + a(x^\ell h_x, x^\ell Dh_x)Dv + b(x^\ell h_x, x^\ell Dh_x)v = F(h_x,Dh_x,D^2h_x),\end{equation}
where $F(h_x,Dh_x,D^2h_x)$ takes the form
\[F = \sum_i F_i(h_x,Dh_x,D^2h_x)P_i(h_x,Dh_x,D^2 h_x) ,\]
where $F_i$ is the sum of a tame map and a fixed section, and $P_i$ is a non-trivial constant-coefficient polynomial. The dependence on $D^2 h_x$ is not harmful, since it is already controlled from the case $i=0$. Using the case $i=0$ (and the determination of the initial data) one may control $\norm{v}_{H^k_x(Y_T)}$. Differentiating this equation with $D$ and using the same argument as above, one may control $\norm{v}_{H^m_x(Y_T)}$ for $m > k$. For higher $i$, one may just differentiate \eqref{eq:C3:Ihatehighders} and use the same argument.
\end{proof}

\subsection{Composition with bundle maps}
In this subsection we prove \cref{thm:C3:composition}. Let $(r,s)$ and $(p_i,q_i)$ $i = 1,\ldots, n$ be pairs of integers, and suppose $ h_i$ are (distributional) sections of ${}^\sp T^{p_i}_{q_i}M$ defined over $M_x$.
Also, let
\[\psi:\bigoplus_{i=1}^n{}^\sp T^{p_i}_{q_i}M \to {}^\sp T^r_s M\] be a tame bundle map.

First, because of \cref{thm:C3:compare}, we may prove a version of the lemma for the manifold $Y \times [0,T]$ and the space $H^k(Y_T)$, rather than $H^k_x(Y_T)$. Next, since $\psi$ is tame, it vanishes outside a compact set. Let $\{U_j\}$ be a locally finite collection of coordinate charts covering $\supp \psi$ on which $TY_T$ is trivializable. Let $\{\rho_j\}$ be a partition of unity subordinate to $\{U_j\}$. Then
\[\psi(h_1,\ldots,h_n) = \sum_j \rho_j\psi(\tilde{\rho}_jh_1,\ldots, \tilde{\rho}_jh_n).\]

Consider for $m \geq 0$ and $N \in \N$ the usual norms $\norm{u}_{H^m(U_j)}$ for $u$ supported in $u_j$ taking values in $\R^N$. Since there are only finitely many $U_j$, the $H^m(U_j)$ norms are \emph{uniformly} equivalent to the $H^m(Y)$ norms of sections supported inside $U_j$ (with constant depending on the $U_j$ and hence ultimately on $\psi$). Because there is no salient distinction between a finite collection of vector-valued functions and one vector-valued function (taking values in a higher-dimensional space), we may additionally assume that $n=1$. Notice that $\tilde{\rho}_j\psi$ is a smooth function of compact support, which is constant in the fibres outside a compact subset.

Thus we have reduced \cref{thm:C3:composition} to:

\begin{lem}\label{thm:C3:compositionprime}Let $U$ be an open subset of $\R^d$ ($d \geq 1$) and suppose for $N,M \geq 1$ that $\psi: U \times \R^N \to \R^M$ is a smooth function of compact support in $U$ such that $\psi(x,0) = 0$ and $\pa_y \psi(x,y) = 0$ for $y$ large enough. Let $h:U \to \R^N$ be in $H^k(U)$, for $\dim U/2 < k \in \N$. Then $\psi(x,h(x)) \in H^k(U)$, too, and the estimate
\[\label{eq:C3:latenight}\norm{\psi(x,h(x))}_{H^k(U)} \lesssim \norm{h}_{H^k(U)}(1+\norm{h}_{H^k(U)}^\alpha),\] holds for some $\alpha \geq 0$ (depending only on $k$).
\end{lem}
\begin{proof}
This proof is on page 273 of \cite{MosRapi}. We reproduce the proof here for the reader's convenience. First, observe that for any $h \in H^k(U)$, $h \in C^0(U)$ by Sobolev embedding, and thus $\phi(x,h(x))$ is well-defined. Let $\pa_x$ denote the gradient operator in the $x$ direction, and $\pa_y$ in the $y$ direction.
First assume $h$ is smooth. To prove \eqref{eq:C3:latenight}, it suffices to bound $\psi(x,h(x))$ in $L^2$ and $\pa_x^k \psi(x,h(x))$ in $L^2$. Since $\psi(0) = 0$, we have the estimate
\[|\psi(x,y)| \lesssim |y|.\] The $L^2$-bound follows immediately. To treat the $\pa_x^k$ derivatives, observe that
\[\pa_x^k \psi(x,h(x)) = \sum_{\rho + \sigma \leq k} (\pa_x^\sigma \pa_y^\rho \psi)(x,h(x)) \sum C_{\alpha\rho\sigma}\prod_{i=1}^k (\pa_x^i h)^{\alpha_i},\]
where 
\begin{gather*}
\alpha_1 + \ldots + \alpha_k = \rho\\
1\alpha_1 + 2\alpha_2 + \cdots + k\alpha_k + \sigma = k.
\end{gather*}

It suffices to bound each term. If $\rho = 0$, then these terms corresponds to $(\pa_x^k \psi)(x,h(x))$, and since $\pa_x^k \psi$ satisfies the same assumptions as $\psi$ does, the bound
\[\norm{(\pa_x^k \psi)(x,h(x))}_{L^2(U)} \lesssim \norm{h}_{L^2(U)}\] still holds.
Now assume $\rho \geq 1$.
Define $p_i = \frac{k}{i\alpha_i}$, $q = \frac{k}{\sigma}$. Then
\[q_i^{-1} + \sum p_i^{-1} = 1.\]
Write $v_{\rho,\sigma} = (\pa_x^\sigma \pa_y^\rho \psi)(x,h(x))$. We may apply H\"older's inequality to obtain that
\begin{align*}
\norm{\pa_x^k \psi(x,h(x))}^2_{L^2(U)} &= \int_U |\pa_x^k \psi(x,h(x))|^2\ dx\\
&\lesssim \sum\int_U (\pa_x^{\sigma}\pa_y^{\rho}\psi)(x,h(x))^2\prod (\pa_x^i h)^{2\alpha_i}\ dx\\
&\lesssim \norm{h}^2_{L^2(U)} + \sum_{\rho \geq 1} \norm{v_{\rho,\sigma}}^2_{L^{2q}(U)}\prod \norm{\pa_x^i h}_{L^{2k/i}(U)}^{2\alpha_i}.\end{align*}
Here and below, the sum runs over all $\alpha,\rho,\sigma$ with $\rho \geq 1$.
The first factor in each summand of the second term is bounded above uniformly by some constant independent of $h$. The Gagliardo-Niremberg inequality followed by Sobolev embedding implies that
\[\norm{\pa_x^i h}_{L^{2k/i}(U)} \lesssim \norm{\pa_x^k h}_{L^2}^{i/k}\norm{h}_{L^{\infty}}^{1-i/k} \lesssim \norm{h}_{H^k}.\]
Thus,
\begin{align*}\norm{\pa_x^k \psi(x,h(x))}_{L^2(U)}^2 &\lesssim \norm{h}_{L^2(U)}^2 + \sum_{\rho \geq 1} \norm{h}_{H^k(U)}^{\sum 2\alpha_i}\\ &=\norm{h}_{L^2(U)}^2 + \sum_{\rho \geq 1} \norm{h}_{H^k}^{2\rho} \lesssim \norm{h}^2_{H^k(U)} + \norm{h}^{2k}_{H^k(U)},\end{align*}
which is the $L^2$ bound on the $k$th derivatives required for \eqref{eq:C3:latenight}.

Now if $h$ is not smooth, we approximate $h$ with smooth $h_j \to h$ in $H^k$. By Sobolev embedding, $h_j \to h$ in $C^0$ and thus $\psi(x,h_j(x)) \to \psi(x,h(x))$ in $C^0$. We must additionally show that the sequence $\psi(x,h_j(x))$ is Cauchy in $H^k$, since then it will also converge to $\psi(x,h(x))$ in $H^k$, and the estimate \eqref{eq:C3:latenight} holds by taking limits. We need only estimate in $L^2$ the $0$th and $k$th derivatives. Estimating the $0$th derivative is not difficult because of the estimate $|\psi(x,y)-\psi(x,y')| \lesssim C|y-y'|$, and gives for $j_1,j_2 \in \N$,
\[\norm{\psi(x,h_{j_1}(x))-\psi(x,h_{j_2}(x))}_{L^2(U)} \lesssim \norm{h_{j_1}(x)-h_{j_2}(x)}_{L^2(U)}.\] For the $k$th derivative, we need only control in $L^2$ the differences
\begin{equation}\label{eq:C3:diff0}(\pa_x^\sigma \pa_y^\rho \psi)(x,h_{j_1}(x))\prod_{i=1}^k (\pa_x^i h_{j_1})^{\alpha_i} - (\pa_x^\sigma \pa_y^\rho \psi)(x,h_{j_2}(x)) \prod_{i=1}^k (\pa_x^i h_{j_2})^{\alpha_i},\end{equation} for $j_1,j_2$. If $\rho = 0$, then we use the same estimate as before since the big product is not present. Otherwise, $\rho \geq 1$, and we may use the triangle inequality and the usual trick for bounding products to bound the $L^2$ norm of the difference by the sum of the $L^2$ norms of
\begin{equation}\label{eq:C3:diff1}((\pa_x^\sigma \pa_y^\rho \psi)(x,h_{j_1}(x))-(\pa_x^\sigma \pa_y^\rho \psi)(x,h_{j_2}(x)))\prod_{i=1}^k (\pa_x^i h_{j_1})^{\alpha_i}\end{equation}
and terms of the form
\begin{equation}\label{eq:C3:diff2}((\pa_x^\sigma \pa_y^\rho \psi)(x,h_\ell(x))(\pa_x^{i_0}h_j - \pa_x^{i_0} h_{j_2})\prod_{\ell=1}^2\prod_{i=1}^k (\pa_x^i h_{j_\ell})^{\alpha_{i,\ell}},\end{equation}
where $1 \leq i_0 \leq k$, and $\alpha_{i_1} + \alpha_{i_2} = \alpha_i$ for $i \neq i_0$ and $\alpha_{i_1} + \alpha_{i_2} = \alpha_i -1$ for $i = i_0$.

For \eqref{eq:C3:diff1}, the first factor can be estimated using the same argument as above, because $|\pa_x^\sigma \pa_y^\rho \psi(x,y)-\pa_x^\sigma \pa_y^\rho \psi(x,y')| \lesssim |y-y'|$, and the second factor is bounded by $\norm{h_{j_1}}_{H^k}^{\rho}$ by the argument we used above.

For \eqref{eq:C3:diff2} we use a similar argument to how we bounded the smooth case. Set $p_{i_\ell} = \frac{k}{i\alpha_{i,\ell}}$, $r = \frac{k}{i_0}$, and $q = k/\sigma$ as before. Then for the same reason as before, $q^{-1} + r^{-1} + \sum_{i,\ell} p_{i,\ell} ^{-1} = 1$, and so applying H\"older's inequality, the Gagliardo-Niremberg inequality, and Sobolev embedding bounds \eqref{eq:C3:diff0} by
\[\norm{h_{j_1}-h_{j_2}}_{H^k(U)}\left(\max_{\ell = 1,2}\norm{h_{j_\ell}}_{H^k(U)}\right)^{\rho-1}.\] Combining the estimates for $\rho = 0$ and all $\rho \geq 1$, we see the sequence in Cauchy.\end{proof}

\section{Completing the proof of \texorpdfstring{\cref{thm:C3:SPCIVP}}{theorem 3.2.19}}
\label{C:C3:completetheproof}
\todo{make sure the texorpdftostring reference is correct in the title}
In this section we complete the proof of \cref{thm:C3:SPCIVP}. By the discussion in \cref{C:C3:SCIVP}, the only thing we have left to do is complete the proof of \cref{thm:C3:SPexistence}.

Recall from \cref{thm:C3:SPTSI} that we have found a short-pulse metric $g$ on $\mathcal M([0,a)\times[0,b))$ solving $\Ric(g) = 0$ in Taylor at $\mathcal H_1 \un \mathcal H_2$ and at $x=0$, and we are looking for an exact solution to $Q(h) = \widetilde{\Ric}(g+h) = 0$, where $\widetilde{\Ric}(g+h) = 0$ are the reduced Einstein equations associated to a solution $g+h$ in wave gauge with respect to $g$. We will use \cref{thm:C3:bigtheoremI}, so let us begin by putting $Q(h)$ into a more amenable form. Let us denote by $\nabla$ the Levi-Civita connection of $g$.

\begin{prop}\label{thm:C3:amenable}The equation $-2Q(h) = 0$ takes the form
\[(g+h)^{-1}\nabla^2 h + A(h,\nabla h)\nabla h + B(h,\nabla h) h = -2\Ric(g),\]
where $A,B$ are smooth bundle maps
\begin{align*}
A:U \subseteq {}^\sp T^0_2 \mathcal M \oplus {}^\sp T^0_3 \mathcal M \to {}^\sp T^3_2 \mathcal M\\
B: U \subseteq {}^\sp T^0_2 \mathcal M \oplus{}^\sp T^0_3 \mathcal M \to {}^\sp T^2_2 \mathcal M,
\end{align*}
and where $U$ is a neighbourhood of the $0$ section.
\end{prop}
\begin{proof}
By \eqref{eq:C2:hastheform},
\[-2Q(h) = (g+h)^{-1}\nabla^2 h + (g+h)^{-1}\cdot (g+h)^{-1}\cdot \nabla h \cdot \nabla h -2(g+h)^{-1}\cdot (g+h)\cdot \Riem(g),\]
where $\cdot$ indicates a certain contraction. The first term has the desired form. The second does, too, since $(\xi,\eta) \mapsto (g+\xi)^{-1}\cdot (g+\xi)^{-1}\cdot \eta$ is a smooth bundle map already. Now let us treat the third term. We begin by splitting
\[2(g+h)^{-1}\cdot (g+h)\cdot \Riem(g) = 2(g+h)^{-1}\cdot g \cdot \Riem(g)+ 2(g+h)^{-1}\cdot h \cdot \Riem(g).\]
The map $\xi \mapsto 2(g+\xi)^{-1}\cdot \Riem(g)$ is a smooth bundle map, so the second term here is fine. Now let us look at the first term. We may write $(g+h)^{-1} = g^{-1} + \psi(h)$, for a smooth bundle map $\psi$ which vanishes at the zero section. Thus the first term is
\[2g^{-1}\cdot g \cdot \Riem(g) + \psi(h)\cdot g\cdot \Riem(g).\] The exact nature of the contractions means that the first term here is just $2\Ric(g)$ (alternatively it is $2Q(0) = 2\Ric(g)$ by definition). The second term is of the desired form, since a fixed smooth section is certainly a bundle map. Since $\psi$ vanishes at the zero section there exists some smooth $\phi$ for which $\psi(h) = \phi(h)\cdot h$. This then puts the second term into the required form.\end{proof}

We have the analogue of \cref{thm:C2:Extension}:
\begin{lem}\label{thm:C3:SPextension}Fix $u' < u'' < a$ and $v' < v'' < b$, and let $\widetilde{\mathcal M} = [0,1)_x\times \R_t \times \R_y \times S_{0,0,0}$ be given the structure of a doubly-foliated manifold with base $\R\times \R$ with projections $(t,y)$. Then there exists an embedding \[\kappa \: \mathcal M([0,u'']\times [0,v'']) \to \widetilde{\mathcal M} := [0,1)_x\times \R_t \times \R_y \times S_{0,0,0}\] such that $\mathcal H_1 \n \mathcal M ([0,u'']\times [0,v''])$ is a subset of $\{t = -y\}$, and $\mathcal H_2 \n \mathcal M([0,u'']\times [0,v''])$ is a subset of $\{t = y\}$. Furthermore, $\kappa^\ast x = x$, and hence $d\kappa: {}^\sp T\mathcal M \to {}^\sp T\widetilde{\mathcal M}$. Also, there is an extension of $g$ to a short-pulse metric $\widetilde{g}$ over $\widetilde{M}$ (i.e.\ $\kappa^\ast \widetilde{g} = g$) for which the hyperplanes $\{t=\pm y\}$ are null, the hyperplanes $\{t = \text{const.}\}$ are spacelike, $\pa_t$ is timelike and future-oriented, and such that $\widetilde{g} = -dt^2 + dy^2 + x^{-2}\mathring{\slash{k}}$ for $t$ or $y$ large enough (here $\mathring{\slash{k}}$ is any fixed Riemannian metric on $S_{0,0,0}$).\end{lem}
\begin{proof}Using \cref{thm:C3:globaltriv}, we may write 
\[\mathcal M([0,u'']\times[0,v'']) \iso [0,1)\times [0,u'']\times[0,v'']\times S.\] The embedding now comes from defining $t$ and $y$ via $u = 1/2(t-y)$, $v = 1/2(t+y)$. The remainder of the proof is the same as \cref{thm:C2:Extension}, and is omitted. That $d\kappa$ acts between the short-pulse bundles is easily seen by picking a basis of vector fields
\[\pa_u, \pa_v, x\pa_{\theta^1},\ldots, x\pa_{\theta^n}\] ($n = \dim S$).
\end{proof}

We are now ready to prove \cref{thm:C3:SPexistence}.
\begin{proof}[Proof of \cref{thm:C3:SPexistence}.]
Let $P(h) = -2Q(h) +2\Ric(g)$, so that we are trying to solve $P(h) = 2\Ric(g)$ on $\mathcal M([0,u']\times [0,v']) \n\{x \leq \delta\}$, with $h = 0$ in Taylor series on $\mathcal H_1 \un \mathcal H_2$. We will make a sequence of technical reductions until we may apply \cref{thm:C3:bigtheoremI} and conclude almost immediately. Embed $\mathcal M$ into $\widetilde{\mathcal M}$ and extend $g$ to $\tilde{g}$ as in \cref{thm:C3:SPextension}. Notice this also extends $A$ and $B$ to smooth bundle maps.

By \cref{thm:C3:amenable}, $P(h)$ has the form
\[P(h) = (g+h)^{-1}\nabla^2 h + A(h,\nabla h) \nabla h + B(h,\nabla h)h,\] where $A$ and $B$ are smooth bundle maps. Observe that $P$ maps symmetric tensors to symmetric tensors, since (the extended) operator $Q$ does. Let $\psi$ be a smooth cutoff which is $1$ on a large set in the base (for instance $\psi \equiv 1$ on $\{|t|,|y| \leq 2(a+b)\}$), and set $A_1 = \psi A$, $B_1 = \psi B$, and $P_1$ for the operator obtained from $P$ by replacing $A$ with $A_1$ and $B$ with $B_1$. Observe that $P_1$ maps symmetric tensors to symmetric tensors, since
\[P_1 = \psi P + (1-\psi)(g+h)^{-1}\nabla^2 h\] and both summands map symmetric tensors to symmetric tensors.

Set $\mathbf{e} = dt^2+dy^2 + x^{-2}\mathring{\slash{k}}$ (where $\mathring{\slash{k}}$ is the same fibre-metric as $\slash{g}$ is for $t,y$ large). Let $D$ denote the Levi-Civita connection of $\mathbf{e}$. Then $V = \nabla-D$ is a short-pulse tensor, being both a tensor on the ordinary tangent bundle, and a smooth map between sections of ${}^\sp T^0_2\widetilde{\mathcal M}$ and ${}^\sp T^0_3 \widetilde{\mathcal M}$. The connections $D$ and $\nabla$ are by definition equal at infinity, so $V$ is compactly supported. 

Then we may write
\[P_1(h) = (g+h)^{-1}D^2 h + (2V+A_1(h,Dh+Vh))Dh + (DV + V\cdot V + B_1(h,Dh+Vh))h,\]
Set $A_2(\xi,\eta) = 2V + A_1(\xi,\eta + V\xi)$, $B_2(\xi,\eta) = DV + V\cdot V + B_1(\xi,\eta + V\xi)$.

The map $h \mapsto (g+h)^{-1}$ is slightly problematic from our perspective, since it only maps small \emph{symmetric} tensors $h$ to Lorentzian cotensors (and additionally is not the constant cometric $-\pa_t^2 + \pa_y^2 + x^2\slash{k}^{-1}$ outside a compact subset of the base). To remedy this, replace it with $\phi(h) = (g+\psi/2(h+h^\ast))^{-1}$, where $h^\ast_{ab} = h_{ba}$, and replace $P_1$ with
\[P_2(h) = \phi(h)D^2 h + A_2(h,Dh)Dh + B_2(h,Dh)h.\] Notice that $P_2(h)$ maps symmetric tensors to symmetric tensors (since $P_2(h) - P_1(h) = ((g+\phi h)^{-1}-g^{-1})D^2h$ whenever $h$ is symmetric) and $P_2(h) = P_1(h) = P(h)$ on $\{\psi \equiv 1\}$ whenever $h$ is symmetric.

Let us define a function $\rho_1$ on $-\infty < u < a,-\infty < v < b$ inside of $\tilde{\mathcal M}$ by
\[\rho_1= \begin{cases}2\Ric(g) & 0 \leq u < a, 0 \leq v < b\\
0 & \text{otherwise}.\end{cases}\]
\begin{figure}[htbp]
\centering
\includegraphics{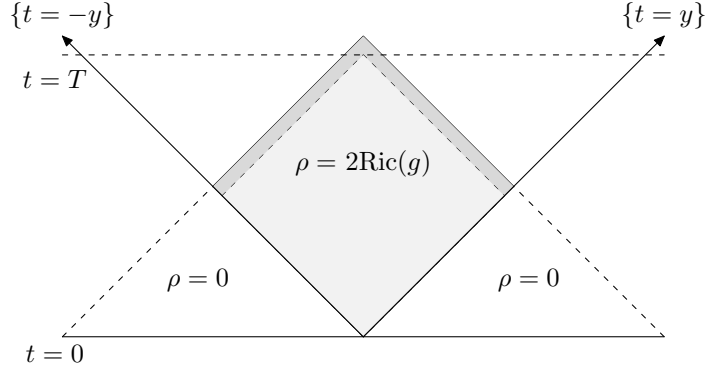}
\caption{The proof of \cref{thm:C3:SPexistence}. The dark shaded region corresponds to $\mathcal M([0,u'']\times[0,v''))$. The light shaded region corresponds to $\mathcal M([0,u']\times[0,v'])$ (i.e.\ the domain of existence of the solution $h$ of the characteristic initial value problem). Compare \cref{fig:C2:shorttime}.}
\end{figure}

Since by assumption $\Ric(g) = 0$ in Taylor series at $\mathcal H_1 \un \mathcal H_2$ and at $x=0$, $\rho_1$ is smooth on its domain of definition. Multiply $\rho_1$ by a smooth cutoff which is identically $1$ on $\mathcal M([0,u'],[0,v'])$, with support inside $-\infty < u < a, -\infty < v < b$, to obtain a function $\rho$. Then $\rho$ is smooth of compact support on all of $\tilde{ \mathcal M}$ and $0$ in Taylor series at $x=0$. Now let us consider the equation
\begin{equation}\label{eq:C3:fineq}P_2(h) = \rho\end{equation} posed on $\tilde{\mathcal M}\n \{t \geq 0\}$ with trivial Cauchy data on $\{t = 0\}$.
We will solve this using \cref{thm:C3:bigtheoremI}.
Shrinking the domain of $\phi$ to a small neighbourhood of the zero section, the triple $(\phi,A_2,B_2)$ is admissible in the sense of \cref{def:C3:admissible}, so by \cref{thm:C3:bigtheoremI} applied to \eqref{eq:C3:fineq} for time $T = u'+v'+2$, we have found $\delta > 0$ and a short-pulse tensor $h$ solving $P_2(h) = \rho$ for $0 \leq t \leq T$, and $x \leq \delta$, and in particular in $\mathcal M([0,u']\times [0,v'])\n \{x < \delta\}$. 

As in the proof of \cref{thm:C2:existence}, since $P_2(0) = 0 = \rho$, for $t \leq y$ and $t \leq -y$, and $\{t=\pm y\}$ are null for $g+0$, $h$ is $0$ in Taylor series at $(\mathcal H_1 \un \mathcal H_2)\n \mathcal M([0,u']\times [0,v'])$ for $x > 0$, and thus $h$ is $0$ in Taylor series for all $x \geq 0$. By definition $P_2(h) = -2\Ric(g)$ on $\mathcal M([0,u']\times [0,v'])$. If we can show that $h$ is symmetric, then this implies that $P_2(h) = P_1(h) = P(h)$ on $\mathcal M([0,u']\times [0,v'])$, and thus $P(h) = -2\Ric(g)$, i.e.\ $Q(h) = 0$.

Symmetry follows from the usual local well-posedness theory for ordinary quasilinear hyperbolic equations, as we now explain. For each $x > 0$, we claim that the set $\mathcal S$ of $t_0 \in [0,T]$ for which $h_x|_{\{t = t_0\}}$ and $\pa_t h_x|_{\{t = t_0\}}$ are symmetric is open, closed, and non-empty. Indeed, closed is obvious, and non-empty is true since $0 \in \mathcal S$. For open, since $P_2$ maps symmetric tensors to symmetric tensors, the usual local well-posedness theory provides for each $t_0 \in \mathcal S$ an $\epsilon > 0$ and a symmetric solution $h'_x$ to $P_2(h'_x)= \rho_x$ for $t \in (t_0-\epsilon,t_0+\epsilon)\n[0,T]$, with initial data $(h'_x,\pa_t h'_x)|_{\{t = t_0\}} = (h_x,\pa_t h_x)|_{\{t = t_0\}}$. Since solutions of $P_2(h_x) = \rho_x$ with given data are unique, it follows that $h_x = h'_x$ is symmetric for $t \in (t_0-\epsilon,t_0+\epsilon)\n[0,T]$.\footnote{The reason why we can deduce symmetry immediately for the ordinary local well-posedness problem but not the short-pulse problem is how one interprets $P$. For the ordinary version, one may write $P(h) = (g+h)^{-1}D^2 h + F(h,\nabla h)$, for some smooth $F$ mapping a symmetric tensor $h$ to another symmetric tensor. Thus, we can interpret $P$ as acting on the bundle of symmetric tensors and obtain a solution in that bundle. In the short-pulse case, we have split up $F$ in order to pull out factors $h$ and $Dh$ to write $P$ in a different fashion, which was necessary to prove \cref{thm:C3:bigtheoremI}. In this splitting one must exit the space of symmetric tensors. For instance, $A(h,\nabla h)$ does not necessarily map the covariant derivative of a symmetric tensor to a symmetric tensor, even if $h$ is symmetric; only $A(h,\nabla h)\nabla h$ is a symmetric tensor. One can also not recover symmetry by looking at the adjoint, since it is not true that $P(h^\ast) = (Ph)^\ast$.}\end{proof}

\section{The formation of trapped surfaces}
\label{C:C3:FormTrapped}
Using the outline outline discussed in \cref{C:C3:ansatz}, we can use \cref{thm:C3:SPCIVP} to prove the long-time existence portion of \cref{thm:C3:ChristMain}. Since the range of $v$ is shrunk when applying \cref{thm:C3:SPCIVP}, we will need, for technical reasons, to slightly extend our data to $v > 1$.

Let us set (for $\epsilon > 0$ small) $\tilde{M} = [0,1)_u\times [0,1+\epsilon)_v \times S^2$ and $\mathcal {\tilde{M}} = [0,1)_x\times \tilde{M}$. Observe that $\tilde{\mathcal M}$ is a parametrized doubly-foliated manifold. Set $N = \pa_u$. We will regard $\tilde{\mathcal M}$ as an extension of the manifold $\mathcal M$ introduced in \cref{C:C3:ansatz}. Set $\tilde{\mathcal H}_1 = \{v = 0\}$ and $\tilde{\mathcal H}_2 = \{u = 0\}$.

Let us extend the short-pulse tensor $\T$ from $[0,1]$ to $[0,1+\epsilon]$, and use it to extend the initial data of the short-pulse ansatz from $M_\delta$ to $\tilde{M}_\delta = \{0 \leq u < 1, 0 \leq \und{u} \leq \delta(1+\epsilon)\}$.

We may also extend the blowdown map $\beta: \mathcal M \to \mathcal R$ to a blowdown map
\[\beta: \tilde{\mathcal M} \to \tilde{\mathcal R} = \{(p,\delta) \in \R^4\times [0,1)\: p \in \tilde{M}_\delta\}.\]

Now, let us set up a characteristic initial value problem for a short-pulse metric on $\mathcal {\tilde{M}}$ which comes from the pulling back the short-pulse ansatz via $\beta$ by providing characteristic initial data for a short-pulse metric $k$. The short-pulse tensor $\T$, extended to be constant in $x$, is naturally a short-pulse fibre tensor defined over $\tilde{\mathcal H}_2$. Let us provide initial data for \cref{thm:C3:SPCIVP} as follows:

\begin{romanumerate}
\item $\widetilde{L} = \pa_v$, a vector field on $\tilde{\mathcal H}_2$;
\item $\slash{k} = x^{-2}\mathring{\slash{g}}$, a fibre-metric on $S_{0,0}$
\item $\hat{\slash{k}} = x^{-2}(1-u)^2\mathring{\slash{g}}$ over $\tilde{\mathcal H}_1$, and $\hat{\slash{g}} = x^{-2}\mathring{\slash{g}}\exp(x\T)$, metrics over $\tilde{\mathcal H}_1 \un \tilde{\mathcal H}_2$;
\item $f_1 = -2$, $f_2 = 2x^2$, functions on $S_{0,0}$;
\item $W = 0$, a fibre-vector field on $S_{0,0}$.
\end{romanumerate}

Assuming the initial data are regular in the sense of \cref{def:C3:regular} (at least for $x$ small), \cref{thm:C3:SPCIVP} says that for any $u^\ast < 1$ and $v^\ast < 1+\epsilon$, there is short-pulse metric $k$ solving $\Ric(k) = 0$ on $\mathcal M([0,u^\ast]\times[0,v^\ast])\n \{0 \leq x \leq \delta\}$ in a double null gauge with $x,u,v,N$, and inducing this initial data. In particular, we may choose $v^\ast = 1$. The point is the following clear lemma:
\begin{lem}\label{thm:C3:relation}Set $g_\delta = x^2\beta_\ast k_x|_{\{v \leq 1\}}$. If $k$ solves the short-pulse characteristic initial value problem for $\{0 \leq x \leq x_0\}$ with the previous initial data, then $g_\delta$ solves the characteristic initial value problem with the short-pulse data in \cref{C:C3:ansatz} for $\{0 < \delta \leq x_0^2\}$.\end{lem}

Thus, if the initial data are regular, we have proved the existence part of \cref{thm:C3:ChristMain}. Assuming this for now, we have the corollary of \cref{thm:C3:SPCIVP} and \cref{thm:C3:SPTSI}, which establishes the notion of smoothness in the parameter.
\begin{cor}\label{thm:C3:smoothnes}The blowdown map $\beta$ desingularizes the family $g_\delta$, i.e.\ $\beta^\ast g_\delta$ form a smooth family of short-pulse metrics, whose Taylor series at $x=0$ are in principle computable. In other words, if one writes in canonical coordinates
\[g_{x^2} = -2\Omega^2(du\otimes d\und{u} + d\und{u}\otimes du) + \slash{g}_{ij}(d\theta^i-f^i)\otimes(d\theta^j-f^j),\]
then $\Omega,\slash{g}_{ij}, f_{ij}$ are smooth functions of $u,\und{u}/\delta,\theta^i$ for which there exists an algorithm to compute them to arbitrarily high order in $x$.\end{cor}

Now let us show regularity.
\begin{lem}[Lemma~2.2 in \cite{ChrForm}]\label{thm:C3:isreg}The initial data are regular, at least if $x$ is small enough.\end{lem}
\begin{proof}
With $\slash{k} = e^{\Phi}\hat{\slash{k}}$, recall we need to show that it is possible to solve
\begin{align*}
2N^2\Phi + (N\Phi)^2 + \Tr(\hat{\slash{k}}^{-1}\Lie_N\hat{\slash{k}})N\Phi + \frac{1}{2}|\Lie_N\hat{\slash{k}}|^2_{\hat{\slash{k}}} + N\Tr(\hat{\slash{k}}^{-1}\Lie_N\hat{\slash{k}})&= 0\\
2{\widetilde{L}}^2\Phi + (\widetilde{L}\Phi)^2 + \Tr(\hat{\slash{k}}^{-1}\Lie_{\widetilde{L}}\hat{\slash{k}}){\widetilde{L}}\Phi + \frac{1}{2}|\Lie_{\widetilde{L}}\hat{\slash{k}}|^2_{\hat{\slash{k}}} + {\widetilde{L}} \Tr(\hat{\slash{k}}^{-1}\Lie_{\widetilde{L}}\hat{\slash{k}})&= 0,\end{align*}
along $\tilde{\mathcal H}_1$ and $\tilde{\mathcal H}_2$, respectively, with initial data $\Phi|_{S_{0,0}} =0$, and $N\Phi_{S_{0,0}}$ and ${\widetilde{L}}\Phi|_{S_{0,0}}$ chosen so that $\slash{\tr}\Lie_N \slash{k} = 2f_1= -4$, $\slash{\tr}\Lie_{\widetilde{L}} \slash{k} = 2f_2= 4x^2$. One easily checks that $\Phi \equiv 0$ solves the equation along $\tilde{\mathcal H}_1$ (alternatively it is the same as Minkowski data there, so must satisfy it).

For the second equation, observe that $\Tr(\hat{\slash{k}}^{-1}\Lie_{\widetilde{L}}\hat{\slash{k}}) = 0$. Indeed, by Jacobi's identity,
\begin{align*}
\Tr(\hat{\slash{k}}^{-1}\Lie_{\widetilde{L}}\hat{\slash{k}}) &= \Tr((\hat{\slash{k}}^{-1}\mathring{\slash{g}})(\Lie_{\widetilde{L}}(\mathring{\slash{g}}^{-1}\hat{\slash{k}})))\\
&= \Lie_{\widetilde{L}} |\det \mathring{\slash{g}}^{-1}\slash{k}|\\
&= \Lie_{\widetilde{L}} |\det e^{\mathbf{T}}|\\
&= \Lie_{\widetilde{L}} \exp(\Tr(\mathbf{T}))\\
&= \Lie_{\widetilde{L}} \exp(0)= 0\end{align*}
 Thus the second equation reduces down to\[2{\widetilde{L}}^2\Phi + ({\widetilde{L}}\Phi)^2 + \frac{1}{2}|\Lie_{\widetilde{L}}\hat{\slash{k}}|^2_{\hat{\slash{k}}} = 0,\] with initial data $\Phi|_{S_{0,0}} = 0$, ${\widetilde{L}}\Phi|_{S_{0,0}} = 2x^2$.

Let us set $\Psi = {\widetilde{L}}\Phi/x^2$, which turns the equation into
\begin{equation}\label{eq:C3:regODE}{\widetilde{L}}\Psi + \frac{x^2}{2}\Psi^2 + \frac{1}{4x^2}|\Lie_{\widetilde{L}}\hat{\slash{k}}|^2_{\hat{\slash{k}}},\end{equation}
with $\Psi|_{S_{0,0}} = 2$.
From the definition of $\exp(x\T)$ via a power series, for instance, it is clear that $|\Lie_{\widetilde{L}}\hat{\slash{k}}|^2_{\hat{\slash{k}}} \lesssim x^2$ uniformly, so the last term is bounded by $C$, for some $C > 0$. Thus if there exists $B > 0$ and $0 \leq T \leq 1+\epsilon$ such that $|\Psi| \leq B$ for $0 \leq v \leq T$, it follows from \eqref{eq:C3:regODE}
\[|\Psi(v)| \leq 2 + \frac{x^2 B^2}{2}+TC.\] In particular, if $B := 8 + 4TC$, then if $x$ is small enough,
\[\sup_{0 \leq v \leq T}|\Psi(v)| \leq B \Rightarrow \sup_{0 \leq v \leq T}|\Psi(v)| \leq \frac{1}{2}B.\] Starting with $|\Psi|_{S_{0,0}}| = 2$, we may use a bootstrap argument to show that $\Psi$ exists on all of $0 \leq v \leq 1+\epsilon$ and hence so does $\Phi$.\end{proof}

Next, we show in the following two lemmas the condition, for a fixed $\delta_0 > 0$, that $g_{\delta_0}$ has no trapped surfaces on $\tilde{\mathcal H}_1 \un \tilde{\mathcal H}_2$ is true so long as the mean curvature $\frac{1}{2}\slash{\tr}\Lie_{L} \slash{g}_{\delta_0}$\footnote{Recall that $L$ is the outgoing null vector field associated with $g_{\delta_0}$, i.e.\ $L = \widetilde{L}/\delta_0$.} is strictly positive on $\tilde {\mathcal H}_2$. In particular, this reduces the question about trapped surfaces on the initial hypersurfaces to the easier question about trapped fibred spheres on $\tilde {\mathcal H}_2$.\footnote{Notice in particular that there are no trapped surfaces in $\widetilde{\mathcal H}_1$ unconditional on the the mean curvatures of the fibred sphere.} We first treat the case of a codimension-two surfaces which like $\tilde{\mathcal H}_2$, and then the case of one in $\tilde{\mathcal H}_1$.\footnote{The codimension is taken with respect to the ambient manifold $\tilde{M}_{\delta_0}$; they are of course codimension $1$ in $\tilde{\mathcal H}_i$ ($i = 1,2$).}

Fix $\delta_0 > 0$. Let $P$ denote a compact codimension two spacelike submanifold in $\tilde{\mathcal H}_2 \n \{\delta = \delta_0\}$, and let $H^{L}_P$ denote the mean curvature of the embedding of $P$ into $\widetilde{M}_{\delta_0}$ with the metric $g_{\delta_0}$, in the direction $L$.\footnote{$L$ is necessarily $g_{\delta_0}$-orthogonal to $P$ since $P \subseteq \tilde{\mathcal H}_2$.}

\begin{lem}\label{thm:C3:Pno}At a point $p$ of $P$,
\[H^{L}_P = \frac{1}{2} \slash{\tr}\Lie_L \slash{g}_{\delta_0}.\] In particular, if the latter is positive along $\tilde{\mathcal H}_2$, then $P$ is not trapped.\end{lem}

Now, let $Q$ denote a compact codimension-two spacelike submanifold in $\tilde{\mathcal H}_1 \n \{\delta = \delta_0\}$, and let $X$ denote the vector field other than $N$ which is future-directed null vector, and $g_{\delta_0}$-orthogonal to $Q$.\footnote{$N$ is necessarily $g_{\delta_0}$-orthogonal to $Q$ since $Q \subseteq \tilde{\mathcal H}_1$.} Denote by $H_{P}^X$ the mean curvatures of the embedding of $P$ into $\widetilde{M}_{\delta_0}$ in the direction $X$.
\begin{lem}\label{thm:C3:Qno}The mean curvature $H_Q^X > 0$ for at least one point in $Q$. In particular, $Q$ is not trapped.\end{lem}

\begin{proof}[Proof of \cref{thm:C3:Pno}.]
Let $\slash{m}$ denote the metric induced by $g_{\delta_0}$ on $P$.
Let $p \in P$ be arbitrary, and let $\{E_1, E_2\}$ denote a local orthonormal frame for $\slash{g}_{\delta_0}$ near $p$. Then there are smooth functions $a^1, a^2$ such that $V_i = a^iL + E_i$, $i = 1,2$, are tangent to $P$ and are $\slash{m}$ orthonormal. Using the formulas for the connection in \cref{C:A1:connection}, it is easy to check that
\[H^{L}_P = \sum_{i=1}^2 g_{\delta_0}(\nabla_{V_i} L,V_i) = \sum_{i=1}^2 g_{\delta_0}(\nabla_{E_i} L, E_i).\] The last expression is just the mean curvature $\frac{1}{2}\slash{\tr} \Lie_L \slash{g}_{\delta_0}$.
\end{proof}

\begin{proof}[Proof of \cref{thm:C3:Qno}]
Observe that the Minkowski metric $\mathring{g}$ is in a double-null gauge with $u,\und{u}$ and $\pa_u$ on $M_{\delta_0}$, and we may write
\[\mathring{g}-2(du\otimes d\und{u} + d\und{u}\otimes du) + (1-u+\und{u})^2\mathring{\slash{g}},\] where $\mathring{\slash{g}}$ of course acts on the vectors tangent to the $S^2$-fibres. We will first show that $g_{\delta_0}$ and $\mathring{g}$ agree to zeroth order on $\tilde{\mathcal H}_1$, which then implies that $X$ is also $\mathring{g}$-null and $\mathring{g}$-orthogonal to $Q$. We will then show that the mean curvature \[\frac{1}{2}\slash{\tr}\Lie_{L} \slash{g}_{\delta_0}\] agrees with its value in Minkowski space, and then use this to show that mean curvature $H^X_Q$ agrees with the mean curvature of the embedding of $Q$ into Minkowski spacetime in the direction $X$. Since the Minkowski spacetime has no trapped surfaces (by the contrapositive to the Penrose incompleteness theorem, \cref{thm:C1:Penrose}), $H^X_Q > 0$ for at least one point in $Q$.

Let us start by showing that $g_{\delta_0}$ and $\mathring{g}$ agree to zeroth order. Observe that $\Omega \equiv 1$ along $\tilde{\mathcal H}_1$, and we have shown in the proof of \cref{thm:C3:isreg} by finding the conformal factor that $\slash{g}_{\delta_0} = (1-u)^2\mathring{\slash{g}}$. From this, observe that \eqref{eq:C2:two}$=0$ becomes
\[0 = \Lie_N [N,L] -2[N,L],\]
and thus $[N,L] = 0$ along $\tilde{\mathcal H}_1$, since it is zero on $\tilde{\mathcal H}_1 \n \tilde{\mathcal H}_2$. In particular $L = \pa_{\und{u}}$ on $\tilde{\mathcal H}_1$. Thus, along $\tilde{\mathcal H}_1$,
\[g_{\delta_0} \equiv \mathring{g}.\] In particular, $X$ is also $\mathring{g}$-orthogonal to $Q$. Denote by $\mathring{\slash{\tr}}$ the trace with respect to $\mathring{\slash{g}}$, and $\rho = 1-u$. Now observe that \eqref{eq:C2:one}$=0$ becomes
\begin{equation}0 = \Lie_N(\Lie_{\pa_{\und{u}}}\slash{g}_{\delta_0})-\frac{1}{\rho}\Lie_{\pa_{\und{u}}}\slash{g}_{\delta_0} - \frac{1}{2\rho}\left(\mathring{\slash{\tr}}\Lie_{\pa_{\und{u}}}{\slash{g}_{\delta_0}}\right)\mathring{\slash{g}} + \frac{2}{\rho}\Lie_{\pa_{\und{u}}}\slash{g}_{\delta_0} + 2\mathring{\slash{g}}.\end{equation}
Taking $\mathring{\slash{\tr}}$ of this yields
\[N(\mathring{\slash{\tr}}\Lie_{\pa_{\und{u}}}\slash{g}_{\delta_0}) + 4 = 0,\] and hence
\[\mathring{\slash{\tr}}\Lie_{\pa_{\und{u}}}\slash{g}_{\delta_0} = 4\rho + C,\] for some constant $C$. Thus
\[\slash{\tr}\Lie_{\pa_{\und{u}}}\slash{g}_{\delta_0} = 4\rho^{-1} + C\rho^{-2}.\]
By assumption, 
\[f_1 := \frac{1}{2}\slash{\tr}\Lie_{\pa_{\und{u}}}\slash{g}_{\delta_0} = 2\] on $\{\und u = 0, \ u = 1-\rho = 0\}$, and thus $C = 0$ and $\slash{\tr}\Lie_{\pa_{\und{u}}}\slash{g}_{\delta_0} = 4\rho^{-1}$. Observe that this coincides with the value for the induced metric in Minkowski spacetime.

Write $X = L + cN + \Sigma$, where $L = \pa_{\und{u}}$ and $\Sigma \in TS$. For $X$ to be null, $c = |\Sigma|^2/2$ (norm taken with respect to $\slash{g}_{\delta_0} = \mathring{\slash{g}}$). A vector $aN + \Theta$, for $\Theta \in TS$, is orthogonal to $X$ (i.e.\ is tangent to $Q$) if \[a = \frac{1}{2}\mathring{\slash{g}}(\Sigma,\Theta) = \frac{1}{2}\mathring{\slash{g}}_{\delta_0}(\Sigma,\Theta).\] Write $a = \phi(\Theta)$ for this linear map for short.

Now, let $\phi(\Theta_i)N + \Theta_i$, $i=1,2$ be two arbitrary tangent vectors to $Q$.
Then, using the formulas in \cref{C:A1:connection}, and recalling that $[N,L] = 0$ on $\tilde{\mathcal H}_1$ because $N = \pa_u$, $L = \pa_{\und{u}}$ there, ($\nabla^1$ and $\slash{\nabla}$ denoting the Levi-Civita connections with respect to $g_{\delta_0}$ and $\slash{g}_{\delta_0} = (1-u)^2\mathring{\slash{g}})$, respectively)
\[g_{\delta_0}(\nabla^1_{\phi(\Theta_i)N + \Theta_i} L,\phi(\Theta_j)N + \Theta_j) = \chi^L(\Theta_i,\Theta_j) + \phi(\Theta_i)(\Theta_j \Omega^2) - \phi(\Theta_j)(\Theta_i \Omega^2),\]
where $\chi^L$ denotes the second fundamental form of the embedding of a fibred-$S^2$ into $\widetilde{M}_{\delta_0}$ under the metric $g_{\delta_0}$. In particular, the last two terms vanish when taking the trace with respect to $\slash{g}_{\delta_0}$, and so the trace of the left-hand side is just the mean curvature
\[\frac{1}{2}\slash{\tr}\Lie_{\pa_{\und{u}}}\slash{g}_{\delta_0} = \frac{2}{\rho},\]
which agrees with the mean curvature of the embedding into Minkowski space. By the same argument, with $\nabla^2$ denoting the Levi-Civita connection of $\mathring{g}$, it follows that the traces with respect to of $\slash{g}_{\delta_0} = (1-u)^2\mathring{\slash{g}}$ of
\[(\Theta_1,\Theta_2) \mapsto g_{\delta_0}(\nabla^2_{\phi(\Theta_i)N + \Theta_i} L,\phi(\Theta_j)N + \Theta_j)\]
is also $\frac{2}{\rho}$.

Since $g_{\delta_0}$ and $\mathring{g}$ agree on $\tilde{\mathcal H}_1$, the second fundamental forms in the direction $N$ associated to $\mathring{\slash{g}}$ and $(1-u)^2\mathring{\slash{g}}$ agree. Write $\chi^N$ for this tensor. Again using the formulas in \cref{C:A1:connection}, one deduces that
\begin{align*}g_{\delta_0}(\nabla^1_{\phi(\Theta_i)N + \Theta_i} |\Sigma|^2/2N + \Sigma,&\phi(\Theta_j)N + \Theta_j) = \phi(\Theta_i)\chi^N(\Sigma,\Theta_j) - \phi(\Theta_j)\chi^N(\Sigma,\Theta_j)\\
&+\frac{|\Sigma|^2}{2}\chi^N(\Theta_i,\Theta_j) + \phi(\Theta_i)\slash{g}_{\delta_0}([N,\Sigma],\Theta_j) + \slash{g}_{\delta_0}(\slash{\nabla}_{\Theta_i}\Sigma,\Theta_j),\end{align*}
with the same formula holding for
\[\mathring{g}(\nabla^2_{\phi(\Theta_i)N + \Theta_i} |\Sigma|^2/2N + \Sigma,\phi(\Theta_j)N + \Theta_j)\]
(since the metrics on the fibred spheres are the same in both cases).
Thus, the traces of both with respect to $\slash{g}_{\delta_0}$ agree. 

Writing $X = L + (|\Sigma|^2/2 N + \Sigma)$ and using the previous two statements, it follows that the trace with respect to $\slash{g}_{\delta_0} = (1-u)^2\mathring{\slash{g}}$ of
\[(\Theta_i,\Theta_j) \mapsto g_{\delta_0}(\nabla^1_{\phi(\Theta_i)N + \Theta_i} X,\phi(\Theta_j)N + \Theta_j)\] and
\[(\Theta_i,\Theta_j) \mapsto \mathring{g}(\nabla^2_{\phi(\Theta_i)N + \Theta_i} X,\phi(\Theta_j)N + \Theta_j)\] 
coincide. Write $\slash{m}$ for the metric induced by $\mathring{g} = g_{\delta_0}$ on $Q$. If $\Theta_1,\Theta_2$ are $\slash{g}_{\delta_0}$-orthonormal, then $\Phi_i = \phi(\Theta_i)N + \Theta_i$, $i = 1,2$, are also orthonormal with respect to $\slash{m}$. It follows that the previous traces agree with the traces with respect to $\slash{m}$ of
\[(\Phi_i,\Phi_j) \mapsto g_{\delta_0}(\nabla^1_{\Phi_i} X,\ \Phi_j)\] 
and
\[(\Phi_i,\Phi_j) \mapsto \mathring{g}(\nabla^2_{\Phi_i} X,\Phi_j),\]
respectively, and so they coincide. Thus the mean curvature $H_Q^X$ agrees with its value in Minkowski spacetime.
\end{proof}

Now, \cref{thm:C3:relation} implies that the signs of the mean curvatures $\frac{1}{2}\slash{\tr}\Lie_{\widetilde{L}}\slash{k_x}, \ \frac{1}{2}\slash{\tr}\Lie_{N}\slash{k_x}$ are the same as the signs of $\frac{1}{2}\slash{\tr}\Lie_{L}\slash{g}_{x^2},\ \frac{1}{2}\slash{\tr}\Lie_{N}\slash{g}_{x^2}$, respectively so to complete the proof of \cref{thm:C3:ChristMain}, we just need to prove:
\begin{prop}\label{thm:C3:computetrapped}Set 
\[U = \{(u,t,\theta) \in [0,1)\times [0,1]\times S^2 \: 1-\frac{1}{4}\mathbf{E}(t,\theta) < u \leq u^\ast\}.\]
If $K$ is any compact subset of $U$, then for small enough $x$, $\slash{\tr}\Lie_{\widetilde{L}}\slash{k},\ \slash{\tr}\Lie_{N}\slash{k} < 0$ on all of $K$.

If $\sup_{\theta \in S^2} \mathbf{E}(1,\theta) < 4$, then choosing $x$ smaller, $\slash{\tr}\Lie_{\widetilde{L}}\slash{k} > 0$ on $\tilde{\mathcal H}_2\n \{ 0 \leq v \leq 1\}$.\end{prop}

To prove this, we will compute the Taylor series of $k$ up until order $2$. This is an easy although somewhat lengthy computation which follows the steps outlined in the proof of \cref{thm:C3:SPTSI} in \cref{C:C3:SPTaylor}.

We use the notation $O(x^j)$ to denote a quantity equal to $x^j$ times a smooth section of an appropriate bundle.

\begin{proof}We will show
\begin{subequations}
\label{eq:C3:MC}
\begin{align}
\slash{\tr}\Lie_N \slash{k} &= -\frac{4}{1-u} + O(x)\\
\slash{\tr}\Lie_{\widetilde{L}} \slash{k} &= \frac{4x^2}{1-u} - \frac{x^2}{(1-u)^2}\mathbf{E}(v,\bullet) + O(x^3).\end{align}
\end{subequations}
This is sufficient to complete the proof, as we now indicate. 

On $\tilde{\mathcal H}_2$, $u= 0$, so 
\[\slash{\tr}\Lie_{\widetilde{L}} \slash{k} = x^2(4-\mathbf{E}(v,\theta)) + O(x^3).\] Since $4-\sup_{\theta \in S^2} \mathbf{E}(1,\theta) > 0$, if $0 \leq v \leq 1$, then $4-\sup_{\theta \in S^2}\mathbf{E}(v,\theta) > 0$. Thus $\slash{\tr}\Lie_{\widetilde{L}} \slash{k} > 0$ for $x$ small.

Next,
\[\frac{4x^2}{1-u} - \frac{x^2}{(1-u)^2}\mathbf{E}(v,\bullet) < 0\] on $U$, so taking $x$ small $\slash{\tr}\Lie_{\widetilde{L}} \slash{k} < 0$ on $K$. Similarly, $\slash{\tr}\Lie_{N} \slash{k} < 0$ everywhere provided $x$ is small enough.

Let us now solve in Taylor series enough to deduce \eqref{eq:C3:MC}. \Cref{thm:C3:SPTSI} shows that Taylor series are unique. Following the notation in that section, set $\slash{h} = x^2\slash{k}$, $\ovelrine{L} = x^{-1}(\pa_v - \widetilde{L})$, and $\omega = \log \Omega$. Write $\slash{h}_i$, $\overline{L}_i$, $\omega_i$ for the coefficient of $x^i$ in the expansion of $\slash{h}$, $\overline{L}$, $\omega$ at $x=0$, respectively, and let $\slash{h}_{(i)}$, $\ovelrine{L}_i$, $\omega_{(i)}$ denote their Taylor expansions, respectively, up to order $x^i$. Observe that $\pa_v$ coincides with $\mathring{L}$ defined by $[N,\mathring{L}] = 0$, $\mathring{L}|_{\widetilde{\mathcal H}_2} = \pa_v$. Since $\widetilde{L} = \pa_v + x\overline{L}$, knowing $\overline{L}_{(1)}$ is sufficient to know $\widetilde{L}$ up to $x^3$ times a smooth vector field. Thus,
\begin{align*}
\slash{\tr}_{\slash{k}}\Lie_N \slash{k} &= \slash{\tr}_{\slash{h}_{(0)}}\Lie_N \slash{h}_{(0)} + O(x)\\
\slash{\tr}_{\slash{k}}\Lie_{\widetilde{L}} \slash{k} &= \slash{\tr}_{\slash{h}_{(2)}}\Lie_{\pa_v} \slash{h}_{(2)} + O(x^3).\end{align*}

It therefore suffices to show that
\begin{align}\label{eq:C3:trace1}\slash{\tr}_{\slash{h}_{(0)}}\Lie_N \slash{h}_{(0)} &= -\frac{4}{1-u}\\
\label{eq:C3:trace2}\slash{\tr}_{\slash{h}_{(2)}}\Lie_{\pa_v} \slash{h}_{(2)} &= \frac{4x^2}{1-u} - \frac{x^2}{(1-u)^2}\mathbf{E}(v,\bullet) + O(x^3).\end{align} We first need to figure out the Taylor series of the initial data at $x=0$. For the rest of the proof, let $\mathring{L} := \pa_v$.

\textbf{Determining the Initial Data: }Let us start with $\mathcal{\widetilde{H}}_1$. To do this, we generally follow the procedure used in \cref{thm:C2:toporder} used to solve exactly on the initial hypersurfaces. We have already found in \cref{thm:C3:isreg} that the conformal factor is $1$. Thus $\slash{h} = (1-u)^2\mathring{\slash{g}}$ along $\mathcal{\widetilde{H}}_1$, and in particular \eqref{eq:C3:trace1} holds. Thus \eqref{eq:C3:two}$=0$ becomes
\[0 = \Lie_N Z -2Z,\] with $0$ initial data. Thus $Z \equiv 0$, and hence $\overline{L} \equiv 0$ is, too. By assumption $\omega \equiv 0$, so we have found everything.

Now let us find the data on $\mathcal{\widetilde{H}}_2$. By assumption, we are given $\widetilde{L} = \pa_v$, $\Omega = 1$ so $\ovelrine{L} \equiv 0$, $\omega \equiv 0$. As mentioned in the proof of \cref{thm:C3:isreg}, 
\[2{\widetilde{L}}^2\Phi + ({\widetilde{L}}\Phi)^2 + \frac{1}{2}|\Lie_{\widetilde{L}}\hat{\slash{h}}|^2_{\hat{\slash{h}}} = 0,\] with initial data ${\widetilde{L}}\Phi|_{S_{0,0}} = 2x^2$.

We may compute a Taylor series for $|\Lie_{\widetilde{L}}\hat{\slash{h}}|^2_{\hat{\slash{h}}}$ using the power series for $\exp$, namely,
\[\exp(x\T) = 1 + x\T + O(x^2).\] Thus,
\[|\Lie_{\widetilde{L}}\hat{\slash{h}}|^2_{\hat{\slash{h}}} = x^2\Tr(\pa_v \T^2) + O(x^3).\] Setting $\Psi = x^{-2}{\widetilde{L}}\Phi$ turns the ODE into
\[{\widetilde{L}}\Psi + \frac{x^2}{2}({\widetilde{L}}\Psi)^2 + \frac{1}{4}\Tr(\pa_v \T^2) = O(x),\] with initial data $\Psi|_{S_{0,0}} = 2$. Thus,
\[\Psi = 2-\frac{1}{2}\mathbf{E} + O(x).\]
Hence,
\[\Phi =2x^2 v - \frac{x^2}{2}\int_0^v \mathbf{E}(s)\ ds + O(x^3),\]
and thus
\begin{align*}
\slash{h} &= e^{\Phi}\mathring{\slash{g}}\exp(x\T)\\
&= \left(1 + 2x^2 v - \frac{x^2}{2}\int_0^v \mathbf{E}(s)\ ds + O(x^3)\right)\mathring{\slash{g}}\left(1 + x\T + \frac{x^2}{2}\T^2\right)\\
&= \mathring{\slash{g}} + x\mathring{\slash{g}}\T + \left(2x^2 v - \frac{x^2}{2}\int_0^v \mathbf{E}(s)\ ds\right)\mathring{\slash{g}} +\frac{x^2}{2} \mathring{\slash{g}}\T^2 + O(x^3).\end{align*}

This gives us all the initial data. 

\textbf{Computing to top order: }Observe that \eqref{eq:C3:one} $=0$ means that $\slash{h}_0 = (1-u)^2\mathring{\slash{g}}$, \eqref{eq:C3:six} $=0$ means that $\omega_0 = 0$, and \eqref{eq:C3:three} $=0$ means that $\overline{L}_0 = 0$. Already this means that $\Lie_N \slash{h}_0 = -2(1-u)\mathring{\slash{g}}$, and hence $\slash{\tr}\Lie_N \slash{h}_0 = -4(1-u)^{-2}$, which is \eqref{eq:C3:trace2}

To solve for the higher order terms, we must linearize. Let $P$ denote the operator as in \cref{thm:C3:linearization}, giving the right-hand side of \eqref{eq:C3:one}, \eqref{eq:C3:three} and \eqref{eq:C3:six} In our case, the components of its linearization, $DP_{(\slash{h}_0,\overline{L}_0,\omega_0)}$, acting on a triple $(\slash{h}', \overline{L}', \omega')$, are:
\begin{subequations}
\begin{align}
\begin{split}
&\Lie_N\Lie_{\mathring{L}}\slash{h}' + \frac{1}{1-u}\Lie_{\mathring{L}} \slash{h}'- \frac{1}{2(1-u)}(\mathring{\slash{\tr}}\Lie_{\mathring{L}} \slash{h}')\mathring{\slash{g}}
\end{split}\\
\begin{split}
&(1-u)^2\Lie_{\mathring{L}}\mathring{\slash{g}}([N,\overline{L}'],\bullet)
\end{split}\\
\begin{split}
&N\mathring{L}\omega' + \frac{1}{4(1-u)^2}N\mathring{\slash{\tr}}\Lie_{\mathring{L}}\slash{h}' + \frac{1}{4(1-u)^3}\mathring{\slash{\tr}}\Lie_{\mathring{L}}\slash{h}'.\end{split}\end{align}
\end{subequations}

Here $\mathring{\slash{\tr}}$ denotes the trace with respect to $\mathring{\slash{g}}$.

\textbf{Computing to first order: }
Observe that $P(\slash{h}_0,\ovelrine{L}_0,\omega_0) = O(x^2)$. Thus the equation
\[(-x^{-1}P(\slash{h}_0+x\slash{h}_1,\ovelrine{L}_0+x\overline{L}_1,\omega_0+x\omega_1))|_{x=0}=DP_{(\slash{h}_0,\ovelrine{L}_0,\omega_0)}(\slash{h}_1,\overline{L}_1,\omega_1)\] becomes
\begin{align*}
\begin{split}
 0&= \Lie_N\Lie_{\mathring{L}}\slash{h}_1 + \frac{1}{1-u}\Lie_{\mathring{L}}\slash{h}_1- \frac{1}{2(1-u)}(\mathring{\slash{\tr}}\Lie_{\mathring{L}} \slash{h}_1)\mathring{\slash{g}}
\end{split}\\
\begin{split}
0 &= (1-u)^2\Lie_{\mathring{L}}\mathring{\slash{g}}([N,\ovelrine{L}_1],\bullet)
\end{split}\\
\begin{split}
0 &= N\mathring{L}\omega_1 + \frac{1}{4(1-u)^2}N\mathring{\slash{\tr}}\Lie_{\mathring{L}}\slash{h}_1 + \frac{1}{4(1-u)^3}\mathring{\slash{\tr}}\Lie_{\mathring{L}}\slash{h}_1.\end{split}\end{align*}

 The initial data for these equations is trivial for $\overline{L}_1$ and $\omega_1$ on $\mathcal{\widetilde{H}}_1\un \mathcal{\widetilde{H}}_2$ and for $\slash{h}_0$ on $\mathcal{\widetilde{H}}_1$ (since $\overline{L}_0$ and $\omega_0$, $\slash{h}_0$ already had the correct initial data exactly). On $\mathcal{\widetilde{H}}_2$, the initial data for $\slash{h}_1$ is $\mathring{\slash{g}}\T$, which we computed above.

Let us make the ansatz that $\Lie_{\mathring{L}}\slash{h}_1$ is $\mathring{\slash{g}}$ tracefree. This is consistent since $\Lie_{\mathring{L}}\slash{h}_1 = \mathring{\slash{g}}\pa_v \T$ initially, and $\T$ is $\mathring{\slash{g}}$ tracefree. Then the first equation becomes a transport equation for $\Lie_{\mathring{L}}\slash{h}_1$, which we can then solve to obtain \[\Lie_{\mathring{L}}\slash{h}_1 = (1-u)\mathring{\slash{g}}\pa_v\T ,\] and hence $\slash{h}_1 = (1-u)\mathring{\slash{g}}\T$, which is still $\mathring{\slash{g}}$-tracefree. Thus, using the other two equations, $\overline{L}_1 = 0$ and $\omega_1 = 0$.

\textbf{Computing to second order: }Of the second-order components, we only need to compute $\slash{h}_2$.
Let us use a shorthand and denote $\T' = \pa_v \T = \Lie_{\mathring{L}}\T$. 

\begin{claim}The first component of $P$ applied to
$(\slash{h}_0 + x\slash{h}_1, \overline{L}_0 + x\overline{L}_1, \omega_0 + x\omega_1)$ is $2x^2\mathring{\slash{g}}$.\end{claim}
\begin{proof}
Recall that $\T$ is $\mathring{\slash{g}}$-symmetric and tracefree. We compute the following quantities in order:
\begin{align*}
\slash{h}_{(1)} &= (1-u)^2\mathring{\slash{g}} + x(1-u)\mathring{\slash{g}}\T\\
\Lie_N \slash{h}_{(1)} &= \mathring{\slash{g}}(-2(1-u)-x\T)\\
\Lie_{\mathring{L}}\slash{h}_{(1)} &= x(1-u)\mathring{\slash{g}}\T'\\
\slash{h}^{-1}_{(1)} &= (1-u)^{-2}\left(1-x(1-u)^{-1}\T\right)\mathring{\slash{g}}^{-1} + O(x^2)\\
\slash{h}^{-1}_{(1)}\Lie_N \slash{h}_{(1)} &= (1-u)^{-2}(-2(1-u) + x\T) + O(x^2)\\
\slash{h}^{-1}_{(1)}\Lie_{\mathring{L}} \slash{h}_{(1)} &= (1-u)^{-2}(x(1-u)\T' - x^2\T\T') + O(x^3)\\
\slash{\tr}\Lie_N\slash{h}_{(1)} &= \frac{-4}{1-u} + O(x^2)\\
\slash{\tr}\Lie_{\mathring{L}}\slash{h}_{(1)} &= -x^2(1-u)^{-2}\Tr(\T\T') + O(x^3)\\
\Lie_{\mathring{L}} \slash{h}_{(1)} \times\Lie_N \slash{h}_{(1)} &= (1-u)^{-2}\mathring{\slash{g}}(-2x(1-u)^2\T' + x^2(1-u)\T'\T) + O(x^3)\\
\Lie_N \slash{h}_{(1)}\times \Lie_{\mathring{L}} \slash{h}_{(1)} &= (1-u)^{-2}\mathring{\slash{g}}(-2x(1-u)^2\T' + x^2(1-u)\T\T') + O(x^3)\\
x^2\slash{\Ric}(\slash{h}_{(1)}) &= 	x^2\mathring{\slash{g}}.\end{align*}
Using these and the facts $\Omega = 1 + O(x^2)$ and $L = \mathring{L} + x\overline{L} = \mathring{L} + O(x^3)$ (coming from $\omega_0 = \omega_1 = 0$, $\overline{L}_0 = \overline{L}_1 = 0$) we may compute the first component of $P$ acting on $(\slash{h}_{(1)},\overline{L}_{(1)},\omega_{(1)})$ as

\begin{equation}\label{eq:C3:secondorder}x^2\mathring{\slash{g}}\left(2+\frac{1}{2(1-u)}\Tr(\T\T')-\frac{1}{2(1-u)}(\T'\T+\T\T')\right).\end{equation}

Observe $\T^2 + \det \T = 0$, since $\T$ is a fibrewise tracefree linear map on a $2$-dimensional vector space.
Thus,
\[\T\T' + \T'\T = \pa_v(\T^2) = -\pa_v \det \T\] is in fact a scalar, and so $\T\T' + \T'\T = \frac{1}{2}\Tr(\T'\T+\T'\T) = \Tr(\T\T')$. So \eqref{eq:C3:secondorder} simplifies to just $2x^2\mathring{\slash{g}}$.\end{proof}

Thus the equation for $\slash{h}_2$,
\[(-x^{-2}P_1(\slash{h}_1,\ovelrine{L}_1,\omega_01)|_{x=0} = D{P_1}_{(\slash{h}_0,\ovelrine{L}_0,\omega_0)}(\slash{h}_2,\overline{L}_2,\omega_2)\]
becomes
\begin{equation}\label{eq:C3:label}-2\mathring{\slash{g}} = \Lie_N\Lie_{\mathring{L}}\slash{h}_2 + \frac{1}{1-u}\Lie_{\mathring{L}}\slash{h}_2- \frac{1}{2(1-u)}(\mathring{\slash{\tr}}\Lie_{\mathring{L}} \slash{h}_2)\mathring{\slash{g}}.\end{equation}

 The initial data for $\slash{h}_2$ is
\[\left(2v - \frac{1}{2}\int_0^v \mathbf{E}(s)\ ds\right)\mathring{\slash{g}} + \frac{ 1}{2}\mathring{\slash{g}}\T^2.\] Observe $\T^2 = \frac{1}{2}\Tr(\T^2)$, since $\T^2$ is a scalar. Thus on $\widetilde{\mathcal H}_2$
\[\Lie_{\mathring{L}}\slash{h}_2 = \left(2 - \frac{ 1}{2}\mathbf{E} + \frac{ 1}{4}\pa_v \Tr(\T^2)\right)\mathring{\slash{g}}.\]

 The initial data for $\Lie_{\mathring{L}} \slash{h}_2$ is thus a scalar multiple of $\mathring{\slash{g}}$. Thus we make the ansatz $\Lie_{\mathring{L}} \slash{h}_2 = \kappa\mathring{\slash{g}}$ for $\kappa$ a scalar. From \eqref{eq:C3:label}, the equation for $\kappa$ is just $-2 = N\kappa$, and so
\[\Lie_{\mathring{L}} \slash{h}_2 = \left(2(1-u)- \frac{ 1}{2}\mathbf{E} + \frac{ 1}{4}\Tr(\T^2)\right)\mathring{\slash{g}},\]
and so
\[\slash{h}_2 = \left(2(1-u)v- \frac{ 1}{2}\int_0^v\mathbf{E}(s)\ ds + \frac{ 1}{4}\Tr(\T^2)\right)\mathring{\slash{g}}.\]
Let us denote by $S$ the big factor out front.

\textbf{Completing the computation:} We have shown that
\[\slash{h}_{(2)} = (1-u)^2\mathring{\slash{g}} + x(1-u)\mathring{\slash{g}}\T + x^2F\mathring{\slash{g}}.\]

In particular,
\begin{align*}
\slash{h}_{(2)}^{-1} &= (1-u)^{-2}\left(1-x(1-u)^{-1}\T\right)\mathring{\slash{g}}^{-1} + O(x^2)\\
\Lie_{\pa_v}\slash{h}_{(2)} &= x(1-u)\mathring{\slash{g}}\T' + x^2\pa_v S\mathring{\slash{g}}\\
\slash{\tr}\Lie_{\mathring{L}}\slash{h}_{(2)} &= (1-u)^{-2}\Tr(x(1-u)\T' + x^2\pa_v S- x^2\T\T') + O(x^3).\end{align*}
The tensor $\T'$ is tracefree, so the first argument of the trace drops out. By definition,
\[\pa_v S - \T\T' = 2(1-u) - \frac{1}{2}\mathbf{E} + \frac{1}{2}\Tr(\T\T')-\T\T',\]
and so its trace is just
\[4(1-u) - \mathbf{E}.\] Putting it all together,
\[\slash{\tr}\Lie_{\mathring{L}}\slash{h}_{(2)} = \frac{4x^2}{1-u} - \frac{x^2\mathbf{E}}{(1-u)^2} + O(x^3),\]
which is \eqref{eq:C3:trace2}.\end{proof}
\chapter{Real blow up and continuation in formal series}
\label{C:C4}

\section{The blown-up space}
\label{C:C4:blowup}
Although the manifold $\mathcal M$ of \cref{C:C3:ansatz} and \cref{C:C3:FormTrapped} successfully desingularizes
\[\mathcal R = \{0 \leq u < 1, 0 \leq u \leq \delta, \delta \geq 0\} = \{(p,\delta) \in \R^4\times[0,1) \: p \in M_\delta\},\] it is not enough to desingularize the space\footnote{We use the usual definitions $u = 1-(|y|-t)/2$, $\und{u} = (t+|y|)/2$.}
\[\mathcal S = \{0 \leq u \leq 1, 0 \leq \und{u} \leq \delta, \delta \geq 0\} \subseteq \R^4_{(t,y)} \times [0,1)_\delta\] because of the line of cone points $\{u = 1,\ \und{u} = 0\}$. We start by desingularizing this region to a mwc $\widebar{\mathcal M}$. We will perform an iterated blowup, starting with $\R^4_{(t,y)} \times [0,\delta)$, and successively blowing up and taking submanifolds until we have the desired resolved space. Set $\rho = 1-u$.

We first need to resolve the biggest singularity, so let us start by blowing up the point $\{t=0,y=0,\delta = 0\}$ (which is the same as the point $\{\rho = \und{u} = \delta = 0\}$), replacing it with a hemisphere $S^4_+$. Call this new manifold $\mathcal M_0$, and the introduced face $\mathbf{sf}$ (``sphere face''). The second biggest singularity is the line of cone tips, which form the (closure of the lift of the) line $\{t=0,y=0, \delta > 0\}$ (which is the same as the line $\{\rho = \und{u} = 0, \delta > 0\}$). Let us blow it up, replacing it with a cylinder $S^3\times [0,1)$. Call this new manifold $\mathcal M_1$, and the introduced face $\mathbf{cf}$ (``cylinder face'').

One may check that the lift of the region 
\[\{0 \leq \und{u} \leq 1, 0 \leq \rho \leq 1\}\] to $\mathcal M_1$ defines a submanifold which is also a mwc. Call this space $\mathcal M_2$. Let us redefine $\mathbf{sf}$ and $\mathbf{cf}$ to be their intersections with $\mathcal M_2$ (see \cref{fig:C4:M2}). The manifold $\mathcal M_2$ now has a number of additional faces: the face corresponding to $\und{u} = 0$ (really the closure of the lift of $\{\und{u} = 0, \rho > 0, \delta > 0\}$) and the face corresponding to $\rho = 1$ (really the closure of the lift of $\{\rho = 1, \delta > 0\}$). Call these faces $\mathbf{lf}$ and $\mathbf{rf}$, respectively (``left face'' and ``right face,'' respectively) since these correspond to the initial hypersurfaces for the characteristic initial value problem.\footnote{In the previous two chapters, these were called $\mathcal H_1$ and $\mathcal H_2$, respectively. In this and subsequent chapters, we will be less concerned with the characteristic initial value problem, since it has been solved, and will choose to use more traditional notation for mwcs and blowups.} The lift of $\{\und{u} = 1\}$ and (the closure of) the lift of $\{\rho = 0, \und{u} > 0, \delta > 0\}$ are also faces, but will not appear in the final desingularized region so we do not name them.

The region $\mathcal M_2^\circ = \{0 < \rho \leq 1 ,0 \leq \und{u} \leq 1, \ \delta > 0\}$\footnote{We regard the initial surfaces $\mathbf{lf}$ and $\mathbf{rf}$ as well as $\{\und{u} = 1\}$ as ``artificial'' faces, so they are allowed to remain in the interior.} has the structure of a trivial $S^2$-bundle over
\[\{(u,\und{u},\delta)\in \R\times \R\times [0,1)\: 0 \leq u < 1, \ 0 \leq \und{u} \leq \delta, \ 0 < \delta \leq 1\},\] with projection taking a point in $\mathcal M_2^{\circ}$ to $(u,\und{u},\delta)$. This is the same fibre bundle we have been thinking of when we consider a single slice $\mathcal R\n \{\delta = \delta_0 > 0\}$ as a doubly-foliated manifold. We claim that this fibre bundle structure extends to $\mathcal M_2$. Let us denote by $\mathcal N_2$ the putative base of this bundle.

\begin{proof}[Proof of claim]Consider $\R_t\times \R_r^{\geq 0} \times [0,1)_\delta$, and blow up first $\{t = r = \delta = 0\}$ replacing it with a quarter sphere $S^2_{++}$, and then the (closure of the lift of) the line $\{t = r = 0, \delta > 0\}$, replacing it with a half-cylinder cylinder $[0,1)\times S^1_+$. Call the result of this iterated blowup $\mathcal N_1$. Set $\rho = (r-t)/2$, $\und{u} = (r+t)/2$. Like above, we may define $\mathcal N_2$ to be the intersection of $\mathcal N_1$ with $\{0 \leq \und{u} \leq 1, 0 \leq \rho \leq 1\}$. Let us think of $S^2$ as the unit sphere in $\R^3$. Observe that the smooth map
\begin{gather*}(\R_t\times \R_r^{\geq 0}\times [0,1)_\delta)\times S^2_{\theta} \to \R_t\times \R^3_y \times [0,1)_\delta\\
(t,r,\delta,\theta) \mapsto (t,r\theta,\delta)\end{gather*}
restricts to a diffeomorphism 
\[\{0 \leq u < 1, \ 0 \leq \und{u} \leq \delta\} \times (0, 1]\times S^2 \to \mathcal M_2^{\circ}\]
that provides a trivialization of fibre bundle mentioned above associated to the projection map $(u,\und{u},\delta)$ (although we of course need to change coordinates $u = 1-(r-t)/2$, $\und{u} = (r+t)/2$ on the base).

We claim the map extends to a diffeomorphism $\mathcal N_2 \times S^2 \to \mathcal M_2$. Indeed, the only problematic points are those with $r = 0$, but we blew up $\{r = 0\}$ and $r \geq |t|$ on $\mathcal N_2$, so in effect this diffeomorphism introduces ``partial polar coordinates'' on $\mathcal M_2$. This can also be checked carefully in projective coordinates, but we omit this computation for brevity.

\begin{figure}[htbp]
\centering
\includegraphics{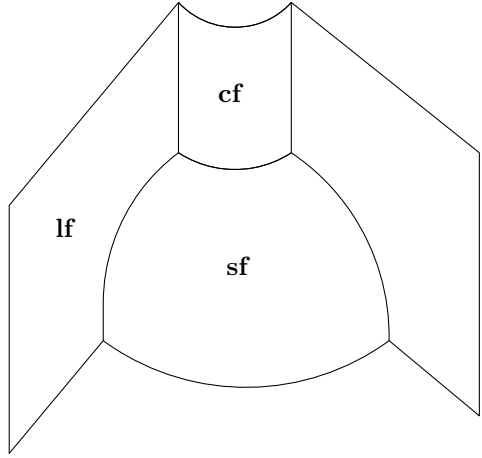}
\caption{A stylized view of $\mathcal M_2$ at a certain angular slice $\theta$. $\mathbf{rf}$ is not drawn, and is a face meeting $\mathbf{lf}$ transversely.}
\label{fig:C4:M2}
\end{figure}

Using this diffeomorphism, we may realize $\mathcal M_2$ as a trivial $S^2$ bundle over $\mathcal N_2$ extending the one on $\mathcal M_2^\circ$.\end{proof}

We may also introduce ``partial'' coordinate systems by introducing coordinates on the base $\mathcal N_2$, but not on the fibres. Let us introduce $(\rho, v_1 = \und{u}/\rho, \eta_1 = \delta/\rho,\theta)$, which are valid away from $\mathbf{cf}$ and (the closure of) the lift of $\{\rho = 0, \und{u} > 0, \delta > 0\}$. The set $\{v_1 = \eta_1\}$, intersects $\mathbf{lf}$ and $\{\eta_1 = 0\}$ at $\{v_1 = \eta_1 = 0\}$, so to resolve this intersection, we blow up the intersection, obtaining a face diffeomorphic to $[0,1]^2 \times S^2$. Call this manifold $\mathcal M_3$. Let us call the introduced face $\mathbf{bf}$ ( for ``bottom face'').

\begin{figure}[htbp]
\centering
\includegraphics{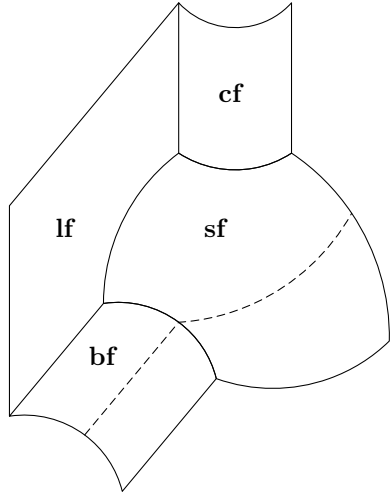}
\caption{A stylized view of $\mathcal M_4$ at a certain angular slice $\theta \in S^2$. The face $\mathbf{rf}$ is not drawn, and is a face meeting $\mathbf{lf}$ and $\mathbf{bf}$ transversely. The dashed line indicates where $\mathbf{ff}$ intersects $\mathbf{sf}$ and $\mathbf{bf}$.}
\label{fig:C4:M4}
\end{figure}

We are really interested in the region $\{\und{u} \leq \delta\}$, which is the same as $\{v_1 \leq \eta_1\}$. The intersection of this with $\mathcal M_3$ is a new mwc, $\mathcal M_4$ and the (closure of the) lift of $\{v_1 = \eta_1 > 0\}$ is a new face, which we call $\mathbf{ff}$ for ``far face'' (see \cref{fig:C4:M4}). Inside $\mathcal M_4$, the coordinates $(\rho, v = v_1/\eta_1 = \und{u}/\delta, \eta_1,\theta)$ form valid coordinates near $\mathbf{bf}$.

The manifold $\mathcal M_4$ is fully desingularized, but it is not desingularized enough with respect to the putative solution metric $g$ to the EVEs with short-pulse data.

\begin{figure}[htbp]
\centering
\includegraphics{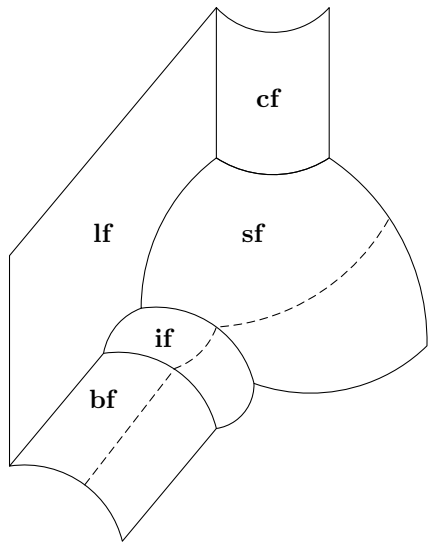}
\caption{Figure~\ref{fig:C1:blowup}, reproduced for the reader's convenience. A stylized view of $\protect\widebar{\mathcal M}$ at an angular slice $\theta \in S^2$. The dashed line is where $\mathbf{ff}$ intersects $\mathbf{bf}$, $\mathbf{if}$ and $\mathbf{sf}$. $\mathbf{rf}$ is not drawn, and is a face meeting $\mathbf{lf}$, $\mathbf{bf}$ transversely.}
\label{fig:C4:Mbar}
\end{figure}

Let us redefine the smooth structure by requiring the square root of a bdf of $\mathbf{sf}$ and the square root of a bdf of $\mathbf{bf}$ to be smooth. Now blow up $\mathbf{sf}\n\mathbf{bf}$, and call the introduced face $\mathbf{if}$ for ``intermediate face.'' We let $\widebar{\mathcal M}$ be this final manifold (see \cref{fig:C4:Mbar}).

Now observe that none of the blowups we performed involved the spherical variables. Thus $\widebar{\mathcal M}$ remains a trivial $S^2$ bundle, over a manifold $\widebar{\mathcal N}$.

Let us list several useful partial coordinate charts on $\widebar{\mathcal M}$ by introducing coordinates only on the base $\widebar{\mathcal N}$. The first set of coordinates on $\widebar{\mathcal M}$ we use is valid near $\mathbf{bf}$, and covers all of $\mathbf{if}$ except for $\mathbf{if}\n\mathbf{sf}$. They are $(\xi = \sqrt{\rho}, v,\eta = \sqrt{\eta_1}/\xi,\theta)$ (here $\theta \in S^2$). In these coordinates, $\xi$ is a bdf of $\mathbf{if}$, $v$ is a bdf of $\mathbf{lf}$, $(1-v)$ is a bdf of $\mathbf{ff}$, $\eta$ is a bdf of $\mathbf{bf}$. 
Near $\mathbf{if}\n\mathbf{sf}$, we may use the coordinates $(\varpi = \sqrt{\eta_1}, \tau = \xi/\varpi, v, \theta)$. In these coordinates, $\varpi$ is a bdf of $\mathbf{if}$, $v$ is a bdf of $\mathbf{lf}$, $\tau$ is a bdf of $\mathbf{sf}$.
Near $\mathbf{sf}\n\mathbf{cf}$, we may use the coordinates $(\sigma = t/r, s= r/\delta, \vepsilon = \sqrt{\delta}, \theta)$, where $\sigma$ is a bdf of $\mathbf{lf}$, $1-\sigma$ is a bdf of $\mathbf{ff}$, $s$ is a bdf of $\mathbf{cf}$ and $\vepsilon$ is a bdf of $\mathbf{sf}$.

We end this section by showing how $\mathcal M$, the desingularization of $\mathcal R$ introduced in \cref{C:C3:ansatz}, embeds naturally into $\widebar{\mathcal M}$. The identity map \[\mathcal M^\circ = \widebar{\mathcal M}^\circ = \{0 \leq u < 1, \ 0 \leq \und{u} \leq \delta, \ 0 < \delta \leq 1\}\] extends to an embedding of $\mathcal M \embeds \widebar{\mathcal M}$ taking a point $(x,u,v,\theta)$ in $\mathcal M$ (in coordinates defined in \cref{C:C3:ansatz}) to
\[(x,u,v,\theta) \mapsto (\xi,v,\eta,\theta) = (\sqrt{1-u}, v, x/(1-u),\theta) \in \widebar{\mathcal M}\n \{u < 1\},\] in the first coordinate system mentioned above. This embedding will be very useful since it allows us to apply the results of \cref{C:C3} to $\widebar{\mathcal M}$ where appropriate. Observe also that $(\mathcal M,u,v)$, when considered as a doubly-foliated manifold (in the sense of \cref{def:C3:doublefoliated}), and $\widebar{\mathcal M}$, when considered as a trivial fibre bundle, have the same fibres. Moreover, the trivialization of $\widebar{\mathcal M}$ extends the trivialization of $\mathcal M$.

\section{The short-pulse bundle and short-pulse metrics}
\label{C:C4:DNG}
In order to clarify what sort of metric we should be solving for on $\widebar{\mathcal M}$, we express the Minkowski metric $\mathring{g}$ in each of the coordinate systems above.

Observe that $\mathring{g}$ is not a metric on $\R^4\times [0,1)$, but on a fixed $\delta$ slice. Therefore, we treat $\mathring{g}$ as a section, not of $\Sym^2(T^\ast\widebar{\mathcal M})$, but of $\Sym^2(T^\ast\widebar{\mathcal M}/\vspan\{ d\delta\})$. Since we are working in the quotient, $d\delta = 0$. This lets us simplify some formulae. For instance, $d\delta = d(\xi^4\eta^2) = 0$, and so $2d\xi/\xi + d\eta/\eta = 0$.

Thus,
\begin{align*}
\mathring{g} &= \xi^6\eta^2[4d\xi/\xi\otimes dv + 4dv\otimes d\xi/\xi + (1+\eta^2\xi^2 v)(\xi\eta)^{-2}\mathring{\slash{g}}]\\
\mathring{g} &= \tau^4\varpi^6[2 d\tau/\tau\otimes dv + 2dv \otimes d\tau/\tau + (1+\varpi^2v)^2\varpi^{-2}\mathring{\slash{g}}]\\
\mathring{g} &= \vepsilon^4\sigma^2[ds\otimes ds + s(d\sigma/\sigma\otimes ds + ds\otimes d\sigma/\sigma) + (1+s^2) d\sigma/\sigma \otimes d\sigma/\sigma + \mathring{\slash{g}}]
\end{align*}
in each of the respective coordinate systems,
where $\mathring{\slash{g}}$ is the round metric on the fibred $S^2$.

This suggests what should be the correct notion of short-pulse bundle $^\sp T\widebar{\mathcal M}$. Recall that $\widebar{\mathcal M} \iso \widebar{\mathcal N} \times S^2$.
\begin{defn}Let $^\sp T\widebar{\mathcal M}$ be the vector bundle over $\widebar{\mathcal M}$ whose smooth sections are smooth sections $X$ of $T\widebar{\mathcal M}$ satisfying:
\begin{romanumerate}
\item $d\delta X = 0$;
\item at $\mathbf{bf}$, $X$ is tangent to the fibres of the projection onto the $S^2$ factor of $\widebar{\mathcal M}$;
\item at $\mathbf{if}$, $X$ is tangent to $\mathbf{if}$ and tangent to the projection of the fibres onto the $S^2$ factor;
\item at $\mathbf{sf}$, $X$ is tangent to $\mathbf{sf}$;
\item at $\mathbf{cf}$, $X$ is tangent to $\mathbf{cf}$.
\end{romanumerate}
\end{defn}
\begin{rk}Such a vector bundle exists by Swan's theorems.\end{rk}
We see that with these definitions, $\mathring{g}$ is a section of $\bm{w}C^{\infty}(\widebar{\mathcal M};{}^\sp T\widebar{\mathcal M})$, where $\bm{w}$ is any product of the following powers of bdfs: the sixth power of a bdf of $\mathbf{if}$, the square of a bdf of $\mathbf{bf}$, the fourth power of a bdf of $\mathbf{sf}$, and the square of a bdf of $\mathbf{cf}$. We will think of $\bm{w}$ as a weight.

Also, because the trivialization of $\widebar{\mathcal M}$ extends that of $\mathcal M$, $^\sp T\widebar{\mathcal M}$, defined in the sense of \cref{def:C3:spbundle} is precisely the pullback of $^\sp T\widebar{\mathcal M}$ to $\mathcal M$ under the embedding $\mathcal M \embeds \widebar{\mathcal M}$.

Let us now give a coordinate description of $^\sp T\widebar{\mathcal M}$. In $(\xi,v,\eta,\theta)$ coordinates, a basis of smooth sections is:
\[\xi\pa_\xi - 2\eta\pa_\eta, \pa_v, \xi\eta\pa_{\theta^1},\xi\eta\pa_{\theta^2}\] (where $\pa_{\theta^i}$ ($i = 1,2$) is a local basis of $TS^2$). In $(\varpi,v,\tau,\theta)$ coordinates:
\[\varpi\pa_{\varpi} - 2\tau\pa_\tau, \pa_v, \varpi\pa_{\theta^1},\varpi\pa_{\theta^2},\]
and in $(s,\sigma,\vepsilon,\theta)$ coordinates:
\[\pa_s, \sigma\pa_\sigma, \pa_{\theta^1}, \pa_{\theta^2}.\]

Let us also denote by $^\sp TS^2$ the subbundle of $^\sp T\widebar{\mathcal M}$ consisting of those vector which are tangent to the fibres (diffeomorphic to $S^2$), and write $\iota: {}^\sp TS^2\to {}^\sp T\widebar{\mathcal M}$ for the inclusion map.

We may of course form the tensor products $^\sp T_q^p\widebar{\mathcal M} = (^\sp T\widebar{\mathcal M})^{\otimes p} \otimes (^\sp T^\ast\widebar{\mathcal M})^{\otimes q}$.

We redefine the following for our extension $\widebar{\mathcal M}$.
\begin{defn}We call sections of ${^\sp T^p_q \widebar{\mathcal M}}$ \emph{type $(p,q)$ short-pulse tensors}.\footnote{Unfortunately this conflicts with the notion of \emph{the} short-pulse tensor $\mathbf{T}$ of the initial data in the short-pulse ansatz. Since we will always use $\mathbf{T}$ to denote this tensor, and will never use the terms in the same context, this will not be a cause for confusion.} We call type $(1,0)$ short-pulse tensors \emph{short-pulse vector fields} and type $(0,1)$ short-pulse tensors \emph{short-pulse one forms}. We similarly call sections of the tensor product ${^\sp T^p_q S^2} $ \emph{type $(p,q)$ short-pulse fibre tensors}. We call type $(1,0)$ short-pulse tensors \emph{short-pulse fibre vector fields} and type $(0,1)$ short-pulse tensors \emph{short-pulse fibre one forms}.\end{defn} 

\begin{defn}Let $g$ be a symmetric type $(0,2)$ short-pulse tensor, i.e.\ $g$ is a section of $\Sym^2({^\sp T^\ast\widebar{\mathcal M}}))$. Suppose that $g$ is non-degenerate and is of Riemannian (resp. Lorentzian) signature. We will call $g$ a \emph{short-pulse Riemannian (resp. Lorentzian) metric}. We similarly call a non-degenerate section (no signature assumptions) of $\Sym^2({{^\sp T^\ast S}})$ a \emph{short-pulse fibre metric}, or a \emph{short-pulse pseudo-Riemannian metric}. \end{defn}

The analogues of \cref{thm:C3:algebroid}, and hence also \cref{thm:C3:connectionisnice} and \cref{thm:C3:curvatureisnice}, are valid in this setting, too, namely:
\begin{lem}\label{thm:C4:geoisnice}Suppose $X_1,X_2$ are short-pulse vector fields. Then $[X_1,X_2]$, with the commutator interpreted as between sections of $T\widebar{\mathcal M}$, is in fact a short-pulse vector field. 

In particular, if $g$ is a short-pulse metric, then its Levi-Civita connection $\nabla$ extends to a map $\nabla: C^\infty(\mathcal M;\spT) \to C^\infty(\widebar{\mathcal M};{}^\sp T^1_1 \widebar{\mathcal M})$.

Furthermore, $\Riem(g)$ and $\Ric(g)$ are both short-pulse tensors.\end{lem}
\begin{proof}The proof is the same as that of its analogues in \cref{C:C3:SPCharValue}, and is omitted.\end{proof}

Let $\bm{w}$ be any product of bdfs as above. Then:
\begin{cor}\label{thm:C4:geoisniceII}Let $\bm{w}^{-1}g$ be a short-pulse metric, i.e.\ $g = \bm{w}h$, for $h$ a non-degenerate section of $\Sym^2(^\sp T^\ast \widebar{\mathcal M})$.
Then $\bm{w}^{-1}\Riem(g)$ is a section of $^\sp T^0_4 \widebar{\mathcal M}$ (i.e.\ $\Riem(g) = \bm{w}k$ for $k$ a section of $^\sp T^0_4 \widebar{\mathcal M}$) and $\Ric(g)$ is a section of $\Sym^2(^\sp T^\ast \widebar{\mathcal M})$.\end{cor}
\begin{proof}The statement about $\Ric(g)$ follows immediately from the statement about $\Riem(g)$ by contracting.

Let $h = \bm{w}^{-1}g$ be a short-pulse metric. Observe that away from $\mathbf{cf}$, we may take $\bm{w} = \delta \rho$. Using the formulae for change of curvature under a conformal change, and that $h$ is a metric on the level sets of $\delta$, we find that away from $\mathbf{cf}$ in $\widebar{\mathcal M}^\circ$, 
\[\Riem(g) = \Riem(\delta \rho h) = \delta \rho\hspace{-3pt} \left(\Riem(h)\hspace{-1.5pt} - \hspace{-1.5pt} h\owedge\left(\frac{1}{2}\nabla (d\rho/\rho) - d\rho/\rho \otimes \frac{1}{4}d\rho/\rho + \frac{1}{8}\norm{d\rho/\rho}^2_{h}h\right)\hspace{-1pt}\right),\]
where $\nabla$ is the Levi-Civita connection of $h$, and where $\owedge$ denotes the Kulkarni-Nomizu product. It follows from \cref{thm:C4:geoisnice} that $\Riem(g)$ is a section of $\bm{w}\Sym^2(^\sp T^\ast \widebar{\mathcal M})$.

Near $\mathbf{cf}$, we may take $\bm{w} = \vepsilon^4\sigma^2 = r^2$. Like above,
\[\Riem(g) = \Riem(r^2 h) = r^2\left(\Riem(h) - h\owedge\left(\nabla (dr/r) - dr/r \otimes dr/r + \frac{1}{2}\norm{dr/r}^2_{h}h\right)\right).\] As above, it follows from \cref{thm:C4:geoisnice} that $\Riem(g)$ is a section of $\bm{w}\Sym^2(^\sp T^\ast \widebar{\mathcal M})$ near $\mathbf{cf}$, too.
\end{proof}

We introduce two important vector fields.
\begin{defn}Let $N$, $\mathring{L}$ be the vector fields which in $\widebar{\mathcal M}^\circ$ defined by 
\[N = \pa_u = \pa_t - \pa_r, \ \mathring{L} = \pa_{\und{u}} = \pa_t + \pa_r,\]
where $\pa_r$ is given by $r\pa_r := y^i\pa_{y^i}$.\end{defn}
The vector fields $N$ and $\mathring{L}$ are of course the future-directed null geodesic generators of the level sets of $\und{u}$ and $u$, respectively, for the Minkowski metric $\mathring{g}$.

The one-forms $du$ and $d\und{u}$ are not short-pulse one-forms, nor are $N$ and $\mathring{L}$ short-pulse vector fields. To see the correct scaling, we introduce more weights. Let $\bm{w}_1$ be a product of the following powers of bdfs: the square of a bdf of $\mathbf{if}$, the square of a bdf of $\mathbf{sf}$ and a bdf of $\mathbf{cf}$, and let $\bm{w}_2$ be a product of: the square of a bdf of $\mathbf{bf}$, the fourth power of a bdf of $\mathbf{if}$, the square of a bdf of $\mathbf{sf}$, and a bdf of $\mathbf{cf}$. Notice that $\bm{w}_1\bm{w}_2$ is a valid choice of weight $\bm{w}$.

Then:
\begin{lem}\label{thm:C4:rescalings}The one-forms $\bm{w}_1^{-1}du$ and $\bm{w}_2^{-1}d\und{u}$ defined a priori on $\widebar{\mathcal M}^\circ$ extend to short-pulse one-forms. In fact in coordinates we have (recall that the representations are only unique up to adding on $d\delta = 0$)
\begin{align}
\begin{split}
du &= \xi^2(-2d\xi/\xi) = \varpi^2\tau^2(-2d\varpi/\varpi-2d\tau/\tau) = \vepsilon^2\sigma(1/2ds + (s-1)/2 d\sigma/\sigma)\\
d\und{u} &= \xi^4\eta^2dv = \varpi^4\tau^2dv = \vepsilon^2\sigma(1/2ds + (s+1)/2 d\sigma/\sigma),
\end{split}
\end{align}
in the three different coordinate systems, respectively.

Additionally, $\bm{w}_1N$ and $\bm{w}_2\mathring{L}$ are both short-pulse vector fields. In fact in coordinates, we have
\begin{align}
\label{eq:C4:rescalings}
\begin{split}
N &= \xi^{-2}(\eta\pa_\eta - 1/2\xi\pa_\xi) = \varpi^{-2}\tau^{_2}(1/2\varpi\pa_{\varpi} - \tau\pa_\tau) = \vepsilon^{-2}\sigma^{-1}((1-s)\pa_s + \sigma\pa_\sigma)\\
\mathring{L} &= \xi^{-4}\eta^{-2}\pa_v = \varpi^{-4}\tau^{-2}\pa_v = \vepsilon^{-2}\sigma^{-1}((1+s)\pa_s + \sigma\pa_\sigma),
\end{split}
\end{align}
in the different coordinate systems, respectively.
\end{lem}
\begin{proof}Just compute. The statement about the vector fields and one-forms are in fact equivalent if we notice that $-2\gradt_{\mathring{g}}u = \mathring{L}$ and $-2\gradt_{\mathring{g}}\und{u} = N$.\end{proof}

Now we are in a position to introduce our notion of a double-null gauge.
For the rest of this thesis, we will only work in a specific double-null gauge, so instead of saying ``$g$ is in double null gauge with...'', we will simply say ``$g$ is in a double-null gauge.''
\begin{defn}\label{def:C4:DNGI}Let $\bm{w}$ be a weight as above, and suppose $\bm{w}^{-1}g$ is a short-pulse Lorentzian metric. We say $g$ is in a \emph{double-null gauge} if over $\widebar{\mathcal M}^\circ$, and for any fixed $\delta_0$, $g|_{\widebar{\mathcal M}^\circ\n \{\delta = \delta_0\}}$ is in a double-null gauge with $u, \ \und{u},\ N$ (in the sense of \cref{def:C2:doublenullgauge}). In other words:
\begin{romanumerate}
\item there exists a nonvanishing timelike vector field, $T$, defined on all of $\widebar{\mathcal M}^\circ$, such that $u,\und{u}$ are increasing towards the future according to $T$;
\item $g^{-1}(du,du) = g^{-1}(d\und{u},d\und{u}) = 0$;
\item setting $-2\Omega^{-2} := g^{-1}(2du,2d\und{u})$, it follows that $N = -2\Omega^2\grad_g \und{u}$.
\end{romanumerate}\end{defn}

This definition looks slightly suspect, since it is a condition on the interior. Nonetheless, if one rescales, then these conditions make sense uniformly up to the boundary, as we will now indicate. Indeed:
\begin{lem}\label{thm:C4:DNGII}Let $g$ be as above, and set $h = \bm{w}^{-1}g$. Then $g$ is in a double-null gauge if and only if on all of $\widebar{\mathcal M}$, even up to the boundary, the following hold:
\begin{romanumerate}
\item there exists a non-vanishing timelike short-pulse vector field $T$ defined on all of $\widebar{\mathcal M}$ for which $\bm{w}_1^{-1}duT, \bm{w}_2^{-1}d\und{u}T > 0$;
\item 
$h(\bm{w}_1^{-1}du,\bm{w}_1^{-1}du) = h(\bm{w}_2^{-1}d\und{u},\bm{w}_2^{-1}d\und{u}) = 0$;
\item setting $-2\Omega^{-2} := h^{-1}(2\bm{w}_1^{-1}du,2\bm{w}_2^{-1}d\und{u})$, it follows that $\bm{w}_1N = -2\Omega^2h(\bm{w}_2^{-1}d\und{u},\bullet)$.
\end{romanumerate}

In this case, $\Omega^{2} \in \R$ defined in (iii) coincides with the $\Omega^2$ defined in (iii) of \cref{def:C4:DNGI}.
\end{lem}
\begin{rk}Notice that all contractions make sense since we are working with short-pulse tensors.\end{rk}
\begin{proof}
For the moment let us ignore the fact that $\Omega > 0$ may not be well-defined, although $\Omega^2 \in \R$ is. We will show that we may take $\Omega > 0$ in \cref{thm:C4:properties}, below. Condition (ii) clearly implies condition (ii) of \cref{def:C4:DNGI}. Conversely, condition (ii) of \cref{def:C4:DNGI} implies that (ii) is true over $\widebar{\mathcal M}^\circ$, and so by continuity on all of $\widebar{\mathcal M}$. Similarly, the quantity $\Omega^2$ defined by each condition (iii) coincides and the conditions are equivalent.

Condition (i) implies condition (i) of \cref{def:C4:DNGI}. For the converse, observe that \cref{thm:C2:DNGeo} it implies that $g^{-1}(du,d\und{u}) < 0$ on $\widebar{\mathcal M}^\circ$, and hence $h^{-1}(\bm{w}_1^{-1}du,\bm{w}_2^{-1}dv) \leq 0$ on all of $\widebar{\mathcal M}$. Since $\bm{w}_1^{-1}du,\bm{w}_2^{-1}dv$ are linearly-independent and null covectors, it follows that $h^{-1}(\bm{w}_1^{-1}du,\bm{w}_2^{-1}dv) \neq 0$, and hence in fact $h^{-1}(\bm{w}_1^{-1}du,\bm{w}_2^{-1}dv) < 0$. Thus
\[L' = -2h^{-1}(\bm{w}_1^{-1}du,\bullet), \ N' = -2h^{-1}(\bm{w}_2^{-1}d\und{u},\bullet)\]
are linearly independent null vector fields, and $T = L'+N'$ is a non-vanishing timelike short-pulse vector field, for which by definition $\bm{w}_1du T$, $\bm{w}_2d\und{u}T > 0$.
\end{proof}

Let us introduce some notation For $g$, $h$ as in \cref{def:C4:DNGI} and \cref{thm:C4:DNGII}, respectively, denote as usual $\slash{g}= \iota^\ast g$, $\slash{h} = \iota^\ast h$, the pullbacks of $g,h$ respectively to the fibres. Let us also define $L = -2\Omega^2\gradt_g u$, so that $\bm{w}_2L = -2\Omega^2h^{-1}(\bm{w}_1du,\bullet)$ is a short-pulse vector field.

Using \cref{thm:C4:DNGII}, we have the analogue of \cref{thm:C2:DNGeo}, whose proof is the same.

\begin{lem}\label{thm:C4:properties}The fibre metrics $\slash{g} = \bm{w}\slash{h}$ and $\slash{h}$ are short-pulse Riemannian fibre metrics. Additionally, $\Omega > 0$ defined by (iii) of either \cref{def:C4:DNGI} or \cref{thm:C4:DNGII} is well-defined. Lastly, $\vspan\{\bm{w}_2L,\bm{w}_1N\}$ is transverse to $^\sp TS^2$.\end{lem}
We will also introduce:
\begin{defn}Let $g$ be in a double-null gauge, and $L = -2\Omega^2\gradt \und{u}$ as above. Let us define a vector field $\overline{L}$, so that $\bm{w}_2\ovelrine{L} = \bm{w}_2(L-\mathring{L})$ is a short-pulse vector field.\end{defn}
The importance of $\overline{L}$ is that $\overline{L}u = \overline{L}\und{u} = \overline{L}\delta = 0$, and so $\bm{w}_2\overline{L}$ is a short-pulse fibre vector field.

We have the analogy of \cref{thm:C2:Reconstruct} and \cref{thm:C3:ReconstructSP}, whose proof is the same:
\begin{lem}\label{thm:C4:Reconstruct}
Let $\Omega$ be a section of the trivial line bundle, $\bm{w}_2\overline{L}$ a short-pulse fibre vector field, and $\bm{w}^{-1}\slash{g}$ a short-pulse fibre Riemannian metric. Then there exists a unique metric $g$ for which $\bm{w}^{-1}g$ is a short-pulse Lorentzian metric in a double null gauge such that $\overline{L} = -2\Omega^2\grad_g\und{u}-\mathring{L}$, and $\slash{g}$ is the restriction of $g$ to fibre-tangent vectors.\end{lem}

We end this section with a simple, but useful observation showing that the short-pulse double-null gauges of \cref{C:C3:SPCharValue} and the double-null gauges introduced this section are essentially the same object when considered over $\mathcal M \subseteq \widebar{\mathcal M}$. Observe that in $\mathcal M$, $x$ is a bdf of $\mathbf{bf}$ (since $x = \eta \rho$ and $\rho > 0$), and so $\bm{w}_1 = 1$ and $\bm{w}_2 = \bm{w} = x^2$ are valid choices of the weights.
\begin{lem}
If $\bm{w}^{-1}g$ is a section of $\Sym^2(^\sp T^\ast \widebar{\mathcal M})$ of Lorentzian signature (i.e.\ is a short-pulse Lorentzian metric) and is in a double-null gauge (in the sense of this section), then over $\mathcal M$, $x^{-2}g$ is a section of $\Sym^2(^\sp T^\ast \mathcal M)$ of Lorentzian signature which is in a short-pulse double-null gauge with $(u,v,N)$. Conversely, if $h \in \Sym^2(^\sp T^\ast \widebar{\mathcal M})$ is a short-pulse metric which over $\mathcal M$ is in in a short-pulse double-null gauge with $(u,v,N)$ (in the sense of \cref{def:C3:DNG}) then $\bm{w}h$ is in a double-null gauge.\end{lem}

\section{Continuation as a formal series}
\label{C:C4:existence}
\subsection{Setup and statement of the Theorem}
In this section we state the precise version of the existence/uniqueness and degeneracy portion of \cref{thm:C1:bigthm} and prove most of it. Let us fix short-pulse data, as described in \cref{C:C3:ansatz}. In particular, we fix the short-pulse tensor $\mathbf{T}$.
First, we define $\digamma$.
Let us work in coordinates $(\xi,v,\eta,\theta)$, which when restricted to $\{\xi = 0\}$ gives us coordinates $(v,\eta,\theta)$ on $\mathbf{if}\setminus \mathbf{sf}$.
\begin{defn}
Define $\digamma:\mathbf{if}\setminus \mathbf{sf} \to \R$ by
\[\digamma(v,\eta,\theta) = 4-2\eta^2\int_0^v \mathbf{E}(t,\theta)\ dt,\]
where
\[\mathbf{E}(t,\theta) = \frac{1}{2}\int_0^v |\pa_t \mathbf{T}(t,\theta)|^2\ dt\]
is the energy of $\mathbf{T}$.\end{defn}

We will often use theorems from \cref{C:A3}. In order to do this, we observe:
\begin{lem}\label{thm:C4:rect}For every point $(v,\eta,\theta) \in \mathbf{if}\n\{\digamma > 0\}$, there is a small neighbourhood $U \subseteq S^2$ of $\theta$ for which $[0,v]\times[0,\eta]\times U \subseteq \{\digamma > 0\}$.\end{lem}
\begin{proof}Notice that $\digamma$ is continuous in $\theta$, and decreasing in both $v$,$\eta$. In particular on a small neighbourhood $U$ of $\theta$, $\digamma(v,\eta,\theta') > 0$ for $\theta' \in U$, and thus $\digamma(v',\eta',\theta') > 0$ for $(v',\eta',\theta') \in [0,v]\times[0,\eta]\times U$.\end{proof}
The advantage is that the results of \cref{C:A3}, which are only stated for rectangles, also apply to $\mathbf{if}\n\{\digamma > 0\}$.

Before we can state the main theorem, we will need to slightly extend the definition of polyhomogeneity to include \emph{log smooth} functions. The reason for this will be primarily notational rather than technical. Roughly speaking a log smooth function $u$ is one which possesses an asymptotic expansion
\[u \sim \sum_{(k,p) \in \N^2} u_{k,p}t^k{\log}^p t,\]
where $t$ is a bdf of some face (see \cref{def:C4:logsmooth}, below for a precise definition). We also require that for each $k$, there are only finitely many $p$ such that $u_{k,p} \neq 0$. We will also sometimes assume that $u_{0,p} = 0$ for $p \geq 1$, i.e.\ there are no ${\log}$s to top order. One may not strictly realize $\log$-smooth functions as polyhomogeneous function. Indeed the sets \begin{align*}
E_{{\log},0} &:= \{(0,0)\} \un (\N^{\geq 1}\times \N^{\geq 0})\\
E_{\log} &:= \N^{\geq 0}\times \N^{\geq 0}
\end{align*} are not a valid index sets because the sets $\{(z,p)\: \Re(z) < N\}$ are infinite for $N \geq 1$. The advantage of log-smoothness will be that it allows us to work with polyhomogeneity without having to keep track of exactly which index set we use, i.e.\ without keeping track of how many ${\log}$s we allow in our expansion. This will be convenient when find a solution to the (nonlinear!) Einstein equations, since the nonlinear interactions mean it becomes cumbersome to keep track of precisely how many ${\log}$s are present in our series. See \cref{thm:C4:algebra} and \cref{rk:C4:algebra} for a more detailed explanation. The difference between $E_{\log,0}$ and $E_{\log}$ is that functions with expansions according to $E_{\log,0}$ have well-defined restrictions to a corresponding face.

Formally, we (re)define:
\begin{defn}\label{def:C4:logsmooth}Let $X$ be a mwc, and let $\mathcal E = (E_1,\ldots, E_n)$ be an index family of ``generalized'' index sets, i.e.\ each $E_i$ is either an index set or is $E_{\log,0}$ or $E_{\log}$. Then we define $\mathcal A^{\mathcal E}_{\mathrm{phg}}(X)$ to be the union
\[\mathcal A^{\mathcal E}_{\mathrm{phg}}(X) := \bigcup_{\mathcal E'} \mathcal A^{\mathcal E'}_{\mathrm{phg}}(X),\]
where the union runs over all (proper) index families $\mathcal E' = (E_1',\ldots, E_n')$ where $E_i' \subseteq E_i$ is a (proper) index set.
If $\mathcal F_i$ is a face corresponding to an index set $E_i = E_{\log,0}$ or $E_i = E_{\log}$, then we say $u \in \mathcal A^{\mathcal E}_{\mathrm{phg}}(X)$ is ``${\log}$-smooth'' at $F_i$.\end{defn}
The same definition extends to the partially polyhomogeneous spaces and for polyhomogeneous sections of vector bundles.

As mentioned in the introduction, we do not expect our solution to continue past $\{\digamma = 0\}$. One checks that the flow of the restriction of $\bm{w}_1 N$ to $\mathbf{if}$ flows from $\mathbf{bf}$ to $\mathbf{sf}$, and likewise the flow on $\mathbf{sf}$ flows from $\mathbf{if}$ to $\mathbf{cf}$. Since $\bm{w}_1 N$ is future-directed for $g$, we expect that the metric does not last long enough to have any sort of expansion on $\mathbf{sf} \un \mathbf{cf}$. Thus, for the rest of this thesis, we will not have cause to deal with all of $\widebar{\mathcal M}$. In fact, in this chapter we will be able to work in coordinates $(\xi,v,\eta,\theta)$, since these cover most of the region of interest. Let $\widebar{\mathcal M}'$ denote the region covered by these coordinates which lies inside $\{\digamma > 0\}$, which is diffeomorphic to
\[([0,1]_\xi \times [0,1]_v \times [0,\infty)_\eta \times S^2_\theta)\n \{\digamma > 0\},\] and write $^{\sp}T\widebar{\mathcal M'}$ for the short-pulse bundle over $\widebar{\mathcal M'}$.\footnote{We have implicitly chosen an extension of $\digamma$ from $\{\xi = 0\}$ to $\{\xi > 0\}$ by requiring it to be constant in $\xi$. Since the extension theorem is a statement only about solutions at the boundary, any other extension works equally well.} In $\widebar{\mathcal M'}$, we fix the values
\begin{romanumerate}
\item $\bm{w} = \xi^6\eta^2$;
\item $\bm{w}_1 = \xi^2$;
\item $\bm{w}_2 = \eta^2\xi^4$.
\end{romanumerate}

\begin{figure}[htbp]
\centering
\includegraphics{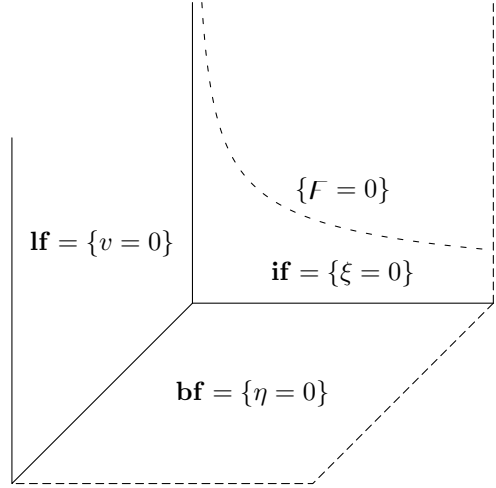}
\caption{A stylized view of $\protect\widebar{\mathcal M'}$ at a certain angular slice $\theta$. $\mathbf{rf} = \{\xi = 1\}$ is not shown, and is the ``front'' face meeting $\mathbf{lf}$ and $\mathbf{bf}$ transversely.}
\end{figure}

We record a special case of \eqref{eq:C4:rescalings}:
\begin{align*}
\bm{w}_1N &= \xi^2N = \eta\pa_\eta - \frac{1}{2}\xi\pa_\xi\\
\bm{w}_2\mathring{L} &= \xi^4\eta^2\mathring{L} = \pa_v.
\end{align*}

\begin{rk}\label{rk:C4:LieDer}The advantage of fixing $\bm{w}_1$ and $\bm{w}_2$ in this way means that $\Lie_{\bm{w}_1 N}$ and $\Lie_{\bm{w}_2 \mathring{L}}$ are well-defined differential operators on fibre-tensors (cf.\ \cref{thm:C3:LieDer}). Indeed, denoting $\pi:\widebar{\mathcal M} \to \widebar{\mathcal N}$ the projection off the fibres, then $\bm{w}_1 N$ and $\bm{w}_2 L$ are both $\pi$-related to vector fields on the base. The proof of \cref{thm:C3:LieDer} now goes through verbatim.\end{rk}

For the rest of this section, let $\mathcal E =(0,0,0,E_{\log,0})$, and $\mathcal E_{\emptyset} = (\emptyset,\ldots,\emptyset)$, corresponding to the ordering $(\mathbf{rf},\mathbf{lf},\mathbf{rf},\mathbf{if})$ of the boundary faces of $\Mp$.

For a metric $g$ in a double-null gauge, we let $\slash{g} = \iota^\ast g$ denotes its restriction to fibre-tangent vectors. As usual, we continue to denote by $\mathring{\slash{g}}$ the metric on the round sphere.

The precise version of the existence/uniqueness and degeneracy statement of \cref{thm:C1:bigthm} is:
\begin{thm}\label{thm:C4:continuation}There exists a section $g \in \bm{w}\mathcal A_{\mathrm{phg}}^\mathcal E(\widebar{\mathcal M'};\Sym^2(^\sp T^\ast\widebar{\mathcal M'}))$ of Lorentzian signature that is in a double-null gauge and which is a solution to $\Ric(g) = 0$ in series at $\mathbf{rf}$, $\mathbf{lf},\mathbf{bf},\mathbf{if}$, i.e.\ \[\Ric(g) \in \mathcal A_{\mathrm{phg}}^{\mathcal E_{\emptyset}}(\widebar{\mathcal M'};\Sym^2(^\sp T^\ast\widebar{\mathcal M'})),\] such that $g$ has short-pulse data. Any two such solutions are the same up to an element of $\bm{w}\mathcal A_{\mathrm{phg}}^{\mathcal E_{\emptyset}}(\widebar{\mathcal M'}; \Sym^2(^\sp T^\ast \widebar{\mathcal M'}))$.

Moreover, 
\[\det (\xi^{-4}\mathring{\slash{g}}^{\vspace{-5pt}-1}\slash{g})|_{\mathbf{if}\n \{\digamma > 0\}} = \frac{1}{16}\digamma^2.\] In particular, $\xi^{-4}\slash{g}$ degenerates as $\digamma \to 0$.
\end{thm}
In more pedestrian terms, \cref{thm:C4:continuation} states that there exists a unique series solution $g$ in each of the faces with the powers of the bdfs in the series expansions given by $\mathcal E$, and $\slash{g}$ cannot continue as a fibre metric past any extension of $\{\digamma = 0\}$ into $\Mp$.
\begin{rk}Notice that that the theorem says nothing about what happens at $\mathbf{ff}$. This is because $\mathbf{ff}$ is an ``artificial face,'' and the solution will always continue past it (for instance by extending the data slightly as we did in \cref{C:C3:FormTrapped}).\end{rk}

We focus on proving the existence portion of \cref{thm:C4:continuation}, the degeneracy of $\slash{g}$ coming as a side effect of the proof. In fact, we will be able to obtain reasonably precise information about components of the metric in a double-null gauge at $\mathbf{if}$. This examination has a slightly different flavour to the rest of this section, and consequently we reserve \cref{C:C4:topo} to do it.

Observe that $\overline{\mathcal M'}\n\{\xi > 0\} = \mathcal M$. Since
\[\bm{w}\mathcal A^{(0,0,0)}_{\mathrm{phg},(\mathbf{rf},\mathbf{lf},\mathbf{bf})}(\mathcal M;\Sym^2(^\sp T^\ast \widebar{\mathcal M'})) = x^2C^{\infty}(\mathcal M; \Sym^2(^\sp T^\ast \mathcal M)),\]
\cref{thm:C3:SPTSI} proves most of \cref{thm:C4:continuation}, in a sense which will be made precise in the next subsection. Thus, we need only worry about what happens near $\mathbf{if}$.\footnote{To be pedantic, \cref{thm:C3:SPTSI} only allows us to consider the region $\mathcal M \n \{v < 1\}$. But we can always arbitrarily extend the data as in \cref{C:C3:FormTrapped} to get up to $\{v = 1\}$. We will often implicitly use this extension for the rest of this section.}

An important feature of polyhomogeneous functions which are smooth or $\log$-smooth at at faces is that they are closed under composition with smooth functions.
\begin{lem}\label{thm:C4:algebra}
Let $X$ be a manifold with corners and faces $\mathcal F = (F_1,\ldots,F_n)$, and $\mathcal E' = (E_1,\ldots,E_n)$, where each $E_i$ is one of $0$, or $E_{\log,0}$.

If $U \subseteq \R$ is open and $f\:U \to \R$ is smooth, and $a \in \mathcal A^{\mathcal E'}_{\mathrm{phg}}(X)$ takes values in $U$, then $f(a) \in \mathcal A^{\mathcal E'}_{\mathrm{phg}}(X)$.

More generally, if $V$, and $W$ are vector bundles, and $U \subseteq V$ is open, and $f:U \to W$ is smooth, and if $a\in \mathcal A^{\mathcal E'}_{\mathrm{phg}}(X;V)$ takes values in $U$, then $f(a) \in \mathcal A^{\mathcal E'}_{\mathrm{phg}}(X;W)$.

In the special case that $f\: V^\ast\times V \to \R$ is tensor contraction, then the conclusion also holds if we allow some of the $E_i$ to be $E_{\log}$. In particular $\mathcal \mathcal A^{\mathcal E'}_{\mathrm{phg}}(X)$ is an algebra, and the tensor contraction map is well-defined $\mathcal A^{\mathcal E'}_{\mathrm{phg}}(X;V^\ast \otimes V) \to \mathcal A^{\mathcal E'}_{\mathrm{phg}}(X)$.

\end{lem}
\begin{rk}\label{rk:C4:algebra}
Notice that this lemma is not true if we replace $E_{\log,0}$ or $E_{\log}$ with any proper subindex set other than $0$. This is the main point of allowing $E_{\log,0}$ and $E_{\log}$ as ``generalized'' index sets: we allow $\log$s, but our class of polyhomogeneous functions are still closed under nonlinear operations.
\end{rk}
\begin{rk}This lemma is ``obvious'' if $\mathcal E' = (0,\ldots,0)$, since in this case $\phg{\mathcal E'}(X) = C^{\infty}(X)$. We present a proof to convince the reader that the argument may be generalized to the more unfamiliar log-smooth setting.\end{rk}
\begin{proof}
It suffices to treat the case that $X = [0,1)^k\times [0,1)^{n-k}\times \R^m$, $0 \leq k \leq n \in \N$, $m \in \N$, and $V$ and $W$ are the trivial bundles $X \times \R^M$, $X \times \R^K$, respectively. We will use coordinates $(t_1,\ldots,t_n,s)$ for $X$.

Let us first handle the case of general smooth maps $f$.

It suffices to prove the following stronger claim: let $J \subseteq \{1,\ldots,n\}$ be a subset, and let $\mathcal E_J$ be an associated index family, each of whose components is $0$ or $E_{\log,0}$. We will prove the lemma not only for the space of polyhomogeneous functions, but also the space of partially polyhomogeneous functions $\mathcal A_{\mathrm{phg},J}^{\mathcal E_J}(X)$, by induction on $\#J$.

If $\#J=0$, then $f(a) \in \mathcal A(X)$ if $a \in \mathcal A(X)$ by the chain rule.

Now assume the claim for $\#J = j-1$. We prove it for $\#J = j$. For simplicity, we will only treat the case $M=K =1$. Fix some $i \in J$. 
Let $\hat{\mathcal E}_J$ denote the index family obtained from $\mathcal E_J$ by removing the index set corresponding to $E_i$, and let $\hat{J} = J\setminus \{i\}$. We need to show that $f(a)$ has a polyhomogeneous expansion at $\{t_i = 0\}$.

Consider first the case $E_i = 0$. Fix $a \in \phgi{J}{\mathcal E_J}(X)$, and suppose $a$ takes values in $U$. Write for each $N$
\[a = \sum_{j=0}^{N-1} a_jt_i^j + t_i^N \mathcal A_{\mathrm{phg},\hat{J}}^{\hat{\mathcal E}_J}(X),\]
where $a_j \in \phgi{\hat{J}}{\hat{\mathcal E}_J}(X)$.
Write $r = (t_1,\ldots,\hat{t_i},\ldots, t_n)$ for a coordinate on $\{t_i = 0\}$.\footnote{The notation $\hat{t_i}$ means we are removing this coordinate.} We may Taylor expand $f$ around any $y \in U$, writing
\begin{equation}\label{eq:algebra:eq1}f(x) = \sum_{j=0}^{N-1} (x-y)^jf_j(y) + (x-y)^N R_N(x,y)),\end{equation} where $j!f_j = (\pa_x^j f)(y) \in C^{\infty}(U)$ and $R_N(x,y) \in C^{\infty}(U \times U)$. Decompose
\[a = a_0 + b_N + t_i^NS_N,\] where
\[b_N = \sum_{j=1}^{N-1}a_jt_i^j,\]
$a_0,\ldots,a_{N-1} \in \phgi{\hat{J}}{\hat{\mathcal E}_J}(X)$ and $S_N \in \mathcal A_{\mathrm{phg},\hat{J}}^{\hat{\mathcal E}_J}(X)$.\footnote{This is the point of the proof where we require $E_i = 0$ or $E_i = E_{\log,0}$: the restriction of $a$ to $\{t_i = 0\}$, $a_0$, must be well-defined!} Choosing $y = a_0$ in \eqref{eq:algebra:eq1}, it follows that for all $N'$
\begin{align*}
f(a) = \sum_{j=0}^{N-1} (b_{N'} + t_i^{N'}S_{N'})^jf_j(a_0) + (b_{N'} + t_i^{N'}S_{N'})^N R_N(a,a_0).
\end{align*}
By induction, $f_j(a_0) \in \mathcal A_{\mathrm{phg},\hat{J}}^{\hat{\mathcal E}_J}(\{t_i = 0\})$, 
$R_N(a,a_0) \in A_{\mathrm{phg},\hat{J}}^{\hat{\mathcal E}_J}(X)$, and, provided $N' \geq N$, for each $j \geq 1$
\[(b_N'+ t_i^{N'}S_{N'})^j = \sum_{\ell = j}^{N-1} c_{N',j}t_i^\ell + t_i^NA_{\mathrm{phg},\hat{J}}^{\hat{\mathcal E}_J}(X),\] where
$c_{N',j} \in A_{\mathrm{phg},\hat{J}}^{\hat{\mathcal E}_J}(\{t_i = 0\})$ (we have used the inductive hypothesis applied to the multiplication map to deal with the products). Using the inductive hypothesis again to handle the products, it follows that
\[f(a) = \sum_{j=0}^{N-1}t_i^{j}F_j + t_i^NA_{\mathrm{phg},\hat{J}}^{\hat{\mathcal E}_J}(X),\] where $F_j \in A_{\mathrm{phg},\hat{J}}^{\hat{\mathcal E}_J}(\{t_i = 0\})$. This shows that $f(a)$ has the desired polyhomogeneous expansion.

The same proof works if $E_i = E_{\log,0}$.

Let us now treat the case that $f$ is tensor contraction. For simplicity, we only treat the case that $V = X \times \R$ is the trivial bundle. As above, we prove the claim more generally for $\mathcal A_{\mathrm{phg},J}^{\mathcal E_J}(X)$, by induction on $\#J$. The case $\#J = 0$ is clear since $\mathcal A(X)$ is an algebra. The inductive step is also clear: if
\begin{align*}
a &= \sum_{j < N,p\in \N}^{N-1} a_{j,p}t_i^j\log^p t_i + t_i^N \mathcal A_{\mathrm{phg},\hat{J}}^{\hat{\mathcal E}_J}(X)\\
b &= \sum_{j < N, p\in \N}^{N-1} b_{j,p}t_i^j\log^p t_i + t_i^N \mathcal A_{\mathrm{phg},\hat{J}}^{\hat{\mathcal E}_J}(X),
\end{align*}
where 
$a_j, b_j \in \phgi{\hat{J}}{\hat{\mathcal E}_J}(X)$, then $ab$ has a similar expansion by induction.
\end{proof}

We obtain the important corollary, which extends \cref{thm:C4:geoisnice} and \cref{thm:C4:geoisniceII} to the log-smooth setting.
\begin{cor}\label{thm:C4:ricisnice}Suppose $g \in \bm{w}\phgd(\Mp;\Sym^2(^\sp T^\ast \Mp)$ is a short-pulse metric. Then its connection is a map
\[\nabla \: \phgd(\Mp;{}^\sp T\Mp) \to \phgd(\Mp;{}^\sp T_1^1\Mp).\]
Moreover, $\Riem(g) \in \bm{w}^{-1}\phgd(\Mp;{}^\sp T^0_4(\Mp))$ and $\Ric(g) \in \phgd(\Mp;\Sym^2(^\sp T^\ast \Mp))$.\end{cor}
\begin{proof}[Proof sketch.]
Let $h = \bm{w}^{-1}g \in \phgd(\Mp;\Sym^2(^\sp T^\ast \Mp)$.

We first show that if $X$ is a short-pulse vector field, then \[\nabla X \in \phgd(\Mp;{}^\sp T_1^1 \Mp).\] This is essentially \cref{thm:C4:geoisnice}, except one now also needs to use \cref{thm:C4:algebra} to deal with the contractions appearing in the Koszul formula, together with one additional fact: if $Y$ is a short-pulse vector field, then $Y$ tangent to $\mathbf{if}$, and so $\phgd(\Mp;\R)$ is stable under derivation via $Y$.\footnote{Both these facts were implicitly used for the $C^{\infty}$ category during the proof of \cref{thm:C4:geoisnice}.}
Using that $\phgd(\Mp;{}^\sp T\Mp) = \phgd(\Mp;\R)\otimes C^{\infty}(\Mp; {}^\sp T\Mp)$, \cref{thm:C4:algebra} implies the desired mapping property of $\nabla$.

Using \cref{thm:C4:algebra} to analyze tensor contradictions again, it is clear that \[\Riem(h) \in \phgd(\Mp;{}^\sp T^0_4(\Mp)), \ \Ric(h) \in \phgd(\Mp;\Sym^2(^\sp T^\ast(\Mp))).\]

To complete the proof of the corollary, use the previous results, together with \cref{thm:C4:algebra} and the formulae for curvature under a conformal change given in the proof of \cref{thm:C4:geoisniceII}.
\end{proof}
We can now start discussing the proof of \cref{thm:C4:continuation} in detail. As in the proofs of \cref{thm:C2:TSI} and \cref{thm:C3:SPTSI}, we will first reduce the full Einstein equations in a double-null gauge to three equations for $(\slash{g} = \iota^\ast g,\overline{L} = L-\mathring{L},\omega = \log \Omega)$ which we will solve.

\subsection{Reduction to simpler equations}
In this subsection, we reduce \cref{thm:C4:continuation} to a couple statements about solutions to tractable equations for components of $g$ in a double-null gauge.

Our first step is to deal with the apparent overdetermined nature of the Einstein equations.

In a double-null gauge, $\Ric(g)$ splits up naturally into several tensors:
\begin{romanumerate}
\item $2\Omega^2\Ric_{-,-}$, a section of $\Sym(^\sp T^\ast S^2)$;
\item $4\Omega^2\Ric_{\bm{w}_1N,-}$, a section of $^\sp T^\ast S^2$;
\item $-4\Omega^2\Ric_{\bm{w}_2L,-}$, a section of $^\sp T^\ast S^2$;
\item $-2\Ric_{\bm{w}_1N,\bm{w}_1N}$, a section of the trivial line bundle;
\item $-2\Ric_{\bm{w}_2L,\bm{w}_2L}$, a section of the trivial line bundle;
\item $-\frac{1}{2}\Ric_{\bm{w}_1N,\bm{w}_2L}$, a section of the trivial line bundle.
\end{romanumerate}
Here a $-$ indicates an argument in $^\sp TS^2$.

As usual, when solving $\Ric(g) = 0$, looking at some of these tensors is redundant, namely:
\begin{prop}\label{thm:C4:integrability}Suppose $g \in \bm{w}\mathcal A_{\mathrm{phg}}^\mathcal E(\widebar{\mathcal M'};\Sym^2(^\sp T^\ast\widebar{\mathcal M'}))$ satisfies that $2\Omega^2\Ric_{-,-}$, $-4\Omega^2\Ric_{\bm{w}_2L,-}$ and $-\frac{1}{2}\Ric_{\bm{w}_1N,\bm{w}_2L}$ are all rapidly vanishing at $\mathbf{lf}$, $\mathbf{rf}$, $\mathbf{bf}$, $\mathbf{if}$, i.e.\ suppose
\begin{align*}
2\Omega^2\Ric_{-,-} &\in \mathcal A_{\mathrm{phg}}^{\mathcal E_\emptyset}(\widebar{\mathcal M'};\Sym^2(^\sp T^\ast S^2))\\
-4\Omega^2\Ric_{\bm{w}_2L,-} &\in \mathcal A_{\mathrm{phg}}^{\mathcal E_\emptyset}(\widebar{\mathcal M'};{}^\sp T^\ast S^2)\\
-\frac{1}{2}\Ric_{\bm{w}_1N,\bm{w}_2L} &\in \mathcal A_{\mathrm{phg}}^{\mathcal E_\emptyset}(\widebar{\mathcal M'};\R),
\end{align*}
and that
$4\Omega^2\Ric_{\bm{w}_1N,-}$, $-2\Ric_{\bm{w}_1N,\bm{w}_1N}$ are identically $0$ on $\mathbf{lf}$ and $-2\Ric_{\bm{w}_2L,\bm{w}_2L}$ is identically $0$ on $\mathbf{rf}$. Then all components of $\Ric(g)$ are rapidly vanishing, i.e.\ also
\begin{align}
\label{eq:integrability:4eq10}
\begin{split}
 -4\Omega^2\Ric_{\bm{w}_1N,-} &\in \mathcal A_{\mathrm{phg}}^{\mathcal E_\emptyset}(\widebar{\mathcal M'};{}^\sp T^\ast S^2)\\
 -2\Ric_{\bm{w}_1N,\bm{w}_1N} &\in \mathcal A_{\mathrm{phg}}^{\mathcal E_\emptyset}(\widebar{\mathcal M'};\R)\\
 -2\Ric_{\bm{w}_2L,\bm{w}_2L} &\in \mathcal A_{\mathrm{phg}}^{\mathcal E_\emptyset}(\widebar{\mathcal M'};\R),
\end{split}
\end{align}
so that
\[\Ric(g) \in \mathcal A_{\mathrm{phg}}^{\mathcal E_{\emptyset}}(\widebar{\mathcal M'};\Sym^2(^\sp T^\ast\widebar{\mathcal M'})).\]
\end{prop}
\begin{proof}
By \cref{thm:C3:integrability}, we have that all three quantities in \eqref{eq:integrability:4eq10} are rapidly vanishing at all faces in $\xi > 0$. However, all quantities are also $\log$-smooth, i.e.\ are in $\phgd(\Mp;\bullet)$ (where $\bullet$ is either $^\sp T^\ast S^2$ or $\R$ where appropriate). The only way this is possible is if they are in fact rapidly vanishing all these faces all the way to $\xi = 0$,\footnote{\label{fn:C4:ifn}Technically, this requires the use of \cref{thm:A4:betterphg} to pass from rapidly vanishing at a collection of faces to having the appropriate index set there, too. The main difficulty is that although the respective expansions at $\mathbf{bf}$ and $\mathbf{lf}$ are by assumption trivial, it is \emph{not} immediate that the expansions at $\mathbf{bf}$ and $\mathbf{lf}$ of \emph{the coefficients in the log-smooth expansion} at $\mathbf{if}$ are also trivial. Since this and similar results are intuitive and we make frequent use of them in this chapter, we only mention this technicality once.} i.e.\ are in \[\phg{(\emptyset,\emptyset,\emptyset,E_{\log,0})}(\Mp;\bullet).\] 

It thus suffices to show that each quantity is rapidly vanishing at $\mathbf{if} = \{\xi = 0\}$, i.e.\ is in
\[\phg{(0,0,0,\emptyset)}(\Mp;\bullet),\]
which we will do by showing that all coefficients in their series expansions at $\mathbf{if}$ are $0$.

If $g \in \bm{w}\mathcal A_{\mathrm{phg}}^\mathcal E(\widebar{\mathcal M'};\Sym^2(^\sp T^\ast\widebar{\mathcal M'}))$ is in a double-null gauge, then by \cref{thm:C4:algebra}, the associated $0 < \Omega^2 \in \mathcal A_{\mathrm{phg}}^\mathcal E(\widebar{\mathcal M'};\R)$, and so, again by \cref{thm:C4:algebra}, $\Omega^{\pm 1} \in \mathcal A_{\mathrm{phg}}^\mathcal E(\widebar{\mathcal M'};\R).$

Observe that $\Ric_{\bm{w}_1N,\xi\eta (-)}$ is in $\phgd(\Mp;T^\ast S^2)$, not $\phgd(\Mp;{}^\sp T^\ast S^2)$, and 
\[\Ric_{\bm{w}_1N,\xi\eta (-)} \in \phg{\mathcal E_{\emptyset}}(\Mp;T^\ast S^2)\text{ if and only if } \Ric_{\bm{w}_1N,-} \in \phg{\mathcal E_{\emptyset}}(\Mp;{}^\sp T^\ast S^2).\]

Using our assumptions, we record schematically the contracted Bianchi identities from \cref{C:A1:computebianchi} (or see them written in a slightly different form in \eqref{eq:C2:onei}--\eqref{eq:C2:threei}). Under the assumptions and also using \cref{thm:C4:algebra}, we may write these equations as
\begin{equation}\label{eq:integrability:4eq1}\Lie_{\bm{w}_2 L}\Ric_{\bm{w}_1N,\xi\eta(-)} + (\bm{w}_2H^L)\Ric_{\bm{w}_1N,\xi\eta(-)} \in \mathcal A_{\mathrm{phg}}^{\mathcal E_\emptyset}(\widebar{\mathcal M'};{}^\sp T^\ast S^2)\end{equation}
(commuting $\bm{w}_1$ through the equation by using the fact that since $\bm{w}_2L \bm{w}_1 = 0$).

If we can prove $\Ric_{\bm{w}_1N,\xi\eta(-)} \in \mathcal A_{\mathrm{phg}}^{\mathcal E_\emptyset}(\widebar{\mathcal M'},^\sp T^\ast S^2)$, then we may further write (commuting $\bm{w}_2$ using that $\bm{w}_1N \bm{w}_2 = 0$)
\begin{align}
\label{eq:integrability:4eq2}
\bm{w}_2L\Ric_{\bm{w}_1N\bm{w}_1N} + (\bm{w}_2H^L)\Ric_{\bm{w}_1N\bm{w}_1N} &\in \mathcal A_{\mathrm{phg}}^{\mathcal E_\emptyset}(\widebar{\mathcal M'};\R)\\
\label{eq:integrability:4eq3}
\bm{w}_1N\Ric_{\bm{w}_2L\bm{w}_2L} + (\bm{w}_1H^N)\Ric_{\bm{w}_2L\bm{w}_2L} &\in \mathcal A_{\mathrm{phg}}^{\mathcal E_\emptyset}(\widebar{\mathcal M'};\R).
\end{align}

Let us look at \eqref{eq:integrability:4eq1}. Set $\slash{h} = \bm{w}^{-1}\slash{g}$, which is a short-pulse fibre Riemannian metric.
Then
\begin{align*}
 \bm{w}_2H^L &= \bm{w}_2\Tr(\slash{g}^{-1}\Lie_L\slash{g})\\
 &= \bm{w}_2\Tr(\bm{w}^{-1}\slash{h}^{-1}\Lie_L \bm{w}\slash{h})\\
 &=\Tr(\slash{h}^{-1}\bm{w}^{-1}\Lie_{\bm{w}_2 L}\bm{w}\slash{h})\\
 &= \Tr(\slash{h}^{-1}\Lie_{\bm{w}_2L}\slash{h}),
\end{align*}
since $\bm{w}_2L\bm{w} = 0$. We could bring the $\bm{w}_2$ into the Lie-derivative since $\slash{h}$ is a fibre tensor, and $\bm{w}_2$ is constant on the fibres. Thus, by \cref{thm:C4:algebra}, $\bm{w}_2H^L \in \phgd(\Mp;\R)$.

In particular, we may expand
\[\bm{w}_2H^L = \sum_{k < M, p \in \N} H_{k,p}\xi^k\log^p\xi + \xi^{M} \phgi{(\mathbf{rf},\mathbf{lf},\mathbf{bf})}{(0,0,0)}(\Mp;\R),\] where each $H_{k,p}$ is smooth and does not depend on $\xi$. Also, for each $k$ there are finitely many $p$ for which $H_{k,p}$ is nonzero, and $H_{0,p} = 0$ for $p > 0$. Likewise, $\Ric_{\bm{w}_1N,\xi\eta(-)}$ admits the expansion
\[\Ric_{\bm{w}_1N,\xi\eta(-)} = \sum_{k < M, p \in \N} R_{k,p}\xi^k\log^p\xi + \xi^{M} \phgi{(\mathbf{rf},\mathbf{lf},\mathbf{bf})}{(0,0,0)}(\Mp;{}^\sp T^\ast S^2).\]

Since $\bm{w}_2L$ commutes with $\xi$, \eqref{eq:integrability:4eq1} becomes
\[\sum_{k < M, p \in \N} \left(\Lie_{\bm{w}_2L}R_{k,p} + \sum_{i+j = k, q+r = p} H_{j,r}R_{i,q}\right)\xi^k\log^p \xi \in \xi^{M} \phgi{(\mathbf{rf},\mathbf{lf},\mathbf{bf})}{(0,0,0)}(\Mp;T^\ast S^2).\]
This implies that for each $k,p$,
\[\Lie_{\bm{w}_2L}R_{k,p} + \sum_{i+j = k, q+r = p} H_{j,r}R_{i,q} = 0.\]
Since $\bm{w}_2L = \pa_v$, and at $\mathbf{lf}$ the assumptions mean that $R_{k,p}|_{\{v = 0\}} \equiv 0$, we can show inductively that each $R_{k,p} = 0$, and thus $\Ric_{\bm{w}_1N,\xi\eta(-)} \in \mathcal A_{\mathrm{phg}}^{(0,0,0,\emptyset)}(\widebar{\mathcal M'},T^\ast S^2)$.

Returning to \eqref{eq:integrability:4eq2}, the exact same argument now works to establish
$\Ric_{\bm{w}_1N,\bm{w}_1N} \in \phg{(0,0,0,\emptyset)}(\Mp;\R)$.

Dealing with \eqref{eq:integrability:4eq3} is slightly trickier because $\bm{w}_1N$ differentiates in the $\xi$ direction and is singular in the $\eta$ direction. First, notice that $\bm{w}_1H^N \in \phgd(\Mp;\R)$, for the same reason that $\bm{w}_2H^L$ is. We may therefore expand
\begin{align*}
\bm{w}_1H^N &= \sum_{k < M, p \in \N} H_{k,p}\xi^k\log^p\xi + \xi^{M} \phgi{(\mathbf{rf},\mathbf{lf},\mathbf{bf})}{(0,0,0)}(\Mp;\R)\\
\Ric_{\bm{w}_2L,\bm{w}_2L} &=\sum_{k < M, p \in \N} R_{k,p}\xi^k\log^p\xi + \xi^{M}\phgi{(\mathbf{rf},\mathbf{lf},\mathbf{bf})}{(0,0,0)}(\Mp;\R),
\end{align*}
where $R_{k,p}$, $H_{k,p}$ are smooth and do not depend on $\xi$, and $R_{0,p}, H_{0,p}$ are zero if $p \geq 1$.

Arguing similarly to above, but now observing that if $\bm{w}_1N = \eta\pa_\eta - 1/2\xi\pa_\xi$ hits a factor $x^k\log^p \xi$ it generates another term, from \eqref{eq:integrability:4eq3} we arrive at
\begin{equation}\label{eq:integrability:4eq5}\Lie_{\bm{w}_1 N} R_{k,p} - \frac{(p+1)}{2}R_{k,p+1} + \sum_{i+j = k, q+r = p} H_{j,r}R_{i,q} = 0.\end{equation}
Since $\bm{w}_1N = \eta\pa_\eta$ on $\{\xi = 0\}$, this is a linear ODE with singular derivative $\eta\pa_\eta$. We already know that $R_{k,p} \in \eta^\infty C^{\infty}(\{\xi = 0\})$. Indeed, this follows from the fact that
\[\Ric_{\bm{w}_2L,\bm{w}_2L} \in \phg{(\emptyset,\emptyset,\emptyset,E_{\log,0})}(\Mp;\R).\]

Now we may apply a theorem on ODEs with singular derivatives from \cref{C:A3:Sing}. Using \cref{thm:A3:GoursatSing} (or rather the special case where there is no $v$ dependence), we may show by induction that there is a unique solution to \eqref{eq:integrability:4eq5} which is rapidly vanishing, and thus must coincide with $R_{k,p}$. However, it is certainly that $R_{k,p} = 0$ (for all $k,p$) is also a solution in $\eta^{\infty}C^{\infty}(\{\xi = 0\})$, and thus $R_{k,p} = 0$ for all $k,p$. So $\Ric_{\bm{w}_2L,\bm{w}_2L} \in \phg{(0,0,0,\emptyset)}(\Mp;\R)$.\end{proof}
Let us now fix $g \in \bm{w}\mathcal A_{\mathrm{phg}}^\mathcal E(\widebar{\mathcal M'};\Sym^2(^\sp T^\ast\widebar{\mathcal M'}))$ in a double-null gauge, and set
$V = \xi^3\eta\overline{L}$ (recall $\overline{L} = L-\mathring{L}$) a section of $TS^2$ over $\Mp$ (rescaled so it is \emph{not} a section of $^\sp TS^2$) and $\slash{k} = \xi^{-4}\slash{g}$, a section of $\Sym^2(T^\ast S^2)$ over $\Mp$ (rescaled so it is \emph{not} a section of $\Sym^2(^\sp T^\ast S^2)$). Let $\bm{w}_3$ be a new weight, a product of a bdf of $\mathbf{if}$ and $\mathbf{bf}$, which in $\Mp$ we fix to be $\bm{w}_3 = \xi\eta$. In light of \cref{thm:C4:integrability}, we only need to solve:
\begin{align}
\label{eq:C4:addon}
\begin{split}2\Omega^2\Ric_{\bm{w}_3(-),\bm{w}_3(-)} &\in \mathcal A_{\mathrm{phg}}^{\mathcal E_\emptyset}(\widebar{\mathcal M'};\Sym^2( T^\ast S^2)),\\
-4\Omega^2\Ric_{\bm{w}_2L,\bm{w}_3(-)} &\in \mathcal A_{\mathrm{phg}}^{\mathcal E_\emptyset}(\widebar{\mathcal M'}; T^\ast S^2),\\
-\frac{1}{2}\Ric_{\bm{w}_1N,\bm{w}_2L} &\in \mathcal A_{\mathrm{phg}}^{\mathcal E_\emptyset}(\widebar{\mathcal M'};\R)\end{split}\end{align}
(notice $-$ is now an argument in $TS^2$, and we have removed the $^\sp$ from the scaling of the sections).

Write $\phi = (\slash{k},V,\omega= \log \Omega)$, a section of $X :=\Sym^2(T^\ast S^2)\oplus TS^2 \oplus \R$ over $\Mp$. Keep the notation from \cref{C:C2:RicciDNG}, except all slash and other tensor operations are taken with respect to $\slash{k}$. Then we may express the Ricci components in \eqref{eq:C4:addon} in terms of $\phi$. Let us make the shorthand
\[Z = Z_{\phi} = \slash{k}([\bm{w}_1N,V]+1/2V,\bullet).\]
We record, by rescaling \eqref{eq:C2:one}, \eqref{eq:C2:three}, \eqref{eq:C2:six},
\begin{subequations}
\begin{align}
\label{eq:C4:one}
\begin{split}
2\Omega^2\Ric_{\bm{w}_3(-),\bm{w}_3(-)} &= \Lie_{\bm{w}_1 N}\Lie_{\bm{w}_2 \mathring{L}} \slash{k}+ \bm{w}_3\Lie_{\bm{w}_1 N}\Lie_V \slash{k} - \frac{\bm{w}_3}{2}\Lie_{[\bm{w}_1 N,V]}\slash{k} + \frac{\bm{w}_3}{4}\Lie_V\slash{k}\\
&+\left(\frac{1}{4}\slash{\tr}(\Lie_{\bm{w}_1 N}\slash{k})-1\right)\left(\Lie_{\bm{w}_2\mathring{L}+\bm{w}_3 V}\slash{k}\right)\\
&+\frac{1}{4}\slash{\tr}\Lie_{\bm{w}_2\mathring{L}+\bm{w}_3 V}\slash{k}\left(\Lie_{\bm{w}_1 N}\slash{k} - 2\slash{k}\right)\\
&-\frac{1}{2}\left[\left(\Lie_{\bm{w}_2\mathring{L}+\bm{w}_3 V}\slash{k}\right)\times \Lie_{\bm{w}_1 N}\slash{k} + \Lie_{\bm{w}_1 N}\slash{k}\times \left(\Lie_{\bm{w}_2\mathring{L}+\bm{w}_3 V}\slash{k}\right)\right]\\
&+2\bm{w}_3^2\Omega^2\slash{\Ric}(\slash{k}) - 2\bm{w}_3^2\slash{\Hess}\Omega^2 + 2\bm{w}_3^2\Omega^2\slash{d}\omega\otimes\slash{d}\omega\\
&- \frac{1}{4\Omega^2}Z\otimes Z.
\end{split}\\
\label{eq:C4:two}
\begin{split}
-4\Omega^2\Ric_{\bm{w}_2 L,\bm{w}_3(-)} &= \Lie_{\bm{w}_2 \mathring{L}+\bm{w}_3 V}Z -2\Omega^2\bm{w}_3\slash{\div}\left(\Lie_{\bm{w}_2\mathring{L}+\bm{w}_3 V}\slash{k}\right)\\
&+2\Omega^2\bm{w}_3\slash{d}\slash{\tr}\Lie_{\bm{w}_2\mathring{L}+\bm{w}_3 V}\slash{k}+\frac{1}{2}\left(\slash{\tr}\Lie_{\bm{w}_2\mathring{L}+\bm{w}_3 V}\slash{k}\right)Z\\
&-\bm{w}_3\left(\slash{\tr}\Lie_{\bm{w}_2\mathring{L}+\bm{w}_3 V}\slash{k}\right)\slash{d}\Omega^2\\
&+4\Omega^2\bm{w}_3 \slash{d}(\bm{w}_2\mathring{L}+\bm{w}_3 V)\omega - 2((\bm{w}_2\mathring{L} + \bm{w}_3 V)\omega) Z
\end{split}\\
\label{eq:C4:three}
\begin{split}
-\frac{1}{2}\Ric_{\bm{w}_2 L,\bm{w}_1 N} &= (\bm{w}_1 N)(\bm{w}_2\mathring{L} + \bm{w}_3 V)\omega + \frac{1}{4}(\bm{w}_1 N)\slash{\tr}\Lie_{\bm{w}_2\mathring{L}+\bm{w}_3 V}\slash{k}\\
&- \frac{\bm{w}_3}{8}\slash{\tr}\Lie_{[\bm{w}_1 N,V]+\frac{1}{2}V}\slash{k}-\frac{\bm{w}_3}{2}[\bm{w}_1 N,V]\omega -\frac{\bm{w}_3}{4}V\omega - \frac{\bm{w}_3^2}{2}\slash{\Delta}\Omega^2\\
&+ \frac{1}{8\Omega^2}\left(|[\bm{w}_1 N,V]|^2 +\frac{1}{4}|V|^2 + \slash{k}([\bm{w}_1N,V],V)\right)\\
&+\frac{1}{8}\slash{k}(\Lie_{\bm{w}_2\mathring{L}+\bm{w}_3 V}\slash{k},\Lie_{\bm{w}_1 N}\slash{k}) - \frac{1}{4}\slash{\tr}\Lie_{\bm{w}_2\mathring{L}+\bm{w}_3 V}\slash{k}.
\end{split}
\end{align}
\end{subequations}

Let us denote by $P = (P^1,P^2,P^3)$ the non-linear differential operator, taking sections of $X$ (with non-degenerate $\slash{k})$ to section of $Y = \Sym^2(T^\ast S^2)\oplus T^\ast S^2\oplus \R$, whose action on a section $\phi = (\slash{h},V,\omega)$ is given setting $P^1\phi$, $P^2\phi$, $P^3\phi$ to be the right-hand side of \eqref{eq:C4:one}, \eqref{eq:C4:two} and \eqref{eq:C4:three}, respectively. It is important to observe:
\begin{lem}If $\phi \in \phgd(\Mp;X)$, then $P\phi \in \phgd(\Mp;Y)$.\end{lem}
\begin{proof}[Proof sketch.]Use \cref{thm:C4:algebra}, together with the fact that $\phgd(\Mp;\R)$ is stable under derivation by smooth short-pulse vector fields (cf.\ \cref{thm:C4:ricisnice}).\end{proof}

From \cref{thm:C4:Reconstruct}, one may from any $\phi$ recover the metric $g$ which induces the components $\slash{h}$, $V$, $\omega$. It is also not hard to see, going through the proof of \cref{thm:C4:Reconstruct} carefully, that $\phi \in \phgd(\Mp,X)$ if and only if $g \in \bm{w}\mathcal A_{\mathrm{phg}}^\mathcal E(\widebar{\mathcal M'};\Sym^2(^\sp T^\ast\widebar{\mathcal M'}))$. In light of \cref{thm:C4:integrability}, provided $g$ has the short-pulse data, and the data for $\phi$ are given accordingly and solve $\Ric_{NN}, \ \Ric_{N,-} = 0$ on $\mathbf{lf}$ and $\Ric_{LL} = 0$ on $\mathbf{rf}$, then
$P\phi \in \phg{\mathcal E_{\emptyset}}(\Mp;Y)$ if and only if $\Ric(g) \in \phg{\mathcal E_{\emptyset}}(\Mp;\Sym^2(^\sp T^\ast \Mp))$. That $\phi$ may be arranged to have these Ricci components $0$ is due to \cref{thm:C3:toporder}, where it was shown how to construct the initial data set for $\phi$ (really $\phi$ rescaled slightly differently), combined with \cref{thm:C3:isreg}, where it was shown that the constraints were satisfied (and hence the data exist on all of $\mathbf{rf}$ and $\mathbf{lf}$).\footnote{We have glossed over one subtle one small issue: we have $\Ric_{NN}$, $\Ric_{N,-} = 0$ only on $\mathbf{lf} \n \{\xi > 0\}$ from \cref{thm:C3:isreg}. However, once we know that $\phi \in \phgd(\Mp;X)$, we know that once rescaled appropriately, they must continue as polyhomogeneous sections of rescaled bundles to $\{\xi = 0\}$, and thus the rescaled versions are $0$ up to $\{\xi = 0\}$ and we may actually apply \cref{thm:C4:integrability}.}

In other words, it suffices to prove:
\begin{thm}
\label{thm:C4:reduction}
There exists a solution $\phi \in \phgd(\Mp;X)$ with
\[P\phi \in \phg{\mathcal E_{\emptyset}}(\Mp,Y),\]
and the associated metric $g$ in a double-null gauge has short-pulse data and solves \[\Ric_{\bm{w}_1N,\bm{w}_1N} = 0, \ \Ric_{\bm{w}_1N,-} = 0 \text{ on }\mathbf{lf}\] and \[\Ric_{\bm{w}_2L,\bm{w}_2L} = 0\text{ on }\mathbf{rf}.\] Any two solutions are the same up to an element of $\phg{\mathcal E_{\emptyset}}(\Mp;X)$.

Moreover, \[\det \mathring{\slash{g}}^{\vspace{-5pt}-1}\slash{k} = \frac{1}{16}\digamma^2.\]
\end{thm}

Thus, for the rest of the proof of \cref{thm:C4:continuation}, and indeed the rest of this section, we will focus mainly on $\phi$ instead of $g$.

By \cref{thm:C3:SPTSI}, there is already a unique $\phi$ in series at all faces other than $\mathbf{if}$. Thus to complete the log-smoothness/polyhomogeneity at $\mathbf{if}$, it suffices to show two things. First, we show that the coefficients in the Taylor series expansions of $\phi$ at $\mathbf{bf}$ and $\mathbf{lf}$ are themselves log-smooth at $\mathbf{bf}\n\mathbf{if}$ and $\mathbf{lf}\n \mathbf{if}$, respectively. Let $\phi \in C^{\infty}(\Mp\n \{\xi > 0\};X)$ be the unique series solution (modulo rapidly vanishing functions) to $P\phi = 0$ with short-pulse data which is provided by \cref{thm:C3:SPTSI}, i.e.\[P\phi \in x^\infty v^\infty u^\infty C^{\infty}(\Mp \n \{\xi > 0\};Y)\] (recall that $\Mp \n \{\xi > 0\} = \mathcal M$). Write
\begin{align*}
\phi &\sim \sum v^k \phi^{\mathbf{lf}}_k\end{align*}
and
\begin{align*}
\phi &\sim \sum x^k \tilde{\phi}^{\mathbf{bf}}_k = \sum \eta^k\phi^{\mathbf{bf}}_k
\end{align*}
for the Taylor-series expansions of $\phi$ at $\mathbf{lf}\n \{ \xi > 0\}$ and $\mathbf{rf}\n \{\xi > 0\}$, respectively, where the coefficients $\phi^{\mathbf{lf}}_k \in C^{\infty}(\mathbf{lf}\n \{\xi > 0\})$ and $\phi^{\mathbf{bf}}_k \in C^{\infty}(\mathbf{bf}\n \{\xi > 0\})$ are unique.
\begin{prop}
\label{thm:C4:first}
For all $k \in \N$ in fact \[\phi_k^{\mathbf{lf}} \in \phgi{(\mathbf{rf},\mathbf{bf},\mathbf{if})}{(0,0,E_{\log,0})}(\mathbf{lf};X), \phi_k^{\mathbf{bf}} \in \phgi{(\mathbf{rf},\mathbf{lf},\mathbf{if})}{(0,0,E_{\log,0})}(\mathbf{bf};X).\footnote{Really we mean $\phi_{k}^{\mathbf{lf}}$ and $\phi_k^{\mathbf{bf}}$ \emph{extend} uniquely to log-smooth functions.}\]
The top-order behaviour at $\mathbf{lf}$ and $\mathbf{bf}$, respectively, is given by
\[\phi_0^{\mathbf{lf}} = (\mathring{\slash{g}},0,0) + \xi \phgi{(\mathbf{rf},\mathbf{bf},\mathbf{if})}{(0,0,E_{\log})}(\mathbf{lf};X),\ \phi_0^{\mathbf{bf}} = (\mathring{\slash{g}},0,0) + \xi \phgi{(\mathbf{rf},\mathbf{lf},\mathbf{if})}{(0,0,E_{\log})}(\mathbf{bf};X).\]
Furthermore, writing $\slash{k}_i^{\mathbf{bf}}$ for first component of $\phi_i^{\mathbf{bf}}$, and $V_i^{\mathbf{bf}}$ for the second,
\begin{align*}\slash{k}_1^{\mathbf{bf}} &= \mathring{\slash{g}}\mathbf{T},\\
\slash{k}_2^{\mathbf{bf}} &= \mathring{\slash{g}}S_0 + \xi \phgi{(\mathbf{rf},\mathbf{lf},\mathbf{if})}{(0,0,E_{\log})}(\mathbf{bf};\Sym^2(T^\ast S^2)),\\
V_1^{\mathbf{bf}} &= 0 + \xi \phgi{(\mathbf{rf},\mathbf{lf},\mathbf{if})}{(0,0,E_{\log})}(\mathbf{bf};TS^2)\end{align*}
where 
\[S_0 = -\frac{1}{2}\int_0^v \mathbf{E}(s)\ ds + \frac{1}{4}\Tr(\mathbf{T}^2).\]
\end{prop}

Second, we show:
\begin{prop}
\label{thm:C4:second}
There are unique coefficients $\phi^{\mathbf{if}}_{k,p} \in \phg{(0,0)}(\mathbf{if}\n\{\digamma > 0\})$ (subject to the condition $\phi^{\mathbf{if}}_{0,p} = 0$ for $p \geq 1$ and $\phi^{\mathbf{if}}_{k,p} = 0$ for each $k$ fixed and $p$ large enough), smooth on $\mathbf{if} \n \{\digamma > 0\}$ and satisfying the compatibility conditions at $\mathbf{bf}\n \mathbf{if}$ and $\mathbf{lf}\n \mathbf{lf}$, such that the asymptotic sum\footnote{The asymptotic sum is not unique, so we really mean \emph{any} asymptotic sum.}
\[\phi \sim \sum_{k \geq 0,p \geq 0}\xi^k\log^p \xi\phi^{\mathbf{if}}_{k,p}\]
solves $P\phi \in \phgi{(\mathbf{lf},\mathbf{bf},\mathbf{if})}{(0,0,\emptyset)}(\Mp;X)$.

Moreover, denoting by $\slash{k}^{\mathbf{if}}_{0,0}$ the first component of $\phi^{\mathbf{if}}_{k,p}$,
\[\det \mathring{\slash{g}}^{\vspace{-5pt}-1}\slash{k}^{\mathbf{if}}_{0,0} = \frac{1}{16}\digamma^2.\]
\end{prop}

Recall that the compatibility conditions here mean, writing
\[\phi^{\mathbf{if}}_{k,p} \sim \sum \eta^j\phi^{\mathbf{if}}_{k,p,j}\]
and
\[\phi^{\mathbf{bf}}_j \sim \sum \xi^k\log^p \xi \phi^{\mathbf{bf}}_{j,k,p}\]
(where $\phi^{\mathbf{bf}}_j$ is given by \cref{thm:C4:first}), that
\[\phi^{\mathbf{if}}_{k,p,j} = \phi^{\mathbf{bf}}_{j,k,p},\]
and similarly for the expansion at $\mathbf{lf}$.

Let us show how these two propositions imply \cref{thm:C4:reduction}. \begin{proof}[Proof of \cref{thm:C4:reduction}]
Let us first show existence.
We know that the coefficients $\phi^{\mathbf{lf}}_j$ and $\phi^{\mathbf{bf}}_k$ satisfy the compatibility conditions $\pa_\eta^k \phi^{\mathbf{lf}}_j = \pa_v^j \phi^{\mathbf{bf}}_k$ in $\xi > 0$ because they are the expansions of an actual function in $C^{\infty}(\Mp\n\{\xi > 0\};X)$. However, since both extend to polyhomogeneous functions towards $\xi = 0$, they continue to satisfy the compatibility conditions for all $\xi$. \Cref{thm:C4:second} also provides us with $\phi^{\mathbf{if}}_{z,p}$ satisfying the compatibility conditions. Combining this with the coefficients $\phi_j^{\mathbf{rf}}$ provided by the series solution at $\mathbf{rf}$ (which exists by \cref{thm:C3:SPTSI}), Borel's lemma now provides for the existence of $\phi \in \phgd(\Mp;Y)$ such that $P\phi$ is vanishing in series at $\mathbf{rf}$, $\mathbf{if} \n \{\digamma > 0\}$ and also for all $\xi_0 > 0$, is vanishing in series at $\mathbf{lf}\n \{\xi \geq \xi_0\}$ and $\mathbf{bf}\n \{\xi \geq \xi_0\}$.

However, since $P\phi \in \phgd(\Mp;Y)$, the last two vanishing properties hold only if $P\phi$ is vanishing rapidly at $\mathbf{lf}$ and $\mathbf{bf}$, too. 

Thus, $P\phi$ is vanishing at all faces, i.e.\ $P\phi \in \phg{\mathcal E_{\emptyset}}(\Mp;Y)$.

By construction of the coefficients $\phi^{\mathbf{lf}}_k$ and $\phi^{\mathbf{rf}}_k$, \[\Ric_{\bm{w}_1N,\bm{w}_1N} = 0, \ \Ric_{\bm{w}_1N,-} = 0 \text{ on }\mathbf{lf}, \ \ \Ric_{\bm{w}_2L,\bm{w}_2L} = 0\text{ on }\mathbf{rf}\]
and has short-pulse data. This completes the proof of existence.

For uniqueness, any solution $\phi$ comes with its own $\phi^{\mathbf{lf}}_k$, etc., and the propositions show that these are unique.

The statement about $\det \mathring{\slash{g}}^{\vspace{-5pt}-1}\slash{k}$ follows immediately from \cref{thm:C4:second}.
\end{proof}

We will devote the next subsection to proving \cref{thm:C4:first}, and the subsequent subsection to showing \cref{thm:C4:second}, given the existence/uniqueness of the top order. We will save the discussion of the top order for the last section.

\subsection{log-smoothness of expansions}
We prove \cref{thm:C4:first}. 
The proposition is really two statements: one for the coefficients $\phi_k^{\mathbf{lf}}$, and one for the coefficients $\phi_k^{\mathbf{bf}}$. We treat the second one first.

\begin{proof}[Proof of \cref{thm:C4:first} for the case of the expansions at $\mathbf{bf}$]
We first compute the top order expressions explicitly, and then we will linearize the equations to find transport equations for the higher order $\phi_j^{\mathbf{bf}}$ which we can solve. Write $\phi_j^{\mathbf{bf}} = \phi_j = (\slash{k}_j,V_j,\omega_j)$. From the proof of \cref{thm:C3:computetrapped}, we know the values of $\phi_0$, $\phi_1$, and $\slash{k}_2$ (although we need to rescale). They are:
\begin{align*}
\phi_0 &= (\mathring{\slash{g}},0,0)\\
\phi_1 &= (\mathring{\slash{g}} \mathbf{T},0,0)\\
\slash{k}_2 &= \mathring{\slash{g}} S,\end{align*}
where
\[S = 2\xi^2 v - \frac{1}{2}\int_0^v \mathbf{E}(s)\ ds + \frac{1}{4}\Tr(\mathbf{T}^2).\]

In particular, the proposition is true for $j = 0,1$, and is true of $\slash{k}_2$. For $j \geq 0$ set
\[\phi_{(j)} = \phi_0 + \cdots + \eta^j\phi_j.\]
We will prove the proposition by induction on $j$.

Assume the inductive hypothesis for $j-1$. We prove it for $j$. We may write over $\xi > 0$
\begin{equation}\label{eq:second:eq:10}P(\phi_{(j-1)}+\eta^{j}\phi_j) = P(\phi_{(j-1)}) +\eta^j T_{\phi_{(j-1)},j}\phi_j + \eta^{j+1}C^{\infty}(\{\xi > 0\}),\end{equation}
where $T_{\phi_{(j-1)},j}$ is a linear operator which depends on $\phi_{(j-1),j}$. However the dependence on $\phi_{(j-1)}$ is smooth, and so \begin{equation}\label{eq:second:eq:11}T_{\phi_{(j-1)},j}\phi_j =T_{\phi_{0},j}\phi_j + \eta C^{\infty}(\{\xi > 0\}).\end{equation}

By the inductive hypothesis, \[\phi_{(j-1)} \in \phgd(\Mp).\] Thus,
\[P\phi_{(j-1)} \in \phgd(\Mp),\] too.
However by definition of the $\phi_i$, 
\[P(\phi_{(j-1)}) \in \eta^jC^{\infty}(\{\xi > 0\}),\] and so in fact
\begin{equation}\label{eq:first:eq5}P\phi_{(j-1)} \in \eta^{j}\phgd(\Mp).\end{equation}

Also, by definition
\begin{equation}\label{eq:second:eq13}P(\phi_{(j-1)}+\eta^{j}\phi_j) \in \eta^{j+1}C^{\infty}(\{\xi > 0\}).\end{equation} Using \eqref{eq:second:eq:10}--\eqref{eq:second:eq13}, it follows that
\[T_{\phi_0,j}\phi_{j} = -(\eta^{-j}P(\phi_{(j-1)})|_{\mathbf{bf}} \in \phgi{(\mathbf{rf},\mathbf{lf},\mathbf{if})}{(0,0,E_{\log,0})}(\mathbf{bf}),\]
at least over $\mathbf{bf} \n \{\xi > 0\}$.
We may record the components of $T_{\phi_0,j} = (T^1,T^2,T^3)$, which are obtained by performing a computation and discarding terms of order $\eta^{j+1}$ and higher, acting on $\phi_j$:
\begin{subequations}
\begin{align}
\label{eq:secondsolve:eqa}
T^1(\phi_j) &= \Lie_{\bm{w}_1 N}\Lie_{\bm{w}_2\mathring{L}}\Lie_{\bm{w}_2\mathring{L}} \slash{k}_j + (j-1)\Lie_{\bm{w}_2\mathring{L}}\slash{k}_j -\frac{1}{2}\left(\slash{\tr}\Lie_{\bm{w}_2\mathring{L}}\slash{k}_j\right)\mathring{\slash{g}}\\
\label{eq:secondsolve:eqb}
T^2(\phi_j) &= \mathring{\slash{g}}\left(\Lie_{\bm{w}_2\mathring{L}}\left((\Lie_{\bm{w}_1 N}V_j) + \left(j+\frac{1}{2}\right)V_j\right),\bullet\right)\\
\label{eq:secondsolve:eqc}
T^3(\phi_j) &= \bm{w}_1 N \bm{w}_2\mathring{L} \omega_j + j\bm{w}_2\mathring{L}\omega +\frac{1}{4}\bm{w}_1 N\slash{\tr}\Lie_{\bm{w}_2\mathring{L}} \slash{k}_j + \frac{j-1}{4}\slash{\tr}\Lie_{\bm{w}_2\mathring{L}}\slash{k}_j.
\end{align}
\end{subequations}
Here, all tensorial operations are taken with respect to $\mathring{\slash{g}}$. 

In the case $j=2$, we already know $\slash{k}_2 \in \phgi{(\mathbf{rf},\mathbf{lf},\mathbf{if})}{(0,0,E_{\log,0})}(\mathbf{bf})$. If $j \geq 3$, then we may take $\mathring{\slash{\tr}}$ of \eqref{eq:secondsolve:eqa} to obtain
\[\Lie_{\bm{w_1}N}(\Lie_{\bm{w}_2\mathring{L}}\slash{k}_j) + (j-2)\slash{\tr}(\Lie_{\bm{w}_2\mathring{L}}\slash{k}_j) \in \phgi{(\mathbf{rf},\mathbf{lf},\mathbf{if})}{(0,0,E_{\log,0})}(\mathbf{bf}).\] Now use \cref{thm:C4:secondsolve}, below, to conclude
\[\Lie_{\bm{w}_2\mathring{L}}\slash{k}_j \in \phgi{(\mathbf{rf},\mathbf{lf},\mathbf{if})}{(0,0,E_{\log,0})}(\mathbf{bf}).\] Plugging this back into \eqref{eq:secondsolve:eqa} means that
\[\Lie_{\bm{w}_1 N}\Lie_{\bm{w}_2\mathring{L}}\slash{k}_j + (j-1)\Lie_{\bm{w}_2\mathring{L}}\slash{k}_j \in\phgi{(\mathbf{rf},\mathbf{lf},\mathbf{if})}{(0,0,E_{\log,0})}(\mathbf{bf}),\] and thus again by \cref{thm:C4:secondsolve}, $\Lie_{\bm{w}_2\mathring{L}}\slash{k}_j \in \phgi{(\mathbf{rf},\mathbf{lf},\mathbf{if})}{(0,0,E_{\log,0})}(\mathbf{bf})$ and thus $\slash{k}_j \in \phgi{(\mathbf{rf},\mathbf{lf},\mathbf{if})}{(0,0,E_{\log,0})}(\mathbf{bf})$ by integrating and using the fact that $\slash{k}_j|_{\{v = 0\}}$ is $0$ for $j \geq 1$ since the data is Minkowskian. The same argument now applies to \eqref{eq:secondsolve:eqc}, yielding $\omega_j \in \phgi{(\mathbf{rf},\mathbf{lf},\mathbf{if})}{(0,0,E_{\log,0})}(\mathbf{bf})$.

For \eqref{eq:secondsolve:eqb}, lower it first via $\mathring{\slash{g}}$ to obtain a transport equation for \[\Lie_{\bm{w}_1 N} V_j + (j+1/2)V_j.\] Since $\bm{w}_1 N$ is tangent to $\mathbf{lf}$, where the data for $V_j$ is $0$, it follows by integrating that
\[\Lie_{\bm{w}_1 N} V_j + (j+1/2)V_j \in \phgi{(\mathbf{rf},\mathbf{lf},\mathbf{if})}{(0,0,E_{\log,0})}(\mathbf{bf}).\] Now apply \cref{thm:C4:secondsolve} to to conclude that
and $V_j \in \phgi{(\mathbf{rf},\mathbf{lf},\mathbf{if})}{(0,0,E_{\log,0})}(\mathbf{bf})$.
\end{proof}

\begin{lem}\label{thm:C4:secondsolve}Fix $0 < k \in \N$. Consider the following ODE for a section $f$ of $T_q^p S^2$, $(p,q) \in \N\times \N$, over $\mathbf{bf}$:
\begin{equation}\label{eq:secondsolve:eq1}\Lie_{\bm{w}_1N}f + \frac{k}{2}f = R \in \phgi{(\mathbf{rf},\mathbf{lf},\mathbf{if})}{(0,0,E_{\log,0})}(\mathbf{bf};T_q^pS^2),\end{equation}
with smooth data for $f$ at $\{\xi = 1\}$. Then \[f \in \phgi{(\mathbf{rf},\mathbf{lf},\mathbf{if})}{(0,0,E_{\log,0})}(\mathbf{bf};T_q^pS^2).\]
\end{lem}

\begin{proof}
On $\mathbf{bf}$, $\bm{w}_1 N = -\frac{1}{2}\xi\pa_\xi$. Thus one may check that the solution $f$ is given by
\begin{equation}\label{eq:secondsolve:form}f(\xi) = 2\xi^k\int_\xi^1 R(t)\ t^{-k-1}\ dt + \xi^kf(1) =: 2S(\xi) + \xi^k f(1),\end{equation} where we have suppressed the dependence on $v,\theta$.
By ODE theory, $f$ is certainly smooth away from $\{\xi = 0\}$, so we just need to check polyhomogeneity near $\{\xi = 0\}$. The second term is harmless, so we focus on $S$.

We start by showing that $f$ has a $\log$-smooth expansion at $\xi=0$. We first show that for each $N$ large enough $S$ has an expansion of the form
\begin{equation}\label{eq:secondsolve:form2}S(\xi) = \xi^k S_{N} + \sum_{(n,p), k < n < N}\xi^n\log^p \xi S_{n,p} + \xi^N H_N,\end{equation}
where $S_{N}$, $S_{n,p}$ are smooth and do not depend on $\xi$,  and $H_N \in \phgi{(\mathbf{lf})}{(0)}(\mathbf{bf}\n\{\xi \leq 1/2\})$. Notice that this is not immediately the same as $\log$-smoothness because we are allowing the top order coefficient $S_{k,N}$ to depend on $N$. However, comparing the coefficient of $\xi^k$ for various $N$, it follows that they must be the same. So let us focus on proving \eqref{eq:secondsolve:form2}.

Expand for all $N$ sufficiently large
\[R(\xi) = \sum_{(n,p), n < N} \xi^n\log^p \xi R_{n,p} + \xi^N G_N,\] where $R_{n,p}$ are smooth and do not depend on $\xi$ and $G_N \in \phgi{(\mathbf{lf})}{(0)}(\mathbf{bf}\n\{\xi \leq 1/2\})$. Let us see the contribution of the terms $\xi^m\log^p \xi R_{n,p}$ to \eqref{eq:secondsolve:form2}:
\begin{align*}
\xi^k\int_\xi^1 t^m\log^p t R_{n,p}\ t^{-k-1}\ dt = c_{m,p,k}\xi^kR_{n,p} + \xi^m q_{m,p}(\log \xi)R_{n,p},
\end{align*}
for some polynomial $q_{m,p}$ of degree at most $p+1$ and constant $c_{m,p,k}$. For $m > 0$, these terms belong in the expansion \eqref{eq:secondsolve:form2}. Since $k > 0$, the terms coming from $m=0$ have no $\log$s in them, so they also belong in the expansion.

Thus, we need only show that if $N$ is large that
\begin{equation}\label{eq:secondsolve:expr}\xi^k\int_\xi^1 t^N G_N t^{-k-1}\ dt \in \xi^kC^{\infty}(\{\xi = 1\}) + \xi^N\phgi{(\mathbf{lf})}{(0)}(\mathbf{bf}\n\{\xi \leq 1/2\}),\end{equation}
where we have extended an element of $C^{\infty}(\{\xi = 1\})$ to the entire space by requiring it to be constant in $\xi$.

Expand \[G_N = \sum_{0 \leq m < M} v^mG_{N,m} + v^M H_{N,M},\]
where $G_{N,m} \in \mathcal A(\mathbf{bf}\n \mathbf{lf} \n \{\xi \leq 1/2\})$ and $H_{N,M} \in \mathcal A(\mathbf{bf}\n \{\xi \leq 1/2\})$. It suffices to show that
\begin{align}
\label{eq:secondsolve:expr1}\xi^k\int_\xi^1 t^N G_{N,m} t^{-k-1}\ dt \in \xi^kC^{\infty}(\mathbf{bf}\n \mathbf{lf} \n \{\xi = 1\}) + \xi^N\mathcal A(\mathbf{bf}\n\mathbf{lf}\n\{\xi \leq 1/2\})\\
\label{eq:secondsolve:expr2}\xi^k\int_\xi^1 t^N H_{M,N} t^{-k-1}\ dt \in \xi^kC^{\infty}(\mathbf{bf} \n \{\xi = 1\}) + \xi^N\mathcal A(\mathbf{bf}\n\{\xi \leq 1/2\}).
\end{align}
The same proof essentially works to prove both, so we focus on proving the second.

If $N-k-1 \geq 0$, then the left-hand side of \eqref{eq:secondsolve:expr2} is equal to
\[\xi^k\int_0^1 t^N H_{M,N} t^{-k-1}\ dt - \xi^k\int_0^\xi t^N H_{M,N} t^{-k-1}\ dt.\] The first term is in $\xi^k L^{\infty}(\mathbf{bf} \n \{\xi = 1\})$, while the second is in $\xi^N L^{\infty}(\mathbf{bf}\n \{\xi \leq 1/2\})$.

$v\pa_v$ and $\pa_\theta$ commute through the first integral, which shows that the first term is in $\xi^k C^{\infty}(\mathbf{bf} \n \{\xi = 1\})$, as desired.

For the second term, its form is stable under differentiation via $v\pa_v$ and $\pa_\theta$, so it suffices to show its $\xi\pa_\xi$ derivatives are all bounded. Differentiating once $\xi\pa_\xi$ gives an expression of the form
\[\xi^k\int_0^\xi t^N A t^{-k-1}\ dt  + \xi^N B,\] where $A,B \in \mathcal A(\mathbf{bf}\n\{\xi \leq 1/2\})$. The form of this expression is stable under further differentiation via $\xi\pa_\xi$, and using the previous argument is certainly in 
\[\xi^k L^{\infty}(\mathbf{bf} \n \{\xi = 1\}) + \xi^N L^{\infty}(\mathbf{bf}\n \{\xi \leq 1/2\}).\] This shows \eqref{eq:secondsolve:expr2}.

We are not quite done, since we are also required that the coefficients in the expansion of $f$ at $\mathbf{lf}$  are \emph{themselves} $\log$-smooth at $\{\xi = 0\}$. This is a similar argument. Expand
\[R(\xi) = \sum_{0 \leq n < N} v^n R_n + v^N G_N,\]
where $R_n \in \phgi{(\mathbf{if})}{E_{\log,0}}(\mathbf{bf}\n\mathbf{lf}\n\{\xi \leq 1/2\})$ and $G_N \in \phgi{(\mathbf{if})}{E_{\log,0}}(\mathbf{bf}\n\{\xi \leq 1/2\})$.
Thus
\[f(\xi) = \sum_{0 \leq n < N} 2v^n \xi^k \int_{\xi}^1 R_n(t) t^{-k-1}\ dt + 2v^N\xi^k \int_{\xi}^1 G_N(t) t^{-k-1}\ dt + \xi^k f(1).\]

Arguing similar to above, we deduce that
\begin{align*}
\xi^k \int_{\xi}^1 R_n(t) t^{-k-1}\ dt \in \phgi{(\mathbf{if})}{E_{\log,0}}(\mathbf{bf}\n\mathbf{lf}\n\{\xi \leq 1/2\})\\
\xi^k \int_{\xi}^1 G_N(t) t^{-k-1}\ dt \in \phgi{(\mathbf{if})}{E_{\log,0}}(\mathbf{bf}\n\{\xi \leq 1/2\}),
\end{align*}
which shows that $f$ has the desired expansion at $\mathbf{lf}$.
\end{proof}

We now prove the (trickier) part of \cref{thm:C4:first}: the case of $\phi_k^{\mathbf{lf}}$.

\begin{proof}[Proof of \cref{thm:C4:first} for the case of the expansions at $\mathbf{lf}$.]
We will examine the top order expressions, and then we will linearize the equations to find transport equations for the higher order $\phi_j^{\mathbf{lf}}$ which we can solve. First, though, let us find the initial data for the equations on $\mathbf{lf}\n \mathbf{rf}$. Write $\phi_j^{\mathbf{lf}} = \phi_j = (\slash{k}_j,V_j,\omega_j)$. Since $\mathring{L} = L$ on $\mathbf{rf}$ by definition, the initial data for each term $V_j$ is $0$. Similarly, $\omega = 0$ on $\mathbf{rf}$, so each $\omega_j = 0$. The data for each $\slash{k}_j$ are harder to analyze. Recall from \cref{C:C3:ansatz} that the data are given on $\mathbf{rf}$ by
\[\slash{g} = e^{\Phi}\mathring{\slash{g}}e^{\eta \mathbf{T}} = e^{\Phi}\hat{\slash{g}},\] (on $\mathbf{rf}$ $\xi = 1$, and so $\eta = x$) for the conformal factor $e^{\Phi}$, where $\Phi$ satisfies
\[2L^2\Phi + (L\Phi)^2 + \Tr(\hat{\slash{g}}^{-1}\Lie_L\hat{\slash{g}})L\Phi + \frac{1}{2}|\Lie_L\hat{\slash{g}}|^2_{\hat{\slash{g}}} + L \Tr(\hat{\slash{g}}^{-1}\Lie_L\hat{\slash{g}})= 0\] in order for $\Ric_{LL} = 0$ on $\mathbf{rf}$ (see \cref{thm:C2:toporder}, \cref{thm:C3:toporder} or \cref{thm:C3:isreg}). Also recall from the proof of \cref{thm:C3:isreg} that $\Phi = \eta^2\Psi$ for some $\Psi$ smooth in a neighbourhood of $v=0$. In particular,
\begin{equation}\label{eq:first:eq1}j!\slash{k}_j|_{\{\xi = 1\}} = \mathring{\slash{g}}\pa_v^j (e^{\eta^2\Psi}e^{\eta T})|_{\{v = 0\}}.\end{equation} For $j \geq 1$, this is in $\eta C^{\infty}(\mathbf{rf}\n \mathbf{lf})$. Denote by $\mathring{\slash{\tr}}$ the trace with respect to $\mathring{\slash{g}}$. We also claim that for $j \geq 1$.
\[\mathring{\slash{\tr}}\slash{h}_j|_{\{\xi = 1\}} = \pa_v^j e^{\eta^2 \Psi}\Tr(e^{\eta T)}|_{\{\xi = 1\}} \in \eta^2 C^{\infty}(\mathbf{rf}\n \mathbf{lf}).\] Indeed, looking at \eqref{eq:first:eq1}, if even a single $\pa_v$ derivative hits the conformal factor, we pick up the desired $\eta^2$ factor. Now we examine what happens when all derivatives hit the factor $e^{\eta T}$.
\[\pa_v^j \Tr(e^{\eta T}) = \sum C_\alpha \Tr(e^{\eta T}(\eta \pa_v^\alpha T)^{\beta_\alpha}),\]
where the sum runs over all $\sum \alpha\beta_\alpha = j$. The only way to pick up a factor less than $\eta^2$ is if $\beta_\alpha = 1$ for a single $\alpha$, and $\beta_\alpha = 0$ for all others. In this case, the term is
\[\Tr(e^{\eta T}\eta \pa_v^{j} T).\] Since $T(0) = 0$, this is just $\eta \pa_v^j \Tr(T) = 0$, so this term vanishes.

Now let us find the behaviour of the top order $\phi_0$ and $\phi_1$. We know from \cref{thm:C2:toporder} that $\phi_0$ and $\phi_1$ are uniquely specified.
Following the procedure described there, it is easy to see that $\phi_0^ = (\mathring{\slash{g}},0,0)$ (after all, it must by assumption coincide with Minkowski data). Finding $\phi_1$ is trickier. Using the procedure, $\slash{k}_1$ must be taken so that $\Ric_{\eta\xi(-),\eta\xi(-)} = 0$ on $\mathbf{lf}$. Using \eqref{eq:C4:one}, this equation becomes
\begin{equation}\label{eq:first:eq2}0 = \Lie_{\bm{w}_1N}\slash{k}_1 - \slash{k}_1 - \frac{1}{2}(\mathring{\slash{\tr}}\slash{k})\mathring{\slash{g}} + 2\bm{w}_3^2\mathring{\slash{g}}.\end{equation}
Taking the trace of this,
\[0 = \Lie_{\bm{w}_1N}\mathring{\slash{\tr}}\slash{k}_1 - 2\mathring{\slash{\tr}}\slash{k}_1 = -4\bm{w}_3^2.\]
The initial data for $\mathring{\slash{\tr}}\slash{k}_1$ is in $\eta^2 C^{\infty}(\mathbf{lf})$, so \cref{thm:C4:firstsolve}, below, implies that \[\mathring{\slash{\tr}}\slash{k}_1 \in \eta^2\phgi{(\mathbf{rf},\mathbf{bf},\mathbf{if})}{(0,0,E_{\log,0})}(\mathbf{lf})\]
(observe that $\bm{w}_3^2 = \xi^2 \eta^2$, so there is no coefficient of $\xi^0$ and hence no $\log$ term to worry about).
Thus \eqref{eq:first:eq2} becomes
\begin{equation}\label{eq:first:eq2.1}\Lie_{\bm{w}_1N}\slash{k}_1 - \slash{k}_1 \in \eta^2\phgi{(\mathbf{rf},\mathbf{bf},\mathbf{if})}{(0,0,E_{\log,0})}(\mathbf{lf})\end{equation} with initial data in $\eta C^{\infty}(\mathbf{lf})$. \Cref{thm:C4:firstsolve}(and \cref{rk:C4:firstsolve} then imply that \[\slash{k}_1 \in \eta\phgi{(\mathbf{rf},\mathbf{bf},\mathbf{if})}{(0,0,E_{\log,0})}(\mathbf{lf}).\]

Using the procedure, $\omega_1$ must be taken so that $\Ric_{\bm{w}_2 L,\bm{w}_1 N} = 0$ on $\mathbf{lf}$. Using \eqref{eq:C4:three}, this becomes
\[\bm{w}_1N \omega_1 + \frac{1}{4}\bm{w}_1N\mathring{\slash{\tr}}\slash{k}_1 -\frac{1}{4}\mathring{\slash{\tr}}\slash{k}_1 = 0\] Thus by \cref{thm:C4:firstsolve}, $\omega_1 \in \eta\phgi{(\mathbf{rf},\mathbf{bf},\mathbf{if})}{(0,0,E_{\log,0})}(\mathbf{lf})$.

Next, let us determine $V_1$, which must be chosen so that \eqref{eq:C4:two} is $0$ on $\mathbf{lf}$. On $\mathbf{lf}$,
\[\pa_v Z = \mathring{\slash{g}}([\bm{w}_1N,V_1] + V_1/2,\bullet).\] Thus the equation becomes
\[\mathring{\slash{g}}([\bm{w}_1N,V_1] + V_1/2,\bullet) - 2\xi\eta \slash{\div} \slash{k}_1 + 2\xi\eta\slash{d}\slash{\tr}\slash{k}_1 = 0,\]
(all slash tensorial operations taken with respect to $\slash{k}_0 = \mathring{\slash{g}}$). Raising the equation via $\mathring{\slash{g}}$, \cref{thm:C4:firstsolve} applies again to show that $V_1 \in \eta\phgi{(\mathbf{rf},\mathbf{bf},\mathbf{if})}{(0,0,E_{\log,0})}(\mathbf{lf})$.

This proves the proposition for the cases $j=0,1$.
Now for $j \geq 0$ set
\[\phi_{(j)} = \phi_0 + \cdots + v^j \phi_j.\]

We will now prove a more precise version of the proposition by induction on $j$, namely we will show:
\begin{romanumerate}
\item $\slash{k}_j \in \eta\phgi{(\mathbf{rf},\mathbf{bf},\mathbf{if})}{(0,0,E_{\log,0})}(\mathbf{lf})$ for $j \geq 1$;
\item $\mathring{\slash{\tr}}\slash{k}_j \in \eta^2\phgi{(\mathbf{rf},\mathbf{bf},\mathbf{if})}{(0,0,E_{\log,0})}(\mathbf{lf})$ for $j \geq 1$;
\item $\omega_j,V_j \in \eta\phgi{(\mathbf{rf},\mathbf{bf},\mathbf{if})}{(0,0,E_{\log,0})}(\mathbf{lf})$ for all $j$;
\item $\mathring{\slash{g}}([\bm{w}_1N,V_j] + V_j/2,\bullet) \in \eta\phgi{(\mathbf{rf},\mathbf{bf},\mathbf{if})}{(0,0,E_{\log,0})}(\mathbf{lf})$ for all $j$;
\end{romanumerate}
(of course $\phi_0 = (\mathring{\slash{g}},0,0) \in\phgi{(\mathbf{rf},\mathbf{bf},\mathbf{if})}{(0,0,E_{\log,0})}(\mathbf{lf})$ as we have just computed).

Observe that we have shown all of these for $j=0,1$.
Assume the inductive hypothesis for $j-1$. We prove it for $j$. We may write over $\xi > 0$
\begin{equation}\label{eq:first:eq10}P(\phi_{(j-1)}+v^{j}\phi_j) = P(\phi_{(j-1)}) +v^{j-1} T_{\phi_{(j-1)},j}\phi_j + v^jC^{\infty}(\{\xi > 0\}),\end{equation}
where $T_{\phi_{(j-1)},j}$ is a linear operator which depends on $\phi_{(j-1),j}$. However the dependence on $\phi_{(j-1)}$ is smooth, and so \begin{equation}\label{eq:first:eq11}T_{\phi_{(j-1)},j}\phi_j =T_{\phi_{0},j}\phi_j + vC^{\infty}(\{\xi > 0\}).\end{equation}

By the inductive hypothesis, \[\phi_{(j-1)} \in \phgi{(\mathbf{rf},\mathbf{bf},\mathbf{if})}{(0,0,E_{\log,0})}(\mathbf{lf}).\] Thus,
\[P\phi_{(j-1)} \in \phgi{(\mathbf{rf},\mathbf{bf},\mathbf{if})}{(0,0,E_{\log,0})}(\mathbf{lf}),\] too. However, looking at the precise form of $P$, we may use the inductive hypothesis and the description of $\phi_0$ to see that there are no terms of order $\eta^0$, and so
\[P\phi_{(j-1)} \in \eta\phgi{(\mathbf{rf},\mathbf{bf},\mathbf{if})}{(0,0,E_{\log,0})}(\mathbf{lf}).\] Looking at the precise form of the first component of $P$, $P^1$, we in fact conclude that
\[P^1\phi_{(j-1)} \in \eta^2\phgi{(\mathbf{rf},\mathbf{bf},\mathbf{if})}{(0,0,E_{\log,0})}(\mathbf{lf}).\]

However by definition of the $\phi_i$, $P(\phi_{(j-1)}) \in v^{j-1}C^{\infty}(\{\xi > 0\})$, and so in fact
\begin{equation}\label{eq:first:eq12}P\phi_{(j-1)} \in v^{j-1}\eta\phgi{(\mathbf{rf},\mathbf{lf},\mathbf{if})}{(0,0,E_{\log,0})}(\mathbf{lf})\end{equation}
and
\begin{equation}\label{eq:first:eq13}P^1\phi_{(j-1)} \in v^{j-1}\eta^2\phgi{(\mathbf{rf},\mathbf{bf},\mathbf{if})}{(0,0,E_{\log,0})}(\mathbf{lf}).\end{equation}

Similarly, by definition
\begin{equation}\label{eq:first:eq14}P(\phi_{(j-1)}+v^j\phi_j) \in v^jC^{\infty}(\{\xi > 0\}).\end{equation} Thus, using \eqref{eq:first:eq10}--\eqref{eq:first:eq14}
\[T_{\phi_0,j}\phi_{j} = -(v^{j-1}P(\phi_{(j-1)})|_{\mathbf{bf}} \in \eta\phgi{(\mathbf{rf},\mathbf{lf},\mathbf{if})}{(0,0,E_{\log,0})}(\mathbf{lf}),\]
at least over $\mathbf{lf} \n \{\xi > 0\}$.
Write the components of $T_{\phi_0,j}$ \[T_{\phi_0,j} = (T^1,T^2,T^3).\] We additionally know from \eqref{eq:first:eq13} that
\[T^1\phi_j \in \eta^2\phgi{(\mathbf{rf},\mathbf{lf},\mathbf{if})}{(0,0,E_{\log,0})}(\mathbf{lf}).\]

We may record the components of $T_{\phi_0,j} = (T^1,T^2,T^3)$, which are obtained by performing a computation and discarding terms of order $v^{j}$ and higher, acting on $\phi_j$:
\begin{subequations}
\begin{align}
\label{eq:firstsolve:eqa}
T^1(\phi_j) &= j\Lie_{\bm{w}_1 N}\slash{k}_j - j\slash{k}_j -\frac{j}{2}(\slash{\tr}\slash{k}_j)\mathring{\slash{g}}\\
\label{eq:firstsolve:eqb}
\begin{split}
T^2(\phi_j) &= j\mathring{\slash{g}}([\bm{w}_1 N,V_j],\bullet)+\frac{j}{2}\mathring{\slash{g}}(V_j,\bullet)\\
&-2j\xi\eta\slash{\div}(\slash{k}_j) + 2j\xi\eta\slash{d}\slash{\tr}\Lie_{\bm{w}_2\mathring{L}}\slash{k}_j+4j\xi\eta\slash{d}\omega_j
\end{split}\\
\label{eq:firstsolve:eqc}
T^3(\phi_j) &= j(\bm{w}_1 N)\omega_j + \frac{j}{4}(\bm{w}_1N)\mathring{\slash{\tr}}\slash{k}_j -\frac{j}{4}\slash{\tr}\slash{k}_j.
\end{align}
\end{subequations}
Here, all tensorial operations are taken with respect to $\mathring{\slash{g}}$.

Let us look at \eqref{eq:firstsolve:eqa}. Taking its trace with respect to $\mathring{\slash{g}}$ gives
\[(\bm{w}_1 N)\slash{\tr}\slash{k}_j - 2\slash{\tr}\slash{k}_j \in \eta^2\phgi{(\mathbf{rf},\mathbf{bf},\mathbf{if})}{(0,0,E_{\log,0})}(\mathbf{lf}),\]
with initial data in $\eta^2 C^\infty$.
\Cref{thm:C4:firstsolve} shows that for some $R \in C^{\infty}(S^2;\R)$, $\slash{\tr}\slash{k}_j - R\eta^2\log\xi \in \eta^2\phgi{(\mathbf{rf},\mathbf{bf},\mathbf{if})}{(0,0,E_{\log,0})}(\mathbf{lf})$.

However, by compatibility, this would mean that the coefficient of $v^j$ in the expansion at $\mathbf{lf}$ of the coefficient of $\eta^2$ in the expansion of $\slash{k}$ at $\mathbf{bf}$ would have a $\log \xi$ in it. But by the part of \cref{thm:C4:first} concerning the expansion at $\mathbf{bf}$, this does not happen.\footnote{This ``magical cancellation'' is related to the apparent lower-order behaviour of some of the terms in the expression for $P^1$. This, in fact, goes back to equation \eqref{eq:C3:secondorder}, where the terms with a $1/(1-u)$ in the denominator cancelled. Tracking carefully, this is why no $\log$s appear in the expansions at $\mathbf{bf}$. In principle, one could inductively show that these terms cancel without recourse to compatibility, but this is quite delicate.}

Equation \eqref{eq:firstsolve:eqa} becomes
\[\Lie_{\bm{w}_1 N}\slash{k}_j - \slash{k}_j \in \eta^2\phgi{(\mathbf{rf},\mathbf{bf},\mathbf{if})}{(0,0,E_{\log,0})}(\mathbf{lf}),\]
with data in $\eta C^{\infty}(\mathbf{lf})$, and hence using \cref{thm:C4:firstsolve} and \cref{rk:C4:firstsolve} again, \[\slash{k}_j \in \eta\phgi{(\mathbf{rf},\mathbf{bf},\mathbf{if})}{(0,0,E_{\log,0})}(\mathbf{lf}),\] too. Then one treats \eqref{eq:firstsolve:eqb} and \eqref{eq:firstsolve:eqc} the third component and second component in kind in the same way as when solving for the first-order behaviour, obtaining $\omega_j \in \eta\phgi{(\mathbf{rf},\mathbf{bf},\mathbf{if})}{(0,0,E_{\log,0})}(\mathbf{lf})$ and
$V_j \in \eta\phgi{(\mathbf{rf},\mathbf{bf},\mathbf{if})}{(0,0,E_{\log,0})}(\mathbf{lf})$.
\end{proof}

\begin{lem}\label{thm:C4:firstsolve}Fix $k,j \in \Z$ and $2j +k > 0$. Consider the following ODE for a section $f$ of $T_q^p S^2$ ($(p,q) \in \N\times \N$) over $\mathbf{lf}$:
\begin{equation}\label{eq:firstsolve:eq1}\Lie_{\bm{w}_1N}f + \frac{k}{2}f = R \in \eta^j \phgi{(\mathbf{rf},\mathbf{bf},\mathbf{if})}{(0,0,E_{\log,0})}(\mathbf{lf};T_q^pS^2),\end{equation}
with data for $f$ at $\{\xi = 1\}$ satisfying $f|_{\{\xi = 1\}} \in \eta^j C^{\infty}(\{\xi = 1\})$. Then \[f \in \eta^j\phgi{(\mathbf{rf},\mathbf{bf},\mathbf{if})}{(0,0,E_{\log,0})}(\mathbf{lf};T_q^pS^2).\]
If $2j + k = 0$, then one modifies the conclusion as follows: let $\widebar{R} \in C^{\infty}(S^2_;T^p_q S^2)$ denote the coefficient of $\xi^0$ which appears when expanding at $\{\xi = 0\}$ the coefficient of $\eta^j$ of the expansion of $R$ at $\{\eta = 0\}$. Then
\[f + 2\eta^j\log\xi \widebar{R} \in \phgi{(\mathbf{rf},\mathbf{bf},\mathbf{if})}{(0,0,E_{\log,0})}(\mathbf{lf};T_q^pS^2).\]
\end{lem}
\begin{rk}\label{rk:C4:firstsolve}In particular if $\widebar{R} = 0$, then $f \in \phgi{(\mathbf{rf},\mathbf{bf},\mathbf{if})}{(0,0,E_{\log,0})}(\mathbf{lf};T_q^pS^2)$. Notice that this condition is automatically met if $R \in \eta^{j'} \phgi{(\mathbf{rf},\mathbf{bf},\mathbf{if})}{(0,0,E_{\log,0})}(\mathbf{lf};T_q^pS^2)$ for any $j' > j$.\end{rk}
\begin{proof}On $\mathbf{lf}$, $\bm{w}_1N = \eta\pa_\eta - \frac{1}{2}\xi\pa_\xi$.

We first find a solution $f$ to \eqref{eq:firstsolve:eq1} in series at all boundary faces in $\mathbf{lf}$, and then turn it into an actual solution by perturbing it by an element of $\phgi{(\mathbf{rf},\mathbf{bf},\mathbf{if})}{(\emptyset,\emptyset,\emptyset)}(\mathbf{lf})$ which is $0$ at $\{\xi = 1\}$. Since solutions to ODEs are unique, this must coincide with $f$.

\textbf{Step 1: Series solution.} We start with $\mathbf{rf} = \{\xi = 1\}$. We may write \[\ R \sim \sum R_i(1-\xi)^i\] where $R_i \in \eta^j C^{\infty}(\{\xi = 1\}$ and we are looking for coefficients
\[f_i \in \eta^j C^{\infty}(\{\xi = 1\}\] with $f_0 = f(1,\bullet)$ so that
\[f \sim \sum f_i(1-\xi)^i\] solves \eqref{eq:firstsolve:eq1} in series at $\{\xi = 1\}$. Thus, $f_i$ are required to solve
\begin{equation}\label{eq:firstsolve:eq5.1}\frac{i+1}{2}f_{i+1} + \frac{k}{2}f^{\mathbf{rf}}_{i-1} + \eta\pa_\eta f_i = R_i\end{equation} for $i \geq 0$. These equations may be solved recursively for $f_{i+1}$.

Next let us solve at $\mathbf{bf} = \{\eta = 0\}$. Again write
\[R \sim \sum R_i\eta^{i+j},\]
where $R_i \in \phgi{(\mathbf{rf},\mathbf{if})}{(0,E_{\log,0})}$, and expand the data
\[f(1,\bullet) \sim \sum \tilde{f}_i \eta^{i+j},\]
where $\tilde{f}_i \in C^{\infty}(\{\eta = 1-\xi = 0\})$. We are looking for
\[f_i \in \phgi{(\mathbf{rf},\mathbf{if})}{(0,E_{\log,0})}\] so that
\[f \sim \sum f_i\eta^{i+j}\] solves \eqref{eq:firstsolve:eq1} in series at $\{\eta = 0\}$ with $f_i|_{\{\xi = 1\}} = \tilde{f}_i$. Thus, $f_i$ are required to solve
\[-\frac{1}{2}\xi\pa_\xi f_i + (k/2+j+i)f_i = R_i.\]

Notice that the assumptions mean that $k/2+j+i \geq 0$.
Assume first that $k/2+j+i > 0$. Then a version of \cref{thm:C4:secondsolve}\footnote{applied to $\mathbf{bf}\n \mathbf{lf}$ rather than $\mathbf{bf}$ although the proof is manifestly the same.} shows that $f_i \in \phgi{(\mathbf{rf},\mathbf{bf})}{(0,E_{\log,0})}$. If $k/2+j+i = 0$, then $k/2+j = 0$ and $i = 0$. Let us expand $R_{0} = R_{0,0} + \xi\tilde{R}_{0}$, where $R_{0,0} = \widebar{R}$ and
$\tilde{R}_0 \in \phgi{(\mathbf{rf},\mathbf{if})}{(0,E_{\log})}$.
Write $f_i = -2\log\xi R_{0,0} + \widebar{f}_i$ Then $\widebar{f}_i$ satisfies
\[-\frac{1}{2}\pa_\xi \widebar{f}_i = \tilde{R}_i,\]
with initial data $\widebar{f}_i|_{\{\xi = 1\}} = \tilde{f}_i$. Thus, by integrating, $\widebar{f}_i \in \phgi{(\mathbf{rf},\mathbf{if})}{(0,E_{\log,0})}$.

Before moving on to $\mathbf{if} = \{\xi = 0\}$, let us verify the compatibility conditions between the series at $\mathbf{bf}$ and the series at $\mathbf{rf}$. This is the same argument as in \cref{thm:C3:compat}. Write $f_i^{\mathbf{rf}}$ and $f_i^{\mathbf{bf}}$ for the coefficients we found above for the expansions at $\mathbf{rf}$ and $\mathbf{bf}$, respectively,and find asymptotic sums $f^{\mathbf{rf}}$ and $f^{\mathbf{bf}}$ respectively for each. We may expand
\[f^{\mathbf{bf}} \sim f^{\mathbf{bf}}_0 + (1-\xi)f^{\mathbf{bf}}_1 + (1-\xi)^2f^{\mathbf{bf}}_2 + \cdots\]
and
\[f^{\mathbf{rf}} \sim f_0^{\mathbf{rf}} + (1-\xi)f_1^{\mathbf{rf}} + (1-\xi)^2f_2^{\mathbf{rf}}+\cdots,\]
where the $f^{\mathbf{rf}}_i$ and $f^{\mathbf{bf}}_i$ do not depend on $\xi$. Write $f^{\mathbf{bf}}_{(n)}$ and $f^{\mathbf{rf}}_{(n)}$ for the sum of the first $n$ terms in the above series, respectively. The compatibility condition is equivalent to $f_i^{\mathbf{rf}} - f^{\mathbf{bf}}_i \in \eta^{\infty}C^{\infty}(\{\xi = 1\})$. We will use induction on $i$. Observe that $f^{\mathbf{bf}}_0 - f^{\mathbf{rf}}_0 \in \eta^{\infty}C^{\infty}(\{\xi = 1\})$, since by construction $f^{\mathbf{bf}}_{\{\xi = 1\}} \sim f^{\mathbf{rf}}$ in Taylor series at $\{\eta = 1-\xi = 0\}$.

Now we use induction, so assume that the claim is true for all $\ell \leq i-1$ and we prove it for $i$. Write $f^{\mathbf{bf}}-f^{\mathbf{bf}}_{(i-1)} = (1-\xi)^{i}f^{\mathbf{bf}}_{i} + (1-\xi)^{i+1}\eta^{j}C^{\infty}(\{\xi = 1\})$. Then
\begin{equation}\label{eq:firstsolve:eq5}(\Lie_{\bm{w}_1 N} + k/2)f^{\mathbf{bf}}_{(i-1)} + \frac{i}{2}f^{\mathbf{bf}}_i \xi(1-\xi)^{i-1} + (1-\xi)^iC^{\infty}(\{\xi = 1\}).\end{equation}

Since all terms in $f^{\mathbf{bf}}_{(i-1)}$ and $f^{\mathbf{rf}}_{(i-1)}$ are the same in Taylor series $\eta = 0$, we may expand
\begin{equation}\label{eq:firstsolve:eq6}(\Lie_{\bm{w}_1 N} + k/2)f^{\mathbf{bf}}_{(i-1)} = (\Lie_{\bm{w}_1 N} + k/2)f^{\mathbf{rf}}_{(i-1)} + \sum_{\ell=0}^{i-1} \kappa_\ell(1-\xi)^{\ell},\end{equation} where the $\kappa_\ell$ do not depend on $\xi$ and are rapidly vanishing at $\{\eta = 0\}$.

Also, by assumption, $(\Lie_{\bm{w}_1 N} + k/2)f^{\mathbf{bf}}-R$ vanishes in series at $\{\eta = 0\}$. In particular, we may expand
\begin{equation}\label{eq:firstsolve:eq7}(\Lie_{\bm{w}_1 N} + k/2)f^{\mathbf{bf}} = R + \sum_{\ell=0}^{i} (1-\xi)^\ell \lambda_\ell + (1-\xi)^{i+1}C^{\infty}(\{\xi = 1\}),\end{equation} where each $\lambda_\ell$ does not depend on $\xi$ and is rapidly vanishing at $\{\eta = 0\}$.

Equating the coefficients of $(1-\xi)^{i-1}$ and taking $\xi \to 1$, it follows from \eqref{eq:firstsolve:eq5}--\eqref{eq:firstsolve:eq7} that
\[\lambda_{i-1} + R_{i-1} = \kappa_{i-1} + \eta\pa_\eta f^{\mathbf{rf}}_{i-1} + \frac{i}{2}f^{\mathbf{bf}}_i + \frac{k}{2}f_i\]
(here, $R_{i-1}$ denotes the coefficient of $(1-\xi)^{i-1}$ of $R$ in its expansion at $\xi = 1$). From \eqref{eq:firstsolve:eq5.1}, it follows that
\[\lambda_{i-1} = \kappa_{i-1} + \frac{i}{2}(f^{\mathbf{bf}}_i-f^{\mathbf{rf}}_i),\]
which shows $f^{\mathbf{bf}}_i = f^{\mathbf{rf}}_i$ up to a rapidly vanishing function at $\{\eta = 0\}$, as desired.

Now we move on to obtaining a series expansion at $\{\xi = 0\}$. As usual, expand
\[R \sim \sum R_{i,p}\xi^{i+j}\log^p \xi,\] where $R_{i,p} \in \eta^jC^{\infty}(\{\xi = 0\})$, and we are looking for $f_{i,p} \in \eta^jC^{\infty}(\{\xi = 0\})$ such that
\[f \sim \sum f_{i,p}\xi^{i}\log^p \xi\] solves \eqref{eq:firstsolve:eq1} in series at $\{\xi = 0\}$. Thus, $f_{i,p}$ are required to solve
\begin{equation}\label{eq:firstsolve:eq2}\eta\pa_\eta f_{i,p} + \frac{k-i}{2}f_{i,p} - \frac{(p+1)}{2}f_{i,p+1} = R_{i,p}.\end{equation}

We already know that $f^{\mathbf{bf}}$ solves \eqref{eq:firstsolve:eq1} in series at $\eta = 0$. Thus, if we expand
\[f^{\mathbf{bf}} \sim \sum \tilde{f}_{i,p}\xi^{i}\log^p \xi,\] then the $\tilde{f}_{i,p}$ solve \eqref{eq:firstsolve:eq2} modulo an error $E_{i,p}$ rapidly vanishing at $\{\eta = 0\}$.

Let $F_{i,p} = f_{i,p}-\tilde{f}_{i,p}$. Now we may apply a theorem on ODEs with singular derivatives from \cref{C:A3:Sing}. Using \cref{thm:A3:GoursatSing} (or rather the special case that there is no $v$ dependence), it is easy to inductively find $F_{i,p}$ which are rapidly vanishing at $\eta = 0$ and solve \[\eta\pa_\eta F_{i,p} + \frac{k-i}{2}F_{i,p} - \frac{(p+1)}{2}F_{i,p+1} = -E_{i,p} \in \eta^\infty C^\infty\] (of course a term $E_{i,p}$ is $0$ if $R_{i,p} = 0$ and $\tilde{f}_{i,q} = 0$ for all $q \geq p$). This gives us $f_{i,p}$ and shows for free that the compatibility conditions are satisfied. However, we need to be certain that the terms $f_{0,p}$ have the correct form, i.e.\ if $2j + k > 0$ then $f_{0,p} = 0$ for $p \geq 1$ and if $2j+k = 0$ then $f_{0,1} = -2\eta^j\widebar{R}$ and $f_{0,p} = 0$ for $p \geq 2$.

Let us first look at the case $2j+k > 0$. In this case, $\tilde{f}_{0,p}$ and $R_{0,p}$ are $0$ for $p \geq 1$, and so $E_{0,p} \equiv 0$ for $p \geq 1$ and thus $F_{0,p} = 0$, and hence $f_{0,p} = 0$, too, so the condition is satisfied.

Not let us look at the case $2j+k = 0$. In this case $f^{\mathbf{bf}} + 2\widebar{R}\eta^j\log \xi \in \phgi{(\mathbf{bf},\mathbf{if})}{(0,E_{\log,0})}$, so it follows that $\tilde{f}_{0,1} = -2\eta^j\widebar{R}$ and $\tilde{f}_{0,p} = 0$ for $p \geq 2$. Since $R_{0,p} = 0$ for $p \geq 1$, this means like above that $f_{0,p} = 0$, too. For $p = 1$, the equation for $\tilde{f}_{i,1}$ reads
\[\eta\pa_\eta (-2\eta^j\widebar{R}) + k/2(-2\eta^j\widebar{R}) = E_{0,1}.\]

However in this case, $j = -k/2$, so the left-hand side if $0$. This means $E_{0,1} \equiv 0$, so $F_{0,1} = 0$ and hence $f_{0,1} = \tilde{f}_{0,1} = -2\eta^j\widebar{R}$, as desired.

The expansions having been found and the compatibility conditions having been verified, we may use Borel's lemma to find a solution
\[\tilde{f} \in \eta^j\phgi{(\mathbf{rf},\mathbf{bf},\mathbf{if})}{(0,0,E_{\log,0})}(\mathbf{lf})\] with the correct initial data, solving
\[\Lie_{\bm{w}_1N}\tilde{f} + \frac{k}{2}\tilde{f} - R = -S \in \phgi{(\mathbf{rf},\mathbf{bf},\mathbf{if})}{(\emptyset,\emptyset,\emptyset)}(\mathbf{lf}).\] 
\textbf{Step 2: Actual solution.} To find $f$, we just need to find $g \in \phgi{(\mathbf{rf},\mathbf{bf},\mathbf{if})}{(\emptyset,\emptyset,\emptyset)}(\mathbf{lf})$ with $0$ data solving
\[\Lie_{\bm{w}_1N}g + \frac{k}{2}g = S \in \phgi{(\mathbf{rf},\mathbf{bf},\mathbf{if})}{(\emptyset,\emptyset,\emptyset)}(\mathbf{lf}).\]
If this holds, then $f = \tilde{f}+g$.

One may write down an explicit formula for $g$:
\[g(\xi,\eta) = 2\xi^k\int_\xi^1 S(t,t^{-2}\eta\xi^2)t^{-k-1}\ dt.\]
By assumption for all $A,B,C \in \N$, $S(v,\eta) \in \xi^A\eta^B(1-\xi)^C \mathcal A(\mathbf{lf})$. It is easy to see that for all $A,B,C \in \N$
\[g(\xi,\eta) \in \xi^A\eta^B(1-\xi)^C L^{\infty}_{\mathrm{loc}}(\mathbf{lf}).\] Differentiating via $\xi\pa_\xi, \eta\pa_\eta$ or $\pa_{\theta}$, one finds that, for $a,b,c \in \N$
\[(\xi\pa_\xi)^a(\eta\pa_\eta)^b\pa_\theta^c g = \xi^k\int_1^\xi S_{a,b,c}(t,t^{-2}\eta\xi^2)t^{-k-1}\ dt + T_{a,b,c}(\xi,\eta),\] where
\[S_{a,b,c},T_{a,b,c} \in \phgi{(\mathbf{rf},\mathbf{bf},\mathbf{if})}{(\emptyset,\emptyset,\emptyset)}(\mathbf{lf}).\] Thus, one concludes that for any $A,B,C\in \N$
\[(\xi\pa_\xi)^a(\eta\pa_\eta)^b\pa_\theta^c g \in \xi^A\eta^B(1-\xi)^C L^{\infty}_{\mathrm{loc}}(\mathbf{lf}),\] as well. Thus $g \in \phgi{(\mathbf{rf},\mathbf{bf},\mathbf{if})}{(\emptyset,\emptyset,\emptyset)}(\mathbf{lf})$, as desired.
\end{proof}

\subsection{Expansion at \texorpdfstring{$\mathbf{if}$}{if}}
We prove \cref{thm:C4:second}, modulo the top order which is more complex.
\begin{proof}[Proof of \cref{thm:C4:second}]
We start with some preliminary observations. If $j \geq 1$, and $\psi \in \phgd(\Mp)$, $\psi_j \in \xi^j \phg{(0,0,0,E_{\log})}(\Mp)$ then
\begin{equation}\label{eq:second:eq1}P(\psi + \psi_j) = P(\psi) + R,\end{equation}
where $R \in \xi^j \phg{(0,0,0,E_{\log})}(\Mp)$. This can be proved by first observing that \cref{thm:C4:algebra} shows that $P(\psi+\psi_j) \in \phgd(\Mp)$, after which a formal series argument may be carried out. Furthermore, if $\psi_j$ has the form
\[\psi_j = \sum_{q=0}^p\xi^j\log^q\xi \psi_{j,p} + \tilde{\psi}_j,\] where $\psi_{j,p} \in C^{\infty}(\mathbf{if}\n \{\digamma > 0\})$ and $\tilde{\psi}_j \in \xi^{j+1} \phg{(0,0,0,E_{\log})}(\Mp)$, then, denoting by $\psi_{0,0} = \psi|_{\{\xi = 0\}}$\footnote{In other words the coefficient of $\xi^0\log^0 \xi$ in the expansion.}
\begin{equation}\label{eq:second:eq2}P(\psi + \psi_j) = P(\psi) + \xi^j\log^p \xi T_{j,\psi_{0,0}}\psi_{j,p} + \sum_{q=0}^{p-1} \xi^j\log^q (T_{j,\psi_{0,0}}\psi_{j,q} + L_{q+1,\psi_{0,0}}\psi_{j,q+1}) + R,\end{equation}
various $T$,
and the where $R \in \xi^{j+1} \phg{(0,0,0,E_{\log})}(\Mp)$, and the various $T$, $L$ are linear operators depending on the indicated arguments.\footnote{The operators $T$ do not depend on $p$ since the dependence on $j$ comes only from hitting a of order term $\xi^j\log^p \xi$ with $\xi\pa_\xi$, which generates a term or order $j\xi^{j}\log^p \xi$ (which manifestly only depends on $j$), whereas $L$ does not depend on $j$ because it only comes from when $\xi\pa_\xi$ lands on a term of order$\xi^j\log^q \xi$, and generates a term of the form $q\xi^{j}\log^{q-1}\xi$.} This can also be proven by using \cref{thm:C4:algebra} followed by a formal series argument.

We can now carry on with the main part proof.

It suffices to prove by induction on $k$ that there exist unique coefficients $\phi_{j,p} \in \phg{(0,0)}(\mathbf{if}\n \{\digamma > 0\})$, $j \leq k$ ($\phi_{0,p} = 0$ for $p \geq 1$), such that their sum
\[\phi_{(k)} = \sum_{i \leq k,p} \phi_{i,p}\xi^i\log^p \xi\] solves
\[P\phi_{(k)} \in \xi^{k+1} \phgi{(\mathbf{lf},\mathbf{bf},\mathbf{if})}{(0,0,E_{\log})}(\Mp)\] and the $\phi_{(k)}$ satisfy the compatibility conditions. To prove the existence of $\phi$ from this, take an asymptotic sum of the $\phi_{k,p}$. Then by \eqref{eq:second:eq1}, then for all $k \geq 0$,
\[P\phi = P(\phi_{(k)} + (\phi-\phi_{(k)})) = P(\phi_{(k)}) + \xi^{k+1} \phgi{(\mathbf{lf},\mathbf{bf},\mathbf{if})}{(0,0,E_{\log})}(\Mp) = \xi^{k+1} \phgi{(\mathbf{lf},\mathbf{bf},\mathbf{if})}{(0,0,E_{\log})}(\Mp),\]
and hence
\[P\phi \in \phgi{(\mathbf{lf},\mathbf{bf},\mathbf{if})}{(0,0,\emptyset)}(\Mp).\]

Aruging uniqueness is similar. If $\psi \in \phgd(\Mp)$ solves $P\psi \in \phgi{(\mathbf{lf},\mathbf{bf},\mathbf{if})}{(0,0,\emptyset)}(\Mp)$, denoting by $\psi_{(k)}$ the sum of terms up to order $\xi^k$, then, by \eqref{eq:second:eq1},
\[P\psi = P(\psi_{(k)} + (\psi-\psi_{(k)}) = P\psi_{(k)} + \xi^{k+1} \phgi{(\mathbf{lf},\mathbf{bf},\mathbf{if})}{(0,0,E_{\log})}(\Mp),\] and so
\[P\psi_{(k)} \in \xi^{k+1}\phgi{(\mathbf{lf},\mathbf{bf},\mathbf{if})}{(0,0,E_{\log})}(\Mp),\]
and thus $\psi_{(k)} = \phi_{(k)}$ by uniqueness.

So let us prove the claim by induction.

Write $\phi^{\mathbf{if}}_{j,p} = \phi_{j,p} = (\slash{k}_{j,p},V_{j,p},\omega_{j,p})$.

We handle the base case $k=0$ in \cref{C:C4:topo}. Unfortunately $\slash{k}_{0,0}$ and $\omega_{0,0}$ do not have easy formulae, although $V_{0,0} = 0$. 

Let us assume the inductive hypothesis is true for $k-1$. We prove it for $k$. Using the coefficients $\phi_\ell^{\mathbf{lf}}$ and $\phi_\ell^{\mathbf{bf}}$ of the series solution provided by \cref{thm:C4:first}, together with the coefficients $\phi^{\mathbf{if}}_{i,p}$ for $0 \leq i \leq k-1$, we may use Borel's lemma to find some $\tilde{\phi}_{(k-1)} = \tilde{\phi} \in \phgd(\Mp)$ whose expansion at $\mathbf{lf}$ and $\mathbf{bf}$ have coefficients $\phi_{\ell}^{\mathbf{lf}}$ and $\phi_{\ell}^{\mathbf{bf}}$, respectively, and whose expansion at $\mathbf{if}$ has the coefficient of $\xi^i\log^p \xi$ for $0 \leq i \leq k-1$ equal to $\phi^\mathbf{if}_{i,p}$.\footnote{To be pedantic, in order to apply Borel's lemma we also need a sequence of functions which will be the coefficients of $\xi^i\log^p \xi$ for $i > k$. We may just choose these coefficients arbitrarily provided they satisfy the compatibility conditions.} Notice that $\tilde{\phi} - \phi_{(k-1)} \in \xi^k\phg{(0,0,0,E_{\log})}(\Mp)$, and so by \eqref{eq:second:eq1}
\[P(\tilde{\phi}) \in \xi^k\phg{(0,0,0,E_{\log})}(\Mp).\]
However, since $\tilde{\phi}$ solves $P\tilde{\phi}$ in series at $\mathbf{rf}$, $\mathbf{lf}$ and $\mathbf{bf}$ in fact
\[P(\tilde{\phi}) \in \xi^k\phg{(\emptyset,\emptyset,\emptyset,E_{\log})}(\Mp).\]
Now we can prove the existence of the $\phi^{\mathbf{if}}_{k,q}$. Let $p$ be the largest natural number such that the coefficient of $\xi^k\log^p \xi$ of $P(\tilde{\phi})$ or the coefficient of $\tilde{\phi}$ is nonzero. Consider some $\widebar{\phi}_{k,p} \in C^{\infty}(\mathbf{if}\n \{\digamma > 0\})$ which is rapidly vanishing at $\{\eta = 0\}$ and $\{\xi = 0\}$, and suppose it solves
\begin{equation}\label{eq:second:eq4}P(\tilde{\phi}+ \xi^k\log^p \xi \widebar{\phi}_{k,p}) = \sum_{q=0}^{p-1} \xi^k\log^q \xi c_{k,q} + \xi^{k+1}\phgi{(\mathbf{rf},\mathbf{lf},\mathbf{bf},\mathbf{if})}{(0,0,0,E_{\log})}(\Mp),\end{equation} for some smooth coefficients $c_{k,q}$, i.e.\ it kills the coefficient of $\xi^k \log^p \xi$. Observe that since $P(\tilde{\phi})$ and $\widebar{\phi}_{k,p}$ are both rapidly vanishing at $\{\eta = 0\}$, $\{v = 0\}$, so are the $c_{k,q}$, provided they exist.
Set $\phi^{\mathbf{if}}_{k,p} = \tilde{\phi}_{k,p} + \widebar{\phi}_{k,p}$. Then $\phi^{\mathbf{if}}_{k,p}$ automatically satisfies the compatibility conditions.
Since 
\begin{equation}\label{eq:second:eq100}\tilde{\phi} = \phi_{(k-1)} + \sum_{q=0}^p \xi^k \log^q \tilde{\phi}_{k,q} + \xi^{k+1}\phgi{(\mathbf{rf},\mathbf{lf},\mathbf{bf},\mathbf{if})}{(0,0,0,E_{\log})}(\Mp),\end{equation}
\eqref{eq:second:eq2} and \eqref{eq:second:eq4} show that $\phi^{\mathbf{if}}_{k,p}$ solves
\begin{align*}P(\phi_{(k-1)} + \xi^k\log^p \xi \phi^{\mathbf{if}}_{k,p}) &= P(\tilde{\phi}\hspace{-0.75pt} + \hspace{-0.75pt}\xi^k\log^p \widebar{\phi}_{k,p})+ \sum_{q=0}^{p-1} \xi^k\log^q \xi a_{k,q}\hspace{-0.75pt} + \hspace{-0.75pt} \xi^{k+1}\hspace{-0.25pt}\phgi{(\mathbf{rf},\mathbf{lf},\mathbf{bf},\mathbf{if})}{(0,0,0,E_{\log})}(\Mp)\\
&= \sum_{q=0}^{p-1} \xi^k\log^q \xi b_{k,q} + \xi^{k+1}\phgi{(\mathbf{rf},\mathbf{lf},\mathbf{bf},\mathbf{if})}{(0,0,0,E_{\log})}(\Mp),
\end{align*}
for some coefficients $a_{k,q}$, $b_{k,q}$. The upshot is that adding $\xi^k\log^p \xi\phi^{\mathbf{if}}_{k,p}$ kills the term of order $\xi^k \log^p\xi$ in the expansion. 
So let us solve \eqref{eq:second:eq4}.
From \eqref{eq:second:eq2}, we have that
\begin{align*}
 P(\tilde{\phi}+ \xi^k\log^p \xi \widebar{\phi}_{k,p}) &= P(\tilde{\phi}) + \xi^k\log^p \xi^k T_{j,\tilde{\phi}_{0,0}} \widebar{\phi}_{k,p}\\
 &+\sum_{q=0}^{p-1} \xi^k\log^q \xi d_{k,q} + \xi^{k+1}\phgi{(\mathbf{rf},\mathbf{lf},\mathbf{bf},\mathbf{if})}{(0,0,0,E_{\log})}(\Mp),\end{align*}
for some smooth coefficients $d_{k,q}$. 
Extracting the coefficient of $\xi^k\log^p \xi$ means that it is necessary and sufficient for
\begin{equation}\label{eq:second:eq101}T_{j,\tilde{\phi}_{0,0}} \widebar{\phi}_{k,p} = -(\xi^{-k}\log^{-p}P(\tilde{\phi}))|_{\{\xi = 0\}}.\end{equation}
Since $P(\tilde{\phi}) \in \xi^k\phg{(\emptyset,\emptyset,\emptyset,E_{\log})}(\Mp)$, in fact $(\xi^{-k}\log^{-p}P(\tilde{\phi}))|_{\{\xi = 0\}} \in v^\infty \eta^\infty C^{\infty}(\mathbf{if}\n \{\digamma > 0\})$. We may record the components of $T_{j,\tilde{\phi}_{0,0}} = (T^1,T^2,T^3)$. Set $\widebar{\phi}_{k,p} = (\slash{k}_{k,p},V_{k,p},\omega_{k,p})$. We record:
\begin{align*}
\begin{split}
T^1(\widebar{\phi}_{k,p}) &= \Lie_{\eta \partial_{\eta}}\Lie_{\partial_v}\slash{k}_{k,p} - \frac{j}{2}\Lie_{\partial_v}\slash{k}_{k,p}\\
&+\left(\frac{1}{4}\slash{\tr}(\Lie_{\eta\partial_\eta}\slash{k}_{0,0})-1\right)(\Lie_{\partial_v}\slash{k}_{k,p})\\
&+ \frac{1}{4}\left(-\slash{k}_{0,0}(\slash{k}_{k,p},\Lie_{\eta\partial_{\eta}}\slash{k}_{0,0})+\slash{\tr}\Lie_{\eta\partial_{\eta}}\slash{k}_{k,p} - \frac{j}{2}\slash{\tr}\slash{k}_{k,p}\right)\Lie_{\partial_v}\slash{k}_{0,0}\\
&+\frac{1}{4}(\slash{\tr}\Lie_{\partial_v}\slash{k}_{0,0})\left(\Lie_{\eta\partial_{\eta}}\slash{k}_{k,p} -\left(\frac{j}{2}+2\right)\slash{k}_{k,p}\right)\\
&+ \frac{1}{4}\left(-\slash{k}_{0,0}(\slash{k}_{k,p},\Lie_{\partial_v}\slash{k}_{0,0})+ \slash{\tr}\Lie_{\partial_v}\slash{k}_{k,p}\right)\left(\Lie_{\eta\partial_{\eta}}\slash{k}_{0,0}-2\slash{k}_{0,0}\right)\\
&-\frac{1}{2}\left[\Lie_{\partial_v} \slash{k}_{k,p}\times \Lie_{\eta\partial_{\eta}}\slash{k}_{0,0} + \Lie_{\eta\partial_{\eta}}\slash{k}_{0,0}\times \partial_v \slash{k}_{k,p}\right]\\
&-\frac{1}{2}\left[\Lie_{\partial_v}\slash{k}_{0,0}\times\left(\Lie_{\eta\partial_{\eta}}\slash{k}_{k,p}-\frac{j}{2}\slash{k}_{k,p}\right) +\left(\Lie_{\eta\partial_{\eta}}\slash{k}_{k,p}-\frac{j}{2}\slash{k}_{k,p}\right)\times\Lie_{\partial_v}\slash{k}_{0,0}\right]\\
&+\frac{1}{2}\left[ \Lie_{\partial_v}\slash{k}_{0,0}\times \slash{k}_{k,p} \times \Lie_{\eta\partial_{\eta}}\slash{k}_{0,0} + \Lie_{\eta\partial_{\eta}}\slash{k}_{0,0}\times\slash{k}_{k,p}\times \Lie_{\partial_v}\slash{k}_{0,0}\right]
\end{split}\\
\begin{split}
T^2(\widebar{\phi}_{k,p}) &= \Lie_{\partial_v}\left(\left(\Lie_{\eta\partial_{\eta}} V_{k,p}\right)^{\flat} + \frac{(1-j)}{2}V_{k,p}^{\flat}\right)\\
&+\frac{1}{2}\hspace{-2pt}\left(\slash{\tr}\Lie_{\partial_v}\slash{k}_{0,0}\right)\left(\left(\Lie_{\eta\partial_{\eta}} V_{k,p}\right)^{\flat} + \frac{(1-j)}{2}V_{k,p}^{\flat}\right)-\hspace{-2.6pt}2\partial_v\omega_0\left(\hspace{-1.1pt}\left(\Lie_{\eta\partial_{\eta}} V_{k,p}\right)^{\flat} + \frac{(1-j)}{2}V_{k,p}^{\flat}\right)
\end{split}\\
\begin{split}
T^3(\widebar{\phi}_{k,p}) &= \eta\partial_{\eta}\partial_v \omega_{k,p} - \frac{j}{2}\partial_v \omega_{k,p}\\
&+\frac{1}{4}\eta\partial_{\eta}\slash{\tr}\Lie_{\partial_v}\slash{k}_{k,p} -\frac{1}{4}\eta\partial_{\eta}\slash{k}_{0,0}(\slash{k}_{k,p},\Lie_{\partial_v}\slash{k}_{0,0})\\
&-\left(\frac{j}{8}+\frac{1}{4}\right)\slash{\tr}\Lie_{\partial_v}\slash{k}_{k,p} +\left(\frac{j}{8}+\frac{1}{4}\right)\slash{k}_{0,0}(\slash{k}_{k,p},\Lie_{\partial_v}\slash{k}_{0,0})\\
&-\frac{1}{8}\Tr(\slash{k}_{0,0}^{-1}(\Lie_{\eta\partial_\eta}\slash{k}_{0,0})\slash{k}_{0,0}^{-1}\slash{k}_{k,p}\slash{k}_{0,0}^{-1}(\Lie_{\partial_v}\slash{k}_{0,0}))\\ &-\frac{1}{8}\Tr(\slash{k}_{0,0}^{-1}(\Lie_{\partial_v}\slash{k}_{0,0})\slash{k}_{0,0}^{-1}\slash{k}_{k,p}\slash{k}_{0,0}^{-1}(\Lie_{\eta\partial_{\eta}}\slash{k}_{0,0}))\\
&+\frac{1}{8}\slash{k}_{0,0}(\Lie_{\partial_v}\slash{k}_{k,p},\Lie_{\eta\partial_{\eta}}\slash{k}_{0,0}) + \frac{1}{8}\slash{k}_{0,0}\left(\Lie_{\partial_v}\slash{k}_{0,0},\Lie_{\eta\partial_{\eta}}\slash{k}_{k,p}-\frac{j}{2}\slash{k}_{k,p}\right).
\end{split}
\end{align*}

Here, all tensorial operations are with respect to $\slash{k}_{0,0}$, and the notation $\bullet^\flat$ indicates raising a vector field via $\slash{k}_{0,0}$, i.e.\ $A^\flat = \slash{k}_{0,0}(A,\bullet)$ for any vector field $A$.

Although these equations look horrendous\footnote{In fact they are horrendous!}, they are just the sort of system considered in \cref{C:A3:Sing}. By assumption, the right-hand side for these is in $v^\infty \eta^\infty C^\infty(\mathbf{if} \n \{\digamma > 0\})$, so, after lowering the expression for $T^2$ via $\slash{k}_{0,0}$, we may apply \cref{thm:A3:GoursatSing} to obtain a unique solution in $v^\infty \eta^\infty C^\infty$, which is our $\widebar{\phi}_{k,p}$.

Now we may iterate this argument, finding inductively, for each $0 \leq q \leq p-1$ (starting with $q=p-1$), a section
$\widebar{\phi}_{k,q}$ satisfying
\begin{equation}P\left(\left(\tilde{\phi}+\sum_{r=q+1}^p \xi^k\log^r \xi \widebar{\phi}_{k,r}\right)\hspace{-0.75pt} +\hspace{-0.75pt} \xi^k\log^q \xi \widebar{\phi}_{k,q}\right)\hspace{-0.75pt} = \hspace{-0.75pt}\sum_{r=0}^{q-1} \xi^k\log^r \xi c_{k,r} + \xi^{k+1}\hspace{-0.25pt}\phgi{(\mathbf{rf},\mathbf{lf},\mathbf{bf},\mathbf{if})}{(0,0,0,E_{\log})}(\Mp),\end{equation} for some smooth coefficients $c_{k,r}$ (by induction the argument of $P$ in parentheses already solves the equation up to order $\xi^k \log^q \xi$, so we just need to kill the error term of order $\xi^k\log^q \xi$) and then setting
\[\phi_{k,q}^{\mathbf{if}} = \tilde{\phi_{k,q}}^{\mathbf{if}} + \widebar{\phi}_{k,q}\]
so that \[P(\phi_{(k-1)} + \xi^k\log^p \xi \phi^{\mathbf{if}}_{k,p} + \cdots + \xi^k\log^q \xi \phi^{\mathbf{if}}_{k,q}) = \sum_{r=0}^{q-1}\xi^k \log^r\xi b_{k,r} + \xi^{k+1}\phgi{(\mathbf{rf},\mathbf{lf},\mathbf{bf},\mathbf{if})}{(0,0,0,E_{\log})}(\Mp).\] 

This completes the existence part of the inductive step. We move onto uniqueness. Suppose $\psi_{j,p}$, $j \leq k$, $\psi_{0,p} = 0$ for $p \geq 1$ are coefficients such that their sum
\[\psi_{(k)} = \sum_{i \leq k, p} \psi_{i,p}\xi^i\log^p \xi\] solves
\[P\psi_{(k)} \in \xi^{k+1}\phgi{(\mathbf{lf},\mathbf{bf},\mathbf{if})}{(0,0,E_{\log})}(\Mp),\]
and satisfy the compatibility conditions. We need to show that $\psi_{j,p} = \phi_{j,p}$ for all $j$, $p$. Write
\[\psi_{(k-1)} = \sum_{i \leq k-1, \ p} \psi_{i,p}\xi^i\log^p \xi.\] Then by \eqref{eq:second:eq1}
\[P\psi_{(k-1)} \in \xi^{k}\phgi{(\mathbf{lf},\mathbf{bf},\mathbf{if})}{(0,0,E_{\log})}(\Mp),\] and so by induction, $\psi_{(k-1)} = \phi_{(k-1)}$. Let $p$ be the largest natural number such that the coefficient of $\xi^k\log^p \xi$ of $P(\tilde{\phi})$ or $\tilde{\phi}$ is nonzero, and let $p'$ be the maximum of $p$ and the largest natural number $q$ such that $\psi_{k,q}$ is nonzero. Since the the compatibility conditions are satisfied by assumption, for all $q$ \[\widebar{\psi}_{k,q} := \psi_{k,q} - \tilde{\phi}_{k,q}\] is rapidly vanishing at $\mathbf{lf}$ and $\mathbf{bf}$ for all $q \leq p$. By assumption
\[P(\psi_{(k-1)} + \xi^k\log^{p'} \xi\psi_{k,p'}) = \sum_{q=0}^{p'-1} \xi^k \log^q \xi c_{k,q} + \xi^{k+1}\phgi{(\mathbf{rf},\mathbf{lf},\mathbf{bf},\mathbf{if})}{(0,0,0,E_{\log})}(\Mp)\] for some smooth coefficients $c_{k,q}$.
Using \eqref{eq:second:eq100} and the inductive hypothesis it follows (like above) that
\[P(\tilde{\phi} + \xi^k\log^{p'} \xi\widebar{\psi}_{k,p'}) = \sum_{q=0}^{p'-1} \xi^k \log^q \xi d_{k,q} + \xi^{k+1} + \xi^{k+1}\phgi{(\mathbf{rf},\mathbf{lf},\mathbf{bf},\mathbf{if})}{(0,0,0,E_{\log})}(\Mp)\]
for some smooth coefficients $d_{k,q}$. As above, this means that
\[T_{j,\tilde{\phi}_{0,0}}\widebar{\psi}_{k,p'} = -(\xi^{-k}\log^{-p'} \xi P(\tilde{\phi}))|_{\{\xi = 0\}}.\]
If $p' = p$, this is the same equation as \eqref{eq:second:eq101}, which means by uniqueness that $\widebar{\psi}_{k,p'} = \widebar{\phi}_{k,p'}$, and so $\psi_{k,p'} = \phi_{k,p'}$. If $p' > p$, then the right-hand side is $0$, and so by uniqueness $\widebar{\psi}_{k,p'} = 0$, and so $\psi_{k,p'} = 0 = \phi_{k,p'}$. Arguing as above using induction, it follows for $0 \leq q < p'$ that
\[T_{j,\tilde{\phi}_{0,0}}\widebar{\psi}_{k,q} = -(\xi^{-k}\log^{-q} \xi P(\tilde{\phi}))|_{\{\xi = 0\}},\] too, and thus like above, $\psi_{k,q} = \phi_{k,q}$, too.
\end{proof}

\section{Top order at \texorpdfstring{$\mathbf{if}$}{if}}
\label{C:C4:topo}
In this section we complete the proof of \cref{thm:C4:second} by showing how to find the top order $\phi^{\mathbf{if}}_{0,0} = (\slash{k}_{0,0},V_{0,0},\omega_{0,0})$. The treatment of how to find $\slash{k}_{0,0}$ and $\omega_{0,0}$ will be important since it is how we are going to find the blowup in \cref{C:C6} (although also important, we will see that $V_{0,0}$ is trivial so there is little content to use later).
For this section, we drop the subscript $0,0$ as well as the superscript $\mathbf{if}$ from our notation.

\subsection{The behaviour of \texorpdfstring{$V$}{V} to top order.}

Perhaps unusually, the first equation we analyze at $\mathbf{if}$ is \eqref{eq:C4:two}, which will give us the top order behaviour of $V$. We will show:
\begin{prop}\label{thm:C4:topV}Suppose $\phi \in \phgd(\Mp)$. Then $V = 0$ is the unique choice of $V$ which satisfies \eqref{eq:C4:two} $=0$ on $\mathbf{if}\n\{\digamma > 0\}$ and satisfies the compatibility conditions with the series solutions at $\mathbf{bf}$ and $\mathbf{if}$.\end{prop}
\begin{proof}
At $\mathbf{if}$, \eqref{eq:C4:two} $=0$ takes the form
\begin{equation}\label{eq:C4:V}0 = \Lie_{\pa_v} Z + \frac{1}{2}\slash{\tr}(\Lie_{\pa_v} \slash{k})Z -2(\pa_v\omega)Z.\end{equation} \Cref{thm:C4:first} shows that if the compatibility conditions are to be satisfied, then $V = 0$ to top order at $\mathbf{lf} \un \mathbf{bf}$ and $V = 0$ to first order as well as $\mathbf{bf}$. This gives us the data required to solve \eqref{eq:C4:V} uniquely if $V$ is to satisfy the compatibility conditions, as we now indicate. The equation is a linear transport equation in $Z$, and so $Z=0$ is the unique solution, since $Z=0$ on $\mathbf{lf}$ by compatibility. Thus
\[0 = \Lie_{\bm{w}_1N} V + \frac{1}{2}V = \Lie_{\eta\pa_\eta} V + \frac{1}{2}V.\] The data we are providing means that both $V,\pa_\eta V = 0$ at $\{\eta = 0\}$. Thus we may set $V = \eta W$, and $W$ satisfies
\[0 = \Lie_{\eta\pa_\eta}W +\frac{3}{2}\eta W,\] with $W|_{\{\eta = 0\}} = 0$, and thus $W = 0$ and so $V=0$.

Lastly, we need to check that all the compatibility conditions are satisfied. Let $\widebar{\phi} = (\widebar{k},\widebar{V},\widebar{\omega})$ be a series solution at $\mathbf{bf}$ and $\mathbf{lf}$, provided by \cref{thm:C4:first}. We need that for all $k$
\[0 = \Lie_{\pa_\eta}^k V = \Lie_{\pa_\eta}^k \widebar{V}, \ 0 = \Lie_{\pa_v}^k V = \Lie_{\pa_v}^k \widebar{V}\]
on $\mathbf{if} \n \mathbf{bf}$ and $\mathbf{if}\n\mathbf{lf}$, respectively.

By assumption, with $\widebar{Z} = \widebar{k}([\mathbf{w}_1 N,\widebar{V}] + 1/2V,\bullet)$, the condition that is \eqref{eq:C4:two} rapidly vanishing at $\mathbf{bf}\n\mathbf{if}$ and $\mathbf{lf}\n\mathbf{if}$ on $\mathbf{if}$ is
\[\Lie_{\pa_v} \widebar{Z} + \frac{1}{2}\slash{\tr}(\Lie_{\pa_v} \slash{k})\widebar{Z} -2(\pa_v\widebar{\omega})\widebar{Z} \in \eta^\infty v^\infty C^{\infty}(\mathbf{if}).\]

Arguing as above, this means that $\widebar{Z} \in \eta^\infty v^\infty C^{\infty}(\mathbf{if})$, and thus $\widebar{V} \in \eta^\infty v^\infty C^{\infty}(\mathbf{if})$. This shows that the compatibility conditions are satisfied.
\end{proof}

\subsection{The behaviour of \texorpdfstring{$\slash{k}$}{slash-k} to top order.} \label{C:C4:topk}

Let us now analyze \eqref{eq:C4:one} at $\mathbf{if}$ to obtain the top-order behaviour of $\slash{k}$. This will be the hardest component of $\phi$ to examine. We will show in this subsection:
\begin{prop}\label{thm:C4:topk}Suppose $\phi \in \phgd(\Mp)$. Then there exists a unique smooth $\slash{k}$, defined on $\mathbf{if}\n\{\digamma > 0\}$ such that $\phi$ satisfies \eqref{eq:C4:one} $=0$ on $\{\digamma > 0\}$ and satisfies the compatibility conditions with the series solutions at $\mathbf{bf}$ and $\mathbf{lf}$. The unique solution may not be continued past $\{\digamma = 0\}$.\end{prop}
We split the proof of \cref{thm:C4:topk} across this subsection.

From \cref{thm:C4:topV}, $V = 0$, so \eqref{eq:C4:one} takes the form
\begin{align}
\label{eq:C4:oneII}
\begin{split}0 &= \Lie_{\eta\partial_\eta}\Lie_{\partial v}\slash{k} + \frac{1}{4}\slash{\tr}(\Lie_{\eta\partial_\eta} \slash{k})\Lie_{\partial v}\slash{k} +\frac{1}{4}\slash{\tr}(\Lie_{\partial v}\slash{k})\Lie_{\eta\partial_\eta}\slash{k}\\
&-\frac{1}{2}\left(\Lie_{\partial v}\slash{k}\times \Lie_{\eta\partial_\eta}\slash{k} + \Lie_{\eta\partial_\eta}\slash{k}\times\Lie_{\partial v}\slash{k}\right)-\Lie_{\partial_v} \slash{k} - \frac{1}{2}\slash{\tr}(\Lie_{\partial_v} \slash{k})\slash{k}.
\end{split}
\end{align}
If the compatibility conditions are to be satisfied, then $\slash{k} = \mathring{\slash{g}}$ along $\{v = 0\}$ and $\{\eta = 0\}$ by \cref{thm:C4:first}.

To analyze \eqref{eq:C4:oneII}, it will be convenient to decompose the derivatives of $\slash{k}$ into their trace and tracefree parts. Write
\begin{align*}\Lie_{\pa_v} \slash{k} &= \slash{k}\frac{f}{2} + \slash{k}F\\
\Lie_{\pa_\eta} \slash{k} &= \slash{k}\frac{h}{2} + \slash{k}H,
\end{align*}
where $f$, $h$ are scalars, and $F$, $H$ are tracefree $\slash{k}$-symmetric tensors. Equation \eqref{eq:C4:oneII} implies that $f,h,F,H$ satisfy the system\footnote{Here, and for the rest of this thesis, we may write $\pa_\eta$ and $\pa_v$ for $\Lie_{\pa_\eta}$ and $\Lie_{\pa_v}$ after choosing a trivialization of the fibre bundle over $\mathbf{if}$, i.e.\ choosing coordinates $(\eta,v,\theta)$ in which the projection is $(\eta,v,\theta) \mapsto (\eta,v)$. Indeed, one may interpret a tracefree fibre-tensor $A$ as a family of matrices or vectors in $\R^{2\times 2} = \R^4$, and for $X = \pa_\eta$ or $\pa_v$, in canonical coordinates $\Lie_X A$ is the same as the component-wise derivative $XA$, since the commutators $[X,\Theta]$ with $\Theta$ a basis vector field $\pa_{\theta^a}$ is $0$. Cf.\ \cref{rk:C4:LieDer}.}
\begin{subequations}
\label{eq:C4:fhFH}
\begin{align}
\label{eq:C4:fh}
\begin{split}
0 &= \pa_\eta f + \frac{1}{2}fh - \frac{2}{\eta}f\\
0 &= \pa_v h + \frac{1}{2}fh - \frac{2}{\eta}f
\end{split}\\
\label{eq:C4:FH}
\begin{split}
0 &= \pa_\eta F + \left(\frac{1}{4}h-\frac{1}{\eta}\right)F + \frac{1}{4}fH+ \frac{1}{2}[H,F]\\
0 &= \pa_v H + \left(\frac{1}{4}h-\frac{1}{\eta}\right)F + \frac{1}{4}fH + \frac{1}{2}[F,H].
\end{split}
\end{align}
\end{subequations}

The notation $[\bullet,\bullet]$ indicates the usual commutator of matrices (or type $(1,1)$ tensors).

Set $\tilde{f} = \eta^{-2}f$, $\tilde{F} = \eta^{-1}F$. We may rewrite \eqref{eq:C4:fhFH} as

\begin{subequations}
\label{eq:C4:fthFtH}
\begin{align}
\label{eq:C4:fth}
\begin{split}
0 &= \pa_\eta \tilde{f} + \frac{1}{2}\tilde{f}h\\
0 &= \pa_v h + \frac{\eta^2}{2}\tilde{f}h - 2\eta \tilde{f}
\end{split}\\
\label{eq:C4:FtH}
\begin{split}
0 &= \pa_\eta \tilde{F} + \frac{1}{4}h\tilde{F} + \frac{1}{4}\tilde{f}H + \frac{1}{2}[H,\tilde{F}]\\
0 &= \pa_vH + \frac{\eta}{4}h\tilde{F} + \frac{\eta^2}{4}\tilde{f}H + \frac{\eta}{2}[\tilde{F},H] - \tilde{F}.
\end{split}
\end{align}
\end{subequations}

The combined systems for $(f,F,h,H)$ and $(\tilde{f},\tilde{F},h,H)$ are nonlinear Goursat problems, which are considered in \cref{C:A3:Goursat}. We must provide initial data for $(h,H)$ along $\{v = 0\}$ and for $(\tilde{f},\tilde{F})$ along $\{\eta = 0\}$. 
\begin{lem}\label{thm:C4:initdata}If the compatibility conditions are satisfied, then $\slash{k}|_{\{v = 0\} \un \{\eta = 0\}} = \mathring{\slash{g}}$, $f \in \eta^2 C^{\infty}(\{\digamma > 0\})$, $F \in \eta C^{\infty}(\{\digamma > 0\})$ and the initial data are:
\begin{align*}
h|_{\{v = 0\}} &= 0\\
H|_{\{v = 0\}} &= 0\\
\tilde{f}|_{\{\eta = 0\}} &= -\mathbf{E}\\
\tilde{F}|_{\{\eta = 0\}} &= \pa_v \mathbf{T}.
\end{align*}
\end{lem}
\begin{proof}
From \cref{thm:C4:first}, the data for $\slash{k}$ is as described. It follows that \[\Lie_{\pa_\eta} \slash{k} = \Lie_{\pa_\eta} \mathring{\slash{g}} = 0\] on $\mathbf{lf}$, and so the data for $h$ and $H$ is as described. The data for $\tilde{f}$ and $\tilde{F}$ are harder to provide, since they must come from higher-order terms in the series expansion of $\slash{k}$ at $\mathbf{bf}$.

However, from \cref{thm:C4:first}, we know that on $\mathbf{if}$,
\[\slash{k} = \mathring{\slash{g}}(1+\eta\mathbf{T} + \eta^2S_0) + \eta^3 C^{\infty}(\mathbf{if}),\]
where 
\[S_0 = - \frac{1}{2}\int_0^v \mathbf{E}(s)\ ds + \frac{1}{4}\Tr(\mathbf{T}^2).\]

Thus, if the compatibility conditions are to be satisfied, then at $\mathbf{if}$,
\begin{align*}\Lie_{\pa_v} \slash{k} &= \mathring{\slash{g}}\left(\eta\pa_v\mathbf{T} + \eta^2\left(-\frac{1}{2}\mathbf{E} + \frac{1}{2}\Tr(\mathbf{T}\pa_v\mathbf{T})\right)\right) + \eta^3 C^{\infty}(\mathbf{if})\\
&= \slash{k}\left(\eta \pa_v \mathbf{T} - \frac{\eta^2}{2}\mathbf{E} + \frac{\eta^2}{2}\Tr(\mathbf{T}\pa_v\mathbf{T})- \eta^2\mathbf{T}\pa_v\mathbf{T} + \eta^3C^{\infty}(\mathbf{if})\right).
\end{align*}

Now recall that since $\mathbf{T}$ is fibrewise tracefree linear map on a $2$-dimensional vector space, $\mathbf{T}\pa_v\mathbf{T}$ is a scalar, and so the last two terms of order $\eta^2$ cancel. Thus, if the compatibility conditions are satisfied, we see that $f|_{\{\eta = 0\}}, \pa_\eta f|_{\{\eta = 0\}} = 0$, and so $f \in \eta^2 C^{\infty}(\{\digamma > 0\})$ and $\tilde{f}|_{\{\eta = 0\}} = -\mathbf{E}$, and $F|_{\{\eta = 0\}} = 0$, and so $F \in \eta C^{\infty}(\{\digamma > 0\})$ and $\tilde{F}|_{\{\eta = 0\}} = \pa_v \mathbf{T}$.
\end{proof}

Recalling \cref{thm:C4:rect}, the uniqueness statement of \cref{thm:A3:wp} shows that:
\begin{lem}Solutions to \eqref{eq:C4:oneII}, the system \eqref{eq:C4:fhFH} and the system \eqref{eq:C4:fthFtH}, with the data specified as in \cref{thm:C4:initdata}, are unique on $\{\digamma > 0\}$, provided they exist.\end{lem}

These last two lemmas show that if there is a solution satisfying the compatibility conditions, then it is unique. So we are left with establishing existence on $\{\digamma > 0\}$, together with checking that all the compatibility conditions are indeed satisfied. Let us begin with the latter.
\begin{lem}\label{thm:C4:topkcompat}Smooth solutions to \eqref{eq:C4:oneII} on $\{\digamma > 0\}$ with the given data for $\slash{k}$, $\tilde{f}$, $\tilde{F}$, $h$, $H$ satisfy the compatibility conditions.\end{lem}
\begin{proof}
Let $\widebar{\phi}= (\widebar{\slash{k}},\widebar{V},\widebar{\omega})$ be a series solution at $\mathbf{bf}$ and $\mathbf{lf}$, provided by \cref{thm:C4:first}. It suffices to show that $\slash{k}-\widebar{\slash{k}} \in \eta^\infty v^\infty C^{\infty}(\{\digamma > 0\})$.

Let us define $\tilde{\widebar{f}}$, $\tilde{\widebar{F}}$, $\widebar{h}$, $\widebar{H}$ analogously to $\tilde{f}$, $\tilde{F}$, $h$, $H$, i.e.\ by
\begin{align*}
\pa_v \widebar{k} &= \widebar{\slash{k}}(\eta^2/2\tilde{\widebar{f}} + \eta\tilde{\widebar{F}})\\
\pa_\eta \widebar{k} &= \widebar{\slash{k}}(\widebar{h} + \widebar{H}).
\end{align*}
Since, by \cref{thm:C4:topV}, $\widebar{V} = V$ modulo a vector field rapidly vanishing at $\{v = 0\}$ and $\{\eta = 0\}$, $\tilde{\widebar{f}}$, $\tilde{\widebar{F}}$, $\widebar{h}$, $\widebar{H}$ satisfy the system \eqref{eq:C4:fthFtH} with the same initial data as $\tilde{f}$, $\tilde{F}$, $h$, $H$, but with a rapidly vanishing error. The $\pa_\eta$--equations in both then establish that $\tilde{f} = \tilde{\widebar{f}}$, $\tilde{F} = \tilde{\widebar{F}}$ on $\{v = 0\}$, since these are ODEs with the same data and same coefficients, and likewise the $\pa_v$--equations establish that $h = \widebar{h}$, $H = \widebar{H}$ on $\{\eta = 0\}$.

Using the equation, it now follows that $\pa_\eta \tilde{f} = \pa_\eta \tilde{\widebar{f}}$, $\pa_\eta \tilde{F} = \pa_\eta \tilde{\widebar{F}}$ along $\{v = 0\} \un \{\eta = 0\}$, and similarly $\pa_v h = \pa_v \widebar{h}$, $\pa_v H = \pa_v \widebar{H}$ along $\{v = 0\} \un \{\eta = 0\}$. 

Now differentiate the equations using $\pa_v$ and $\pa_\eta$ and iterate the argument. We conclude that $\tilde{f}$, $\tilde{F}$, $h$, $H$ are equal to $\tilde{\widebar{f}}$, $\tilde{\widebar{F}}$, $\widebar{h}$, $\widebar{H}$ up to a smooth function rapidly decreasing at $\{v = 0\} \un \{\eta = 0\}$.

By definition $\pa_v \slash{k} = \slash{k}(f/2+F)$ and $\pa_v \widebar{\slash{k}} = \widebar{\slash{k}}(\widebar{f}/2+\widebar{F})$. Since $\widebar{f}/2 + \widebar{F}$ and $f/2 + F$ differ by a rapidly vanishing function, it follows that on $\{v = 0\} \un \{\eta = 0\}$, $\pa_v \slash{k} = \pa_v \widebar{\slash{k}}$. Similarly, since $\pa_\eta \slash{k} = \slash{k}(h/2 + H)$ and $\pa_\eta \widebar{\slash{k}} = \widebar{\slash{k}}(\widebar{h}/2 + \widebar{H})$, it follows that $\pa_\eta \slash{k} = \pa_\eta \widebar{\slash{k}}$ on $\{v = 0\} \un \{\eta = 0\}$. Now differentiate the equations for $\pa_v \slash{k}$ and $\pa_\eta \slash{k}$ via $\pa_v$ and $\pa_\eta$ and iterate the argument. We conclude $\slash{k}-\widebar{\slash{k}} \in \eta^\infty v^\infty C^{\infty}(\{\digamma > 0\})$.
\end{proof}

Now we focus on existence. While the systems \eqref{eq:C4:fhFH} and \eqref{eq:C4:fthFtH} are obviously equivalent, it is not clear that every solution $(f,F,h,H)$ to \eqref{eq:C4:fhFH} comes from a solution to \eqref{eq:C4:oneII}. We will for the rest of this thesis work almost exclusively with the systems \eqref{eq:C4:fhFH} or \eqref{eq:C4:fthFtH}, so we need to know that we may consider these equations rather than \eqref{eq:C4:oneII}.
\begin{prop}\label{thm:C4:backtometric}Let $(f,F,h,H)$ be a smooth solution to \eqref{eq:C4:fhFH} on $\{\digamma > 0\}$ (or indeed on any ``rectangle'' $[0,v_0]\times[0,\eta_0]\times U$, $v_0,\ \eta_0 > 0 \subseteq \{\digamma > 0\}$, $U \subseteq S^2$) with the specified data. Then there exists a solution $\slash{k}$ to \eqref{eq:C4:oneII} such that $\pa_v\slash{k} = \slash{k}(f/2 + F)$ and $\pa_\eta \slash{k} = \slash{k}(h/2 + H)$ with data $\slash{k} = \mathring{\slash{g}}$ on $\{v = 0\}\un\{\eta =0\}$.\end{prop}
\begin{proof}
Fix a rectangle $\mathcal R = [0,v_0]\times[0,\eta_0]\times U \subseteq \{\digamma > 0\}$ for $v_0,\eta_0 > 0$ and $U \subseteq S^2$.
Define $\slash{k}$ by requiring it to solve the linear ODE 
\[\pa_v\slash{k} = \slash{k}(f/2 + F)\] with data $\slash{k}|_{\{v = 0\}} = \mathring{\slash{g}}$. We need to verify:
\begin{romanumerate}
\item $\slash{k}_{\{\eta = 0\}} = \mathring{\slash{g}}$;
\item $\pa_\eta \slash{k} = \slash{k}(h/2 + H)$;
\item $\slash{k}$ solves \eqref{eq:C4:oneII};
\item $\slash{k}$ is a Riemannian metric on $TS^2$.
\end{romanumerate}
Property (iv) entails:
\begin{enumerate}[label= (iv)\alph*., leftmargin=2.5\parindent]
\item $\slash{k}$ is non-degenerate;
\item $\slash{k}$ is symmetric;
\item $\slash{k}$ is positive.
\end{enumerate}
Property (iii) follows from (ii) and (iv) by direct substitution. Property (i) is clear since on $\{\eta = 0\}$, $f/2+ F = 0$ and so $\pa_v\slash{k} = 0$. Let us now show (iv)a. Set $K = \mathring{\slash{g}}^{\vspace{-5pt}-1}\slash{k}|_{\mathbf{if}\n\{\digamma > 0\}}$ to be a $(1,1)$ tensor. It suffices to show that $\det K \neq 0$. By Jacobi's formula
\[\pa_v \det K = \Tr(K^{-1}\pa_v K )\det K = f\det K\] if $\det K \neq 0$ at a point $p \in \mathcal R$. In particular, if $\det K > 0$ at $p$ then
\begin{equation}\label{eq:C4:logs}\pa_v \log \det K = f.\end{equation} We know that on $\{v = 0\}$, $\det K > 0$ since $\slash{k} = \mathring{\slash{g}}$. In particular, for any $\eta_1 \in [0,\eta_0]$ and $\theta_0 \in U$ fixed, for $\epsilon > 0$ small enough, $\det K > 0$ on $[0,\epsilon]\times\{\eta_1\}\times\{\theta_0\}$. Now \eqref{eq:C4:logs}, combined with an open-closed argument shows that $\det K > 0$ on $[0,v_0]\times\{\eta_1\}\times \{\theta_0\}$. Since $\eta_1,\ \theta_0$ were arbitrary, it follows that $\det K > 0$ on all of $\mathcal R$.

Let us now show (ii).
Since $\slash{k}$ is non-degenerate, we may write $\pa_\eta \slash{k} = \slash{k}(h'/2 + H')$ for some scalar $h'$ and tracefree $H'$. Using that $\Lie_{\pa_v}\Lie_{\pa_\eta}\slash{k} = \Lie_{\pa_\eta}\Lie_{\pa_v}\slash{k}$, we deduce that
\[\pa_\eta(f/2+F) - \pa_v(h'/2+H') = [f/2+F,h'/2+H'] = [F,H'].\]
Looking at the system \eqref{eq:C4:fhFH}, it is clear that
\[\pa_\eta(f/2+F) - \pa_v(h/2+H) = [F,H].\] Thus, taking the trace and tracefree parts, 
\begin{align*}
\pa_v (h-h') &= 0\\
\pa_v(H-H') &= [F,H-H'].
\end{align*}
Since $F$ is smooth, these are both ODEs. By assumption $h' = 0$, $H' = 0$ on $\{v = 0\}$ since $\slash{k} = \mathring{\slash{g}}$ on $\{v = 0\}$. Thus $h-h'$ and $H-H'$ are $0$ on $\{v = 0\}$ and solve an ODE for which the pair $(0,0)$ is a solution. Thus $h = h'$ and $H = H'$. This proves (ii).

Next we show (iv)b. Denote $k^\ast_{ab} = k_{ba}$. We need to show that $\slash{k} = \slash{k}^\ast$. Observe that $\slash{k}^\ast$ also satisfies \eqref{eq:C4:oneII}. We may write
\begin{equation}\label{eq:C4:Ihh}\pa_v \slash{k}^\ast = \slash{k}^\ast(f'/2 + F'), \ \pa_\eta \slash{k}^\ast = \slash{k}^\ast(h'/2 + H'),\end{equation} for $f'$, $h'$ scalars and $F'$, $H'$ tracefree. Since $\slash{k}^\ast = \mathring{\slash{g}}$ initially, $h' = 0$, $H' = 0$ on $\{v = 0\}$, and $f' = 0$, $F' = 0$ on $\{\eta = 0\}$. It suffices to show that 
\begin{subequations}\label{eq:C4:ih8this}
\begin{align}
\label{eq:C4:ih8thisa}
f' \in \eta^2C^{\infty}(\mathcal R), \ \ \tilde{f}' := \eta^{-2}f' = \tilde{f}\text{ on }\{\eta = 0\}\\
\label{eq:C4:ih8thisb}
\ F' \in \eta C^{\infty}(\mathcal R), \ \ \tilde{F}' := \eta^{-1}F' = \tilde{F}\text{ on }\{\eta = 0\}.\end{align}
\end{subequations}Indeed, since $\slash{k}^\ast$ solves \eqref{eq:C4:oneII}, \eqref{eq:C4:ih8this} shows that $\tilde{f}'$, $\tilde{F}'$, $h'$, $H'$ solve the system \eqref{eq:C4:fthFtH} with the same data as $\tilde{f}$, $\tilde{F}$, $h$, $H$, and so by the uniqueness statement of \cref{thm:A3:wp}, $\tilde{f}' = \tilde{f}$, $\tilde{F}' = \tilde{F}$, $h' = H$, $g' = G$. Thus from \eqref{eq:C4:Ihh}, the definition of $\slash{k}$ and property (ii), $\slash{k} = \slash{k}^\ast$, since they satisfy the same linear ODE with the same data.

Let us prove \eqref{eq:C4:ih8this}. By definition it is true that
\begin{equation}\label{eq:C4:Ihhh}\pa_v(\slash{k}^\ast) = (\pa_v\slash{k})^\ast = (f/2+F)\slash{k}^\ast.\end{equation}
Let $\slash{h}$ denote the $(2,0)$ tensor satisfying $\slash{h}\slash{k}^\ast = \delta^a_b$.\footnote{In other words, $\slash{h} = (\mathring{\slash{g}\slash{k}}^\ast)^{-1}\mathring{\slash{g}}$, where the inverse is well-defined since its argument is an invertible $(1,1)$ tensor. Notice that we need to keep track of ``left inverses'' of $\slash{k}^\ast$ and ``right inverses'' of $\slash{k}^\ast$ since $\slash{k}^\ast$ has not yet been proven to be symmetric.} Since $\slash{k}$ is smooth, so is $\slash{h}$. Thus from \eqref{eq:C4:Ihh}-\eqref{eq:C4:Ihhh}
\begin{equation}(f'/2+F') = \slash{h}(f/2+F)\slash{k}^\ast = f/2 +\slash{h}F\slash{k}^\ast\end{equation}
and so taking the trace and tracefree parts shows the regularity statement of \eqref{eq:C4:ih8thisa}-\eqref{eq:C4:ih8thisb}. Thus, we may divide the trace part by $\eta^2$ and the tracefree part by $\eta$ to obtain
\begin{equation}\label{eq:C4:why}\tilde{f'} = \tilde{f}, \ \tilde{F'} = \slash{h}\tilde{F}\slash{k}^\ast.\end{equation}
This shows the rest of \eqref{eq:C4:ih8thisa}.
At $\{\eta = 0\}$, $\tilde{F} = \pa_v \mathbf{T}$ and $\slash{k}= \mathring{\slash{g}}$ is symmetric. In particular $\slash{h} = \mathring{\slash{g}}^{\vspace{-10pt}-1}$. Recall that the short-pulse tensor $\mathbf{T}$ was defined to be $\mathring{\slash{g}}$-symmetric, and so on $\{\eta = 0\}$, $\tilde{F}$ is also $\mathring{\slash{g}}$-symmetric, and so
\[\mathring{\slash{g}}\tilde{F}= (\mathring{\slash{g}}\tilde{F})^\ast = \tilde{F}\mathring{\slash{g}}.\]
Thus, \eqref{eq:C4:why} implies that $\tilde{F'} = \tilde{F}$, which is the last part of \eqref{eq:C4:ih8thisb}.

Finally, we show (iv)c. The set of positive, symmetric matrices $(0,2)$ tensors is open in symmetric $(0,2)$ tensors, and disjoint from the set of symmetric nondegenerate $(0,2)$ tensors of other signatures. The $(0,2)$ tensors $\slash{k}$ form a continuous family of symmetric $(0,2)$ tensors which is at a point (really along all of $\{v = 0\}\un \{\eta= 0\}$) positive and symmetric. It follows that $\slash{k}$ is positive everywhere.
\end{proof}

The combined systems \eqref{eq:C4:fhFH}, \eqref{eq:C4:fthFtH} are upper triangular, in the sense that \eqref{eq:C4:fh} is an equation only on $f,\ h$, and \eqref{eq:C4:fth} is an equation only for $\tilde{f},h$, respectively. It thus makes sense to solve \eqref{eq:C4:fh}/\eqref{eq:C4:fth} first.
\begin{prop}\label{thm:C4:fh}With the given data, $f$, $h$ are given explicitly by the formulae 
\begin{align*}
f &= 2\pa_v \log \digamma = -\frac{4\eta^2\mathbf{E}}{\digamma}\\
h &=2\pa_\eta \log \digamma = -\frac{2\eta\int_0^v \mathbf{E}}{\digamma}.
\end{align*}
Of course, this means 
\[\tilde{f} = f/\eta^2 = -\frac{4\mathbf{E}}{\digamma}.\]

In particular, $\tilde{f},\ h$ exist and are smooth on $\{\digamma > 0\}$, and may not be continued past $\{\digamma = 0\}$.
\end{prop}
\begin{proof}
Recalling that $\digamma = 4-2\eta^2\int_0^v \mathbf{E}$, it is clear that $\tilde{f}$, $h$ have the correct data, and one may easily verify that $f$, $h$ solve \eqref{eq:C4:fh}.

This is slightly mysterious, however, so we provide a formal derivation. Let us assume we have a smooth solution $(\tilde{f},h)$.

The equation
\[0 = \pa_\eta \tilde{f} + \frac{1}{2}\tilde{f}h\] is a first order linear ODE for $\tilde{f}$ in terms of $h$, and therefore by the method of integrating factors has a solution
\begin{equation}\label{eq:fh:eq1}\tilde{f}(v,\eta) = \tilde{f}(v,0)\exp\left(-\frac{1}{2}\int_0^\eta h(v,t)\ dt\right).\end{equation}
Now $\pa_v h = \pa_\eta f$ by definition, and so
\[\pa_v\int_0^\eta h(v,t)\ dt = f(v,\eta) - f(v,0) = \eta^2\tilde{f}(v,\eta).\] Thus
\[\int_0^\eta h(v,t)\ dt = \int_0^v \eta^2\tilde{f}(s,\eta) \ ds + \eta^2\tilde{f}(0,\eta) = \int_0^v \eta^2\tilde{f}(s,\eta) \ ds\]
(the term $\tilde{f}(0,\eta)$ vanishes using \eqref{eq:fh:eq1} because $\tilde{f}(0,0) = 0$)
and so we may write \eqref{eq:fh:eq1} as
\begin{equation}\label{eq:fh:eq2}\tilde{f}(v,\eta) = \tilde{f}(v,0)\exp\left(-\frac{\eta^2}{2}\int_0^v \tilde{f}(s,\eta)\ dt\right).\end{equation}

Let us assume that $\mathbf{E}(v) > 0$ on $(v_0,1]$, and $\mathbf{E}(v) = 0$ on $[0,v_0]$. Uniqueness of the solution shows that $\tilde{f}= h = 0$ on $[0,v_0]\times[0,\infty)$, so the region $v > v_0$ is the only interesting region. On this region, we may differentiate \eqref{eq:fh:eq2} using $\pa_v$ to obtain
\begin{equation}\label{eq:fh:eq3}\pa_v \tilde{f}(v,\eta) = \pa_v(\log \tilde{f}(v,0))\tilde{f}(v,\eta) -\frac{\eta^2}{2}\tilde{f}^2.\end{equation} This is a Bernoulli differential equation, and so we may transform it into a linear equation by setting $\Phi = \tilde{f}^{-1}$ (let us ignore the issue of dividing by zero since we are only providing a formal derivation; one may simply check directly that the resulting formulae are indeed solutions). This transforms \eqref{eq:fh:eq3} into
\[\pa_v \Phi = - \pa_v(\log \tilde{f}(v,0))\Phi + \frac{\eta^2}{2},\] which has solution (obtained using the method of integrating factors)
\[\Phi(v,\eta) = \frac{\Phi(v_1,\eta)\tilde{f}(v_1,0)}{\tilde{f}(v,0)} +\frac{\eta^2}{2\tilde{f}(v,0)}\int_{v_1}^v\tilde{f}(s,0)\ ds,\] for any $v_1 > v_0$. Thus
\[\tilde{f}(v,\eta) = \frac{4\eta^2\tilde{f}(v,0)}{4\Phi(v_1,\eta)\tilde{f}(v_1,0)+2\eta^2\int_{v_1}^v \tilde{f}(s,0)\ ds} = -\frac{4\mathbf{E}}{4\tilde{f}(v_1,0)/\tilde{f}(v_1,\eta) -2\eta^2\int_{v_1}^v \mathbf{E}}.\]

From \eqref{eq:fh:eq2}, 
\[\tilde{f}(v_1,0)/\tilde{f}(v_1,\eta) = \exp\left(\frac{\eta^2}{2}\int_0^{v_1} f(v,t)\ dt\right),\]
and so $\tilde{f}(v_1,0)/\tilde{f}(v_1,\eta) \to 1$ as $v_1 \to 0$. Therefore we recover \[f = 2\pa_v\log \digamma= -\frac{4\eta^2\mathbf{E}}{\digamma}.\] From the ODE
\[\pa_v h = \pa_\eta f = 2\pa_v \pa_\eta \log \digamma,\] we recover
\[h = -2\pa_\eta \log\digamma = \frac{2\eta\int_0^v \mathbf{E}}{\digamma}.\]
\end{proof}

The formulae for $f$ and $h$ can easily be used to proved the formula for $\det \mathring{\slash{g}}^{-1}\slash{k}$.
\begin{cor}\label{thm:C4:www}The following formula is valid:
\[\det \mathring{\slash{g}}^{\vspace{-5pt}-1}\slash{k} = \frac{1}{16}\digamma^2.\]
\end{cor}
\begin{proof}
Set $K = \det \mathring{\slash{g}}^{\vspace{-5pt}-1}\slash{k}$.
From \eqref{eq:C4:logs} and \cref{thm:C4:fh},
\[\log \det K = \int_0^v f(s,\eta) \ ds = 2\log\digamma(v,\eta) - 2\log 4.\]
Thus $\det K = \frac{1}{16}\digamma^2$.
\end{proof}

Finally we treat the existence of $\tilde{F}$ and $H$.
\begin{prop}\label{thm:C4:FH}With $\tilde{f}$, $h$ given as above, the solution $(\tilde{F},H)$ to \eqref{eq:C4:FtH} with the given data exists on $\{\digamma > 0\}$ and is smooth.\end{prop}
\begin{proof}
By \cref{thm:C4:rect}, \cref{thm:A3:globalsolution} and \cref{thm:A3:globalsmooth}, it suffices to show that if $\mathcal R_0 = [0,v_0]\times[0,\eta_0]\times U \subseteq \{\digamma > 0\}$, for $v_0 > 0, \eta_0 > 0$ and $U \subseteq S^2$ compact, then there exists $M > 0$ depending only on $v_0$, $\eta_0$, $U$, such that if $0 < v_1 < v_0$ and $0 < \eta_1 < \eta_0$ and $(\tilde{F},H)$ is any solution to \eqref{eq:C4:FtH} on $\mathcal R_1 = [0,v_1]\times[0,\eta_1]\times U \subseteq \mathcal R_0$ with the given data, then \begin{equation}\label{eq:FH:eq1}|\tilde{F}| + |H| \leq M.\end{equation} Indeed, then there would be a solution on $\mathcal R_0$ for any choice of $v_0$, $\eta_0$ and compact $U \subseteq S^2$, and hence on all of $\{\digamma > 0\}$. Here, we define as our norm for a $(1,1)$ tensor $A$ by $|A| := |A|_{\mathring{\slash{g}}} = \sqrt{\Tr(A^\ast A)}$, where the adjoint of $A$ is taken with respect to $\mathring{\slash{g}}$.

We will do this in two steps. Observe that by \cref{thm:C4:backtometric}, $\slash{k}$ exists on $\mathcal R_1$. We first show that there is some $M'$, not depending on $\mathcal R_1$ and only on $\mathcal R_0$, such that
\begin{equation}\label{eq:FH:eq2}\Tr(\tilde{F}^2) + \Tr(H^2) \leq M'.\end{equation} Since $\tilde{F}$ and $H$ are $\slash{k}$-symmetric, the left hand side is the same as
\[|\tilde{F}|_{\slash{k}} + |H|_{\slash{k}},\] where, for a $(1,1)$ tensor $A$, $|A|_{\slash{k}} := \sqrt{\Tr(A^\ast A)}$, where here the adjoint is taken with respect to $\slash{k}$. In the second step we show that on $\mathcal R_1$, $\slash{k}$ and $\mathring{\slash{g}}$ are uniformly equivalent metrics, with constants depending only on $\mathcal R_0$, which in turn implies that $|A| \lesssim|A|_{\slash{k}}$ uniformly, which means \eqref{eq:FH:eq2} implies \eqref{eq:FH:eq1}. 

Let us prove \eqref{eq:FH:eq2}. The key observation is that if $A,B$ are $(1,1)$ tensors, then \[\Tr(A[A,B]) = 0.\] We thus derive from \eqref{eq:C4:FtH} the following equations for $\Tr(\tilde{F}^2)$ and $\Tr(H^2)$:
\begin{align}
\label{eq:FH:eq3}
\begin{split}
\pa_\eta \Tr(\tilde{F}^2) &= -\frac{h}{2} \Tr(\tilde{F}^2) - \frac{\tilde{f}}{2}\Tr(\tilde{F}H)\\
\pa_v \Tr(H^2) &= -\frac{\eta^2\tilde{f}}{2}\Tr(H^2) - \frac{\eta h}{2}\Tr(\tilde{F}H).\end{split}
\end{align}

Since $\tilde{F},H$ are $\slash{k}$-symmetric, in fact by the Cauchy-Schwarz inequality
\[|\Tr(\tilde{F}H)| = |\Tr(\tilde{F}^\ast H)| \leq \frac{1}{2}|\tilde{F}|^2_{\slash{k}} + \frac{1}{2}|H|^2_{\slash{k}} = \frac{1}{2}(\Tr(\tilde{F}^2) + \Tr(H^2)).\]

Thus we derive the differential inequalities
\begin{align}
\label{eq:FH:eq4}
\begin{split}
\pa_\eta \Tr(\tilde{F}^2) &\leq \left(\frac{1}{4}\tilde{f}-\frac{1}{2}h\right)\Tr(\tilde{F}^2) + \frac{\tilde{f}}{2}\Tr(H^2)\\
\pa_v \Tr(H^2) &\leq \left(\frac{1}{4}\eta h - \frac{1}{2}\eta^2\tilde{f}\right)\Tr(H^2) + \frac{\eta h}{2}\Tr(\tilde{F}^2).\end{split}
\end{align}

By \cref{thm:C4:fh}, on $[0,v_0]\times[0,\eta_0]\times U$, the coefficients are all uniformly bounded by some $C$. Since $\Tr(\tilde{F}^2), \Tr(H^2) \geq 0$ (since $\tilde{F}$ and $H$ are $\slash{k}$-symmetric), we deduce that there exists some $B$ depending only on the data for $\tilde{F}$ and $H$ and a constant $C$ so that
\begin{align*}
\Tr(\tilde{F}^2)(v,\eta) &\leq B + C\int_0^\eta \Tr(\tilde{F}^2)(v,t) + \Tr(H^2)(v,t)\ dt\\
\Tr(H^2)(v,\eta) &\leq B + C\int_0^v \Tr(\tilde{F}^2)(s,\eta) + \Tr(H^2)(s,\eta)\ ds.
\end{align*} Now we may use a result from \cref{C:A3:Goursat}. Namely, we use a two-dimensional version of Gronwall's inequality, \cref{thm:A3:Gronwall}, to immediately conclude \eqref{eq:FH:eq2}.

Next let us show that $\slash{k}$ and $\mathring{\slash{g}}$ are uniformly equivalent. Let $X \in TS^2$ be a fixed vector. We first show that $\slash{k}(X,X) \lesssim \mathring{\slash{g}}(X,X)$ where the implied constant does not depend on $X$. Let us differentiate
\begin{align*}
\pa_v (\slash{k}(X,X)) &= (\pa_v \slash{k})(X,X)\\
&= (\slash{k}(f/2+F))(X,X) = \slash{k}(X,(f/2+F)X).\end{align*}
The last equality follows from, denoting $A = (f/2+F)$, \begin{align*}(\slash{k}A)(X,X) &= (\slash{k}A)_{ab}X^aX^b\\
&= \slash{k}_{ac}A^{c}_bX^a X^b\\
&= \slash{k}_{ac}X^a (AX)^c = \slash{k}(X,AX).
\end{align*}
Now notice that the operator norm of $F$ with respect to $\slash{k}$ is at most $|F|_{\slash{k}}$,\footnote{This is just the fact that on any Hilbert space, the operator norm is at most the Hilbert-Schmidt norm.} and so by \eqref{eq:FH:eq2}, the operator norm is bounded above. By \cref{thm:C4:fh}, we know that $|f|$ is bounded above on $[0,v_0]\times[0,\eta_0]\times U$, and so we conclude that
\[\pa_v \slash{k}(X,X) \lesssim \slash{k}(X,X)\] where the constant depends only on $\mathcal R_0$. Gronwall's inequality now implies that
\[\slash{k}(X,X) \lesssim \mathring{\slash{g}}(X,X),\] since $\slash{k} = \mathring{\slash{g}}$ on $\{v = 0\}$. Again the constant depends only on $\mathcal R_0$. This is one half of equivalence. To show the other, we need to put a lower bound on the lowest eigenvalue of $K = \mathring{\slash{g}}^{\vspace{-5pt}-1}\slash{k}$, which is a $\mathring{\slash{g}}$-symmetric matrix. The upper bound we have just proven gives an upper bound on its largest eigenvalue, so to obtain a lower bound on its lowest eigenvalue, we need only lower bound its determinant. But this is true by \cref{thm:C4:www} since $\digamma^2 > 0$ and is continuous on $\mathcal R_0$.\end{proof}

We observe the following important remark, which we will extend in the next chapter.
\begin{rk}$\digamma f$ and $\digamma h$ continue as smooth functions to all of $\mathbf{if}$.\end{rk}

One has a similar, albeit much weaker, result about $F$ and $H$. As in the proof of \cref{thm:C4:FH}, denote, for a $(1,1)$ tensor $A$, $|A|_{\slash{k}} = \sqrt{\Tr(A^\ast A)}$, where the adjoint is taken with respect to $\slash{k}$.

\begin{prop}\label{thm:C4:nonblowup}Fix $v_0 > 0$ and let $U$ be an open subset of $S^2$ on which $\inf_{\theta \in U} E(v_0) > 0$. Then $\digamma|\tilde{F}|_{\slash{k}}$ and $\digamma|H|_{\slash{k}}$ are uniformly bounded on compact subsets of $\{\digamma \geq 0\}\n \{v \geq v_0, \theta \in U\}$, even as $\digamma \to 0$.\end{prop}
In particular, the eigenvalues of $\tilde{F}$ and $H$ grow at most like $\frac{1}{\digamma}$ as $\digamma \to 0$.

As we do not need this proposition to prove either \cref{thm:C4:continuation} or \cref{thm:C1:blowupvague}, we delay its proof until the end of \cref{C:C6}, since the techniques to prove it will be similar to those developed in \cref{C:C5}, but will require some techniques from \cref{C:C6}, as well.

\subsection{The behaviour of \texorpdfstring{$\omega$}{omega} to top order.}
Let us finally analyze \eqref{eq:C4:three} at $\mathbf{if}$ to obtain the top-order behaviour of $\slash{\omega}$. We will prove:
\begin{prop}\label{thm:C4:topomega}Suppose $\phi \in \phgd(\Mp)$. Then there exists a unique smooth $\omega$, defined on $\mathbf{if}\n\{\digamma > 0\}$ such that $\phi$ satisfies \eqref{eq:C4:one}$=0$ on $\{\digamma > 0\}$ and satisfies the compatibility conditions with the series solutions at $\mathbf{bf}$ and $\mathbf{lf}$.\end{prop}
\begin{proof}
By \cref{thm:C4:topV} and \cref{thm:C4:topk}, we know $V=0$ and the uniqueness of $\slash{k}$. With $f$, $h$, $F$, $H$ as in the previous subsection,
\eqref{eq:C4:three} becomes
\begin{align*}0 &= \eta\pa_\eta \pa_v \omega + \frac{1}{4}\eta\pa_\eta f -\frac{1}{4} \frac{\eta}{16}fh + \frac{\eta}{8}\Tr(FH) - \frac{1}{4}f.\end{align*}
Using \eqref{eq:C4:fh}, one can substitute $\eta\pa_\eta f$ and rewrite this as
\begin{align}\label{eq:topomega:eq1}0 &= \pa_\eta \pa_v \omega - \frac{1}{16}fh + \frac{1}{8}\Tr(FH) + \frac{1}{4\eta}f.\end{align}
Recall that $f/\eta \in C^{\infty}(\{\digamma > 0\})$, so the last three terms are smooth in $\{\digamma > 0\}$. Either by integrating out both derivatives or using \cref{thm:A3:globalsmooth}, it follows that for given data for $\omega$ on $\{v = 0\}\un \{\eta = 0\}$ there exists a unique smooth solution on $\{\digamma > 0\}$. From \cref{thm:C4:first}, in order to satisfy the compatibility conditions, $\omega$ is $0$ to top order at $\mathbf{lf}$ and $\mathbf{bf}$, so this shows uniqueness and the existence of a solution. To show that this solution satisfies all the compatibility conditions, let $\widebar{\phi} = (\widebar{\slash{k}},\widebar{V},\widebar{\omega})$ be a series solution at $\mathbf{bf}$ and $\mathbf{lf}$, provided by \cref{thm:C4:first}. By \cref{thm:C4:topk} and \cref{thm:C4:topV}, $\slash{k}$ and $\widebar{\slash{k}}$, $V$ and $\widebar{V}$, respectively, are the same modulo a rapidly vanishing error, and so $\widebar{\omega}$ satisfies \eqref{eq:topomega:eq1} modulo a rapidly vanishing error. Since $\widebar{\omega} = \omega = 0$ on $\{v = 0\}\un \{\eta = 0\}$, either by integrating or using a similar argument to the proof of \cref{thm:C4:topkcompat}, it follows that $\omega = \widebar{\omega}$ modulo a rapidly vanishing error. This is sufficient to show that the compatibility conditions hold.
\end{proof}
\chapter{Behaviour at \texorpdfstring{$\digamma = 0$}{digamma = 0} for generic commutative data}
\label{C:C5}
\section{Preliminaries}
\label{C:C5:intro}
Let us continue to fix short-pulse data, i.e.\ the short-pulse tensor $\mathbf{T}$, and continue to denote its energy by $\mathbf{E}(v) = \frac{1}{2}\int_0^v |\pa_s \mathbf{T}|^2(s)\ ds$. Let $g$ be the formal series solution with short-pulse data given by $\mathbf{T}$, which is provided by \cref{thm:C4:continuation}, and let $\slash{k} = \xi^{-4}\slash{g}$ be the associated rescaled fibre metric, and $\omega = \log\Omega$.
In \cref{C:C4:topo}, we began to examine the behaviour of $\omega$ and the derivatives of $\slash{k}$ on $\mathbf{if}\n\{\digamma > 0\}$. Keeping the notation from that subsection, let us denote
\[\pa_v \slash{k} = \slash{k}(f/2 + F), \ \pa_\eta \slash{k} = \slash{k}(h/2 + H),\] where $f$, $h$ are scalars, and $F$, $H$ are tracefree. In this chapter we seek to obtain more detailed information about how $f$, $h$, $F$, $H$ and $\omega$ behave near the boundary of the region $\{\digamma > 0\}$. While this is easy for $f$, $h$, since they are given by explicit formulae in \cref{thm:C4:fh} it will be more difficult for $F$, $H$ and $\omega$.\footnote{Since $V = \xi^3\eta\overline{L} \equiv 0$, we already know what happens to it, so we will ignore it in this chapter.}

We will only be able to obtain our most detailed information about $F$, $H$ and $\omega$ under more restrictive assumptions on $\mathbf{T}$, which we call \emph{commutative}.
\begin{defn}We will call short-pulse data \emph{commutative} if $\mathbf{T}$ factors as $\psi(v)\mathbf{T}_0(\theta)$, for a real-valued smooth function $\psi$ on $[0,1]$ with $\psi(0) = 0$, not depending on $\theta \in S^2$, and a fixed $\mathring{\slash{g}}$-symmetric tracefree $(1,1)$ tensor $\mathbf{T}_0$ on $S^2$. We will also call short-pulse tensors admitting such a factorization \emph{commutative}.\end{defn}
The upshot of commutative data, and the reason for its name, comes from a simplification to the equation for $F$ and $H$, \eqref{eq:C4:FH}. Since the data for $H$ is just $0$, and the data for $F/\eta$ is $\pa_v \mathbf{T} = \pa_v \psi \mathbf{T}_0$, in the case of commutative data, we may write $F = F'\mathbf{T}_0$ and $H = H'\mathbf{T}_0$, where now $(F',H')$ satisfy the equation
\begin{align}
\label{eq:C5:FH}
\begin{split}0 &= \pa_\eta F' + \left(\frac{1}{4}h-\frac{1}{\eta}\right)F' + \frac{1}{4}fH'\\
0 &= \pa_v H' + \left(\frac{1}{4}h-\frac{1}{\eta}\right)F' + \frac{1}{4}fH',\end{split}
\end{align}
with initial data $H'|_{\mathbf{lf}} = 0$, $(F'/\eta)|_{\mathbf{bf}} = \pa_v \psi$. This is the same equation as \eqref{eq:C4:FH}, except the commutator drops out (hence the term ``commutative'') and the resulting equation is linear. We also recall \eqref{eq:topomega:eq1} for $\omega$:
\begin{align}\label{eq:C5:omega}0 &= \pa_\eta \pa_v \omega - \frac{1}{16}fh + \frac{1}{8}F'H'\Tr(\mathbf{T}_0)^2 + \frac{1}{4\eta}f\end{align}
with data $\omega|_{\mathbf{lf}} = 0$, $\omega_{\mathbf{bf}} = 0$.
\begin{rk}Due to topological reasons, any tracefree $(1,1)$ tensor $\mathbf{T}_0$ must vanish somewhere. In particular the energy $\mathbf{E}$ associated to a commutative short-pulse tensor must vanish on $[0,1]\times \{\theta_0\}$ for at least one $\theta_0$. Thus, even though we will prove a curvature blow-up theorem for such data (\cref{thm:C1:blowupvague}), Christodoulou does not (and consequently we do not) prove the formation of trapped surfaces for such data in \cref{thm:C1:ChrVague}. However, the work of Klainerman--Luk--Rodnianski \cite{KlaLukRodFull} deals with the anisotropic formation of trapped surfaces, so we may appeal to the theorem proved there to show that commutative data also form trapped surfaces, even though it is not the fibred spheres which become trapped.\end{rk}

In order to undertake the most refined analysis, we will also need to impose two more restrictions: one on $\psi$ and one on $\mathbf{T}_0$. For $\psi$, we will require that $0$ is a \emph{simple} zero of $\psi$, i.e.\ $\pa_v \psi(0) \neq 0$. Since $\mathbf{T}_0$ must vanish somewhere, we will impose the condition that all zeroes of $\mathbf{T}_0$ are also \emph{simple}. This means that in any local trivialization $\R^2 \to U \subseteq S^2$ of the bundle of tracefree $(1,1)$ tensors over $S^2$ (a rank $2$ bundle), if we consider $\mathbf{T}_0: U \to \R^2$ as a smooth map, then the differential $(d\mathbf{T}_0)_p$ must be non-singular at any point $p$ where $\mathbf{T}_0(p) = 0$. An equivalent way of stating this is that $\mathbf{T}_0$ is transverse to the zero section of the bundle of tracefree $(1,1)$ tensors.

\begin{rk}\label{rk:C5:nongeneric}While we make these two further restrictions, all the work in this chapter will also apply to general commutative data, provided we restrict ourselves to examining regions away from $\{v = 0, \eta = \infty\}$ and the zeroes of $\mathbf{T}_0$. Indeed, for the second restriction on the region, the equations are only parametrized by $\theta \in S^2$, so we can always restrict to an open set disjoint from the zeroes of $\mathbf{T}_0$. For the first restriction on the region, a backwards domain of dependence of $\{\digamma = 0, \ \eta < \infty\}$, only intersects $\{v = 0\}$ for finite values of $\eta$.\end{rk}

Even if we do not restrict our attention to the regions mentioned in the previous remark, the restrictions on $\psi$ and $\mathbf{T}_0$ are not substantial. Indeed,
\begin{lem}The set of commutative short-pulse tensors $\psi\mathbf{T}_0$ with $0$ a simple $0$ of $\psi$ and all zeroes of $\mathbf{T}_0$ simple is \emph{generic} in the set of all commutative short-pulse tensors, in the sense that such tensors form an open dense set (in the $C^\infty$ topology).\end{lem}
\begin{proof}
Let us first show openness. Fix a commutative short-pulse tensor $\T = \psi \T_0$ satisfying the restrictions. Let us assume that $\pa_v \psi(0) = 1$. Write $\psi(v) = v\tilde{\psi}(v)$, for $\tilde{\psi} \in C^{\infty}([0,1])$ and $\tilde{\psi}(0) \neq 0$. It suffices to show that if 
\[\mathbf{S} = \chi \mathbf{S}_0 \] is a a factorization of any other commutative short-pulse tensor, into a smooth function $\chi \in C^{\infty}([0,1])$ with $\chi(0) = 0$ and a tracefree tensor $\mathbf{S}_0$, then $\tilde{\chi} = \chi/v$ is not $0$ at $0$ and $\mathbf{S}_0$ is transverse to the zero section. Fix $\theta_0 \in S^2$ for which $\mathbf{T}_0(\theta_0) \neq 0$. If $\mathbf{T}$ and $\mathbf{S}$ are sufficiently close in the $C^\infty$ topology, then $\tilde{\psi}(0)\mathbf{T}_0(\theta_0)$ and $\tilde{\chi}(0)\mathbf{S}_0(\theta_0)$ are sufficiently close tensors. In particular, $\tilde{\chi}(0) \neq 0$. Dividing $\tilde{\chi}$ by a constant and multiplying $\mathbf{S}_0$ by the same constant, we may assume that $\tilde{\chi}(0) = \tilde{\psi}(0)$. Thus, $\mathbf{T}_0$ and $\mathbf{S}_0$ are also close in the $C^\infty$ topology. Since the set of tensors with simple zeroes is certainly open, it follows that if $\mathbf{T}$ and $\mathbf{S}$ are sufficiently close, then $\mathbf{S}_0$ only has simple zeroes.

Next, let us show density. If $\mathbf{T} = \psi\mathbf{T}_0$ is an arbitrary commutative short-pulse tensor, then certainly $\psi + \epsilon v$ has simple zeroes for $\epsilon$ small and converges to $\psi$ in the $C^\infty$ topology as $\epsilon \to 0$. Thus to show density, it suffices to show that the set of tensors transverse to the zero section is dense. But this is a consequence of the parametric transversality theorem.\footnote{The use of the parametric transversality theorem in this argument is essentially the same as its use in proving the density of Morse functions.}
\end{proof}

In order to state the theorem giving the precise behaviours of $f$, $h$, $F'$, $H'$ and $\omega$ at $\{\digamma = 0\}$, it will be necessary to desingularize $\mathbf{if}\n \{\digamma = 0\}$ via real blowup. In fact, we will also be able to examine the corner $(v=0,\eta = \infty)$, which means that we will not be able to only consider $\mathbf{if}$ as $[0,\infty)_v \times [0,\infty)_\eta \times S^2_{\theta}$, as we did in \cref{C:C4:existence}. Rather, we will also need to restrict the coordinates $(\varpi,\tau,v,\theta)$, which are valid near $\mathbf{if}\n \mathbf{sf}$, to $\mathbf{if}$ ($ = \{\varpi = 0\}$), and will consider $\mathbf{if}$ as
\[([0,1]_v \times [0,\infty)_\eta \times S^2_{\theta})\sqcup([0,1]_v \times [0,\infty)_\tau \times S^2_{\theta})/\sim,\]
where $\sim$ is the relation $\eta = 1/\tau$. We will perform the blowups and state the theorem in the next section. In the subsequent sections, we prove the theorem, devoting a section for different regions of the blowup space.

\section{Blowups and the statement of the main theorem}
\label{C:C5:blowup}
\begin{figure}
\centering
\begin{subfigure}{.5\textwidth}
 \centering
 \includegraphics{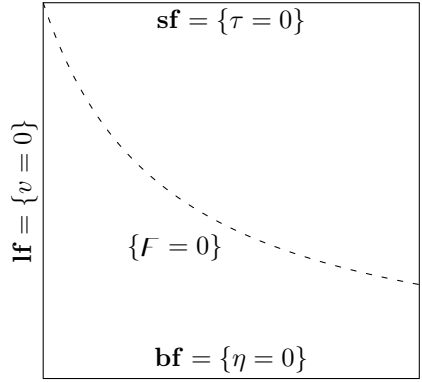}
 \caption{A view of $\mathbf{if}$ at an angle $\theta$ for which $\mathbf{T}_0(\theta) \neq 0$.}
\end{subfigure}%
\begin{subfigure}{.5\textwidth}
 \centering
 \includegraphics{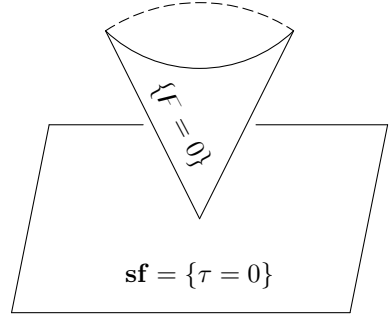}
 \caption{An upside-down view of $\mathbf{if}$ at some $v > 0$ near a $\theta$ for which $\mathbf{T}_0(\theta) = 0$.}
\end{subfigure}
\caption{The reason for blowing up. The closure of $\{\digamma = 0\}$ intersects the boundary at the corner $\mathbf{lf} \n \mathbf{sf}$ in $\{\mathbf{T}_0 > 0\}$, and meets $\mathbf{sf}$ in a cone point at $\mathbf{T}_0(\theta) = 0$.}
\end{figure}
We perform the blowups necessary to give the precise behaviours of $f$, $h$, $F'$, $H'$. For the rest of this chapter, let $\mathbf{lf}$, $\mathbf{sf}$, $\mathbf{bf}$ denote their intersections with $\mathbf{if}$, rather than the entire face in $\widebar{\mathcal M}$.

Let us blow up the sphere $S = \{v = \tau = 0\} = \mathbf{sf} \n \mathbf{lf}$, forming $M_1 = [\mathbf{if},S]$. This introduces a front face $\mathbf{frf}$, which separates $\mathbf{lf}$ from $\mathbf{sf}$. Two sets of projective coordinates will be useful for us to cover this blowup: $(\lambda = v/\tau, \tau,\theta)$ and $(v,\zeta = \tau/v,\theta)$, where $\theta$ parameterizes some coordinate chart on $S^2$. The first is valid near $\mathbf{lf}$, and the second near $\mathbf{sf}$. Because $\mathbf{T}_0$ has zeroes, this is not quite enough. We will need to blow up the zero set $\nu_2 = \mathbf{sf} \n \{\mathbf{T}_0 = 0\}$, which is diffeomorphic to a finite disjoint union of intervals $[0,1]$. So let us set $M_2 = [M_1,\nu_2]$. This introduces several new faces $\mathbf{nf}_i \iso S^2_+ \times [0,1]$ ($S^2_+ \iso D^2$ denoting the upper hemisphere of $S^2$), one for each zero of $\mathbf{T}_0$, which intersect the lift of $\mathbf{sf}$ transversely.

In order to simplify the computations, we will use special coordinates near each $\mathbf{nf}_i$. The manifold $\nu_2$ may be described entirely in $(v,\zeta,\theta)$ coordinates. We choose special coordinates $\theta$ around the zeroes of $\mathbf{T}_0$ in order to choose coordinates on $M_2$. In any stereographic chart, $\mathbf{T}_0$ takes the form
\[\mathbf{T}_0 = \begin{pmatrix} a & b\\ b & -a \end{pmatrix}\] for smooth functions $a,b \: U \subseteq \R^2 \to \R$. Let us assume that $z$ is a zero of $\mathbf{T}_0$ and that $z = 0$ in $U$. Since the zeroes of $\mathbf{T}_0$ are simple, we may change coordinates again so that $a(\theta^1,\theta^2) = \sqrt{2}\theta^1$, and $b(\theta^1,\theta^2) = \sqrt{2}\theta^2$ (the normaliztion will be made for convenience, below).\footnote{Observe carefully it is \emph{not} true in these coordinates that $\mathbf{T}_0 = \begin{psmallmatrix} \sqrt{2}\theta^1 & \sqrt{2}\theta^2\\
\sqrt{2}\theta^2 & -\sqrt{2}\theta^1\end{psmallmatrix}$. Indeed, changing coordinates would involve conjugating the matrix by the Jacobian. Instead we have identified the rank two bundle of tracefree symmetric tensors with $\R^2$ and have trivialized it as an abstract bundle.} The advantage is that in these coordinates, $\Tr(\mathbf{T}_0^2) = 2|\theta|^2 = 2(\theta^1)^2 + 2(\theta^2)^2$. Now, choose coordinates $(v, \zeta, \alpha^1 = \theta^1/\zeta, \alpha^2 = \theta^2/\zeta)$ near $\mathbf{nf}_i$. In these coordinates, $\mathbf{nf}_i = \{\zeta = 0\}$. One could also introduce coordinates which are valid near the lift of $\mathbf{sf}$, but we will not need this.

We will check below that $\mathbf{sgf}$ (for \emph{singular face}), the lift of the set $\{\digamma = 0\}$ intersects the boundary faces of $M_2$ transversely, and thus we set $M = M_2 \n \overline{\{\digamma \geq 0\}}$. We will also check that $\mathbf{sf} \n M = \emptyset$, which is why in the previous paragraph we didn't introduce coordinates near $\mathbf{sf}$.

\begin{figure}[htbp]
\centering
\includegraphics{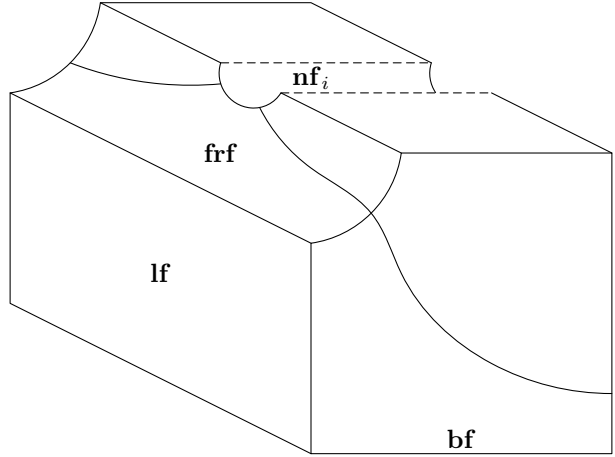}
\caption{A stylized view of the boundary faces of $M$. The directions in $S^2$ are pictured as one dimensional and going into the page. The curved lines in the middle of the figure represent the intersection of the closure of $\{\digamma = 0\}$ with the boundary faces of $M_2$. $\{\digamma = 0\}$ itself continues as a surface inside the interior. $M$ is the portion of $M_2$ below this surface.}
\end{figure}

Let us verify the transversality of (the closure of the lift of) $\{\digamma = 0\}$ to the boundary faces of $M_2$. Is is convenient to introduce the notation
\begin{align*}
R^2 &:= \frac{1}{2}\Tr(T_0^2)\\
\frac{1}{2}\int_0^v \mathbf{E}(t)\ dt &= \frac{1}{2}R^2\int_0^v\int_0^s |\pa_t \psi|^2(t)\ dtds =: R^2 J.
\end{align*}
Observe that $0 \leq J \in v^2 C^{\infty}([0,1])$, $(J/v^2)(0) =: \kappa^2 > 0$. Thus we may write $\sqrt{J} = \kappa v(1+vK)$ where $\kappa > 0$, $K$ is smooth, and $1+vK > 0$. With this notation
\begin{gather*}
  \frac{1}{2}\int_0^v \mathbf{E}(t)\ dt =\kappa^2 v^2 (1+vK)^2 R^2\\
  \mathbf{E} = 4\kappa^2 v (1+vK)(1+2vK + v^2\pa_v K)R^2.
\end{gather*}

In $(v,\zeta,\theta)$ coordinates, this mean that
\[\digamma = 4(1-\kappa^2\zeta^{-2}(1+vK)^2R^2).\] Thus, away from the zeroes of $\mathbf{T}_0$ (which are the zeroes of $R$) the closure of the lift of $\{\digamma = 0\}$ is
\[\{\zeta - \kappa (1+vK)R = 0\}.\] This is transverse to $\mathbf{frf}$, and the closure of $\{\digamma \geq 0\}$ does not intersect $\{\zeta = 0\} = \mathbf{sf}$. In $(v,\zeta,\alpha)$ coordinates,
\[\digamma = 4(1-\kappa^2(1+vK)^2|\alpha|^2),\] and so the closure of the lift is $\{\digamma = 0\}$ is
\[\{1-\kappa (1+vK)|\alpha| = 0\}\] Since $|\alpha|$ is smooth away from $\alpha = 0$ (which is certainly disjoint from the $\{\digamma = 0\}$), this intersects $\mathbf{nf}_i$ transversely for each $i$, and the closure of $\{\digamma \geq 0\}$ does not intersect $\mathbf{sf}$, which is at $|\alpha| = \infty$ in these coordinates.

$M$ has a number of boundary faces, which we order $\mathbf{lf}$, $\mathbf{bf}$, $\mathbf{frf}$, $\mathbf{sgf}$ followed by the faces $\mathbf{nf_i}$. We continue to regard $\mathbf{ff} = \{v = 1\}$ as artificial.

We are now almost in a position to state the theorem about the behaviours on $M$ of $f$, $h$, $F'$, $H'$ and $\omega$.
For an integer $n$, let $n'$ denote the index set \[n' = n \un \{n+1,n+2,\ldots\}\times \{1\},\] i.e.\ the index set consisting of those functions which have expansions in integral powers, but may have up to one logarithm appearing in lower-order terms, and \[n'' = n \un \{n,n+1,\ldots,\}\times \{1\} \un \{n+1,n+2,\ldots\}\times \{2\},\] i.e.\ the index set consisting of those functions which have expansions in integral powers, but may have a logarithm appearing in the highest power, and up to two logarithms in lower-order terms.

\begin{thm}\label{thm:C5:behaviour}$f$, $h$, $F'$, $H'$ and $\omega$ are polyhomogeneous with the following index families:
\begin{align*}
\mathcal E_f &= (1,2,-1,-1,0,\ldots,0)\\
\mathcal E_h &= (2,1,1,-1,1,\ldots,1)\\
\mathcal E_{F'} &= (0,1,-1,-1',-1,\ldots,-1)\\
\mathcal E_{H'} &= (1,0,1,-1',0,\ldots,0)\\
\mathcal E_{\omega} &= (2,2,0,0'',0,\ldots,0).\end{align*}
In other words,
\[f \in \phg{\mathcal E_f}(M), \ h \in \phg{\mathcal E_h}(M), \ F' \in \phg{\mathcal E_{F'}}(M), \ H' \in \phg{\mathcal E_{H'}}(M), \ \omega \in \phg{\mathcal E_\omega}(M).\]
\end{thm}

For most of the proof, it will be convenient to instead work with $\tilde{h} = -\eta^2 h$ and $\tilde{H} = -\eta^2 H'$. This corresponds to setting
\[\pa_\tau \slash{k} = -\eta^2 \pa_\eta \slash{k} = \slash{k}(\tilde{h}/2 + \tilde{H}\mathbf{T}_0).\]
Since $\eta^2 = \tau^{-2} = \zeta^{-2}v^{-2}$,
\begin{align*}
\mathcal E_{\tilde{h}} = (2,3,-1,-1,-1,\ldots,-1)\\
\mathcal E_{\tilde{H}} = (1,2,-1,-1',-2,\ldots,-2),
\end{align*}
it is equivalent to show
\[\tilde{h} \in \phg{\mathcal E_{\tilde{h}}}(M), \ \tilde{H} \in \phg{\mathcal E_{\tilde{H}}}(M).\]

The statements for $f$, $\tilde{h}$, are easy to prove since we have explicit formulae for them from \cref{thm:C4:fh}.
\begin{proof}[Proof of \cref{thm:C5:behaviour} for $f$, $h$] We cover $M$ by several coordinate charts $U$ and show that $f$, $\tilde{h}$ are polyhomogeneous at the faces which lie inside $U$ with the correct index families, respectively. We will write down the formulae for $f$, $\tilde{h}$ derived from \cref{thm:C4:fh}.

Near $\mathbf{bf}$, where coordinates $(v,\eta,\theta)$ are valid,
\begin{align}
\label{eq:C5:fhve}
\begin{split}
f &= -\frac{4\eta^2 v\kappa^2(1+vK)[1+2vK + v^2\pa_v K]R^2}{1-v^2\eta^2\kappa^2(1+vK)^2R^2}\\
\tilde{h} &= \frac{\eta^3 v^2\kappa^2(1+vK)^2R^2}{1-v^2\eta^2\kappa^2(1+vK)^2R^2}.\end{split}\end{align}Since $v$ is a bdf of $\mathbf{lf}$ and $\eta$ is a bdf of $\mathbf{bf}$, the desired polyhomogeneity follows.

Let us now turn to $(\lambda,\tau,\theta)$ coordinates. In these coordinates, we may use the explicit formulae to write
\begin{align}
\label{eq:C5:fhlt}
\begin{split}
f &= -\frac{4\kappa^2 \lambda(1+\tau\lambda K)[1+2\tau\lambda K+\tau\lambda^2\pa_\lambda K]R^2}{\tau(1-\lambda^2 \kappa^2(1+\tau\lambda K)^2R^2)}\\
\tilde{h} &= \frac{\kappa^2\lambda^2(1+\tau\lambda K)^2R^2}{\tau(1-\lambda^2 \kappa^2(1+\tau\lambda K)^2R^2)}.
\end{split}
\end{align}
Since $\lambda$ is a bdf of $\mathbf{lf}$ and $\tau$ is a bdf of $\mathbf{frf}$, the desired polyhomogeneity is clear. In $(v,\zeta,\theta)$ coordinates, $f$, $\tilde{h}$ are given by
\begin{align}
\label{eq:C5:fhvz}
\begin{split}
f &= -\frac{4\kappa^2(1+vK)[1+2vK+v^2\pa_v K]R^2}{v(\zeta+\kappa(1+vK)R)(\zeta-\kappa(1+vK)R)}\\
 \tilde{h} &= \frac{\kappa^2(1+vK)^2R^2}{\zeta v(\zeta+\kappa(1+vK)R)(\zeta-\kappa(1+vK)R)},
\end{split}
\end{align}
 and $v$ is a bdf of $\mathbf{frf}$ and $\zeta - \kappa(1+vK)R$ is a bdf of $\mathbf{sgf}$. The desired polyhomogeneity is again clear. Finally, in $(v,\zeta,\alpha)$ coordinates around the $i$th zero of $\mathbf{T}_0$,
\begin{align}
\label{eq:C5:fhalpha}
\begin{split}
f &= -\frac{4\kappa^2(1+vK)[1+2vK+v^2\pa_v K]|\alpha|^2}{v(1+\kappa(1+vK)|\alpha|)(1-\kappa(1+vK)|\alpha|)}\\
\tilde{h} &= \frac{\kappa^2(1+vK)^2|\alpha|^2}{\zeta v(1+\kappa(1+vK)|\alpha|)(1-\kappa(1+vK)|\alpha|)},
\end{split}
\end{align}
and $v$ is a bdf of $\mathbf{frf}$, $1-\kappa(1+vK)|\alpha|$ is a bdf of $\mathbf{sgf}$ and $\zeta$ is a bdf of $\mathbf{nf}_i$. The desired polyhomogeneity is clear.
\end{proof}
\begin{rk}
We point out that the formulae \eqref{eq:C5:fhlt}-\eqref{eq:C5:fhalpha} are also important to prove the remainder of \cref{thm:C5:behaviour}.\end{rk}

Dealing with $F'$, $\tilde{H}$ and $\omega$ is more difficult, although treating $\omega$ will be significantly easier than dealing with $F'$ and $\tilde{H}$. We split the proof into four propositions, the first of which follows immediately from the results of \cref{C:C4:topo}, and the remaining three we split over the next three sections. Each proposition will consider polyhomogeneity in regions where the coordinates introduced above are valid: $(v,\eta,\theta)$, $(\lamda,\tau,\theta)$, $(v,\zeta,\theta)$, $(v,\zeta,\alpha)$, respectively (the spherical variables act like parameters in the first three charts so we mostly ignore them). We start with a proposition about the behaviour away from the new faces.
\begin{prop}\label{thm:C5:prop1} It holds that $F',\tilde{H},\omega \in C^{\infty}(M^\circ)$. Moreover, setting $U_1 = \{\digamma > 0, \eta < \infty\} \subseteq M$, 
\begin{romanumerate}
\item $F' \in \phgi{(\mathbf{lf},\mathbf{bf})}{(0,1)}(U_1),\\
\tilde{H} \in \phgi{(\mathbf{lf},\mathbf{bf})}{(1,2)}(U_1);$
\item $\omega \in \phgi{(\mathbf{lf},\mathbf{bf})}{(2,2)}(U_1)$.
\end{romanumerate}
\end{prop}
\begin{proof}
The smoothness in $M^\circ$ follows from \cref{thm:C4:FH} for $F'$ and $\tilde{H}$, and \cref{thm:C4:topomega}. Now we show that the index families are as described. From the specification of the initial data, we know $F'|_{\{\eta = 0\}} = 0$, so this shows that $F'$ has the correct index family. $\tilde{H} = \eta^2 H' \in \eta^2 C^{\infty}(U_1)$ and $H'|_{\{v = 0\}} = 0$ by assumption, so this shows that $\tilde{H}$ has the correct index family.

From \eqref{eq:C5:omega} and the above,
\[\pa_\eta \pa_v \omega \in \phgi{(\mathbf{lf},\mathbf{bf})}{(1,1)}(U_1).\] Since $\omega = 0$ on $\{v = 0\} \un \{\eta = 0\}$, it follows that $\omega$ has the desired index family.
\end{proof}

Observe that there exists $\lambda_0 > 0$ sufficiently small so that
\[U_2 = \{0 \leq \lambda = v/\tau = v\eta \leq \lambda_0, \ 0 \leq \tau \leq 1\} \subseteq M\] is disjoint from $\mathbf{sgf}$ and each $\mathbf{nf}_i$. Indeed, if $\lambda_0$ is sufficiently small, then $U_2 \subseteq M_2$ is disjoint from $\mathbf{sf}$, so cannot intersect any $\mathbf{nf}_i$, and $\mathbf{sgf}$ does not intersect $\mathbf{lf}$, so since $\lambda$ is a bdf of $\mathbf{lf}$, $U_2$ cannot intersect $\mathbf{sgf}$ either. With this definition of $U_2$, we have:

\begin{prop}\label{thm:C5:prop2}The following are true: 
\begin{romanumerate}
\item $F' \in \phgi{(\mathbf{lf},\mathbf{frf})}{(0,-1)}(U_2),\\
\tilde{H} \in \phgi{(\mathbf{lf},\mathbf{frf})}{(1,-1)}(U_2);$
\item $\omega \in \phgi{(\mathbf{lf},\mathbf{frf})}{(2,0)}(U_2).$
\end{romanumerate}
\end{prop}

Choose any $\zeta_0 > 1/\lambda_0$. Let $W \subseteq S^2$ be any open set whose closure is disjoint from the zeroes of $\mathbf{T}_0$. Then the region 
\[U_3 = \{\zeta = \tau/v = 1/(v\eta) \leq \zeta_0, \ \theta \in W, \zeta - \kappa(1+vK) R \geq 0\} \subseteq M\] contains part of $\mathbf{frf}$, corresponding to $\{v = 0\}$, and part of $\mathbf{sgf}$, corresponding to $\{\zeta - \kappa(1+vK) R = 0\}$. Then:
\begin{prop}\label{thm:C5:prop3} The following are true in $U_3$:
\begin{romanumerate}
\item $F' \in \phgi{(\mathbf{frf},\mathbf{sgf})}{(-1,-1')}(U_3),\\
\tilde{H} \in \phgi{(\mathbf{frf},\mathbf{sgf})}{(-1,-1')}(U_3);$
\item $\omega \in \phgi{(\mathbf{frf},\mathbf{sgf})}{(0,0'')}(U_3).$
\end{romanumerate}
\end{prop}

Finally, we work in a neighbourhood near the zeroes of $\mathbf{T}_0$. Let $(v,\zeta,\alpha)$ be coordinates near some face $\mathbf{nf}_i$. These are valid provided $|\alpha|\zeta = |\theta| \leq \delta_i$ for $\delta_i$ sufficiently small (recall that a zero of $\mathbf{T}_0$ is $\{\theta = 0\}$ in these coordinates) and $1-\kappa(1+vK)|\alpha| \geq 0$, the latter being a bdf of $\mathbf{sgf}$. With $\zeta_0$ as above, the region
\[U_4 = \{\zeta \leq \zeta_0, |\alpha|\zeta \leq \delta_i, 1-\kappa(1+vK)|\alpha| \geq 0\} \subseteq M\]
contains part of $\mathbf{frf}$, corresponding to $\{v = 0\}$, part of $\mathbf{sgf}$, corresponding to $\{1-\kappa(1+vK)|\alpha| = 0\}$, and part of $\mathbf{nf}_i$ corresponding to $\{\zeta = 0\}$.

\begin{prop}\label{thm:C5:prop4}The following are true in $U_4$:
\begin{romanumerate}
\item $F' \in \phgi{(\mathbf{frf},\mathbf{sgf},\mathbf{nf}_i)}{(-1,-1',-1)}(U_4),\\
\tilde{H} \in \phgi{(\mathbf{frf},\mathbf{sgf},\mathbf{nf}_i)}{(-1,-1',-2)}(U_4);$
\item $\omega \in \phgi{(\mathbf{frf},\mathbf{sgf},\mathbf{nf}_i)}{(0,0'',0)}(U_4).$
\end{romanumerate}
\end{prop}

Since $M$ is covered by the union of $U_1$, $U_2$, $U_3$ and all possible $U_4$ for each face $\mathbf{nf}_i$, these propositions together prove \cref{thm:C5:behaviour}.

The general outline for proving \cref{thm:C5:prop2}, \cref{thm:C5:prop3}, and \cref{thm:C5:prop4} is the same, although the details will differ (sometimes by a little and other times significantly). The outline is as follows:
\begin{romanumerate}
\item Establish weighted $L^\infty$ bounds.
\item Commute with b-vector fields and establish conormality, i.e.\ show that $F'$, $\tilde{H}$ and $\omega$ in $\mathcal A$.
\item Upgrade conormality to full polyhomogeneity.
\end{romanumerate}

\subsection{Specialized notation and polyhomogeneity}
Before moving on, we introduce some specialized notation and state some theorems which will be useful during the rest of the chapter. Let $X$ be a mwc with a tuple $\mathcal F$ of $n$ boundary faces, and let $\rho = (\rho_1,\ldots, \rho_n)$ be a tuple of bdfs. 

For $m \in \N_0$ let $\Diff_b^m$ denote the set of $\mathrm{b}$-differential operators of order $m$, i.e.\ those which admit a local description in a chart $U \subseteq (\R^+_x)^k\times \R_y^{j}$:
\[\sum_{|\beta| \leq N} a_\beta (x_1\pa_{x_1})^{\beta_1}\cdots (x_k \pa_{x_k})^{\beta_k}(\pa_{y_1})^{\beta_{k+1}}\cdots (\pa_{y_j})^{\beta_{k+j}},\] where $a_\beta$ are smooth functions.

For $\alpha \in \R^n$ and $N \in \N$, we define the $L^\infty$ $b$-Sobolev spaces
\[W_b^{N,\alpha}(X) = \{u \: X \to \R \text{ measurable }\: L\rho^{-\alpha} \in L^{\infty}(X),\text{ for all } L \in \Diff^N_b(X)\}.\]
Of course,
\[W_{b}^{N,\alpha}(X) = \rho^{\alpha}W_b^{N}(X).\]
Observe that $\bigcap_N W_b^N(X) \subseteq \mathcal A(X)$, with equality if $X$ is compact.\footnote{In this section, we will only work over compact spaces so that this equality holds. One can of course define $W_{b,\loc}^N(X)$ and the intersection of these will be $\mathcal A(X)$, regardless of compactness.}

We will often deal with the Cauchy problem, so for some proofs it will be convenient at first only to measure regularity tangent to the level sets of sets ``parallel'' to Cauchy hypersurfaces surfaces. For this reason, if $x$ is a smooth function with everywhere non-singular derivative, we define the space
\[W_{b,x}^{N,\alpha}(X) = \{u \: X \to \R \text{ measurable }\: L\rho^{-\alpha} \in L^{\infty}(X),\text{ for all } L \in \Diff^N_{b,x}(X)\},\]
where $\Diff^N_{b,x}(X) \subseteq \Diff^N_{b}(X)$ is the subset of those $\mathrm{b}$-differential operators which only consist of derivatives tangent to the level sets of $X$. If $X$ admits a product decomposition $[0,1)_x\times Y$, then
\[W_{b,x}^{N,\alpha}(X) = L^{\infty}_xW_b^{N,\alpha}(Y).\] The advantage of the space $W_{b,x}^{N,\alpha}$ is that we do not need to choose an explicit product decomposition.

We now give a useful equivalent description of polyhomogeneity, which can be found in \cite{MelCalc}.

Let $E$ be an index set, and define the polynomial for $N \in \R$
\[b(E,N;t) = \prod_{(z,p) \in E, \ \Re(z) < N} (t-z)\]
(notice that the multiplicity of the root $z$ is equal to $q+1$ for $q$ the largest integer such that $(z,q) \in E)$, and $b(E,N;t) \equiv 1$ if $\Re(z) \geq N$ for all $(z,p) \in E$.)

For each face $F_i \in \mathcal F$, we may choose a ``radial'' vector field $\nu_i$, i.e.\ a vector field tangent to the faces of $X$ such that $\nu_i/\rho_i$ is transverse to $F_i$. In local coordinates $(x,y)$ in which $x_1$ is a bdf of $F_i$, $\nu_i = a(x,y)(x_1\pa_{x_1})$, where $a$ is smooth and nonvanishing on $\{x_1 = 0\}$.

We have:
\begin{thm}[\cite{MelCalc}]\label{thm:C5:melrose}Suppose $\mathcal E = (E_1,\ldots,E_n)$ is an index family associated to $\mathcal F$. Let $(\nu_1,\ldots,\nu_n)$ be some choice of radial vector fields for the faces $F_i$. Then $u \in \phgd(X)$ if and only if for all tuples $\alpha = (N_1, \ldots, N_n)$
\[\prod_{i=1}^N b(E_i,N_i;\nu_i)u \in \rho^{\alpha}\mathcal A(X).\]\end{thm}

Because of the exact index sets we will use, we give some of these polynomials special names. Define the sequences of polynomials $a_k(t)$ and $p_k(t)$ for $k = 0,1,2,3,\ldots$ by
\begin{align*}
1, t+1, (t+1)(t), (t+1)(t)(t-1), \ldots\\
1, t, t(t-1), t(t-1)t,-2), \ldots.
\end{align*}
The $a_k$, $p_k$ should be thought of as the polynomials $b$ associated to the index sets $-1$ and $0$, respectively. Indeed, if $N \leq -1$ or $N \leq 0$, respectively, $b(-1,N;t) = a_0(t)$, $b(0,N;t) = p_0(t)$, respectively, and if $j < N \leq j+1$, for $j \geq -1$ or $j \geq 0$, respectively, then $b(-1,N;t) = a_{j+2}(t)$, $b(0,N;t) = p_{j+1}(t)$.

We record an elementary but useful observation about applying the product rule with these polynomials.
\begin{lem}\label{thm:C5:product}Suppose $A,B \in C^{\infty}((0,\infty))$. Then for all $j \in \R$, and $k,\ell \in \R$ with $k+\ell = j$,
\[(x\pa_x-j)(AB) = ((x\pa_x-k)A)B + A((x\pa_x-\ell)B).\]

In particular, for all $m$, there are constants $c_{m,i,1}$ and $c_{m,i,2}$ such that
\begin{align*}
a_m(x\pa_x)(AB) &= \sum_{i=0}^m c_{m,i,1}p_i(A)a_{m-i}(B)\\
p_m(x\pa_x)(AB) &= \sum_{i=0}^m c_{m,i,2}p_i(A)p_{m-i}(B).
\end{align*}
The same statements remains true if $A$, $B$ depend on some parameters.
\end{lem}

\todo{check reference in texorpdfstring}\section{The proof of \texorpdfstring{\cref{thm:C5:prop2}}{proposition~5.2.4}}
We use coordinates $(\lambda,\tau,\theta)$, for $\theta$ in some compact neighbourhood, to cover $U_2$. Let $U \subseteq U_2$ denote an arbitrary such coordinate chart, which is compact. For the rest of this section, we work only in $U$. In these coordinates, \eqref{eq:C5:FH} becomes
\begin{align}
\label{eq:C5:FHlt}
\begin{split}
(\tau\pa_\tau - \lambda\pa_\lambda)F' + \frac{\tau f}{4}\tilde{H} + \frac{\tau \tilde{h}}{4}F' + F' &= 0\\
\pa_\lambda \tilde{H} + \frac{\tau f}{4}\tilde{H} + \frac{\tau \tilde{h}}{4}F' +F'&= 0.\end{split}\end{align}

and \eqref{eq:C5:omega} becomes
\begin{equation}\label{eq:C5:omegalt}(\tau\pa_\tau -\lambda\pa_\lambda)\pa_\lambda \omega + \frac{\tau^2}{16}f\tilde{h} - \frac{\tau^2}{4}F'\tilde{H}R^2 - \frac{\tau}{4}f = 0\end{equation}

From \cref{thm:C5:prop1}, we know smoothness of $F'$, $\tilde{H}$ and $\omega$ away from the boundary $\{\tau = 0\}$ (=$\{\eta = \infty\})$. Thus we will treat \eqref{eq:C5:FHlt} as an initial value problem with data for $F'$ posed on $\{\tau = 1\}$, and data for $\tilde{H}$ given by $\tilde{H}|_{\{\lambda = 0\}} \equiv 0$ coming from the requirement $H|_{\mathbf{lf}} \equiv 0$, and \eqref{eq:C5:omegalt} as an initial-value problem with data for $\omega$ posed on $\{\tau = 1\}$ and given on $\{\lambda = 0\}$ by $\omega \equiv 0$, coming form the requirement that $\omega|_{\mathbf{lf}} = 0$. The spherical variables act like a parameter, so we will usually ignore them unless we are trying to establish regularity in the spherical directions.\begin{figure}[htbp]
\centering
\includegraphics{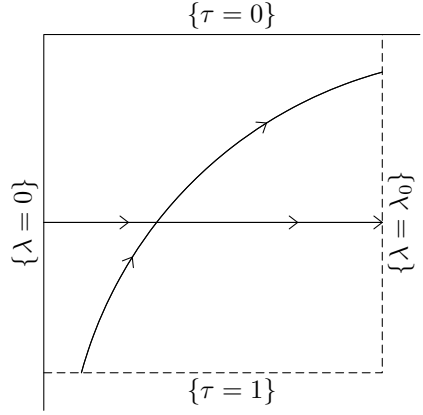}
\caption{A view of $U = U_2$ at some $\theta$. Indicated is a sample forwards-directed integral curve of $\pa_\lambda$ and a backwards-directed integral curve of $\tau\pa_\tau - \lambda\pa_\lambda$.}
\end{figure}

Let us first work on proving part (i) of \cref{thm:C5:prop2}.
We will need to to expand the type of equations we consider to include an inhomogeneous right-hand side. Such equations will arise when commuting \eqref{eq:C5:FHlt} with derivatives which we will use to establish regularity. Let us consider for $j,\ k\in \N_0$ the system
\begin{align}
\label{eq:C5:meqlt}
\begin{split}
(\tau\pa_\tau - \lambda\pa_\lambda)A + a_{11}B +a_{12}A + A = R_A\\
\pa_\lambda B + a_{21}B +a_{22}A +A = R_B,\end{split}\end{align}
where $a_{i1} \in \lambda L^{\infty}(U)$, $a_{i2} \in \lambda^2 L^{\infty}(U)$, $R_A \in \lambda^{k+\max(j,1)}\tau^{k-1}L^{\infty}(U)$, $R_B \in \lambda^{k+j}\tau^{k-1}L^{\infty}(U)$.
and $A|_{\{\tau = 1\}} \in \lambda^{k+j}L^{\infty}(\{\tau = 1\})$, and $B|_{\{\lambda = 0\}} = 0$.
Set
\begin{align*}
\mathcal E &= \sum_{i=1}^2 \norm{a_{i1}}_{\lambda L^{\infty}(U)} + \norm{a_{i2}}_{\lambda^2 L^{\infty}(U)}\\
\mathcal I &= \norm{A|_{\{\tau = 1\}}}_{\lambda^{k+j}L^{\infty}(\{\tau = 1\})}+\norm{R_A}_{\lambda^{k+\max(j,1)}\tau^{k-1}L^{\infty}(U)} + \norm{R_B}_{\lambda^{k+j}\tau^{k-1}L^{\infty}(U)}.
\end{align*}
Then:
\begin{prop}\label{thm:C5:initialreglt}Suppose $A,B \in C^{1}(\{\tau > 0\})$ solves \eqref{eq:C5:meqlt}.\hspace{-2pt} Then $A \in \lambda^{k+j}\tau^{k-1}L^{\infty}(U)$, $B \in \lambda^{k+1+j}\tau^{k-1}L^{\infty}(U)$, and 
\[\norm{A}_{\lambda^{k+j}\tau^{k-1}L^{\infty}(U)} + \norm{B}_{\lambda^{k+1+j}\tau^{k-1}L^{\infty}(U)} \leq C\mathcal I,\]
where $C = C(\mathcal E,\lambda_0,j,k)$ is increasing in all its arguments.
\end{prop}
\begin{rk}
The reason for $\max(1,j)$ appearing is technical. If it were just $j$, then if $j=0$ one obtains a logarithmic loss in the estimate. This is analogous to how the ODE $(\tau\pa_\tau + 1)A = \tau^{-1-j}$ has a solution in $O(\tau^{-1-j})$ except if $j = 0$. Since in our applications either $R_A \equiv 0$ or $j>0$, this will not be a cause for concern.\end{rk}
\begin{proof}We will let $C$ denote a constant depending only on $\lambda_0$, $\mathcal E$, and explicit numerical constants, which may vary from line to line.
On a first reading, it is helpful to suppose $k=j = 0$, and $R_A = R_B = 0$. 

We will recast the system as an integral system.
Define $\phi$, $\psi$ by requiring them to solve in $U$
\begin{align*}
(\tau\pa_\tau -\lambda\pa_\lambda)\phi &= \phi\left(a_{12}+1\right)\\
\pa_\lambda \psi &= \psi a_{21},\end{align*}
with $\phi(\lambda,1) \equiv 1$, $\psi(0,\tau) \equiv 1$. Then

 \begin{align*}
 \phi(\lambda,\tau) &= \tau\exp\left(\int_1^{\tau} a_{12}(\lambda\tau/t,t)\ dt/t\right)\\
 \psi(\lambda,\tau) &= \exp\left(\int_0^\lambda a_{21}(s,\tau)\ ds\right).\end{align*} 
 
 Since
 \begin{align*}
 \left|\int_1^{\tau} a_{12}(\lambda\tau/t,t)\ dt/t\right| &\leq \lambda^2\mathcal E \leq C\mathcal E\\
 \left|\int_0^{\lambda} a_{21}(s,\tau)\ ds\right| &\leq \lambda^2 \mathcal E \leq C\mathcal E\end{align*}
 we have
 \begin{align*}
(\phi/\tau), (\tau/\phi) &\leq e^{C\mathcal E}\\
\psi, \psi^{-1} &\leq e^{C\mathcal E}.\end{align*}
 
 Now using the method of integrating factors,
 \begin{align}
 \label{eq:C5:bootstraplt}
\begin{split}
A(\lambda,\tau) &= \phi^{-1}\left(A(\lambda \tau,1) - \int_{1}^{\tau} \left(\phi a_{11}B\right)(\lambda\tau/t,t)\ dt/t + \int_1^{\tau} \left(\phi R_A\right)(\lambda\tau/t,t)\ dt/t\right)\\
B(\lambda,\tau) &= \psi^{-1}\left(- \int_0^\lambda \left(\psi\left(a_{22}A + A\right)\right)(s,\tau)\ ds + \int_0^{\lambda} (\psi R_B)(s,\tau)\ ds\right).\end{split}\end{align}

Observe that the integrals are taken along the flow of $\tau\pa_\tau - \lambda\pa_\lambda$ and $\pa_\lambda$, respectively.

We will use a bootstrap argument, and find constants $K_1, K_2, N_1$, depending on $\lambda_0$ and $\mathcal E$ (which we will choose later) such that for any $0 < \tau^\ast \leq 1$ there exists $N_2$ such that for any $\epsilon > 0$ such that for $\tau^\ast \leq \tau \leq 1$,
\begin{align*}
|A| &\leq \epsilon \mathcal I e^{N_2(1+\lambda-\tau)} + K_1\mathcal I e^{N_1 \lambda^2}\lambda^j(\lambda \tau)^{k}\tau^{-1}\\
|B| &\leq\epsilon \mathcal I e^{N_2(1+\lambda-\tau)}+ K_1K_2\mathcal I e^{N_1\lambda^2}\lambda^j(\lambda\tau)^{k}\lambda\tau^{-1}.\end{align*}
(observe that $e^{N_2(1-\tau)}$ is increasing in the direction in which we are propagating).
Since $K_1,K_2, N_1$ do not depend on $\tau^\ast$ or $\epsilon$, taking $\epsilon \to 0$, following by $\tau^\ast \to 0$ shows the desired bound. We need only check that the first bound is satisfied for $\tau$ close to $1$, the second bound is satisfied for $\lambda$ close to $0$, and if $0 < \lambda_1 \leq \lambda_0$ and $\tau^\ast \leq \tau_1 < 1$ and the bound is satisfied for $\tau_1 \leq \tau < 1$ and $0 < \lambda_1 \leq \lambda$, then $A$ and $B$ in fact satisfy an improved bound on the region.
 
By continuity, for any $\epsilon > 0$ the appropriate bounds are satisfied near the initial surface. Observe that if $(\lambda,\tau) \in \{0 < \lambda \leq \lambda_1, \ \tau_1 \leq \tau < 1\}$, then the arguments of the integrands of \eqref{eq:C5:bootstraplt} lie in the same region. Therefore we may plug the bootstrap assumption into \eqref{eq:C5:bootstraplt} (and use that $\tau/t \leq 1$ if $\tau \leq t$) to see that if $0 < \lambda \leq \lambda_1$ and $\tau_1 \leq \tau < 1$, then
\begin{align*}
|A|(\lambda,\tau) &\leq C\tau^{-1}|A|(\lambda\tau, 1) + \epsilon \mathcal I C\tau^{-1}\int_{\tau}^1 \lambda e^{N_2(\lambda\tau/t +1- t)}\ dt/t\\
&+ C K_1K_2 \mathcal I \tau^{-1}\int_{\tau}^1t(\lambda\tau/t)e^{N_1 (\lambda\tau/t)^2}(\lambda\tau/t)^j (\lambda\tau)^k (\lambda\tau/t) t^{-1}\ dt/t\\
&+C \mathcal I \tau^{-1}\int_{\tau}^1 t(\tau\lambda)^{k+\max(1,j)}t^{-1-\max(1,j)}\ dt/t\\
&\leq C\lambda^{k+j}\tau^{k+j-1}\mathcal I + \epsilon\tau^{-1} \mathcal I C \int_{\tau}^1 e^{N_2(1-t)}\ dt\\
&+ C K_1K_2\mathcal I\lambda^{j+k}\tau^{k-1} \int_{\tau}^1 (\lambda\tau)^2 t^{-3}e^{N_1(\lambda\tau/t)^2}\ dt\\
&+ C \mathcal I \lambda^{k+\max(1,j)}\tau^{k+\max(1,j)-1}\int_{\tau}^1 t^{-1-\max(1,j)}\ dt\\
&\leq \frac{C\epsilon\mathcal I}{N_2\tau^\ast}e^{N_2(1-\tau)}+ C\mathcal I \lambda^{k+j}\tau^{k-1}\left(\tau^j - \frac{K_1 K_2}{2N_1}+ \frac{K_1K_2}{2N_1}e^{N_1\lambda^2} + \frac{\lambda^{\max(1,j)-j}}{\max(1,j)}\right)\\
&\leq \frac{C\epsilon\mathcal I }{N_2\tau^\ast}e^{N_2(1-\tau)}+ C\mathcal I \lambda^{k+j}\tau^{k-1}\left(C- \frac{K_1 K_2}{2N_1}+ \frac{K_1K_2}{2N_1}e^{N_1\lambda^2}\right)\\
\end{align*}
and
\begin{align*}
|B|(\lambda,\tau) &\leq C\epsilon \mathcal I \int_0^{\lambda}e^{N_2(s+1-\tau)}\ ds\\
&+ C\mathcal I K_1 \int_0^{\lambda} \tau^{k-1} s^{k+j}e^{N_1s^2} \ ds+ C\mathcal I \int_0^{\lambda} s^{k+j}\tau^{k-1}\ ds\\
&\leq \frac{C\epsilon \mathcal I}{N_2} e^{N_2(\lambda+1-\tau)} + C\mathcal I\lambda^{k+j+1}\tau^{k-1}(K_1e^{N_1\lambda^2}+1)\\
&\leq \frac{C\epsilon \mathcal I}{N_2} e^{N_2(\lambda+1-\tau)} + C\mathcal I\lambda^{k+j+1}\tau^{k-1}e^{N_1\lambda^2}K_1\end{align*}

Replacing $C$ with $\max(C,1)$, we will improve over the bootstrap assumption provided:
\begin{romanumerate}
\item $N_2 > C(\tau^\ast)^{-1}$;
\item $K_2 > \max(C,2)$;
\item $N_1 > CK_2$;
\item $K_1 > \max(2N_1C/K_2,2)$.
\end{romanumerate}
The first condition deals with the first term of each bound. The second condition bounds $CK_1 < K_1K_2$. The third bounds $C\frac{K_1K_2}{2N_1} < K_1$, and the fourth bounds $C-\frac{K_1K_2}{2N_1} < 0$.
\end{proof}

With additional regularity assumptions, we have additional regularity of the solution.
\begin{cor}\label{thm:C5:corlt}Let $a_{ij}$ ($i,j=1,2$), $R_A$, $R_B$ be as in \eqref{eq:C5:meqlt}. Fix $N \geq 1$ and suppose additionally that $a_{i1} \in \lambda W_b^N(U)$, $a_{i2} \in \lambda^2 W_b^N(U)$, $R_A \in \lambda^{k+\max(1,j)}\tau^{k-1}W_b^N (U)$, $R_B \in \lambda^{k+j}\tau^{k-1}W_b^N(U)$. If $A,B \in C^{N+1}(\{\tau > 0\})$, solve \eqref{eq:C5:meqlt} with data $A|_{\{\tau = 1\}} \in \lambda^{k+j}W_b^N(\{\tau = 1\})$ and $B|_{\{\lambda = 0\}} = 0$, then $A \in \lambda^{k+j}\tau^{k-1}W_b^N (U)$, $B \in \lambda^{k+1+j}\tau^{k-1}W_b^N(U)$, with analogous bounds to those in \cref{thm:C5:initialreglt}.

In particular, if the assumptions are satisfied for all $N$ (i.e.\ with $\mathcal A(U)$ instead of $W_b^N(U)$), then $A \in \lambda^{k+j}\tau^{k-1}\mathcal A (U)$, $B \in \lambda^{k+1+j}\tau^{k-1}\mathcal A(U)$.\end{cor}
\begin{proof}

We first show $A \in \lambda^{k+j}\tau^{k-1}W_{b,\tau}^{N'}(U), B \in \lambda^{k+j+1}\tau^{k-1}W_{b,\tau}^{N'}U)$, for $N' \leq N$, by induction on $N'$. The case $N' = 0$ is the previous proposition.

Now let us show is for $N'$, assuming it is true for $N'' < N'$.

For all $P \in \Diff^{N'}_{b,\tau}(U)$ by definition $PB|_{\{\lambda = 0\}} = 0$,\footnote{Recall that the derivatives are $\lambda\pa_\lambda$, not $\pa_\lambda$.} and by assumption \[PA \in \lambda^{k+j} L^{\infty}(\{\tau = 1\}).\]

Choosing $P$ as a product of coordinate vector fields $\lambda\pa_\lambda,\pa_{\theta^i}$, and commuting $P$ with \eqref{eq:C5:meqlt}, we see that $PA,PB$ satisfies and equation of the form \eqref{eq:C5:meqlt}, except with some additional terms on the right-hand side of the form
\[(Q_1 a)(Q_2 S),\] where $a \in \lambda W_b^N(U)$, $S$ is $A$ or $B$, and $Q_1 \in \Diff^{L_1}_{b,\tau}(U)$, $Q_2 \in \Diff^{L_2}_{b,\tau}(U)$ for $L_1 + L_2 = N$, $L_2 < N$. The first factor is in $\lambda L^{\infty}(U)$ by assumption, and the second is in $\lambda^{k+j}\tau^{k-1} L^{\infty}(U)$, by the inductive hypothesis. 

Thus their product is in $\lambda^{k+\max(j,1)}\tau^{k-1} L^{\infty}(U)$, and we may apply \cref{thm:C5:initialreglt}. This completes the inductive step.

Now we complete the proof by showing by induction on $i$ that \[(\tau\pa_\tau)^i A \in \lambda^{k+j}\tau^{k-1}W_{b,\tau}^{N-i}(U), \ (\tau\pa_\tau)^i B \in \lambda^{k+j+1}\tau^{k-1}W_{b,\tau}^{N-i}(U).\] The base case $i = 0$ was treated above. Now let us show it for $i$, assuming it is true for $i' < i$. Commuting $(\tau\pa_\tau)^i$ with \eqref{eq:C5:meqlt}, we see that $(\tau\pa_\tau)^i A, (\tau\pa_\tau)^i B$ satisfies and equation of the form \eqref{eq:C5:meqlt}, except with additional terms on the right hand side of the from
\[[(\tau\pa_\tau)^{i-i'}a][(\tau\pa_\tau)^{i'} S],\] where $a \in \lambda W_b^N(U)$, $S$ is $A$ or $B$, and $0 \leq i' < i$. By the inductive hypothesis, \[(\tau\pa_\tau)^{i'}S \in \lambda^{k+\max(1,j)}\tau^{k-1}W_{b,\tau}^{N-i'}(U) \subseteq  \lambda^{k+\max(1,j)}\tau^{k-1}W_{b,\tau}^{N-i}(U).\]Now we may argue by induction on the amount of regularity, as in the base case and conclude, provided that the data is in the right space, i.e.\ $(\tau\pa_\tau)^i A$ is in $\lambda^{k+j}W_b^{N-i}(\{\tau = 1\})$ and $(\tau\pa_\tau)^i B|_{\{\lambda = 0\}} = 0$. The second of these is true because $\tau\pa_\tau$ is tangent to $\{\lambda = 0\}$. To show the first, use the equation for $(\tau\pa_\tau)^{i-1} A$ to express $(\tau\pa_\tau)^i A$ in terms of lower-order $\tau$ derivatives of $A$ and $B$, to see that \[(\tau\pa_\tau)^i A \in \lambda^{k+\max(j,1)}W_b^{N-i}(\{\tau = 1\}),\] whici is the correct space.\end{proof}

We need one more lemma:
\begin{lem}\label{thm:C5:finiteseriesI}For $k \in \N_0$,
$a_k(\tau\pa_\tau)F'|_{\{\tau = 1\}} \in \lambda^{k}C^{\infty}(\{\tau = 1\})$.
\end{lem}
\begin{proof}
It suffices to prove that in a neighbourhood of $\{\tau = 1\}$ $F^j = \pa_\lambda^j F'|_{\{\lambda = 0\}}$ and $\tilde{H}^j = \pa_{\lambda}^j \tilde{H}|_{\{\lambda = 0\}}$ have finite expansions in integral powers $\tau^\ell$, \emph{without remainder},
\begin{align*}
F^j &= \sum_{k=-1}^{j-1} F^j_k\tau^k\\
\tilde{H}^j &= \sum_{k=-1}^{j-2} \tilde{H}^j_k \tau^k,\end{align*}
where $F^j_k,\tilde{H}^j_k \in C^{\infty}(\{\lambda = 0, \ \tau = 1\})$
(except for $j=0$ for $\tilde{H}^j$, which we know by assumption is $\equiv 0$). Indeed, this implies that $a_k(\tau\pa_\tau)\pa_{\lambda}^j F'|_{\{\lambda = 0\}} = 0$ for $k \geq j+1$, and since we know $F'$ is smooth for $\tau > 0$, it follows that $a_k(\tau\pa_\tau)F'|_{\{\tau = 1\}}$ is smooth in $\lambda$ and has its first $k-1$ $\pa_\lambda$-derivatives vanishing at $\lambda = 0$.

Since we know smoothness for $\tau > 0$, $F^j$, $\tilde{H}^j$ are well-defined, and commuting \eqref{eq:C5:fhlt} with $\pa_\lambda^j$ shows that they satisfy
\begin{subequations}
\begin{align}
\label{eq:C5:eq11a}
\tau\pa_\tau F^j + (1-j)F^j + \sum_{\ell=1}^j {\binom{j}{\ell}}\left[ \pa_\lambda^{\ell}\left(\frac{\tau}{4}f\right)\tilde{H}^{j-\ell} + \pa_\lambda^{\ell}\left(\frac{\tau}{4}\tilde{h}\right)F^{j-\ell}\right] &= 0\\
\label{eq:C5:eq11b}
\tilde{H}^j + F^{j-1}+ \sum_{\ell=1}^j {\binom{j}{\ell}}\left[ \pa_\lambda^{\ell}\left(\frac{\tau}{4}f\right)\tilde{H}^{j-\ell} + \pa_\lambda^{\ell}\left(\frac{\tau}{4}\tilde{h}\right)F^{j-\ell}\right] &= 0,\end{align}
\end{subequations}
except if $j=0$, in which case there is no equation available for $\tilde{H}^0$.
We will show the desired expansions by induction on $j$.
Observe that by \eqref{eq:C5:fhlt} $\tau f$ is $\lambda$ times a smooth function of $\lambda, \lambda\tau$, $\theta$. Indeed, $K$ is a function of $\theta$ and $v = \tau\lambda$, and this is the only dependence on $\tau$. Therefore, at $\lambda = 0$, $\pa_{\lambda}^{\ell}(\tau f)$ has a finite expansion in powers of $\tau$ with powers $\tau^k$, $0 \leq k \leq \ell-1$ and coefficients in $C^{\infty}$. Similarly, $\pa_{\lambda}^{\ell}(\tau \tilde{h})$, has an expansion at $\{\lambda = 0\}$, except with powers $\tau^k$, $0 \leq k \leq \ell-2$.

For $j=0$, \eqref{eq:C5:eq11a} reads $\tau\pa\tau F^0 + F^0 = 0$, which means $F^0 = \tau^{-1}F^0|_{\tau = 1}$, which proves the base case for $F^0$. Now for higher $j$, by induction the right-hand side of \eqref{eq:C5:eq11a} has the desired series expansion, i.e.\ has a finite expansion in powers $\tau^k$ for $-1 \leq k \leq j-2$.
Equation \eqref{eq:C5:eq11b} is just an identity and by induction we have the desired expansion. Equation \eqref{eq:C5:eq11a} is of the form
\[(\tau\pa_\tau + (1-j))F^j = R,\] where $R$ is a series of the desired form, with largest power of $\tau^k$, that with $k = j-2$. Thus we may solve and obtain
\[F^j(\tau) = \tau^{j-1}F^j(1) + \tau^{j-1}\int_{1}^{\tau} t^{1-j}R dt/t,\] which we may integrate to show that $F^j$ has a representation as a series of the desired form (no $\log$s appear because the largest power of $k$ in the integrand is $t^{-2}$.
\end{proof}

We may now prove part (i) of \cref{thm:C5:prop2}.
\begin{proof}[Proof of \cref{thm:C5:prop2}:(i)]It suffices to show that for all $i,k \in \N_0$
\begin{align*}
p_i(\lambda\pa_\lambda)a_k(\tau\pa_\tau) F' \in \lambda^{i}\tau^{k-1}\mathcal A(U), \ p_i(\lambda\pa_\lambda)a_k(\tau\pa_\tau)\tilde{H} \in \lambda^{i+1}\tau^{k-1}\mathcal A(U)
\end{align*}
While this apparently only proves that $\tilde{H}$ is smooth at $\{\lambda = 0\}$, the assumption that $\tilde{H}|_{\{\lambda = 0\}} = 0$ means it actually has index set $1$. We will in fact prove the stronger claim that for $i \leq k$
\begin{align*}
p_i(\lambda\pa_\lambda)a_k(\tau\pa_\tau) F' \in \lambda^{k}\tau^{k-1}\mathcal A(U),\ p_i(\lambda\pa_\lambda)a_k(\tau\pa_\tau)\tilde{H} \in \lambda^{k+1}\tau^{k-1}\mathcal A(U)
\end{align*}
and for $i \geq k+1$, setting $i = k+j$ that
\begin{align*}
p_i(\lambda\pa_\lambda)a_k(\tau\pa_\tau) F' \in \lambda^{k+j}\tau^{k-1}\mathcal A(U), \ p_i(\lambda\pa_\lambda)a_k(\tau\pa_\tau)\tilde{H} \in \lambda^{k+1+j}\tau^{k-1}\mathcal A(U).\end{align*}

Let us prove this by induction on $i+k$. If $i+k = 0$, then the initial data falls under the assumptions of \cref{thm:C5:corlt}, and so we may conclude.
Now if we commute $p_i(\lambda\pa_\lambda)a_k(\tau\pa_\tau)$ with \eqref{eq:C5:FHlt} and use \cref{thm:C5:product}, we obtain that $p_i(\lambda\pa_\lambda)a_k(\tau\pa_\tau) F'$ and $p_i(\lambda\pa_\lambda)a_k(\tau\pa_\tau)\tilde{H}$ satisfy \eqref{eq:C5:FHlt}, except there is a right-hand side consisting of sums of terms of the form
\begin{gather}\label{eq:C5:types}\begin{split}[p_{i-i'}(\lambda\pa_\lambda)p_{k-k'}(\tau\pa_\tau)(\tau f)][p_{i'}(\lambda\pa_\lambda)a_{k'}(\tau\pa_\tau)(\tilde{H})],\\
[p_{i-i'}(\lambda\pa_\lambda)p_{k-k'}(\tau\pa_\tau)(\tau \tilde{h})][p_{i'}(\lambda\pa_\lambda)a_{k'}(\tau\pa_\tau)(F')]\end{split}\end{gather}
for $0 \leq i'+k' < i+k$.
Now observe
$p_i(\lambda\pa_\lambda)a_k(\tau\pa_\tau) H'|_{\{\lambda = 0\}} = 0$ and by \cref{thm:C5:finiteseriesI}
$a_k(\tau\pa_\tau) F' \in \lambda^{k}C^{\infty}(\{\tau = 1\})$ and so
\[p_i(\lambda\pa_\lambda)a_k(\tau\pa_\tau) F' \in \lambda^{\max(k,i)}C^{\infty}(\{\tau = 1\}).\]
This means, treating the cases $i \leq k$ and $i \geq k$ separately, that the initial data verify the hypotheses necessary to apply \cref{thm:C5:corlt}.
In order to verify the remaining hypotheses, and conclude the inductive step, we just need both types of terms in \eqref{eq:C5:types} to be in $\lambda^{\max(k+1,i)}\tau^{k-1}\mathcal A(U)$.

Let us look at the first type term. The second factor is by induction in \[\tau^{k'-1}\lambda^{\max(k',i')}\mathcal A(U).\] By \eqref{eq:C5:fhlt}, $\tau f$ is not just a smooth function of $\lambda,\tau$; it is in fact $\lambda$ times a smooth function of $\lambda$ and $\tau\lambda$. Thus $p_{k-k'}(\tau\pa_\tau)(\tau f) \in \lambda^{k-k'+1}\tau^{k-k'}C^{\infty}(U)$. And so the first factor is in
\[\tau^{k-k'}\lambda^{\max(i-i',k'-k'+1)}C^{\infty}(U),\] and thus the first type of term is in
\[\tau^{k-1}\lambda^{\max(i-i',k-k'+1)+\max(i',k')}\mathcal A(U).\] Since $\tau\tilde{h}$ is $\lambda^2$ times a smooth function of $\lambda\tau$, we similarly have that the second type of term is in
\[\tau^{k-1}\lambda^{\max(i-i',k-k'+2)+\max(i',k')}\mathcal A(U).\]

The exponent of $\tau$ matches, so we just need to ensure that the exponent of $\lambda$ is at least as big as $\max(i,k+1)$. Set $x = i'-k'$, $y = i-(k+1)$. Use the identity $2\max(a,b) = a + b + |a-b|$, for $a,b\in \R$ to rewrite $2\max(i,k+1) = i+k+1+|y|$, and (twice) the exponents of the $\lambda$ in the first and second type of terms as
\[i + k + 2 + |y-x| + |x|, \ i+k+2 + |(y-1)-x| + |x|,\] respectively. Using the reverse triangle inequality establishes that the exponents are big enough, which allows us to conclude the proof.\end{proof}

We now turn our attention to proving part (ii) of \cref{thm:C5:prop2}. Let us consider for $j,k \in \R$ the equations:
\begin{equation}\label{eq:C5:wlti} \pa_\lambda w = S\end{equation}
and
\begin{equation}\label{eq:C5:wltii}(\tau\pa_\tau - \lambda\pa_\lambda) w = S,\end{equation}
where in both cases $S \in \lambda^j \tau^k L^{\infty}(U)$.
We start with a general proposition.
\begin{prop}\label{thm:C5:omegalthelper}Suppose $w \in C^1(\{\tau > 0\})$ solves \eqref{eq:C5:wlti} with data $w|_{\{\lambda = 0\}} = 0$. Then $w \in \lambda^{j+1}\tau^k L^{\infty}(U)$ with the bound
\[\norm{w}_{\lambda^{j+1}\tau^k L^{\infty}(U)} \leq C\norm{S}_{\lambda^{j}\tau^k L^{\infty}(U)},\]
where $C = C(j,k,\lambda_0)$ is increasing in all its arguments.

Now fix $N \geq 1$. If instead $w \in C^{N+1}(\{\tau > 0\})$ and $S \in \lambda^j\tau^k W_b^N(U)$, then $w \in \lambda^{j+1}\tau^{k} W_b^N(U)$ with an analogous bound. In particular, if the hypotheses hold for all $N$, then $w \in \lambda^{j+1}\tau^k \mathcal A(U)$.

If $w$ solves \eqref{eq:C5:wltii} with data $w|_{\{\tau = 1\}} \in \lambda^{j}L^{\infty}(U)$ and $k\neq j$, then $w \in \lambda^{j}\tau^{\min(j,k)} L^{\infty}(U)$ with the bound
\[\norm{w}_{\lambda^{j}\tau^{\min(j,k)} L^{\infty}(U_2)} \leq C\left(\norm{w}_{\lambda^j L^{\infty}(\{\tau = 1\})} + \norm{S}_{\lambda^{j}\tau^kL^{\infty}(U_2)}\right),\]
where $C = C(j,k,\lambda_0)$ is increasing in all its arguments. If instead $S \in \lambda^j\tau^k W_b^N(U)$, then $w \in \lambda^j\tau^{\min(j,k)} W_b^N(U)$ with analogous bounds. 

In particular, if the hypotheses hold for all $N$, then $w \in \lambda^{j}\tau^{\min(j,k)} \mathcal A(U)$.\end{prop}
\begin{proof}
Suppose $w$ solves \eqref{eq:C5:wlti} with the specified data. The first statement is obvious by integrating. Since $\tau\pa_\tau, \lambda\pa_\lambda, \pa_{\theta^i}$ ($i = 1,2$) commute through the equation the second statement is also true.

Now suppose $w$ solves \eqref{eq:C5:wltii} with the specified data. Then $w$ admits the representation formula
\[w(\lambda,\tau) = w(\lambda\tau,1) + \int_{1}^\tau S(\lambda\tau/t,t)\ dt/t.\]
Since $k-j \neq 0$, no $\log$s appear when naively bounding the integral, and the first statement follows. Since $\tau\pa_\tau, \lambda\pa_\lambda, \pa_{\theta^i}$ ($i = 1,2$) all commute through the equation the second statement is also clear.
\end{proof}

We need one more lemma:
\begin{lem}\label{thm:C5:finiteseriesomegaI}For $k \in \N_0$, $p_k(\tau\pa_\tau)\omega|_{\{\tau = 1\}} \in \lambda^{k+2}C^{\infty}(\{\tau = 1\})$.
\end{lem}
\begin{proof}
The proof is similar to that of \cref{thm:C5:finiteseriesI}. It suffices to show that for $j \geq 1$, in a neighbourhood of $\{\tau = 1\}$, $\omega^j = \pa_\lambda^j| \omega_{\{\lambda = 0\}}$ has a finite expansion in integral powers $\tau^\ell$ without remainder
\[\omega^j = \sum_{\ell=0}^{j-2} \omega_\ell^j \tau^j,\]
where $\omega_\ell^j \in C^{\infty}(\{\tau = 1, \lambda = 0\})$. Indeed, this implies that for $j \geq 1$, $p_k(\tau\pa_\tau)\pa_\lambda^j \omega|_{\{\lambda = 0\}} = 0$ for $k \geq j-1$, and hence the all but the zeroth of the first $k+1$ $\pa_\lambda$-derivatives of $p_k(\tau\pa_\tau)\omega|_{\{\tau = 1\}}$ vanish at $\lambda = 0$. However, we already know for all $k$ that $p_k(\tau\pa_\tau)\omega \in \lambda^2 C^{\infty}(\{\tau = 1\})$ from \cref{thm:C5:prop1}, and so the first $k+1$ $\pa_\lambda$-derivatives of $p_k(\tau\pa_\tau)\omega|_{\{\tau = 1\}}$ vanish.

Commuting \eqref{eq:C5:omegalt} with $\pa_\lambda^{j-1}$ shows that for $j \geq 1$, $\omega^j$ satisfies
\begin{equation}\label{eq:C5:seriesomegai}(\tau\pa_\tau + (1-j))\omega^j = \pa_\lambda^{j-1}\left.\left(\frac{\tau^2}{16}f\tilde{h} - \frac{\tau^2}{4}F'\tilde{H}R^2 - \frac{\tau}{4}f\right)\right|_{\{\lambda = 0\}}.\end{equation}

Now by \eqref{eq:C5:fhlt}, $\tau^2 f\tilde{h}$ is $\lambda^3$ times a smooth function of $\lambda,\theta, \lambda\tau$, and thus $\pa_\lambda^{j-1} (\tau^2 f\tilde{h})|_{\{\lambda = 0\}}$ has a series expansion in powers of $\tau^\ell$ for $0 \leq \ell \leq j-4$. Similarly, $\tau f$ is $\lambda$ times a smooth function of $\lambda$, $\tau\lambda$, $\theta$, and so 
$\pa_\lambda^{j-1}(\tau^2 f)|_{\{\lambda = 0\}}$ has an expansion in powers of $\tau^\ell$ for $0 \leq \ell \leq j-2$.

From the proof of \cref{thm:C5:finiteseriesI}, we also have that for $s \geq 0$, $\pa_\lambda^{s}F'|_{\{\lambda = 0\}}$ and $\pa_\lambda^s H|_{\{\lambda = 0\}}$ have expansions in $\tau^\ell$ for $-1 \leq \ell \leq s-1$ and $-1 \leq \ell \leq s-2$, respectively. Thus, since $R^2$ depends only on $\theta$, $\pa_\lambda^{j-1}|_{\{\lambda = 0\}}(\tau^2 F'\tilde{H})R^2$ has a finite expansion in powers of $\tau^{\ell}$ for $0 \leq \ell \leq j-2$.

Putting this all together, it follows from \eqref{eq:C5:seriesomegai} that
\[(\tau\pa_\tau + (1-j))\omega^j = \sum_{\ell = 0}^{j-2} R_\ell \tau^\ell,\] for some coefficients $R_{\ell} \in C^{\infty}(\{\tau = 1, \lambda = 0\})$. Solving this equation as in the proof of \cref{thm:C5:finiteseriesI}, it follows that $\omega^j$ has the desired expansion.
\end{proof}

We may now prove part (ii) of \cref{thm:C5:prop2}.
\begin{proof}[Proof of \cref{thm:C5:prop2}:(ii)]
It suffices to show that for all $j,k \in \N_0$
\begin{equation}\label{eq:C5:omegaltsf}p_j(\lambda\pa_\lambda)p_k(\tau\pa_\tau)\omega \in \lambda^{\max(j,2)} \tau^k \mathcal A(U),\end{equation} the $\max(j,2)$ showing that there are no terms $\lambda^0$ or $\lambda^1$ in the polyhomogeneous expansion of $\omega$ at $\{\lambda = 0\}$.

We need to examine the right-hand side of \eqref{eq:C5:omegalt}. Observe from \eqref{eq:C5:fhlt} that $\tau^2 f\tilde{h}$ is $\lambda^3$ times a smooth function of $\tau \lambda$, $\theta$, and $\lambda$. Therefore
\begin{equation}\label{eq:C5:omegalti}p_j(\lambda\pa_\lambda)p_k(\tau\pa_\tau)(\tau^2 f\tilde{h}) \in \lambda^{\max(j,k+3)}\tau^k \mathcal A(U).\end{equation}
Similarly, $\tau f$ is $\lambda$ times a smooth function of $\tau \lambda$, $\theta$ and $\lambda$ and so
\begin{equation}\label{eq:C5:omegaltii}p_j(\lambda\pa_\lambda)p_k(\tau\pa_\tau)(\tau f) \in \lambda^{\max(j,k+1)}\tau^k \mathcal A(U).\end{equation}
Now, $R^2$ is only a function of $\theta$, and the proof of \cref{thm:C5:prop2}:(i) in fact shows that
\[p_j(\lambda\pa_\lambda)a_k(\tau\pa_\tau) F' \in \lambda^{\max(j,k)}\tau^{k-1}\mathcal A(U), \ p_j(\lambda\pa_\lambda)a_k(\tau\pa_\tau)\tilde{H} \in \lambda^{\max(j,k)+1}\tau^{k-1}\mathcal A(U).\]
Observe that for a function $X$, $p_j(\tau\pa_\tau)(\tau X) = \tau a_j(\tau\pa_\tau)(X)$.
Thus by \cref{thm:C5:product},
\[p_j(\lambda\pa_\lambda)p_k(\tau\pa_\tau) (\tau^2 F'\tilde{H} )\] is a sum of terms of the form
\[\tau^2 (p_{j-j'}(\lambda\pa_\lambda)a_{k-k'}(\tau\pa_\tau)F')(p_{j'}(\lambda\pa_\lambda)a_{k'}(\tau\pa_\tau)\tilde{H})\]
where $0 \leq j' \leq j$, $0 \leq k' \leq k$, and such a term is in
\[\lambda^{\max(j-j',k-k')+\max(j',k')+1}\tau^{k}\mathcal A(U).\]
Using the reverse triangle inequality as in the proof of \cref{thm:C5:prop2}:(i) shows that $\max(j-j',k-k')+\max(j',k') \geq \max(j,k)$. Thus
\begin{equation}\label{eq:C5:omegaltiii}p_j(\lambda\pa_\lambda)p_k(\tau\pa_\tau) (\tau^2 F'\tilde{H}) \in \lambda^{\max(j,k+1)}\tau^{k}\mathcal A(U).\end{equation}

Commuting $p_j(\lambda\pa_\lambda)p_k(\tau\pa_\tau)$ with \eqref{eq:C5:omegalt} and using \eqref{eq:C5:omegalti}--\eqref{eq:C5:omegaltiii} shows that
\[\pa_\lambda(\tau\pa_\tau - \lambda\pa_\lambda) p_j(\lambda\pa_\lambda)p_k(\tau\pa_\tau)\omega \in \lambda^{\max(j,k+1)}\tau^{k}\mathcal A(U).\]

Now $(\tau\pa_\tau - \lambda\pa_\lambda) p_j(\lambda\pa_\lambda)p_k(\tau\pa_\tau)\omega = 0$ on $\{\lambda = 0\}$ since $\omega$ is smooth away from $\{\tau = 0\}$ and vanishes on $\{\lambda = 0\}$. Thus by the first part of \cref{thm:C5:omegalthelper}, 
\[(\tau\pa_\tau - \lambda\pa_\lambda) p_j(\lambda\pa_\lambda)p_k(\tau\pa_\tau)\omega \in \lambda^{\max(j,k+1)+1}\tau^{k}\mathcal A(U).\] \Cref{thm:C5:finiteseriesomegaI} shows that $p_j(\lambda\pa_\lambda)p_k(\tau\pa_\tau)\omega \in \lambda^{\max(j,k+2)}C^{\infty}(\{\tau = 1\})$, and therefore we may apply the second part of \cref{thm:C5:omegalthelper} to conclude \eqref{eq:C5:omegaltsf}.
\end{proof}

\todo{check reference in texorpdfstring}\section{The proof of \texorpdfstring{\cref{thm:C5:prop3}}{proposition~5.2.5}}
\label{C:C5:vz}
Let us set 
\begin{equation}\label{eq:C5:sigmadefi}\sigma = \zeta - \kappa(1+vK)R,\end{equation}
a bdf of $\mathbf{sgf}$.
We will use coordinates $(v,\zeta,\theta)$, for $\theta$ in some compact neighbourhood, to cover $U_3$. Let $U \subseteq U_3$ denote an arbitrary such coordinate chart, which is compact. For the rest of this section, we work only in $U$. In these coordinates \eqref{eq:C5:FH} becomes
\begin{align}
\label{eq:C5:FHvz}
\begin{split}
0 &= \pa_\zeta F' + \frac{v}{4}f\tilde{H} + \frac{v}{4}\tilde{h}F' + \frac{F'}{\zeta}\\
0 &= (\zeta\pa_\zeta-v\pa_v) \tilde{H} - \frac{v}{4}f\tilde{H} - \frac{v}{4}\tilde{h}F' - \frac{F'}{\zeta}\end{split}\end{align}
and \eqref{eq:C5:omega} becomes
\begin{equation}\label{eq:C5:omegavz}0 = \pa_\zeta(\zeta\pa_\zeta - v\pa_v)\omega - \frac{v^2}{16}f\tilde{h} + \frac{v^2}{4}F'\tilde{H}R^2 + \frac{v}{4}f.\end{equation}
(observe that $\zeta > 0$ on $U$, since $\sigma \geq 0$ implies $\zeta \geq \kappa (1+vK)R$).
From \cref{thm:C5:prop1}, we know smoothness away from the boundaries $\{v=0\}$ and $\{\sigma = 0\}$. Since $\zeta_0 \geq 1/\lambda_0$, \cref{thm:C5:prop2} implies that we know
\begin{align*}F',\tilde{H} \in \phgi{(\mathbf{frf})}{(-1)}(U_3 \n \{\zeta_0 \leq \zeta \leq 1/\lambda_0,\ 0 \leq \zeta v \leq 1\})\\
\omega \in \phgi{(\mathbf{frf})}{(0)}(U_3 \n \{\zeta_0 \leq \zeta \leq 1/\lambda_0,\ 0 \leq \zeta v \leq 1\}).\end{align*} Combining these yields
\begin{align*}F',\tilde{H} \in \phgi{(\mathbf{frf})}{(-1)}(U_3 \n \{\zeta_0 \leq \zeta \leq 1/\lambda_0\})\\
\omega \in \phgi{(\mathbf{frf})}{(0)}(U_3 \n \{\zeta_0 \leq \zeta \leq 1/\lambda_0\}).\end{align*} 
Thus we will treat \eqref{eq:C5:FHvz} and \eqref{eq:C5:omegavz} as initial value problems with data for $F'$,\ $\tilde{H}$ and $\omega$ given on $\{\zeta = \zeta_0\}$. The spherical variables continue to act like a parameter, so we will usually ignore them unless we are trying to establish regularity in the spherical directions. \begin{figure}[htbp]
\centering
\includegraphics{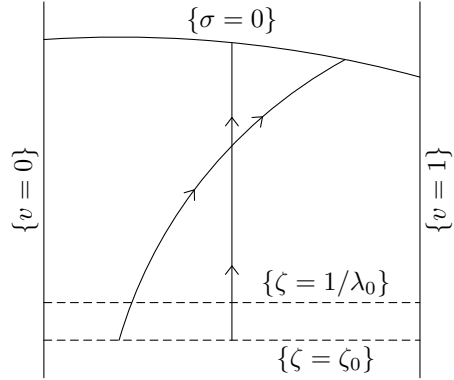}
\caption{A view of $U = U_3$ at an angle $\theta$. Indicated is a sample backwards-directed integral curve of $\pa_\zeta$ and $\zeta\pa_\zeta - v\pa_v$.}
\end{figure}

Propagating the polyhomogeneity at $\mathbf{frf}$ to all of $\{\sigma > 0\}$ and then establishing polyhomogeneity at $\mathbf{sgf}$ will require different arguments, to we split \cref{thm:C5:prop3} into two different propositions:
\begin{prop}\label{thm:C5:prop3i}For all $\epsilon > 0$ sufficiently small depending only on $\zeta_0$ and the specific choice \eqref{eq:C5:sigmadefi} of the bdf $\sigma$ of $\mathbf{sgf}$:
\begin{romanumerate}
\item$F',\tilde{H} \in \phgi{(\mathbf{frf})}{(-1)}(U \n \{\sigma \geq \epsilon\});$
\item $\omega \in \phgi{(\mathbf{frf})}{(0)}(U \n \{\sigma \geq \epsilon\}).$
\end{romanumerate}
\end{prop}
and
\begin{prop}\label{thm:C5:prop3ii}For all $0 < \epsilon < 1$ sufficiently small depending only on $\zeta_0$ and the specific choice \eqref{eq:C5:sigmadefi} of the bdf $\sigma$ of $\mathbf{sgf}$:
\begin{romanumerate}
\item$F',\tilde{H} \in \phgi{(\mathbf{frf},\mathbf{sgf})}{(-1,-1')}(U \n \{\sigma \leq \epsilon\});$
\item $\omega \in \phgi{(\mathbf{frf},\mathbf{sgf})}{(0,0'')}(U \n \{\sigma \leq \epsilon\}).$
\end{romanumerate}
\end{prop}

The splitting of \cref{thm:C5:prop3} into two makes sense in light of hyperbolic nature of \eqref{eq:C5:FHvz} and \eqref{eq:C5:omegavz}. Indeed, one may compute
\begin{align*}
\pa_\zeta \sigma &= 1\\
(\zeta\pa_\zeta-v\pa_v)\sigma &= \kappa(vK+v^2\pa_v K)R + \zeta = \kappa(1+2vK+v^2\pa_v K)R +\sigma\\
&= \mathbf{E}/P+\sigma,
\end{align*}
where $P = 4\kappa v(1+vK)R > 0$. Since $\mathbf{E} \geq 0$ (with equality only if $v = 0$) and $(\mathbf{E}/P)(0) = \kappa R > 0 \mathbf{E}/P > 0$ on $U$, and hence $(\zeta\pa_\zeta - v\pa_v)\sigma > 0$ on $U$.
\begin{rk}If one takes a different bdf instead of $\sigma$, then this may only be true for $\sigma$ sufficiently small rather than on all of $U$. The choice \eqref{eq:C5:sigmadefi} is in some sense ``lucky.''\end{rk}

Thus $\sigma$ is increasing along the flows of $\pa_\zeta$ and $\zeta\pa_\zeta -v\pa_v$. Notice that we are propagating backwards along the flow of $\pa_\zeta$ and $\zeta\pa_\zeta - v\pa_v$, so a backwards domain of dependence is forwards along the flow. In particular, $\{\sigma \geq \epsilon\}$ forms a domain of dependence for \eqref{eq:C5:FHvz} and \eqref{eq:C5:omegavz}, and so what happens in $\{\sigma \leq \epsilon\}$ cannot affect what happens in $\{\sigma \geq \epsilon\}$ (as we are propagating backwards along the flow).

Let us now specify how small $\epsilon$ needs to be in \cref{thm:C5:prop3i} and \cref{thm:C5:prop3ii}. For \cref{thm:C5:prop3i}, $\epsilon$ should be small enough so $\epsilon < \inf \sigma|_{\{\zeta = \zeta_0\}}$ and $\inf \zeta|_{\{\sigma = \epsilon\}} \leq \zeta_0$. This ensures that every backwards-directed integral curve of $\pa_\zeta$ or $\zeta\pa_\zeta - v\pa_v$ which starts at $\{\zeta = \zeta_0\}$ eventually intersects $\{\sigma = \epsilon\}$, and conversely, every forwards-directed integral curve starting at $\{\sigma = \epsilon\}$ intersects $\{\zeta = \zeta_0\}$. For \cref{thm:C5:prop3ii}, we need to additionally take $\epsilon < 1$ (for technical reasons). Notice that every forwards-directed integral curve starting at $\{\sigma = 0\}$ eventually intersects $\{\sigma = \epsilon\}$, so the latter is a Cauchy hypersurface of \eqref{eq:C5:FHvz} and \eqref{eq:C5:omegavz}.

We will prove \cref{thm:C5:prop3i} and \cref{thm:C5:prop3ii} over the next two subsections, devoting one subsection to each proposition.

\subsection{Behaviour away from \texorpdfstring{$\{\sigma = 0\}$}{\{sigma = 0\}}.} We first prove part (i) of \cref{thm:C5:prop3i}.
For technical purposes, we will need to localize to domains of dependence of \eqref{eq:C5:FHvz}. For our purposes:
\begin{defn}A \emph{domain of dependence} for \eqref{eq:C5:FHvz} is subset $\mathcal D \subseteq U$ for which every forwards-directed integral curve of $\pa_\zeta$ or $\zeta\pa_\zeta-v\pa_v$ which starts in $\mathcal D$ does not exit $\mathcal D$ until it intersects $\{\zeta = \zeta_0\}$.\end{defn}
We also define:
\begin{defn}For $p \in U$, the \emph{backwards domain of dependence} of $p$, denoted $B(p)$, is the smallest domain of dependence containing $p$.\end{defn}
One may explicitly describe $B(p)$. If $p = (v_1,\zeta_1)$, then
\[B(p) = \{(v,\zeta), \: \zeta_0 \leq \zeta \leq \zeta_1, \ v_1\zeta_1/\zeta \leq v \leq v_1\}.\]
Indeed, one can explicitly check that $B(p)$ is a domain of dependence, and is the smallest set containing the forwards-directed integral curve of $\zeta\pa_\zeta$ starting from $p$ and all the forwards-directed integral curves of $\zeta\pa_\zeta - v\pa_v$ starting at any point of the previous curve. As remarked above, observe that if $\epsilon < 1$, then $\sigma(p) \geq \epsilon$ implies $\sigma(q) \geq\epsilon$ for all $q \in B(p)$, since every backwards-directed integral curve of $\zeta\pa_\zeta$ and $\zeta\pa_\zeta - v\pa_v$ starting at a point $q'$ with $\sigma(q') \leq \epsilon$ remains in the region $\{\sigma \leq \epsilon\}$.

We will need to to expand the type of equations we consider to include an inhomogeneous right-hand side. Such equations will arise when commuting our original equation with derivatives which we will use to establish regularity. Let us consider for $k \in \N_0$ the system\begin{align}
\label{eq:C5:mevz}
\begin{split}
\pa_\zeta &A + a_{11}B +a_{12}A = R_A\\
(\zeta\pa_\zeta-v\pa_v) &B + a_{21}B + a_{22}A = R_B,\end{split}\end{align}
where $a_{ij} \in L^{\infty}(\{\sigma \geq \epsilon\})$, $R_A,R_B \in v^{k-1}L^{\infty}(\{\sigma \geq \epsilon\})$
and \[A|_{\{\zeta = \zeta_0\}},B|_{\{\zeta = \zeta_0\}} \in v^{k-1}L^{\infty}(\{\zeta = \zeta_0\}).\]
Set
\begin{align*}
\mathcal E &= \sum_{i,j=1}^2 \norm{a_{ij}}_{ L^{\infty}(\{\sigma \geq \epsilon\})}\\
\mathcal I &=\norm{A|_{\{\zeta = \zeta_0\}}}_{v^{k-1}L^{\infty}(\{\zeta = \zeta_0\})} + \norm{B|_{\{\zeta = \zeta_0\}}}_{v^{k-1}L^{\infty}(\{\zeta = \zeta_0\})}\\
&+\norm{R_A}_{v^{k-1}L^{\infty}(\{\sigma \geq \epsilon\})} + \norm{R_B}_{v^{k-1}L^{\infty}(\{\sigma \geq \epsilon\})}.
\end{align*}
Then:
\begin{prop}\label{thm:C5:initialvz}Suppose $A,B \in C^{1}(\{v > 0\})$ solves \eqref{eq:C5:mevz}. Then $A,B \in v^{k-1}L^{\infty}(\{\sigma \geq \epsilon\})$ and
\[\norm{A}_{v^{k-1}L^{\infty}(\{\sigma \geq \epsilon\})} + \norm{B}_{v^{k-1}L^{\infty}(\{\sigma \geq \epsilon\})} \leq C\mathcal I,\]
where $C = C(k,\mathcal E, \zeta_0,\sup_{\sigma \geq \epsilon} \zeta^{-1})$ (the fourth argument is finite) is increasing in all its arguments.
\end{prop}
\begin{proof}
We will let $C$ denote a constant depending only on $\zeta_0$, $\sup \zeta^{-1}$, $\mathcal E$, and explicit numerical constants, which may vary from line to line.
On a first reading, it is helpful to suppose $k=0$ and $R_A = R_B = 0$.
We recast the system as an integral system.
Define $\phi,\psi$ by requiring them to solve in $U$
\begin{align*}
\pa_\zeta \phi &= \phi a_{12}\\
(\zeta\pa_\zeta-v\pa_v)\psi &= \psi a_{21}\end{align*}
with $\phi|_{\{\zeta = \zeta_0\}},\psi|_{\{\zeta = \zeta_0\}} \equiv 1$. Then

 \begin{align*}
 \phi(v,\zeta) &= \exp\left(\int_{\zeta_0}^{\zeta} a_{12}(v,t)\ dt\right)\\
 \psi(v,\zeta) &= \exp\left(\int_{\zeta_0}^{\zeta} a_{21}(v\zeta/s,s)\ ds/s\right),\end{align*} and so
 $\frac{1}{C} \leq \phi,\psi \leq C$.
 
 Now using the method of integrating factors,
 \begin{align}
\label{eq:C5:bootstrapvz}
\begin{split}
A(v,\zeta) &= \phi^{-1}\left(A(v,\zeta_0) - \int_{\zeta_0}^{\zeta} \left(\phi a_{11}B\right)(v,t)\ dt + \int_{\zeta_0}^{\zeta} \left(\phi R_A\right)(v,t)\ dt\right)\\
B(v,\zeta) &= \psi^{-1}\left(B(v\zeta/\zeta_0,\zeta_0)- \int_{\zeta_0}^{\zeta} \left(\psi a_{22} A\right)(v\zeta/s,s)\ ds/s + \int_{\zeta_0}^{\zeta} (\psi R_B)(v\zeta/s,s)\ ds/s\right).
\end{split}
\end{align}

Observe that the integrals are taken along the flow of $\pa_\zeta$ and $\zeta\pa_\zeta-v\pa_v$, respectively.

We will use a bootstrap argument, and find numerical constants $K,\ N > 0$ large that for any $p \in U$ with $\sigma(p) \geq \epsilon$ and $v(p) > 0$, and for and $\delta > 0$ small, then inside $B(p)$
\[|A|,|B| \leq K(\mathcal I+\delta)v^{k-1}e^{N(\zeta_0-\zeta)}\]
(observe that $e^{N(\zeta_0-\zeta)}$ is increasing in the direction in which we are propagating). Then taking $\delta \to 0$ will prove the proposition, since $K,\ N$ do not depend upon $p$. We need only check that if the bootstrap assumptions are satisfied in a neighbourhood of $B(p)\n \{\zeta = \zeta_0\}$, and, for $\zeta(p) \leq \zeta_1 < \zeta_0$, if they are satisfied in $B(p)\n\{\zeta_1 \leq \zeta < \zeta_0\}$, then $A,B$ in fact satisfy an improved estimate on the region. 
 
By continuity, for any $p$, the bounds are satisfied near $\{\zeta = \zeta_0\}$, provided $K > 1$.\footnote{This is essentially the reason for localizing to $B(p)$. We don't know that $A,B \in v^{k-1}C^0(\{\sigma \geq \epsilon\})$ yet, so in order to perform a bootstrap argument we need to restrict to a domain in the interior where we know continuity.} Observe that if $q = (v,\zeta) \in B(p) \n \{\zeta_1 \leq \zeta < \zeta_0\}$, then the arguments of the integrands of \eqref{eq:C5:bootstrapvz} lie in the same region (and in particular in $\{\sigma \geq \epsilon\})$. Therefore we may plug in the bootstrap assumptions into \eqref{eq:C5:bootstrapvz} (and use that $\zeta/s \leq 1$ if $\zeta \leq s$) to see that if $q \in B(p) \n \{\zeta_1 \leq \zeta < \zeta_0\}$
\begin{align*}
|A|(v,\zeta) &\leq C|A|(v,\zeta_0) + CK(\mathcal I + \delta) \int_{\zeta}^{\zeta_0} v^{k-1}e^{N(\zeta_0-t)}\ dt + C\mathcal I \int_{\zeta}^{\zeta_0}v^{k-1}\ dt\\
&\leq C(\mathcal I+\delta) v^{k-1}\left(C- \frac{K}{N} + \frac{K}{N}e^{N(\zeta_0-\zeta)}\right).
\end{align*}

and
\begin{align*}
|B|(v,\zeta) &\leq C|B|(v\zeta/\zeta_0,\zeta_0) + CK(\mathcal I + \delta)\int_{\zeta}^{\zeta_0} (v\zeta/s)^{k} s (v\zeta)^{-1} e^{N(\zeta_0-s)}\ ds/s\\ &+ C\mathcal I\int_{\zeta}^{\zeta_0} (v\zeta/s)^k s(v\zeta)^{-1}\ ds/s\\
&\leq C\mathcal I v^{k-1} + CK(\mathcal I + \delta)\zeta^{-1}v^{k-1}\int_{\zeta}^{\zeta_0}e^{N_1(\zeta_0-s)}\ ds + C\mathcal I v^{k-1}\zeta^{-1}\int_{\zeta}^{\zeta_0} 1\ ds\\
&\leq C(\mathcal I+\delta) v^{k-1}\left(1 + \frac{\zeta_0}{\zeta} - \frac{K}{N\zeta} + \frac{K}{N\zeta}e^{N(\zeta_0-\zeta)}\right)\\
&\leq C(\mathcal I+\delta) v^{k-1}\left(C - \frac{K}{N\zeta} + \frac{K}{N\zeta}e^{N(\zeta_0-\zeta)}\right).
\end{align*}

We will improve over the bootstrap assumption provided $N > C(1+\sup \zeta^{-1})$, and $K > CN(1+\zeta^{-1})$.
\end{proof}

With additional regularity assumptions, we have more regularity of the solution.
\begin{cor}\label{thm:C5:corvz}Let $a_{ij}$ ($i,j, = 1,2$) $R_A, R_B$, be as in \eqref{eq:C5:mevz}. Fix $N \geq 1$ and suppose $a_{ij} \in W_b^N(\{\sigma \geq \epsilon\})$, $R_A,R_B \in v^{k-1}W_b^N(\{\sigma \geq \epsilon\})$. If $A,B \in C^{N+1}(\{\sigma \geq \epsilon, v > 0\})$ solve \eqref{eq:C5:mevz} with data \[A|_{\{\zeta = \zeta_0\}},\ B|_{\{\zeta = \zeta_0\}} \in v^{k-1}W_b^N(\{\zeta = \zeta_0\}),\]
then $A,B \in v^{k-1}W_b^N (\{\sigma \geq \epsilon\})$, with analogous bounds to those in \cref{thm:C5:initialvz}. In particular, if the assumptions are satisfied for all $N$ (i.e.\ with $\mathcal A$ instead of $W_b^N$), then $A,B \in v^{k-1}\mathcal A (\{\sigma \geq \epsilon\})$.\end{cor}
The proof is the almost the same as that of \cref{thm:C5:corlt}, so is omitted. 

We may use the corollary to prove part (i) of \cref{thm:C5:prop3i}.
\begin{proof}[Proof of \cref{thm:C5:prop3i}:(i)]It suffices to show that for all $k$
\begin{align*}
a_k(v\pa_v) F', a_k(v\pa_v)\tilde{H} \in v^{k-1}\mathcal A(\{\sigma \geq \epsilon\}).
\end{align*}
We can easily prove this by induction on $k$. First, recall from \cref{thm:C5:prop2} that $F',\tilde{H} \in \mathcal A_{\mathrm{phg},(\mathbf{frf})}^{(-1)}(\{\zeta = \zeta_0\})$. Thus case $k=0$ is handled by \cref{thm:C5:corvz}. Now for higher $k$, commuting $a_k(v\pa_v)$ with \eqref{eq:C5:FHvz}, we see that $a_k(v\pa_v) F',a_k(v\pa_v)\tilde{H}$ satisfies the same equation as $F',\tilde{H}$, but with a right-hand side involving terms of the form
\[(p_{k-k'}(v\pa_v)\phi)(a_{k'}(v\pa_v) A),\] where $\phi$ is smooth, $A$ is one of $F',\tilde{H}$, and $0 \leq k' < k$. The first factor is in $v^{k-k'}C^{\infty}(\{\sigma \geq \epsilon\})$. The second by induction is in $v^{k'-1}\mathcal A(\{\sigma \geq \epsilon\})$. Thus the right-hand side is in $v^{k-1}\mathcal A(\{\sigma \geq \epsilon\})$. By assumption, the data for $a_k(v\pa_v) F'$ and $a_k(v\pa_v) \tilde{H}$ is in $\mathcal A_{\mathrm{phg}}^{(k-1)}(\{\zeta = \zeta_0\})$. Thus \cref{thm:C5:corvz} applies and we may conclude the inductive step.\end{proof}

We now turn our attention to proving part (ii) of \cref{thm:C5:prop3i}. Let us consider for $k \in \R$ the equations:
\begin{equation}\label{eq:C5:wvzi} \pa_\zeta w = S\end{equation}
and
\begin{equation}\label{eq:C5:wvzii} (\zeta\pa_\zeta -v\pa_v)w = S,\end{equation}
where in both cases $S \in v^kL^{\infty}(\{\sigma \geq \epsilon\})$.
We start with a general proposition.
\begin{prop}\label{thm:C5:omegavzhelper}Fix $k \in \R$. Suppose $w \in C^1(\{\sigma \geq \epsilon, v > 0\})$ solves \eqref{eq:C5:wvzi} or \eqref{eq:C5:wvzii}
with data $w|_{\{\zeta = \zeta_0\}} \in v^k L^{\infty}(\{\zeta = \zeta_0\})$. Then $w \in v^k L^{\infty}(\{\sigma \geq \epsilon\})$ with the bound
\[\norm{w}_{v^k L^{\infty}(\{\sigma \geq \epsilon\})} \leq C(\norm{w|_{\{\zeta = \zeta_0\}}}_{v^kL^{\infty}(\{\zeta = \zeta_0\})} + \norm{S}_{v^k L^{\infty}(\{\sigma \geq \epsilon\})}),\]
where $C = C(k,\zeta_0,\sup_{\{\sigma \geq \epsilon\}} 1/\zeta)$ is increasing in all its arguments.

Now fix $N \geq 1$. If instead $w \in C^{N+1}(\{v > 0, \sigma \geq \epsilon\})$ and $S \in v^k W_b^N(\{\sigma \geq \epsilon\})$ and $w|_{\{\zeta = \zeta_0\}} \in v^kW_b^N(\{\zeta = \zeta_0\})$, then $w \in v^k W_b^N(\{\sigma \geq \epsilon\})$ with an analogous bound. In particular, if the hypotheses hold for all $N$, then $w \in v^{k} \mathcal A(\{\sigma \geq \epsilon\})$.\end{prop}
\begin{proof}
The proof is the essentially the same as that of \cref{thm:C5:omegalthelper}, with the integrals appearing taken along the flows as in the proof of \cref{thm:C5:initialvz}, so is omitted. We only mention that since $\{\sigma \geq \epsilon\}$ is disjoint from $\{\zeta = 0\}$, we do not need to worry about $\log$s appearing.
\end{proof}
We can now prove 
part (ii) of \cref{thm:C5:prop3i}.
\begin{proof}[Proof of \cref{thm:C5:prop3i}:(ii)]
It suffices to show that for all $k \in \N_0$
\begin{equation}\label{eq:C5:sfvz}p_k(v\pa_v)\omega \in v^{k}\mathcal A(U).\end{equation}
We need to examine the right-hand side of \eqref{eq:C5:omegavz}. Observe from \eqref{eq:C5:fhvz} that $v^2f \tilde{h}$ is in $vC^{\infty}(\{\sigma \geq \epsilon\})$ and $vf \in C^{\infty}(\{\sigma \geq \epsilon\})$. From \cref{thm:C5:prop3i}:(i),
\[vF',v\tilde{H} \in \phgi{(\{v = 0\}}{(0)}(\{\sigma \geq \epsilon\}) = C^{\infty}(\{\sigma \geq \epsilon\}).\] It follows from \eqref{eq:C5:omegavz} that
\[\pa_\zeta(\zeta\pa_\zeta - v\pa_v)\omega \in C^{\infty}(\{\sigma \geq \epsilon\}),\] and so commuting $p_k(v\pa_v)$ that
\[\pa_\zeta(\zeta\pa_\zeta - v\pa_v)p_k(v\pa_v)\omega \in v^k C^{\infty}(\{\sigma \geq \epsilon\}).\] Since $\omega|_{\{\zeta = \zeta_0\}} \in C^{\infty}(\{\zeta = \zeta_0\})$ because we know $\omega \in \phgi{(\mathbf{frf})}{(0)}(U_3 \n \{\zeta_0 \leq \zeta \leq 1/\lambda_0\}$, it follows that
\[p_k(v\pa_v)\omega|_{\{\zeta = \zeta_0\}} \in v^{k}C^{\infty}(\{\zeta = \zeta_0\})\] and so we may apply \cref{thm:C5:omegavzhelper} to conclude \eqref{eq:C5:sfvz}.
\end{proof}

\subsection{Behaviour near \texorpdfstring{$\{\sigma = 0\}$}{\{sigma = 0\}}.}
We begin by showing \cref{thm:C5:prop3ii}.

Define
\[\beta = (\zeta\pa_\zeta - v\pa_v) \sigma = \kappa(1+2vK + v^2\pa_v K)R + \sigma.\] Recall that we have shown that $\beta > 0$ at the beginning of this section. Since $U$ is compact, $\beta$ is bounded uniformly away from $0$. Let us introduce $\Phi = \sigma \tilde{H}$ and $\Psi = -\sigma/\beta F'$. From \eqref{eq:C5:fhvz}, 
\[vf = \frac{-2\beta}{\sigma}(1+\sigma C^{\infty}(U)), \ v\tilde{h} = \frac{2}{\sigma}(1+\sigma C^{\infty}(U)).\]

We focus first on proving part (i) of \cref{thm:C5:prop3ii}.

Observe that $\zeta^{-1}$ is smooth for $\sigma \geq 0$, and so \eqref{eq:C5:FHvz} becomes
\begin{align}
\label{eq:C5:PPvz}
\begin{split}
\sigma(\zeta\pa_\zeta-v\pa_v)\Phi &= \frac{\beta}{2}(1+\sigma C^{\infty}(U))\Phi -\frac{\beta}{2}(1+C^{\infty}(U))\Psi\\
\sigma\pa_\zeta \Psi &= -\frac{1}{2}(1+\sigma C^{\infty}(U))\Phi + \frac{1}{2}(1+\sigma C^{\infty}(U))\Psi
\end{split}.\end{align}

In order to get a sense of this equation, one should think of it as essentially the system of ODEs \[\sigma\pa_\sigma \begin{pmatrix} \Upsilon_1\\ \Upsilon_2\end{pmatrix} = \begin{pmatrix} 1/2 & -1/2\\-1/2 & 1/2\end{pmatrix}\begin{pmatrix} \Upsilon_1\\ \Upsilon_2\end{pmatrix}.\] This is in line as thinking of $\sigma(\zeta\pa_\zeta-v\pa_v)$ and $\sigma\pa_\zeta$ as ``warped'' versions of $\sigma\pa_\sigma$. The presence of $\beta$ accounts for the fact that $\sigma(\zeta\pa_\zeta - v\pa_v)\sigma = \beta\sigma$ rather than $\sigma$. This linear ODE has its behaviour as $\sigma \to 0$ determined by the eigenvalues of $\begin{psmallmatrix} 1/2 & -1/2\\-1/2 & 1/2\end{psmallmatrix}$, which are $0$ and $1$. Thus one expects solutions to be bounded, and solutions to have series in $\sigma^n$ for $n \geq 0$ and $\sigma^n \log \sigma$ for $n \geq 1$. In other words, one expects polyhomogeneity at $\{\sigma = 0\}$ with index set $0'$.

As above, we will need to to add an inhomogeneous right-hand side to \eqref{eq:C5:PPvz}, and also extend from a pair of scalar functions to a pair of vector-valued function on $\R^d$ ($d \geq 1$). Consider the following $2d\times 2d$ system, for $k \in \N_0$
\begin{subequations}
\label{eq:C5:mevzs}
\begin{align}
\label{eq:C5:mevzsa}
\sigma(\zeta\pa_\zeta-v\pa_v)A &= \frac{\beta}{2}(1+\sigma E_{11})A -\frac{\beta}{2}(1+\sigma E_{12})B + \sigma R_A\\
\label{eq:C5:mevzsb}
\sigma\pa_\zeta B &= -\frac{1}{2}(1+\sigma E_{21})A + \frac{1}{2}(1+\sigma E_{22})B + \sigma R_B,\end{align}
\end{subequations}
Where $E_{ij} \in L^{\infty}(\{0 \leq \sigma \leq \epsilon\};M_d(\R))$, and for some $R_A,R_B \in v^{k-1}L^{\infty}(\{0 \leq \sigma \leq \epsilon\};\R^d)$ (here $M_d(\R)$ denotes the space of $d\times d$ matrices). Suppose further $A|_{\{\sigma = \epsilon\}}, B|_{\{\sigma = \epsilon\}} \in v^{k-1}L^{\infty}(\{\sigma = \epsilon\};\R^d)$.

To define the space $L^{\infty}(\{0 \leq\sigma \leq \epsilon\};M_d(\R))$ and $L^{\infty}(\{0 \leq\sigma \leq \epsilon\};\R^d)$, we use any norm on $\R^d$, and the associated operator norm on $M_d(\R)$ to simplify our computations. Henceforth, we will omit $M_d(\R)$ and $\R^d$ from out notation for clarity. 

Let
\begin{align*}
\mathcal I &= \max(\norm{A|_{\{\sigma = \epsilon\}}}_{v^{k-1}L^{\infty}(\{\sigma = \epsilon\})}, \norm{B|_{\{\sigma = \epsilon\}}}_{v^{k-1}L^{\infty}(\{\sigma = \epsilon\})})\\
\mathcal R &= \norm{R_A}_{v^{k-1}L^{\infty}(\{0 \leq \sigma \leq \epsilon\})} + \norm{R_B}_{v^{k-1}L^{\infty}(\{0 \leq \sigma \leq \epsilon\})}\\
\mathcal E &= \sum_{ij} \norm{E_{ij}}_{L^{\infty}(\{0 \leq \sigma \leq \epsilon\})}.\end{align*}

\begin{prop}\label{thm:C5:initialvzs}Suppose $A,B \in C^1(\{0 < \sigma \leq \epsilon, v > 0\})$ solves \eqref{eq:C5:mevzs}. Then $A,B \in v^{k-1}L^{\infty}(\{0 \leq \sigma \leq \epsilon\})$ and
\[\norm{A}_{v^{k-1}L^{\infty}(\{0 \leq \sigma \leq \epsilon\})} + \norm{B}_{v^{k-1}L^{\infty}(\{0 \leq \sigma \leq \epsilon\})} \leq C(\mathcal I + \mathcal R),\] where $C = C(k,\mathcal E, \sup \beta^{-1})$ is increasing in all its arguments.
\end{prop}

Before proving \cref{thm:C5:initialvzs}, we will need some overhead. Fix $0 < \vepsilon \leq \epsilon$, and $p = (v,\zeta,r,u)$ with $\sigma(p) \leq \epsilon$. Let $q_1 = (v_1,\zeta_1,r_1,u_1)$ be the unique point for which the integral curve of $\zeta\pa_\zeta$ starting from $p$ intersects $\{\sigma = \vepsilon\}$.\footnote{The integral curve intersects $\{\sigma = \vepsilon\}$ precisely once because $\sigma$ is strictly increasing along the flow for $\{\sigma \leq \epsilon\}$.} Up to scaling the parameter, an integral curve starting at $p$ is the curve
\[t \mapsto (v,t,\theta).\] Thus $q_1 = (v,\mathbf{i}(v,\theta),\theta)$ for some $\mathbf{i} = \mathbf{i}_\vepsilon$.

Now, let $q_2 = (v_2,\zeta_2,\theta_2)$ be the unique point for which the forwards-directed integral curve of $\zeta\pa_\zeta-v\pa_v$ starting from $p$ intersects $\{\sigma = \vepsilon\}$. Up to scaling the parameter, an integral curve starting at $p$ is the curve
\[t \mapsto (v\zeta/t,t,\theta).\] Thus $q = (v\zeta/\mathbf{j}(v\zeta,\theta),\mathbf{j}(v\zeta,\theta),\theta)$, for some $\mathbf{j} = \mathbf{j}_\vepsilon$.

Since \[\sigma(v,\mathbf{i},\theta) = \sigma(v\zeta/\mathbf{j},\mathbf{j},\theta) = \vepsilon,\]
the implicit function theorem shows that $\mathbf{i},\ \mathbf{j}$ are smooth and satisfy $\pa_\zeta \mathbf{i} = 0$, $(\zeta\pa_\zeta-v\pa_v)\mathbf{j} = 0$.
The following lemma will be useful:
\begin{lem}\label{thm:C5:integration1}Let $f:(0,\infty) \to \R$ be a smooth function, and let $f'$ denote its derivative. Then for $0 < \vepsilon \leq \epsilon$

\begin{gather*}
\int_{\mathbf{i}_{\vepsilon}(v)}^{\zeta} f'(\sigma(v,t))\ dt = f(\sigma(v,\zeta))-f(\vepsilon)\\
\int_{\mathbf{j}_{\vepsilon}(v\zeta)}^{\zeta} \beta f'(\sigma(v\zeta/s,s))\ ds/s = f(\sigma(v,\zeta))-f(\vepsilon).\end{gather*}\end{lem}
\begin{proof}
Observe that since $\pa_\zeta \mathbf{i} = 0$, $(v\pa_v-\zeta\pa_\zeta)\mathbf{j} = 0$,
\begin{align*}
\pa_\zeta \int_{\mathbf{i}_{\vepsilon}(v)}^{\zeta} f'(\sigma(v,t))\ dt &= f'(\sigma(v,\zeta))\\
(\zeta\pa_\zeta-v\pa_v) \int_{\mathbf{j}_{\vepsilon}(v\zeta)}^{\zeta} \beta f'(\sigma(v\zeta/s,s))\ ds/s &= \beta f'(\sigma(v,\zeta)) + \int_{\mathbf{j}}^{\zeta} (\zeta\pa_\zeta-v\pa_v)f(\sigma(v\zeta/s,s))\ ds/s\\
&= \beta f'(\sigma(v,\zeta). \end{align*}
On the other hand,
\begin{align*}
\pa_\zeta ( f(\sigma(v,\zeta))-f(\vepsilon)) &= f'(\sigma(v,\zeta))\pa_\zeta \sigma = f'(\sigma(v,\zeta))\\
(\zeta\pa_\zeta-v\pa_v) (f(\sigma(v,\zeta))-f(\vepsilon)) &= f'(\sigma(v,\zeta))(\zeta\pa_\zeta-v\pa_v) \sigma = \beta f'(\sigma(v,\zeta)).\end{align*}
Thus the left-hand and right-hand sides in the lemma statement differ by a function constant on the integral curves of $\pa_\zeta$ and $\zeta\pa_\zeta-v\pa_v$, respectively. Since both sides are $0$ if $\sigma(v,\zeta) = \vepsilon$, and both sets of integral curves intersect $\{\sigma = \vepsilon\}$ at some point, it follows that the constant is $0$.\end{proof}We can now prove \cref{thm:C5:initialvzs}. \begin{proof}[Proof of \cref{thm:C5:initialvzs}.] Fix $0 < \vepsilon \leq \epsilon$ and
set $\mathbf{i} = \mathbf{i}_\vepsilon$, $\mathbf{j} = \mathbf{j}_{\vepsilon}$. First of all, conjugating by $v^{1-k}$, we may suppose $k=1$, since this only changes the $\sigma C^{\infty}(U)$ terms. We will let $C$ denote a numerical constant, depending only on $k$, which may change from line to line (in particular for this proof $C$ does not depend on $\mathcal I$, $\mathcal E$, $\epsilon$ or $\sup \beta^{-1}$). On a first reading, it is helpful to suppose $E_{ij} = R_A = R_B = 0$. 
We recast the system as an integral system.
Define $\phi,\psi$ by requiring them to solve in $\{0 \leq \sigma \leq \vepsilon\}$.
\begin{align*}
\sigma(\zeta\pa_\zeta-v\pa_v)\phi &= -\phi\frac{\beta}{2}\\
\sigma\pa_\zeta \psi &= -\psi\frac{1}{2}
\end{align*}
with $\phi|_{\{\sigma = \vepsilon\}},\psi|_{\{\sigma = \vepsilon\}} \equiv 1$. Then,

\begin{gather*}
\phi(v,\zeta) = \exp\left(\int_{\mathbf{j}}^{\zeta} -\frac{\beta}{2\sigma}(v\zeta/s,s)\ ds/s\right)\\
\psi(v,\zeta) = \exp\left(\int_{\mathbf{i}}^{\zeta} -\frac{1}{2\sigma}\ dt\right).
\end{gather*}
Thus by \cref{thm:C5:integration1},
\begin{align*}
\phi &= \frac{\vepsilon^{1/2}}{\sigma^{1/2}}\\
\psi &= \frac{\vepsilon^{1/2}}{\sigma^{1/2}}.
\end{align*}
 
Now using the method of integrating factors,
\begin{align}
\label{eq:C5:bootstrapvzs}
\begin{split}
A(v,\zeta) &= \frac{\sigma^{1/2}}{\vepsilon^{1/2}}\left(A(v\zeta/\mathbf{j},\mathbf{j})+ \int_{\mathbf{j}}^{\zeta} \frac{\vepsilon^{1/2}\beta}{2\sigma^{3/2}}(-B+\sigma E_{11}A + \sigma E_{12}B)(v\zeta/s,s)\ ds/s \right.\\
&\left.+ \int_{\mathbf{j}}^{\zeta} (\vepsilon^{1/2}\sigma^{-1/2} R_A)(v\zeta/s,s)\ ds/s\right)\\
B(v,\zeta) &= \frac{\sigma^{1/2}}{\vepsilon^{1/2}}\left(B(v,\mathbf{i}) + \int_{\mathbf{i}}^{\zeta}\frac{\vepsilon^{1/2}}{2\sigma^{3/2}}(-A - \sigma E_{21}A + \sigma E_{22}B)(v,t)\ dt\right.\\
&\left. + \int_{\mathbf{i}}^{\zeta} (\vepsilon^{1/2}\sigma^{-1/2}R_B)(v,t)\ dt\right).
\end{split}\end{align}

Observe that the integrals are along the flow of $\pa_\zeta$ and $\zeta\pa_\zeta-v\pa_v$, respectively.

Set $\mathcal O = \sup |\beta^{-1}|$. For $p \in U^\z$, write
\[\mathcal I_{\vepsilon} = \max(\norm{A}_{L^{\infty}(B(p)\n \{\sigma = \vepsilon\})},\norm{B}_{L^{\infty}(B(p)\n \{\sigma = \vepsilon\})}).\] We will use a bootstrap argument to show that there exist constants $K_1,K_2 >1,\kappa > 0$ such that in all of $U \n \{\sigma \geq \vepsilon\}$ such that for any $p \in U^\z$, $\sigma(p) \leq \vepsilon$ any $\delta > 0$, \[|A|,|B| \leq (\mathcal I_\vepsilon+\vepsilon K_1\mathcal O \mathcal R+\delta)(1+\vepsilon K_2 \mathcal E + 2\delta)\sigma^{-\kappa\vepsilon\mathcal E-\delta}.\]
Taking $\delta \to 0$ shows that 
\begin{equation}\label{eq:C5:ttt}|A|,|B| \leq (\mathcal I_\vepsilon+\vepsilon K_1\mathcal O\mathcal R)(1+\vepsilon K_2\mathcal E)\sigma^{-\kappa\vepsilon\mathcal E}.\end{equation} With $\vepsilon = \epsilon$, this is apparently weaker than what we want to prove, since there is mild blow-up in $\sigma$ as $\sigma \to 0$. However, the exponent decreases in $\vepsilon$, and so we may iterate this bound on a sequence of initial surfaces $\{\sigma = \vepsilon_i\}$ for $\vepsilon_0 = \epsilon$ and $\vepsilon_i \to 0$ to derive the conclusion. We will show how to do this after we prove the bound. We need only check that the bootstrap assumption is satisfied in a neighbourhood of $B(p) \n \{\sigma = \vepsilon\}$, and for $\sigma(p) \leq \sigma_0 \leq \vepsilon$, that if they are satisfies in $B(p) \n \{\sigma_0 \leq \sigma < \vepsilon\}$, then $A,B$ in fact satisfy an improved estimate on the region.

Certainly the assumptions are satisfied in a small neighbourhood of $B(p) \n \{\sigma = \vepsilon\}$, for any $\kappa, K_1,K_2 >0 $. Observe that if $q = (v,\zeta) \in B(p) \n \{\sigma_0 \leq \sigma \leq \vepsilon\}$, then the arguments of the integrands in \eqref{eq:C5:bootstrapvzs} lie in the same region. There we may plug in the bootstrap assumptions into \eqref{eq:C5:bootstrapvzs} to see that if $q \in B(p) \n \{\sigma_0 \leq \sigma \leq \vepsilon\}$ (using \cref{thm:C5:integration1})
\begin{align*}
|\phi A|(v,\zeta) &\leq \mathcal I_\vepsilon + \int_{\zeta}^{\mathbf{j}} \frac{\vepsilon^{1/2}\beta}{2\sigma^{3/2+\kappa\vepsilon\mathcal E+\delta}(v\zeta/s,s)}(1+2\vepsilon \mathcal E)(\mathcal I_\vepsilon+\vepsilon K_1\mathcal O \mathcal R+\delta)(1+\vepsilon K_2 \mathcal E + 2\delta)\ ds/s\\
&+ \int_{\zeta}^{\mathbf{j}} \beta\vepsilon^{1/2}\sigma^{-1/2} \mathcal R/\beta\ ds/s\\
&\leq \mathcal I_\vepsilon + \frac{\vepsilon^{1/2}(1+2\vepsilon \mathcal E)(\mathcal I_\vepsilon+\vepsilon K_1\mathcal O \mathcal R+\delta)(1+\vepsilon K_2 \mathcal E + 2\delta)}{1+2\kappa\vepsilon\mathcal E + 2\delta}\left(\sigma^{-1/2-\kappa\vepsilon\mathcal E}-\vepsilon^{-1/2-\kappa\vepsilon\mathcal E}\right)\\
&+ 2\mathcal \vepsilon^{1/2}\mathcal O\mathcal R(\vepsilon^{1/2}-\sigma^{1/2}).\end{align*}
Thus, 
\begin{align*}
|A|(v,\zeta) &\leq\frac{(1+2\vepsilon \mathcal E)(\mathcal I_\vepsilon+\vepsilon K_1\mathcal O \mathcal R+\delta)(1+\vepsilon K_2 \mathcal E + 2\delta)}{1+2\kappa\vepsilon\mathcal E + 2\delta}\sigma^{-\kappa\vepsilon\mathcal E-\delta}\\
&+ \phi^{-1}\left(2\vepsilon \mathcal O \mathcal R + \mathcal I_\vepsilon - \frac{(1+2\vepsilon \mathcal E)(\mathcal I_\vepsilon+\vepsilon K_1\mathcal O \mathcal R+\delta)(1+\vepsilon K_2 \mathcal E +2 \delta)}{1+2\kappa\vepsilon\mathcal E + 2\delta}\vepsilon^{-1/2-\kappa\vepsilon\mathcal E-\delta}\right).\end{align*}
A similar bound holds for $B$.

Since $\vepsilon < 1$, we will improve over the bootstrap assumptions provided $\kappa > 2$, $K_2 \geq 2\kappa$, and $K_1 \geq 2$.

Now we will show how \eqref{eq:C5:ttt} implies the desired bounds. We will apply \eqref{eq:C5:ttt} on a sequence of initial surfaces: $\{\sigma = \vepsilon_i\}$ for $\vepsilon_i = 2^{-i}\epsilon$ for $i =0,1,2,\ldots$. Write
\[\mathcal I_i = \mathcal I_{\vepsilon_i}.\]
Introduce the shorthand
\[\Pi_i = (1+K_2\mathcal E 2^{-i}\epsilon)(2^{-i}\epsilon)^{-\kappa 2^{-i}\epsilon\mathcal E}.\] Then the bound \eqref{eq:C5:ttt} applied to $\vepsilon_i$ implies that
\[\mathcal I_i \leq (\mathcal I_{i-1}+\epsilon 2^{-i+1}K_1\mathcal O \mathcal R)\Pi_{i-1}.\] Iterating, we see that
\begin{align*}
\mathcal I_{i} &\leq \cdots \leq \mathcal I_{j}\prod_{k=j}^{i-1}\Pi_k + \sum_{k=j}^{i-1} \epsilon 2^{-k}K_1\mathcal O \mathcal R\prod_{\ell=k}^{i-1}\Pi_{\ell}\\
&\leq \cdots \leq \mathcal I_0\prod_{k=0}^{i-1}\Pi_k + \sum_{k=0}^{i-1} \epsilon 2^{-k}K_1\mathcal O \mathcal R\prod_{\ell=k}^{i-1}\Pi_{\ell}.\end{align*}
Now

\[\prod_{\ell=0}^{\infty} \Pi_{\ell} \leq e^{2K_2\epsilon\mathcal E}2^{10\kappa\epsilon\mathcal E}\epsilon^{10\kappa\epsilon\mathcal E} \leq Ce^{C \mathcal E},\]
since $\epsilon^{-\epsilon} \lesssim 1$ for $0 <\epsilon \leq 1$. Thus
\[\mathcal I_i \leq Ce^{C\mathcal E}(\mathcal I_0 + \epsilon K_1 \mathcal O \mathcal R) \leq Ce^{C\mathcal E}(\mathcal I+ \epsilon K_1 \mathcal O \mathcal R).\]
Applying \eqref{eq:C5:ttt} again shows that on $B(p)\n \{2^{-i-1}\epsilon \leq \sigma \leq 2^{-i}\epsilon\}$
\[|A|,|B| \leq (\mathcal I_i + 2^{-i}\epsilon K_1 \mathcal O \mathcal R)(1+2^{-i}\epsilon K_2\mathcal E)(2^{-i-1})^{-\kappa \epsilon 2^{-i}\mathcal E}) \leq Ce^{C\mathcal E}(\mathcal I + \epsilon K_1 \mathcal O \mathcal R).\] Taking the union over all $i$ proves the proposition.
\end{proof}

With additional regularity assumptions, we have more regularity of the solution, although we may need to take $\epsilon$ slightly smaller.
\begin{cor}\label{thm:C5:corvzs}Let $E_{ij}$ ($i,j = 1,2$), $R_A,R_B$ be as in \eqref{eq:C5:mevzs}. Fix $N \geq 1$ and suppose $E_{ij} \in W_b^N(\{\sigma \leq \epsilon\})$, $R_A,R_B \in v^{k-1}W_b^N(\{\sigma \leq \epsilon\})$. If $A,B \in C^{N+1}(\{0 < \sigma \leq \epsilon, v > 0\})$ solve \eqref{eq:C5:mevz} with data \[A|_{\{\sigma = \epsilon\}},\ B|_{\{\sigma = \epsilon\}} \in v^{k-1}W_b^N(\{\sigma = \epsilon\}),\]
then $A,B \in v^{k-1}W_b^N (\{\sigma \leq \epsilon\})$, with analogous bounds to those in \cref{thm:C5:initialvz} (except now the constant is allowed to depend on $\beta,\beta^{-1}$ and their first $N$ derivatives).

If the assumptions are true for all $N$, then in particular $A,B \in v^{k-1}\mathcal A(\{0 \leq \sigma \leq \epsilon\})$.
\end{cor}
\begin{rk}It is important to remark that $\epsilon$ does not depend on $N$.\end{rk}
\begin{proof}
Unlike in the proofs of \cref{thm:C5:corlt} and \cref{thm:C5:corvz}, we may not find a basis of vector fields tangent to the level sets of $\sigma$ which commute with both $\zeta\pa_\zeta - v\pa_v$ and $\pa_\zeta$, so we will need to adopt a different approach.
Setting $\sigma = \sigma(v,\zeta,\theta)$, $(\sigma, v,\theta)$ become new coordinates on $U \n \{\sigma \leq \epsilon\}$. Consider the vector fields
\[\nu_0 = \sigma\pa_\sigma, \ \nu_1 = v\pa_v, \ \nu_{1+i} = \pa_{\theta^i}\ (i = 1,2).\]
Observe that in these coordinates $\zeta\pa_\zeta - v\pa_v = \beta\pa_\sigma - v\pa_v$ and $\pa_\zeta = \pa_\sigma$.
We record the commutator formulae
\begin{align*}[\nu_i,\sigma(\zeta\pa_\zeta-v\pa_v)] &= \nu_i\log \beta \sigma(\zeta\pa_\zeta- v\pa_v) + \sigma(\nu_i\log \beta-\delta_{i0}) \nu_1\\
[\nu_i,\sigma\pa_\zeta] &= 0,
\end{align*}
where $\delta_{ij}$ denotes the Kronecker $\delta$.
For $\tau \in \N^4$ a multiindex, set $\nu^\tau = \nu_0^{\tau_0}\nu_1^{\tau_1}\nu_2^{\tau_2}\nu_3^{\tau_3}$, and for $N' \leq N$, set $A_{N'}$ to be the vector formed by concatenating each $\nu^\tau A$, where the length of $\tau$ is $N'$.\footnote{Notice that the $\nu_i$ commute with each other, so this covers all derivatives.}
Define $B_{N'}$ similarly.

Let us first show how to prove the corollary for $N=1$. Commuting a $\nu_i$ through \eqref{eq:C5:mevzs} yields
\begin{align}
\label{eq:C5:induct1}
\begin{split}
\sigma(\zeta\pa_\zeta-v\pa_v)\nu_iA &= \frac{\beta}{2}(1+\sigma E_{11})\nu_iA
- \sigma(\nu_i\log \beta-\delta_{i0}) \nu_1 A-\frac{\beta}{2}(1+\sigma E_{12})\nu_iB\\
& +\frac{1}{2}\nu_i(\sigma\beta E_{11})A -\frac{1}{2}\nu_i(\sigma\beta E_{12})B+ \nu_i(\sigma R_A)\\
\sigma\pa_\zeta \nu_iB &= -\frac{1}{2}(1+\sigma E_{21})\nu_iA + \frac{1}{2}(1+\sigma E_{22})\nu_iB\\
&- \frac{1}{2}\nu_i(\sigma \beta E_{21})A - \frac{1}{2}\nu_i(\sigma\beta E_{22})B + \nu_i(\sigma R_B).
\end{split}
\end{align}
From \cref{thm:C5:initialvzs}, $A, \ B \in v^{k-1}L^{\infty}(\{\sigma \leq \epsilon\})$. Thus all terms in the second line of both equations in \eqref{eq:C5:induct1} are in $\sigma v^{k-1}L^{\infty}(\{\sigma \leq \epsilon\})$.

This means \eqref{eq:C5:induct1} almost has the same form as \eqref{eq:C5:mevzs}, except for the presence of the $\sigma(\nu_i\log\beta)\nu_1 A$ term in the first line. However, we can combine the equations \eqref{eq:C5:induct1} for different $i$ into a single equation for $(A_1,B_1)$, and incorporate the term $\sigma(\nu_i\log \beta)$ into the $\sigma E_{11}$ term. This equation takes the form
\begin{align}
\label{eq:C5:induct2}
\begin{split}
\sigma(\zeta\pa_\zeta-v\pa_v)A_1 &= \frac{\beta}{2}(1+\sigma E^1_{11})A_1
-\frac{\beta}{2}(1+\sigma E^1_{12})B_1 + \sigma R_A^1\\
\sigma\pa_\zeta B_1 &= -\frac{1}{2}(1+\sigma E_{21}^1)A_1 + \frac{1}{2}(1+\sigma E_{22}^2)B_1 + \sigma R_B^1,
\end{split}
\end{align}
where $E_{ij}^1$ are linear maps in $L^{\infty}(\{\sigma \leq \epsilon\})$, $E_{ij}^1$ for $(i,j) \neq (1,1)$ is just $E_{ij}$ arranged in a block-diagonal fashion, and $E_{11}^i$ has additional components $\nu_i\log \beta - \delta_{i0}$ in the column corresponding to $\tau = (1)$, and $R_A^1,\ R_A^2$ are terms involving $\nu_i$-derivatives of $R_A$, $R_B$ of order at most $1$, and of $A$, $B$ of order at most $0$.

Provided we can establish that $A_1,B_1 \in v^{k-1}L^{\infty}(\{\sigma = \epsilon\})$, then we can apply \cref{thm:C5:initialvzs} and conclude the $N=1$ case. This is clear for the components $\nu_i A$, $\nu_i B$, for $i \geq 1$ since these involve just tangential derivatives of the initial data for $A$, $B$, respectively. However, using $\sigma\pa_\sigma = \sigma\pa_\zeta$ and \eqref{eq:C5:mevzsb} shows that $\sigma\pa_\sigma B \in v^{k-1}L^{\infty}(\{\sigma = \epsilon\})$, and using $\sigma\pa_\sigma = \beta^{-1}((\zeta\pa_\zeta - v\pa_v) + \nu_1)$ and \eqref{eq:C5:mevzsa} shows that $\sigma\pa_\sigma A \in v^{k-1}L^{\infty}(\{\sigma = \epsilon\})$. Thus we have established the required regularity of the initial data.

For $N > 1$, we use induction. The inductive step is similar to what we have already carried out, so we only provide a sketch. Using induction, one may establish an equation for $(A_{N'},B_{N'})$ for $N' \leq N$ of the following form
\begin{align}
\label{eq:C5:induct3}
\begin{split}
\sigma(\zeta\pa_\zeta-v\pa_v)A_{N'} &= \frac{\beta}{2}(1+\sigma E^{N'}_{11})A_{N'}
-\frac{\beta}{2}(1+\sigma E^{N'}_{12})B_{N'} + \sigma R_A^{N'}\\
\sigma\pa_\zeta B_{N'}&= -\frac{1}{2}(1+\sigma E_{21}^{N'})A_{N'} + \frac{1}{2}(1+\sigma E_{22}^{N'})B_{N'} + \sigma R_B^{N'},
\end{split}
\end{align}
where $E_{ij}^{N'}$ are linear maps in $W_b^{N-N'}(\{\sigma \leq \epsilon\})$,and $R_A^{N'},\ R_B^{N'}$ are terms involving only $\nu_i$-derivatives of $R_A$, $R_B$ of order at most $N'$, and of $A$, $B$ of order at most $N'-1$.

Thus we can use \cref{thm:C5:initialvzs} provided the data $A_{N'}$ and $B_{N'}$ have sufficient regularity. We use a similar argument to the above. Indeed, by assumption $A_{N'-1}, B_{N'-1} \in v^{k-1}W_b^{N-(N'-1)}(\{\sigma = \epsilon\})$, and so all components $\nu_i A_{N'-1}, \ \nu_i A_{N'-1} \in v^{k-1}L^{\infty}(\{\sigma = \epsilon\})$ for $i \geq 1$. If $i=0$, one uses \eqref{eq:C5:induct3} for $A_{N'-1}$ and $B_{N'-1}$ to write $\nu_0A_{N'-1}$ and $\nu_0B_{N'-1}$ in terms of $A_{N'-1}$, $\nu_1A_{N'-1}$ and $B_{N'-1}$, respectively, to establish that $\nu_0 A_{N'-1}, \ \nu_0 B_{N'-1} \in v^{k-1}L^{\infty}(\{\sigma = \epsilon\})$, too.
\end{proof}

We argue polyhomogeneity at $\mathbf{sgf} = \{\sigma = 0\}$ in two steps. In the first, we show partial polyhomogeneity, propagating the polyhomogeneity at $\mathbf{frf} = \{v = 0\}$ from $\{\sigma \geq \epsilon\}$ given by \cref{thm:C5:prop3i} to all of $\{\sigma \geq 0\}$, and then improving this to full polyhomogeneity by treating the behaviour as $v \to 0$ essentially as a parameter.
\begin{prop}\label{thm:C5:prop3iihelper}It holds that $F', \ \tilde{H} \in \sigma^{-1}\phgi{(\{v=0\})}{(-1)}(\{\sigma \leq \epsilon\})$.\end{prop}
\begin{proof}In light of \cref{thm:C5:prop3i}:(i), we already know that for all $\epsilon > 0$ sufficiently small,
\begin{equation}\label{eq:prop3iihelperstart}F', \ \tilde{H} \in \sigma^{-1}\phgi{(\{v=0\})}{(-1)}(U\n \{\sigma \geq \epsilon/2\}).\end{equation}

Let us return to the coordinates $(\sigma,v,\theta)$ of \cref{thm:C5:corvzs}, and let $\nu = \nu_1 = v\pa_v$ in these coordinates. Observe carefully that this is \emph{not} $v\pa_v$ in the original coordinates, since this vector field is not tangent to the level sets of $\sigma$. Using that $\zeta\pa_\zeta - v\pa_v = \beta\pa_\sigma - \nu$ and $\pa_\zeta = \pa_\sigma$, we may rewrite \eqref{eq:C5:PPvz} as
\begin{align}
\label{eq:C5:PPvzrw}
\begin{split}
\sigma(\pa_\sigma - \beta^{-1}\nu)\Phi &= \frac{1}{2}(1+\sigma C^{\infty}(U))\Phi - \frac{1}{2}(1+\sigma C^{\infty}(U))\Psi\\
\sigma\pa_\sigma \Psi &=-\frac{1}{2}(1+\sigma C^{\infty}(U))\Phi + \frac{1}{2}(1+\sigma C^{\infty}(U))\Psi.
\end{split}
\end{align}

It suffices to show that for all $k \in \N_0$, $a_k(\nu)\Phi, a_k(\nu)\Psi \in v^{k-1}\mathcal A( \{\sigma \leq \epsilon\})$. We will prove this by induction on $k$. First, observe that, because of \eqref{eq:prop3iihelperstart}, $a_k(\nu)\Phi, a_k(\nu)\Psi \in v^{k-1}\mathcal A(\{\sigma = \epsilon\})$. The case $k=0$ is handled by \cref{thm:C5:corvzs}. Recall that from \eqref{eq:C5:fhve}, $vf$ and $v\tilde{h}$ are both smooth. Now for higher $k$, we can commute $a_k(\nu)$ through \eqref{eq:C5:PPvzrw} and use \cref{thm:C5:product} to see that $a_k(\nu)\Phi,a_k(\nu)\Psi$ satisfy \eqref{eq:C5:PPvzrw}, except with a right-hand side involving terms coming from the commutator of $a_k(\nu)$ with the $\sigma C^{\infty}$ error terms and the commutator with $\beta^{-1}\nu$. The extra terms resulting from the first commutator are of the form
\[\sigma(p_{k-k'}(\nu)\phi)(a_{k'}(\nu) S),\]
 where $\phi$ is smooth, $S$ is one of $\Phi,\Psi$, and $0 \leq k' < k$. The extra terms coming from the commutator with $\beta^{-1}\nu$ are of the form
\[\sigma(p_{k-k'}(\nu)\beta^{-1})(a_{k'}(\nu) \nu\Phi).\]
Since $\phi, \beta^{-1}$ is smooth, the first factor of both is in $v^{k-k'}C^{\infty}(\{\sigma \leq \epsilon\})$. By induction, the second factors are both in $v^{k'-1}\mathcal A(\{\sigma \leq \epsilon\})$ (since applying $\nu$ does not affect this estimate). Thus, the error terms are in $\sigma v^{k-1}\mathcal A(U)$, and we use \cref{thm:C5:corvzs} to conclude.
\end{proof}

Now we can prove part (i) of \cref{thm:C5:prop3ii}.
\begin{proof}[Proof of \cref{thm:C5:prop3ii}:(i).] It suffices to prove $\Phi,\Psi \in \mathcal A_{\mathrm{phg},(\{v=0\},\{\sigma = 0\})}^{(-1,0')}(\{\sigma \leq \epsilon\})$. We continue to work in $(\sigma,v,\theta)$ coordinates (valid for $\epsilon$ small). Define $A$, $B$ via $B-A = \Phi$, $B+A = \Psi$. From \cref{thm:C5:prop3iihelper}, $A,B \in \phgi{(\{v = 0\})}{(-1)}(\{\sigma \leq \epsilon\})$. Using $A$ and $B$ diagonalizes to top order the coefficient of the zeroth order term in \eqref{eq:C5:PPvz}. Explicitly, $A$, $B$ satisfy
\begin{align}
\label{eq:C5:ABvz}
\begin{split}
\sigma\pa_\sigma A = \sigma L_{11}A + \sigma L_{12}B\\
(\sigma\pa_\sigma -1)B = \sigma L_{21}A + \sigma L_{22}B,
\end{split}
\end{align}
where $L_{ij} \in \Diff_b^1(\{\sigma \leq \epsilon\})$ ($i,j = 1,2$). 

We will exhibit a series expansions at $\{\sigma = 0\}$ for all $\ell$
\begin{align}
\label{eq:C5:inductiveclaim}
\begin{split}
A &= A_{0,0} + \sigma A_{1,0} + \sum_{i = 2}^{\ell-1} \sigma^i(\log\sigma A_{i,1} + A_{i,0}) + A_{[\ell]}\\
B &= \sum_{i=1}^{\ell-1} \sigma^i(\log\sigma B_{i,1} + B_{i,0}) + B_{[\ell]},\end{split}\end{align}
where $A_{i,j}, B_{i,j} \in \mathcal A_{\mathrm{phg},(\{v = 0\})}^{(-1)} (\{\sigma = 0\})$, and $A_{[\ell]}, B_{[\ell]} \in \sigma^{\ell^-}\mathcal A_{\mathrm{phg},(\{v=0\})}^{(-1)}(\{\sigma \leq \epsilon\})$. 

One this is proven, we may use \cref{thm:C5:melrose} to conclude that \[A,B \in A_{\mathrm{phg},(\{v=0\},\{\sigma = 0\})}^{(-1,0')}(\{\sigma \leq \epsilon\}).\] Indeed, let $b_k(X)$, $k= 0,1,2,\ldots$ be the sequence of polynomials
\[1, X, X(X-1)^2, X(X-1)^2(X-2)^2, \ldots.\]
Then for all $k$ \[b_k(\sigma\pa_\sigma)A, \ b_k(\sigma\pa_\sigma)B \in \sigma^{k}\phgi{(\{v = 0\})}{(-1)}(\{\sigma \leq \epsilon\}),\] and so for all $j,k$, denoting by $\nu = v\pa_v$ (in $(\sigma,v,\theta)$ coordinates),
\[a_j(\nu)b_k(\sigma\pa_\sigma)A, \ a_j(\nu)b_k(\sigma\pa_\sigma)B \in v^{j-1}\sigma^{k}\phgi{(\{v = 0\})}{(-1)}(\{\sigma \leq \epsilon\}).\] Alternatively, we could expand $A$, $B$ as series in powers $v^k$ with coefficients in $\mathcal A(\{v = 0\})$ and remainder in $v^N \mathcal A(\{\sigma \leq \epsilon\})$, and run almost the same argument as the one which we will use to establish \eqref{eq:C5:inductiveclaim} in order to establish that the coefficients and remainders have series expansions at $\{\sigma = 0\}$ analogous to \eqref{eq:C5:inductiveclaim}.

We will show \eqref{eq:C5:inductiveclaim} by induction on $\ell$, together with the claim that
\[A_{(\ell-1)} := A - A_{[\ell]}, B_{(\ell-1)} := B-B_{[\ell]}\] solve \eqref{eq:C5:ABvz} modulo $\sigma^{\ell^-}\mathcal A_{\mathrm{phg},(\{v=0\})}^{(-1)}(\{\sigma \leq \epsilon\})$.

The base case for our inductive step will be to do the cases $\ell=1,2$. Equation \eqref{eq:C5:ABvz} implies that $\sigma\pa_\sigma A \in \mathcal \sigma A_{\mathrm{phg},(\{v=0\})}^{(-1)}(\{\sigma \leq \epsilon\})$. Thus by \cref{thm:C5:phglem}, below, $A = A_{0,0} + A_{[1]}$, where $A_{0,0} \in \mathcal A_{\mathrm{phg},(\{v=0\})}^{(-1)}(\{\sigma = 0\})$, and $A_{[1]} \in \sigma\mathcal A_{\mathrm{phg},(\{v=0\})}^{(-1)}(\{\sigma \leq \epsilon\})$. It also implies that \[(\sigma\pa_\sigma -1)B \in \sigma\mathcal A_{\mathrm{phg},(\{v=0\})}^{(-1)}(\{\sigma \leq \epsilon\}) \subseteq \sigma^{1^-} \mathcal A_{\mathrm{phg},(\{v=0\})}^{(-1)}(\{\sigma \leq \epsilon\}),\] and so by \cref{thm:C5:phglem}, $B \in \sigma^{1^-} \mathcal A_{\mathrm{phg},(\{v=0\})}^{(-1)}(\{\sigma \leq \epsilon\})$. This proves the case $\ell = 1$.

Denote by $L_{ij}(0)$, ($i,j = 1,2$) the restriction of $L_{ij}$ to $\{\sigma = 0\}$. Then $L_{11}A_0 = L_{11}(0)A_0 + \sigma\mathcal A_{\mathrm{phg},(\{v=0\})}^{(-1)}(\{\sigma \leq \epsilon\}) $. Thus
\[\sigma\pa_\sigma (A_{[1]} - \sigma L_{11}(0)A_0) \in \sigma^{2^{-}}\mathcal A_{\mathrm{phg},(\{v=0\})}^{(-1)}(\{\sigma \leq \epsilon\}),\]
so since we already know $A_{[1]} \in \sigma\mathcal A_{\mathrm{phg},(\{v=0\})}^{(-1)}(\{\sigma \leq \epsilon\})$, \cref{thm:C5:phglem} implies that
\[A_{[1]} = \sigma L_{11}(0)A_0 + A_{[2]},\] where $A_{[2]} \in \sigma^{2^-}\mathcal A_{\mathrm{phg},(\{v=0\})}^{(-1)}(\{\sigma \leq \epsilon\})$. Similarly,
\[(\sigma\pa_\sigma -1)(B- \sigma\log\sigma L_{21}(0)A_0) \in \sigma^{2^-}\mathcal A_{\mathrm{phg},(\{v=0\})}^{(-1)}(\{\sigma \leq \epsilon\}),\] and so again by \cref{thm:C5:phglem}
\[B = \sigma\log\sigma L_{21}(0)A_0 + \sigma B_{1,0}+B_{[2]},\] where $B_{1,0} \in \mathcal A_{\mathrm{phg},(\{v=0\})}^{(-1)}(\{\sigma = 0\})$, and $B_{[2]} \in \sigma^{2^-}\mathcal A_{\mathrm{phg},(\{v=0\})}^{(-1)}(\{\sigma \leq \epsilon\})$. This completes the base case.

Now assume the claim for $\ell$. We prove it for $\ell + 1$. By the inductive hypothesis
\[\sigma\pa_\sigma A_{(\ell-1)} = \sigma L_{11}A_{(\ell-1)} + \sigma L_{12}B_{(\ell-1)} + \sigma^{\ell^-}\mathcal A_{\mathrm{phg},(\{v=0\})}^{(-1)}(\{\sigma \leq \epsilon\}).\] However, $A_{(\ell-1)}$, $B_{(\ell-1)}$ are finite series, so the remainder must be a finite series as well. Thus in fact the error is of the form $-\sigma^{\ell}\log\sigma R_{\ell,1} + -\sigma^{\ell}R_{\ell,0} + \sigma^{(\ell+1)^-}\mathcal A_{\mathrm{phg},(\{v=0\})}^{(-1)}(\{\sigma \leq 0\})$, where $R_{\ell,1} \in \mathcal A_{\mathrm{phg},(\{v=0\})}^{(-1)}(\{\sigma = 0\})$, and similarly for $B$. Since the pair $A,B$ solve the equations exactly, we deduce that
\begin{align*}
\sigma\pa_\sigma A_{[\ell]} &= \sigma^{\ell}\log\sigma R_{\ell,1} + \sigma^{\ell}R_{\ell,0} + \sigma L_{11}A_{[\ell]} + \sigma L_{12}B_{[\ell]}\\
&=\sigma^{\ell}\log\sigma R_{\ell,1} + \sigma^{\ell}R_{\ell,0} + \sigma^{(\ell+1)^-}\mathcal A_{\mathrm{phg},(\{v=0\})}^{(-1)}(\{\sigma \leq \epsilon\}),
\end{align*} and similarly for $B$, since induction we already know $A_{[\ell]} \in \sigma^{\ell^-}\mathcal A_{\mathrm{phg},(\{v=0\})}^{(-1)}$. Thus
\[\sigma\pa_\sigma (A_{[\ell]} - \ell^{-1} \sigma^{\ell}\log\sigma R_{\ell,1} - \ell^{-1} \sigma^{\ell} R_{\ell,0} + \ell^{-1}\sigma^{\ell} R_{\ell,1}) \in \sigma^{(\ell+1)^-}\mathcal A_{\mathrm{phg},(\{v=0\})}^{(-1)}(\{\sigma \leq \epsilon\}),\]
and so, setting, \[A_{\ell,1} = -\ell^{-1} R_{\ell,1},\ A_{\ell,0} = -\ell^{-1} R_{\ell,0} + \ell^{-1} R_{\ell,1} \in \mathcal A_{\mathrm{phg},(\{v=0\})}^{(-1)}(\{\sigma =0\}),\]
it follows from \cref{thm:C5:phglem} that
\[A_{[\ell]} = \sigma^{\ell}(\log\sigma A_{\ell,1} + A_{\ell,0}) + A_{[\ell+1]},\] where $A_{[\ell+1]} \in \sigma^{(\ell+1)^-}\mathcal A_{\mathrm{phg},(\{v=0\})}^{(-1)}(\{\sigma \leq \epsilon\})$. A similar thing is true for $B$. 
Thus setting $A_{(\ell)} = A_{(\ell-1)} + \sigma^{\ell}(\log\sigma A_{\ell,1} + A_{\ell,0})$ and similarly for $B$ completes the inductive step.
\end{proof}

\begin{lem}\label{thm:C5:phglem}Choose coordinates $(\sigma,v,\theta)$ on $\{\sigma \leq \epsilon\}$. Fix $k \in \Z$, and fix $x,y,z \in \R$, $x < y$, $y \neq z$, and suppose $u \in \sigma^x \phgi{(\{v = 0\})}{(k)}(\{\sigma \leq \epsilon\})$, and $u$ solves $(\sigma\pa_\sigma -z)u = \sigma^y w \in \sigma^y\phgi{(\{v = 0\})}{(k)}(\{\sigma \leq \epsilon\})$ with data $u(\epsilon) \in \phgi{(\{v = 0\})}{(k)}(\{\sigma = \epsilon\})$.

If  $y < z$, then in fact $u \in \sigma^y \phgi{(\{v = 0\})}{(k)}(\{\sigma \leq \epsilon\})$. If $y > z$, then we distinguish two sub-cases. If $x \leq z$, then $u = \sigma^z u_1 + u_2$, where $u_1 \in \phgi{(\{v = 0\})}{(k)}(\{\sigma = 0\})$, and $u_2 \in \sigma^y\phgi{(\{v = 0\})}{(k)}(\{\sigma \leq \epsilon\})$. If instead $x > z$, then the same conclusion holds with $u_1 \equiv 0$. 
\end{lem}
\begin{rk}The assumption $x < y$ is only needed to ensure that the conclusion is not vacuous, since in the case $x > y$, the conclusion of the lemma is weaker than the hypothesis $u \in \sigma^x \phgi{(\{v = 0\})}{(k)}(\{\sigma \leq \epsilon\})$.\end{rk}
\begin{rk}One could extend this lemma to deal with the case $y=z$ and even extend it to prove \cref{thm:C5:melrose}. Since this would involve some notational overhead, we content ourselves with only this lemma.\end{rk}
\begin{proof}Conjugating by $v^k\sigma^z$, we may assume $k=0$, $z = 0$. One has the representation formula
\begin{equation}\label{eq:phglem:eq1}u(\sigma) = \int_{\epsilon}^\sigma \tau^{y-1}w(\tau)\ d\tau + u(\epsilon)\end{equation}
(we suppress the variables $v$, $\theta$). 

Let us first assume $y > z$, i.e.\ $y > 0$.

Set
\begin{align*}
I &:= \int_{0}^\sigma \tau^{y-1}w(\tau)\ d\tau\\
J &:= \int_{0}^\epsilon \tau^{y-1}w(\tau)\ d\tau
\end{align*}
so that $u = I -J + u(\epsilon)$.

Since $y > 0$, one easily checks that $ I \in \sigma^yL^{\infty}(\{\sigma \leq \epsilon\})$. The vector fields
$v\pa_v$ or $\pa_{\theta^i}$ essentially commute with the integral in the definition of $I$, so all $v\pa_v$ and $\pa_{\theta^i}$ derivatives of $I$ are also in $\sigma^yL^{\infty}(\{\sigma \leq \epsilon\})$. From this and the fact that $\sigma\pa_\sigma I = \sigma^y w$, it follows that $I \in \sigma^y\mathcal A(\{\sigma \leq \epsilon\})$. The differential operator $p_j(v\pa_v)$ also commutes through the integral for all $j$, and thus we conclude $I \in \sigma^y\phgi{(\{v = 0\})}{(0)}(\{\sigma \leq \epsilon\})$.

Since $J$ does not depend on $\sigma$, a similar argument shows that we may interpret $J \in \phgi{(\{v = 0\})}{(0)}(\{\sigma = 0\})$. Since $u(\epsilon)$ does not depend on $\sigma$, we may also interpret $u(\epsilon) \in \phgi{(\{v = 0\})}{(0)}(\{\sigma = 0\})$.

Thus $u = \sigma^y\phgi{(\{v = 0\})}{(0)}(\{\sigma \leq \epsilon\}) + \phgi{(\{v = 0\})}{(0)}(\{\sigma = 0\})$. If $x > z = 0$, then the second term cannot be present, since $u \in \sigma^x\phgi{(\{v = 0\})}{(0)}(\{\sigma \leq \epsilon\})$, and we derive $u \in \sigma^y\phgi{(\{v = 0\})}{(0)}(\{\sigma \leq \epsilon\})$. Otherwise, we obtain the conclusion if $x \leq z$.

Now assume $y < z$, i.e.\ $y < 0$. Certainly \[u(\epsilon) \in \phgi{(\{v = 0\})}{(0)}(\{\sigma = \epsilon\}) \subseteq \sigma^0 \phgi{(\{v = 0\})}{(0)}(\{\sigma \leq \epsilon\}) \subseteq \sigma^y \phgi{(\{v = 0\})}{(0)}(\{\sigma \leq \epsilon\}).\] We focus our attention on
\[I_2 := \int_{\epsilon}^{\sigma} \tau^{y-1}w(\tau)\ d\tau,\]
so that $u = I_2 + u(\epsilon)$. 
We show also that $I_2 \in \sigma^y \phgi{(\{v = 0\})}{(0)}(\{\sigma \leq \epsilon\})$. Certainly $|I_2| \lesssim \sigma^y + 1 \lesssim \sigma^y$. Now commuting through first $v\pa_v,\pa_{\theta^i}$, and then applying $\sigma\pa_\sigma$ as above, it follows that
\[I_2 \in \sigma^y\mathcal A(\{\sigma \leq \epsilon\}).\] Commuting through $p_j(v\pa_v)$ for all $j$ shows that $I_2 \in \sigma^y \phgi{(\{v = 0\})}{(0)}(\{\sigma \leq \epsilon\})$. Thus
\[u = I_2 + u(\epsilon) \in \sigma^y \phgi{(\{v = 0\})}{(0)}(\{\sigma \leq \epsilon\}),\] as desired.
\end{proof}

We now turn our attention to proving part (ii) of \cref{thm:C5:prop3ii}. Fix $k, \ell \in \R$, and $\ell \neq x \in \R$. Consider the equations for a vector-valued function $w$:
\begin{equation}\label{eq:C5:wvzsi} (\sigma\pa_\zeta - x) w - \sigma E w = S\end{equation}
and
\begin{equation}\label{eq:C5:wvzsii} (\sigma(\zeta\pa_\zeta -v\pa_v) -\beta x)w - \sigma E w = S,\end{equation}
where in both cases $S \in v^k\sigma^{\ell}L^{\infty}(\{\sigma \leq \epsilon\})$ and $E \in L^{\infty}(\{\sigma \leq \epsilon\})$.
We start with a general proposition. 
\begin{prop}\label{thm:C5:omegavzshelper}Suppose $w \in C^1(\{\sigma \geq \epsilon, v > 0\})$ solves \eqref{eq:C5:wvzsi} or \eqref{eq:C5:wvzsii} with data $w|_{\{\sigma = \epsilon\}} \in v^k L^{\infty}(\{\sigma = \epsilon\})$. Then $w \in v^k\sigma^{\min(x,\ell)} L^{\infty}(\{\sigma \leq \epsilon\})$ with the bound
\[\norm{w}_{v^k \sigma^{\min(x,\ell)}L^{\infty}(\{\sigma \leq \epsilon\})} \leq C(\norm{w|_{\{\sigma = \epsilon\}}}_{v^kL^{\infty}(\{\sigma = \epsilon\})} + \norm{S}_{v^k \sigma^{\ell}L^{\infty}(\{\sigma \leq \epsilon\})}),\]
where 
\[C = C(\norm{E}_{L^{\infty}(\{\sigma \leq \epsilon\})},\sup \beta^{-1},k,\ell).\]

Now fix $N \geq 1$. If instead $w \in C^{N+1}(\{\sigma \leq \epsilon, v > 0\})$ and $S \in v^k\sigma^{\ell} W_b^N(\{\sigma \leq \epsilon\})$, $E \in W_b^N(\{\sigma \leq \epsilon\})$, and $w|_{\{\sigma = \epsilon\}} \in v^kW_b^N(\{\sigma = \epsilon\})$, then $w \in v^k \sigma^{\min(x,\ell)}W_b^N(\{\sigma \leq \epsilon\})$ with an analogous bound. In particular, if the hypotheses hold for all $N$, then $w \in v^{k}\sigma^{\min(x,\ell)} \mathcal A(\{\sigma \leq \epsilon\})$.\end{prop}
\begin{proof}
The cases of $w$ solving either of \eqref{eq:C5:wvzsi} or \eqref{eq:C5:wvzsii} is similar, so we will only treat the (slightly) harder case of \eqref{eq:C5:wvzsii}. Conjugating by $\sigma^x v^k$, we may assume $k = x = 0$.

Let $\phi$ solve
\[\sigma(\zeta\pa_\zeta-v\pa_v)\phi = -\sigma E,\]
with data $\phi(\epsilon) = 1$. Then
\[\phi = \exp\left(\int_{\mathbf{j}}^\zeta - E(v\zeta/t,t)\ dt/t\right).\]
Since 
\[\int_{\mathbf{j}}^\zeta |E(v\zeta/t,t)|\ dt/t = \int_{\mathbf{j}}^\zeta |E(v\zeta/t,t)|\beta/\beta\ dt/t \lesssim \int_{\mathbf{j}}^\zeta \beta\ dt/t \lesssim 1\]
by \cref{thm:C5:integration1}, it follows that $1 \lesssim \phi \lesssim 1$.
We use the method of integrating factors to see that
\[(\phi w)(v,\zeta) = (\phi w)(v\zeta/\mathbf{j},\mathbf{j}) +\int_{\mathbf{j}}^\zeta (\phi S)(v\zeta/t,t)\ dt/t,\]
and thus
\[|w(v,\zeta)| \lesssim |w(v\zeta/\mathbf{j},\mathbf{j})| + \int_{\zeta}^{\mathbf{j}} \sigma^{-1}|S(v\zeta/t,t)|\beta\ dt/t.\]

To show the first part of the proposition, it therefore suffices to bound the integral. By assumption
\[|S(v\zeta/t,t)| \leq \norm{S}_{ \sigma^{\ell}L^{\infty}(\{\sigma \leq \epsilon\})}\sigma(v\zeta/t,t)^{\ell-1}.\] Since $\ell \neq 0$, it follows from \cref{thm:C5:integration1} that 
\[\int_{\zeta}^{\mathbf{j}} \sigma(v\zeta/t,t)^{\ell-1}\beta\ dt/t \leq \epsilon^{\ell}-\sigma^{\ell} \leq 1 - \sigma^{\ell}.\]

Putting this all together completes the proof of the first part of the proposition.

The proof of higher regularity uses the first part in the same manner as the proof of \cref{thm:C5:corvzs} uses \cref{thm:C5:initialvzs}, so we only provide an outline of the proof, and focus on the case $N=1$ as an example.

Let us return to $(\sigma,v,\theta)$ coordinates, and recall the vector fields $\nu_0,\ldots, \nu_3$ as in the proof of \cref{thm:C5:corvzs}. Then commuting through $\nu_i$,
\begin{equation}\label{eq:C5:imtired}\sigma(\zeta\pa_\zeta-v\pa_v)\nu_i w - \sigma (E \nu_i + (\delta_{i0}-\nu_i\log \beta) \nu_1)w = \nu_i S + \nu_i(\sigma E) w.\end{equation}
All the terms on the right-hand side are in $\sigma^{\min(0,\ell)}L^{\infty}(\{\sigma \leq \epsilon\}$. \hspace{-2pt}Writing $w_1 = (\nu_0w,\ldots,\nu_3w)$, we may combine \eqref{eq:C5:imtired} over all $i$ to find an equation for $w_1$:
\[\sigma(\zeta\pa_\zeta - v\pa_v)w_1 = \sigma E^1w_1 + S_1,\] where $E^1 \in L^{\infty}(\{\sigma \leq \epsilon\}))$ and $S_1 \in v^{k-1}\sigma^{\min(1,\ell)}$. The first part of the proposition, using \eqref{eq:C5:wvzii} to treat the data, now gives $w_1 \in\sigma^{\min(\ell,0)}L^{\infty}(\{\sigma \leq \epsilon\}))$, which shows the $N=1$ case.
\end{proof}

We will argue polyhomogeneity $ \{\sigma = 0\}$ in two steps. First we propagate the polyhomogeneity at $\{v = 0\}$ all the way to $\{\sigma = 0\}$, and then prove full polyhomogeneity.

We may rewrite \eqref{eq:C5:omegavz} as
\begin{equation}\label{eq:C5:omegavzs}
(\sigma\pa_\zeta -1)\sigma(\zeta\pa_\zeta - v\pa_v)\omega = \frac{v^2\sigma^2}{16}f\tilde{h} - \frac{v^2\sigma^2}{4}F'\tilde{H}R^2 - \frac{v\sigma^2}{4}f =: S.\end{equation}

Observe that by \cref{thm:C5:prop3}:(i) and \eqref{eq:C5:fhvz}, $S \in \phgi{(\{v = 0\},\{\sigma = 0\})}{(0,0'')}(U)$.

\begin{prop}\label{thm:C5:omegaphgI}It holds that
 \[\omega \in \sigma^{0^-}\phgi{(\{v = 0\})}{(0)}(\{\sigma \leq \epsilon\}).\]
\end{prop}
\begin{proof}
In light of \cref{thm:C5:prop3i}:(ii), we already know that
\[\omega \in \phgi{(\{v = 0\})}{(0)}(\{\sigma \geq \epsilon/2\}).\]
Let us return to $(\sigma,v,\theta)$ coordinates and let $\nu$ be the vector field is $v\pa_v$. Recall that $\pa_\zeta = \pa_\sigma$ and $\zeta\pa_\zeta - v\pa_v = \beta\pa_\sigma - \nu$.

 It suffices to show that for all $k \in \N_0$
\[p_k(\nu)\omega \in \sigma^{0^-}\mathcal A(\{\sigma \leq \epsilon\}).\] We will use induction on $k$ in the same way as in the proof of \cref{thm:C5:prop3iihelper}.

Let us start with $k=0$. We may apply \cref{thm:C5:omegavzshelper} to \eqref{eq:C5:omegavzs} twice to conclude that \[\omega \in \sigma^{0^-}\mathcal A(\{\sigma \leq \epsilon\}).\]
For higher $k$, we can commute $p_k(\nu)$ through \eqref{eq:C5:omegavzs} and use \cref{thm:C5:product} to see that $p_k(\nu)\omega$ satisfies an equation of the form
\[(\sigma\pa_\zeta-1) \sigma(\zeta\pa_\zeta - v\pa_v)p_k(\nu)\omega = p_k(\nu)S + R,\]
where $R$ is an error term consisting of sums of terms of the form
\[(\sigma\pa_\zeta-1)(p_{k-k'}(\nu)(\beta))(p_{k'}(\nu)(\sigma\pa_\sigma \omega))\] for $0 \leq k' < k$. By induction this is in $v^k\sigma^{0-}\mathcal A(\{\sigma \leq \epsilon\})$, so we may again apply \cref{thm:C5:omegavzshelper} twice and conclude.\end{proof}

We can now prove part (ii) of \cref{thm:C5:prop3ii}.
\begin{proof}[Proof of \cref{thm:C5:prop3ii}:(ii).]
We continue to work in $(\sigma,v,\theta)$ coordinates. We may write \eqref{eq:C5:omegavzs} schematically as
\begin{equation}\label{eq:C5:omegavzsc}(\sigma\pa_\sigma-1)\sigma\pa_\sigma \omega = \sigma L \omega + S/\beta,\end{equation} where $L \in \Diff^2_b(\{\sigma \leq \epsilon\})$.

We will exhibit a series expansion at $\{\sigma = 0\}$ for all $\ell$
\[\omega = \omega_{0,1}\log \sigma + \omega_{0,0} + \sum_{i=1}^{\ell-1}\sum_{j=0}^2\sigma^i\log^j \sigma \omega_{i,j} + \omega_{[\ell]},\]
where $\omega_{i,j} \in \phgi{(\{v = 0\})}{(0)}(\{\sigma = 0\})$ and $\omega_{[\ell]} \in \sigma^{\ell^-}\phgi{(\{v = 0\})}{(0)}(\{\sigma = 0\})$.
We may expand for all $\ell$
\begin{equation}\label{eq:C5:Rseries}S/\beta = S_{0,0} + \sum_{i=1}^{\ell-1}\sum_{j=0}^2\sigma^i\log^j \sigma S_{i,j}+ S_{[\ell]},\end{equation}
where $S_{i,j} \in \phgi{(\{v = 0\})}{(0)}(\{\sigma = 0\})$ and $S_{[\ell]} \in \sigma^{\ell^-}\phgi{(\{v = 0\})}{(0)}(\{\sigma = 0\})$ (and $S_{1,2} \equiv 0$).

We show this by induction on $\ell$, together with the claim that $\omega_{(\ell-1)} := \omega - \omega_{[\ell]}$ solves \eqref{eq:C5:omegavzs} modulo $\sigma^{\ell^-}\phgi{(\{v = 0\})}{(0)}(\{\sigma \leq \epsilon\})$. The base case is $\ell = 1$. For this case observe that \cref{thm:C5:omegaphgI} and \eqref{eq:C5:omegavzsc} implies that
\[(\sigma\pa_\sigma-1)(\sigma\pa_\sigma)(\omega +\log \sigma S_{0,0}) \in \sigma^{1^-}\phgi{(\{v = 0\})}{(0)}(\{\sigma \leq \epsilon\}),\] and so by \cref{thm:C5:phglem},
\[\omega = \log \sigma S_{0,0} + \omega_{0,0} + \sigma^{1^-}\phgi{(\{v = 0\})}{(0)}(\{\sigma \leq \epsilon\}).\]
The inductive step follows in almost the same way as the inductive step in the proof of \cref{thm:C5:prop3ii}:(i), finding the $\omega_{\ell,j}$ such that
\[(\sigma\pa_\sigma-1)(\sigma\pa_\sigma)\left(\omega_{[\ell]} - \sigma^{\ell}\sum_{j=0}^2 \log^j \sigma \omega_{\ell,j}\right) \in \sigma^{(\ell+1)^-}\phgi{(\{v = 0\})}{(0)}(\{\sigma \leq \epsilon\})\] and then invoking \cref{thm:C5:phglem}. We omit the details.
\end{proof}

\todo{check reference in texorpdfstring}\section{The proof of \texorpdfstring{\cref{thm:C5:prop4}}{proposition~5.2.6}} We work in $(v,\zeta,\alpha)$ coordinates to cover $U_4$. Recall that $\zeta|\alpha| \leq \delta_i$ on $U_4$. Set \begin{equation}\label{eq:C5:sigmadefii}\sigma = 1-\kappa(1+vK)|\alpha|,\end{equation} a bdf of $\mathbf{sgf}$. It will be convenient to rescale $\tilde{H}$ by setting $\bar{H} = \zeta\tilde{H}$. \eqref{eq:C5:FH} takes the form
\begin{align}
\label{eq:C5:FHa}
\begin{split}
(\zeta\pa_\zeta - (\alpha^1\pa_{\alpha^1} + \alpha^2 \pa_{\alpha^2}))F' + &\frac{v f}{4}\bar{H} + \frac{v\zeta}{4}\tilde{h}F' + F' = 0\\
(\zeta\pa_\zeta-v\pa_v - (\alpha^1\pa_{\alpha^1} + \alpha^2 \pa_{\alpha^2})\bar{H} - &\frac{vf}{4}\bar{H} - \frac{v\zeta}{4}\tilde{h}F' - F' - \bar{H} = 0
\end{split}
\end{align}
and
\eqref{eq:C5:omega} takes the form
\begin{equation}\label{eq:C5:omegavza} (\zeta\pa_\zeta- (\alpha^1\pa_{\alpha^1} + \alpha^2 \pa_{\alpha^2})(\zeta\pa_\zeta - v\pa_v- (\alpha^1\pa_{\alpha^1} + \alpha^2 \pa_{\alpha^2}))\omega - \frac{v^2\zeta}{16}f\tilde{h} + \frac{v^2\zeta^2}{4}F'\bar{H}|\alpha|^2 + \frac{v\zeta}{4}f = 0.\end{equation}
From \cref{thm:C5:prop1}, we know smoothness away from the boundaries $\{v=0\}$, $\{\sigma = 0\}$, $\{\zeta = 0\}$, and from \cref{thm:C5:prop3} we know that
\begin{gather*}F',\bar{H} \in \phgi{(\mathbf{frf})}{(-1)}(U_4 \n \{\zeta_0 \leq \zeta \leq 1/\lambda_0, \ 0 \leq v\zeta \leq 1\})\\
\omega \in \phgi{(\mathbf{frf})}{(0)}(U_4 \n \{\zeta_0 \leq \zeta \leq 1/\lambda_0, \ 0 \leq v\zeta \leq 1\})\end{gather*}
Combining these shows
\begin{align}\label{eq:C5:initialvaluea}\begin{split}F',\bar{H} \in \phgi{(\mathbf{frf})}{(-1)}(U_4 \n \{\zeta_0 \leq \zeta \leq 1/\lambda_0\})\\
\omega \in \phgi{(\mathbf{frf})}{(0)}(U_4 \n \{\zeta_0 \leq \zeta \leq 1/\lambda_0\}).\end{split}\end{align}
Because of the appearance of the radial vector field $\alpha^1\pa_{\alpha^1} + \alpha^2\pa_{\alpha^2}$ in \eqref{eq:C5:FHa}, it will be convenient for technical purposes to blow up and set $W' = [U_4,\{\alpha = 0\}]$, and cover it by projective coordinates $r = \pm \alpha^i$, $u = \alpha^j/\alpha^i$, for $i \neq j$. To be explicit, we will assume without loss of generality in this subsection that $r = \alpha^1$ and $u = \alpha^2/\alpha^1$. We may assume that in such a chart $|u| \leq 2$. Let $U$ denote the range of such a chart. Since $\sigma \geq 0$ implies that $|\alpha|$ is bounded, $U$ is compact.
On $U$, $\sigma = 1-\kappa(1+vK)r \langle u \rangle$, where $\langle u\rangle = \sqrt{1+u^2}$ denotes the Japanese brackets. On $U$, we may rewrite \eqref{eq:C5:FHa} as
\begin{align}
\label{eq:C5:FHU}
\begin{split}
(\zeta\pa_\zeta - r\pa_r)F' + \frac{v f}{4}\bar{H} + \frac{v\zeta}{4}\tilde{h}F' + F' &= 0\\
(\zeta\pa_\zeta-v\pa_v - r\pa_r)\bar{H} - \frac{vf}{4}\bar{H} - \frac{v\zeta}{4}\tilde{h}F' - F' - \bar{H} &= 0.
\end{split}
\end{align}
and \eqref{eq:C5:omegavza} as
\begin{equation}\label{eq:C5:omegaU}(\zeta\pa_\zeta- r\pa_r)(\zeta\pa_\zeta - v\pa_v- r\pa_r)\omega - \frac{v^2\zeta}{16}f\tilde{h} + \frac{v^2\zeta^2}{4}F'\bar{H}r^2\langle u \rangle^2 + \frac{v\zeta}{4}f = 0.\end{equation}

Translating \eqref{eq:C5:initialvaluea} to the blowup, we have
\begin{align}\begin{split}\label{eq:C5:boredothis}F',\bar{H} \in \phgi{(\mathbf{frf},\{r= 0 \})}{(-1,0)}(U \n \{\zeta_0 \leq \zeta \leq 1/\lambda_0\})\\
\omega \in \phgi{(\mathbf{frf},\{r= 0 \})}{(0,0)}(U \n \{\zeta_0 \leq \zeta \leq 1/\lambda_0\}).\end{split}\end{align}
Thus, we will treat \eqref{eq:C5:FHU} and \eqref{eq:C5:omegaU} as initial-value problems with data for $F'$, $\bar{H}$, $\omega$ given on $\{\zeta = \zeta_0\}$. \begin{figure}[htbp]
\centering
\includegraphics{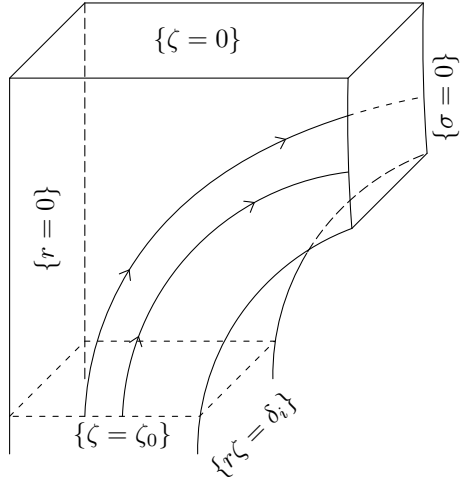}
\caption{A view of $U$ at some $u$ fixed. The horizontal axis is $r$, the vertical axis is $\zeta$, and the axis into the page is $v$. Indicated is a sample backwards-directed integral curve of $\zeta\pa_\zeta - r\pa_r$ and $\zeta\pa_\zeta - v\pa_v - r\pa_r$.}
\end{figure}

For this section, we will primarily work in $U$. Since the blowdown $W' \to U_4$ is a diffeomorphism away from $\alpha = 0$, any statement about $U$ which is true away from $\alpha = 0$ will carry over automatically to $U_4$. We will need to work a little harder to pass from statements about $U$ which are true near $\alpha = 0$ to $U_4$.
In $U$, the variable $u$ acts like a parameter, so we will usually ignore it unless we are explicitly differentiating $\pa_u$.
Propagating polyhomogeneity at $\mathbf{frf}$ to all of $\{\sigma > 0\}$ in $U_4$ and then establishing polyhomogeneity at $\mathbf{sgf}$ will require different arguments. Thus, we break \cref{thm:C5:prop4} into three propositions:
\begin{prop}\label{thm:C5:prop4i}
For all $\epsilon > 0$ sufficiently small depending only on $\zeta_0$ and the specific choice \eqref{eq:C5:sigmadefii} of the bdf $\sigma$ of $\mathbf{sgf}$:
\begin{romanumerate}
\item $F', \bar{H} \in \phgi{(\{v = 0\},\{r=0\},\{\zeta = 0\})}{(-1,0,-1)}(\{\sigma \geq \epsilon\} \n U);$
\item $\omega \in \phgi{(\{v = 0\},\{r=0\},\{\zeta = 0\})}{(0,0,0)}(\{\sigma \geq \epsilon\} \n U).$
\end{romanumerate}
\end{prop}
\begin{prop}\label{thm:C5:prop4ii}
For all $\epsilon > 0$ sufficiently small depending only on $\zeta_0$ and the specific choice \eqref{eq:C5:sigmadefii} of the bdf $\sigma$ of $\mathbf{sgf}$:
\begin{romanumerate}
\item $F', \bar{H} \in \phgi{(\{v = 0\},\{\zeta = 0\})}{(-1,-1)}(\{\sigma \geq \epsilon\} \n U_4);$
\item $\omega \in \phgi{(\{v = 0\},\{\zeta = 0\})}{(-1,0)}(\{\sigma \geq \epsilon\} \n U_4).$
\end{romanumerate}
\end{prop}
and
\begin{prop}\label{thm:C5:prop4iii}For all $0 < \epsilon < 1$ sufficiently small depending only on $\zeta_0$ and the specific choice \eqref{eq:C5:sigmadefii} of the bdf $\sigma$ of $\mathbf{sgf}$:
\begin{romanumerate}
\item $F', \bar{H} \in \phgi{(\{v = 0\},\{\zeta = 0\},\{\sigma = 0\})}{(-1,-1,-1')}(\{\sigma \leq \epsilon\} \n U);$
\item $\omega \in \phgi{(\{v = 0\},\{\zeta = 0\},\{\sigma = 0\})}{(0,0,0'')}(\{\sigma \leq \epsilon\} \n U).$
\end{romanumerate}
\end{prop}

We briefly mention how to derive \cref{thm:C5:prop4} from these propositions. As mentioned above, the blowdown $W' \to U_4$ is a diffeormoprhism away from $0$. Since $\alpha \neq 0$ near $\sigma = 0$, \cref{thm:C5:prop4iii} immediately implies that
\begin{align*}F', \bar{H} \in \phgi{(\{v = 0\},\{\zeta = 0\},\{\sigma = 0\})}{(-1,-1,-1')}(\{\sigma \leq \epsilon\} \n U_4)\\
\omega \in \phgi{(\{v = 0\},\{\zeta = 0\},\{\sigma = 0\})}{(0,0,0'')}(\{\sigma \leq \epsilon\} \n U_4).\end{align*}
 It is now clear (after rescaling $\bar{H}$ back to $\tilde{H}$) that \cref{thm:C5:prop4ii} and \cref{thm:C5:prop4iii} imply \cref{thm:C5:prop4}. However, as we mentioned, for technical reasons it is convenient to blow up, so we will first prove \cref{thm:C5:prop4i} and then use it to prove \cref{thm:C5:prop4ii}.

The splitting of \cref{thm:C5:prop4} into regions $\{\sigma \leq \epsilon\}$ and $\{\sigma \geq \epsilon\}$ makes sense in light of the hyperbolic nature of \eqref{eq:C5:FHU} and \eqref{eq:C5:omegaU}. Indeed, one may compute
\begin{align}
\label{eq:C5:sigders2}
\begin{split}
(\zeta\pa_\zeta - r\pa_r)\sigma &= \kappa(1+vK)r\langle u \rangle = 1-\sigma\\
(\zeta\pa_\zeta - v\pa_v - r\pa_r)\sigma &= \kappa(1+2vK + v^2\pa_v K)r\langle u \rangle\\
&= \frac{1+2vK + v^2\pa_v K}{1+vK}(1-\sigma)
\end{split}
\end{align}
If $\epsilon < 1$, then the first line is strictly positive on $\{\sigma \leq \epsilon\}$. Since $1+vK > 0$, and in the beginning of the previous section we established $1+2vK + v^2\pa_v K > 0$, the second line is also strictly positive on $\{\sigma \leq \epsilon\}$.

Thus $\sigma$ is increasing along the flows of $\zeta\pa_\zeta- r\pa_r$ and $\zeta\pa_\zeta -v\pa_v-r\pa_r$ for $\sigma \leq \epsilon$. In particular, $\{\sigma \geq \epsilon\}$ forms a domain of dependence for \eqref{eq:C5:FHvz} and what happens in the region $\{\sigma \leq \epsilon\}$ cannot affect what happens in the region $\{\sigma \geq \epsilon\}$ (as we are propagating backwards along the flow).

Let us now specify how small $\epsilon$ needs to be in \cref{thm:C5:prop4i}, \cref{thm:C5:prop4ii} and \cref{thm:C5:prop3ii}. As we've mentioned, we need to take $\epsilon  1$. For \cref{thm:C5:prop4i} and \cref{thm:C5:prop4ii}, $\epsilon$ should be small enough so $\epsilon < \inf \sigma|_{\{\zeta = \zeta_0\}}$ and $\inf \zeta|_{\{\sigma = \epsilon\}} \leq \zeta_0$. This ensures that every backwards-directed integral curve of $\zeta\pa_\zeta-r\pa_r$ or $\zeta\pa_\zeta - v\pa_v-r\pa_r$ which starts at $\{\zeta = \zeta_0\}$ eventually intersects $\{\sigma = \epsilon\}$, and conversely, every forwards-directed integral curve starting at $\{\sigma = \epsilon\}$ intersects $\{\zeta = \zeta_0\}$. Since $\epsilon < 1$, every forwards-directed integral curve starting at $\{\sigma = 0\}$ eventually intersects $\{\sigma = \epsilon\}$, so the latter is a Cauchy hypersurface of \eqref{eq:C5:FHU} and \eqref{eq:C5:omegaU}.

We will prove \cref{thm:C5:prop4i}, \cref{thm:C5:prop4ii} and \cref{thm:C5:prop4iii} over the next two subsections. The first subsection will be devoted to proving \cref{thm:C5:prop4i} and \cref{thm:C5:prop4ii}. The proof of \cref{thm:C5:prop4i} will follow the same outline as the proofs of \cref{thm:C5:prop2} and \cref{thm:C5:prop3i}, with a few additional difficulties. The second subsection will be devoted to proving \cref{thm:C5:prop4iii}. Because the proof is almost identical to that of \cref{thm:C5:prop3ii}, since in the latter $v$ is treated mostly as a parameter, and in the former $v,r$ are treated as parameters, we omit most of the details.

\subsection{Behaviour away from \texorpdfstring{$\{\sigma = 0\}$}{\{sigma = 0\}}.}
For most of this subsection, we work in $U$, only returning to $U_4$ when we prove \cref{thm:C5:prop4ii}.

Let us first work on proving part (i) of \cref{thm:C5:prop4i} and \cref{thm:C5:prop4ii}. The proof of the former is essentially a hybrid between the proofs of \cref{thm:C5:prop2}:(i) and \cref{thm:C5:prop3i}:(i).

As in the proof of \cref{thm:C5:prop3i}:(i), for technical purposes, we will need to localize to domains of dependence of \eqref{eq:C5:FHU}. For our current purposes, we need to redefine:
\begin{defn}A \emph{domain of dependence} for \eqref{eq:C5:FHU} is subset $\mathcal D \subseteq U$ for which every forwards-directed integral curve of $\zeta\pa_\zeta-r\pa_r$ or $\zeta\pa_\zeta-v\pa_v-r\pa_r$ which starts in $\mathcal D$ does not exit $\mathcal D$ until it intersects $\{\zeta = \zeta_0\}$.\end{defn} 
Notice that we are propagating backwards along the flow of $\pa_\zeta$ and $\zeta\pa_\zeta - v\pa_v$, so the backwards domain of dependence is forwards along the flow. We also define:
\begin{defn}For $p \in U$, the \emph{backwards domain of dependence} of $p$, denoted $B(p)$, is the smallest domain of dependence containing $p$.\end{defn}
One may explicitly describe $B(p)$. If $p = (v_1,\zeta_1,r_1)$, then
\[B(p) = \{(v,\zeta,r) \: \zeta_0 \leq \zeta \leq \zeta_1,\ r\zeta = r_1\zeta_1,\ v_1\zeta_1/\zeta \leq v \leq v_1\}.\] Indeed, one can explicitly check that $B(p)$ is a domain of dependence, and is the smallest set containing the forwards-directed integral curve of $\zeta\pa_\zeta - r\pa_r$ starting from $p$ and all the forwards-directed integral curves of $\zeta\pa_\zeta - v\pa_r - r\pa_r$ starting at any point of the previous curve. As remarked above, observe then $\sigma(p) \geq \epsilon$ implies $\sigma(q) \geq\epsilon$ for all $q \in B(p)$, since every backwards-directed integral curve of $\zeta\pa_\zeta-r\pa_r$ and $\zeta\pa_\zeta - v\pa_v-r\pa_r$ starting at a point $q'$ with $\sigma(q') < \epsilon$ remains in the region $\{\sigma < \epsilon\}$.

We will need to to expand the type of equations we consider to include an inhomogeneous right-hand side. Such equations will arise when commuting our original equation with derivatives which we will use to establish regularity. Let us consider for $k,j,\ell \in \N_0$ the system
\begin{align}
\label{eq:C5:mea}
\begin{split}
(\zeta\pa_\zeta - r\pa_r)A + a_{11}B +a_{12}A + A = R_A\\
(\zeta\pa_\zeta-v\pa_v - r\pa_r ) B +a_{21}B + a_{22}A -A-B = R_B,\end{split}\end{align}
where $a_{ij} \in r^2 L^{\infty}(\{\sigma \geq \epsilon\})$, \[R_A \in v^{\ell-1}r^{k+\max(j,1)}\zeta^{k-1}L^{\infty}(\{\sigma \geq \epsilon\}), \ \ R_B \in v^{\ell-1}r^{k+j}\zeta^{k-1}L^{\infty}(\{\sigma \geq \epsilon\}),\]
and $A|_{\{\zeta = \zeta_0\}},B|_{\{\zeta = \zeta_0\}} \in v^{\ell-1}r^{k+j}L^{\infty}(\{\zeta = \zeta_0\})$.
Set
\begin{align*}
\mathcal E &= \sum_{ij} \norm{a_{ij}}_{r^2L^{\infty}(\{\sigma \geq \epsilon\})}\\
\mathcal I &= \norm{A|_{\{\zeta = \zeta_0\}}}_{v^{\ell-1}r^{k+j}L^{\infty}(\{\zeta = \zeta_0\})} + \norm{B|_{\{\zeta = \zeta_0\}}}_{v^{\ell-1}r^{k+j}L^{\infty}(\{\zeta = \zeta_0\})}\\
&+\norm{R_A}_{v^{\ell-1}r^{k+\max(j,1)}\zeta^{k-1}L^{\infty}(\{\sigma \geq \epsilon\})} + \norm{R_B}_{v^{\ell-1}r^{k+j}\zeta^{k-1}L^{\infty}(\{\sigma \geq \epsilon\})}.\end{align*}
Then:
\begin{prop}\label{thm:C5:initiala}Suppose $A,B \in C^{1}(\{\zeta > 0, v > 0\})$ solves \eqref{eq:C5:mea}. Then $A ,B \in v^{\ell-1}r^{k+j}\zeta^{k-1}L^{\infty}(\{\sigma \geq \epsilon\})$ and 

\[\norm{A}_{v^{\ell-1}r^{k+j}\zeta^{k-1}L^{\infty}(\{\sigma \geq \epsilon\})} + \norm{B}_{v^{\ell-1}r^{k+j}\zeta^{k-1}L^{\infty}(\{\sigma \geq \epsilon\})} \leq C\mathcal I,\]
where $C = C(\mathcal E, \zeta_0,\sup_{\sigma \geq \epsilon} r)$ (the third argument is finite) is increasing in all its arguments.\end{prop}
\begin{rk}
As in \cref{thm:C5:initialreglt}
The reason for $\max(1,j)$ appearing is technical, as it was in \cref{thm:C5:initialreglt}.\end{rk}
\begin{proof}
Let $C$ denote a numerical constant depending only on $\mathcal E$, $\zeta_0$, $j$, and $\sup r < \infty$ which may change from line to line.

On a first reading, it is helpful to suppose $k=j = \ell=0$, and $R_A = R_B = 0$.
We recast the system as an integral system.
Define $\phi,\psi$ by requiring them to solve in $U$
\begin{align*}
(\zeta\pa_\zeta - r\pa_r)\phi &= \phi\left(a_{12}+1\right)\\
(\zeta\pa_\zeta-v\pa_v-r\pa_r) \psi &= \psi \left(a_{21}-1\right),\end{align*}
with $\phi(\lambda,1) \equiv 1$, $\psi(0,\tau) \equiv 1$. Then
\begin{align*}
\phi(v,r,\zeta) &= \exp\left(\int_{\zeta_0}^{\zeta} (1+a_{12})(v,r\zeta/s,s)\ ds/s\right)\\
\psi(v,r,\zeta) &= \exp\left(\int_{\zeta_0}^{\zeta} \left(-1+a_{21}\right)(v\zeta/s,r\zeta/s,s)\ ds/s\right).\end{align*}

Thus,
\begin{align*}(\phi/\zeta), (\zeta/\phi) \leq C\\
\zeta\psi, 1/(\zeta \psi) \leq C.
\end{align*}

 Now using the method of integrating factors,
\begin{align}
\label{eq:C5:boostrapa}
\begin{split}
A(v,r,\zeta)\hspace{-2pt} &=\hspace{-2pt} \phi^{-1}\hspace{-2pt}\left(A(v,r\zeta/\zeta_0,\zeta_0) -\hspace{-2pt}\int_{\zeta_0}^{\zeta} \hspace{-2pt}\left(a_{11}B\phi\right)(v,r\zeta/s,s)\ ds/s + \int_{\zeta_0}^{\zeta} (\phi R_A)(v,r\zeta/s,s)\ ds/s\right)\\
B(v,r,\zeta)\hspace{-2pt} &=\hspace{-2pt} \psi^{-1}\hspace{-2pt}\left(B(v\zeta/\zeta_0,r\zeta/\zeta_0,\zeta_0) + \hspace{-2pt}\int_{\zeta_0}^{\zeta} \hspace{-2pt}\left(-(a_{22}+1) A\psi\right)(v\zeta/s,r\zeta/s,s)\ ds/s \right.\\
&\ \ \ \ \ \ \ \ \ \left.+ \int_{\zeta_0}^{\zeta} (\psi R_B)(v\zeta/s,r\zeta/s,s)\ ds/s\right).\end{split}\end{align}

As always, the integrals are taken along the flow of $\zeta\pa_\zeta - r\pa_r$ and $\zeta\pa_\zeta - v\pa_v -r\pa_r$, respectively.

We will use the same strategy as usual, i.e.\ a bootstrap argument, and find numerical constants $K_1,K_2 > 1$, $N> 0$, large, such that for all $p \in \{\sigma \geq \epsilon\}$ with $\zeta(p),v(p),r(p) > 0$ and $\delta > 0$, then inside $B(p)$
\begin{align*}
|A| &\leq K_1(\mathcal I+\delta)v^{\ell-1}\zeta^{k-1} r^{k+j}e^{N r^2}\\
|B| &\leq K_1K_2(\mathcal I + \delta)v^{\ell-1}\zeta^{k-1} r^{k+j}e^{Nr^2}.\end{align*}
Taking $\delta \to 0$ will prove the proposition, since $K_1, K_2, N$ do not depend on $\delta$ or $p$. We need only check that if the bootstrap assumptions are satisfied in a neighbourhood of $B(p) \n \{\zeta = \zeta_0\}$ and for $\zeta(p) \leq \zeta_1 < \zeta_0$, if they are satisfied in $B(p) \n \{\zeta_1 \leq \zeta < \zeta_0\}$, then $A$, $B$ in fact
satisfy an improved estimate on the region.
By continuity, these bounds are verified in a small neighbourhood of $B(p) \n \{\zeta = \zeta_0\}$ if $K_1, K_2 > 1$. Observe that if $q = (v,\zeta) \in B(p) \n \{\zeta_1 \leq \zeta < \zeta_0\}$,then the arguments of the integrands of \eqref{eq:C5:boostrapa} lie in the same region (and in particular in $\{\sigma \geq \epsilon\})$. Therefore we may plug in the bootstrap assumptions into \eqref{eq:C5:boostrapa} (and use that $\zeta/s \leq 1$ if $\zeta \leq s$) to see that if $q \in B(p) \in \{\zeta_1 \leq \zeta < \zeta_0\}$ \begin{align*}
|A|(v,r,\zeta) & \leq \zeta^{-1}v^{\ell-1}(r\zeta)^{k+j}C\mathcal I\\
&+ C\zeta^{-1}K_1(\mathcal I + \delta)K_2\int_{\zeta}^{\zeta_0}(r\zeta/s)^2v^{\ell-1}(r\zeta)^k(r\zeta/s)^j s^{-1}e^{N(r\zeta/s)^2}s\ ds/s\\
&+ C\zeta^{-1}\mathcal I \int_{\zeta}^{\zeta_0} v^{\ell-1}(r\zeta)^{k+\max(j,1)}s^{-1-\max(j,1)}s\ ds/s\\
&\leq v^{\ell-1}r^{k+j}\zeta^{k+j-1}C\mathcal I\\
&+ CK_1(\mathcal I + \delta)K_2r^{k+j}\zeta^{k-1}v^{\ell-1}\int_{\zeta}^{\zeta_0}(r\zeta)^2 s^{-3}e^{N(r\zeta/s)^2}\ ds\\
&+ C\mathcal I r^{k + \max(j,1)}\zeta^{k-1}v^{\ell-1}\int_{\zeta}^{\zeta_0} \zeta^{\max(j,1)}s^{-1-\max(j,1)}\ ds\\
&\leq C(\mathcal I + \delta) r^{k+j}\zeta^{k-1} v^{\ell-1}\left(\zeta^j + 1- \frac{K_1K_2}{2N} + \frac{K_1K_2}{2N}e^{N r^2}\right)\\
&\leq C(\mathcal I + \delta) r^{k+j}\zeta^{k-1} v^{\ell-1}\left(C- \frac{K_1K_2}{2N} + \frac{K_1K_2}{2N}e^{N r^2}\right)\\
\end{align*}
and
\begin{align*}
|B|(v,r,\zeta) &\leq C\mathcal I \zeta(v\zeta)^{\ell-1}(r\zeta)^{k+j}\\
&+C\zeta K_1(\mathcal I + \delta)\int_{\zeta}^{\zeta_0}s^{-1}v^{\ell-1}(\zeta/s)^{\ell}\zeta^{-1} (r\zeta)^k(r\zeta/s)^{j} e^{N(r\zeta/s)^2}\ ds/s\\
&+ C\zeta\mathcal I\int_{\zeta}^{\zeta_0} s^{-1}v^{\ell-1}(\zeta/s)^{\ell}\zeta^{-1} (r\zeta)^k(r\zeta/s)^{j} \ ds/s\\
&\leq C \mathcal I v^{\ell-1}r^{k+j}\zeta^{\ell+k+j}\\
&+ C(\mathcal I+\delta)( K_1e^{Nr^2}+1)v^{\ell-1}r^{k+j}\zeta^{k}\int_{\zeta}^{\zeta_0} s^{-2}\ ds \\
&\leq C(\mathcal I+\delta) v^{\ell-1}r^{k+j}\zeta^{k-1} \left(\zeta^{\ell+j+1} + \frac{\zeta}{2}(\zeta^{-1}-\zeta_0^{-1})(K_1e^{Nr^2}+ 1)\right)\\
&\leq C(\mathcal I+\delta)v^{\ell-1}r^{k+j}\zeta^{k-1} e^{Nr^2}(K_1+C).
\end{align*}
Replacing $C$ with $\max(C,1)$, we will improve over the bootstrap assumption provided
\begin{romanumerate}
\item $K_2 > \max(2,10C^2)$;
\item $N > CK_2$;
\item $K_1 > \max(2NC/K_2,2)$.
\end{romanumerate}
The first condition bounds $C(K_1+C) < K_1K_2$ (since $K_1 > 2$) the second condition bounds $\frac{CK_1K_2}{2N} < K_1$, and the third bounds $C-\frac{K_1K_2}{2N} < 0$.
\end{proof}

With additional regularity assumptions, we have additional regularity of the solution.
\begin{cor}\label{thm:C5:cora}Let $a_{ij}$, $R_A,R_B$ be as in \eqref{eq:C5:mea}. Fix $N \geq 1$ and suppose $a_{ij} \in r^2 W_b^N(\{\sigma \geq \epsilon\})$, $R_A,R_B \in v^{\ell-1}r^{k+\max(j,1)}\zeta^{k-1}W_b^{N}(\{\sigma \geq \epsilon\})$, $R_B \in v^{\ell-1}r^{k+j}\zeta^{k-1} W^N_b(\{\sigma \geq \epsilon\})$, If $A,B \in C^{N+1}(\{\sigma \geq \epsilon, v > 0\})$ solves \eqref{eq:C5:mea} with data $A|_{\{\zeta = \zeta_0\}},B|_{\{\zeta = \zeta_0\}} \in v^{\ell-1}r^{k+j}W_b^N\{\zeta = \zeta_0\})$, then $A,B \in v^{\ell-1}r^{k+j}\zeta^{k-1} W_b^N(\{\sigma \geq \epsilon\})$, with analogous bounds to those in \cref{thm:C5:initiala}. In particular, if the assumptions are satisfied for all $N$, then \[A,B \in v^{\ell-1}r^{k+j}\zeta^{k-1} \mathcal A(\{\sigma \geq \epsilon\}.\]\end{cor}
The proof is the almost the same as that of \cref{thm:C5:corlt}, so is omitted. 

We need one more lemma, the analogue of \cref{thm:C5:finiteseriesI}.

\begin{lem}\label{thm:C5:finiteseriesII}For $k \in \N_0$, $a_k(\zeta\pa_\zeta)F'|_{\{\zeta = \zeta_0\}},\ a_k(\zeta\pa_\zeta)\bar{H}|_{\{\zeta = \zeta_0\}} \in r^{k}\mathcal A_{\mathrm{phg},(\{v = v_0\},\{r = 0\})}^{(-1,0)}(\{\zeta = \zeta_0\})$.
\end{lem}
\begin{proof}
Recall that polyhomogeneity with index set $0$ is the same as smoothness, and recall \eqref{eq:C5:boredothis} implies that $F'$ and $\bar{H}$ are well-defined and polyhomogeneous on \[\{1/\lambda_0 \leq \zeta \leq \zeta_0, \sigma > 0\}.\] It suffices to prove that for $1/\lambda_0 \leq \zeta \leq \zeta_0$, $F^j = \pa_r^j F'|_{\{r=0\}}$ and $\bar{H}^j = \pa_{r}^j \bar{H}|_{\{r= 0\}}$ have finite expansions in integral powers $\zeta^\ell$, \emph{without remainder}
\begin{align*}
F^j &= \sum_{k=-1}^{j-1} F^j_k\zeta^k\\
\bar{H}^j &= \sum_{k=-1}^{j-1} \bar{H}^j_k \zeta^k,\end{align*}
where \[F^j_k,\bar{H}^j_k \in \mathcal A_{\mathrm{phg},(\{v = 0\})}^{(-1)}(\{\zeta = \zeta_0, \ r = 0\}).\] Indeed, this implies that \[a_k(\zeta\pa_\zeta)\pa_{r}^j F',a_k(\zeta\pa_\zeta)\pa_{r}^j \bar{H} = 0, \ \ k \geq j+1\] at $r=0$, and so \[a_k(\zeta\pa_\zeta) F'|_{\{\zeta = \zeta_0\}},a_k(\zeta\pa_\zeta) \bar{H}|_{\{\zeta = \zeta_0\}}\] have their first $k-1$ $\pa_r$-derivatives equal to $0$ at $r=0$. 

Commuting $\pa_r^j$ with \eqref{eq:C5:FHU} yields
\begin{align}
\label{eq:C5:finiteseriesII}
\begin{split}
\zeta\pa_\zeta F^j + (1-j)F^j &= -\sum_{\ell=2}^j {\binom{j}{\ell}}\left[\left(\pa_r^{\ell}\frac{vf}{4}\right)\bar{H}^{j-\ell} + \left(\pa_r^{\ell} \frac{v\zeta \tilde{h}}{4}\right) F^{j-\ell}\right]\\
(\zeta\pa_\zeta-v\pa_v)\bar{H}^j-F^j - (j+1)\bar{H}^j &= \sum_{\ell=2}^j {\binom{j}{\ell}}\left[\left(\pa_r^{\ell}\frac{vf}{4}\right)\bar{H}^{j-\ell} + \left(\pa_r^{\ell} \frac{v\zeta \tilde{h}}{4}\right) F^{j-\ell}\right]
\end{split}\end{align}
(the sum starts at $2$ since by \eqref{eq:C5:fhalpha} $vf$ and $v\zeta\tilde{h} \in r^2C^{\infty}(U)$).
Also by \eqref{eq:C5:fhalpha}, $v\zeta \tilde{h}$ and $vf$ are both $r^2$ times smooth functions of $r$, $\zeta r$, $v$ $u$. Indeed, the only dependence on $\zeta$ comes from the dependence of $K$ on \[\theta = (\theta^1,\theta^2) = (\zeta r, \zeta r u).\] Thus, at $r=0$, $\pa^{\ell}_r (v\zeta \tilde{h})$ and $\pa^{\ell}_r (vf)$ have finite expansions in powers of $\zeta^k$, $0 \leq k \leq \ell-2$.

Let us look at $j=0$. In this case, $\zeta\pa_\zeta F^0 + F^0 = 0$, and so $F^0 = \zeta_0/\zeta F^0|_{\{\zeta = \zeta_0\}}$, which is a series of the desired form. Also, $\bar{H}^0$ satisfies $(\zeta\pa_\zeta - v\pa_v)\bar{H}^0 - F^0 - \bar{H}^0 = 0$. It will be useful, both for the base case and the inductive step, to return to $(\lambda,\tau,\theta)$ coordinates, so we provide some overheard. To be consistent with the blowup of $r=0$, we will also need to blow up $|\theta| = 0$, choosing projective coordinates $s = \theta^1, u = \theta^2/\theta^1 = \alpha^2/\alpha^1$ (remember we are working without loss of generality in a coordinate chart in which $r = \alpha^1, u = \alpha^2/\alpha^1$). Then at $s=0$, $\zeta\pa_\zeta - v\pa_v = -\lambda\pa_\lambda$. Moreover, setting $\tilde{H}^j = \pa_s^j \tilde{H}|_{\{s = 0\}}$, we have that $\bar{H}^j = \lambda^{-j-1}\tilde{H}^j = \zeta^{j+1}\tilde{H}^j$, since $\pa_s = \lambda\pa_r$. In particular, \[(\zeta\pa_\zeta - v\pa_v - (j+1))\bar{H}^j = -\lambda^{-j}\pa_{\lambda}\tilde{H}^j.\]

Since $\tilde{H}|_{\mathbf{lf}} \equiv 0$ and $\pa_s$ is tangential to $\{\lambda = 0\}$, $\tilde{H}^j|_{\{\lambda = 0\}} = 0$, for all $j$.

Thus, the equation for $\bar{H}^0$ becomes
\[\pa_\lambda \tilde{H}^0 = -F^0 = \tau^{-1}A(\lambda\tau,u),\] where $A \in C^{\infty}(\{\lambda \leq \lambda_0, s = 0\})$, since $v\zeta F^0 \in C^{\infty}(\{\lambda \leq \lambda_0\})$. Thus
\begin{align*}
\tilde{H}^0 &= \int_0^{\lambda} \tau^{-1}A(\mu\tau)\ d\mu\\
&= \frac{1}{\tau^2}\int_0^{\lambda\tau}A(t)\ dt \ \ \ (\text{substituting }t = \mu\tau).\end{align*}
Since the integrands are a smooth function of $t$, $\int_0^x A(t)\ dt = xS(x)$, where $S$ is smooth. Thus
\[\tilde{H}^0 = \lambda/\tau S(\lambda\tau,u).\] Hence, $\bar{H}^0 = \lambda^{-1}\tilde{H}^0 = \tau^{-1}S(\lambda\tau,u)$, or in $(v,\zeta)$ coordinates, $\bar{H}^0 = \zeta^{-1}v^{-1}S(v,u)$, where $S$ is smooth, which is a series of the desired form.

For higher $j$, we have by induction the right-hand sides of \eqref{eq:C5:finiteseriesII} are finite series in $\zeta$ with polyhomogeneous coefficients, with powers $\zeta^k$, $0 \leq k \leq j-3$.

The equation for $F^j$ is
\[(\zeta\pa_\zeta + (1-j))F^j = R,\] where $R$ is the right-hand side. Integrating as in the proof of \cref{thm:C5:finiteseriesI}, we know that this means $F^j$ has a series of the desired form. Now for $\bar{H}^j$, we switch to $(v,\lambda,s)$ coordinates, where the equation reads
\[\pa_\lambda \tilde{H}^j = -\lambda^j F^j + \lambda^j R.\] The right-hand side of this equation we know is itself a finite series
\[\lambda^j\sum_{i=-1}^{j-1} \zeta^i v^{-1}A_i(v) = \sum_{i=-1}^{j-1} \lambda^{j-i-1}\tau^{-1}A_i(\lambda\tau),\]where $A_i$ are smooth. Thus
\begin{align*}
\tilde{H}^j &= \sum_{i=-1}^{j-1} \int_0^{\lambda} \mu^{j-i-1}\tau^{-1}A_i(\mu\tau)\ d\mu\\
&=\sum_{i=-1}^{j-1} \tau^{i-1-j}\int_0^{\tau\lambda} t^{j-i-1}A_i(t)\ dt \ \ \ (\text{substituting }t = \mu\tau).\end{align*}
Since the integrand is a non-negative integral power of $t$ times a smooth function,
\[\int_0^x t^{j-i-1}A_i(t) = x^{j-i}S_i(x),\] for some smooth $S_i$. Thus
\[\tilde{H}^j = \sum_{i=-1}^{j-1} \tau^{-1}\lambda^{j-i}S_i(\lambda\tau),\] and so
\[\bar{H}^j = \sum_{i=-1}^{j-1} \tau^{-1}\lambda^{-i-1}S_i(\lambda\tau) = \sum_{i=-1}^{j-1}v^{-1}\zeta^{-i}S_i(v),\] which is a series of the desired form.
\end{proof}

We may now prove \cref{thm:C5:prop4i}. The proof is very similar to that of \cref{thm:C5:prop2}, so we only indicate what needs to be changed.
\begin{proof}[Proof of \cref{thm:C5:prop4i}]
It suffices to show that for all $\ell,i,k$
\begin{align*}
a_\ell(v\pa_v)p_i(r\pa_r)a_k(\zeta\pa_\zeta)F', a_\ell(v\pa_v)p_i(r\pa_r)a_k(\zeta\pa_\zeta) \bar{H}, \in v^{\ell-1}r^{i}\zeta^{k-1}A(\{\sigma \geq \epsilon\}).
\end{align*}
We will in fact prove the stronger claim that for $i \leq k$
\begin{align*}
a_\ell(v\pa_v)p_i(r\pa_r)a_k(\zeta\pa_\zeta)F', a_\ell(v\pa_v)p_i(r\pa_r)a_k(\zeta\pa_\zeta) \bar{H}, ' \in v^{\ell-1}r^{\max(i,k)}\zeta^{k-1}A(\{\sigma \geq \epsilon\}).
\end{align*}
One argues by induction on $i + k + \ell$ using, \cref{thm:C5:cora} and \cref{thm:C5:finiteseriesII}, in almost the same manner as in the proof of \cref{thm:C5:prop2}. We omit the details.
\end{proof}

We can now prove \cref{thm:C5:prop4ii}.
\begin{proof}[Proof of \cref{thm:C5:prop4ii}]
Let $U_4'$ the image of $U$ in $U_4$ under the blowdown. It suffices to show
\begin{align}
\label{eq:C5:phgrel}
\begin{split}
a_j(\zeta\pa_\zeta)a_k(v\pa_v)F' \in \zeta^{j-1}v^{k-1}\mathcal A(U_4'\n\{\sigma \geq \epsilon\})\\
a_j(\zeta\pa_\zeta)a_k(v\pa_v) \bar{H} \in \zeta^{j-1}v^{k-1}\mathcal A(U_4'\n\{\sigma \geq \epsilon\}).
\end{split}
\end{align}
In our projective coordinate patch, $U$, the vector fields $\pa_{\alpha^1}$ and $\pa_{\alpha^2}$ are spanned over $C^{\infty}(U)$ by $\pa_r$ and $r^{-1}\pa_u$. Thus, it suffices to show that for all $j$
\begin{equation}\label{eq:C5:blowdown}r^{-j}\pa_u^j F',r^{-j}\pa_u^j \bar{G}' \in \mathcal A_{\mathrm{phg},(\{v = 0\}, \{r = 0\}, \{\zeta = 0\})}^{(-1,0,-1)}(\{\sigma \geq \epsilon\}\n U).\end{equation}
Indeed, since polyhomogeneity with index set $0$ is the same as smoothness, it follows that\[A_{\mathrm{phg},(\{v = 0\}, \{r = 0\}, \{\zeta = 0\})}^{(-1,0,-1)}(\{\sigma \geq \epsilon\}\n U)\] is stable under differentiation via $\pa_r$, and so \eqref{eq:C5:blowdown} would imply that for all $j,k$,$\ell_1,\ell_2$,
\begin{align*}
\begin{split}
a_j(\zeta\pa_\zeta)a_k(v\pa_v)\pa_{\alpha_1}^{\ell_1}\pa_{\alpha_2}^{\ell_2} F' \in \zeta^{j-1}v^{k-1}\mathcal A(\{\sigma \geq \epsilon\} \n U)\\
a_j(\zeta\pa_\zeta)a_k(v\pa_v)\pa_{\alpha_1}^{\ell_1}\pa_{\alpha_2}^{\ell_2} \bar{H} \in \zeta^{j-1}v^{k-1}\mathcal A(\{\sigma \geq \epsilon\} \n U),
\end{split}
\end{align*}
which is the same thing as \eqref{eq:C5:phgrel}.

To prove \eqref{eq:C5:blowdown}, it suffices for\footnote{Technically, this requires the use of \cref{thm:A4:betterphg}. The quoted identity is sufficient to show that the first few terms in the polyhomogeneous expansion at $\{r = 0\}$ vanish. However, it is not immediately obvious that the first few coefficients in the expansions at $\{r = 0\}$ of the expansions at $\{\zeta = 0\}$ and $\{v = 0\}$ also vanish. Cf. \cref{fn:C4:ifn}.}

\[\pa_u^j \pa_r^i F'|_{\{r = 0\}}, \pa_u^j \pa_r^i \bar{H}|_{\{r = 0\}} \equiv 0, \ \ \ 0 \leq i \leq j-1.\] In the notation of the proof of \cref{thm:C5:finiteseriesII}, since $\pa_u$ is tangent to $\{r = 0\}$, 
\[\pa_u^j \pa_r^i F'|_{\{r = 0\}} = \pa_u^{j}F^i, \ \pa_u^j \pa_r^i \bar{H}|_{\{r = 0\}} = \pa_u^{j} \bar{H}^i.\] Certainly $\pa_u^{j}F^i, \pa_u^{j}\bar{H}^i = 0$ for all $0 \leq i \leq j-1$ if and only if \begin{equation}\label{eq:C5:aaad}\pa_u^{i+1}F^i, \pa_u^{i+1}\bar{G}^i = 0\end{equation} so we have reduced \eqref{eq:C5:phgrel} to \eqref{eq:C5:aaad}.

From \eqref{eq:C5:finiteseriesII} in the proof of \cref{thm:C5:finiteseriesII}, $\pa_u^{i+1}F^i, \pa_u^{i+1}\bar{H}^i$ satisfy equations of the form
\begin{subequations}
\label{eq:C5:FHil}
\begin{gather}
\label{eq:C5:Fil}
\zeta\pa_\zeta \pa_u^{i+1}F^i + (1-i)\pa_u^{i+1}F^i = \pa_u^{i+1}R_{F^i}\\
\label{eq:C5:Hil}
(\zeta\pa_\zeta-v\pa_v)\pa_u^{i+1}\bar{H}^i-\pa_u^{i+1}F^i - (i+1)\pa_u^{i+1}\bar{H}^i = \pa_u^{i+1}R_{\bar{H}^i},
\end{gather}
\end{subequations}
where $R_{F^i}$, $R_{\bar{H}^i}$ are sums of terms of the form
\[\pa_{r}^{i-\ell}(\phi)\Phi^{\ell},\] for $0 \leq \ell < i$, where $\phi$ is either $vf,v\zeta\tilde{h}$, and $\Phi^{\ell}$ is either $F^{\ell}$ or $\bar{H}^{\ell}$. We will prove \eqref{eq:C5:aaad} by induction. For $i=0$, \eqref{eq:C5:Fil} reads
\[\zeta\pa_\zeta \pa_u F^0 + \pa_uF^i = 0,\] and thus \begin{equation}\label{eq:C5:blowdown2}\pa_uF^0 = \frac{\zeta_0}{\zeta}\pa_u F^0|_{\{\zeta = \zeta_0\}}.\end{equation}

We already know that
\[F' \in \phgi{(\{v = 0\})}{(-1)}(U_4 \n \{\zeta = \zeta_0\}),\] and thus since $\pa_u = r\pa_{\alpha_2}$(at least in the projective chart $r = \alpha_1$, $u = \alpha^2/\alpha^1$ which covers $U$),
\[\pa_u F^0|_{\{\zeta = \zeta_0\}} = \pa_u F'|_{\{\zeta = \zeta_0,r = 0\}} = 0.\]
Therefore from \eqref{eq:C5:blowdown2} that $\pa_uF^0 \equiv 0$, too, which is \eqref{eq:C5:aaad} in the case $i = 0$ for $F$.Similarly, for $i=0$, \eqref{eq:C5:Hil} reads
\[(\zeta\pa_\zeta - v\pa_v)\pa_u \bar{H}^0 - \pa_u \bar{H}^0 = 0,\] and so from \eqref{eq:C5:FHil}
\begin{equation}\label{eq:C5:Hbarz}\pa_u\bar{H}^0 = \frac{\zeta}{\zeta_0}\pa_u \bar{H}^0(v\zeta/\zeta_0,\zeta_0).\end{equation} Since
\[\bar{H} \in \phgi{(\{v = 0\})}{(-1)}(U_4 \n \{\zeta = \zeta_0\}),\] too, in particular for all $i$, $\pa_u^{i+1}\bar{H}^i|_{\{\zeta = \zeta_0\}} \equiv 0$, too, and thus by \eqref{eq:C5:Hbarz} $\pa_u \bar{H}^0 \equiv 0$. This completes the proof of \eqref{eq:C5:aaad} in the case $i = 0$.

Now for higher $i$, observe that $\pa_u^{i+1}R_{F^i}$, $\pa_u^{i+1}R_{\bar{H}^i}$ are sums of terms of the form
\begin{equation}\label{eq:C5:ada}[\pa_u^{i+1-k}\pa_r^{i-\ell}(\phi)][\pa_u^{k}\Phi^{\ell}],\end{equation}{eq:C5:ada} for $0 \leq \ell < i, 0 \leq k < i+1$, where $\phi$ and $\Phi^\ell$ are as above. By induction, the second factor is $0$ whenever $k \geq \ell+1$. Since $\phi$ is smooth on the blowdown $U_4$,
\[\pa_u^{i+1-k}\pa_r^{i-\ell}\phi = \pa_r^{i-\ell}\pa_u^{i+1-k}\phi = \pa_r^{i-\ell}(r^{i+1-k}\pa_{\alpha^2}^{i+1-k})\] is $0$ on $\{r = 0\}$ whenever $i+1-k \geq i-\ell+1$. Since one of the conditions $i+1-k \geq i-\ell+1$ and $k \geq \ell + 1$ these two conditions is always met, \eqref{eq:C5:ada} is exactly $0$, and so $\pa_u^{i+1}R_{F^i} = \pa_u^{i+1}R_{\bar{G}^i} = 0$. Thus, like in the base case we may argue
\[\pa_u^{i+1}F^i = \frac{\zeta^{i-1}}{\zeta_0^{i-1}}F^i|_{\{\zeta = \zeta_0\}} = 0,\]
and
\[\pa_u^{i+1}\bar{G}^i = \frac{\zeta^{i+1}}{\zeta_0^{i+1}}\bar{G}^i(v\zeta/\zeta_0,\zeta_0) = 0.\] This completes the induction, and so \eqref{eq:C5:aaad} holds for all $i$.
\end{proof}

We may now prove part (ii) of \cref{thm:C5:prop4i} and \cref{thm:C5:prop4ii}. The proof closely follows that of \cref{thm:C5:prop2}:(ii). Fix $j,k,\ell \in \R$, and consider the equations:
\begin{equation}\label{eq:C5:wvzaii}(\zeta\pa_\zeta-r\pa_r) w = S\end{equation}
and
\begin{equation}\label{eq:C5:wvzai}(\zeta\pa_\zeta -v\pa_v - r\pa_r) w = S,\end{equation}
where in both cases $S \in v^jr^k \zeta^\ell L^{\infty}(\{\sigma \geq \epsilon\})$.
We start with a general proposition.
\begin{prop}\label{thm:C5:omegavzahelper}Suppose $w \in C^1(\{\sigma \geq \epsilon, v > 0\})$ solves \eqref{eq:C5:wvzaii} with data $w|_{\{\zeta = \zeta_0\}} \in v^jr^k L^{\infty}(\{\zeta = \zeta_0\})$.
Assume also that $\ell-k \neq 0$.
Then $w \in v^jr^k\zeta^{\min(\ell,k)} L^{\infty}(\{\sigma \geq \epsilon\})$ with the bound
\[\norm{w}_{v^jr^k\zeta^{\min(k,\ell)} L^{\infty}(\{\sigma \leq \epsilon\})} \leq C(\norm{w|_{\{\zeta = \zeta_0\}}}_{v^jr^kL^{\infty}(\{\zeta=\zeta_0\})} + \norm{S}_{v^jr^k\zeta^\ell L^{\infty}(\{\sigma \geq \epsilon\})}),\]
where $C = C(j,k,\ell,\sup_{\{\sigma \geq \epsilon\}} r)$ (the last argument is finite) is increasing in all its arguments.

Now fix $N \geq 1$. If instead $w \in C^{N+1}(\{\sigma \geq \epsilon, v > 0\})$ and $S \in v^jr^k\zeta^\ell W_b^N(\{\sigma \geq \epsilon\})$ and $w|_{\{\zeta = \zeta_0\}} \in v^jr^k\zeta^\ell W_b^N(\{\zeta = \zeta_0\})$, then $w \in v^jr^k\zeta^{\min(k,\ell)} W_b^N(\{\sigma \geq \epsilon\})$ with an analogous bound. In particular, if the hypotheses hold for all $N$, then $w \in v^jr^k\zeta^{\min(k,\ell)} \mathcal A(\{\sigma \geq \epsilon\})$.

Suppose instead $w \in C^1(\{\sigma \geq \epsilon, v > 0\})$ solves \eqref{eq:C5:wvzai}
with data $w|_{\{\zeta = \zeta_0\}}\hspace{-2pt} \in v^jr^k \hspace{-2pt} L^{\infty}(\{\zeta = \zeta_0\})$.
Assume also that $\ell-j-k \neq 0$.
Then $w \in v^jr^k\zeta^{\min(\ell,k)} L^{\infty}(\{\sigma \geq \epsilon\})$ with the bound
\[\norm{w}_{v^jr^k\zeta^{\min(j+k,\ell)} L^{\infty}(\{\sigma \geq \epsilon\})} \leq C(\norm{w|_{\{\zeta = \zeta_0\}}}_{v^jr^kL^{\infty}(\{\zeta = \zeta_0\})} + \norm{S}_{v^jr^k\zeta^\ell L^{\infty}(\{\sigma \geq \epsilon\})}),\]
where $C = C(j,k,\ell,\sup_{\{\sigma \geq \epsilon\}} r,\zeta_0)$ is increasing in all its arguments.

Now fix $N \geq 1$. If instead $w \in C^{N+1}(\{\sigma \geq \epsilon, v > 0\})$ and $S \in v^jr^k\zeta^\ell W_b^N(\{\sigma \geq \epsilon\})$ and $w|_{\{\zeta = \zeta_0\}} \in v^jr^k\zeta^\ell W_b^N(\{\zeta = \zeta_0\})$, then $w \in v^jr^k\zeta^{\min(j+k,\ell)} W_b^N(\{\sigma \geq \epsilon\})$ with an analogous bound. In particular, if the hypotheses hold for all $N$, then $w \in v^jr^k\zeta^{\min(j+k,\ell)} \mathcal A(\{\sigma \geq \epsilon\})$.
\end{prop}
\begin{proof}
The proof is nearly the same as that of \cref{thm:C5:omegalthelper}, and is omitted.
\end{proof}

We need one more lemma:
\begin{lem}\label{thm:C5:finiteseriesomegaII}For $k \in \N_0$, 
\[p_k(\zeta\pa_\zeta)\omega|_{\{\zeta=\zeta_0\}}, \ p_k(\zeta\pa_\zeta)(\zeta\pa_\zeta)\omega|_{\{\zeta = \zeta_0\}} \in r^{k+2}C^{\infty}(\{\zeta=\zeta_0\}).\]
\end{lem}
\begin{proof}
The proof is a hybrid of the proof of \cref{thm:C5:finiteseriesomegaI} and \cref{thm:C5:finiteseriesII}. It suffices to show that for $j \geq 0$, in a neighbourhood of $\{\zeta=\zeta_0\}$, $\omega^j = \pa_r^j \omega |_{\{r= 0\}}$ has a finite expansion in integral powers $\zeta^\ell$ without remainder
\[\omega^j = \sum_{\ell=0}^{j-2} \omega_\ell^j \zeta^j,\]
where $\omega_\ell^j \in C^{\infty}(\{\zeta=\zeta_0\})$.
We know $\omega$ is smooth in a neighbourhood of $\{\zeta = \zeta_0\}$, even to the boundaries, so commuting $\pa_r$ with \eqref{eq:C5:omegavza}, we obtain the equation
\begin{equation}\label{eq:C5:seriesomegaai}(\zeta\pa_\zeta-j)(\zeta\pa_\zeta - v\pa_v-j)\omega^j = \pa_r^{j}\left.\left(\frac{v^2\zeta}{16}f\tilde{h} - \frac{v^2\zeta^2r^2\langle u \rangle^2}{4}F'\bar{H} - \frac{v\zeta}{4}f\right)\right|_{\{r=0\}}.\end{equation}

Using \eqref{eq:C5:fhalpha}, one sees that both $\pa_r^{j}( v^2\zeta f\tilde{h})|_{\{r = 0\}}$ and $\pa_r^{j}(v\zeta f)|_{\{r = 0\}}$ have expansions in powers of $\zeta^\ell$ for $0 \leq \ell \leq j-2$.

From the proof of \cref{thm:C5:finiteseriesII}, we also have that for $i \geq 0$, $\pa_r^{i}F'|_{\{r = 0\}}$ and $\pa_r^i \bar{H}|_{\{r = 0\}}$ have expansions in $\zeta^\ell$ for $-1 \leq \ell \leq i-1$, with coefficients in 
\[\mathcal A_{\mathrm{phg},(\{v = 0\})}^{(-1)}(\{\zeta = \zeta_0, \ r = 0\}).\] Thus $\pa_r^j(\zeta^2 v^2r^2\langle u \rangle^2 F'\bar{H})|_{\{r = 0\}}$ has a finite expansion in powers of $\zeta^{\ell}$ for $0 \leq \ell \leq j-2$, with smooth coefficients (the $v^2$ factor means the index set of the polyhomogeneity at $\{v=0\}$ is $0$, and hence the coefficients are actually smooth).

Putting this all together, it follows from \eqref{eq:C5:seriesomegai} that
\begin{equation}(\zeta\pa_\zeta-j)(\zeta\pa_\zeta - v\pa_v-j)\omega^j = \sum_{\ell = 0}^{j-2} S_\ell \zeta^\ell,\end{equation} for some coefficients $S_{\ell} \in C^{\infty}(\{\zeta=\zeta_0, r= 0\})$. Using the same argument as in the proof of \cref{thm:C5:finiteseriesII}, one deduces that
\[(\zeta\pa_\zeta - v\pa_v-j)\omega^j = \sum_{\ell = 0}^{j-2} S'_\ell \zeta^\ell\]
for some other coefficients $S'_{\ell} \in C^{\infty}(\{\zeta=\zeta_0, r= 0\})$. We may use coordinates $(\lambda,\tau,s,u)$ as in the proof of \cref{thm:C5:finiteseriesII} to rewrite this equation as
\[-\pa_\lambda\tilde{\omega}^j = \sum_{\ell = 0}^{j-2} S'_\ell(\lambda\tau)\lambda^{j-\ell},\]
where $\tilde{\omega}^j = \pa_\lambda^j \omega^j|_{\{s = 0\}}$. Since $\omega|_{\mathbf{lf}} \equiv 0$, the same argument as in the proof of \cref{thm:C5:finiteseriesII} shows that that $\omega^j$ has the desired series expansion.
\end{proof}

We may now prove part (ii) of \cref{thm:C5:prop4i} and \cref{thm:C5:prop4ii}.
\begin{proof}[Proof of \cref{thm:C5:prop2}:(ii)]
It suffices to show that for all $j,k \in \N_0$
\begin{equation}\label{eq:C5:omegavzasf}p_j(v\pa_v)p_k(r\pa_r) p_\ell(\zeta\pa_\zeta)\omega \in v^jr^{\max(k,2)}\zeta^{\ell} \mathcal A(U),\end{equation} the $\max(k,2)$ showing that there are no terms $r^0$ or $r^1$ in the polyhomogeneous expansion of $\omega$ at $\{\lambda = 0\}$. The proof is nearly the same as that of \cref{thm:C5:prop1}:(ii), so we only provide a sketch.

We need to examine the right-hand side of \eqref{eq:C5:omegaU}. Using \eqref{eq:C5:fhalpha}, one sees that
\begin{align}\label{eq:C5:omegavzai}p_j(v\pa_v)p_k(r\pa_r)p_{\ell}(\zeta\pa_\zeta)(v^2\zeta f \tilde{h}) \in v^jr^{\max(k,\ell+4)}\zeta^{\ell}A(U)\\
\label{eq:C5:omegavzaii}p_j(v\pa_v)p_k(r\pa_r)p_{\ell}(\zeta\pa_\zeta)(v \zeta f) \in v^jr^{\max(k,\ell+2)}\zeta^{\max(1,\ell)} \mathcal A(U).\end{align}
The proof of \cref{thm:C5:prop4i}:(i) in fact shows that
\[a_j(v\pa_v)p_k(r\pa_r)a_\ell(\zeta\pa_\zeta) F', a_j(v\pa_v)p_k(r\pa_r)a_\ell(\zeta\pa_\zeta)\bar{H} \in v^{j-1}r^{\max(k,\ell)}\zeta^{\ell-1}\mathcal A(\{\sigma \geq \epsilon\}).\]
Therefore
\begin{equation}\label{eq:C5:omegavzaiv}p_j(v\pa_v)p_k(r\pa_r)p_\ell(\zeta\pa_\zeta)(v^2\zeta^2 r^2 F'\bar{H}) \in v^jr^{\max(k,\ell+2)}\zeta^{\ell}\mathcal A(\{\sigma \geq \epsilon\}).\end{equation}

Commuting $p_j(v\pa_v)p_k(r\pa_r)p_\ell(\zeta\pa_\zeta)$ with \eqref{eq:C5:omegaU} and using \eqref{eq:C5:omegavzai}--\eqref{eq:C5:omegavzaiv} shows that
\[(\zeta\pa_\zeta-r\pa_r)(\zeta\pa_\zeta-v\pa_v-r\pa_r) p_j(v\pa_v)p_k(r\pa_r)p_\ell(\zeta\pa_\zeta)\omega \in v^jr^{\max(k,\ell+2)}\zeta^{\ell}\mathcal A(U).\] Notice that regardless of $j,k,\ell \geq 0$ \[\ell-\max(k,\ell+2), \ell-\max(k,\ell+2)-j < 0.\]
Thus we may conclude by \cref{thm:C5:finiteseriesomegaII} and \cref{thm:C5:omegavzahelper}.
\end{proof}

\begin{proof}[Proof of \cref{thm:C5:prop4ii}:(ii)]
As in the proof of \cref{thm:C5:prop4ii}:(i), it suffices to show that for all $j$
\[r^{-j}\pa_u^j \omega \in \phgi{(\{v = 0\},\{r = 0\},\{\zeta = 0\})}{(0,0,0)}(\{\sigma \geq \epsilon\}),\]
which will be implied provided for $i \geq 0$ $\pa_u^{i+1}\omega^i = 0$, where $\omega^i = \pa_r^i \omega|_{\{r = 0\}}$ is as in the proof of \cref{thm:C5:finiteseriesomegaII}. From \eqref{eq:C5:seriesomegaai}, and \cref{thm:C5:prop4ii}:(i), $\pa_u^{i+1}\omega^i$ satisfies
\[(\zeta\pa_\zeta-i)(\zeta\pa_\zeta - v\pa_v-i)\pa_u^{i+1}\omega^i = \pa_u^{i+1}\pa_r^{i}S|_{\{r = 0\}},\]
where $S$ is smooth on the blowdown $U_4 \n \{\sigma \geq \epsilon\}$. Therefore, the right-hand side is $0$ and so in fact
\[(\zeta\pa_\zeta-i)(\zeta\pa_\zeta - v\pa_v-i)\pa_u^{i+1}\omega^i \equiv 0.\]
Integrating, it follows that
\[(\zeta\pa_\zeta - v\pa_v-i)\pa_u^{i+1}\omega^i = \frac{\zeta^i}{\zeta_0^i}((\zeta\pa_\zeta - v\pa_v-i)\pa_u^{i+1}\omega^i)|_{\{\zeta = \zeta_0\}}.\]

However, $\omega^i$ is smooth on the blowdown in a neighbourhood of $\{\zeta = \zeta_0\}$, since $\omega \in \phgi{(\{v = 0\})}{(0)}(U_4 \n \{1/\lambda_0 \leq \zeta \leq \zeta_0\})$. In particular, the right-hand side is just $0$. Therefore
\[(\zeta\pa_\zeta - v\pa_v-i)\pa_u^{i+1}\omega^i = 0,\]
and so integrating
\[\pa_u^{i+1}\omega^i = \frac{\zeta^i}{\zeta_0^i}\pa_u^{i+1}\omega^i(v\zeta/\zeta_0,\zeta_0),\] which is again $0$ since $\omega^i$ is smooth on the blowdown in a neighbourhood of $\{\zeta = \zeta_0\}$. This completes the proof.
\end{proof}

\subsection{Behaviour near \texorpdfstring{$\{\sigma = 0\}$}{\{sigma = 0\}}.}

For this subsection, we work only in $U$. The proof of \cref{thm:C5:prop4iii} mirrors very closely that of \cref{thm:C5:prop3ii}, so we only give the indications of what needs to be changed. The main difference is that in the current setting the vector field $-r\pa_r$ plays the role of $\pa_\zeta$. The intuition behind this is that the forwards-directed integral curves of $r\pa_r$ are transverse to $\{\sigma = 0\}$ and cross from $\{\sigma > 0\}$ to $\{\sigma < 0\}$ in the current setting, while in the proof of \cref{thm:C5:prop3ii}, it is the backwards-directed integral curves of $\pa_\zeta$ that have the same property.

Define
\begin{gather*}
\beta = (\zeta\pa_\zeta-v\pa_v-r\pa_r)\sigma\\
\gamma = (\zeta\pa_\zeta-r\pa_r)\sigma,
\end{gather*}
which are both positive for $\{\sigma \leq \epsilon\}$ by \eqref{eq:C5:sigders2}. Since $U \n \{\sigma \leq \epsilon\}$ is compact, $\beta$, $\gamma$ are bounded uniformly away from $0$.  

 From \eqref{eq:C5:sigders2},
\[vf = \frac{-2\beta}{\sigma}(1+\sigma C^{\infty}(U)), \ v\zeta \tilde{h} = \frac{-2\gamma}{\sigma}(1+\sigma C^{\infty})(U).\]

We first focus on proving part (i) of \cref{thm:C5:prop4iii}. Let us introduce $\Phi = -\sigma/\beta \bar{H}$ and $\Psi = \sigma/\gamma F'$. 

We rewrite \eqref{eq:C5:FHU}
\begin{align}
\label{eq:C5:PPU}
\begin{split}
\sigma(\zeta\pa_\zeta - v\pa_v-r\pa_r)\Phi &= \frac{\beta}{2}(1+\sigma C^{\infty})\Phi - \frac{\beta}{2}(1+\sigma C^{\infty})\Psi\\
\sigma(\zeta\pa_\zeta -r\pa_r)\Psi &= -\frac{\gamma}{2}(1+\sigma C^{\infty})\Phi + \frac{\gamma}{2}(1+\sigma C^{\infty})\Psi.
\end{split}
\end{align}

We will need to to add an inhomogeneous right-hand side to \eqref{eq:C5:PPvz}, and also extend from a pair of scalar functions to a pair of vector-valued function on $\R^d$ ($d \geq 1$). Fix $0 < \epsilon < 1$ small enough so that $\{\sigma \leq \epsilon\} \subseteq U$, and consider the following $2d\times 2d$ system, for $k,\ell \in \N_0$
\begin{align}
\label{eq:C5:mevu}
\begin{split}
\sigma(\zeta\pa_\zeta-v\pa_v-r\pa_r)A &= \frac{\beta}{2}(1+\sigma E_{11})A -\frac{\beta}{2}(1+\sigma E_{12})B + \sigma R_A\\
\sigma(\zeta\pa_\zeta-r\pa_r) B &= -\frac{\gamma}{2}(1+\sigma E_{21})A + \frac{\gamma}{2}(1+\sigma E_{22})B + \sigma R_B,\end{split}\end{align}
Where $E_{ij} \in L^{\infty}(\{0 \leq \sigma \leq \epsilon\};M_d(\R))$, and for some $R_A,R_B \in v^{k-1}\zeta^{\ell-1}L^{\infty}(\{0 \leq \sigma \leq \epsilon\};\R^d)$ (here $M_d(\R)$ denotes the space of $d\times d$ matrices). Suppose further $A|_{\{\sigma = \epsilon\}}, B|_{\{\sigma = \epsilon\}} \in v^{k-1}\zeta^{\ell-1}L^{\infty}(\{\sigma = \epsilon\};\R^d)$.

To define the space $L^{\infty}(\{0 \leq\sigma \leq \epsilon\};M_d(\R))$ and $L^{\infty}(\{0 \leq\sigma \leq \epsilon\};\R^d)$, we use any norm on $\R^d$, and the associated operator norm on $M_d(\R)$ to simplify our computations. Henceforth, we will omit $M_d(\R)$ and $\R^d$ from out notation for clarity. 

Let
\begin{align*}
\mathcal I &= \max(\norm{A|_{\{\sigma = \epsilon\}}}_{v^{k-1}\zeta^{\ell-1}L^{\infty}(\{\sigma = \epsilon\})}, \norm{B|_{\{\sigma = \epsilon\}}}_{v^{k-1}\zeta^{\ell-1}L^{\infty}(\{\sigma = \epsilon\})})\\
\mathcal R &= \norm{R_A}_{v^{k-1}\zeta^{\ell-1}L^{\infty}(\{0 \leq \sigma \leq \epsilon\})} + \norm{R_B}_{v^{k-1}\zeta^{\ell-1}L^{\infty}(\{0 \leq \sigma \leq \epsilon\})}\\
\mathcal E &= \sum_{ij} \norm{E_{ij}}_{L^{\infty}(\{0 \leq \sigma \leq \epsilon\})}.\end{align*}

\begin{prop}\label{thm:C5:initialvu}Suppose $A,B \in C^1(\{0 < \sigma \leq \epsilon, v > 0\})$ solves \eqref{eq:C5:mevu}. Then $A,B \in v^{k-1}\zeta^{\ell-1}L^{\infty}(\{0 \leq \sigma \leq \epsilon\})$ and
\[\norm{A}_{v^{k-1}\zeta^{\ell-1}L^{\infty}(\{0 \leq \sigma \leq \epsilon\})} + \norm{B}_{v^{k-1}\zeta^{\ell-1}L^{\infty}(\{0 \leq \sigma \leq \epsilon\})} \leq C(\mathcal I + \mathcal R),\] where $C = C(k,\ell,\mathcal E, \sup \beta^{-1},\sup \gamma^{-1})$ is increasing in all its arguments.
\end{prop}

Before proving \cref{thm:C5:initialvu}, we will need some overhead. Fix $0 < \vepsilon \leq \epsilon$, and $p = (v,\zeta,r,u)$ with $\sigma(p) \leq \epsilon$. Let $q_1 = (v_1,\zeta_1,r_1,u_1)$ be the unique point for which the integral curve of $\zeta\pa_\zeta-r\pa_r$ starting from $p$ intersects $\{\sigma = \vepsilon\}$.\footnote{The integral curve intersects $\{\sigma = \vepsilon\}$ precisely once because $\sigma$ is strictly increasing along the flow for $\{\sigma \leq \epsilon\}$ since $\epsilon < 1$.} Up to scaling the parameter, an integral curve starting at $p$ is the curve
\[t \mapsto (v,\zeta r/t,t,u).\] Thus \[q_1 = (v,\zeta/\mathbf{i}(v,\zeta r,u),\mathbf{i}(v,\zeta r,u),u)\] for some $\mathbf{i} = \mathbf{i}_\vepsilon$.

Now, let $q_2 = (v_2,\zeta_2,\theta_2)$ be the unique point for which the integral curve of $\zeta\pa_\zeta-v\pa_v-r\pa_r$ starting from $p$ intersects $\{\sigma = \epsilon\}$. Up to scaling the parameter, an integral curve starting at $p$ is the curve
\[t \mapsto (tv/r,\zeta r/t,t,u)\] Thus 
\[q_2 = (\mathbf{j}(v/r,\zeta r,u)v/r,\zeta r/\mathbf{j}(v/r,\zeta r,u),\mathbf{j}(v/r,\zeta r,u),y)\] for some $\mathbf{j} = \mathbf{j}_\vepsilon$.

Since \[\sigma(v,r\zeta/\mathbf{i},\mathbf{i},u) = \sigma(v/r\mathbf{j},r\zeta/\mathbf{j},\mathbf{j},u) = \vepsilon,\]
the implicit function theorem shows that $\mathbf{i}, \mathbf{j}$ are smooth and satisfy $(\zeta\pa_\zeta-r\pa_r) \mathbf{i} = 0$, $(\zeta\pa_\zeta-v\pa_v-r\pa_r)\mathbf{j} = 0$.

\Cref{thm:C5:integration1} has an analogue in the current setting, whose proofs is the same.
\begin{lem}Let $f:(0,\infty) \to \R$ be a smooth function, and let $f'$ denotes its derivative. Then for $0 < \vepsilon \leq \epsilon$
\begin{align*}
\int_{\mathbf{i}_{\vepsilon}(r\zeta)}^{r} \gamma f'(\sigma(v,r\zeta/s,s))\ ds/s &= f(\sigma(v,r,\zeta))-f(\vepsilon)\\
\int_{\mathbf{j}_{\vepsilon}(r\zeta,v/r)}^{r} \beta f'(\sigma(sv/r,r\zeta/s),s)\ ds/s &= f(\sigma(v,r,\zeta))-f(\vepsilon).\end{align*}\end{lem}

We can now sketch the proof of \cref{thm:C5:initialvu}, which follows the same steps as the proof of \cref{thm:C5:initialvzs}.
\begin{proof}[Proof sketch of \cref{thm:C5:initialvu}]
Fix $0 < \vepsilon \leq \epsilon$, and set
Set $\mathbf{i} = \mathbf{i}_{\vepsilon}$ and $\mathbf{j} = \mathbf{j}_{\vepsilon}$.
Conjugating, we may suppose $k = \ell = 1$. Using the method of integrating factors,
\begin{align*}
A(v,r,\zeta) &= \frac{\sigma^{1/2}}{\vepsilon^{1/2}}\left(A(\mathbf{j}v/r,r\zeta/\mathbf{j},\mathbf{j})+ \int_{\mathbf{j}}^{r} \frac{\vepsilon^{1/2}\beta}{2\sigma^{3/2}}(-B+\sigma E_{11}A - \sigma E_{12}B)(sv/r,s,r\zeta/s)\ ds/s\right.\\
&\left. + \int_{\mathbf{j}}^{r} (\vepsilon^{1/2} \sigma^{-1/2} R_A)(sv/r,r\zeta/s,s)\ ds/s\right)\\
B(v,r,\zeta) &= \frac{\sigma^{1/2}}{\vepsilon^{1/2}}\left(B(v,r\zeta/\mathbf{i},\mathbf{i}) + \int_{\mathbf{i}}^{r} \frac{\vepsilon^{1/2}\gamma}{2\sigma^{3/2}}(-A = \sigma E_{21}A + \sigma E_{22}B)(v,r\zeta/s,s)\ ds/s\right.\\
&\left. + \int_{\mathbf{i}}^{r} (\vepsilon^{1/2} \sigma^{-1/2} R_B)(v,s,r\zeta/s)\ ds/s\right).
\end{align*}
The rest of the proof proceeds in an identical fashion to that of \cref{thm:C5:initialvzs}.
\end{proof}

With additional regularity assumptions, we have more regularity of the solution, although we may need to take $\epsilon$ slightly smaller. The proof is almost the same as that of \cref{thm:C5:corvzs}
\begin{cor}\label{thm:C5:corvU}Let $E_{ij}$ ($i,j, = 1,2$), $R_A,R_B$ be as in \eqref{eq:C5:mevzs}. Fix $N \geq 1$ and suppose $E_{ij} \in W_b^N(\{\sigma \leq \epsilon\})$, $R_A,R_B \in v^{k-1}\zeta^{\ell-1}W_b^N(\{\sigma \leq \epsilon\})$. If $A,B \in C^{N+1}(\{0 < \sigma \leq \epsilon, v > 0\})$ solve \eqref{eq:C5:mevz} with data \[A|_{\{\sigma = \epsilon\}},\ B|_{\{\sigma = \epsilon\}} \in v^{k-1}\zeta^{\ell-1}W_b^N(\{\sigma = \epsilon\}),\]
then $A,B \in v^{k-1}\zeta^{\ell-1}W_b^N (\{\sigma \leq \epsilon\})$, with analogous bounds to those in \cref{thm:C5:initialvu} (except now the constant is allowed to depend on $\beta,\beta^{-1},\gamma,\gamma^{-1}$ and its first $N$ derivatives).

If the assumptions are true for all $N$, then in particular $A,B \in v^{k-1}\zeta^{\ell-1}\mathcal A(\{0 \leq \sigma \leq \epsilon\})$.
\end{cor}
\begin{proof}
Setting $\sigma = \sigma(v,\zeta,r,u)$, $(\sigma,v,\zeta,u)$ becomes new coordinates on $\{\sigma \leq \epsilon\}$. Consider the vector fields
\[\nu_0 = \sigma\pa_\sigma, \ \nu_1 = v\pa_v,\ \nu_2 = \zeta\pa_\zeta, \ \nu_3 = \pa_u.\]
Observe that in these coordinates $\zeta\pa_\zeta-r\pa_r = \gamma\pa_\sigma + \zeta\pa_\zeta$ and $\zeta\pa_\zeta - v\pa_v - r\pa_r = \beta\pa_\sigma + \zeta\pa_\zeta - v\pa_v$.

We record the commutator formulae
\begin{align*}
[\nu_i,\sigma(\zeta\pa_\zeta-r\pa_r)] &= (\nu_i \log \gamma)\sigma(\zeta \pa_\zeta - r\pa_r) - \sigma( \nu_i\log \gamma - \delta_{i0})\nu_2\\
[\nu_i,\sigma(\zeta\pa_\zeta-v\pa_v-r\pa_r)] &= (\nu_i \log \beta)\sigma(\zeta \pa_\zeta -v\pa_v- r\pa_r) +\sigma (\nu_i\log \beta -\delta_{i0})(\nu_1 - \nu_2).\end{align*}

Using these identities, one may proceed using nearly verbatim as in the proof of \cref{thm:C5:corvzs}, using \cref{thm:C5:initialvu} instead of \cref{thm:C5:initialvzs}.\end{proof}

To complete the proof of \cref{thm:C5:prop4iii}, one first establishes partial polyhomogeneity i.e.:
\begin{prop}It holds that \label{thm:C5:prop4iiihelper}
\[F',\bar{H} \in \sigma^{-1}\phgi{\{v = 0\}, \{\zeta = 0\}}{(-1,-1)}(\{\sigma \leq \epsilon\})\]
\end{prop}
\begin{proof}
In light of \cref{thm:C5:prop4i}, we already know that
\[F',\bar{H} \in \sigma^{-1}\phgi{\{v = 0\}, \{\zeta = 0\}}{(-1,-1)}(\{\sigma \geq \epsilon/2\}.\]

We use coordinates $(\sigma,v,\zeta,u)$. Write $\nu_1 = v\pa_v$ and $\nu_2 = \zeta\pa_\zeta$ in these coordinates. We may rewrite \eqref{eq:C5:PPU} as
\begin{align}
\begin{split}
\sigma(\pa_\sigma - \beta^{-1}(\nu_1-\nu_2))\Phi &= \frac{1}{2}(1+\sigma C^{\infty}(U))\Phi - \frac{1}{2}(1+\sigma C^{\infty}(U))\Psi\\
\sigma(\pa_\sigma +\gamma^{-1}\nu_2)\Psi &=-\frac{1}{2}(1+\sigma C^{\infty}(U))\Phi + \frac{1}{2}(1+\sigma C^{\infty}(U))\Psi.
\end{split}
\end{align}
It suffices to show that for all $k,\ell \in \N_0$, \[a_k(\nu_1)a_\ell(\nu_2)\Phi, \ a_k(\nu_1)a_\ell(\nu_2)\Psi \in v^{k-1}\zeta^{\ell-1}\mathcal A(\sigma \leq \epsilon).\] One proceeds by induction on $k + \ell$ using \cref{thm:C5:cora} in the exact same way as the proof of \cref{thm:C5:prop3iihelper}. We omit the details.\end{proof}

\Cref{thm:C5:prop4iii}:(i) now follows from \cref{thm:C5:prop4iiihelper} in the same way that \cref{thm:C5:prop3ii}:(i) follows from \cref{thm:C5:prop3iihelper}. We omit the details.

We now turn our attention to proving part (ii) of \cref{thm:C5:prop4iii}. Fix $j,k,\ell \in \R$, and $\ell \neq x \in \R$. Consider the equations for a vector-valued function $w$:
\begin{equation}\label{eq:C5:wai} (\sigma(\zeta\pa_\zeta -r\pa_r) -\gamma x) w - \sigma E w = S\end{equation}
and
\begin{equation}\label{eq:C5:waii} (\sigma(\zeta\pa_\zeta -v\pa_v-r\pa_r) -\beta x)w -\sigma E = S,\end{equation}
where in both cases $S \in v^j\zeta^k\sigma^{\ell}L^{\infty}(\{\sigma \leq \epsilon\})$ and $E \in L^{\infty}(\{\sigma \leq \epsilon\})$. We start with a general proposition.
\begin{prop}\label{thm:C5:omegaashelper}Suppose $w \in C^1(\{\sigma > 0, \ v > 0\})$ solves \eqref{eq:C5:wai} or \eqref{eq:C5:waii} with data $w|_{\{\sigma = \epsilon\}} \in v^j\zeta^k L^{\infty}(\{\sigma = \epsilon\})$. Then $w \in v^j\zeta^k\sigma^{\min(x,\ell)} L^{\infty}(\{\sigma \leq \epsilon\})$ with the bound
\[\norm{w}_{v^j\zeta^k \sigma^{\min(x,\ell)}L^{\infty}(\{\sigma \leq \epsilon\})} \leq C(\norm{w|_{\{\sigma = \epsilon\}}}_{v^j\zeta^kL^{\infty}(\{\sigma = \epsilon\})} + \norm{S}_{v^j\zeta^k \sigma^{\ell}L^{\infty}(\{\sigma \leq \epsilon\})}),\]
where \[C = C(\norm{E}_{L^{\infty}(\{\sigma \leq \epsilon\})},\sup \beta^{-1},\sup \gamma^{-1},j,k,\ell).\]

Now fix $N \geq 1$. If instead $w \in C^{N+1}(\{\sigma > 0, \ v > 0\})$, $S \in v^j\zeta^k\sigma^{\ell} W_b^N(\{\sigma \leq \epsilon\})$, $E \in W_b^N(\{\sigma \leq \epsilon\})$, and $w|_{\{\sigma = \epsilon\}} \in v^j\zeta^kW_b^N(\{\sigma = \epsilon\})$, then $w \in v^j\zeta^k \sigma^{\min(x,\ell)}W_b^N(\{\sigma \leq \epsilon\})$ with an analogous bound. In particular, if the hypotheses hold for all $N$, then $w \in v^{j}\zeta^k\sigma^{\min(x,\ell)} \mathcal A(\{\sigma \leq \epsilon\})$.\end{prop}
\begin{proof}
Conjugating by $v^j\zeta^k\sigma^{\ell}$, we may assume $j = k = \ell = 0$. Using the commutator identities established in the proof of \cref{thm:C5:cora}, the proof is now nearly identical to that of \cref{thm:C5:omegavzshelper}, and is omitted.
\end{proof}

To complete the proof, we first propagate polyhomogeneity at $\{v = 0\}$ and $\{\zeta = 0\}$ to $\{\sigma = 0\}$, and then develop full polyhomogeneity.

We may rewrite \eqref{eq:C5:omegaU} as
\begin{equation}\label{eq:C5:omegaas}
\sigma(\zeta\pa_\zeta - r\pa_r-\gamma)\sigma(\zeta\pa_\zeta - v\pa_v-r\pa_r)\omega = \frac{v^2\sigma^2\zeta}{16}f\tilde{h} - \frac{v^2r^2\langle u \rangle^2\sigma^2}{4}F'\bar{H} - \frac{v\zeta\sigma^2}{4}f =: S.\end{equation}

Observe that by \cref{thm:C5:prop4}:(i) and \eqref{eq:C5:fhalpha}, $S \in \phgi{(\{v = 0\},\{\zeta = 0\},\{\sigma = 0\})}{(0,0,0'')}(U)$.

\begin{prop}It holds that
\[\omega \in \phgi{(\{v = 0\},\{\zeta = 0\})}{(0,0)}(\{\sigma \leq \epsilon\}).\]
\end{prop}
\begin{proof}
The proof is a hybrid between the proof of \cref{thm:C5:omegaphgI} and \cref{thm:C5:prop4iiihelper}. In light of \cref{thm:C5:prop4i}:(ii), we already know that
\[\omega \in \phgi{(\{v = 0\},\{\zeta = 0\})}{(0,0)}(\{\sigma \geq \epsilon/2\}).\] Let us return to $(\sigma,v,r,u)$ coordinates and the vector fields $\nu_1 = v\pa_v$ and $\nu_2 = \zeta\pa_\zeta$ in those coordinates. Recall that $\zeta\pa_\zeta - r\pa_r = \gamma\pa_\sigma +\nu_2$ and $\zeta\pa_\zeta - v\pa_v -r\pa_r = \beta\pa_\sigma -\nu_1 + \nu_2$.

It suffices to show that for all $j,k \in \N_0$. \[p_j(\nu_1)p_k(\nu_2)\omega \in \sigma^{0^-}\mathcal A(\{\sigma \leq \epsilon\}).\] For $j = k= 0$, this follows by applying \cref{thm:C5:omegaashelper} twice to \eqref{eq:C5:omegaas}. For $j,k > 0$, one now uses induction on $j+k$ and commutes $p_j(\nu_1)p_k(\nu_2)$ through \eqref{eq:C5:omegaas} to see that $p_j(\nu_1)p_k(\nu_2)\omega$ satisfies an equation of the form
\[\sigma(\zeta\pa_\zeta - r\pa_r-\gamma)\sigma(\zeta\pa_\zeta - v\pa_v-r\pa_r)p_j(\nu_1)p_k(\nu_2)\omega = p_j(\nu_1)p_k(\nu_2)S + R,\]
where $R$ is an error term consisting of sums of terms of the form
\[v^{j-j'}\zeta^{k-k'}L_{j,k,j',k'}p_{j'}(\nu_1)p_{k'}(\nu_2)\omega,\]
where $0 \leq j' +k'< j+k$, and $L_{j,k,j',k'} \in \Diff^2_b(\{\sigma \leq \epsilon\})$ (the coefficients of $L_{j,k,j',k'}$ depend on the derivatives of $\beta$ and $\gamma$). By induction this is in $v^j\zeta^k\sigma^{0^-} \mathcal A(\{\sigma \leq \epsilon\})$, so we may again apply \cref{thm:C5:omegaashelper} twice and conclude.
\end{proof}

\begin{proof}[Proof of \cref{thm:C5:prop4iii}:(ii)]
We can rewrite \eqref{eq:C5:omegaas} schematically as
\[(\sigma\pa_\sigma - 1)\sigma\pa_\sigma \omega = \sigma L \omega + S/(\beta\gamma)\] where $L \in \Diff^2_b(\{\sigma \leq \epsilon\})$. The proof of \cref{thm:C5:prop3ii}:(ii) now applies nearly verbatim. We omit the details.
\end{proof}
\chapter{Curvature blowup for commutative data.}
\label{C:C6}
\section{Curvature blowup}
\label{C:C6:prelims}
We are now in a position to state the precise version of \cref{thm:C1:blowupvague} which we will prove. 

First some preliminaries. Fix a short-pulse tensor $\mathbf{T}$ and associated $\digamma$ as defined in \cref{C:C4:existence}. Recall from \cref{C:C4:existence} the manifold $\Mp = \overline{\mathcal M} \n \{\digamma > 0\}$, which we will think of as
\[([0,1]_\xi \times [0,1]_v \times [0,\infty)_\eta \times S^2_\theta)\n \{\digamma > 0\},\]
as well as the short-pulse bundle $^\sp T\Mp$ from \cref{C:C4:existence} and \cref{C:C4:DNG} and the notion of being in a double-null gauge from \cref{C:C4:DNG}. Let $g$ be a formal solution to $\Ric(g) = 0$ on $\Mp$ in a double-null gauge with some short-pulse data given by $\mathbf{T}$, which exists by \cref{thm:C4:continuation}. Denote by
\[\mathcal K = \mathcal K(g) = |\Riem(g)|_{g}^2= R_{abcd}R^{abcd}\] its Kretschmann scalar.
Then by \cref{thm:C4:ricisnice}, it follows that \[\mathcal K \in \xi^{-12}\eta^{-4}\phgd(\Mp),\] where $\mathcal E = (0,0,0,E_{\log,0})$. Thus
\[\tilde{\mathcal K} := (\xi^{12} \mathcal K)|_{\mathbf{if} \n \{\digamma > 0\}} \in C^{\infty}(\mathbf{if} \n \{\digamma > 0, \eta > 0\}).\]

\begin{thm}\label{thm:C6:blowup}Let $T$ be any tracefree $\mathring{\slash{g}}$-symmetric type $(1,1)$ tensor on $S^2$, and let $U \subseteq \{T \neq 0\}$ be compactly contained. Then for any $v_0 > 0$ there exists a class $\mathcal C$ of commutative short-pulse data, open under a class of $C^\infty$ perturbations, such that, for the formal solution $g$ to $\Ric(g) = 0$ with short-pulse data $\mathbf{T} \in \mathcal C$, whose existence is provided by \cref{thm:C4:continuation}, $\tilde{\mathcal K}(g)|_{\{v \geq v_0, \theta \in U\}}$ blows up as $\digamma \to 0$.

More precisely, for $\delta_1,\delta_2 > 0$ let $\mathcal C_{\delta_1,\delta_2}$ denote the set of commutative short pulse tensors $\mathbf{T} = \psi \mathbf{T}_0$ satisfying 
\[\supp \psi \subseteq [0,\delta_1]\] and
\[\norm{\mathbf{T}_0-T}_{C^0(S^2)} < \delta_2, |\norm{\pa_v \psi}^2_{L^2([0,1])}-1| < \delta_2^2\]
(using $\mathring{\slash{g}}$ to define the $C^0$ norm). Then for $\delta_1,\delta_2$ sufficiently small depending only on $\norm{T}_{C^0(S^2)}$, with $g$ the formal solution with short-pulse data given by $\mathbf{T}$, $\tilde{\mathcal K}(g)|_{\{v \geq v_0, \theta \in U\}}$ blows up as $\digamma \to 0$. In fact, 
\[\tilde{\mathcal K}(g)|_{\{v \geq v_0, \theta \in U\}} \gtrsim \digamma^{-3}\] (where the implied constant depends on $\mathbf{T}$).
\end{thm}
Observe that the class $\mathcal C_{\delta_1,\delta_2}$ consists of data having a very sharp ``pulse'' near $v=0$, and is open under perturbations supported in $[0,\delta_1]$.

Before beginning the proof, we need some notation and to state some results about polyhomogeneity of metric components at $\{\digamma = 0\}$ for non-generic commutative data, which are of independent interest. These are essentially translations of \cref{thm:C5:behaviour} to the current setting (more specifically translations of \cref{thm:C5:prop3} since we are away from $\{v = 0\}$ and the zeroes of $\mathbf{T}_0$) Fix $\mathbf{T} = \psi \mathbf{T}_0$, an arbitrary commutative short-pulse tensor. Set $v^\ast = \sup_v \{\pa_v \psi(t) = 0\text{ on }[0,v]\}$.\footnote{Observe that $v^\ast < 1$ if we assume that $\mathbf{T}_0$ is not trivial, i.e.\ is not identically $0$.} Then on $(v^\ast,1]\times W$, the energy $\mathbf{E}$ is nonzero. Let $W = \{\mathbf{T}_0 \neq 0\}$. Thus the function $\gamma:\mathbf{if}\n \{\theta \in W\} \to \R$ defined by
\[\gamma(v,\theta) = \frac{\sqrt{2}}{\sqrt{\int_0^v \mathbf{E}(s)\ ds}}\] is well-defined and smooth. The dependence of $\gamma$ on $\theta$ will be treated parametrically, and we will often omit writing this dependence explicitly.
Let us set $\sigma = \gamma(v) - \eta$, a bdf of $\{\digamma = 0\}$ for $v > v^\ast, \theta \in W$. Observe that $\gamma, -\pa_v \gamma > 0$ on $\{v > v^\ast, \theta \in U\}$, and $\pa_v \sigma = \pa_v \gamma$, $\pa_\eta \sigma = -1$. Let us set $X = \{v > v^\ast, \sigma \geq 0, \theta \in W\}$.

Let $g$ denote the solution to $\Ric(g) = 0$ in series with data given by $\mathbf{T}$. Recall from \cref{C:C4:existence} our notation for the metric components $\slash{k} = \xi^{-4}\slash{g}$, $\omega = \log \Omega$, $V = \xi^{3}\eta(L - \mathring{L})$, and from \cref{C:C4:topo} \[\pa_v \slash{k} = \slash{k}(f/2 + F), \ \pa_\eta \slash{k} = \slash{k}(h/2 + H),\]
as well as $\tilde{f} = \eta^{-2}f$ and $\tilde{F} = \eta^{-1}F$.
Recall also from \cref{C:C5:blowup} the index sets 
\begin{align*}n' &= n \un \{n+1,n+2,\ldots\}\times \{1\}\\
n'' &= n \un \{n,n+1,\ldots,\}\times \{1\} \un \{n+1,n+2,\ldots\}\times \{2\}\times\{n+1,n+2,\ldots\}.\end{align*}
Then:
\begin{thm}\label{thm:C6:phg}
\begin{align*}
f,h \in \phgi{(\{\sigma = 0\})}{(-1)}(X)\\
F,H \in \phgi{(\{\sigma = 0\})}{(-1')}(X)\\
\omega \in \phgi{(\{\sigma = 0\})}{(0'')}(X).
\end{align*}
\end{thm}
If $\mathbf{T}$ satisfied the assumption of the ``generic'' commutative data considered in \cref{C:C5:intro}, then this would follow immediately from \cref{thm:C5:behaviour}. The proof of the theorem in general follows from techniques introduced in \cref{C:C5:vz}, except there will be a few details by having to localize away from $\{v = v^\ast\}$. We postpone the details until \cref{C:C6:kernel}. For the moment, however, the skeptical reader may suppose $\mathbf{T}_0$ has simple zeroes and replace $\mathcal C_{\delta_1,\delta_2}$ with its intersection with generic short-pulse data. Because they have explicit formulae, the proof for $f$ and $h$ is easy, however.
\begin{proof}[Proof of \cref{thm:C6:phg} for $f$, $h$]
From \cref{thm:C4:fh}, $f,h$ have the explicit formulae
\begin{align}
\label{eq:C6:fhs}
\begin{split}
f &= \frac{4}{\sigma}\frac{\eta^2 \pa_v \gamma}{\gamma^2 + \eta\gamma}\\
h &= -\frac{4}{\sigma}\frac{\eta}{\gamma + \eta},
\end{split}
\end{align}
and so the statement about $f$, $h$ follows.

\end{proof}
\begin{rk}All sets and quantities above, $X$, $\digamma$, $\sigma$, $f$, $h$, $X$, etc.\ ultimately depend only on the choice of short-pulse tensor $\mathbf{T}$. If we wish to emphasize this dependence, we will write $\digamma[\mathbf{T}]$, $\sigma[\mathbf{T}]$, etc.\end{rk}

It is important to examine the top-order behaviour of the metric quantities. From \eqref{eq:C6:fhs},
\begin{align}
\label{eq:C6:fht}
\begin{split}
f &= \frac{2\pa_v \gamma}{\sigma} + C^{\infty}(X)\\
h &= -\frac{2}{\sigma} + C^{\infty}(X).
\end{split}
\end{align}

From \cref{thm:C6:phg}, we may also write
\begin{align*}
F &= \frac{F_{-1}}{\sigma} + \phgi{(\{\sigma = 0\})}{E_{\log}}(X)\\
H &= \frac{H_{-1}}{\sigma} + \phgi{(\{\sigma = 0\})}{E_{\log}}(X),
\end{align*}
for some tracefree $\mathring{\slash{g}}$-symmetric $F_{-1}, \ H_{-1} \in C^{\infty}(X)$.
However, from \eqref{eq:C5:FH}, $\pa_\eta F = \pa_v H$, and thus there is some tracefree $\mathring{\slash{g}}$-symmetric $\Sigma \in C^{\infty}(X)$ for which $F_{-1} = \pa_v \gamma \Sigma$ and $H_{-1} = -\Sigma$. Of course we may write $\Sigma = \Sigma' \mathbf{T}_0$ for some smooth function $\Sigma'$. Thus
\begin{align}
\label{eq:C6:FHt}
\begin{split}
F &= \frac{\pa_v \gamma\Sigma}{\sigma} + \phgi{(\{\sigma = 0\})}{E_{\log}}(X)\\
H &= \frac{-\Sigma}{\sigma} + \phgi{(\{\sigma = 0\})}{E_{\log}}(X).
\end{split}
\end{align}

Using \eqref{eq:topomega:eq1}, \eqref{eq:C6:fht} and \eqref{eq:C6:FHt} and the above, we may write
\begin{equation}\label{eq:C6Llastomega}\pa_\eta\pa_v \omega = -\frac{\pa_v \gamma}{4\sigma^2} + \frac{\pa_v \gamma}{8\sigma^2}\Tr(\Sigma^2) + \sigma^{-1}\phgi{(\{\sigma = 0\})}{E_{\log}}(X).\end{equation}
Using \cref{thm:C6:phg}, we know a priori that
\[\omega = \omega_{0,1}\log \sigma + \omega_{0,0} + \sigma\phgi{(\{\sigma = 0\})}{E_{\log}}(X),\]
where $\omega_{0,1},\omega_{0,0} \in C^{\infty}(X)$. It follows that \[\omega_{0,1} = \frac{\Tr(\Sigma^2)-2}{8}.\]
Thus \begin{equation}\label{eq:C6:Omega^2}\Omega^2 = C\sigma^{\frac{1}{4}(\Tr(\Sigma^2)-2)},\end{equation} where $1 \lesssim C \lesssim 1$ locally uniformly on $X$.

We may now outline the proof of \cref{thm:C6:blowup}. We will start with an examination of $\tilde{\mathcal K}[\mathbf{T}]$.
\begin{prop}\label{thm:C6:Kretsch}
Fix commutative short-pulse data $\mathbf{T} = \psi \mathbf{T}_0$, and let $X = X[\mathbf{T}]$. Let $g = g[\mathbf{T}]$ be the solution metric. Then
\begin{equation}\label{eq:C6:curveform}\tilde{\mathcal K}(g)|_{X} = \frac{(\pa_v \gamma)^2}{2\eta^2 \Omega^4\sigma^4}\left(\Tr((1-\Sigma^2)^2) + o(1)\right),\end{equation}
where $o(1)$ indicates a quantity converging to $0$ as $\sigma \to 0$ uniformly on compact subsets of $X$.\footnote{$\tilde{\mathcal K}$ should have scaling $\eta^{-4}$ at $\eta \to 0$, even though \eqref{eq:C6:curveform} suggests that it has scaling $\eta^{-2}$. This scaling is hidden in the $o(1)$ term, since $\{\digamma = 0\}$ is disjoint from $\{\eta = 0\}$.}
\end{prop} Proving this will require a detailed examination of the Kretschmann scalar in a double-null gauge and the leading-order asymptotics of various metric components, so we postpone the proof until \cref{C:C6:Kretsch}

From \cref{thm:C6:blowup} it is easy to prove:
\begin{cor}\label{thm:C6:blowupcond}
Fix non-trivial commutative short-pulse data $\mathbf{T} = \psi \mathbf{T}_0$, and let $v^\ast = v^\ast[\mathbf{T}]$, and $U \subseteq \{\mathbf{T}_0 \neq 0\}$ be compact, and fix $\epsilon > 0$ arbitrary (provided $v^\ast + \epsilon < 1$). Let $\Sigma = \Sigma[\mathbf{T}]$, and suppose that $\Sigma|_{\{\sigma = 0,\ v \geq v^\ast + \epsilon\}}$ does not have eigenvalue $\pm 1$.\footnote{Since $\Sigma$ is tracefree, having eigenvalue $+1$ is equivalent to having eigenvalue $-1$.} Then
\[\tilde{\mathcal K}(g)|_{\{v \geq v^\ast + \epsilon, \theta \in U\}} \gtrsim \sigma^{-3}.\]
In particular $\tilde{\mathcal K}(g)$ blows up as $\sigma \to 0$.
\end{cor}
\begin{proof}
Notice that $\Tr((1-\Sigma^2)^2)$ is only zero if $\pm 1$ is an eigenvalue of $\Sigma$. Since $\Sigma$ is smooth, and does not have $\pm 1$ as an eigenvalue, it follows that $\Tr((1-\Sigma^2)^2) \gtrsim 1$ on $\{v \geq v^\ast + \epsilon, \theta \in U\}$. We know from \eqref{eq:C6:Omega^2} that $\Omega^{4} \lesssim \sigma^{\frac{1}{2}(\Tr(\Sigma^2)-2)}$. Now $\frac{1}{2}(\Tr(\Sigma^2)-2) \geq -1$, and so $\Omega^4\sigma^4 \leq \sigma^3$ for $\sigma \leq 1$. The conclusion now follows from \cref{thm:C6:Kretsch}.
\end{proof}

So, to prove \cref{thm:C6:blowup}, it suffices to show for $0 < \delta_1,\delta_2$ small enough and $\mathbf{T} \in \mathcal C_{\delta_1,\delta_2}$, $v^\ast(\mathbf{T}) < v_0$, that $\Sigma[\mathbf{T}]|_{\{\sigma = 0\}}$ does not have eigenvalue $\pm 1$ for $v \geq v_0$ and $\theta \in U \subseteq \{\mathbf{T}_0 \neq 0\}$ compact.

While this may seem easy to do since the equation for $F,H$, \eqref{eq:C5:FH}, is linear in $F$ and $H$, the coefficients of the equation, as well as the position of the boundary $\{\sigma = 0\}$ depend on $f$, $h$ which depend nonlinearly on $\mathbf{T}_0$. Thus, changing the data slightly for $F,H$ will create a corresponding change in the data for $f,h$, which could a priori have the effect of cancelling out the change in the data. Thus, we need to understand the linear equation for $F,H$ (or rather $\tilde{F} = F/\eta,H$) for fixed $f,h$, but with data not necessarily coupled to that of $\mathbf{T}_0$.

Therefore, let us look look at \eqref{eq:C4:FtH} for commutative short-pulse data, so that the commutator drops out. Write $\tilde{F} = \tilde{F}'\mathbf{T}_0$ and $H = H' \mathbf{T}_0$. This is the equation
\begin{align}
\begin{split}0 &= \pa_\eta \tilde{F}' + \frac{1}{4\eta}fH' + \frac{1}{4}h\tilde{F}'\\
0 &= \pa_v H + \frac{f}{4}H' + \left(\frac{\eta}{4}h-1\right)\tilde{F}'.\end{split}
\end{align}
with data $H|_{\{v = 0\}} = 0$ and $\tilde{F}|_{\{\eta = 0\}} = \pa_v \mathbf{T} = \pa_v \psi \mathbf{T}_0$.

Motivated by this, let us introduce the equation for $\Psi,\Phi$, which we pose for $\theta \in U \subseteq S^2$
\begin{align}
\label{eq:C6:PP}
\begin{split}0 &= \pa_\eta \Psi + \frac{1}{4\eta}f\Phi + \frac{1}{4}h\Psi\\
0 &= \pa_v \Phi + \frac{f}{4}\Phi + \left(\frac{\eta}{4}h-1\right)\Psi.\end{split}
\end{align}
We will have reason to specify several different initial data for this equation, so we do not make any specific assumption at this point.

We may use the results of \cref{C:A3:kernel} to find the integral kernel of the system \eqref{eq:C6:PP} (i.e.\ the forwards fundamental solution). Namely, we will show in \cref{C:C6:kernel}:
\begin{prop}
\label{thm:C6:ffs}
Fix short-pulse data $\mathbf{T}$ and $U \subseteq S^2$. Let $f = f[\mathbf{T}]$, $h = h[\mathbf{T}]$, $\digamma = \digamma[\mathbf{T}]$ and let $(A(v,\eta;s,t,\theta) = A[\mathbf{T}](v,\eta;s,t,\theta), B[\mathbf{T}](v,\eta;s,t,\theta) = B(v,\eta;s,t,\theta))$ be a solution to \eqref{eq:C6:PP} on \[\{v \geq s,\ \eta \geq t, \ \theta \in U, \ \digamma(v,\eta,\theta), \digamma(s,t,\theta) > 0\}\] with data
\begin{align*}
A(v,t;s,t,\theta) &\equiv 0\\
B(s,\eta;s,t,\theta) &= \left(-\frac{\eta h(s,\eta)}{4}+1\right)\sqrt{\frac{\digamma(s,t,\theta)}{\digamma(s,\eta,\theta)}}
\end{align*}
(the second set of arguments $(s,t,\theta)$ should be thought of as parameters). Then $A,B$ exist and are smooth and any solution $(\Psi,\Phi)$ to \eqref{eq:C6:PP} with data $\Phi_{\{v = 0\}} = 0$, $\Psi|_{\{\eta = 0\}} = \Psi_0$ admits the representation formula
\begin{align}
\label{eq:C6:repform}
\begin{split}
\Psi(v,\eta,\theta) &= \frac{2}{\sqrt{\digamma(v,\eta)}}\Psi_0(v) + \int_0^v A(v,\eta;s,0,\theta)\Psi_0(s)\ ds\\
\Phi(v,\eta,\theta) &= \int_0^v (-B(v,\eta;s,\eta,\theta) + B(v,\eta;s,0,\theta))\Psi_0(v).
\end{split}
\end{align}
\end{prop}
Like $F$ and $H$, $\sigma A$ and $\sigma B$ remain bounded even as $\sigma \to 0$. Indeed:
\begin{prop}
\label{thm:C6:boundkernel}
Fix commutative short-pulse data $\mathbf{T} = \psi \mathbf{T}_0$, and $U \subseteq\{\mathbf{T}_0 \neq 0\}$ compact.
Suppose $v^\ast = v^\ast[\mathbf{T}] < 1$, and set $\gamma = \gamma[\mathbf{T}]$, $\sigma = \sigma[\mathbf{T}]$, etc., as above. Suppose $v_0 > v^\ast$, and $v^\ast < v_1 < v_0$. Then
\[\sup_{\{v \geq v_0, \sigma \geq 0, 0 \leq s \leq v_1 , \theta \in U\}} \sigma |A[\mathbf{T}]|(v,\eta;s,0,\theta) + \sigma |B[\mathbf{T}]|(v,\eta;s,0,\theta) \leq C,\]
where $C$ depends only on 
\begin{equation}\label{eq:C6:dependence}(v_0-v_1)^{-1}, \ \sup_{\theta \in U}\max\left(\mathbf{E}(1),\mathbf{E}(v_1)^{-1},\int_0^1 \mathbf{E}(s)\ ds, \left(\int_0^{v_1} \mathbf{E}(s)\ ds\right)^{-1}\right),\end{equation}
and is increasing in both its arguments. Observe that since $v_1 > v^\ast$ and $U$ is compact, \eqref{eq:C6:dependence} is finite.
\end{prop}
We will prove this as well in \cref{C:C6:kernel}.
We may combine these results to prove \cref{thm:C6:blowup}:
\begin{proof}[Proof of \cref{thm:C6:blowup}]
First choose any $0 < v_1 < v_0$, and fix $\delta_1,\delta_2 > 0$ (how small will become clear during the proof). Notice that for $\mathbf{T} \in \mathcal C_{\delta_1,\delta_2}$ and $\delta_1$ sufficiently small, $v^\ast = v^\ast[\mathbf{T}] \leq \delta_1 < v_1$, and $U \subseteq \{\mathbf{T}_0 \neq 0\}$. Set $Y = \{v \geq v_0, \ \theta \in U, \ \digamma \geq 0\}$. Then from \cref{thm:C6:ffs}, for any $\mathbf{T} = \psi \mathbf{T}_0 \in \mathcal C_{\delta_1,\delta_2}$:
\begin{equation}\label{eq:C6:gotobed}F(v,\eta,\theta) = \frac{2}{\sqrt{\digamma(v,\eta)}} \pa_v \psi(v) \mathbf{T}_0(\theta) + \int_0^vA(v,\eta;s,0,\theta)\pa_v \psi(s)\mathbf{T}_0(\theta)\ ds,\end{equation}
where $A = A[\mathbf{T}]$ is defined in \cref{thm:C6:ffs}. Notice that for $v \geq \delta_1$, $\pa_v \psi(v) \equiv 0$, and so the first term in \eqref{eq:C6:gotobed} vanishes on $Y$.

Let us write $A'(v,\eta;s,\theta) := \sigma A(v,\eta;s,0,\theta)$. It follows from \eqref{eq:C6:FHt} that on $Y$,
\[\pa_v \gamma(v,\eta,\theta)\Sigma(v,\eta,\theta) = \int_0^v A'(v,\eta;s,\theta)\pa_v \psi(s)\mathbf{T}_0(\theta)\ ds + o(1),\]
where the $o(1)$ is uniform in $Y$ as $\sigma \to 0$ (since $Y$ is compact).

Thus
\begin{align*}
|\Sigma(v,\eta,\theta)| &\leq |\pa_v \gamma|^{-1}|\mathbf{T}_0(\theta)|\int_0^{\delta_1} |A'(v,\eta;s,\theta)\pa_v \psi(s)|\ ds + o(1)\\
&\leq \norm{|\pa_v \gamma|^{-1}}_{L^{\infty}(Y)}\norm{\mathbf{T}_0}_{L^{\infty}(U)}\norm{A'(v,\eta;s,0)}_{L^{\infty}(Y)}\sqrt{\delta_1}\norm{\pa_v \psi}_{L^2([0,\delta_1])} + o(1)\\
&\leq \norm{|\pa_v \gamma|^{-1}}_{L^{\infty}(Y)}\norm{\mathbf{T}_0}_{L^{\infty}(U)}\norm{A'(v,\eta;s,0)}_{L^{\infty}(Y)}\sqrt{\delta_1}\sqrt{1+\delta_2^2}+ o(1).
\end{align*}

It follows that if $\norm{A'(\cdot,\cdot;s)}_{L^{\infty}(Y)}, \norm{|\pa_v \gamma|^{-1}}_{L^{\infty}(Y)} \lesssim 1$ uniformly for $\mathbf{T} \in \mathcal C_{\delta_1,\delta_2}$, $\theta \in U$, and $\delta_1 > 0$ small enough, then $|\Sigma|_{\{\sigma = 0\}}|$ can be made as small as desired for $\delta_1$ sufficiently small (and $\delta_2 \lesssim 1$). In particular, $\Sigma$ it will not have eigenvalue $\pm 1$, and thus by \cref{thm:C6:blowupcond}, $\tilde{K}$ blows up at least like $\sigma^{-3} \sim \digamma^{-3}$ as $\sigma \to 0$ (recall that $\sigma$ and $\digamma$ are bdfs of the same compact hypersurface).

To control $\pa_v \gamma$, notice first that
\[\frac{\mathbf{E}(v_1)}{\left(\int_0^{1}\mathbf{E}(s)\ ds\right)^{3/2}} \lesssim |\pa_v\gamma|\] for $v \geq v_1$.
To control $A'$, we use \cref{thm:C6:boundkernel}.

Therefore, to control both, it suffices to show that
\[\mathbf{E}(1),\ \int_0^1 \mathbf{E}(s)\ ds,\ (\mathbf{E}(v_1))^{-1},\ \left(\int_0^{v_1} \mathbf{E}(s)\ ds\right)^{-1}\] are all uniformly bounded for $\mathbf{T} \in \mathcal C_{\delta_1,\delta_2}$, $\theta \in U$, and $\delta_2 > 0$ small enough. By definition for $v \geq \delta_1$,
\begin{align*}\mathbf{E}(v) &= \frac{1}{2}\Tr(\mathbf{T}_0^2)\norm{\pa_v \psi}_{L^2([0,\delta_1])}^2.\end{align*}
Thus
\begin{equation}\label{eq:C6:Ebound}\frac{1}{2}\Tr(\mathbf{T}_0^2)\left(1-\delta_2^2\right)\lesssim \mathbf{E}(v) \lesssim \frac{1}{2}\Tr(\mathbf{T}_0^2)(1+\delta_2^2).\end{equation}
If $\delta_2$ is small enough then $1 \lesssim \Tr(\mathbf{T}_0^2) \lesssim 1$ on $U$ uniformly for $\mathbf{T} \in \mathcal{C}_{\delta_1,\delta}$, and so \eqref{eq:C6:Ebound} shows that $\mathbf E(1), (\mathbf E(v_1))^{-1} \lesssim 1$.

Now
\[\int_0^1 \mathbf{E}(s)\ ds \leq \frac{1}{2}\Tr(\mathbf{T}_0^2)\int_0^1 \norm{\pa_v \psi}^2_{L^2([0,\delta_1])}\ ds \leq \frac{1}{2}\Tr(\mathbf{T}_0^2)(1+\delta_2^2),\]
and so similarly $\int_0^1 \mathbf{E}(s) \lesssim 1$ if $\delta_2$ is small enough. Finally, using \eqref{eq:C6:Ebound},
\[\int_0^{v_1} \mathbf{E}(s)\ ds\geq \int_{\delta_1}^{v_1} \mathbf{E}(s)\ ds \geq \frac{1}{2}\Tr(\mathbf{T}_0^2)\left(1-\delta_2^2\right)(v_1-\delta_1),\]
and so $\left(\int_0^{v_1} \mathbf{E}(s)\ ds\right)^{-1} \lesssim 1$ for $\delta_2$ small enough.
\end{proof}

\section{Analyzing the Kretschmann scalar}
\label{C:C6:Kretsch}
We start by describing the Kretschmann scalar in a double-null gauge.
\begin{lem}\label{thm:C6:Kdesc}Let $g$ be a short-pulse metric on a subset of $\overline{\mathcal M}$ in a double-null gauge (with definitions as in \cref{C:C4}). Recall
\begin{align*}
N &= -2\Omega^2 \grad \und{u} = \pa_u = \xi^{-2}\left(\eta\pa_\eta - 1/2\xi\pa_\xi\right)\\
L &= -2\Omega^2 \grad u = \pa_{\und{u}} = \xi^{-4}\eta^{-2}(\pa_v + \xi\eta V).
\end{align*}
Then
\begin{align}
\begin{split}\label{eq:C6:Kdesc}\mathcal K(g) &= R_{abcd}R_{ijkl}\slash{g}^{ai}\slash{g}^{bj}\slash{g}^{ck}\slash{g}^{dl} - \frac{4}{\Omega^2}R_{Nbcd}R_{Ljkl}\slash{g}^{bj}\slash{g}^{ck}\slash{g}^{dl} - \frac{1}{\Omega^4}R_{NLcd}R_{NLkl}\slash{g}^{ck}\slash{g}^{dl}\\
&+ \frac{2}{\Omega^4}R_{NbNd}R_{LkLl}\slash{g}^{bj}\slash{g}^{dl} + \frac{2}{\Omega^4}R_{NbLd}R_{NlLj}\slash{g}^{bj}\slash{g}^{dl} - \frac{1}{\Omega^6}R_{NLNd}R_{LNLl}\slash{g}^{dl} + \frac{1}{4\Omega^{16}}R^2_{NLNL}.\end{split}\end{align}
Here, roman letters denote arguments in the fibred $TS^2$, so $R_{abcd}$ denotes a type $(4,0)$ fibre tensor, $R_{Nbcd}$ denotes a type $(3,0)$ fibre-tensor, etc.
\end{lem}
The proof is a routine computation, so is omitted.

We may now prove \cref{thm:C6:Kretsch}.
\begin{proof}[Proof of \cref{thm:C6:Kretsch}]
Let $\slash{k}_0 = \slash{k}|_{\mathbf{if}}$. We start off by recording the curvature components to top order at $\mathbf{if}$ and at $\{\sigma = 0\}$ (at least over $\mathbf{if}\n X)$). We will derive the formulae from \cref{C:A1:ComputeCurvature}. We will use the notation $q \in o(\sigma^x)$, for $x \in \R$ to indicate a quantity $q$ such that $\sigma^{-x}q \to 0$ as $\sigma \to 0$ uniformly on compact subsets of $X$. We will show
\begin{align}
\label{eq:C6:curv}
\begin{split}
(\xi^{-2}\slash{k}_0^{-1}R_{abcd})|_{\mathbf{if}\n X} &= -\pa_v\gamma \frac{\eta^{-1}}{8\sigma^2\Omega^2}\left(\slash{k}(1+\Sigma)\owedge(\slash{k}(1+\Sigma))\right) + o(\sigma^{-2})\\
(\xi^{-1}R_{abcL})|_{\mathbf{if}\n X} &= 0\\
(\xi R_{abcN})|_{\mathbf{if}\n X} &= 0\\
(\xi^{-2}R_{abNL})|_{\mathbf{if}\n X}&= 0+ o(\sigma^{-2})\\
(\xi^4\slash{k}_0^{-1}R_{aLbL})|_{\mathbf{if}\n X} &= \frac{1}{4}\eta^{-4}(\pa_v \gamma)^2\sigma^{-2}(1-\Sigma^2) + o(\sigma^{-2})\\
(\slash{k}_0^{-1}R_{aNbN})|_{\mathbf{if}\n X} &= \frac{1}{4}\eta^2\sigma^{-2}(1-\Sigma^2) + o(\sigma^{-2})\\
(\slash{k}_0^{-1}\xi^2R_{LaNb})|_{\mathbf{if}\n X} &= -\frac{1}{4}\eta^{-1}(\pa_v\gamma)\sigma^{-2}(1-\Sigma^2) + o(\sigma^{-2})\\
(\xi^{-5}R_{LNLa})|_{\mathbf{if}\n X} &= 0\\ (\xi^{-3}R_{NLNc})|_{\mathbf{if}\n X} &= 0\\
(\xi^{-2}R_{abNL})|_{\mathbf{if}\n X} & = 0 + o(\sigma^{-2})\\
(\xi^6 R_{NLNL})|_{\mathbf{if}\n X} &=\Omega^2\left(\frac{1}{4}\eta^{-1}((\pa_v\gamma)\sigma^{-2}(2-\Tr(\Sigma^2))^2 + o(\sigma^{-2})\right)
\end{split}
\end{align}
Before we derive these formulas, let us show how to use them to prove \cref{thm:C6:Kretsch}. Since $\Sigma$ is $\mathring{\slash{k}}$ symmetric, we may plug in the formulae \eqref{eq:C6:curv} into the formula given by \cref{thm:C6:Kdesc} to obtain
\begin{align}
\label{eq:C6:adone}
\begin{split}
(\eta^2 \xi^{12}\Omega^4\sigma^4 \mathcal K)|_{\mathbf{if} \n X} &= \frac{(\pa_v \gamma)^2}{64}\left|\slash{k}(1+\Sigma)\owedge(\slash{k}(1+\Sigma))\right|_{\slash{k}}^2 + \frac{(\pa_v \gamma)^2}{4}\Tr((1-\Sigma^2)^2)\\
&+\frac{(\pa_v\gamma)^2}{16}(2-\Tr(\Sigma^2)) + o(1).
\end{split}
\end{align}

Let us now simplify the norm of the Kulkarni-Nomizu product in the first term. Since $1+\Sigma$ is $\slash{k}$-symmetric, we may at at point $p$ choose a $\slash{k}$-orthonormal basis $E_1,E_2$ of $T_p S^2$ which diagonalizes $1+\Sigma$. Write $\slash{k}(1+\Sigma) = \Upsilon$, a symmetric $(1,1)$ tensor, for which $\Upsilon_{ij} = 0$ if i$i \neq j$, and $\Upsilon_{ii} = \lambda_i$, the eigenvalue of $1+\Sigma$ associated to $E_i$.

Since $T_pS^2$ is two-dimensional, up to symmetry the only non trivial component of $\Upsilon\owedge\Upsilon$ is
\[(\Upsilon\owedge\Upsilon)_{1212}.\] So in this orthonormal basis,
\[|(\Upsilon\owedge\Upsilon)|^2_{\slash{k}} = 4(\Upsilon\owedge\Upsilon)_{1212}^2.\]
Using the definition of the Kulkarni-Nomizu product,
\[(\Upsilon\owedge\Upsilon)_{1212} = 2\Upsilon(E_1)\Upsilon(E_2) = 2\lambda_1\lambda_2 = (\lambda_1+\lambda_2)^2-(\lambda_1^2+\lambda_2^2).\]
We recognize this expression as
\[(\Tr(1+\Sigma))^2 - \Tr((1+\Sigma)^2),\]
which because $\Tr(\Sigma) = 0$ simplifies to
\[2-\Tr(\Sigma^2).\]
Thus
\[\left|\slash{k}(1+\Sigma)\owedge(\slash{k}(1+\Sigma))\right|_{\slash{k}}^2 = 4(2-\Tr(\Sigma)^2)^2,\]
and so \eqref{eq:C6:adone} simplifies to
\begin{equation}\label{eq:C6:fff}(\eta^2 \xi^{12}\Omega^4\sigma^4 \mathcal K)|_{\mathbf{if} \n X} = \frac{(\pa_v \gamma)^2}{8}(2-\Tr(\Sigma^2)) + \frac{(\pa_v \gamma)^2}{4}\Tr((1-\Sigma^2)^2) + o(1).\end{equation}
Since $\Sigma^2$ is tracefree and $S^2$ is 2-dimensional, $\Sigma^2$ is a scalar, and so $\Tr(\Sigma^2) = 2\Sigma^2$ and thus \eqref{eq:C6:fff} simplifies to
\[(\eta^2 \xi^{12}\Omega^4\sigma^4 \mathcal K)|_{\mathbf{if} \n X} = \frac{(\pa_v \gamma)^2}{2}\Tr((1-\Sigma^2)^2),\]
which is \eqref{eq:C6:curveform}

Let us now derive the formulae in \eqref{eq:C6:curv}. Let $\chi^N = \frac{1}{2}\Lie_N \slash{g}$ and $\chi^L = \frac{1}{2}\Lie_L \slash{g}$ denote the second fundamental forms of $\slash{g}$ in the $N$ and $L$ directions, respectively. We record first their top order behaviour, as well as that of $[N,L]$ at $\mathbf{if} = \{\xi = 0\}$. Recall from \cref{C:C4:topo} that $V|_{\mathbf{if}} = 0$.

Using \cref{thm:C4:algebra} to handle the non-linear operations, we may conclude
\begin{align}
\label{eq:C6:comp}
\begin{split}
\chi^L &= \frac{\eta^{-2}}{2}\pa_v \slash{k} = \frac{\eta^{-2}}{2}\slash{k}_0\left(\frac{f}{2} + F\right) + \xi \eta^{-2}\phgi{(\mathbf{if})}{E_{\log}}(\Mp;\Sym^2 T^\ast S^2)\\
\chi^N &= \frac{\xi^2}{2}\eta\pa_\eta \slash{k} - \xi^2 \slash{k} = \xi^2\slash{k}_0\left(\eta \frac{h}{4} -1 + \eta H\right) + \xi^3\phgi{(\mathbf{if})}{E_{\log}}(\Mp;\Sym^2 T^\ast S^2)\\
[N,L] & = [N,\xi^{-3}\eta^{-1}V] \in \xi^{-4}\eta^{-1}\phgi{(\mathbf{if})}{E_{\log}}(\Mp;\Sym^2 TS^2).
\end{split}
\end{align}
We now begin the computation, using the previous formulae and \cref{C:A1:ComputeCurvature}.

\begin{romanumerate}
\item $R_{abcd}$: Notice that $\slash{R}(\slash{g})_{abcd} = \xi^4 \slash{R}(\slash{k})_{abcd} \in \xi^4\phgi{(\mathbf{if})}{E_{\log}}(\Mp)$. Thus,
\[R_{abcd} = \frac{1}{2\Omega^2}(\chi^L\owedge \chi^N) + \xi^3\eta^{-2}\phgi{(\mathbf{if})}{E_{\log}}(\Mp).\]

Thus,
\[(\xi^{-2}\slash{k}_0^{-1}R_{abcd})|_{\mathbf{if}\n X} = -\pa_v\gamma \frac{\eta^{-1}}{8\sigma^2\Omega^2}\left(\slash{k}(1+\Sigma)\owedge(\slash{k}(1+\Sigma))\right) + o(\sigma^{-2}).\]

\item $R_{abcL}\text{ and }R_{abcN}$: It is clear from \eqref{eq:C6:comp} that
\[(\xi^{-1}R_{abcL})|_{\mathbf{if}}, (\xi R_{abcN})|_{\mathbf{if}} = 0.\]

\item $R_{abNL}$: We have
\[R_{abNL} = (\chi^L \times \chi^N - \chi^N\times \chi^L) + \xi^3\eta^{-2}\phgi{(\mathbf{if})}{E_{\log}}(\Mp),\]
where $\times$ is taken with respect to $\slash{g}$.
Since on $X$,
\begin{align*}4\xi^2\eta^2\slash{h}^{-1}(\chi^L\times \chi^N -\chi^N\times \chi^L)&= [f/2+F,\eta h/2 -1 + \eta H]\\
&= [F,H]\\
&= -\frac{\pa_v\gamma\eta}{\sigma^2}[\Sigma,\Sigma] + o(\sigma^{-2}) \in o(\sigma^{-2}),\end{align*}
it follows that
\[(\xi^{-2}R_{abNL})|_{\mathbf{if}\n X} \in o(\sigma^{-2}).\]

\item $R_{aLbL}$:We have
\[R_{aLbL} = \chi^L \times \chi^L - \Lie_L \chi^L + \xi^{-3}\eta^{-4}\phgi{(\mathbf{if})}{E_{\log}}(\Mp).\]

On the other hand, on $X$,
\begin{align*}
4\eta^4\xi^4\slash{k}_0^{-1}(\chi^L \times \chi^L - \Lie_L \chi^L) &= (f/2+F)^2 - 2(\pa_v f/2 + \pa_v F)-2(f/2+F)^2\\
&=(\pa_v \gamma)^2\sigma^{-2}(1-\Sigma^2 + o(\sigma^{-2})).
\end{align*}
Thus
\[(\xi^4\slash{k}_0^{-1}R_{aLbL})|_{\mathbf{if}\n X}= \frac{1}{4}\eta^{-4}(\pa_v \gamma)^2\sigma^{-2}(1-\Sigma^2) + o(\sigma^{-2}).\]

\item $R_{aNbN}$: In the same fashion as above, 
\[R_{aNbN} = \chi^N \times \chi^N - \Lie_N \chi^N + \xi\phgi{(\mathbf{if})}{E_{\log}}(\Mp),\]
and on $X$
\[4\slash{k}_0^{-1}( \chi^N \times \chi^N - \Lie_N \chi^N) = \eta^2\sigma^{-2}(1-\Sigma^2 + o(\sigma^{-2})).\]
Thus
\[(\slash{k}_0^{-1}R_{aNbN})|_{\mathbf{if}\n X} = \frac{1}{4}\eta^2\sigma^{-2}(1-\Sigma^2) + o(\sigma^{-2}).\]

\item $R_{LaNb}$: Notice that $\slash{\Hess}_{\slash{g}} = \slash{\Hess}_{\slash{k}}$, and thus $\slash{\Hess}_{\slash{g}}\Omega^2 \in \phgi{(\mathbf{if})}{E_{\log}}(\Mp)$.
Thus
\[R_{LaNb} = (\chi^L \times \chi^N) -\frac{1}{2}(\Lie_N \chi^L + \Lie_L \chi^N) + \xi^{-1}\eta^{-2}\phgi{(\mathbf{if})}{E_{\log}}(\Mp).\]

In the same fashion as above, on $X$,
\[4\xi^2\eta^2\slash{k}_0^{-1}((\chi^L \times \chi^N) -\frac{1}{2}(\Lie_N \chi^L + \Lie_L \chi^N)) = -\eta\pa_v \gamma\sigma^{-2}(1-\Sigma^2) + o(\sigma^{-2}).\]
Thus
\[(\slash{k}_0^{-1}\xi^2R_{LaNb})|_{\mathbf{if}\n X} = -\frac{1}{4}\eta^{-1}(\pa_v\gamma)\sigma^{-2}(1-\Sigma^2) + o(\sigma^{-2}).\]

\item $R_{LNLa}\text{ and }R_{NLNc}$: Since $\slash{\nabla}^{\slash{g}} = \slash{\nabla}^{\slash{k}}$, \[\slash{\nabla}^{\slash{g}} N \in \xi^{-2}\phgi{(\mathbf{if})}{E_{\log}}(\Mp), \ \slash{\nabla}^{\slash{g}} N \in \xi^{-4}\phgi{(\mathbf{if})}{E_{\log}}(\Mp).\]
It follows that
\[(\xi^{-5}R_{LNLa})|_{\mathbf{if}\n X}, (\xi^{-3}R_{NLNc})|_{\mathbf{if}\n X} = 0.\]

\item$ R_{NLNL}$: Since $N = \xi^{-2}\eta\pa_\eta + \xi \phgi{(\mathbf{if})}{E_{\log}}(\Mp)$ and $L = \xi^{-4}\eta^{-2}v + \xi^{-3}\eta^{-2}\phgi{(\mathbf{if})}{E_{\log}}(\Mp)$,
it follows that
\[\xi^6\eta^2(LN+NL)\log \Omega = 2\pa_v\eta \pa_\eta \omega + \xi\phgi{(\mathbf{if})}{E_{\log}}(\Mp).\]

However, from \eqref{eq:C6Llastomega}, on $X$,
\[2\pa_v\eta\pa_\eta \omega = \frac{\eta\pa_v \gamma}{4}\sigma^{-2}(\Tr(\Sigma^2)-2) + o(\sigma^{-2}).\]

Thus
\[(\xi^6 R_{NLNL})|_{{\mathbf{if}\n X}} = \Omega^2\left(\frac{1}{4}\eta^{-1}(\pa_v\gamma)\sigma^{-2}(2-\Tr(\Sigma^2)) + o(\sigma^{-2})\right).\]
\end{romanumerate}
\end{proof}

\section{The integral kernel and other loose ends}
\label{C:C6:kernel}
We start by proving \cref{thm:C6:ffs}.
\begin{proof}[Proof of \cref{thm:C6:ffs}]
This is an application the results of \cref{C:A3:kernel}. In the notation of \cref{thm:A3:ffs}, we need to find the $E_{ij}$ ($i,j = 1,2$) and $I$, $J$, and then apply \cref{thm:A3:repform}.

For clarity, we do not include the parameter $\theta \in S^2$ in our notation. In our case, $I$ satisfies
\[\pa_\eta I(\eta;s,t) + \frac{h(s,\eta)}{4}I(\eta;s,t) = 0,\]
with data $I(t;s,t) = 1$. Therefore using \cref{thm:C4:fh},
\[I = \frac{\sqrt{\digamma(s,t)}}{\sqrt{\digamma(s,\eta)}}.\]
Similarly,
\[\pa_v J(v;s,t) + \frac{f(v,t)}{4}J(v;s,t) = 0,\]
with data $J(s;s,t) = 1$, and therefore
\[J = \frac{\sqrt{\digamma(s,t)}}{\sqrt{\digamma(v,t)}}.\]
Now, $(E_{11},E_{21})$ solve \eqref{eq:C6:PP} with data \begin{align*}E_{11}(v,t;s,t) &= 0\\
E_{21}(s,\eta;s,t) &= -\left(\frac{\eta h}{4}-1\right)\frac{\sqrt{\digamma(s,t)}}{\sqrt{\digamma(s,\eta)}}.\end{align*} Thus by definition, $(E_{11},E_{21}) = (A,B)$. In particular by \cref{thm:A3:ffssmooth}, $A$, $B$ are smooth in the appropriate region.

By assumption $\Phi|_{\{v = 0\}} = 0$, and thus by \cref{thm:A3:repform}
\begin{align*}
\Psi(v,\eta) = \frac{\sqrt{\digamma(v,0)}}{\sqrt{\digamma(v,\eta)}}\Psi_0(s) + \int_{0}^v A(v,\eta;s,0)\Psi_0(s)\ ds\\
\Phi(v,\eta) = \int_0^v (-B(v,\eta;s,\eta) + B(v,\eta;s,0))\Psi_0(s)\ ds.
\end{align*}
Since $\digamma(s,0) = 4$, this shows \eqref{eq:C6:repform}. 
\end{proof}

We now start to prove \cref{thm:C6:boundkernel}. Until otherwise noted, let hereafter $v^\ast$, $v_1$, $v_0$, $U$, denote the same quantities as in the statement of \cref{thm:C6:boundkernel}. Also set
\[\eta^\ast = \sup_{\theta\in U} \gamma(v_0) = \sup \{\eta \: v \geq v_0, \ \sigma(v,\eta) \geq 0\}.\]
We will prove \cref{thm:C6:boundkernel} by showing two general boundedness statements about arbitrary solutions $\Psi$, $\Phi$ to \eqref{eq:C6:PP}, one to show boundedness on a region $\{\sigma \leq \epsilon, \eta \leq \gamma(v_0)\}$, and another to show boundedness all the way to $\sigma = 0$. We will then apply this to the particular case of $A$ and $B$.

\begin{figure}[htbp]
\centering
\includegraphics{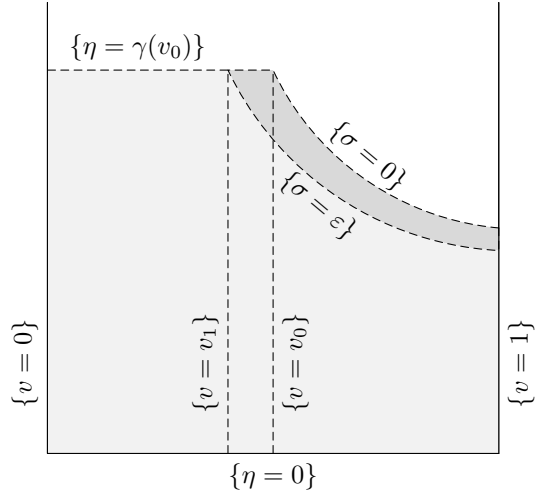}
\caption{The regions in the proof of \cref{thm:C6:boundkernel}. The light shaded region is the region in which the estimates from \cref{thm:C6:helperi} are valid, and the dark shaded region is the region in which the estimates from \cref{thm:C6:helperii} are valid. In the case of this figure, it is assumed that $v_2 = 0$.}
\label{Fig:C6:regions}
\end{figure}

We start with the first:
\begin{prop}\label{thm:C6:helperi}Fix $0 \leq v_2 \leq v_1$, $\epsilon > 0$ and let $(\Psi,\Phi)$ be a bounded solution to \eqref{eq:C6:PP} on $Y = Y_{\epsilon,v_2} = \{v \geq v_2, \ \eta \leq \gamma(v_0), \ \sigma \geq \epsilon, \ \theta \in U\}$ with data $\Psi|_{\{\eta = 0\}} = \Psi_0$ and $\Phi|_{\{v = v_2, \eta \leq \gamma(v_0)\}} = \Phi_{v_2}$.
Then
\[\norm{\Psi}_{L^{\infty}(Y)} + \norm{\Phi}_{L^{\infty}(Y)} \leq C(\norm{\Psi_0}_{L^{\infty}(\{\eta = 0\})} + \norm{\Phi_{v_2}}_{L^{\infty}(\{v = v_2, \eta \leq \gamma(v_0)\})}),\]
where
\[C = C(\eta^\ast, \norm{\tilde{f}}_{L^{\infty}(Y)}, \norm{h}_{L^{\infty}(Y)})\]
is increasing in all its arguments.
Observe that $C$ does not depend on $v_2$.
\end{prop}
\begin{proof}
Observe that if $(v,\eta,\theta) \in Y$, then $[v_1,v]\times [0,\eta]\times U \subseteq Y$. Therefore we may use the results of \cref{C:A3:Goursat}. More precisely, we use a two-dimensional version of Gronwall's inequality, \cref{thm:A3:Gronwall}, which lets us derive the desired bound from the inequalities
\begin{align*}
|\Psi(v,\eta,\theta)| &\leq \norm{\Psi_0}_{L^{\infty}(\{\eta = 0\}}) +\\
&\frac{1}{4}(\norm{\tilde{f}}_{L^{\infty}(Y)}+\norm{h}_{L^{\infty}(Y)})\int_0^\eta |\Phi(v,t,\theta)| + |\Psi(v,t,\theta)| dt\\
|\Psi(v,\eta,\theta)| &\leq \norm{\Phi_{v_2}}_{L^{\infty}(\{v = v_2, \eta \leq \gamma(v_0)\})}\\
&+\frac{1}{4}(\eta^\ast\norm{\tilde{f}}_{L^{\infty}(Y)}+\eta^\ast\norm{h}_{L^{\infty}(Y)}+ 1)\int_{v_2}^v |\Phi(s,\eta,\theta)| + |\Psi(s,\eta,\theta)| dt.
\end{align*}
\end{proof}

Before stating the theorem about boundedness to $\{\sigma = 0\}$, we need to discuss domains of dependence.
\begin{lem}\label{thm:C6:suffsmall}For $\epsilon > 0$ sufficiently small, \[Z = Z_\epsilon = \{\sigma \leq \epsilon,\ v \geq v_1, \eta \leq \gamma(v_0), \theta \in U\}\]is a domain of dependence according to \eqref{eq:C6:PP} and \[H = H_\epsilon = \{\sigma = \epsilon,\ v \geq v_1, \ \eta \leq \gamma(v_0), \ \theta \in U\}\] is a Cauchy hypersurface for $Z$, i.e.\ every backwards-directed integral curve of $\pa_v$ or $\pa_\eta$ starting in $Z$ stays in $Z$ until it intersects $H$.\end{lem} \begin{proof}Since $\gamma(1,\theta) > 0$ for $\theta \in U$ and $\gamma$ is strictly decreasing in $v$, we may find $0 < \epsilon$ so that
\begin{gather*}\epsilon < \inf_{\theta \in U} \gamma(1)\\
\epsilon < \inf_{\theta \in U}\gamma(v_1) -\gamma(v_0).\end{gather*} 

Let us show that for this choice of $\epsilon$, $H$ is a Cauchy hypersurface for $Z$.
We first show that $\sigma(v_1,\eta) > \epsilon$ for any $\eta \leq \gamma(v_0)$.
Since $\gamma$ is decreasing as a function of $v$, so is $\sigma$. Then by definition,
\[\sigma(v_1,\eta) = \gamma(v_1)-\eta > \epsilon + \gamma(v_0)-\eta = \epsilon + \sigma(v_0,\eta) \geq \epsilon.\]

Suppose $p = (v,\eta,\theta) \in Z$. We will now show that a backwards-directed integral curve of $\pa_v$ starting at $p$ remains in $Z$ until it intersects $H$. Such an integral curve is of the form $t \mapsto (v-t,\eta,\theta)$ for $t > 0$. 

Since $\sigma$ is increasing along such an integral curve, $\sigma(p) \leq \epsilon$, and $\lim_{v' \to v^\ast} \sigma(v',\eta) = \infty$, there exists $t^\ast$ such that $\sigma(v-t^\ast,\eta,\theta) = \epsilon$. We need to show that $v-t \geq v_1$ for all $0 \leq t \leq t^\ast$. Indeed, if instead for some $t$ it holds that $v-t < v_1$, then \[\epsilon \geq \sigma(v-t,\eta,\theta) \geq \sigma(v_1,\eta,\theta) > \epsilon,\] a contradiction.

Next we will show that a backwards directed integral curve of $\pa_\eta$ starting at $p$ remains in $Z$ until it intersects $H$. Such an integral curve is of the form $t \mapsto (v,\eta-t,\theta)$ for $t > 0$.

Since $\sigma$ is increasing along such an integral curve, and $\sigma(p) \leq \epsilon$, \[\sigma(v,\eta-t,\theta) \geq \sigma(1,0,\theta) = \gamma(1) > \epsilon\], there exists $t^\ast$ such that $\sigma(v,\eta-t^\ast,\theta) = \epsilon$. Since certainly \[\sigma(v,\eta-t,\theta) \leq \sigma(v,\eta-t^\ast,\theta) = \epsilon\] for $0 \leq t \leq t^\ast$, and $v \geq v_1$ by definition, this means that a backwards-directed integral curve remains in $Z$ until it intersects $H$.
\end{proof}

From \eqref{eq:C6:fht}, we know that to top order $f$ is $\frac{2 \pa_v \gamma}{\sigma}$ and $h = -\frac{2}{\sigma}$. We will let $\bar{f}$ and $\bar{h}$ denote the ``regular'' parts of $f$ and $h$, respectively,which by \eqref{eq:C6:fhs} are explicitly given by
\begin{align}
\label{eq:C6:bfh}
\begin{split}
\bar{f} := \frac{1}{\pa_v \gamma}f - \frac{2}{\sigma} = -\frac{6\eta + 2\sigma}{\gamma^2+\gamma\eta} \in C^{\infty}(\{\sigma \geq 0, v > v^\ast, \theta \in U\})\\
\bar{h} := h+ \frac{2}{\sigma}=-\frac{2}{\gamma + \eta}\in C^{\infty}(\{\sigma \geq 0, v > v^\ast, \theta \in U\}).
\end{split}
\end{align}

We now have the statement about boundedness towards $\sigma = 0$, together with a statement about polyhomogeneity which we will use later to show \cref{thm:C6:phg}.
\begin{prop}\label{thm:C6:helperii}Fix $\epsilon > 0$ sufficiently small as specified by \cref{thm:C6:suffsmall}. Let $(\Psi,\Phi)$ be a solution to \eqref{eq:C6:PP} on with data $\Psi|_{H_\epsilon} = \Psi_\epsilon$ and $\Phi_{H_\epsilon} = \Phi_\epsilon$. Assume $\Psi,\Phi \in C^1(Z\n \{\sigma > 0\})$. Then choosing $\epsilon > 0$ perhaps smaller (depending only on some numerical constant), the following bound holds:
\[\norm{\sigma\Psi}_{L^{\infty}(Z)} + \norm{\sigma \Phi}_{L^{\infty}(Z)} \leq C(\epsilon\norm{\Psi_\epsilon}_{L^{\infty}(H_\epsilon)} + \epsilon\norm{\Phi_\epsilon}_{L^{\infty}(H_\epsilon)}),\]
where
\[C = C\left(\eta^\ast,1/\eta^\ast, \norm{\bar{f}}_{L^{\infty}(Z)}, \norm{\bar{h}}_{L^{\infty}(Z)}, \sup_{Z} |\pa_v \gamma|, \sup_Z |\pa_v \gamma|^{-1}\right)\]
is increasing in all its arguments. If $\Psi,\Phi \in C^{\infty}(Z \n \{\sigma > 0\})$, then in fact, \[\Psi, \Phi \in \phgi{(\{\sigma = 0\})}{(-1')}(Z).\]
\end{prop}

\begin{proof}
Set $S = -\frac{\sigma\eta}{\pa_v \gamma}\Psi$ and $T = \sigma \Phi$. Then $S$, $T$ satisfy
\begin{align}
\label{eq:C6:sobored}
\begin{split}
-\sigma\pa_v T &= \frac{-\pa_v \gamma }{2}\left(1+\frac{\sigma\bar{f}}{2}\right)T - \frac{-\pa_v \gamma}{2}\left(1-\frac{\sigma\bar{h}}{2} + \frac{2\sigma}{\eta}\right)S\\
-\sigma\pa_\eta S &=-\frac{1 }{2}\left(1+\frac{\sigma\bar{f}}{2}\right)T + \frac{1}{2}\left(1+\frac{\sigma\bar{h}}{2} - \frac{\sigma}{\eta}\right)S.
\end{split}
\end{align}

It suffices to show 
\[\norm{T}_{L^{\infty}(Z)} + \norm{S}_{L^{\infty}(Z)} \leq C'(\norm{T}_{L^{\infty}(H_\epsilon)} + \epsilon\norm{S}_{L^{\infty}(H_\epsilon)}),\]
where
\[C' = C'\left(\norm{\bar{f}}_{L^{\infty}(Z)}, \norm{\bar{h}}_{L^{\infty}(Z)}, \sup_{Z} |\pa_v \gamma|^{-1}\right)\]
is increasing in all its arguments.

The proof of \cref{thm:C5:initialvzs} now applies nearly verbatim to give the desired bounds. Applying the rest of the propositions in \cref{C:C5:vz} required to prove \cref{thm:C5:prop3}:(i) nearly verbatim shows the polyhomogeneity. Indeed, the vector fields are written $-\sigma\pa_v$ and $-\sigma\pa_\eta$ so that they propagate \emph{towards} $\{\sigma = 0\}$, just like in \cref{C:C5:vz}, and $-\pa_v \gamma > 0$, which is why we insist on writing \eqref{eq:C6:sobored} in that form. We may need to shrink $\epsilon$ so that $\epsilon < 1$, for technical reasons.
\end{proof}

We are finally in a position to prove \cref{thm:C6:boundkernel}.
\begin{proof}[Proof of \cref{thm:C6:boundkernel}]
Fix $\epsilon > 0$ sufficiently small as in the \cref{thm:C6:helperii}, and let $Y = Y_{\epsilon,0}$ be as in \cref{thm:C6:helperi}, and $Z = Z_\epsilon$ be as in \cref{thm:C6:helperii}. Combining \cref{thm:C6:helperi} and \cref{thm:C6:helperii} and applying them to the pair $(A(v,\eta;s,0),B(v,\eta;s,0))$ shows that
\[\sup_{\{v \geq v_0, \sigma > 0, s \in [0,v_1],\theta \in U\}} \sigma |A|(v,\eta;s,0,\theta) + \sigma |B|(v,\eta;s,0,\theta) \leq C,\]
where
\begin{align*}C = C\left(\eta^\ast,1/\eta^\ast, \norm{\bar{f}}_{L^{\infty}(Z)}, \norm{\bar{h}}_{L^{\infty}(Z)}, \sup_{Z} |\gamma|^{-1}, \sup_{Z} |\pa_v \gamma|, \sup_Z |\pa_v \gamma|^{-1},\right.\\
\left.\norm{\tilde{f}}_{L^{\infty}(Y)}, \norm{h}_{L^{\infty}(Y)},\sup_{\{0 \leq s \leq v_1\}}\norm{B|_{\{ v = s\}}}_{L^{\infty}(\{v = s, \eta \leq \eta^\ast\})}\right).\end{align*}
We need to show that we may control these quantities in terms of \eqref{eq:C6:dependence}.\footnote{For the remainder of this proof, we say that a quantity $a$ is ``controlled'' by a quantity $b$ if there exists a non-decreasing, potentially nonlinear, function $C$ such that $|a| \lesssim C(b)$. We choose to use the word ``control'' instead of ``bound'' to emphasize that $C$ is not necessarily linear, i.e.\ it may not be true that $a \lesssim b$.}

First, it is clear from its definition that for $v_1 \leq v \leq 1$,
\begin{align*}
\frac{1}{\sqrt{\int_0^{1}\mathbf{E}(s)\ ds}} \lesssim |\gamma| &\lesssim \frac{1}{\sqrt{\int_0^{v_1}\mathbf{E}(s)\ ds}}\\
\frac{\mathbf{E}(v_1)}{\left(\int_0^{1}\mathbf{E}(s)\ ds\right)^{3/2}} \lesssim |\pa_v\gamma| &\lesssim \frac{\mathbf{E}(1)}{\left(\int_0^{v_1}\mathbf{E}(s)\ ds\right)^{3/2}},
\end{align*}
where the implied constants are numerical constants.
Thus $|\gamma|$, $|\gamma|^{-1}$, $|\pa_v \gamma|$, $|\pa_v \gamma|^{-1}$ on $v_1 \leq v \leq 1$ are all controlled by \eqref{eq:C6:dependence}, and in particular so are their respective suprema on $Z$.

Let us now control $\bar{f}$ and $\bar{h}$.

Since $\eta^\ast = \sup_{\theta}\gamma(v_0)$, we have also controlled $\eta^\ast, 1/\eta^\ast$.

We will need to control $\epsilon$. By \cref{thm:C6:suffsmall} and \cref{thm:C6:helperii}, we may take
\[\epsilon = \frac{1}{2}\min\left(c,\inf_{\theta \in U} \gamma(1),\inf_{\theta \in U} \gamma(v_1)-\gamma(v_0)\right),\]
for some numerical constant $c$. Thus $\epsilon$ is bounded by $c$. Let us now control $\epsilon$ from below (equivalently $\epsilon^{-1}$ from above). We have already controlled $\gamma(1)$ from below. To control $\gamma(v_1)-\gamma(v_0)$, we use the mean value theorem
\[\gamma(v_1)-\gamma(v_0) = \pa_v \gamma(v')(v_1-v_0)\] for some $v' \in [v_1,v_0]$. Since $\pa_v \gamma < 0$, and we have controlled $|\pa_v \gamma|$ from below, this means that we have controlled $\gamma(v_1)-\gamma(v_0)$ from below by $v_0-v_1$ and the desired energy quantities.

Let us now control $\bar{f}$ and $\bar{h}$ on $Z$. From \eqref{eq:C6:bfh}, these are both controlled from above by $1/\gamma$, $\eta^\ast$ and $\epsilon$, which we know are controlled by \eqref{eq:C6:dependence}.

Next we control $\tilde{f}$ and $h$ on $Y$.
We first control them on $Y \n \{v > v^\ast\}$. By \eqref{eq:C6:fhs}, $h$ is controlled by $\epsilon^{-1}$ (since $\eta/(\gamma + \eta) \leq 1$). We have already controlled $\epsilon^{-1}$ by \eqref{eq:C6:dependence}. Since by definition \[\pa_v \gamma/\gamma = -\frac{1}{4}\mathbf{E}\gamma^2,\] it follows from \eqref{eq:C6:fhs} that
\[\tilde{f} = -\eta\mathbf{E}\frac{\gamma^2}{\sigma(\gamma+\eta)}.\] The first factor is controlled by \eqref{eq:C6:dependence}. To control the second factor, rewrite it as
\[\frac{\sigma^2 + 2\sigma \eta + \eta^2}{\sigma^2 + 2\sigma\eta} = 1 + \frac{1}{\sigma}\frac{\eta^2}{\sigma + 2\eta}.\] Since $\sigma \geq \epsilon$, this is controlled by $\epsilon^{-1}$ and $\eta$, which are in turn controlled by \eqref{eq:C6:dependence}.
On $Y \n \{v < v^\ast\}$, $\tilde{f},h \equiv 0$, ad so they are certainly controlled.

Now let us control 
\[B(s,\eta;s,0) = \left(-\frac{\eta h(s,\eta)}{4}+1\right)\frac{\sqrt{2}}{\sqrt{\digamma(s,\eta)}}\] for $\eta \leq \gamma(v_0)$ and $0 \leq s \leq v_1$. The first factor is already controlled. If $s \leq v^\ast$, then $\digamma \equiv 4$, so the second factor is controlled in this case. Now we control the second factor if $s \geq v^\ast$. By definition
\[\digamma(s,\eta)^{-1} = \frac{\gamma^{2}}{\gamma^2-\eta^2}.\] We have already shown how to control this quantity.
Thus all quantities are controlled by \eqref{eq:C6:dependence}.\end{proof}

Now let us prove \cref{thm:C6:phg}. Observe that it is sufficient to show that for any $U \subseteq \{\mathbf{T}_0 \neq 0\}$ and $v_0 > v^\ast$,
\begin{gather*}\tilde{F},H \in \phgi{(\{\sigma = 0\})}{(-1')}(X\n \{v \geq v_0, \ \theta \in U\})\\
\omega \in \phgi{(\{\sigma = 0\})}{(0'')}(X\n \{v \geq v_0, \ \theta \in U\}).\end{gather*}
Let $v^\ast < v_1 < v_0$ be arbitrary, and set
\[\eta^\ast = \sup_{\theta \in U} \gamma(v_0),\]
as in the statement of \cref{thm:C6:boundkernel}. Choose any $v^\ast < v_0 < v_1$. Then the previous discussion still applies to the quantities $v^\ast$, $v_1$, $v_0$, $U$, $\eta^\ast$.

\begin{proof}[Proof of \cref{thm:C6:phg}]Let $\epsilon > 0$ be be provided by \cref{thm:C6:suffsmall}, and construct the associated Cauchy hypersurface $H_\epsilon$ and region $Z_\epsilon$.

\Cref{thm:C4:FH} implies that $\tilde{F},H \in C^{\infty}(\{\digamma > 0\})$ and so in particular $\tilde{F},H \in C^{\infty}(H_\epsilon)$, and we need only show
\[\tilde{F},H \in \phgi{(\{\sigma = 0\})}{(-1')}(Z_\epsilon).\] But to do so, we need only use \cref{thm:C6:helperii}.

Let us now turn our attention to $\omega$. As the case for $\tilde{F}$ and $H$, we already know $\omega C^{\infty}(\{\digamma > 0\})$, so in particular $\omega \in C^{\infty}(H_\epsilon)$ and we need only show \[\omega \in \phgi{(\{\sigma = 0\})}{(0'')}(Z_\epsilon).\]

From \eqref{eq:topomega:eq1} we may deduce
\begin{equation}\label{eq:C6:omega}(-\sigma\pa_\eta -1)(-\sigma\pa_v)\omega = \frac{fh}{16} + \frac{1}{8}\Tr(FH) + \frac{1}{4\eta}f :=S.\end{equation}

From the parts of this proposition dealing with $f$, $h$, $F$, and $H$, we know \[S \in \phgi{\{\sigma = 0\})}{(0'')}(Z_\epsilon).\] This is analogous to \eqref{eq:C5:omegavzs}, and so the steps in \cref{C:C5:vz} in the proof of \cref{thm:C5:prop3}:(ii) apply nearly verbatim.
\end{proof}

We may also use the techniques of the proof of \cref{thm:C5:initialvzs} and \cref{thm:C6:helperii} to prove \cref{thm:C4:nonblowup}. Let $v_0 > 0$ and $U \subseteq S^2$ be open and such that $\inf_{\theta \in S^2} \mathbf{E}(v_0) > 0$. Let $X = \{v \geq v_0, \theta \in S^2\}$. Notice that $\gamma$ and hence $\sigma$ are well-defined on $X$. We will show will prove is: \begin{prop}Let $|\bullet|_{\slash{k}}$ denote the norm according to $\slash{k}$. Then \[ |F|_{\slash{k}}, |H|_{\slash{k}} \in \sigma^{-1}L^{\infty}_{\loc}(X).\]
\end{prop}
\begin{proof}As above, we may restrict to \[Y = X\n \{v \geq v_0, \ \theta \in U\}.\] Define $S = -\frac{\sigma\eta}{\pa_v \gamma} \tilde{F} = -\frac{\sigma}{\pa_v \gamma} F$ and $T = \sigma H$. Using \eqref{eq:C6:fht} (which is valid even in the non-commutative setting) and \eqref{eq:C4:FH}, $\Phi,\Psi$ satisfy the equation
\begin{align}
\label{eq:nonblowup:eq1}
\begin{split}
-\sigma \pa_v S &= -\frac{\pa_v \gamma}{2}(1+\sigma C^{\infty}(Y))S - \frac{-\pa_v \gamma}{2}(1+\sigma C^{\infty}(Y)T + \frac{-\pa_v \gamma}{2}[S,T]\\
-\sigma \pa_\eta T &= -\frac{1}{2}(1+\sigma C^{\infty}(Y)S + \frac{1}{2}(1+\sigma C^{\infty}(Y)T + \frac{1}{2}[T,S].
\end{split}
\end{align}

Since $S$, $T$ are $\slash{k}$-symmetric, it suffices to show that \[A := \Tr(S^2),\ B := \Tr(T^2) \in L^{\infty}_{\loc}(Y).\] From \eqref{eq:nonblowup:eq1} we obtain as in the proof of \cref{thm:C4:FH}
\begin{align}
\label{eq:nonblowup:eq2}
\begin{split}
-\sigma \pa_v A &= -\pa_v \gamma(1+\sigma C^{\infty}(Y))\ A - (-\pa_v \gamma)(1+\sigma C^{\infty}(Y)) \Tr(\Phi \Psi)\\
-\sigma \pa_\eta B &= -(1+\sigma C^{\infty}(Y))\Tr(\Phi \Psi) + (1+\sigma C^{\infty}(Y))B.
\end{split}
\end{align}

 Since $-\pa_v\gamma > 0$, we may use the Cauchy-Schwarz inequality
 \[\Tr(\Phi\Psi) = \Tr(\Phi^\ast \Psi) \leq \frac{1}{2}(A+B)\] on \eqref{eq:nonblowup:eq1} to obtain, for $\sigma$ sufficiently small, the differential inequalities
\begin{align}
\label{eq:C6:CS}
\begin{split}
-\sigma \pa_v A &\geq \frac{|\pa_v \gamma|}{2}(1+\sigma L^{\infty}(Y))A - \frac{|\pa_v \gamma|}{2}(1+\sigma L^{\infty}(Y))B\\
-\sigma \pa_\eta B &\geq -\frac{1}{2}(1+\sigma L^{\infty}(Y))A +\frac{1}{2}(1+\sigma L^{\infty}(Y))B.
\end{split}
\end{align}

Using \cref{thm:C6:suffsmall}, we may find $\epsilon > 0$ small enough and an associated Cauchy hypersurface $H_\epsilon$ and region $Z_\epsilon$. 

From \cref{thm:C4:FH}, we know that $A,B \in C^\infty(\{\digamma > 0 \})$, and so in particular $A,B \in L^{\infty}(H_\epsilon)$, and we need only show $A,B \in L^{\infty}(Z_\epsilon)$. If the inequality in \eqref{eq:C6:CS} were an equality, then this would follow from \cref{thm:C5:initialvzs} in the same manner as \cref{thm:C6:helperii}. However, the proof of \cref{thm:C5:initialvzs} was via integrating factors and then bounding the right-hand side above, anyway, and so the presence of the inequality does not affect the validity of the result.\footnote{The inequality may appear the wrong way, but this is because we are propagating from $\{\sigma > 0\}$ to $\{\sigma = 0\}$, so the inequality must be backwards when using the vector field $-\sigma\pa_v$ and $-\sigma \pa_\eta$, which propagate in the opposite direction. Alternatively, the inequality is the correct direction if considered propagating using $\sigma\pa_v$ and $\sigma \pa_\eta$ which propagate in the ``correct'' direction.}
\end{proof}

\appendix
\chapter{Connections and curvature in a double-null gauge}
\label{C:A1:ComputeRicci}

\section{Organization}
In this appendix we compute the Ricci curvature of a metric in a double-null foliation, as defined in \cref{C:C2}. Let us introduce some notation. Let $(M,g,u,v)$ be a double-null foliation, and let $\slash{g}$, $L$, $N$, $\Omega$ are the quantities associated to it, i.e.\ $\slash{g} = \iota^\ast g = g|_{TS}$, $-2\Omega^{-2} = g^{-1}(2du,2dv)$, $N = -2\Omega^2\grad v$, $L = -2\Omega^2 \grad u$. Let $S$ be the manifold to which all fibres are diffeomorphic, and write $TS$ for the subbundle of $TM$ of fibre-tangent vectors. 

We will use slash notation to indicate operations performs with respect to $\slash{g}$ rather than $g$. For example $\slash{\tr}$ denotes contraction with $\slash{g}$ and $\slash{\divg}$ is that $\slash{g}$-divergence. We will denote by $\slash{d}$ the fibre exterior derivative.

We will denote by capital Greek letters $\Theta,\Sigma, \Psi$, etc., vectors in $TS$. $\nabla$ will denote the connection of $g$, while $\slash{\nabla}$ will denote the connection of $\slash{g}$. This is the same as the projection of the connection of $g$ onto the fibres $S_{u,v}$. More generally, we will denote by $\slash{\nabla}$ to be the projection of $\nabla$ to $TS$, even when taking an argument which is not in $TS$. We will denote for this section
\begin{align*}
\chi^L(\Theta,\Psi) &= \frac{1}{2}\Lie_L \slash{g}(\Theta,\Psi) = g(\nabla_\Theta L,\Psi) = g(\nabla_\Psi L,\Theta)\\
\chi^N(\Theta,\Psi) &= \frac{1}{2}\Lie_N \slash{g}(\Theta,\Psi) = g(\nabla_\Theta N,\Psi) = g(\nabla_\Psi N,\Theta)\end{align*}
the second fundamental forms in the directions $L,N$, and
\begin{align*}
H^L &= \slash{\tr}\chi^L\\
H^N &= \slash{\tr}\chi^N\end{align*}
the mean curvature in the same respective directions.

In the second subsection, we will state some tensor identities which will be useful. Since these are all tensor identities, they are easily proved by picking a local geodesic frame, and by extending fibre-tangent vector fields defined over a single fibre $S_{u,v}$ for some $u,v$, to be $0$ along the flow of $L$ and $N$. In the third subsection, we will compute the connection $\nabla$ of $g$. In the fourth subsection, we will compute some components of the type $(1,3)$ curvature operator of $g$ which we will think of as an operator
\[R\:  (TM)^{\otimes 3} \to T M.\]
We will compute it acting on vectors spanned by $L,\ N,\ TS$.

We will not compute all of them, only enough to obtain the rank $(0,4)$ curvature tensor. Since all our formulae are symmetric in $L,N$, we will in fact compute fewer formulae explicitly. Recall we use the sign convention for curvature given by
\[R(X,Y,Z) = (\nabla_{Y}\nabla_X - \nabla_X\nabla_Y + \nabla_{[X,Y]})Z.\]

In the fifth subsection, we will compute the all components of the curvature tensor $R(X,Y,Z,W)$, which for us is defined by
\[R(X,Y,Z,W) = g(R(X,Y,Z),W).\] Again, we will compute it in acting between $L$, $N$ and $TS$.

In the sixth subsection, we will compute the Ricci curvature
\[R(X,Y) = \tr R(X,\cdot,Y,\cdot).\]

In the seventh subsection, we will compute integrability conditions on the Ricci curvature implied by the contracted Bianchi identities.

\section{Tensor calculations}
\label{C:A1:tensoroperations}
Recall that for symmetric tensors $\phi,\psi$ on $TS$,
\[(\phi\times \psi)_{ab} := \slash{g}^{cd}\phi_{ac}\psi_{bd}.\]
Stated in a coordinate-free fashion, this is
\[(\phi\times \psi)(\Theta,\Sigma) = \tr(\phi(\Theta,\cdot)\otimes \psi(\Sigma,\cdot)).\]
\begin{lem}For $\Theta,\Psi \in TS$, \[\slash{g}(\slash{\nabla}_{\Theta} L,\slash{\nabla}_{\Psi} N) = (\chi^L\times \chi^N)(\Theta,\Psi)).\]\end{lem}

The following tensor appears often:
\[T(\Theta,\Phi) = \slash{g}([L,\slash{\nabla}_{\Theta}N],\Phi) - \chi^N(\Phi,[L,\Theta]),\]
and similar tensors with $N$ and $L$ swapped, or with all $L$ or all $N$.
\begin{lem}The following formula holds: $T = \Lie_L \chi^N - 2\chi^N\times \chi^L$, and similarly if $L$ replaces $N$, or with $N$ everywhere instead of $L$, or $L$ everywhere instead of $N$.\end{lem}

The following identity also often appears:
\begin{lem}For $\Theta,\Phi,\Psi \in TS$, \[\Lie_N\chi^L(\Theta,\Psi) - \Lie_L\chi^N(\Theta,\Phi) = \frac{1}{2}\left(\slash{g}(\slash{\nabla}_{\Theta}[N,L],\Phi) + \slash{g}(\slash{\nabla}_\Phi [N,L],\Theta)\right).\]\end{lem}
We also have:
\begin{lem}The following formula holds: $\slash{\tr}(\Lie_L\chi^N) = LH^N + 2\slash{g}(\chi^L,\chi^N)$, and similarly if $L$ replaces $N$, or with $N$ everywhere instead of $L$, or $L$ everywhere instead of $N$.\end{lem}

We will denote
\[\slash{R}(\Theta,\Psi,N) = \slash{\nabla}_{\Psi}\slash{\nabla}_{\Theta} N  - \slash{\nabla}_\Theta\slash{\nabla}_{\Psi} - \slash{\nabla}_{[\Theta,\Psi]}N.\]
\begin{lem}
For $\Theta,\Phi,\Psi \in TS$, \[\slash{g}(\slash{R}(\Theta,\Phi,N),\Psi) = (\slash{\nabla}_{\Phi})\chi_N(\Psi,\Theta)-(\slash{\nabla}_{\Theta})\chi^N(\Phi,\Psi).\] In particular, the 1-form
\[\Theta \mapsto \slash{\tr} \slash{g}(\slash{R}(\Theta,\cdot),N),\cdot)\]
is equal to
\[\slash{\div}\chi^N - \slash{d}H^N,\]
and similarly if $L$ replaces $N$.\end{lem}

\section{The connection}
\label{C:A1:connection}
Recall that $N = -2\Omega^2 \grad v$. Thus,
\[g(\nabla_{\Theta} N,L)  = g(\nabla_{\Theta} (-2\Omega^2 \grad v), L),\] where we can now use the fact the the Hessian is symmetric, i.e.\ \[g(\nabla_{\Theta} \grad v,L) = g(\nabla_L \grad v, \Theta),\] to conclude that
\[g(\nabla_{\Theta} N,L) = g(\nabla_L,\Theta) - 2\Theta \Omega^2,\]
and similarly
\[g(\nabla_{\Theta}L,N) = g(\nabla_N,\Theta) - 2\Theta \Omega^2.\]
By metric compatibility of $\nabla$ and the above,
\[g(\nabla_L N,\Theta) + g(\nabla_N L ,\Theta) = 2\Theta \Omega^2.\] Since $\nabla$ is torsion free,
\[g(\nabla_L N,\Theta) - g(\nabla_L N,\Theta) = \slash{g}([L,N],\Theta).\]
The previous two identities imply
\begin{align*}
g(\nabla_L N,\Theta) &= \Theta(\Omega^2) + \frac{1}{2}\slash{g}([L,N],\Theta)\\
g(\nabla_N L,\Theta) &= \Theta(\Omega^2) + \frac{1}{2}\slash{g}([N,L],\Theta).\end{align*}

Now we may compute:
\begin{align*}
\nabla_{\Theta}\Psi &= \slash{\nabla}_{\Theta}\Psi  + \frac{1}{2\Omega^2}\chi^L(\Theta,\Psi)N + \frac{1}{2\Omega^2}\chi^N(\Theta,\Psi) L\\
\nabla_{\Theta} N &= \slash{\nabla}_{\Theta} N+ \left(\frac{1}{4\Omega^2}\slash{g}([N,L],\Theta) + \Theta \log(\Omega)\right)N\\
\nabla_{\Theta} L & =\slash{\nabla}_{\Theta} L + \left(\frac{1}{4\Omega^2}\slash{g}([L,N],\Theta) + \Theta \log(\Omega)\right)L\\
\nabla_{N} \Theta &= \slash{\nabla}_{\Theta} N + [N,\Theta] + \left(\frac{1}{4\Omega^2}\slash{g}([N,L],\Theta) + \Theta \log(\Omega)\right)N\\
\nabla_{N} N &= (2N\log \Omega)N\\
\nabla_{N} L & = \slash{\grad}\Omega^2 + \frac{1}{2}[N,L]\\
\nabla_{L} \Theta &= \slash{\nabla}_{\Theta} L + [L,\Theta] + \left(\frac{1}{4\Omega^2}\slash{g}([L,N],\Theta) + \Theta \log(\Omega)\right)L\\
\nabla_{L} N  &= \slash{\grad}\Omega^2 + \frac{1}{2}[L,N]\\
\nabla_{L} L &= (2L \log \Omega)L
\end{align*}

\section{The curvature tensor I}
We may compute
\begin{align*}
R(N,L,N)& = 2N\log \Omega \slash{\grad}\Omega^2 - (N\log \Omega)[N,L]  - [N,\slash{\grad}\Omega^2]+ \frac{1}{2}[N,[N,L]]\\ &+\frac{3}{2}\slash{\nabla}_{[N,L]}N-\slash{\nabla}_{\slash{\grad}\Omega^2}N\\
 &+\left(LN\log\Omega + NL\log\Omega - 2|\slash{d}\Omega|^2 + \frac{3}{8\Omega^2}|[N,L]|^2\right)N\\
R(N,\Theta,N) &=2N\log \Omega \slash{\nabla}_{\Theta}N - \slash{\nabla}_{\slash{\nabla}_\Theta N}N - [N,\slash{\nabla}_\Theta N] + \slash{\nabla}_{[N,\Theta]}N\\
&+\left(\Theta N \log \Omega +\frac{N\log\Omega}{2\Omega^2}\slash{g}([N,L],\Theta)-\frac{3}{4\Omega^2}\chi^N([N,L],\Theta)\right.\\
& \qquad\qquad- \left.\frac{x}{4\Omega^2}\slash{g}([N,[N,L]],\Theta)-\slash{\nabla}_\Theta N \log \Omega)\right)N\\
R(N,\Theta,L) &= \slash{\nabla}_{\Theta}\slash{\grad}\Omega^2 + \frac{1}{2}\slash{\nabla}_{\Theta}[N,L] - \slash{\nabla}_{\slash{\nabla}_{\Theta} L}N - [N,\slash{\nabla}_\Theta L] + \slash{\nabla}_{[N,\Theta]}L\\
 &\qquad\qquad+ \frac{1}{4\Omega^2}\slash{g}([N,L],\Theta)\slash{\grad}\Omega^2+\frac{1}{8\Omega^2}\slash{g}([N,L],\Theta)[N,L]\\
 &\qquad\qquad -\Theta \log \Omega \slash{\grad} \Omega^2 -\frac{1}{2}\Theta \log\Omega [N,L]\\
&+ \left(\frac{1}{2\Omega^2}\chi^N(\Theta,\slash{\grad}\Omega^2) + \frac{3}{4\Omega^2}\chi^N(\Theta,[N,L]) - \frac{1}{2\Omega^2}N\log\Omega \slash{g}([N,L],\Theta) \right.\\
&\qquad\qquad+ \left.\frac{1}{4\Omega^2}\slash{g}([N,[N,L]],\Theta) -\Theta N \log\Omega\right)L\\
R(N,L,\Theta) &= \slash{\nabla}_{\slash{\nabla}_\Theta N}L - \slash{\nabla}_{\slash{\nabla}_\Theta L}N + [L,\slash{\nabla}_\Theta N] -[N,\slash{\nabla}_\Theta L] + \slash{\nabla}_{[N,\Theta]}L - \slash{\nabla}_{[L,\Theta]}N\\ - [[N,L],\Theta]
&+ \frac{1}{2\Omega^2}\slash{g}([N,L],\Theta)\slash{\grad}\Omega^2-\Theta\log\Omega[N,L] + \slash{\nabla}_{[N,L]}\Theta\\
&+\left(-\slash{\nabla}_\Theta L \log \Omega + \Theta L\log\Omega - \frac{1}{2\Omega^2}L\log\Omega \slash{g}([N,L],\Theta)\right.\\
&\qquad\qquad \left.+ \frac{3}{4\Omega^2}\chi^L([N,L],\Theta) + \frac{1}{4\Omega^2}\slash{g}([L,[N,L]],\Theta)\right)N\\
&+\left(\slash{\nabla}_\Theta N \log \Omega - \Theta N\log\Omega - \frac{1}{2\Omega^2}N\log\Omega \slash{g}([N,L],\Theta)\right.\\
&\qquad\qquad \left.+ \frac{3}{4\Omega^2}\chi^N([N,L],\Theta) + \frac{1}{4\Omega^2}\slash{g}([N,[N,L]],\Theta)\right)L\\
R(\Theta,\Phi,\Psi) &= \slash{R}(\Theta,\Phi,\Psi) + \frac{1}{2\Omega^2}\chi^L(\Theta,\Psi)\slash{\nabla}_\Phi N + \frac{1}{2\Omega^2}\chi^N(\Theta,\Psi)\slash{\nabla}_\Phi L\\
&\qquad\qquad - \frac{1}{2\Omega^2}\chi^L(\Phi,\Psi)\slash{\nabla}_\Theta N - \frac{1}{2\Omega^2}\chi^N(\Phi,\Psi)\slash{\nabla}_\Theta L\\
&+\left(\frac{1}{2\Omega^2}\slash{g}(\slash{R}(\Theta,\Phi,L),\Psi)-\frac{\Phi\log\Omega}{2\Omega^2}\chi^L(\Theta,\Psi) + \frac{1}{8\Omega^4}\chi^L(\Theta,\Psi)\slash{g}([N,L],\Phi)\right.\\
&\qquad\qquad \left. +\frac{\Theta\log\Omega}{2\Omega^2}\chi^L(\Phi,\Psi) - \frac{1}{8\Omega^4}\chi^L(\Phi,\Psi)\slash{g}([N,L],\Theta) \right)N\\
&\hspace{-1.3pt}+\left(\frac{1}{2\Omega^2}\slash{g}(\slash{R}(\Theta,\Phi,N),\Psi)\hspace{-1.1pt}-\frac{\Phi\log\Omega}{2\Omega^2}\chi^N(\Theta,\Psi)\hspace{-1.1pt} - \frac{1}{8\Omega^4}\chi^N(\Theta,\Psi)\slash{g}([N,L],\Phi)\right.\\
&\qquad\qquad \left. +\frac{\Theta\log\Omega}{2\Omega^2}\chi^N(\Phi,\Psi) + \frac{1}{8\Omega^4}\chi^N(\Phi,\Psi)\slash{g}([N,L],\Theta) \right)L.
\end{align*} We may also compute the formulae which are symmetric in $L$, $N$, but we will not record these.

\section{The curvature tensor II}
\label{C:A1:ComputeCurvature}From the previous formulae, and using the lemmas about tensors proved in \cref{C:A1:tensoroperations} if necessary, we may compute:
\begin{align*}R(\Theta,\Phi,\Psi,\Sigma) &= \slash{R}(\Theta,\Phi,\Psi,\Sigma) + \frac{1}{2\Omega^2}\left((\chi^L\owedge \chi^N)(\Theta,\Psi,\Phi,\Sigma)\right)\\
R(\Theta,\Phi,\Psi,N) &= \hspace{-0.7pt} -\slash{g}(\slash{R}(\Theta,\Phi,N),\Psi)+\hspace{-0.7pt}(\Phi\log\Omega)\chi^N(\Theta,\Psi) + \hspace{-0.7pt} \frac{1}{4\Omega^2}\chi^N(\Theta,\Psi)\slash{g}([N,L],\Phi)\\
&-(\Theta\log\Omega)\chi^N(\Phi,\Psi) - \frac{1}{4\Omega^2}\chi^N(\Phi,\Psi)\slash{g}([N,L],\Theta)\\
R(\Theta,\Phi,\Psi,L) &= -\slash{g}(\slash{R}(\Theta,\Phi,L),\Psi)+(\Phi\log\Omega)\chi^L(\Theta,\Psi) - \frac{1}{4\Omega^2}\chi^L(\Theta,\Psi)\slash{g}([N,L],\Phi)\\
&-(\Theta\log\Omega)\chi^L(\Phi,\Psi)+ \frac{1}{4\Omega^2}\chi^L(\Phi,\Psi)\slash{g}([N,L],\Theta)\\
R(\Theta,\Phi,N,L) &= \frac{1}{2}(\slash{g}(\nabla_\Theta [N,L],\Phi)-\slash{g}(\nabla_\Phi [N,L],\Theta))+(\chi^L\times \chi^N-\chi^N\times \chi^L)(\Theta,\Phi)\\
 &+ (\Phi\log\Omega)\slash{g}([N,L],\Theta) - (\Theta\log\Omega)\slash{g}([N,L],\Phi)\\
R(\Theta,N,\Psi,N) &= (2N\log \Omega)\chi^N(\Theta,\Psi) +  \chi^N\times\chi^N(\Theta,\Psi)-\mathcal L_N \chi^N(\Theta,\Psi)\\
R(\Theta,N,\Psi,L) &= \slash{\Hess}_{\Omega^2}(\Theta,\Psi) +\frac{1}{8\Omega^2}\slash{g}([N,L],\Theta)\slash{g}([N,L],\Psi)\\
&+(\chi^L\times\chi^N)(\Theta,\Psi)
 + \frac{1}{2}(\Psi\log\Omega)\slash{g}([N,L],\Theta) - \frac{1}{2}\Theta\log\Omega\slash{g}([N,L],\Psi)\\
 &-2\Omega^2\Theta\log\Omega\Psi\log\Omega-\frac{1}{2}(\mathcal L_N\chi^L(\Theta,\Psi)+\mathcal L_L \chi^N(\Theta,\Psi))\\
 &+ \frac{1}{4}(\slash{g}(\slash{\nabla}_{\Theta}[N,L],\Psi) - \slash{g}({\nabla}_{\Psi}[N,L],\Theta)) \\
R(\Theta,N,N,L) &= 2\Omega^2\Theta N \log\Omega + (N\log\Omega)\slash{g}([N,L],\Theta)-\frac{3}{2}\chi^N([N,L],\Theta)\\
&-\frac{1}{2}\slash{g}([N,[N,L]],\Theta) - 2\Omega^2\slash{\nabla}_\Theta N \log \Omega\\
R(\Theta,L,\Psi,L) &= (2L\log \Omega)\chi^L(\Theta,\Psi) +\chi^L\times\chi^L\Theta,\Psi)-\mathcal L_L \chi^L(\Theta,\Psi)\\
R(\Theta,L,N,L) &= -2\Omega^2\Theta L \log\Omega + (L\log\Omega)\slash{g}([N,L],\Theta)-\frac{3}{2}\chi^L([N,L],\Theta)\\
&-\frac{1}{2}\slash{g}([L,[N,L]],\Theta) + 2\Omega^2\slash{\nabla}_\Theta L \log \Omega\\
R(N,L,N,L) &= -2\Omega^2(LN\log\Omega + NL\log\Omega) + |\slash{d}\Omega^2|^2 -\frac{3}{4}|[N,L]|^2.
\end{align*}

\section{The Ricci curvature}
One can now compute the Ricci curvature (for instance by using an orthonormal frame for $TS$ according to $\slash{g}$)
\begin{align*}
\Ric(\Theta,\Psi) &= \slash{\Ric}(\Theta,\Psi) + \frac{1}{2\Omega^2}(\chi^L(\Theta,\Psi)H^N + \mathcal L_N \chi^L(\Theta,\Psi)\\
&+ \chi^N(\Theta,\Psi)H^L + \mathcal L_L \chi^N(\Theta,\Psi))\\
&-\frac{1}{\Omega^2}((\chi^L\times\chi^N)(\Theta,\Psi) + (\chi^N\times\chi^L)(\Theta,\Psi))\\
&-\frac{1}{\Omega^2}\Hess_{\Omega^2}(\Theta,\Psi) - \frac{1}{8\Omega^4}\slash{g}([N,L],\Theta)\slash{g}([N,L],\Psi) +2\Theta\log\Omega\Theta\log\Psi\\
\Ric(\Theta,N) &= -\Theta H^N+\slash{\divg} \chi^N(\Theta) + \frac{1}{2\Omega^2}\chi^N(\Theta,[N,L])\\
&+\frac{1}{4\Omega^2}H^N\slash{g}([N,L],\Theta))
-\Theta N \log\Omega - \frac{1}{2\Omega^2}N \log\Omega \slash{g}([N,L],\Theta)\\ &+\frac{1}{4\Omega^2}\slash{g}([N,[N,L]],\Theta)+(\Theta\log\Omega) H^N\\
\Ric(\Theta,L) &=-\Theta H^L+\slash{\divg} \chi^L(\Theta) - \frac{1}{2\Omega^2}\chi^L(\Theta,[N,L])\\
&-\frac{1}{4\Omega^2}H^L\slash{g}([N,L],\Theta))-\Theta L \log\Omega + \frac{1}{2\Omega^2}L \log\Omega \slash{g}([N,L],\Theta) \\&-\frac{1}{4\Omega^2}\slash{g}([L,[N,L]],\Theta)+(\Theta\log\Omega) H^L\\
\Ric(N,N) &= 2(N\log\Omega) H^N - |\chi^N|^2 - NH^N\\
\Ric(L,L) &= 2(L\log\Omega) H^L - |\chi^L|^2 - LH^L\\
\Ric(L,N) &= \slash{\Delta}\Omega^2 - \frac{1}{4\Omega^2}|[L,N]|^2 - \slash{g}(\chi^L,\chi^N)\\
&- \frac{1}{2}(NH^L+LH^N) -(LN\log \Omega + NL\log\Omega).
\end{align*}

\section{The contracted Bianchi identities}
\label{C:A1:computebianchi}
The contracted Bianchi identities,
\[\div \Ric = \frac{1}{2}d\tr \Ric,\]
or in abstract index notation
\[g^{\beta\gamma}\nabla_{\beta}R_{\gamma \alpha} = \nabla_{\alpha} g^{\beta\gamma}R_{\beta\gamma},\]
imply certain integrability conditions on the Ricci curvature. We compute these conditions in this subsection.
First, the scalar curvature is
\[\tr \Ric = \slash{\tr}\Ric_{-,-} -\frac{1}{\Omega^2}\Ric_{NL}.\]
Now we compute enough of $\nabla \Ric$ to compute $\div \Ric$.
\begin{align*}
(\nabla \Ric)(\Theta,\Phi,\Psi) &= (\slash{\nabla}\slash{\Ric})(\Theta,\Phi,\Psi) - \frac{1}{2\Omega^2}\chi^L(\Theta,\Phi)\Ric(N,\Psi)\\
&- \frac{1}{2\Omega^2}\chi^N(\Theta,\Phi)\Ric(L,\Psi) -\frac{1}{2\Omega^2}\chi^L(\Theta,\Psi)\Ric(N,\Phi)\\
&- \frac{1}{2\Omega^2}\chi^N(\Theta,\Psi)\Ric(L,\Phi)\\
(\nabla \Ric)(N,L,\Psi) &= (\Lie_N \Ric(L,-))(\Psi) - \Ric(\slash{\grad}\Omega^2 + 1/2[N,L],\Psi)\\
&- \Ric(L,\slash{\nabla}_\Psi N) - \left(\frac{1}{4\Omega^2}\slash{g}([N,L],\Psi) + \Psi\log\Omega\right)\Ric(N,L)\\
(\nabla \Ric)(L,N,\Psi) &= (\Lie_L \Ric(N,-))(\Psi) - \Ric(\slash{\grad}\Omega^2 - 1/2[N,L],\Psi)\\
&- \Ric(N,\slash{\nabla}_\Psi L) - \left(-\frac{1}{4\Omega^2}\slash{g}([N,L],\Psi) + \Psi\log\Omega\right)\Ric(N,L)\\
(\nabla \Ric)(\Theta,\Phi,N) &= (\slash{\nabla}_{\Theta} \Ric(N,-))(\Phi) - \frac{1}{2\Omega^2}\chi^L(\Theta,\Phi)\Ric(N,N)\\ &- \frac{1}{2\Omega^2}\chi^N(\Theta,\Phi)\Ric(N,L)-\Ric(\Phi,\slash{\nabla}_{\Theta} N)\\
&- \left(\frac{1}{4\Omega^2}\slash{g}([N,L],\Theta)+ \Theta\log\Omega\right)\Ric(N,\Phi)\\
(\nabla \Ric)(N,L,N) &= N\Ric(N,L) - \Ric(\slash{\grad}\Omega^2 + 1/2[N,L],N) - 2N\log\Omega \Ric(N,L)\\
(\nabla\Ric)(L,N,N) &= L\Ric(N,N) - 2\Ric(\slash{\grad}\Omega^2-1/2[N,L],N)\\
(\nabla \Ric)(\Theta,\Phi,L) &= (\slash{\nabla}_{\Theta} \Ric(L,-))(\Phi) - \frac{1}{2\Omega^2}\chi^N(\Theta,\Phi)\Ric(L,L)\\ &- \frac{1}{2\Omega^2}\chi^L(\Theta,\Phi)\Ric(N,L)-\Ric(\Phi,\slash{\nabla}_{\Theta} L)\\
&- \left(-\frac{1}{4\Omega^2}\slash{g}([N,L],\Theta) + \Theta\log\Omega\right)\Ric(L,\Phi)\\
(\nabla \Ric)(L,N,L) &= L\Ric(N,L) - \Ric(\slash{\grad}\Omega^2 - 1/2[N,L],L) - 2L\log\Omega \Ric(N,L)\\
(\nabla\Ric)(N,L,L) &= N\Ric(L,L) - 2\Ric(\slash{\grad}\Omega^2-1/2[N,L],L).
\end{align*}

Let $\Ric_{-,-}$ denote the two fibre two-form obtained by plugging in vectors $\Theta,\Phi \in TS$ into $\Ric$, let $\Ric_{N,-}$, $\Ric_{L,-}$ denote the fibre $1$-forms obtained by plugging in $N$, $L$, respectively, and a vector $\Theta \in TS$, into $\Ric$, and let $\Ric_{NN}$, $\Ric_{LL}$, and $\Ric_{NL}$ denote the functions obtained by plugging in $(N,N)$, $(L,L)$, and $(N,L)$ respectively into $\Ric$.

Plugging in each of $\Psi,N,L$ for $X$ in the equation $\div \Ric(X) = \frac{1}{2}d\tr\Ric\cdot X$ then gives
\begin{align*}
\Lie_L\Ric_{N,-}+\Lie_N \Ric_{L,-} &+H^L\Ric_{N,-} + H^N\Ric_{L,-} = 2\Omega^2\slash{\div}\Ric_{-,-}-\Omega^2\slash{d}\slash{\tr}\Ric_{-,-}\\
&\hspace{4.7cm}+2\Ric_{-,-}\cdot \slash{\grad}\Omega^2+\slash{d}\Ric_{NL}.\\
L\Ric_{NN} +H^L\Ric_{NN}&=-\Omega^2N\slash{\tr}\Ric_{-,-}-2\Omega^2\slash{g}(\chi^N, \Ric_{-,-})-H^N\Ric_{NL}\\
&+2\Omega^2\slash{\div}\Ric_{N,-}+2\Ric_{N,-}\cdot (\slash{\grad}\Omega^2) - \Ric_{N,-}\cdot [N,L].\\
N\Ric_{LL} + H^N\Ric_{LL} &=-\Omega^2L\slash{\tr}\Ric_{-,-}-2\Omega^2\slash{g}(\chi^L, \Ric_{-,-})-H^L\Ric_{LN}\\
&+2\Omega^2\slash{\div}\Ric_{L,-}+2\Ric_{L,-}\cdot (\slash{\grad}\Omega^2) + \Ric_{L,-}\cdot [N,L].
\end{align*}

\chapter{Real blowup and polyhomogeneity}
\label{C:A4}

In this appendix give a brief primer on real blowup and polyhomogeneity. Some classic references for this are \cite{MelTapsit}, \cite{MelDiff} and \cite{MelCalc}. For a gentle introduction, see \cite{GriBasi}.

\section{Real blowup}
\label{C:A4:blowup}

Let $X$ be a manifold with corners (mwc), i.e.\ a manifold modelled locally on the product $[0,\infty)^{n-k}\times \R^k$, for $k \leq n \in \N$. We always assume that the boundary faces of $X$ are embedded, as opposed to immersed. Let $Y$ be a ``reasonable'' submanifold of $X$. Here ``reasonable'' means that $Y$ is a $p$-submanifold, i.e.\ there are coordinates $(x,y) \in U \subseteq [0,\infty)^n\times \R^k$ for $X$ in which $Y$ is the intersection $U \n S$, for some coordinate subspace $S$. We seek to define the blowup $[X,Y]$. In pedestrian terms, this means we take polar coordinates around a point in $Y$ in the directions transverse to $Y$. In more sophisticated language, one glues in the cosphere bundle $SN^\ast Y$ in place of $Y$.

The simplest case is when $X = \R^n_x$ and $Y = \{0\}$. In this case, $[X,Y] \iso S^{n-1}\times [0,\infty)$. Observe that $[X,Y]$ comes equipped with a natural ``blowdown'' map $\beta:[X,Y] \to \R^n$, which is the usual map from polar coordinates back to $\R^n$. Rather than use spherical coordinates, it is convenient to cover the space $[X,Y]$ in $2n$ ``projective'' coordinate charts, each valid in the region $\pm x^i > 0$ for some $i$. These coordinates are given by (without loss of generality fixing $i= 1$ and $x^1 \geq 0$) defining $\xi^j = x^j/x^1$ for $j \neq 1$ and using $(x^1,\xi^2,\ldots,\xi^n) \in [0,\infty)\times \R^{n-1}$ as coordinates. Here, $x^1$ is the radial variable, and $\xi^k$ are the angular variables. Naturally, in these coordinates $\beta(x^1,\xi^2,\ldots, \xi^n) = (x^1,x^1x^2,\ldots,x^1x^n)$. Blowing up a point introduces a boundary face to $X$, namely the set $S^{n-1}\times \{0\}$. This face is traditionally called the ``front face'' and is denoted by $\mathbf{ff}$. In the projective coordinates, $\mathbf{ff} = \{x^1 = 0\}$, and $x^1$ is a boundary-defining function (bdf) for it (i.e.\ a function which vanishes on $\mathbf{ff}$ but whose derivative does not). Notice that $\beta:\mathbf{ff} \to \{0\}$.

One may also perform this construction in the case that $X = \R^{n-k}\times [0,\infty)^k$, in which case one replaces $S^{n-1}$ (thought of as the unit sphere in $\R^n$) by its intersection with the space $\R^{n-k}\times [0,\infty)^k$. 

\begin{figure}[htbp]
\centering
\includegraphics{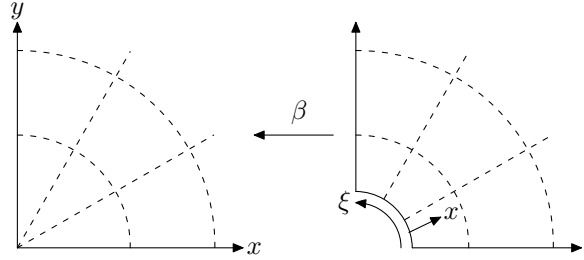}
\caption{Blowing up $\{(0,0)\}$ in $[0,\infty)^2_{(x,y)}$. The dashed lines are curves in $[0,\infty)^2$ and in $[[0,\infty)^2, \{0\}]$. Also labelled are the directions of increase of projective coordinates $(x,\xi= y/x)$.}
\end{figure}

If $X = \R^n \times \Sigma$ and $Y = \{0\} \times \Sigma$ for some manifold $\Sigma$, then one may define
\[[X,Y] = [\R^n,0]\times \Sigma\] i.e.\ by blowing up with the parameter $\Sigma$, and the blowdown map is trivial in the direction of $\Sigma$. This construction is in fact coordinate invariant (this is most easily checked in polar coordinates), so allows one to define $[X,Y]$ for an arbitrary manifold $X$ and $Y$ a $p$-submanifold, and the associated blowdown map $\beta:[X,Y] \to X$, as well as a front face $\mathbf{ff}$ which $\beta$ maps to $Y$.

The purpose of blowups is to ``separate'' the behaviour of a function $f$ on $X$ as it comes in radially to the submanifold $Y$. We illustrate this with a simple example. The function $f(x,y) = \sqrt{x^2+y^2}$ on $\R^2$ is not smooth. However, its lift to $[\R^2,0]$ is just $r$, which is smooth.

A technique we will use in conjunction with blowing up is changing the smooth structure. Given a mwc $X$, one may change the smooth structure by requiring an $n$th root of a bdf to be smooth, for $n \in \N$. In this way, one obtains a new mwc $\widetilde{X}$. While the resulting manifold is diffeomorphic to the original, one regards them as \emph{distinct} since the diffeomorphism is not canonical. However, there is always a smooth map $f:\widetilde{X} \to X$ which raises the new bdf to the power of $n$.

For an example, consider $X = [0,\infty)$ and the boundary face $\{0\}$. Then one may obtain a new smooth manifold $\tilde{X}$ by requiring $x^{1/3}$ to be smooth. Then $\widetilde{X} \iso [0,\infty)$, except the diffeomorphism is not canonical and one instead has a smooth map $\widetilde{X} \to X$ sending $x \mapsto x^3$.

\section{Polyhomogeneity}
\label{C:A4:phg}
We make the definition:
\begin{defn}An index set is a discrete subset $E \subseteq \C \times \N$ such that: 
\begin{romanumerate}
\item $\{(z,p) \in E\: \Re(z) < N\}$ is finite;
\item $(z,p) \in E \Rightarrow (z+k,p) \in E$ for all $k \in \N$;
\item $(z,p) \in E \Rightarrow (z,q) \in E$ for all $q \leq p$.
\end{romanumerate}
\end{defn} 
We make the convention that an integer $n$ stands for the index set $\{n,n+1,\ldots\}\times \{0\}$.

Let us first define polyhomogeneity of a function on $X = [0,1)_t\times \R^k_x$. We say, a little informally, that a function $u:[0,1)\times \R^k \to \C$ is polyhomogeneous with index set $E$ if $u$ admits the expansion
\[u(t) \sim \sum_{(z,p) \in E} u_{z,p}t^z\log^p(t)\] where $u_{z,p} \in C^{\infty}(\R^k)$.
Define $\mathcal A_{\mathrm{phg}}^E(X)$ to be the set of all functions possessing such expansions. In order to make this precise, we need to make sense of the remainder term.

\begin{defn}The \emph{conormal} space $\mathcal A(X)$ is the space of all functions $u$ such that $(t\pa_t)^m\pa_x^n u \in L^{\infty}_{\mathrm{loc}}(X)$ for all integers $m$ and multi-indices $n$.\end{defn}

With this definition, we can define:
\begin{defn} We say $u \in \mathcal A_{\mathrm{phg}}^E(X)$ if for all $(z,p) \in E$ there are $u_{z,p} \in C^{\infty}(\R^k)$ such that $u$ admits the expansion for any $N$
\begin{equation}\label{eq:A4:phgi}u(t) = \sum_{\substack{(z,p) \in E\\ \Re(z) < N}} u_{z,p}t^z \log^p(t) + t^N \mathcal A(X).\end{equation}
\end{defn}

Also observe that for the index set $E = 0$, $A_{\mathrm{phg}}^E(X) = C^{\infty}(X)$, and the index set $E = \emptyset$ corresponds to functions rapidly vanishing at $\{t = 0\}$. Observe that part (ii) of the definition of an index set implies that $C^{\infty}(X) \cdot \mathcal A_{\mathrm{phg}}^E(X) \subseteq \mathcal A_{\mathrm{phg}}^E(X)$.

Now let $X = [0,1)^j_t \times \R^k_x$, and let $\mathcal E = (E_1,\ldots, E_n)$ be a collection of index sets. We call $\mathcal E$ an index family. We seek to define the space $\mathcal A_{\mathrm{phg}}^{\mathcal E}(X)$. Let $\widehat{\mathcal{E}}_i$ denote the index family obtained from $\mathcal E$ after removing $E_i$.

Heuristically, we wish to define $\mathcal A_{\mathrm{phg}}^{\mathcal E}(X)$ inductively by saying that $u \in \mathcal A_{\mathrm{phg}}^{\mathcal E}(X)$ if for each $1 \leq i \leq j$, $u$ admits an asymptotic expansion
\[u \sim \sum_{(z,p) \in E_i} u_{z,p}t^z\log^p(t),\] where each $u_{z,p} \in \mathcal A_{\mathrm{phg}}^{\widehat{\mathcal {E}}_i}(\{t_i = 0\})$. However, addressing the remainder is a slightly tricky. To do this, we will in fact inductively define the ``partially polyhomogeneous spaces.'' 

Let $\mathcal F \subseteq \{1,\ldots,j\}$ and let $\mathcal E$ be an index family, which we think of as a map $\mathcal F$ into index sets (and write $E_i$ for $\mathcal E(i)$). We will define by induction on $n = \#\mathcal F$ the space $\mathcal A_{\mathrm{phg},\mathcal F}^{\mathcal E'}(X)$. Intuitively, this is the space of functions which have polyhomogeneous expansions at the faces in $\mathcal F'$, but are only conormal at all other faces. Of course, in the end, we define
$\mathcal A_{\mathrm{phg}}^{\mathcal E}(X) = \mathcal A_{\mathrm{phg},\{1,\ldots,j\}}^{\mathcal E}(X)$.

Let us first define the conormal space $\mathcal A(X)$ to be the collection of those functions $u$ with $(t\pa_t)^{m}\pa_x^n u \in L^{\infty}_{\mathrm{loc}}(X)$ for all multi-indices $m,n$.

\begin{defn}If $n=0$, then $\mathcal A_{\mathrm{phg},\emptyset}^{\emptyset}(X) = \mathcal A(X)$ (here $\emptyset$ is the empty index family, not the empty index set).

If $n \geq 1$, then we say $u \in \mathcal A_{\mathrm{phg},\mathcal F}^{\mathcal E}(X)$ for all $i \in \mathcal F$, the following is true. Denote by $\widehat{\mathcal F}_i = \mathcal F\setminus \{i\}$, and $\widehat{\mathcal E}_i$ the index family of all index sets in $\mathcal E$ other than $E_i$. Then $u$ is required to have an expansion for all $N$
\[u = \sum_{\substack{(z,p) \in E_i\\ \Re(z) < N}} u_{z,p}t_i^z\log^p(t_i) + t_i^N\mathcal A_{\mathrm{phg},\widehat{\mathcal F}_i}^{\widehat{\mathcal E}_i}(X),\] where each $u_{z,p} \in \mathcal A_{\mathrm{phg},\hat{\mathcal F}_i}^{\widehat{\mathcal E}_i}(\{t_i = 0\})$.
\end{defn}
In other words, at any face $u$ has a polyhomogeneous expansion where all coefficients are polyhomogeneous, and a remainder which is decaying at that face, but polyhomogeneous at the others.

We make a few observations. First, the index family $\mathcal E = (0,\ldots,0)$ corresponds precisely to smooth functions, and if $E_i = \emptyset$ for some $i$, then functions in $\phgi{J}{\mathcal E}$ are rapidly vanishing at $\{t_i = 0\}$. We will often call say that such functions are ``$0$ in (Taylor) series at $\{t_i = 0\}$,'' since in their expansion as a polyhomogeneous series at $\{t_i = 0\}$ all coefficient are identically $0$.

Second, part (ii) of the definition of an index set implies that $C^{\infty}(X) \cdot \mathcal A_{\mathrm{phg,\mathcal F}}^\mathcal E(X) \subseteq \mathcal A_{\mathrm{phg}, \mathcal F}^\mathcal E(X)$. Third, part (iii) of the definition shows that the class of partially polyhomogeneous functions are stable under differentiation via $t_i\pa_{t_i}$ and $\pa_{x_i}$ (of course if $E_i = 0$, then it is also stable under $\pa_{t_i}$).

Parts (ii) and (iii) of the definition of an index set also imply that $\mathcal A_{\mathrm{phg,\mathcal F}}^{\mathcal E}(X)$ is coordinate invariant.\footnote{If $\phi: X \to X$ is a diffeomorphism, it necessarily permutes the boundary faces of $X$. If $\phi$ maps $\{t_{i'} = 0\}$ into $\{t_i = 0\}$, it follows that $\phi_{i}(t,x) = t_{i'}f_i(t,x)$, where $f$ is smooth and nonvanishing on $\{t_i = 0\}$. To  show coordinate invariance, it is necessary to prove a more general statement. Suppose $j' \geq j$, and let $X' = [0,1)^{j'}_s\times \R^{k'}_y$, and $\phi:X' \to X$ be a smooth map for which there is an injection $\sigma: \{1,\ldots, j\} \to \{1, \ldots, j'\}$ for which $\phi_{i}(s,y) = s_{\sigma(i)}f_{\sigma(i)}(s,y)$ and $f_{\sigma(i)}(s,y)$ is smooth and nonvanishing at $\{s_{\sigma(i)} = 0\}$. Let $\mathcal F \subseteq \{1, \ldots, j\}$, and $\mathcal E$ be an associated index family. Let $\mathcal F' = \sigma(\mathcal F) \un \{1,\ldots,j'\}\setminus\sigma(\{1,\ldots,j\}) \subseteq \{1,\ldots, j'\}$ and let $\mathcal E'$ be the associated index family defined by $E'_{\sigma(i)} = E_i$ and $E'_{i'} = 0$ if $i' \not \in \sigma(\{1,\ldots,j\})$. In other words, $\mathcal F'$ corresponds to the union those faces which $\phi$ maps to a face in $\mathcal F$ and the faces which do not get mapped into a face of $X$, and $\mathcal E'$ is smooth at all faces which do not get mapped into a face of $X$. Then if $u \in \mathcal A_{\mathrm{phg,\mathcal F}}^{\mathcal E}(X)$, $\phi^\ast u \in A_{\mathrm{phg,\mathcal F'}}^{\mathcal E'}(X')$. This is proved by induction on $j$. That the index set $E'_{i'} = 0$ if $i' \not \in \sigma(\{1,\ldots,j\})$ comes from the fact that \[\pa_{s_{i'}} \phi^\ast u = \sum_{i=1}^j \phi^\ast (t_i\pa_{t_i} u)\pa_{s_{i'}}\log f_{\sigma(i)} + \sum_{i=1}^k \phi^\ast(\pa_{x_i} u)\pa_{s_{i'}}\phi_{j+i},\] and polyhomogeneity is preserved by taking the $b$-derivatives $t_i\pa_{t_i}$ and $\pa_{x_i}$, and then iterating with higher-derivatives. The diffeomorphism case is the case $j = j'$ and $\sigma$ is a permutation. 

Even more generally, this is a special case of Melrose's famous pullback theorem.}

Thus, if $X$ is a mwc with faces $F_1,\ldots,F_j$, and  $\mathcal F$ is a subcollection of the faces, with $\mathcal E$ an associated index family, we may define the partially polyhomogeneous spaces (and the polyhomogeneous spaces) $\mathcal A_{\mathrm{phg},\mathcal F}^{\mathcal E}(X)$ via a partition of unity.

Let us return to the case $X = [0,1)^j_t \times \R^k_x$. Suppose $u \in \mathcal A_{\mathrm{phg}}^{\mathcal E}(X)$. Fix $i$, and let $u^i_{z,p}$, $(z,p) \in E_i$ denote the coefficients in the polyhomogeneous expansion at $\{t_i = 0\}$. By definition of polyhomogeneity, for $i' \neq i$, $u^i_{z,p}$, is itself polyhomogeneous, and so has an expansion at $\{t_{i'} = 0\}$. Let $u^{i,i'}_{z,p,z',p'}$, $(z',p') \in E_{i'}$, denote the coefficients of this polyhomogeneous expansion, which are polyhomogeneous on $\{t_i = t_{i'} = 0\}$. One can also reverse this process, obtaining other polyhomogeneous functions $u^{i',i}_{z',p',z,p}$ on $\{t_i = t_{i'} = 0\}$.  As one expects, $u^{i,i'}_{z,p,z',p'} = u^{i',i}_{z',p',z,p}$ for any $(z,p) \in E_i$, $(z',p') \in E_{i'}$. In fact more is true, and there is a relationship between the expansions of the remainders in these expansions, too. To state this, let us make some notational definitions. Let $\widehat{\mathcal E}_i$ denote the index family obtained from $\mathcal E$ by removing $E_i$, and for $i' \neq i$ let $\widehat{\mathcal E}_{i,i'}$ denote the index family obtained by removing both $E_i$ and $E_{i'}$. For any $i$ and $N$ we may expand
\[u = \sum_{\substack{(z,p) \in E_i\\ \Re(z) < N}} u^i_{z,p} t_i^z \log^p t_i + t_i^N R^i_N,\] where $R^i_N \in \mathcal A_{\mathrm{phg}}^{\widehat{\mathcal E}_i}(X)$, and for $i' \neq i$ and any $N'$ we may expand
\[u_{z,p}^i = \sum_{\substack{(z',p') \in E_{i'}\\ \Re(z') < N'}} u^{i,i'}_{z,p,z',p'} t_{i'}^{z'} \log^{p'} t_{i'} + t_{i'}^N S^{i,i'}_{z,p,N'},\]
where $R^{i,i
}_{N'} \in \mathcal A_{\mathrm{phg}}^{\widehat{\mathcal E}_{i,i'}}(\{t_i = 0\})$. Also expand for any $N'$
\[R^i_{N} = \sum_{(z',p') \in E_{i'},\Re(z') < N'} R^{i,i'}_{N,z',p'}t_{i'}^{z'}\log^{p'}t_{i'} + t_{i'}^{N'}T^{i,i'}_{N'}\] where
$T^{i,i'}_{N,N'} \in \mathcal A_{\mathrm{phg}}^{\widehat{\mathcal E}_{i,i'}}(X)$. Then:
\begin{lem}\label{thm:A4:compat}For any $i \neq i$, $(z,p) \in E_i$ and $(z',p') \in E_{i'}$ it holds that $u^{i,i'}_{z,p,z',p'} = u^{i',i}_{z',p',z,p}$. For any $N,N' \in \R$ if $\Re(z) < N$, $\Re(z') < N'$ it holds that $S^{i,i'}_{z,p,N'} = R^{i',i}_{N',z,p}$, and $T^{i,i'}_{N,N'} = T^{i',i}_{N',N}$.\end{lem}
We have an important corollary. Let $u$ and $\mathcal E$ be as above. Fix some $i$ and let $E'_i \subseteq E_i$ denote a sub-index set. Let $\mathcal E'$ denote the index family with $E_i$ replaced by $E'_i$. Then:
\begin{cor}\label{thm:A4:betterphg}Suppose $u^i_{z,p} \equiv 0$ for $(z,p) \in E_i\setminus E'_i$. Then $u \in \phg{\mathcal E'}(X)$.\end{cor}
A particular important case is the case $E'_i = \emptyset$.
Of course \cref{thm:A4:compat} and \cref{thm:A4:betterphg} remain true if the polyhomogeneity is replaced by partial polyhomogeneity, and if $X$ is replaced by an arbitrary manifold with corners.
\Cref{thm:A4:betterphg} may seem obvious. However, in the definition of polyhomogeneity one also needs to check that the expansions at faces other than $\{t_i = 0\}$ have their expansions only corresponding to terms in $E'_i$, as well, and hence the need for \cref{thm:A4:compat}. Indeed:
\begin{proof}[Proof of \cref{thm:A4:betterphg}.]For $i' \neq i$, we need to show for all $(z',p') \in E_{i'}$ that $u^{i'}_{z',p'}$, which is a priori in $\phgd{\{t_{i'} = 0\}}$ is actually in $\phg{\mathcal E'}(\{t_{i'} = 0\})$, and for all $N' \in \R$ that $R_{N'}^{i'}$, which is a piori in $\phgd(X)$ is actually in $\phg{\mathcal E'}(X)$. This entails showing that $u^{i',i}_{z',p',z,p} = 0$ and $R^{i',i}_{N',z,p} = 0$ if $(z,p) \not \in E'_i$. But this is true by \cref{thm:A4:compat} since $u^i_{z,p} = 0$.\end{proof}
We will prove \cref{thm:A4:compat} at the end of this section.

We end this primer by reviewing Borel's lemma, which allows one to asymptotically sum a collection of coefficients $u_{z,p}$, $(z,p) \in E$, for some index set $E$, obtaining a polyhomogeneous function. The first version we state is:
\begin{lem}[Borel's lemma I]Let $X = [0,1)\times \R^k$, and let $E$ be an index set. Suppose $u_{z,p} \in C^{\infty}(\R^k)$ are given for each $(z,p) \in E$. Then there exists $u \in \mathcal A_{\mathrm{phg}}^E(X)$ with the expansion
\[u \sim \sum_{(z,p) \in E} u_{z,p}t^z \log^p(t).\]\end{lem}
The proof is routine (although one should use a partition of unity to first reduce to the cast that all $u_{z,p}$ are compactly supported in the same compact set).

One also has its higher-dimensional version. Let $X = [0,1)_t^j \times \R^k$. Let $\mathcal E = (E_1,\ldots, E_j)$ be an index family and for $1 \leq i \leq j$ let $\widehat{\mathcal E}_i$ be the index family obtained by removing the index set $E_i$. Suppose for $1 \leq i \leq j$ and $(z,p) \in E_i$ that we are given functions $u^i_{z,p} \in \mathcal A_{\mathrm{phg}}^{\widehat{\mathcal E}_i}(\{t_i = 0\})$. One wishes to find $u \in \mathcal A_{\mathrm{phg}}^{\mathcal E}(X)$ such that for each $i$
\begin{equation}\label{eq:A4:poly}u \sim \sum_{(z,p) \in E_i} u^i_{z,p}t_i^z\log^p(t_i).\end{equation} However, this is patently not possible in general, since \cref{thm:A4:compat} imposes for $i \neq i$ a condition on the coefficients in the expansions of $u^i_{z,p}$ at $\{t_{i'} = 0\}$ and the expansions of $u^{i'}_{z',p'}$ at $\{t_i = 0\}$. Indeed, denoting by $u^{i,i'}_{z,p,z',p'}$ and $u^{i',i}_{z',p',z,'p}$ the coefficients in the expansions, respectively, it must be true that $u^{i,i'}_{z,p,z',p'} = u^{i',i}_{z',p',z,'p}$ if \eqref{eq:A4:poly} is to hold. We call this the the \emph{compatibility condition} on the coefficients $u^i_{z,p}$.
However, the compatibility condition is the only obstruction to summing coefficients $u^i_{z,p}$ asymptotically. Thus we have
\begin{lem}[Borel's lemma II]\label{thm:A4:Borel}Suppose we are given $u^i_{z,p}$ as above, and the compatibility condition is satisfied. Then there exists $u \in \phgd(X)$ satisfying \eqref{eq:A4:poly} for all $i$.\end{lem}

\begin{rk}If $E_i = 0$ for all $i$, then the compatibility condition is equivalent to $\pa_{t_i'}^m u^i_\ell = \pa_{t_i}^\ell u^{i'}_m$ for all $i,i'$ and $m,\ell$.\end{rk}

\begin{rk}Of course, the asymptotic sum is not unique. However, it is however unique modulo an element of $\phg{(\emptyset,\ldots,\emptyset)}(X)$. Indeed, if $u_1$ and $u_2$ are two asymptotic sums, then by definition for all $i$ the coefficients in the polyhomogeneous expansion of $v = u_1-u_2$ at $\{t_i = 0\}$, $v^i_{z,p}$, are identically $0$. Thus by \cref{thm:A4:betterphg} $v \in \phg{(\emptyset,\ldots,\emptyset)}(X)$.\end{rk}

We carry out the proof at the end of this section.

Of course, these lemmas carry through to an arbitrary mwc via partitions of unity. Lastly, we talk about vector bundles. If $X$ is a mwc and $V$ is a vector bundle over $X$, then we may also analogously define $\mathcal A^{\mathcal E}_{\mathrm{phg}}(X;V)$, or more simply by tensoring
\[\mathcal A^{\mathcal E}_{\mathrm{phg}}(X;V) = C^{\infty}(X;V)\otimes\mathcal A^{\mathcal E}_{\mathrm{phg}}(X).\] All lemmas continue to hold in this setting.

\subsection{Proofs of lemmas}
We first prove \cref{thm:A4:compat}. \begin{proof}[Proof of \cref{thm:A4:compat}.]
It suffices to fix $N,N'$ and show the lemma if $\Re(z) < N$, $\Re(z') < N'$. The notation is cumbersome, even though the idea is clear.
We may expand
\begin{align*}
u &= \sum_{\substack{(z,p) \in E_i \\ \Re(z) < N}}t_i^{z}\log^p t_i \left(\sum_{\substack{(z',p') \in E_{i'} \\ \Re(z') < N'}} u^{i,i'}_{z,p,z',p'}t_{i'}^{z'}\log^{p'} t_{i'} + t_{i'}^{N'}S^{i,i'}_{z,p,N'}\right)\\
&+ t_i^N\left(\sum_{\substack{(z',p') \in E_{i'} \\ \Re(z') < N'}} R^{i,i'}_{N,z',p'}t_{i'}^{z'}\log^{p'} t_{i'} + t_{i'}^{N'}T^{i,i'}_{N,N'}\right)\\
&= \sum_{\substack{(z',p') \in E_{i'} \\ \Re(z') < N'}}t_i^{z}\log^p t_i \left(\sum_{\substack{(z,p) \in E_i \\ \Re(z) < N}} u^{i',i}_{z',p',z,p}t_{i}^{z}\log^{p} t_{i'} + t_{i}^{N}S^{i',i}_{z',p',N}\right)\\
&+ t_{i'}^{N'}\left(\sum_{\substack{(z,p) \in E_{i}\\ \Re(z) < N}} R^{i',i}_{N',z,p}t_{i}^{z}\log^{p} t_{i} + t_{i}^{N}T^{i',i}_{N',N}\right).
\end{align*}
Subtracting these, one obtains
\begin{align}
\label{eq:A4:eeee}
\begin{split}
0 &= \sum_{\substack{(z,p) \in E_i \\ \Re(z) < N}}t_i^{z}\log^p t_i \left(\sum_{\substack{(z',p') \in E_{i'} \\ \Re(z') < N'}} (u^{i,i'}_{z,p,z',p'}-u^{i',i}_{z',p',z,p})t_{i'}^{z'}\log^{p'} t_{i'} + t_{i'}^{N'}(S^{i,i'}_{z,p,N'}-R^{i',i}_{N',z,p})\right)\\
&+ t_i^N\left(\sum_{\substack{(z',p') \in E_{i'} \\ \Re(z') < N'}} (R^{i,i'}_{N,z',p'}-S^{i',i}_{z',p',N})t_{i'}^{z'}\log^{p'} t_{i'} + t_{i'}^{N'}(T^{i,i'}_{N,N'}-T^{i',i}_{N',N})\right)
\end{split}
\end{align}
Denote by $q^i_z = \max \{ p \: (z,p) \in E_i\}$. Fix some $z$ with $\Re(z) < N$, and consider for $0 \leq k \leq q^i_z$ the differential operator
\[P^i_{z,k} = (t_i\pa_{t_i}-z)^{k}\prod_{\substack{w \: \exists (w,p) \in E_i\\w \neq z, \Re(w) < N}} (t_i\pa_{t_i}-z)^{q^i_w+1}.\]
Notice that: $P_{z,k} t_i^w\log^p t_i = 0$ if $(w,p) \in E_i$, $w \neq z$, $\Re(z) < N$,  $P_{z,k} t_i^z\log^p t_i = 0$ if $p < k$, and $P_{z,k} t_i^z \log^p t_i = \frac{p!}{(p-k)!}$ if $p \geq k$. Define $q^{i'}_{z'}$ and $P^{i'}_{z',k'}$ analogously. The key idea is that we can use $P^i_{z,k}$ and $P^{i'}_{z',k'}$ to isolate the correct coefficients in \eqref{eq:A4:eeee} and show that that are $0$.
Also notice that $S_{z,p,N'}^{i,i'}$, $R^{i',i}_{N',z,p}$, $S_{z',p',N}^{i',i}$, $R^{i,i'}_{N,z',p'}$, $T_{N,N'}^{i,i'}$ and $T_{N',N}^{i',i}$ are in $\bm{w} A(X)$, where $\bm{w}$ is a weight $\bm{w} = \prod_{j \neq i,i'}t_j^{\alpha_j}$, for some $\alpha_j \in \R$ (the $\alpha_j$ just need to be large enough to offset whatever growth can come from the polyhomogeneity at other faces $\{t_j = 0\}$).
Fixing some $z$ and applying $P_{z,q^i_z}$ to \eqref{eq:A4:eeee} thus yields
\[0 = t_i^{z}\log^{q^i_z} t_i\left(\sum_{\substack{(z',p') \in E_{i'} \\ \Re(z') < N'}}\hspace{-4pt} (u^{i,i'}_{z,p,z',p'}-u^{i',i}_{z',p',z,p})t_{i'}^{z'}\log^{p'} t_{i'}\hspace{-0.5pt} + \hspace{-0.5pt} t_{i'}^{N'}(S^{i,i'}_{z,p,N'}-R^{i',i}_{N',z,p})\right) + t_i^N\bm{w}\hspace{-0.5pt}\mathcal A(X).\]
Dividing by $t_i^{z}\log^{q^i_z}$ and taking $t_i \to 0$ then shows that
\begin{equation}\label{eq:A4:ffff}0 = \left(\sum_{\substack{(z',p') \in E_{i'} \\ \Re(z') < N'}} (u^{i,i'}_{z,q^i_z,z',p'}-u^{i',i}_{z',p',z,q^i_z})t_{i'}^{z'}\log^{p'} t_{i'} + t_{i'}^{N'}(S^{i,i'}_{z,q^i_z,N'}-R^{i',i}_{N',z,q^i_z})\right).\end{equation}
Fixing some $z'$ and applying $P^{i'}_{z',q^{i'}_{z'}}$ then shows that
\[0 =u^{i,i'}_{z,q^i_z,z',q^{i'}_{z'}}-u^{i',i}_{z',q^{i'}_{z'},z,q^i_z} + t_{i'}^{N'}\bm{w}\mathcal A(X).\] Taking $t_{i'} \to 0$ then shows that
\[0 =u^{i,i'}_{z,q^i_z,z',q^{i'}_{z'}}-u^{i',i}_{z',q^{i'}_{z'},z,q^i_z}.\]
Now we may apply inductively $P^{i'}_{z',p'}$ for each $q^{i'}_{z'} > p' \geq 0$ to \eqref{eq:A4:ffff} to obtain as above
\[0 =u^{i,i'}_{z,q^i_z,z',p'}-u^{i',i}_{z',p',z,q^i_z}\] for all $p'$. Since $z'$ was arbitrary, this is true for all $z'$, too. With this true, it follows from \eqref{eq:A4:ffff} that $S^{i,i'}_{z,q^i_z,N'}-R^{i',i}_{N',z,q^i_z}$. This is part of what we need to show, except we only know now that it is true for when $p = q^{i_z}$. However, we may now apply $P^{i}_{z,k}$ for $q^i_z > p \geq 0$ inductively and use the same argument to obtain $u^{i,i'}_{z,p,z',p'}=u^{i',i}_{z',p',z,p}$ and $S^{i,i'}_{z,p,N'}=R^{i',i}_{N',z,p'}$, which is what we want. Thus \eqref{eq:A4:eeee} shows that
\[0 = \sum_{\substack{(z',p') \in E_{i'} \\ \Re(z') < N'}} (R^{i,i'}_{N,z',p'}-S^{i',i}_{z',p',N})t_{i'}^{z'}\log^{p'} t_{i'} + t_{i'}^{N'}(T^{i,i'}_{N,N'}-T^{i',i}_{N',N})\] applying the same argument as we did to \eqref{eq:A4:ffff}, we see that $R^{i,i'}_{N,z',p'} = S^{i',i}_{z',p',N}$ and $T_{N,N'}^{i,i'} = T_{N',N}^{i',i}$ for all $\Re(z) > N, \Re(z') > N$.\end{proof}

We may now prove the higher-dimensional Borel's lemma. \begin{proof}[Proof of \cref{thm:A4:Borel}]Using a partition of unity, we may suppose that all $u^i_{z,p}$ are compactly supported in the same compact set. For each $n \leq j$, we will prove the theorem is true for the special index family
\[\mathcal E = (E_1,\ldots,E_n,\emptyset,\ldots,\emptyset),\]
i.e.\ if everything is rapidly vanishing at the last $j-n$ faces. The case $n=j$ is Borel's lemma. We will accomplish this by induction on $n$.

The case $n=0$ is vacuous, since $u \equiv 0$ works.

For the $n=1$ case, set
\[v = \sum_{(z,k) \in E_1} t_1^z\log^k(t_1) u^1 \phi(t_1/\epsilon_{z,k}),\] for $\phi$ a smooth cutoff of $\{t_1 = 0\}$ and $\epsilon_{z,k}$ small. \Cref{thm:A4:Borelhelper}, below, implies $v \in \phg{(E_1,\emptyset,\ldots,\emptyset)}$, and has the correct series expansion, which completes the $n=1$ case.

Now assume the claim is true for $n-1$. We prove it for $n$. Set
\[v = \sum_{(z,p) \in E_1} u^1_{z,p}t_1^z \log^p(t_1) \phi(t_1/\epsilon_{z,p}),\] for $\phi$ a smooth cutoff of $\{t_1 = 0\}$ and $\epsilon_{z,k}$ small. \Cref{thm:A4:Borelhelper}, below, implies $v \in \phgd(X)$, and $v$ has the correct expansion at $\{t_1 = 0\}$. It has the correct expansions at $\{t_i = 0\}$ for $i > n$ (since by definition it is polyhomogeneous with empty index set at those faces).
However, $v$ may not have the correct expansion at $\{t_i = 0\}$ for $2 \leq i \leq n$. Write the putative asymptotic sum as $u = v+w$, where $w$ is another polyhomogeneous function. Then $u$ will satisfy \eqref{eq:A4:poly} if for all $1 \leq i \leq j$ and $(z,p) \in E_i$, $w^i_{z,p} = u^i_{z,p} - v^i_{z,p}$. For $i = 1$ and $i > n$, this difference is $0$. For $2 \leq i \leq n$, notice that the compatibility conditions on the coefficients in the expansion of $v$ are satisfied by assumption, since $v$ is polyhomogeneous, and thus the compatibility conditions are satisfied for the coefficients of the putative function $w$. We may now use the inductive hypothesis to find the polyhomogeneous $w$, since $w$ is to be rapidly vanishing at both $\{t_i = 0\}$ for $i > n$ and $i = 1$.\end{proof}

\begin{lem}\label{thm:A4:Borelhelper}Let $X = [0,1)^j\times \R^k$, and let $\emptyset \neq J \subseteq \{1,\ldots,j\}$, and suppose without loss of generality that $1 \in J$. Let $I = J\setminus \{1\}$. Let $\mathcal E$ be the index family corresponding to the faces in $J$, and let $\mathcal E'$ be the index family corresponding to $I$ obtained by deleting $E_1$. Suppose we are given for $(z,p) \in E_1$ functions $v_{z,p} \in \phgi{I}{\mathcal E'}(\{t_1 = 0\})$ which are supported in some compact set $K$.

Set
\begin{equation}\label{eq:A4:idkanymore}v = \sum_{(z,p) \in E_1} t_1^z \log^p (t_1) v_{z,p} \phi(t_1/\epsilon_{z,p}),\end{equation} where $\phi$ is smooth cutoff of $\{t_1 = 0\}$ and $\epsilon_{z,p}$ are small enough. Then $v \in \phgi{J}{\mathcal E}(X)$ and
\[v \sim \sum_{(z,p) \in E_1} t_1^z \log^p (t_1) v_{z,p}\] at $\{t_1 = 0\}$. Moreover, the convergence of \eqref{eq:A4:idkanymore} is in the topology of
\[\prod_{i \in J} t_i^{\rho_i}\mathcal A(X) =: t^{\rho}\mathcal A(X),\]
where $\rho_i$ are any real numbers such that
\[\rho_i <  \inf \{\Re(z) \: (z,p) \in E_i\}.\] \end{lem}
\begin{proof}This is an extension of the usual proof of Borel's lemma. By the standard techniques in proving (ordinary) Borel's lemma, it is clear that the convergence is in $t^{\rho}\mathcal A(X)$, so the hard part is proving the convergence to a (partially) polyhomogenous function rather than just a (weighted) conormal function.

We prove this by induction on $n = \#J$.
Let us treat the base case $n = 1$. Near $t_1 = 0$, we may write
\begin{align*}
v &= \sum_{\substack{(z,p) \in E_1\\ \Re(z) < N}} t_1^z \log^p (t_1) v_{z,p}\\
&+ t_1^N\sum_{\substack{(z,p) \in E_1\\ \Re(z) \geq N} } t_1^{z-N}\log^p(t_1) v_{z,p}\phi(t_1/\epsilon_{z,p})\\
&+ \sum_{\substack{(z,p) \in E_1\\ \Re(z) < N}} t_1^z \log^p (t_1)v_{z,p}(1-\phi(t_1/\epsilon_{z,p})).\end{align*}
If $t_1 \gtrsim \sup_{\Re(z) < N} \epsilon_{z,p}$, then the last sum vanishes. The standard argument shows that if $\epsilon_{z,p}$ are small enough (and not depending on $N$) then second sum converges in the topology of $t_1^N \mathcal A(X)$. This completes the proof of the case $j=1$.
Now assume we have proven it for $n-1$, we prove it for $n$.

Near $t_1 = 0$, we may again write
\begin{align*}v &= \sum_{\substack{(z,p) \in E_1\\ \Re(z) < N}} t_1^z \log^p (t_1) v_{z,p}\\
&+ t_1^N\sum_{\substack{(z,p) \in E_1\\ \Re(z) \geq N} } t_1^{z-N}\log^p(t_1) v_{z,p}\phi(t_1/\epsilon_{z,p})\\
&+ \sum_{\substack{(z,p) \in E_1\\ \Re(z) < N}} t_1^z \log^p (t_1)v_{z,p}(1-\phi(t_1/\epsilon_{z,p}),\end{align*}
where again the last term vanishes if $t_1$ is small enough, so we may ignore it. The first term is fine since the sum is finite. Thus, we only need to show that 
\begin{equation}\label{eq:A4:helper}\sum_{\substack{(z,p) \in E_1\\ \Re(z) \geq N} } t_1^{z-N}\log^p(t_1) v_{z,p}\phi(t_1/\epsilon_{z,p}) \in \phgi{I}{\mathcal E'}(X),\end{equation}
Let $k \in J$, $k \neq 1$ and set $K = J\setminus \{1,k\}$, and let $\mathcal E''$ denote the index family obtained from $\mathcal E'$ by removing $E_k$. Also let $\tilde{K} = J\setminus \{k\}$ and $\tilde{\mathcal E''}$ be the index family obtained from $\mathcal E$ by removing $E_k$ (equivalently obtained from $\mathcal E''$ by adding $E_1$).
Then we may expand
\[v_{z,p} = \sum_{\substack{(w,q) \in E_k\\ \Re(w) < M}} v_{z,p,w,q}t_k^w + t_k^M R_{z,p,M},\]
where $v_{z,p,w,q} \in \phgi{K}{\mathcal E''}(\{t_1 = t_k = 0\})$ and $R_{z,p,M} \in \phgi{K}{\mathcal E''}(\{t_1 = 0\})$.
By induction (and shrinking the $\epsilon_{z,p}$ if necessary)
\[v^k_{w,q} := \sum_{(z,p) \in E_1 }v_{z,p,w,q}t_1^{z}\log^p(t_1)\phi(t_1/\epsilon_{z,p}) \in \phgi{\tilde{K}}{\tilde{\mathcal E''}}(\{t_k = 0\}).\]
In particular, if we denote by $\widetilde{\mathcal {E''}^N}$ the index set obtained from $\widetilde{\mathcal {E''}}$ by replacing $E_1$ with $\{(z,p) \in E_1 \: \Re(z) \geq N\}$, it follows from \cref{thm:A4:betterphg} that
\[\sum_{\substack{(z,p) \in E_1\\ \Re(z) \geq N} }v_{z,p,w,q}t_1^{z}\log^p(t_1)\phi(t_1/\epsilon_{z,p}) \in \phgi{\tilde{K}}{\tilde{\mathcal {E''}^N}}(\{t_k = 0\}),\] and so
\[v^k_{w,q,N} := \sum_{\substack{(z,p) \in E_1\\ \Re(z) \geq N}}v_{z,p,w,q}t_1^{z-N}\log^p(t_1)\phi(t_1/\epsilon_{z,p}) \in t^{-N}\phgi{\tilde{K}}{\tilde{\mathcal {E''}^N}}(\{t_k = 0\}) \subseteq \phgi{K}{\mathcal E''}(\{t_k = 0\})\]
(we can ``forget'' about the polyhomogeneity at $\{t_1 = 0\}$ and replace it with conormality since by construction there are no $(z,p)$ with $\Re(z) < 0$ in the expansion of $v^k_{w,q,N}$ at $\{t_1 = 0\}$).
Similarly, by induction (shrinking the $\epsilon_{z,p}$ if necessary\footnote{It appears that we may need to shrink $\epsilon_{z,p}$ for each $M \in \R$, and so this argument does not close. However, we can always get away with only shrinking those $\epsilon_{z,p}$ with $\Re(z) \geq M$, since changing the value of $\epsilon_{z,p}$ in finitely many terms does not affect the conclusion. Thus, each $\espilon_{z,p}$ is only shrunk a finite number of times, and the argument closes.})
\[R^{k}_M := \sum_{(z,p) \in E_1 }R_{z,p,M}t_1^{z}\log^p(t_1)\phi(t_1/\epsilon_{z,p}) \in \phgi{\tilde{K}}{\tilde{\mathcal E''}}(X)\] and so 
\[R^k_{M,N} := \sum_{\substack{(z,p) \in E_1\\ \Re(z) \geq N}}R_{z,p,M}t_1^{z-N}\log^p(t_1)\phi(t_1/\epsilon_{z,p}) \in \phgi{K}{\mathcal E''}(X).\] We conclude that
\begin{align*}\sum_{\substack{(z,p) \in E_1\\ \Re(z) \geq N} } t_1^{z-N}&\log^p(t_1) v_{z,p}\phi(t_1/\epsilon_{z,p})\\
&= \sum_{\substack{(w,q) \in E_k\\ \Re(w) < M}} t_k^w\log^q(t_k)v^k_{w,q,N} + t_k^MR^k_{M,N}\\
&\in \sum_{\substack{(w,q) \in E_k\\ \Re(w) < M}} t_k^w\log^q(t_k)\phgi{K}{\mathcal E''}(\{t_k = 0\}) + t_k^M\phgi{K}{\mathcal E''}(X).\end{align*} Doing this for all $k$ shows \eqref{eq:A4:helper}.

Now we need to show that $v$ is polyhomogeneous at the other faces. Fix $k \in J$, $k \neq 1$, and let $K,\ \mathcal E''$ as above. Expand
\[v_{z,p} = \sum_{\substack{(w,q) \in E_k\\ \Re(w) < N}} v_{z,p,w,q}t_k^w \log^q(t_k) + t_k^NR_{z,p,N},\]
as above (but this time for all $(z,p) \in E_1$).

Using the definitions above, we may then write for all $N$ \[v = \sum_{\substack{(w,q) \in E_k\\ \Re(w) < N}} v^k_{w,q} t_k^w\log^q(t_k) + t_k^NR^k_N.\]
Since $v^k_{w,q} \in \phgi{\tilde{K}}{\tilde{\mathcal E''}}(\{t_k = 0\})$ and $R^k_N \in \phgi{\tilde{K}}{\tilde{\mathcal E''}}(X)$, writing this expansion for all $k$ shows that $v \in \phgi{J}{\mathcal E}(X)$.\end{proof}

\numberwithin{equation}{chapter}
\chapter{The proof of \texorpdfstring{\cref{thm:C3:SPconversion}}{lemma 3.2.16}}
\label{C:A2:dumblemma}
\todo{make sure reference in texorpdftostring is correct}

\theoremstyle{plain}
\newtheorem*{lemspec}{\Cref{thm:C3:SPconversion}}

We restate and prove \cref{thm:C3:SPconversion}. Recall that we are working in the case $\mathcal R = [0,a)\times[0,b)$ for some $a,b > 0$. Set $\mathcal H_1 = u^{-1}(0)$, $\mathcal H_2 = v^{-1}(0)$, and fix $h \in \Sym^2(\spT)$ which is $0$ on $\mathcal H_1 \un \mathcal H_2$. Then:
\begin{lemspec}Let $(\mathcal M,g,x,u,v,N)$ be a short-pulse double-null gauge, and set $g_x = g+xh$. Fix $a' < a$ and $b' < b$. Then there exists an open set $U \subseteq \mathcal M$ containing $S_{0,0,0}$, $\epsilon > 0$ and a unique smooth map $\phi:U \to \mathcal M$, a diffeomorphism onto its image, the identity on $(\mathcal H_1 \un \mathcal H_2)\n U$ and $\{x = 0\} \n U$, with $\mathcal M([0,a')\times[0,b'))\n \{0 \leq x < \epsilon\}$ in the range of $\phi$, and satisfying $\phi^\ast x = x$, such that $\phi_{\ast}g_x$ is in a double-null gauge with $x,u,v,N$.\end{lemspec}
\begin{proof}
It is slightly inconvenient that $(0,0)$ does not lie in the interior of $[0,a)\times [0,b)$, so, like in \cref{thm:C2:onbd}, we will need to extend $\mathcal M$. For $\delta$ sufficiently small, replace $[0,a)\times [0,b)$ with $(-\delta,a)\times(-\delta,b)$, and using the trivialization \[\mathcal M \iso [0,\delta)\times ([0,a)\times[0,b))\times S_{0,0,0},\] extend $\mathcal M$ to
\[\widetilde{\mathcal M} = [0,\delta)\times ((-\delta,a)\times(-\delta,b))\times S_{0,0,0},\] and extend $(u,v)$ to $\widetilde{\mathcal M}$ so that $(\widetilde{\mathcal M},x,u,v)$ is a doubly-foliated manifold. Extend $\mathcal H_i$ ($i=1,2$) in the same way. Extend $N$, $\slash{g}$, and $\Omega^2$ to $\widetilde{\mathcal M}$, which then extends $g$ by \cref{thm:C3:ReconstructSP}. Extend $h$ to $\widetilde{\mathcal M}$, which then also extends $g_x$.

For this proof, let $T\widetilde{M}$ denote the bundle over $\widetilde{\mathcal M}$ consisting of those vectors tangent to the level sets of $x$. We also interpret $d$ as taking values in $C^\infty(\widetilde{\mathcal M};T^\ast \widetilde{M})$ (i.e.\ we quotient by $\vspan \{dx\}$). 

The metric $g_x$ is certainly Lorentzian is $x$ if small enough, provided we restrict to some compact subset of $\widetilde{\mathcal M}$ (so we may take $U$ to be compactly contained). This also gives $g_x$ the same time-orientation as $g$.

Uniqueness is clear from \cref{thm:C2:diffeo}.

The goal now is to find for $\epsilon$ small enough an open set $U$ containing $S_{0,0}\n \{x \leq \epsilon\}$, a pair of optical functions $u',v'$ for $g_x$ defined over $U$ equal to $u,v$ at $x=0$, such that $(x,u',v')$ turns $(U,g,x,u',v')$ into a double-null foliation, and such that for some $\delta' < \delta$ and $\epsilon$ sufficiently small, $(-\delta',a')\times (-\delta',b')$ in the range of $(x_0,u',v')$ for $x_0 \leq \epsilon$. If this is true, we may construct the diffeomorphism $\phi$ as in \cref{thm:C2:badconversion} and \cref{thm:C2:conversion} to show that $\phi_\ast g_x$ is in double-null gauge with $u,v,\widetilde{M}$, and use \cref{thm:C2:onbd} to replace $\widetilde{\mathcal M}$ with $\mathcal M$.\footnote{However, unlike in \cref{thm:C2:onbd}, there is no need to shrink $U$, since $du = du'$ and $du = dv'$ on $\{x = 0\}$ implies that for small enough $x$, $g(-2du', dv) > 0$ and $d(-2dv',du) > 0$.}, That $\phi_\ast g_x$ is still a short-pulse metric follows from the fact that $\phi|_{\widetilde{M}_0} = \text{id}$, and so for an ordinary vector field $X$, $\phi_\ast X = X+ xR$, for another ordinary vector field $R$ (i.e.\ $xR$ is a short-pulse vector field, regardless of what $X$ was). 

For the rest of the proof, we will only deal with $\widetilde{\mathcal M}$, so for convenience let us replace $\widetilde{\mathcal M}$ in our notation with just $\mathcal M$. The main part of the proof now consists of two steps. In the first we construct $u',v'$ on a large set $U$ using Hamilton-Jacobi theory, and in the second step we use elementary topology to show that $U$ is in fact large enough (i.e.\ $(u',v')$ has large enough range).

We may find canonical coordinates $u,v,\theta^1,\ldots,\theta^n,x$ on $\mathcal M$ for $g$ (not $g_x$) and use them to construct dual coordinates $x,u,v,\theta^1,\ldots,\theta^n,\mu,\nu,\tau^1,\ldots,\tau^n$ on $T^\ast M$, with $\mu$ dual to $u$, $\nu$ dual to $v$, and $\tau^i$ dual to $\theta^i$. Write $\xi = (\mu,\nu,\tau^1,\ldots,\tau^n)$ for a general coordinate on the fibres. Define \[
q(\xi) := -\Omega^2 g_x^{-1}(\xi,\xi),\] a smooth Hamiltonian function on $^\sp T^\ast \mathcal M$, and hence on $T^\ast M$ as well. Notice also that $q(\xi) = 0$ if and only if $\xi$ is a null covector.

At $x=0$, 
\[g^{-1} = -\frac{1}{2\Omega^2}(\pa_u \otimes \pa_v + \pa_v \otimes \pa_u)\] the $\theta^i$ derivatives dropping out since they come with a factor of $x$. Thus at $x=0$ it follows that $q = \mu\nu$, and so the Hamiltonian vector field is
\begin{equation}\label{eq:A2:H}H_q = \nu\pa_u + \mu\pa_v.\end{equation}

Let us find $u'$, first. Let $\mathcal H_1^\ast \subseteq T^\ast M$ be the surface over $\mathcal H_1$ whose points consist of $(p,du_p)$ for $p \in \mathcal H_1$. We will construct $du'$ as the flow via $H_q$ of $\mathcal H_1^\ast$.

Now, fix $\delta' < \delta_1 < \delta$ and $a' < a_1 < a$, and replace $\mathcal H_1^\ast$ with $\mathcal H_1^\ast \n \{-\delta_1 \leq u \leq a_1\}$, turning $\mathcal H_1^\ast$ into a compact set.
\begin{claim}For for $x$ fixed, $\mathcal H_1^\ast\n T^\ast M_x$ is isotropic (i.e.\ the symplectic form pulls back to $0$), $q|_{\mathcal H_1^\ast} = 0$, and for $x$ small enough, $H_q$ is transverse to $\mathcal H_1^\ast$.\end{claim}
\begin{proof}In canonical coordinates the symplectic form is
\[\omega = d\mu\wedge du + d\nu\wedge dv+ \sum d\tau^i \wedge d\theta^i,\] and in coordinates \[\mathcal H_1^\ast = \{(x,u,0,\theta,1,0,0)\}.\] From this, it is clear $\omega$ pulls back to $0$. From the description of $H_q$ above, it is also clear that $H_q$ is transverse to $\mathcal H_1^\ast$ at $x=0$, and thus also for small $x$, $H_q$ remains transverse to $\mathcal H_1^\ast$.
Since $g_x^{-1}(du,du) = 0$ by assumption (even for $x > 0$) , it follows that $q|_{\mathcal H_1^\ast} = 0$.
\end{proof}

Provisionally define for $0 < S \leq \delta_1$ and $0 < T < b$ the map
\[\Phi:\mathcal H_1^\ast\n\{x \leq \epsilon\} \times [-S,T] \to T^\ast \widetilde{M},\]
the flowout from $\mathcal H_1^\ast$ by $H_q$ (which may not a priori exist for long enough times). Throughout the proof we will see how small $\epsilon$ needs to be, and how to choose $S$ and $T$ appropriately.

Set $\Lambda$ to be its image, and $\pi:\Lambda \to \widetilde{\mathcal M}$ the projection map, and set $U_1 = \pi(\Lambda)$. Let us denote $\kappa = \pi\circ \Phi$.
\begin{claim}Suppose $0 < S \leq \delta_1$, $0 < T < b$ and $\epsilon > 0$ is sufficiently small. Then $\Phi$ is well-defined and in fact $\Phi,\kappa$ are diffeomorphisms onto $\Lambda$ and $U_1$, respectively. \end{claim}
\begin{proof}Let us examine how $\Phi$ and $\kappa$ act at $x=0$. In canonical coordinates, fix $p = (0,u,0,\theta,1,0,0) \in \mathcal H_1^\ast$ and $t \in (-\delta,b)$. 
Over $x=0$ it is clear from \eqref{eq:A2:H} that $\Phi$ takes $(p,t)$ to $(0,u,t,\theta,1,0,0)$, which is just $du$ over the point $(0,u,t,\theta)$. Thus $\kappa(p,t) = (0,u,t,\theta)$. So $\Phi|_{\{x = 0\}}$ is defined on $\mathcal H_1^\ast \times [0,b)$, and $\Phi|_{x=0}$, $\kappa|_{x=0}$ are diffeomorphisms. Using ODE theory, it follows that if $S \leq \delta_1$, $T < b$ and $\epsilon$ is small enough, $\Phi$ is defined on $\mathcal H_1^\ast \n\{x \leq \epsilon\}\times [0,T]$. Indeed, since $\mathcal H_1^\ast$ is compact and is contained inside an open subset of $T\widetilde{M}$, the flow of $H_q$ exists for time $T$ and the integral curves of $H_q$ for $x > 0$ stay close to those of $H_q$ at $\{x= 0\}$ (since the data are the same).

Now let us show that $\Phi$ and $\kappa$ are diffeomorphisms if $\epsilon$ is small enough. Since they are diffeomorphisms at $\{x = 0\}$, in particular their Jacobians are non-singular there. Since the domain of $\Phi$, and $\kappa$ is compact, choosing $\epsilon$ small, it it thus clear that if $x_0 \leq \epsilon$, the Jacobians of $\Phi_{\{x = x_0\}}$, $\kappa_{\{x = x_0\}}$ have full rank. Since $\Phi$, $\kappa$ also preserve the level sets of $x$, it follows that $\frac{\pa \Phi}{\pa x}$, $\frac{\pa \kappa}{\pa x}$ are non-zero. Thus the (total) Jacobians of $\Phi$, $\kappa$ have full rank, and thus $\Phi$ and $\kappa$ are local diffeomorphisms. Since $\Phi$, $\kappa$ are proper,\footnote{The preimage of any closed set is a closed subset of a compact set.} they are covering maps. Since $\Phi$ and $\kappa$ are diffeomorphisms at $x=0$, all points in $\{x = 0\}$ are covered once, which means that in fact $\Phi$ and $\kappa$ are global diffeomorphisms.\end{proof}

\begin{claim}There exists a closed one-form $\alpha$ (i.e.\ a section of $T^\ast \tilde{M}$ defined on $U_1$) such that $q(\alpha) = 0$, and $\alpha|_{\mathcal H_1} = du$, and $\alpha|_{x=0} = du$.\end{claim}
\begin{proof}Since $\mathcal H_1^\ast\n T^\ast M_x$ is isotropic, and $H_q$ is a Hamiltonian vector field not tangent to $\mathcal H^\ast_1$, the flowout of $\mathcal H_1^\ast\n T^\ast M_x$ is Lagrangian. Thus $\Lambda \n M_x$ is Lagrangian for each $x$. Since $\Phi$, $\kappa$ are diffeomorphisms, so is $\pi$. Let us define $\alpha$ by $\alpha(p) = \pi^{-1}(p)$. Then $\alpha|_{\mathcal H_1} = du$ by definition. Observe that since we know $\Phi|_{x=0}$, we can deduce that $\alpha(p) = du$ for $x(p) = 0$. Since $\Lambda \n M_x$ is Lagrangian, $d\alpha = 0$. We already know $q(\alpha) = q(du) = 0$ over $\mathcal H_1$. Let $(x(s),\xi(s))$ be a bicharacteristic of $H_q$ starting at $\mathcal H_1^\ast$. By definition, $\alpha(x(s)) = \xi(s)$. Also, $\frac{d}{ds} q(x(s),\xi(s)) = H_q q = 0$, and thus $q(\alpha) = 0$ everywhere.\end{proof}

\begin{claim}There exists a function $u'$ defined over $U_1$ such that $du' = \alpha$ and $u' = u$ over $\mathcal H_1$ and at $x=0$.\end{claim}
\begin{proof}
We may cover $U_1$ by contractible open sets $V_i$, obtaining a collection of functions $u'_i$ defined over $V_i$ with $du'_i = \alpha$. For each $i,j$, there is thus a function $f_{ij}$, only a function of $x$, defined over $V_i \n V_j$ such that $u'_i - u'_j = f_{ij}(x)$. Define
\begin{gather*}
 \mathcal S^+ = \{t \in [0,T]\: u' \text{ exists on a neighbourhood of }\kappa(\mathcal H_1^\ast \times [0,t])\} \\
  \mathcal S^-= \{t \in [-S,0]\: u' \text{ exists on a neighbourhood of }\kappa(\mathcal H_1^\ast \times [t,0])\}
.\end{gather*}
We will show that $\mathcal S^\pm$ are open, closed and non-empty, and thus $u'$ exists everywhere. The proofs for $\mathcal S^{\pm}$ are the same, so we only treat the case of $\mathcal S^+$.

First let us show open. Suppose $t \in \mathcal S^+$. We may cover $\kappa(\mathcal H_1^\ast \times \{t\})$ by a finite subcollection of the open sets $V_i$ (where we require each $V_i$ to actually intersect $\kappa(\mathcal H_1^\ast \times \{t\})$). The function $u'$ is defined on an open subset of each $V_i$ by assumption. On each, we have a $u'_i$ with $u'_i - u'$ only depending on $x$. Modifying this function, we may assume that it is $0$. Thus $u'_i = u'$ where both are defined, and hence $u'_i$ extends $u'$ to all of $V_i$. For different $i,j$, $u'_i-u'_j$ depends only on $x$ (assuming $V_i \n V_j \neq \emptyset)$. However, both are equal to $u'$ somewhere, and hence $u'_i = u'_j$ where both are defined. Thus, we have extended $u'$ to a function on the union of all the $V_i$, which by compactness cover $\kappa(\mathcal H_1^\ast \times (t-\delta,t+\delta))$ if $\delta$ is small enough.

Now let us show non-empty. This is practically the same argument as above. Cover $\mathcal H_1 = \kappa(\mathcal H_1^\ast \times \{0\})$ by a finite subcollection of the open sets $V_i$. By construction, $du'_i = \alpha = du$ over $\mathcal H_1$. In particular, their pull-backs to $\mathcal H_1$ agree, and hence $u'_i - u$ is a function depending only on $x$, which we can take to be $0$. Now continue as above.

Finally let us show closed. Suppose $t_n \in S$, and $t_n \to t$. We wish to show that $t \in S$. This is again the same argument, since any open set intersecting $\kappa(\mathcal H_1^\ast \times \{t\})$ also intersects $\kappa(\mathcal H_1^\ast \times \{t_n\})$ for $n$ large enough.

Thus $du' = \alpha$ and $du' = du$ at $\mathcal H_1$. By construction $\alpha = du$ at $x=0$, and so $du' = du$, and hence $u'=u$ since they agree on $\mathcal H_1$.
\end{proof}

Now let us do the same thing with $v'$, establishing an analogous region $U_2$ for its domain. Set $U' = U_1\n U_2$. We may still assume that $U'$ contains an open neighbourhood of $M_0((-\delta_1,a_1)\times (-\delta_1,b_1))$, so long as we chose $b_1 \leq T$ to define $U_1$, $a_1 \leq T$ to define $U_2$, and $S = \delta_1$ to define both.  Observe that $U'$ still depends on $\epsilon$, and only includes points with $\{x \leq \epsilon\}$.

Since $u'$ and $v'$ are smooth and agree with $u$ and $v$, respectively, at $\{x = 0\}$, it follows that $u' = u+O(x)$, $v' = v + O(x)$, where the $O(x)$ terms represent $x$ times a smooth function. Since $U'$ is compact, if $\epsilon$ is small enough, $(x,u',v')$ is a submersion on $U'$ and $u',v'$ are increasing towards the future. Moreover, $(x,u',v')$ is a proper map onto its image.

Let us now show that the range of $(x,u',v')$ is large enough.
\begin{claim}The range of $U' \subseteq \mathcal M$ under $(x,u',v')$ contains an open neighbourhood of $\{0\}\times(-\delta_1,a_1)\times(-\delta_1,b_1)$.\end{claim}
\begin{proof}If $(\alpha,\beta) \in (-\delta_1,a_1)\times(-\delta_1,b_1)$, then $(0,\alpha,\beta)$ is in the range of $(0,u',v') = (0,u,v)$. Say $(0,u',v')(p) = (0,\alpha,\beta)$. Then an open set of $p$ in $\mathcal M$ is contained in $U'$, and so using the submersion property, an open subset of $(0,\alpha,\beta)$ is contained in the range, too.
\end{proof}
Because of this claim, since $\delta' < \delta_1$, $a' < a_1$ and $b' < b_1$, shrinking $\epsilon$ further, we may assume $[0,\epsilon)\times(-\delta',a')\times(-\delta',b')$ is in the range of $(x,u',v')$. Now, set $U$ to be the preimage of $ [0,\epsilon)\times(-\delta',a')\times(-\delta',b')$ under $(x,u',v')$, so that $U$ is an open subset of $\mathcal M$ containing $S_{0,0}\n\{x \leq \epsilon\}$.

The last thing we need is:
\begin{claim}The quadruple $(U,x,u',v')$ is a doubly-foliated manifold.\end{claim}
\begin{proof}We need to check that $(x,u',v')$, in addition to being a surjective submersion, is proper and has compact fibre the same compact manifold $S$ as the fibres $(x,u,v)$, and that the level sets of $(x,u')$ and $(x,v')$ are connected

Being proper (over a path-connected base) implies that all the fibres are diffeomorphic to $S$. Indeed, this is a general property of proper submersions. 
Let us therefore show that $(x,u',v')$ is proper. We know that $(x,u',v')$ is proper from $U'$ onto its image (since $U'$ is compact). It follows that $(x,u',v')$ is proper from $U$ onto its image  $[0,\epsilon)\times(-\delta',a')\times(-\delta',b')$.

The fibres of $(x,u')$ and $(x,v')$ are connected because we may always lift a path from the base, now that we know $(x,u,v)$ is a submersion with compact fibre. \end{proof}

Thus $(U,g_x,x,u',v')$ is the desired double-null foliation.

\end{proof}

\numberwithin{equation}{section}
\chapter{The Goursat problem}
\label{C:A3}

In this appendix we study several aspects of the characteristic initial value problem for hyperbolic equations in spacetime dimension $2$. This is classically known as the \emph{Goursat problem}. For this appendix, we work in a rectangle $\mathcal R = [a,b]_v \times [c,d]_\eta$. We will use the notation $\mathcal R(a_0,c_0)$ to denote the subrectangle $[a,a_0]\times[c,c_0] \subseteq \mathcal R$ for $a \leq a_0 \leq b$, $c \leq c_0 \leq d$.
\section{The classical Goursat problem}
\label{C:A3:Goursat}
We will work with two sorts of system. The first will be for functions $X,Y:\mathcal R \to \R^N$ ($N \in \N$)
\begin{align}
\label{eq:A3:GoursatI}
\begin{split}
\pa_\eta X = P&(v,\eta,X,Y)\\
\pa_v Y = Q&(v,\eta,X,Y)\\
X(v,c) &= X_0(v)\\
Y(a,\eta) &= Y_0(\eta)\\
\end{split}
\end{align}
where $P,Q:\mathcal R\times \R^{2N} \to \R^N$. We will often suppress the explicit dependence of $P$ and $Q$ on $v$ and $\eta$. The second will be for a function $Z: \mathcal R \to \R^N$
\begin{align}
\label{eq:A3:GoursatII}
\begin{split}
\pa_v\pa_\eta Z = P(v,&\eta,\pa_v Z,\pa_\eta Z, Z)\\
Z(v,c) &= Z_0(v)\\
Z(a,\eta) &= Z_1(\eta),
\end{split}
\end{align}
subject to the constraint $Z_0(a) = Z_1(c)$,
where $P:\mathcal R\times\R^{3N} \to \R^N$. We will often not write explicitly the dependence of $P$ on $v$ and $\eta$.

System \eqref{eq:A3:GoursatII} may be reduced to \eqref{eq:A3:GoursatI} if $P$, $Z_0$, $Z_1$ have high enough regularity (which will always be true in practice). We show how to do this. Assume that $\pa_\eta Z_1$ and $\pa_v Z_0$ are continuous.

Now, let us consider a system of type \eqref{eq:A3:GoursatI} for $X = (X_1,X_2), Y = (Y_1,Y_2) \in \R^{2N}$ described by
\begin{gather}
\label{eq:A3:reduction}
\begin{split}
\pa_\eta (X_1,X_2) = (P(X_1,Y_1,X_2),Y_1)\\
\pa_v (Y_1,Y_2) = (P(X_1,Y_1,X_2),X_1)\\
(X_1,X_2)(v,c) = (\pa_v Z _0(v),Z_0(v))\\
(Y_1,Y_2)(a,\eta) = (\pa_\eta Z_0(\eta),Z_1(\eta)).
\end{split}
\end{gather}
Let us assume that we have a continuous solution $(X,Y)$ of \eqref{eq:A3:reduction}, and that $P$ is continuous. Then $\pa_\eta X_1 = \pa_v Y_1$, and thus there is some potential function $\Phi$, unique up to constant, for which $\pa_v \Phi = X_1$, $\pa_\eta \Phi = Y_1$. The initial data assumptions mean that $\Phi(v,c)-Z_0(v)$ and $\Phi(a,\eta) - Z_1(\eta)$ are constants. The constraint $Z_0(a) = Z_1(c)$ implies that they are the same constant. Thus we can assume that $\Phi$ is chosen so that the constant is $0$.

We claim that setting $Z := \Phi$ gives a solution to \eqref{eq:A3:GoursatII}. Indeed, $\pa_\eta \Phi = Y_1 = \pa_\eta X_2$ and $\Phi|_{\{\eta = c\}} = Z_0 = X_2|_{\{\eta = c\}}$ by assumption, and similarly $\pa_v \Phi = X_1 = \pa_v Y_2$ and $\Phi|_{\{v = a\}} = X_2|_{\{v= a\}}$ by assumption.

Thus, $X_2 = \Phi = Y_2$, and so with $Z = \Phi$,
\[\pa_\eta\pa_v Z = \pa_\eta X_1 = P(X_1,Y_1,X_2) = P(\pa_v Z,\pa_\eta Z,Z).\]

With the reduction carried out, we focus most of our attention on \eqref{eq:A3:GoursatI}. The well-posedness theory of \eqref{eq:A3:GoursatI} is analogous to that of an ODE. The notion of a weak solution is well-defined by integrating, i.e.:
\begin{defn}[Weak solution]Fix $a < a_0 \leq b$ and $c < c_0 \leq c$. We say a pair $(X,Y) \in L^{\infty}(\mathcal R(a_0,c_0))$ is a weak solution of \eqref{eq:A3:GoursatI} on $\mathcal R(a_0,c_0)$ if for all $(\eta,v) \in \mathcal R(a_0,c_0)$
\begin{align}
\label{eq:A3:Goursatweak}
\begin{split}
X(v,\eta) &= X_0(v) + \int_c^{\eta} P(v,t,X(v,t),Y(v,t))\ dt\\
Y(v,\eta) &= Y_0(\eta) + \int_a^v Q(s,\eta,X(s,\eta),Y(s,\eta))\ ds.\end{split}\end{align}
\end{defn}
To state a well-posedness theorem, we will need to assume that $P(v,\eta,x,y)$, $Q(v,\eta,x,y)$ are continuous in $v,\eta$, locally uniformly in $x,y$, and are locally Lipschitz in $x,y$, uniformly in $v,\eta$. This means that $P(\cdot,\cdot,x,y)$ is continuous, uniformly in $x,y$ (i.e.\ the modulus of the continuity in $v,\eta$ does not depend on $x,y$ in a compact set) and $|P(v,\eta,x_1,y_1)-P(v,\eta,x_1,y_1)| \lesssim |x_1-x_2| + |y_1-y_2|$ uniformly in $v,\eta$, if $x_1,y_1,x_2,y_2$ are uniformly bounded (and similarly for $Q$).
\begin{thm}[Well-posedness]\label{thm:A3:wp}Consider the system \eqref{eq:A3:GoursatI}. Assume $P$, $Q$ satisfy the previous assumptions and provide continuous data $X_0$, $Y_0$. Then:
\begin{description}
\item[Local well-posedness]
There exists $\epsilon > 0$ depending only on \[M = \max\left(\sup_{\eta} |Y_0(\eta)|,\sup_{v} |X_0(v)|\right),\] and the Lipschitz constant of $P$, $Q$ in a ball of radius $2M$ centred at the origin such that there is a (weak) solution $(X,Y)$ to \eqref{eq:A3:GoursatI} on $C^0(\mathcal R(\epsilon,\epsilon))$;
 \item[Uniqueness]
for any $a \leq a_0 \leq b$, $c \leq c_0 \leq d$, two solutions in $C^0(\mathcal R(a_0,c_0))$ agree;
\item[Automatic continuity] if $(X,Y) \in L^{\infty}([a,a_0)\times[c,c_0))$ is a weak solution, then in fact $(X,Y) \in C^0(\mathcal R(a_0,c_0))$;
\item[Persistence of regularity]
if $X_0,Y_0$ and $P,\ Q$ are smooth and a (weak) solution $(X,Y)$ exists in $C^0(\mathcal R(a_0,c_0))$, then the solution is a classical solution in $C^{\infty}(\mathcal R(a_0,c_0))$;
\item[Dependence on parameters]If we allow $X_0,Y_0$, $P,Q$, as well as the initial data, to depend in a $C^k$ fashion on a parameter $\sigma \in \Sigma$ (for some manifold $\Sigma$), then all statements above hold with parameter. For the local well-posedness to hold, we require $P,Q$, together with their first $k$ $\sigma$-derivatives to be locally Lipschitz in $x,y$, uniformly in $\sigma,v,\eta$ and continuous in $\sigma,v,\eta$. In fact, if a solution exists in $L^\infty(\mathcal R(a_0,c_0)\times \Sigma)$, then it is automatically in $C^k_\sigma C^0_{v,\eta}\mathcal R(a_0,c_0))$, and hence in $C^\infty(\mathcal R(a_0,c_0)\times \Sigma)$ if $P,Q$, and the data are all smooth.
 \end{description}
 \end{thm}
 
 Just like Gronwall's inequality is a key ingredient in the well-posedness theory of ODEs, we will need a two-dimensional version in order to prove \cref{thm:A3:wp}.
 \begin{prop}[2D Gronwall's inequality]\label{thm:A3:Gronwall}Let $X,\ Y$ be non-negative bounded functions on $\mathcal R(a_0,c_0)$, for $a \leq a_0 \leq b$, $c \leq c_0 \leq d$ which satisfy for some $B \geq 0$, $M >0$ and all $(v,\eta) \in \mathcal R(a_0,c_0)$
 \begin{align*}
     X(v,\eta) &\leq B + M\int_0^\eta X(v,t)+Y(v,t)\ dt\\
     Y(v,\eta) &\leq B + M\int_0^\eta X(s,\eta) + Y(s,\eta)\ ds.
 \end{align*}
 Then there is a constant $C = C(M,a,b,c,d)$, increasing in all its arguments, such that $X,Y \leq CB$ on all of $\mathcal R(a_0,c_0)$.\end{prop}
 \begin{proof}
 Let us assume without loss of generality that $a=c=0$.
 Substituting the second inequality into the first gives
 \begin{equation}\label{eq:A3:Vishesh}X(v,\eta) \leq B + BM \eta+ M\int_0^\eta X(v,t)\ dt + M^2\int_0^\eta\int_0^v X(s,t) + Y(s,t)\ dsdt.\end{equation}
 Let us denote by $C$ a quantity depending only on $a,b,c,d,M$, but not $B$, which is allowed to change from line to line. Denote by $Z(v,\eta)$ the double integral in \eqref{eq:A3:Vishesh}. Then
 \[ X(v,\eta) \leq CB + C\int_0^\eta X(v,t)\ dt + CZ(v,\eta).\] Applying (1D) Gronwall's inequality gives
 \begin{equation}\label{eq:Gronwall:eq1}X(v,\eta) \leq C(B+Z(v,\eta)).\end{equation} Running the same argument gives the same bound for $Y$. Thus
 \[Z(v,\eta) \leq CB\eta v + C\int_0^\eta \int_0^v Z(s,t)\ dsdt.\]
 Plugging this bound back into itself repeatedly gives that for all $n \geq 1$
 \begin{align*}
     Z(v,\eta) &\leq B\left(C\eta v + \cdots + \frac{1}{n!}(C\eta v)^n\right)\\
     &+C^n\int_0^\eta\int_0^v\int_0^{\eta_2}\int_0^{v_2}\cdots\int_0^{\eta_n}\int_0^{v_n} Z(s,t)\ dsdtdv_nd\eta_n\cdots dv_2d\eta_2.
 \end{align*}
 Since $X,Y$ are bounded, so is $Z$. Say $Z \leq A$. Then using the previous bound,
 \[Z(v,\eta) \leq B\left(C\eta v + \cdots + \frac{1}{n!}(C\eta v)^n\right) + \frac{1}{n!}C^nA(v\eta)^n.\]
 Taking $n\to \infty$ gives
 \[Z(v,\eta) \leq B(e^{Cv\eta}-1) \leq CB.\]
 Plugging this into into \eqref{eq:Gronwall:eq1} gives
 \[X(v,\eta) \leq C(B+BC) \leq CB,\] and likewise for $Y$. This completes the proof.\end{proof}
 
 With a two-dimensional version of Gronwall's inequality proven, it is not hard to prove \cref{thm:A3:wp}.
 \begin{proof}[Proof \cref{thm:A3:wp}]Existence locally in $C^0$ follows from the contracting mapping principle. For uniqueness, if $(X_1,Y_1)$ and $(X_2,Y_2)$ are both weak solutions, then for some constant $M$ depending on the $L^{\infty}$ norms of $X_1,X_2,Y_1,Y_2$ and the Lipschitz constant of $P,Q$ in a ball of radius $2M$ centred at the origin,
 \begin{align*}
|X_1-X_2|(v,\eta) \leq M\int_c^\eta |X_1-X_2|(v,t) + |Y_1-Y_2|(v,t)\ dt\\
|Y_1-Y_2|(v,\eta) \leq M\int_a^v |X_1-X_2|(s,\eta) + |Y_1-Y_2|(s,\eta)\ ds,
 \end{align*}
 and hence 2D Gronwall's inequality implies $|X_1-X_2| \equiv 0$ and $|Y_1-Y_2| \equiv 0$.
 
 Now let us prove automatic continuity. Using the weak formulation, it is clear that $X$ is Lipschitz in $\eta$, uniformly in $v$, i.e.\
 \[|X(v,\eta_1)-X(v,\eta_2)| \lesssim |\eta_1-\eta_2|\] for any $c \leq \eta_1,\eta_2 \leq c_0$, where the implied constant does not depend on $\eta_1,\eta_2$ or $v$. Likewise, $Y$ is Lipschitz in $v$, uniformly in $\eta$. Thus, to show continuity, it suffices to show that $X$ is continuous in $v$, uniformly in $\eta$ (i.e.\ the modulus of continuity in $v$ does not depend on $\eta)$ and vice-versa for $Y$. Define
 \[\tilde{P}(v,\eta,x) = P(v,\eta,x,Y(v,\eta)).\] Since $Y$ is continuous in $v$, uniformly in $\eta$, the function $\tilde{P}(\eta,v,x)$ is continuous in $v$, uniformly in $\eta$, $x$. Also, $\tilde{P}$ is locally Lipschitz in $x$, uniformly in $v,\eta$. Since
 \[\pa_\eta X = \tilde{P}(X),\] (at least weakly), $X$ satisfies an ODE in $\eta$ depending continuously on a parameter $v$. Since $X_0$ depends continuously on the same parameter, it follow from ODE theory that $X$ is continuous in $v$, uniformly in $\eta$. The same argument works for $Y$.
 
 Next, let us prove persistence of regularity. This is an extension of the previous argument. The equations show that $\pa_\eta X$ exists and is continuous (in both variables), and similarly for $\pa_v Y$. Thus $\tilde{P}$ is now $C^1$ in $v$. Thus $X$ solves an ODE depending in a $C^1$ fashion on a parameter $v$, which means that $X$ is also $C^1$ in $v$. Thus $X \in C^1(\mathcal R(a_0,c_0))$. The same argument works for $Y$. Iterating, we establish that $X,Y$ are of class $C^{\infty}$.
 
 If we introduce a parameter $\sigma$, then for any $k$ we may (by the contraction mapping principle) obtain local well-posedness in the space $C^k_\sigma C^0(\mathcal R(\epsilon_k,\epsilon_k))$, for $\epsilon_k$ small enough. The previous statements hold with parameter with proofs almost unchanged, so the only thing we need to prove is that if a solution is in $L^\infty(\mathcal R(a_0,c_0))\times \Sigma)$, then it is automatically in $C^k_\sigma C^0_{v,\eta}\mathcal R(a_0,c_0))$. Differentiating \eqref{eq:A3:GoursatI} $\pa_\sigma$,\footnote{For simplicity we assume $\Sigma$ is an compact interval, although the proof in the general case requires only a small change in notation.} we see that
 \begin{align*}
 \pa_\sigma X &= (\pa_\sigma P)(X,Y) + (\pa_x P)(X,Y)\pa_\sigma X + (\pa_yP)(X,Y)\pa_\sigma Y\\
 \pa_\sigma Y &= (\pa_\sigma Q)(X,Y) + (\pa_x Q)(X,Y)\pa_\sigma X + (\pa_yQ)(X,Y)\pa_\sigma Y
 \end{align*}
 (or more appropriately, this holds in a weak sense).
 We already know (by a version of automatic continuity with parameter) that $X,Y \in L^{\infty}_\sigma C^0(\mathcal R(a_0,c_0)$. But this system is linear, so \cref{thm:A3:linearexists}, below, applies, and so there exists a global solution in $\mathcal R(a_0,c_0)$. In particular, solutions of this equation in $\mathcal R(a',c')$, for $a \leq a' \leq a_0$ and $c \leq c' \leq c_0$ are bounded independent of $a',b'$. This allows us to iterate the local well-posedness in $C^1_\sigma C^0$ to show that $X \in C^1_\sigma C^0(\mathcal R(a_0,c_0))$. Iterating this with more $\pa_\sigma$-derivatives proves the claim.
 \end{proof}
 
 We now state and prove a couple corollaries of the statements of the first three parts of \cref{thm:A3:wp} (since we need them to prove the fourth part).
 \begin{cor}\label{thm:A3:globalsolution}Fix continuous initial data $X_0$, $Y_0$ for \eqref{eq:A3:GoursatI}. If there exists $M > 0$ such that for any $a \leq a_0 \leq b$, $c \leq c_0 \leq d$, any bounded (weak) solution $(X,Y)$ on $\mathcal R(a_0,c_0)$ is bounded above by $M$, then there exists a global (weak) solution in $C^0(\mathcal R(a,b))$.
\end{cor}
\begin{proof}
It suffices to show that the set $\mathcal S$ of points $(a_0,c_0)$ for which there is a continuous solution on $\mathcal R(a_0,c_0)$ is open, closed, and non-empty. Non-empty is trivial since $(a_0,c_0) = (0,0)$ is in $\mathcal S$. Closed is also clear: if $(a_0^n,c_0^n) \to (a_0,c_0)$, then by assumption $(X,Y)$ is bounded by $M$ on each rectangle $[a,a+a_0^n]\times[c,c+c_0^n]$, and hence also on $[a,a+a_0)\times[c,c+c_0)$. By automatic continuity, $X,Y$ exist as bounded and continuous solutions on the closure, too, so $(a_0,c_0) \in \mathcal S$. Now for open. Suppose $(a_0,c_0) \in \mathcal S$. It suffices to show that for $\epsilon > 0$ small $(a_0+\epsilon,c_0) \in \mathcal S$ and $(a_0,c_0+\epsilon) \in \mathcal S$, since then we may run local well-posedness again from $(a_0+\epsilon,c_0+\epsilon)$ to obtain $(a_0+\epsilon,c_0+\epsilon) \in \mathcal S$, too (perhaps after shrinking $\epsilon$).

Without loss of generality we show the first. By local well-posedness, we know that there is $\epsilon > 0$ small, depending only on the $C^0$ norms of $Y$  on $\{v = a_0\} \times\{c \leq \eta \leq c_0\}$ and $X$ on $\{\eta = 0\}$,(which by automatic continuity and our assumption we may bound above uniformly by $M$), and the Lipschitz constant of $P,Q$ on a ball of radius $2M$ centred at the origin such that the solution continues to $[a_0,a_0+\epsilon]\times[c,c+\epsilon]$. Let $\mathcal T$ denote the set of $c_1 \leq c_0$ for which the solution continues to $[a_0,a_0+\epsilon]\times[c,c+c_1]$. We have just shown that $\mathcal T$ is non-empty. That it is closed follows in the same way that $\mathcal S$ was closed. Openness is also easy, since we may just iterate local well-posedness with the same $\epsilon$, because by assumption the data are always bounded above by $M$. Thus $\mathcal T = [c,c_0]$, and so the solution continues to $[a,a+a_0+\epsilon]\times[c,c_0]$, and so $(a_0+\epsilon,c_0) \in \mathcal S$.
\end{proof}

An important application of this corollary will be in the case that $P$ and $Q$ are affine linear, i.e.\ there exist continuous functions $f(v,\eta)$, $g(v,\eta)$, and continuous functions $S(v,\hspace{-0.25pt}\eta\hspace{-0.25pt},x\hspace{-0.25pt},y\hspace{-0.25pt})$, $T(v,\eta,x,y)$ linear in $(x,y)$, such that $P = f + S$ and $Q = g + T$. \begin{cor}\label{thm:A3:linearexists}Suppose that $P$ and $Q$ in \eqref{eq:A3:GoursatI} are affine linear. Then there exists a global (weak) solution $(X,Y)$ to \eqref{eq:A3:GoursatI} defined on all of $\mathcal R$.\end{cor}
\begin{proof}Suppose a weak solution $(X,Y)$ exists on a subrectangle $\mathcal R(a_0,c_0)$. Then, since $S,\ T$ are linear, there exists a $C$, not depending on $a_0,c_0$, or $v,\eta$,$x,y$ such that
\begin{align*}
|X(v,\eta)| \leq C + C\int_c^\eta |X(v,t)| + |Y(v,t)|\ dt\\
|Y(v,\eta)| \leq C + C\int_a^v |X(s,\eta)| + |Y(s,\eta)|\ ds.\end{align*}
2D Gronwall's inequality now implies that $|X|,|Y| \leq M$, where $M$ only depends on $C$, and hence not on $a_0,c_0$ or $v,\eta$, $x,y$ such that $|X|,|Y| \leq M$. Thus we may invoke \cref{thm:A3:globalsolution} to conclude.\end{proof}
Finally, we may use all parts of \cref{thm:A3:wp} to conclude:
\begin{cor}\label{thm:A3:globalsmooth}Assume that $P$ and $Q$ and the initial data are smooth.

Assume first that the hypotheses of \cref{thm:A3:globalsolution} are met. Then there exists a unique global smooth solution to \eqref{eq:A3:GoursatI} $(X,Y) \in C^{\infty}(\mathcal R)$. If $P$ and $Q$ and the data depend smoothly on a parameter $\sigma \in \Sigma$, then $(X,Y) \in C^{\infty}(\mathcal R\times \Sigma)$.

If instead $P$ and $Q$ are affine linear and smooth, then there exists a unique global solution $(X,Y)$ to \eqref{eq:A3:GoursatI} in $C^{\infty}(\mathcal R)$. If the data and $f$, $g$, $S$, $T$ depend smoothly in a parameter $\sigma \in \Sigma$, then $(X,Y) \in C^{\infty}(\mathcal R\times \Sigma)$.

Finally, suppose we consider system \eqref{eq:A3:GoursatII}. If the data are smooth and $P$ is smooth and affine linear (in the sense that $P(v,\eta,x,y,z) = f(v,\eta) + S(v,\eta,x,y,z)$ for smooth $f,S$, and $S$ linear in $x,y,z$), and the data are smooth, then there exists a unique global smooth solution $Z \in C^{\infty}(\mathcal R)$. If the data and $f$, $S$ depend smoothly on a parameter $\sigma$, then $Z \in C^{\infty}(\mathcal R\times \Sigma)$.
\end{cor}
\begin{proof}The first and second statement follow from \cref{thm:A3:globalsolution} and \cref{thm:A3:linearexists}, respectively, used in conjunction with the persistence of regularity and dependence on parameters statement of \cref{thm:A3:wp}. The final statement follows from the second statement together with the reduction carried out at the start of this section.\end{proof}

\section{The singular Goursat problem}
\label{C:A3:Sing}
We will have cause as well to treat a  singular (albeit linear) version of \eqref{eq:A3:GoursatI}, where  singularity with respect to a derivative, and a singular version of \eqref{eq:A3:GoursatII}. We will assume $c = 0$, i.e.\ $\mathcal R = [a,b]_v \times [0,d]_\eta$ and consider the problems \begin{align}
\label{eq:A3:GoursatSingI}
\begin{split}
\eta\pa_\eta X + p^1(v,\eta)X(v,\eta) &+ q^1(v,\eta)Y(v,\eta) = F(v,\eta)\\
\pa_v Y + p^2(v,\eta)X(v,\eta) &+ q^2(v,\eta)Y(v,\eta) = G(v,\eta)\\
X(v,0) &= 0\\
Y(a,\eta) &= Y_0(\eta)
\end{split}
\end{align}
for smooth $p^i,q^i$ ($i = 1,2$), $F$, $G$,
and
\begin{align}
\label{eq:A3:GoursatSingII}
\begin{split}
\eta\pa_\eta\pa_v Z + r^1(v,\eta)\pa_v Z(v,\eta) + r^2(v,\eta)\eta\pa_\eta &Z(v,\eta) + r^3(v,\eta) Z(v,\eta) = H(v,\eta)\\
Z(v,0) &= 0\\
Z(a,\eta) &= Z_1(\eta),
\end{split}
\end{align}
subject to the constraint $Z_1(0) = 0$, for smooth $r^i$ $(i = 1,2,3)$, $H$.

While it is possible to develop a general theory for these systems, we will only need it for restrictive assumptions on the data and right-hand sides. For any $0 \leq j \leq \infty$, let $\eta^\infty C^j = \bigcap_{k=1}^{\infty} \eta^k C^j$.

We prove:
\begin{thm}\label{thm:A3:GoursatSing}Consider system \eqref{eq:A3:GoursatSingI}, and suppose additionally that $F,G \in \eta^\infty C^{\infty}(\mathcal R)$, $Y_0 \in \eta^\infty C^{\infty}(\{v = a\})$. Then there exists a unique solution $(X,Y)  \in\eta^\infty C^{\infty}(\mathcal R)$. If $F,G$ and the data depend in a smooth fashion on a parameter $\sigma \in \Sigma$ (for some manifold $\Sigma$), then so does $(X,Y)$.

Now consider system \eqref{eq:A3:GoursatSingII}, and suppose additionally that $H \in \eta^\infty C^{\infty}(\mathcal R)$, $Z_1 \in \eta^\infty C^\infty(\{v = a\})$. Then there exists a unique solution $Z$ with $Z \in \eta^\infty C^{\infty}(\mathcal R)$. If $H$ and the data depend in a smooth fashion on a parameter $\sigma \in \Sigma$ (for some manifold $\Sigma$), then so does $Z$.\end{thm}
\begin{rk}Notice that we say a unique solution in $\eta^\infty C^{\infty}(\mathcal R)$, since we do not preclude the possibility of other solutions with less decay. For a toy case, consider the ODE $\eta\pa_\eta X - NX = 0$, for $N > 0$ (with $X(0) = 0$). Then $X = 0$ is a solution to this equation, but so is $X = \eta^N$. There will not be a unique solution which is $0$ at $\eta = 0$, but there will always be a unique solution which is rapidly vanishing.\end{rk}
\begin{rk}
Of course, in the special case that $Y_0$ does not depend upon $v$, and $q^1,p^2,q^2 \equiv 0$, then \eqref{eq:A3:GoursatSingI} is essentially an ODE, and the theorem applies to the ODE case as well. In fact we prove the theorem first for ODEs in the following proof.\end{rk}
\begin{proof}
The reduction from \eqref{eq:A3:GoursatII} to \eqref{eq:A3:GoursatI} works equally well to reduce \eqref{eq:A3:GoursatSingII} to \eqref{eq:A3:GoursatSingI} if $\pa_\eta$ is replaced with $\eta\pa_\eta$ and $Y_1$ is replaced with $\eta Y_1$ (so that the sense that $\eta Y_1 = \eta\pa_\eta Z$ is kept), so we will only focus on proving the first statement of the theorem.

Let us first treat the ODE case since the key ideas are already present in the simpler setting.

\textbf{ODE case.} Consider the equation \begin{equation}\label{eq:GoursatSing:ODE}\eta\pa_\eta X + pX = F,\end{equation} where $p \in C^{\infty}([0,d])$, $F \in \eta^{\infty}C^{\infty}([0,d])$. We are searching for a unique solution $X \in \eta^{\infty}C^{\infty}([0,d])$.

We recast the equation as an integral equation
\begin{equation}\label{eq:GoursatSing:eq1}X(\eta) = \int_0^\eta t^{-1}F(t) - t^{-1}p(t)X(t)\ dt.\end{equation}

Let us fix $N$ large (how large will become clear during the proof) and set $\tilde{X}(\eta) = \eta^{-N} X(\eta)$ and $\tilde{F}(\eta) = \eta^{-N}\int_0^\eta t^{-1}F(t)\ dt$, which turns \eqref{eq:GoursatSing:eq1} into
\begin{equation}\label{eq:GoursatSing:eq2}\tilde{X}(\eta) = \tilde{F}(\eta) - \eta^{-N}\int_0^\eta t^{N-1} p(t)\tilde{X}(t)\ dt =:  \tilde{F}(\eta)-P\tilde{X}.\end{equation}
Observe that $\tilde{F} \in C^\infty([0,d])$ since $F \in \eta^{\infty}C^\infty([0,d])$ by assumption.

Now, for $C$ depending only $p,d$, the estimate
\[\norm{PA}_{C^0([0,d])} \leq \frac{C}{N}\norm{A}_{C^0([0,d])}\]
holds for any $A \in C^0([0,d])$. Thus, for $N$ large, $P$ is a contraction and we conclude that \eqref{eq:GoursatSing:eq2} has a unique solution in $C^0([0,d])$, and hence for $N$ large there is a unique solution $X \in \eta^N C^0([0,d])$ to \eqref{eq:GoursatSing:ODE}. 

Now for smoothness, if $N$, $M$ are large enough, we may write $X = \eta^{N+M}\tilde{\tilde{X}}$, where $\tilde{\tilde{X}} \in C^0([0,d])$. Hence $\tilde{X} = \eta^M\tilde{\tilde{X}}$ solves \eqref{eq:GoursatSing:eq2}. If $M$ is large enough, we may differentiate \eqref{eq:GoursatSing:eq2} to conclude that $\tilde{X} \in C^1([0,d])$, and hence $X \in \eta^NC^1([0,d])$. Iterating, we conclude that $X \in \eta^NC^{\infty}([0,d])$ for any $N$, and hence $X \in \eta^{\infty}C^{\infty}([0,d])$.

Next, let us briefly mention how to treat the case with parameter. We simply perform the arguments in $C_k^\sigma \eta^N C^0([0,d])$, etc.

We move on to proving the 2D case. For simplicity, we ignore the parameter, only mentioning that treating it requires little change to the argument.

\textbf{Local existence in $C^0$.}
We recast the equations as an integral system
\begin{align}
\label{eq:GoursatSing:eq3}
\begin{split}
X(v,\eta) &= \int_0^\eta t^{-1}(F(v,t) -  p^1(v,t)X(v,t) - q^1(v,t)Y(v,t))\ dt\\
Y(v,\eta) &= Y(a,\eta) +\int_a^v G(s,\eta) - p^2(s,\eta)X(s,\eta) - q^2(s,\eta)Y(s,\eta)\ ds.
\end{split}
\end{align}

Set $\tilde{X} = \eta^{-N}X$, $\tilde{Y} = \eta^{-N}Y$, $\tilde{F} = \eta^{-N}\int_0^\eta t^{-1}F\ dt$, $\tilde{G} = \eta^{-N}\int_0^v G\ ds$. Then
\begin{align}
\label{eq:GoursatSing:eq4}
\begin{split}
\tilde{X}(v,\eta) &= \tilde{F}(v,\eta) - \eta^{-N}\int_0^\eta t^{N-1}p^1(v,t)\tilde{X}(v,t) + t^{N-1}q^1(v,t)\tilde{Y}(v,t)\ dt\\
\tilde{Y}(v,\eta) &= \tilde{Y}(a,\eta) + \tilde{G}(v,\eta) - \int_a^v p^2(s,\eta)\tilde{X}(s,\eta) + q^2(s,\eta)\tilde{Y}(s,\eta)\ ds.
\end{split}
\end{align}
As in the ODE case, $\tilde{F}, \ \tilde{G} \in \eta^{\infty} C^{\infty}(\mathcal R)$.

Now, arguing as in the ODE case to treat the first line, if $N$ is large enough, the contraction mapping principle provides for $\epsilon > 0$ small and a unique solution in
\[C^0([a,\epsilon]\times[0,d]).\]

\textbf{Global existence in $C^0$}.
We seek to upgrade this to $\epsilon = b$, provided $N$ is perhaps larger.
Let us define
\begin{align*}
I(v) = \sup_{0 \leq \eta \leq d} |\tilde{X}(v,\eta)|\\
J(v) = \sup_{0 \leq \eta \leq d} |\tilde{Y}(v,\eta)|
.\end{align*}
Assume $\tilde{X},\tilde{Y}$ exist on $\mathcal R(a_0,d)$ for $a \leq a_0 \leq b$. We show first that $\tilde{X}$, $\tilde{Y}$ are bounded on $\mathcal R(a_0,d)$ by a constant independently of $a_0$, and then that $\tilde{X}$, $\tilde{Y}$ are uniformly continuous, where the modulus of continuity does not depend on $a_0$. This allows us to iterate local well-posedness and conclude.

The first line of \eqref{eq:GoursatSing:eq4} implies that for some $C$ depending only on $p^1$, $q^1$,
\[I(v) \leq \norm{\tilde{F}(v,\cdot)}_{C^0}+ CN^{-1}I(v) + CN^{-1}J(v).\] In particular, if $N$ is large enough, we may rearrange to obtain
\begin{equation}\label{eq:GoursatSing:eq5}I(v) \lesssim 1 + J(v),\end{equation} where the implied constant depends only on $\norm{\tilde{F}(v,\cdot)}_{C^0}$, $p^1$, $q^1$.
Plugging this into the second line of \eqref{eq:GoursatSing:eq4} implies that
\[J(v) \lesssim 1 + J(a) + \int_a^v I(s) + J(s)\ ds \lesssim 1 + J(a) + \int_a^v J(s)\ ds,\] where the implied constants depend only on $p^i,q^i,\tilde{F},\tilde{G}$ $(i = 1,2)$. Thus Gronwall's inequality gives the bound $J(v) \lesssim 1$, independently of $a_0$, and hence by \eqref{eq:GoursatSing:eq5} $I(v) \lesssim 1$, too

Now we show uniform continuity, with the modulus of continuity not depending on $a_0$. For brevity during the proof, we will omit mentioning explicitly that the modulus of continuity does not depend on $a_0$, and by ``uniform continuity'' always mean that the modulus of continuity is independent of $a_0$. 

We first show that $\tilde{X}$ is uniformly continuous. As a start, we will show that $\tilde{X}$ is uniformly continuous in $\eta$, uniformly in $v$. Using the same argument as above, together with the fact that $\tilde{F}$ and $\tilde{G}$ are in $\eta^\infty C^0(\mathcal R)$, and $\tilde{Y}(a,\eta) \in \eta^{\infty}C^0(\{v = a\})$, we may deduce that
\begin{align*}
I(v,\eta_0) := \sup_{\eta \leq \eta_0} |\tilde{X}|(v,\eta) \\
J(v,\eta_0) := \sup_{\eta \leq \eta_0} |\tilde{Y}|(v,\eta)
\end{align*}
are both in $o_{\eta_0}(1)$, where the $o_{\eta_0}(1)$ does not depend on $v$ or $a_0$.

Thus for $0 \leq \eta_1 < \eta_2$, we deduce from \eqref{eq:GoursatSing:eq4} that
\begin{equation}\label{eq:GoursatSing:eq11}|\tilde{X}(v,\eta_2)-\tilde{X}(v,\eta_1)| \leq |\tilde{F}(v,\eta_2)-\tilde{F}(v,\eta_1)| + \left(1-\left(\frac{\eta_1}{\eta_2}\right)^N\right) \frac{C}{N}o_{\eta_2}(1),\end{equation}
where $C$ depends only on $p^1, q^1$. The estimate
\begin{equation}\label{eq:GoursatSing:eq55}L(\eta_1,\eta_2) := \left(1-\left(\frac{\eta_1}{\eta_2}\right)^N\right) \frac{C}{N}o_{\eta_2}(1) \in o_{\eta_2-\eta_1}(1)\end{equation} holds uniformly for $0 \leq \eta_1 \leq \eta_2 \leq d$. Indeed, if not, then there would be $\epsilon > 0$ and sequences $\eta_1^n$ and $\eta_2^n$, in $[0,d]$ such that $\eta_1^n < \eta_2^n$ and $|\eta_2^n-\eta_1^n| \to 0$, but $|L(\eta^n_1,\eta^n_2)| \geq \epsilon$. Taking a subsequence, we may assume $\eta_1^n \to \eta_1$ and $\eta_2^n \to \eta_2$. Then $\eta_1 = \eta_2$. If $0 < \eta_1$, this is a contradiction, since for $n$ large
\[\epsilon \leq |L(\eta_1^n,\eta_2^n)| \lesssim_{\eta_1} (\eta_2^n)^N-(\eta_1^n)^N \to 0,\] where the constant depends only on $\eta_1$.

If $\eta_1 = 0$, then $\eta_2^n \to 0$, and
\[\epsilon \leq |L(\eta^n_1,\eta^n_2)| \in o_{\eta_2}(1),\] a contradiction.

Using \eqref{eq:GoursatSing:eq11} and \eqref{eq:GoursatSing:eq55} we deduce that $\tilde{X}$ is uniformly continuous in $\eta$, uniformly in $v$.

We have left to show that $X$ is uniformly continuous in $v$, uniformly in $\eta$.
From \eqref{eq:GoursatSing:eq4} and the bounds on $\tilde{X}$ and $\tilde{Y}$, it is clear that $\tilde{Y}$ is Lipschitz in $v$, uniformly in $\eta$. For $0 \leq v_1,v_2 \leq a_0$, set $\tilde{X}_{\Delta}(\eta) = \tilde{X}(v_1,\eta)-\tilde{X}(v_2,\eta)$, and similarly for $\tilde{F}_{\Delta}$, $p^1_{\Delta}$, $q^1_{\Delta}$ and $\tilde{Y}_{\Delta}$. From \eqref{eq:GoursatSing:eq4},
\begin{align}
\label{eq:GoursatSing:eq10}
\begin{split}
    \tilde{X}_{\Delta}(\eta) = \tilde{F}_{\Delta}(\eta) - \eta^{-N}\int_0^\eta &t^{N-1}p^1(v_1,t)\tilde{X}_{\Delta}(t) + t^{N-1}p^1_{\Delta}\tilde{X}(v_2,t)\\
    &t^{N-1}q^1(v_1,t)\tilde{Y}_{\Delta}(t) + t^{N-1}q^1_{\Delta}\tilde{Y}(v_2,t)\ dt.
\end{split}
\end{align}
Set
\[C_{\Delta} = \norm{p^1_{\Delta}\tilde{X}(v_2,\cdot)+q^1(v_1,\cdot)\tilde{Y}_{\Delta}(\cdot)+q^1_{\Delta}\tilde{Y}(v_2,\cdot)}_{C^0([0,d])}.\] Since $\tilde{X}(v_2,\cdot)$ and $\tilde{Y}(v_2,\cdot)$ are uniformly bounded by a constant independent of $a_0$, $C_{\Delta} \in o_{|v_1-v_2|}(1)$, where the $o(1)$ does not depend on $v_1,v_2$ or $a_0$. We conclude from \eqref{eq:GoursatSing:eq10} that
\[\norm{\tilde{X}_{\Delta}}_{C^0([0,d])} \leq \norm{\tilde{F}_{\Delta}}_{C^0(\mathcal R)} + \frac{1}{N}\norm{p^1}_{C^0(\mathcal R)}\norm{\tilde{X}_{\Delta}}_{C^0([0,d])} + \frac{1}{N}C_{\Delta}.\]
However, $N$ was already chosen large enough so that $\norm{p^1_{\Delta}}_{C^0(\mathcal R)}/N < 1$, so we may rearrange to obtain that
\[|\tilde{X}_{\Delta}| \in o_{|v_1-v_2|}(1),\] where the $o(1)$ does not depend on $v_1,v_2$ or $a_0$.

With $\tilde{X}$ uniformly continuous, \eqref{eq:GoursatSing:eq4} shows that $\tilde{Y}$ is a (weak) solution to an ODE in $v$ depending on $\eta$ in a uniformly continuous fashion. It follows that $\tilde{Y}$ is uniformly continuous.

Thus $\tilde{X}$ and $\tilde{Y}$ are uniformly continuous, with modulus of continuity not depending on $a_0$. This allows us to iterate well-posedness.

\textbf{Uniqueness.} Global uniqueness follows from the same argument with $I,J$, except applied to the difference of two solutions.

\textbf{Smoothness.} Putting everything together, we have established the existence of a unique solution $(X,Y)$ in $\eta^\infty C^0(\mathcal R)$.

Thus, we may suppose that for $M$ large enough, $\tilde{X},\ \tilde{Y} \in \eta^M C^0(\mathcal R)$. Differentiating \eqref{eq:GoursatSing:eq4} shows that $\tilde{X} \in C_0^v C_1^\eta (\mathcal R)$ and $\tilde{Y} \in C^0_\eta C^1_v (\mathcal R)$.

$\tilde{Y}$ now satisfies an ODE with continuous parameter $\eta$, and thus $\tilde{Y} \in C^1 (\mathcal R)$. Thus $\tilde{X}$ satisfies a singular ODE with $C^1$ dependence on a parameter $v$, and so $\tilde{X} \in C^1(\mathcal R)$ (the right-hand side is no longer $C^{\infty}$ in $\eta$ --it is only $C^1$--but this does not affect the proof). Doing this for all $N$ and shows that $X,Y \in \eta^\infty C^1(\mathcal R)$. Iterating, we deduce that $\tilde{X},\tilde{Y} \in \eta^{\infty} C^{\infty}(\mathcal R)$.

\end{proof}

\section{Riemann's method and the kernel of the solution operator}
\label{C:A3:kernel}
Riemann's method is the classical way to find the forward fundamental solution $E(v,\eta;s,t)$ to \eqref{eq:A3:GoursatII}, and hence the kernel of the solution operator, in the case that $P$ is linear. See \cite{EncRiemm}, for instance.

We will need an extension of this to \eqref{eq:A3:GoursatI}. For this section, we work with $\mathcal R = [a,b]\times [c,d]$ and the operator $L = L(X,Y) = (L^1(X,Y),L^2(X,Y))$ defined by
\begin{gather}
\label{eq:A3:GoursatIII}
\begin{split}
L^1(X,Y) = \pa_\eta X + p^1(v,\eta)X(v,\eta)+q^1(v,\eta)Y(v,\eta)\\
L^2(X,Y) = \pa_v Y + p^2(v,\eta)X(v,\eta)+ q^2(v,\eta)Y(v,\eta),
\end{split}
\end{gather}
where $p^i,q^i \in C^{\infty}(\mathcal R)$, ($i = 1,2$).\footnote{One could get away with less regularity, say $C^0(\mathcal R)$, but we will not need this.} For simplicity, we will only treat the case that $X$ and $Y$ are scalar-valued, although results and proofs carry over to the setting that $X$ and $Y$ are vector-valued.\footnote{Albeit having to change all matrices to block matrices, each block entry taking the place of a scalar.}

We first find the forward fundamental solution of $L$, i.e.\ the kernel \[E(v,\eta;s,t) = \begin{pmatrix} E_{11}(v,\eta;s,t) & E_{12}(v,\eta;s,t)\\
E_{21}(v,\eta;s,t) & E_{22}(v,\eta;s,t)\end{pmatrix}\] for which
\[LE = \mathrm{diag}(\delta(v-s)\delta(\eta-t), \delta(v-s)\delta(\eta-t))\] (we treat $L$ acting on matrices as acting on the columns as vectors), and $\supp E \subseteq \{v \geq s, \eta \geq t\}$.

Let us motivate how to find the solution. We will look at the first column of $E$, since the second is similar. Make the ansatz
\begin{align*}
E_{11}(v,\eta;s,t) &= E_{11}^1 \1_{\{v \geq s\}}\1_{\{\eta \geq t\}} + E_{11}^2 \delta(v-s)\1_{\{\eta \geq t\}}\\
E_{21}(v,\eta;s,t) &+ E_{21}^1 \1_{\{v \geq s\}}\1_{\{\eta \geq t\}} + E_{21}^2 \1_{\{v \geq s\}}\delta(\eta-t),
\end{align*}
where $E_{ij}^k$ are all smooth functions. Then
\begin{align*}
L^1E_{11} &= (\pa_\eta E_{11}^1 + p^1E_{11}^1 + q^1E_{21}^1)\1_{\{v \geq s\}}\1_{\{\eta \geq t\}}\\
&+ (\pa_\eta E_{11}^2 + p^1 E_{11}^2)\delta(v-s)\1_{\{\eta \geq t\}}
+ (q^1E_{21}^2 + E_{11}^1)\1_{\{v \geq s\}}\delta(\eta-t)\\
&+ E_{11}^2\delta(v-s)\delta(\eta-t)\\
L^2E_{21} &= (\pa_vE^1_{21} + p^2E_{11}^1 + q^2E_{21}^1)\1_{\{v \geq s\}}\1_{\{\eta \geq t\}}\\
&+ (\pa_v E_{21}^2 + q^2 E_{21}^2)\1_{\{v \geq s\}}\delta(\eta-t)
+(p^2 E_{11}^2 + E_{21}^1)\delta(v-s)\1_{\{\eta \geq t\}}\\
&+ E_{21}^2\delta(v-s)\delta(\eta-t)
\end{align*}

In order for $L(E_{11},E_{21}) = (\delta,0)$, it suffices for:
\begin{romanumerate}
\item $L(E_{11}^1,E_{21}^1) = 0$;
\item $E_{11}^1(v,t;s,t) = 0$;
\item $E_{21}^1(s,\eta;s,t) = -p^2(s,\eta)E_{11}^2(s,\eta;s,t)$;
\item $E_{21}^2 \equiv 0$;
\item $\pa_\eta E_{11}^2(s,\eta;s,t) + p^1(s,\eta)E_{11}^2(s,\eta;s,t) = 0$;
\item $E_{11}^2(s,t;s,t) = 1$.\end{romanumerate}
The first three conditions specify that the pair $(E_{11}^1,E_{21}^1)$ solve the Goursat problem with certain initial data. The last two specify how to find the initial data for $E_{21}^1$ by solving an ODE along $\{v=s\}$ with data given at $\{v = s,\eta = t\}$.

This motivates:
\begin{thm}[Forwards fundamental solution]\label{thm:A3:ffs}Let $I(\eta;s,t)$, $J(v;s,t)$ denote the functions satisfying
\begin{romanumerate}
\item $\pa_\eta I(\eta;s,t) +p^1(s,\eta)I(\eta;s,t) = 0$;
\item $\pa_v J(v;s,t) + q^2(v,t)J(v;s,t) = 0$
\end{romanumerate}
with initial data $I(t;s,t) = 1$, $J(s;s,t) = 1$.\footnote{Of course if $P^1$ and $Q^2$ are any $\pa_\eta$- or $\pa_v$-primitives for $p^1$ and $q^2$, respectively, then $I(\eta;s,t) = \exp(P^1(s,t)-P^1(s,\eta))$ and $J(v;s,t) = \exp(Q^2(s,t)-Q^2(v,t))$.} Suppose \[E(v,\eta;s,t) = \begin{pmatrix} E_{11}(v,\eta;s,t) & E_{12}(v,\eta;s,t)\\
E_{21}(v,\eta;s,t) & E_{22}(v,\eta;s,t)\end{pmatrix}\] solves $LE = 0$ on $\{v \geq s, \eta \geq t\}$ (and is extended to be $0$ outside this region) with initial data
\begin{romanumerate}
\item $E_{11}(v,t;s,t) = 0$;
\item $E_{21}(s,\eta;s,t) = -p^2(s,\eta)I(\eta;s,t)$;
\item $E_{22}(s,\eta;s,t) = 0$;
\item $E_{12}(v,t;s,t) = -q^1(v,t)J(v;s,t)$.\end{romanumerate}

Let $(X,Y)$ be the unique solution to the Goursat problem
\[L\begin{pmatrix}X \\ Y\end{pmatrix} = \begin{pmatrix}F\\G\end{pmatrix}\] in $\mathcal R$ with initial data $X|_{\{\eta = c\}} = 0$, $Y|_{\{v = a\}} = 0$, and $F,G \in C^{\infty}(\mathcal R)$.
Then 
\begin{equation}\label{eq:ffs:eq1}\begin{pmatrix}X \\ Y\end{pmatrix} = \begin{pmatrix} \int_{c}^{\eta} I(\eta;v,t)F(v,t)\ dt\\ \int_a^v J(v;s,\eta)G(s,\eta)\ ds\end{pmatrix} + \int_a^v\int_c^\eta E(v,\eta;s,t)\begin{pmatrix}F\\ G\end{pmatrix}(s,t)\ dsdt.\end{equation} 
\end{thm}
We abuse terminology and refer to the triple $(E,I,J)$ as the forward fundamental solution of $L$.

In order for this theorem to make sense, we actually need to be able to integrate against $E$. We will in fact show $E$ is smooth in a suitable way, even though this doesn't quite make sense since they are only really defined on different regions varying in $s,t$. What we mean is:
\begin{lem}\label{thm:A3:ffssmooth}Fix $(s_0,t_0) \in [a,b)\times [c,d)$. Then $E$ is jointly smooth in $v,\eta,s,t$ in the region (a mwc with piecewise linear boundary)
\[\mathcal S = \{(s,t) \in \mathcal R(s_0,t_0), (v,\eta) \in [s_0,b]\times[t_0,d]\}.\]

If $p^i$, $q^i$ ($i = 1,2$), depend on a parameter $\sigma \in \Sigma$, then $E$ is also smooth in the parameter $\sigma$.\end{lem}

\begin{proof}
Suppose without loss of generality that $a = c = 0$.
Extend $p^i,q^i$ $(i = 1,2)$ to $[0,\infty)\times [0,\infty)$, and extend $E$ by requiring it to solve the same equation, but now with the extended coefficients on all of $[0,\infty)\times [0,\infty)$. Set $\tilde{E}(v,\eta;s,t) = E(v-s,\eta-t;s,t)$. Then it is easy to check that $\tilde{E}$ satisfies a linear Goursat problem with coefficients smooth in $s,t$ (and if there is a parameter $\sigma$, then it is also smooth in $\sigma$) and initial data given on $\{v = 0\}$ and $\{\eta = 0\}$ which are smooth. Thus by \cref{thm:A3:globalsmooth}, $\tilde{E}$ is smooth on $[0,\infty)^2 \times [0,\infty)^2$, and thus $E$ is smooth on $\mathcal S$.
\end{proof}

\begin{proof}[Proof of \cref{thm:A3:ffs}]
We compute, where $E_{i,\bullet}$ denotes the $i$th row of $E$ ($i = 1,2$),
\begin{align*}
L^1X &= I(\eta;v,\eta)F(v,\eta) - \int_c^\eta p^1(v,\eta)I(\eta;v,t)F(v,t)\ dt+ \int_a^v E_{1,\bullet}(v,\eta;s,\eta)\begin{pmatrix}F\\ G\end{pmatrix}(s,t)\ ds\\
&+\int_a^v L^1E(v,\eta;s,t)\begin{pmatrix}F\\ G\end{pmatrix}(s,t)\ dsdt\\
&+p^1(v,\eta) \int_c^\eta I(\eta;v,t)F(v,t)\ dt + q^1(v,\eta)\int_a^v J(v;s,\eta)G(s,\eta)\ ds\\
&= F(v,\eta),
\end{align*}
and similarly
\[L^2Y = G.\]

The initial data are clearly $0$.
\end{proof}

Equation \eqref{eq:ffs:eq1} allows us to solve the problem $L(X,Y) = (F,G)$ with trivial data. We may extend it to solve the equation with data $X|_{\{\eta = c\}} = X_0$ and $Y|_{\{v = a\}} = Y_0$.
\begin{thm}\label{thm:A3:repform}Let $(E,I,J)$ be the forward fundamental solution to $L$. Let $(X,Y)$ be the unique solution to the Goursat problem
\[L\begin{pmatrix}X \\ Y\end{pmatrix} = \begin{pmatrix}F\\G\end{pmatrix}\] in $\mathcal R$ with initial data $X|_{\{\eta = c\}} = X_0 \in C^{\infty}(\{\eta = c\})$, $Y|_{\{v = a\}} = Y_0 \in C^{\infty}(\{v = a\})$ and $F,G \in C^{\infty}(\mathcal R)$.
Then \begin{align}\label{eq:repform:eq1}\begin{split}\begin{pmatrix}X \\ Y\end{pmatrix} &= \begin{pmatrix} \int_{c}^{\eta} I(\eta;v,t)F(v,t)\ dt\\ \int_a^v J(v;s,\eta)G(s,\eta)\ ds\end{pmatrix} + \int_a^v\int_c^\eta E(v,\eta;s,t)\begin{pmatrix}F\\ G\end{pmatrix}(s,t)\ dsdt\\
&+\begin{pmatrix}I(\eta;v,c)X_0(v) \\J(v;a,\eta)Y_0(\eta)\end{pmatrix} + \begin{pmatrix} \int_a^v E_{11}(v,\eta;s,c)X_0(s)\ ds\\ \int_c^\eta E_{22}(v,\eta;a,t)Y_0(t)\ dt\end{pmatrix}\\
&+\begin{pmatrix}\int_c^\eta (-E_{12}(v,\eta;v,t)+ E_{12}(v,\eta;a,t))Y_0(t)\ dt\\\int_a^v (-E_{21}(v,\eta;s,\eta) + E_{21}(v,\eta;s,c))X_0(s)\ ds\end{pmatrix}.\end{split}
\end{align}
 \end{thm}
\begin{proof}
Using the properties of the forward fundamental solution $(E,I,J)$, it is easy to check that \eqref{eq:repform:eq1} gives the right initial data for $X$ and $Y$. Thus we only need to check that \eqref{eq:repform:eq1} holds for $v > a$ and $\eta > c$.

It suffices to show \eqref{eq:repform:eq1} in the case $(F,G) = (0,0)$, because $L$ is linear and we know by \cref{thm:A3:ffs} that it is true for the case of trivial data and arbitrary $(F,G)$.
Extend $X_0$ to a function defined on all of $\mathcal R$ so that $\pa_\eta X_0 = 0$ and $Y_0$ to a function defined on all of $\mathcal R$ so that $\pa_v Y_0 = 0$. Assume we could solve the problem
\begin{align}
\label{eq:repform:eq2}
\begin{split}
L\begin{pmatrix}\tilde{X} \\ \tilde{Y}\end{pmatrix} &= -L\begin{pmatrix}X_0 \\ Y_0\end{pmatrix}\\
\tilde{X}|_{\{\eta=c\}} &= 0\\
\tilde{Y}|_{\{v = a\}} &= 0.
\end{split}
\end{align}
Then $X = \tilde{X} + X_0$ and $Y = \tilde{Y} + Y_0$.

Using \cref{thm:A3:ffs} we may solve for $\tilde{X}$ and $\tilde{Y}$ via
\begin{equation}\label{eq:repform:eq3}\begin{pmatrix}\tilde{X} \\ \tilde{Y}\end{pmatrix} = -\begin{pmatrix} \int_{c}^{\eta} I(\eta;v,t)L^1(X_0,Y_0)(v,t)\ dt\\ \int_a^v J(v;s,\eta)L^2(X_0,Y_0)(s,\eta)\ ds\end{pmatrix} - \int_a^v\int_c^\eta E(v,\eta;s,t)L\begin{pmatrix}X_0\\ Y_0\end{pmatrix}(s,t)\ dsdt.\end{equation}
Let us analyze the first term. Let $(L^1)^\ast$, $(L^2)^\ast$ denote the formal adjoint to $L^1$ and $L^2$, respectively. Explicitly, $(L^1)^\ast$ and $(L^2)^\ast$ are given, acting on a function $A$, and $B$, respectively, by
\begin{align*}
(L^1)^\ast A &= \left( -\pa_\eta A + A p^1,\ A q^1\right)\\
(L^2)^\ast B &= \left(Bp^2, \  -\pa_v A + Bq^2\right).\end{align*}
Then
\begin{align*}
    \begin{pmatrix} \int_{c}^{\eta} I(\eta;v,t)L^1(X_0,Y_0)(v,t)\ dt\\ \int_a^v J(v;s,\eta)L^2(X_0,Y_0)(s,\eta)\ ds\end{pmatrix} &= \begin{pmatrix}\int_c^\eta ((L^1)^\ast  I(\eta;v,t))(X_0,Y_0)(v,t)\ dt\\
    \int_a^v ((L^2)^\ast  J(v;s,\eta))(X_0,Y_0)(s,\eta)\ ds\end{pmatrix}\\
    &+
    \begin{pmatrix}(1-I(\eta;v,c))X_0(v)\\ (1-J(v;a,\eta))Y_0(\eta)\end{pmatrix}\\
\end{align*}

Now let us analyze the second term in \eqref{eq:repform:eq3}.

 Let $L^\ast$ denote the formal adjoint of $L$, which operates on $2\times 2$ matrices by operating on each row. Explicitly, $L^\ast$ is given, acting on a row vector $(A,B)$, by\
 \[L^\ast (A,B) = (L^1)^\ast \oplus (L^2)^\ast = \left(-\pa_\eta A + Ap^1 + B p^2, \ -\pa_v B +  B q^1 +  B q^2\right).\footnote{In particular, $(L^\ast)^i \neq (L^i)^\ast$, $i = 1,2$.}\]
 Then
\begin{align*}
    \int_a^v\int_c^\eta E(v,\eta;s,t)L\begin{pmatrix}X_0\\ Y_0\end{pmatrix}(s,t)\ dsdt &= \int_a^v\int_c^\eta L^\ast E(v,\eta;s,t)\begin{pmatrix}X_0\\ Y_0\end{pmatrix}(s,t)\ dsdt\\
    &+ \int_a^v (E_{\bullet,1}(v,\eta;s,\eta)-E_{\bullet,1}(v,\eta;s,c))X_0(s)\ ds\\
    &+\int_c^\eta(E_{\bullet,2}(v,\eta;v,t)-E_{\bullet,2}(v,\eta;a,t))Y_0(t)\ dt,
\end{align*}
where $E_{\bullet,i}$ denotes the $i$th column of $E$ ($i = 1,2$). Now $E_{11}(v,\eta;s,\eta) = 0$ and $E_{22}(v,\eta;v,t) = 0$, and so the second line is
\[\begin{pmatrix}-\int_a^v E_{11}(v,\eta;s,c)X_0(s)\ ds \\ \int_a^v (E_{21}(v,\eta;s,\eta) - E_{21}(v,\eta;s,c))X_0(s)\ ds\end{pmatrix}\]
and the third line is
\[\begin{pmatrix} \int_c^\eta (E_{12}(v,\eta;v,t)- E_{12}(v,\eta;a,t))Y_0(t)\ dt\\ -\int_c^\eta E_{22}(v,\eta;a,t)Y_0(t)\ dt\end{pmatrix}.\]
Putting everything together, we deduce that
\begin{align*}\begin{pmatrix}X \\ Y\end{pmatrix} &= \begin{pmatrix}I(\eta;v,c)X_0(v) \\J(v;a,\eta)Y_0(\eta)\end{pmatrix} + \begin{pmatrix} \int_a^v E_{11}(v,\eta;s,c)X_0(s)\ ds\\ \int_c^\eta E_{22}(v,\eta;a,t)Y_0(t)\ dt\end{pmatrix}\\
&+\begin{pmatrix}\int_c^\eta (-E_{12}(v,\eta;v,t)+ E_{12}(v,\eta;a,t))Y_0(t)\ dt\\\int_a^v (-E_{21}(v,\eta;s,\eta) + E_{21}(v,\eta;s,c))X_0(s)\ ds\end{pmatrix}\\
&- \begin{pmatrix}\int_c^\eta ((L^1)^\ast I(\eta;v,t))(X_0,Y_0)(v,t)\ dt\\
    \int_a^v ((L^2)^\ast J(v;s,\eta))(X_0,Y_0)(s,\eta)\ ds\end{pmatrix}\\
&-\int_a^v\int_c^\eta L^\ast E(v,\eta;s,t)\begin{pmatrix}X_0\\ Y_0\end{pmatrix}(s,t)\ dsdt.
\end{align*}
If we can show that the sum of the last two terms are $0$, then we will have completed the proof.
In fact, we claim that for any $(v_0,\eta_0) \in \mathcal R$ and $A,B \in L^1(\mathcal R(v_0,\eta_0))$,
\begin{equation}\label{eq:repform:eq4}K(v_0,\eta_0) := \begin{pmatrix}\int_c^{\eta_0} ((L^1)^\ast  I(\eta;v,t))(A,B)(v,t)\ dt\\
    \int_a^{v_0} ((L^2)^\ast  J(v;s,\eta))(A,B)(s,\eta)\ ds\end{pmatrix} + \int_a^{v_0}\int_c^{\eta_0} L^\ast E(v,\eta;s,t)\begin{pmatrix}A\\ B\end{pmatrix}(s,t)\ dsdt\end{equation}
vanishes. Approximating, it suffices to treat the case $(A,B) \in C_c^\infty(\mathcal R(v_0,\eta_0)^\circ)$.

Then, for $v \geq v_0$, $\eta \geq \eta_0$,
\begin{align*}
&\begin{pmatrix}\int_c^\eta ((L^1)^\ast  I(\eta;v,t))(A,B)(v,t)\ dt\\
\int_a^v ((L^2)^\ast  J(v;s,\eta))(A,B)(s,\eta)\ ds\end{pmatrix} + \int_a^{v}\int_c^{\eta} L^\ast E(v,\eta;s,t)\begin{pmatrix}A\\ B\end{pmatrix}(s,t)\ dsdt \\
&= \begin{pmatrix} \int_{c}^\eta I(\eta;v,t)L^1(A,B)(v,t)\ dt\\
\int_a^v J(v;s,\eta)L^2(A,B)(s,\eta)\ ds\end{pmatrix} + \int_a^{v}\int_c^\eta E(v,\eta;s,t)L\begin{pmatrix}A \\ B\end{pmatrix}(s,t)\ dsdt =: M(v,\eta)
\end{align*}
Hence, if $v \geq v_0$, $\eta \geq \eta_0$, then $K(v,\eta) = M(v,\eta)$.
By \cref{thm:A3:ffs}, $M$ solves the Goursat problem $LM = L(A,B)$ in all of $\mathcal R$ with trivial data. Since $(A,B) \in C_c^\infty(\mathcal R(v_0,\eta_0)^\circ)$, it also has trivial data, and thus by uniqueness, $M = (A,B)$. In particular, for $\eta \geq \eta_0, v \geq v_0$,
\[K(v,\eta) = M(v,\eta) = (A,B)(v,\eta) = 0.\]\end{proof}

\begin{singlespace}
\bibliographystyle{abbrv}
\addcontentsline{toc}{chapter}{Bibliography}
\bibliography{bib}

\begin{bibdiv}
\begin{biblist}

\bib{AnForm}{article}{
      author={{An}, Xinliang},
       title={{Formation of trapped surfaces from past null infinity}},
        date={2012},
     journal={ArXiv e-prints},
      eprint={1207.5271},
}

\bib{AnEmer}{article}{
      author={{An}, Xinliang},
       title={{Emergence of apparent horizon in gravitational collapse}},
        date={2017},
     journal={ArXiv e-prints},
      eprint={1703.00118},
}

\bib{AnScale}{article}{
      author={{An}, Xinliang},
       title={{A scale-critical trapped surface formation criterion: a new
  proof via signature for decay rates}},
        date={2019},
     journal={arXiv e-prints},
      eprint={1903.02967},
}

\bib{AnLukTrap}{article}{
      author={{An}, Xinliang},
      author={{Luk}, Jonathan},
       title={Trapped surfaces in vacuum arising dynamically from mild incoming
  radiation},
        date={2014},
     journal={Advances in Theoretical and Mathematical Physics},
      volume={21},
}

\bib{XinXueExam}{article}{
      author={An, Xinliang},
      author={Zhang, Xuefeng},
       title={Examples of naked singularity formation in higher-dimensional
  {E}instein-vacuum spacetimes},
        date={2018},
     journal={Annales Henri Poincar{\'e}},
      volume={19},
      number={2},
       pages={619\ndash 651},
}

\bib{BasVasWunAsym}{article}{
      author={Baskin, Dean},
      author={Vasy, Andr{\'a}s},
      author={Wunsch, Jared},
       title={Asymptotics of radiation fields in asymptotically {M}inkowski
  space},
        date={2015},
        ISSN={0002-9327},
     journal={American Journal of Mathematics},
      volume={137},
      number={5},
       pages={1293\ndash 1364},
}

\bib{ChoTheo}{article}{
      author={{Choquet-Bruhat}, Yvonne},
       title={Th{\'e}or{\`e}me d'existence pour certains syst{\`e}mes
  d'{\'e}quations aux d{\'e}riv{\'e}es partielles non lin{\'e}aires},
        date={1952},
     journal={Acta Math.},
      volume={88},
       pages={141\ndash 225},
         url={https://doi.org/10.1007/BF02392131},
}

\bib{ChoGene}{book}{
      author={{Choquet-Bruhat}, Yvonne},
       title={{General Relativity and the Einstein Equations}},
   publisher={Oxford University Press},
     address={Oxford},
        date={2009},
        ISBN={978-0-19-923072},
}

\bib{ChoGerGlob}{article}{
      author={Choquet-Bruhat, Yvonne},
      author={Geroch, Robert},
       title={Global aspects of the {C}auchy problem in general relativity},
        date={1969},
     journal={Comm. Math. Phys.},
      volume={14},
      number={4},
       pages={329\ndash 335},
         url={https://projecteuclid.org:443/euclid.cmp/1103841822},
}

\bib{ChrNonl}{article}{
      author={Christodoulou, Demetrios},
       title={{Nonlinear nature of gravitation and gravitational wave
  experiments}},
        date={1991},
     journal={Phys. Rev. Lett.},
      volume={67},
       pages={1486\ndash 1489},
}

\bib{ChrForm91}{article}{
      author={Christodoulou, Demetrios},
       title={{The formation of black holes and singularities in spherically
  symmetric gravitational collapse}},
        date={1991},
     journal={Commun. Pure Appl. Math.},
      volume={44},
      number={3},
       pages={339\ndash 373},
}

\bib{ChrBoun}{article}{
      author={Christodoulou, Demetrios},
       title={Bounded variation solutions of the spherically symmetric
  {E}instein-scalar field equations},
        date={1993},
     journal={Comm. Pure Appl. Math.},
      volume={46},
       pages={1131\ndash 1220},
}

\bib{ChrForm}{book}{
      author={Christodoulou, Demetrios},
       title={{The Formation of Black Holes in General Relativity}},
   publisher={European Mathematical Society},
     address={ETH-Zentrum},
        date={2009},
        ISBN={978-3-03179-068-5},
}

\bib{AbbObse}{article}{
      author={Collaboration, LIGO~Scientific},
      author={Collaboration, Virgo},
       title={Observation of gravitational waves from a binary black hole
  merger},
        date={2016},
     journal={Phys. Rev. Lett.},
      volume={116},
         url={https://link.aps.org/doi/10.1103/PhysRevLett.116.061102},
}

\bib{DafForm}{article}{
      author={Dafermos, Mihalis},
       title={The formation of black holes in general relativity},
        date={2013},
     journal={Ast{\'e}risque},
      volume={352},
}

\bib{DeTExis}{article}{
      author={DeTurck, Dennis},
       title={Existence of metrics with prescribed {R}icci curvature: Local
  theory.},
        date={1981/82},
     journal={Inventiones mathematicae},
      volume={65},
       pages={179\ndash 208},
         url={http://eudml.org/doc/142846},
}

\bib{GriBasi}{incollection}{
      author={Grieser, Daniel},
       title={Basics of the b-calculus},
        date={2001},
   booktitle={Approaches to singular analysis},
   publisher={Springer},
       pages={30\ndash 84},
}

\bib{HawEllLarg}{book}{
      author={Hawking, Stephen},
      author={Ellis, George F.~R.},
       title={The large scale structure of space-time},
      series={Cambridge Monographs on Mathematical Physics},
   publisher={Cambridge University Press},
        date={1973},
        ISBN={9780521099066},
}

\bib{HinNonL}{article}{
      author={Hintz, Peter},
       title={Non-linear stability of the {K}err--{N}ewman--de {S}itter family
  of charged black holes},
        date={2018Apr},
        ISSN={2199-2576},
     journal={Annals of {PDE}},
      volume={4},
      number={1},
       pages={11},
}

\bib{HinVasGlobN}{article}{
      author={{Hintz}, Peter},
      author={Vasy, Andr{\'a}s},
       title={The global non-linear stability of the {K}err-de {S}itter family
  of black holes},
        date={201606},
     journal={Acta Mathematica},
      volume={220},
}

\bib{HinVasGlobA}{article}{
      author={{Hintz}, Peter},
      author={Vasy, Andr{\'a}s},
       title={{A global analysis proof of the stability of Minkowski space and
  the polyhomogeneity of the metric}},
        date={2017},
     journal={ArXiv e-prints},
      eprint={1711.00195},
}

\bib{KlaLukRodFull}{article}{
      author={Klainerman, Sergiu},
      author={Luk, Jonathan},
      author={Rodnianski, Igor},
       title={A fully anisotropic mechanism for formation of trapped surfaces
  in vacuum},
        date={2014},
        ISSN={1432-1297},
     journal={Inventiones mathematicae},
      volume={198},
      number={1},
       pages={1\ndash 26},
         url={https://doi.org/10.1007/s00222-013-0496-6},
}

\bib{KlaNicEvel}{book}{
      author={{Klainerman}, Sergiu},
      author={{Nicol\`o}, Francensco},
       title={{The Evolution Equation in General Relativity}},
   publisher={Birkh\"auser},
     address={Boston},
        date={2003},
        ISBN={978-1-4612-7408-7},
}

\bib{KlaRodEmer}{article}{
      author={Klainerman, Sergiu},
      author={Rodnianski, Igor},
       title={On emerging scarred surfaces for the {E}instein vacuum
  equations},
        date={2010},
        ISSN={1078-0947},
     journal={Discrete \& Continuous Dynamical Systems - {A}},
      volume={28},
       pages={1007},
  url={http://aimsciences.org//article/id/7977f64e-031a-4fc2-aa9a-d8c5afb01da9},
}

\bib{KlaRodForm}{article}{
      author={Klainerman, Sergiu},
      author={Rodnianski, Igor},
       title={On the formation of trapped surfaces},
        date={2012},
     journal={Acta Math.},
      volume={208},
      number={2},
       pages={211\ndash 333},
         url={https://doi.org/10.1007/s11511-012-0077-3},
}

\bib{LiYuCons}{article}{
      author={Li, Junbin},
      author={Yu, Pin},
       title={Construction of {C}auchy data of vacuum {E}instein field
  equations evolving to black holes},
        date={2012},
     journal={Annals of Mathematics},
      volume={181},
}

\bib{LukLoca}{article}{
      author={{Luk}, Jonathan},
       title={{On the local existence for the characteristic initial value
  problem in general relativity}},
        date={2011},
     journal={ArXiv e-prints},
      eprint={1107.0898},
}

\bib{MazElli}{article}{
      author={{Mazzeo}, Rafe},
       title={Elliptic theory of differential edge operators {I}},
        date={1991},
     journal={Comm. Partial Differential Equations},
      volume={16},
      number={10},
       pages={1615\ndash 1664},
}

\bib{MazMelMero}{article}{
      author={Mazzeo, Rafe},
      author={Melrose, Richard~B.},
       title={Meromorphic extension of the resolvent on complete spaces with
  asymptotically constant negative curvature},
        date={1987},
        ISSN={0022-1236},
     journal={Journal of Functional Analysis},
      volume={75},
      number={2},
       pages={260 \ndash  310},
  url={http://www.sciencedirect.com/science/article/pii/0022123687900978},
}

\bib{MazMelAdia}{article}{
      author={Mazzeo, Rafe},
      author={Melrose, Richard~B.},
       title={The adiabatic limit, {H}odge cohomology and {L}eray's spectral
  sequence for a fibration},
        date={1990},
     journal={J. Differential Geom.},
      volume={31},
      number={1},
       pages={185\ndash 213},
         url={https://doi.org/10.4310/jdg/1214444094},
}

\bib{MelCalc}{article}{
      author={Melrose, Richard~B.},
       title={{Calculus of conormal distributions on manifolds with corners}},
        date={1992},
        ISSN={1073-7928},
     journal={International Mathematics Research Notices},
      volume={1992},
      number={3},
       pages={51\ndash 61},
}

\bib{MelTapsit}{book}{
      author={Melrose, Richard~B.},
       title={The {A}tiyah-{P}atodi-{S}inger index theorem},
      series={Research notes in Mathematics},
   publisher={A K Peters, Ltd.},
     address={Wellesley, MA},
        date={1993},
      volume={4},
}

\bib{MelDiff}{book}{
      author={Melrose, Richard~B.},
       title={{Differential analysis on manifolds with corners}},
   publisher={Book, in preparation, available online},
        date={1996},
}

\bib{MelSpec}{article}{
      author={Melrose, Richard~B.},
       title={Spectral and scattering theory for the {L}aplacian on
  asymptotically {E}uclidian spaces},
        date={1997},
     journal={Spectral and Scattering Theory},
      volume={161},
}

\bib{MosRapi}{article}{
      author={{Moser}, J{\"u}rgen},
       title={A rapidly convergent iteration method and non-linear partial
  differential equations - {I}},
        date={1966},
     journal={Ann. Sc. Norm. Sup. Pisa},
      volume={20},
       pages={265\ndash 316},
}

\bib{EncRiemm}{misc}{
      author={Nakhushev, A.~M.},
       title={{R}iemann method},
        date={1999},
        note={Available at
  \url{https://www.encyclopediaofmath.org/index.php/Riemann_method}},
}

\bib{PenGrav}{article}{
      author={{Penrose}, Roger},
       title={Gravitational collapse and space-time singularities},
        date={1965},
     journal={Physical Review Letters},
      volume={14},
       pages={57\ndash 59},
}

\bib{RenRedu}{article}{
      author={{Rendall}, Alan~D.},
       title={Reduction of the characteristic initial value problem to the
  {C}auchy problem and its applications to the {E}instein equations},
        date={1990},
     journal={Proc. R. Soc. Lond. Ser. A Math. Phys. Eng. Sci.},
      volume={427},
      number={1872},
       pages={201\ndash 213},
}

\bib{SchUber}{article}{
      author={{Schwarzschild}, Karl},
       title={{{\"U}ber das Gravitationsfeld eines Massenpunktes nach der
  Einsteinschen Theorie}},
        date={1916},
     journal={Sitzungsberichte der K{\"o}niglich Preu{\ss}ischen Akademie der
  Wissenschaften (Berlin), 1916, Seite 189-196},
}

\bib{YuDyn}{article}{
      author={{Yu}, Pin},
       title={{Dynamical formation of black holes due to the condensation of
  matter field}},
        date={2011},
      eprint={1105.5898},
}

\bib{ZhuElev}{article}{
      author={Zhu, Xuwen},
       title={{The eleven dimensional supergravity equations on edge
  manifolds}},
        date={2018},
     journal={Annales Henri Poincar{\'e}},
      volume={19},
      number={8},
       pages={2347\ndash 2400},
}

\end{biblist}
\end{bibdiv}
\end{singlespace}
\end{document}